\documentclass[11pt]{article}

\usepackage{a4wide}
\usepackage{hyperref}
\hypersetup{colorlinks=true,citecolor=blue,linkcolor=blue,urlcolor=blue}
\usepackage{amsmath}
\usepackage{amssymb}
\usepackage{amsthm}   
\usepackage{enumerate}
\usepackage{thm-restate}
\usepackage[margin=1in]{geometry}
\usepackage{times} 
\usepackage{tikz}

 \usetikzlibrary{arrows.meta,positioning}

\tikzset{
  state/.style={circle, draw, thick, minimum size=16pt, inner sep=1pt},
  specialstate/.style={circle, draw, thick, minimum size=18pt, inner sep=2pt, double},
  >={Stealth[length=2.5mm,width=2mm]}
}

\theoremstyle{plain}
\newtheorem{theorem}{Theorem}
\numberwithin{theorem}{section}
\newtheorem{conjecture}[theorem]{Conjecture}
\newtheorem*{observation*}{Observation}
\newtheorem*{theorem*}{Theorem}
\newtheorem{corollary}[theorem]{Corollary}
\newtheorem{lemma}[theorem]{Lemma}
\newtheorem*{lemma*}{Lemma}
\newtheorem{observation}[theorem]{Observation}
\newtheorem{definition}[theorem]{Definition}
\newtheorem{remark}[theorem]{Remark}
\newtheorem*{remark*}{Remark}

\newtheorem{proposition}[theorem]{Proposition}

\usepackage[UKenglish]{babel}
\usepackage[UKenglish]{isodate}
\usepackage{url}
 
\usepackage{graphicx}

\usepackage[nonumberlist]{glossaries} 
\usepackage{longtable} 
\usepackage{etoolbox}

\newglossarystyle{symbtable}{%
  \setglossarystyle{long}%
  \renewcommand*{\glsgroupheading}[1]{}%
}
\setglossarystyle{symbtable}
\makeglossaries

\newcommand{\newsym}[4][]{%
  \ifglsentryexists{#2}{}{%
    \newglossaryentry{#2}{%
      name={#4},%
      description={#3},%
      sort={#2}%
    }%
  }%
  \glsadd{#2}%
  \ifcsdef{glsfirstpage@#2}{}{%
    \csxdef{glsfirstpage@#2}{\thepage}%
  }%
  \ensuremath{\glsentryname{#2}}%
}

\let\epsilon=\varepsilon
\let\eps=\varepsilon 
\let \phi=\varphi
\renewcommand{\Pr}{\pr}  
\newcommand{\E}{\mathbb{E}}

\def\integers{\mathbb{Z}}
\def\reals{\mathbb{R}} 

\newcommand{\K}{K}

\newcommand{\rr}{r}

\newcommand{\TEnd}{T}

\def\bottleneck{\mathrm{b}}
\def\cbottle{c_{\bottleneck}}
\def\abottle{\alpha_{\bottleneck}}

\def\p{\overline{p}} 
\def\vec{\overline}
\def\b{\mathbf{b}}
\def\n{\mathbf{n}}
\def\B{\mathbf{B}}
\def\s{\mathbf{s}}
\def\e{\mathbf{e}}
\renewcommand{\S}{\mathbf{S}}

\def\calB{\mathcal{B}}
\def\calC{\mathcal{C}}
\def\calD{\mathcal{D}}
\def\calE{\mathcal{E}}
\def\calF{\mathcal{F}}
\def\calI{\mathcal{I}}
\def\calL{\mathcal{L}}
\def\calM{\mathcal{M}}
\def\calN{\mathcal{N}}
\def\calP{\mathcal{P}}
\def\calR{\mathcal{R}}
\def\calS{\mathcal{S}}
\def\calT{\mathcal{T}}

\def\arr{\mathbf{arr}}
\def\romanarr{\mathrm{arr}}

\def\ball{\calB}
\newcommand{\bbar}{{\vec{b}}}
\newcommand{\bbbar}{{\vec{\b}}}

\newcommand{\type}{\textup{type}}
\newcommand{\Failure}{\textup{failure}}
\newcommand{\advancing}{\textup{advancing}}
\newcommand{\filling}{\textup{filling}}
\newcommand{\refilling}{\textup{refilling}}
\newcommand{\stabilising}{\textup{stabilising}}
\newcommand{\initialising}{\textup{initialising}} 

\newcommand{\dist}{\mathrm{dist}}

\def\Po{\mathsf{Po}}

\renewcommand{\Pr}{\mathbb{P}}
\newcommand{\pr}{\mathbb{P}}

\def\Einit{\calE_{\mathrm{init}}}
\def\Edom{\calE_{\mathrm{dom}}}

\newcommand{\CSE}{C_\mathrm{SE}}
\newcommand{\Cnoise}{C_\mathrm{No}}
\newcommand{\CNoise}{\Cnoise}  
\def\L{\mathrm{L}}
\def\LL{\mathrm{LL}}

\newcommand{\CUpsilon}{C_{\Upsilon}}

\def\Wtilde{\widetilde{W}} 
\newcommand{\Cfill}{C_{\mathrm{F}}}  
\newcommand{\CFill}{\Cfill}
\newcommand{\CNoiseBelow}{c_{\mathrm{NB}}}

\def\traje{\mathrm{tr}}
\def\vectraje{\overline{\traje}}

\def\Avect{\overline{A}}
\def\avect{\overline{a}}
\def\abar{\avect}
\def\bvect{\overline{b}}
\def\evect{\overline{\e}}
\def\Gammavect{\overline{\Gamma}}
\def\hvect{\overline{h}}
\def\Lambdavect{\overline{\Lambda}}
\def\mvect{\overline{m}}
\def\muvect{\overline{\mu}}
\def\onevect{\overline{1}}
\def\psivect{\overline{\psi}}
\def\Psivect{\overline{\Psi}}
\def\Pvect{\overline{\Psi}}

\def\xbar{\overline{x}}
\def\xvect{\overline{x}}
\def\yvect{\overline{y}}

\def\zvect{\overline{z}}
\def\zerovect{\overline{0}}
\def\Rvect{\overline{R}}
 
\def\xbarWbar{\avect}
\def\Xvect{\overline{X}}
\def\Avect{\overline{A}}
\def\Pvect{\overline{P}}
\def\arrvect{\overline{\arr}}
\def\romanarrvect{\overline{\romanarr}}

\def\f{F}  
\def\F{F}

\def\In{I}
\def\Out{O}

\def\LoopStarts{\mathrm{LoopStarts}}
\def\Fills{\mathrm{Fills}}

\newcommand{\hls}{\mathrm{hls}}

\def\bup{b^{+}}
\def\bupsome{b^{+}}
\def\bbarup{\overline{\bup}}
\def\bbarupsome{\overline{\bupsome}}
\def\bupset{\mathsf{over}(\bbar)}
\def\bupsomeset{\mathsf{over}_{=}({\bbar})}

 \def\esc{\mathrm{Esc}}
\def\Esc{\esc}
\def\ebar{\overline{e}}
\def\bare{\overline{e}}
\def\bara{\overline{a}}
\def\barr{\overline{r}}
\def\hatS{\hat{S}}
\def\quiet{q}
\def\gammavect{\overline{\gamma}}

\newcommand{\poly}{\mathrm{poly}}

\newcommand\OurHLS{\ensuremath{\Psi^{\mathrm{BB}}(\Psi_0)}}

\def\trans{\mathrm{trans}}
\def\exit{\mathrm{exit}}
\def\back{\mathrm{back}}
\def\backexit{\mathrm{back}\textnormal{-}\mathrm{exit}}

\def\full{\mathrm{full}}
\def\tmin{t_{\mathrm{min}}}
\def\amin{a_{\mathrm{min}}}

\def\Hastad{H\r{a}stad}
\def\fail{\mathrm{fail}}
\def\succ{\mathrm{succ}}

\newcommand{\high}{\mathrm{hi}}
\newcommand{\low}{\mathrm{lo}}
\newcommand{\excess}{\mathrm{excess}}

\def\Ready{\mathrm{Ready}} 

\newcommand{\kMean}{\kappa}

\newcommand{\badRone}{1}
\newcommand{\badRtwo}{2}
\newcommand{\badRfour}{3}
\newcommand{\badLoopStarts}{4}
\newcommand{\badLongLoop}{5}

\title{The Instability of all Backoff Protocols\thanks{For the purpose of Open Access, the authors have applied a CC BY public copyright licence to any Author Accepted Manuscript (AAM) version arising from this submission. 
}
}
\author{Leslie Ann Goldberg and John Lapinskas
}
 
\usepackage{times}
\usepackage{multicol}

\begin{document}

\date{July 2026}
\maketitle

\begin{abstract}

In this paper we prove Aldous's conjecture from 1987 that there is no backoff protocol that is stable for any positive arrival rate. The setting is a communication channel for coordinating requests for a shared resource. Each user who wants to access the resource makes a request by sending a message to the channel. The users don't have any way to communicate with each other, except by sending messages to the channel. The operation of the channel proceeds in discrete time steps. If exactly one message is sent to the channel during a time step then this message succeeds (and leaves the system). If multiple messages are sent during a time step then these messages collide. Each of the users that sent these messages therefore waits a random amount of time before re-sending. A backoff protocol is a randomised algorithm for determining how long to wait --- the waiting time is a function of how many collisions a message has had. Specifically, a backoff protocol is described by a send sequence $\p = (p_0,p_1,p_2,\ldots)$. If a message has had $k$ collisions before a  time step then, with probability~$p_k$, it sends during that time step, whereas with probability $1-p_k$ it is silent (waiting for later).  
The most famous backoff protocol is binary exponential backoff, where $p_k = 2^{-k}$. Under Kelly's model, in which the number of new messages that arrive in the system at each time step is given by a Poisson random variable with mean~$\lambda$,
Aldous proved that binary exponential backoff is unstable for any positive~$\lambda$. 
He conjectured that the same is true for \emph{any} backoff protocol.
We prove this conjecture.
\end{abstract}

\newpage

\setcounter{tocdepth}{2} 
\tableofcontents
\newpage

\section{Introduction}\label{sec:intro}

A \emph{multiple-access channel} is a communication channel for coordinating requests for a shared resource. Each user who wants to access the resource makes a request by sending a message to the channel. Users don't have any way to communicate with each other, except by sending messages to the channel (so they cannot coordinate when to send).  
The operation of the channel proceeds in discrete time steps. If exactly one message is sent to the channel during a  time step, then this message succeeds (and leaves the system). If multiple messages are sent during any time step, then these messages \emph{collide} and must be re-transmitted later. The users that sent these messages 
observe that a  collision has occurred (but don't get any additional information). Each of these users waits some random amount of time before re-sending.  

A \emph{backoff protocol} is a natural randomised algorithm for determining how long to wait, encapsulating the idea that the waiting time of a message before re-sending is a function of how many times it has already collided. Thus, a backoff protocol is described by a \emph{send sequence} $\p = (p_0,p_1,p_2,\ldots)$.
If a message has already had $k$~collisions before a given time step then, with probability~$p_k$, it sends during that time step, whereas with probability~$1-p_k$ it is silent (waiting for later). The most famous backoff protocol is \emph{binary exponential backoff} where $p_k=2^{-k}$. Binary exponential backoff forms the basis of Ethernet~\cite{MB-BEB} and TCP/IP~\cite{TCP}. The name comes from the fact that, each time a message experiences a collision, it ``backs off'' by halving its send rate. 

The most natural model for studying backoff protocols is the model of Kelly~\cite{Kelly}. The number of new messages that arrive in the system at each time step is given by a Poisson random variable with some mean~$\lambda$,   and the central question is whether the resulting stochastic process is \emph{stable} as it evolves over time. Following Kelly, we say that a backoff protocol (given by a send sequence~$\p$) is stable for a \emph{birth rate} $\lambda$ if the resulting Markov process (which we call a ``backoff process'') is positive recurrent. 
Stability implies that the backoff process converges to a unique stationary distribution with a finite expected number of messages waiting in the system (this implies that 
backoff processes are not stable  for $\lambda>1$ since at most one message leaves the system during each step). On the other hand, instability means 
that the messages build up in a very unacceptable way -- instability implies
that the expected return time from \emph{any} state of the system back to that state is infinite and, as $t\rightarrow \infty$, the probability of being in any given state at time~$t$ tends to~$0$.

The earliest example of a backoff protocol is the slotted ALOHA protocol~\cite{ALOHA, Roberts-ALOHA}, which was motivated by radio communication at the University of 
Hawaii.  This is a backoff protocol with a send sequence $\p = (p,p,p,\ldots)$ for some $p\in (0,1)$.
It is easy to show that this protocol is unstable for any positive birth rate~$\lambda$ -- the send probabilities are simply too high (they don't back off at all).
In fact, Kelly and Macphee~\cite{kelly-macphee} give a characterisation of the
collection of protocols with the property that, with probability~$1$, the corresponding backoff process   only has finitely many successful sends (over the course of its infinite evolution). Slotted ALOHA fits into this category. 
Perhaps more surprisingly, \emph{polynomial backoff} with   $p_k = k^{-c}$ (for any constant~$c$) also fits into this category -- the rate of backoff is too low so, with probability~$1$, there are only finitely many successful sends. This is a very strong form of  instability -- 
there are unstable backoff processes which still have
a positive probability of having an infinite number of successful sends (see~\cite{Kelly}).

A very important result was obtained by Aldous in 1987~\cite{Aldous}.  
He showed that binary exponential backoff is unstable for every positive birth rate~$\lambda$. The result is important for two reasons. The first reason is practical -- Ethernet is built on binary exponential backoff, hence it is surprising that this is unstable. The second reason is mathematical, and has to do with the nature of his proof -- in essence, he showed that even binary exponential backoff doesn't back off enough. We will defer any discussion of the practical side to Section~\ref{sec:related-work} -- Ethernet networks don't run pure backoff, instead they are assisted by queues and rely on dropping packets. In this paper we are more interested in the mathematical question.
Consider the evolution of a multiple-access channel whose messages run binary exponential backoff. At any step~$t$, the state of the system can be described by a vector
$\bbar(t) = (b_0(t),b_1(t),\ldots)$ where $b_k(t)$ is the number of messages waiting in the system that have already collided $k$~times. 
We think of~$b_k$ as a ``bin'' that holds $b_k(t)$ messages at time~$t$.
Given the state $\bbar(t)$, the expected number of sends at the next step, from messages that are already in the system,  is $\sum_k p_k b_k(t)$. Informally, we refer to this as the ``expected noise'' at time~$t$. What Aldous shows is that, as $t$ increases, the expected noise grows without bound. In fact, it grows so much that, with positive probability, every step~$t$ has
$\Omega(\log t)$ expected noise.
Furthermore, the probability that this fails for a particular step~$t$ is so small that one can take a union bound, summing the failure probabilities over~$t$ and getting a total failure probability that is less than~$1$.
(We refer to this situation as ``summable failure probability''.)
Aldous's proof shows that, even though binary exponential backoff does back off quite substantially,  the expected number of sends still grows  with every step.
Aldous conjectured that the same holds for every backoff protocol (that is, for every send sequence~$\p$). He believed, therefore, that his argument could be modified to show the instability of every backoff protocol for every positive birth rate. We record his conjecture as follows.

\begin{conjecture}[Aldous's Conjecture]\label{conj:backoff} 
No backoff process is stable for \textit{any}  birth rate $\lambda>0$.
\end{conjecture}

The result of this paper, which we state as Theorem~\ref{thm:goal}, is to finally prove Aldous's conjecture. 

\begin{restatable}{theorem}{thmgoal}\label{thm:goal}
Consider 
any send sequence~$\p$ and any birth rate $\lambda  >0$.
Then the backoff process 
with send sequence~$\p$ and birth rate~$\lambda$ is not positive recurrent.
\end{restatable}

Whereas  Conjecture~\ref{conj:backoff} is correct,  the truth is much more subtle than the arguments that we have already described.
It turns out that there are send sequences that are unstable due to excessive sending but that, nevertheless, go completely silent (with near-zero expected noise) infinitely often (an explicit example can be constructed by interleaving
the send probabilities from slotted ALOHA and binary exponential backoff --   see Section~\ref{sec:context}). 
As we will see, the behaviour of some backoff protocols is very complicated. For example, there are some send sequences for which the expected number of sends at time~$t$ grows, but grows very slowly with~$t$ ---
at rate  $\Omega(\log\log\log t)$ rather than 
rate~$\Omega(\log t)$ as in the protocols covered by Aldous's proof. For
such sequences, the expected noise, conditioned on the current state, goes to near~$0$ infinitely often --- so much more subtle arguments are needed to demonstrate instability.  

This section is organised as follows. Section~\ref{sec:proofsketch} gives a very high-level sketch of our proof, picking out the main ideas. Section~\ref{sec:context} 
gives an explicit example of a backoff protocol that goes silent infinitely often and gives more detail about the academic context of our work. 
There is a lot of very interesting recent work 
\cite{Bender2025SICOMP,
BKKP-adversarial, 
CCDKPP,
CJZ-jamming, ShahShin,
SST,
XZ-SPAA25}  on contention resolution in other settings, for example where there are exactly $n$ messages, or  where messages have strictly more information than they do in  our  model. This will be discussed in 
Section~\ref{sec:related-work}.  
Section~\ref{sec:introfuture} briefly discusses potential future work and Section~\ref{sec:outline} gives an outline of the rest of the paper.

\subsection{Proof sketch}\label{sec:proofsketch}

We start with a slightly more detailed  definition of a ``backoff process''.    The  state of a backoff process after step~$t$ is a vector 
$\bbbar(t) = (\b_0(t),\b_1(t),\ldots)$ where $\b_k(t)$ is the set of messages that have collided $k$ times during steps $1,\ldots,t$. We
think of $\b_k$ as a ``bin'' that holds a certain set of messages $\b_k(t)$ at time~$t$, which we think of as ``balls''. Since the balls-and-bins viewpoint will be relevant for the whole paper, we simply refer to messages as ``balls'' from here onwards. 
Definition~\ref{def:backoff} introduces some further useful notation. We introduce the superscript~$X$ for the backoff process~$X$ because the paper will need to consider \emph{many} different processes, often coupled with each other, so the notation needs to indicate the relevant process.

\begin{definition}\label{def:backoff} 
Fix a send sequence
$\newsym{p}{send sequence, Def~\ref{def:backoff}}{\p = (p_0,p_1,p_2,\ldots)}$  and a birth rate 
$\newsym{lambda}{birth rate, Def~\ref{def:backoff}}{\lambda}>0$. Each~$p_k$   is  a real number in~$(0,1]$.   $W_k$   denotes~$1/p_k$.
The \emph{backoff process} $X$
with send sequence~$\p$ and birth rate~$\lambda$ is defined as follows.
For each $i\geq 0$ and $t\geq 1$, 
$\b^X_i(t)$~is the set of balls that  collide exactly $i$ times during steps $1,\ldots,t$. 
The process starts from the empty state, so  
for each~$i\geq 0$, 
$\b^X_i(0) = \emptyset$.
For each $t\geq 1$, $\n^X(t)$ is the set of balls that are born at time~$t$. Its size is given by $n^X(t)$,  which 
is a Poisson random variable with mean~$\lambda$.
After the births, but just
before the sends at step~$t$, the bin populations are defined as follows.
$\B^X_0(t) = \n^X(t) \cup \b^X_0(t-1)$. For every $i\geq 1$, $\B^X_i(t) = \b^X_i(t-1)$.  
For every $i\geq 0$, the set of balls that send from bin~$i$ at time~$t$,  {$\s^X_i(t)$}, is
formed by including each ball from $\B^X_i(t)$ independently with probability~$p_i$.  The
set of balls that escape from bin~$i$ at time~$t$,  $\e^X_i(t)$, is equal to $\s^X_i(t)$ if 
$|\s_i^X(t)| = |\cup_{k\geq 0} \s_k^X(t)| = 1$. Otherwise,
$\e^X_i(t) =  \emptyset$.
Finally, $\b_0^X(t) = \B_0^X(t) \setminus \s_0^X(t)$ and for all  $i\geq 1$,
$\b_i^X(t) = (\B_i^X(t) \setminus \s_i^X(t)) \cup (\s^X_{i-1}(t) \setminus \e^X_{i-1}(t))$.
For every variable $\b^X_i(t)$, $\B^X_i(t)$, $\s^X_i(t)$, and $\e^X_i(t)$, the roman (non-bold-face) version represents its size, so, for example, $b^X_i(t) = |\b^X_i(t)|$.
\end{definition}

The evolution of a backoff process is a Markov chain since the state $\bbbar^X(t)=(\b_0^X(t),\b_1^X(t),\ldots)$ is a (random) function of $\bbbar^X(t-1)=(\b_0^X(t-1),\b_1^X(t-1),\ldots)$. Following Kelly~\cite{Kelly}, we say that the backoff process~$X$ is \emph{stable} 
if the empty state is positive recurrent, meaning that, after the process starts in the empty state (with no messages waiting) it eventually returns to that state with probability~$1$, and the expected time until the empty state is revisited is finite. 
Note that, except in a pathological case which is easily dealt with (Proposition~\ref{prop:irreducible-aperiodic}),  the Markov chain corresponding to a backoff process is irreducible and 
aperiodic.
This justifies our previous claim that 
stability implies that the backoff process converges to a unique stationary distribution with a finite expected number of messages. On the other hand, instability  implies
that the expected return time from \emph{any} state to itself is infinite and, as $t\rightarrow \infty$, the probability of being in any given state at time~$t$ tends to~$0$.

Before sketching our proof it is useful to consider how full a bin~$i$ can be expected to get in a backoff process. For this, consider briefly an \emph{externally-jammed channel} $Y$  which is the same as the backoff process (Definition~\ref{def:backoff}), except that balls never escape, so 
for all $i$ and $t$,
$\e_i^Y(t) = \emptyset$. As Aldous noted~\cite{Aldous}, 
in the stationary distribution of the externally-jammed channel,  for each~$t$,
the variables in $\{b_i^Y(t)\colon i \ge 0\}$ are independent, and $b_i^Y(t)$ is a Poisson variable with mean $\lambda W_i$. 
This implies that, for a backoff process~$X$, $\E[b_i^X(t)] \leq \lambda W_i$. 
Given this upper bound, we refer to $W_i$ as the \emph{weight} of bin~$i$ and we say that bin~$i$ is \emph{full} if its population is near $\lambda W_i$.
In this sketch, we define a ``low-weight bin'' to be a bin~$j$ with $W_j < j^2$ 
and a ``high-weight bin'' to be a bin~$j$ with $W_j\geq j^2$. We define a ``very high-weight bin'' to be a bin~$j$ with $W_j \geq j^{1000}$. The particular constants are not important --- the labels are just to de-clutter the explanation --- the send sequences that we deal with  
will also have bins of exponential weight (see Corollary~\ref{KM-cor}).  

We have already raised an important issue that we will have to deal with, namely analysing send sequences that have infinitely many ``quiet periods''. 
Our   first step is to split the set of send sequences into two types --- we will use different proofs for the two types.
To do this, we consider the bins corresponding to a send sequence~$\p$ and we show how to partition the bins into three categories.
Given~$\p$, every bin~$j$ is  \emph{covered}, \emph{weakly exposed}, or \emph{strongly exposed}.
These concepts are defined in Definition~\ref{def:covered}.
Here  it is sufficient to describe the consequences of the definition, which indicates what they mean intuitively. To do this
we will use 
two quantities  (Definition~\ref{def:noise})  
$\L(j) = \Theta(\log j)$ and $\LL(j) = \Theta(\log\log j)$, both with appropriately-large constants.

If bin~$j$ is \emph{covered} there is a ``reservoir''
$\calR_j$ of at least~$L(j)$ bins in $[j-1]$. 
The key property of~$\calR_j$ is that, once these bins are (nearly) full, we expect them to stay full for $\Omega(W_j)$ steps, which is long enough to enable bin~$j$ to fill up. The technical definition of~$\calR_j$ is a little bit complicated -- if $W_j$ isn't too large then $\calR_j$ may be equal to $[j-1]$. In other circumstances $\calR_j$ contains $\L(j)$ high weight bins.\footnote{The details are in 
Definition~\ref{def:covered}. Note the text before the definition  that describes the properties of $\calR_j$ using the notation of the actual proof -- 
we use the notation ``$\calR_j$''
in this proof sketch, but it is not used in  the proof itself.} The high-level proof technique of  Aldous (as adapted by~\cite{GL-oldcontention} to accommodate low-weight bins, which lack concentration)
can be made to work whenever $\p$ has the property that  all but finitely many bins are covered.

If bin~$j$ is \emph{weakly exposed} then it is not covered
but there is still a reservoir $\calR_j$ of size at least
$\LL(j)$ of very high-weight bins in
$[j-1]$.   The weight of a weakly-exposed bin is exponentially large in~$j$. 
The behaviour of a backoff process  is complicated regarding its weakly exposed bins. 
Suppose that bin $j$ is weakly exposed, and that the process has reached a state where 
bins in $[j-1]$ are fairly full and bin~$j$ is filling. It is likely that at some time while bin $j$ is filling, most bins in $[j-1]$ go mostly empty. This will drop the expected noise to the point where newly-born balls are likely to escape before reaching bin $j$. 
However, the noise from the bins in $\calR_j$ continues for a long time, and these bins will
provide enough expected noise to  
(repeatedly) restore the populations of the bins in~$[j-1]$, each time enabling further progress towards filling bin~$j$. Essentially, the noise from $\calR_j$ ensures that new balls always arrive much faster than existing balls escape.

What we develop in this work is a new technique for tracking the evolution of the backoff process while this happens, in terms of an algorithm (which we call the high-level state transition algorithm) that serves as an interface between the mundane sends in the backoff process and the more complicated pattern of escapes. We will describe this idea in more detail shortly.

Finally, a \emph{strongly exposed} bin is not covered or weakly exposed. Its weight is exponentially large in~$j$ and there are only $O(\log \log j)$ very high-weight bins in $[j-1]$. In general, even if the backoff process gets to a point where all bins in $[j-1]$ are full, these bins will be likely to go completely empty while bin~$j$ is filling. The send sequence described in Section~\ref{sec:context} has infinitely-many strongly exposed bins.
 
An important aspect of our approach is that we provide one instability proof for the case where there are infinitely many strongly exposed bins (Section~\ref{sec:intro-inf-SE}) and a different (and much more difficult and interesting) proof for the case where there are finitely many strongly exposed bins (Section~\ref{sec:finiteSE}).
Note that  past work, including \cite{GL-oldcontention}, was restricted to the case where every (sufficiently large) bin is covered -- nothing was known about the other cases.
We restrict this proof sketch to the most difficult case of the conjecture, which is the case where bin weights are at most exponentially large.

\subsubsection{Finitely-many strongly exposed bins}\label{sec:finiteSE}

Sections~\ref{sec:finiteSE} -- \ref{sec:bottomline} 
consider a send sequence~$\p$ with finitely-many strongly-exposed bins, so that, for a positive integer~$j_0$, all strongly-exposed bins are in $[j_0]$. 
Our goal is to show,   for any positive~$\lambda$, that
the backoff process~$X$ with send sequence~$\p$ and birth rate~$\lambda$ is not positive recurrent. 
The first step is to identify a finite time $\tau_0$ (depending on~$\p$ and~$j_0$) and a
possible, but unlikely, event $\Einit$ 
such that, conditioned on $\Einit$, 
the distribution  of the 
populations of the
first~$j_0$ bins of the backoff process~$X$ at time~$\tau_0$ is dominated below
by a tuple of independent Poisson distributions.
These independent Poisson distributions 
will have the property that
each bin is fairly full, so that the population of bin~$i$ is a Poisson random variable with mean $\eps\lambda W_i$ for some appropriate $\eps>0$. We write this as
$\dist(\bvect^X_{[j_0]}(\tau_0)) \gtrsim \Po(\eps\lambda W_1,\ldots,\eps\lambda W_{j_0})$, where
$\dist(\bvect^X_{[j_0]}(\tau_0))$ is the distribution of
the tuple
$(b_1^X(\tau_0),\ldots, b^X_{j_0}(\tau_0))$, the symbol~$\gtrsim$ indicates stochastic domination, 
and, for a tuple $\zvect = (z_1,\ldots,z_j)$, 
$\Po(\zvect)$ indicates a $j$-tuple of independent Poisson distributions where the mean of the $i$'th distribution is~$z_i$. 

We will condition on~$\Einit$ for the rest of the proof.
From time~$\tau_0$, the argument continues as follows.
Very roughly,  we consider each $j$ in order, starting with $j=j_0+1$. Bin $j$ is covered or weakly exposed (since all strongly-exposed bins are in~$[j_0]$) so there is a reservoir $\calR_j\subseteq [j-1]$, as we have discussed. In both cases, we will show that
bin $j$ fills in time $W_j \poly(j)$, except with some failure probability that can be summed over~$j$. 
By ``fills'' we (very roughly) mean getting to 
a situation at some time~$t$ where 
$\dist(\bvect^X_{[j]}(t)) \gtrsim \Po(\epsilon \lambda W_1,\ldots,\epsilon \lambda W_{j})$.
The case where $j$ is covered is fairly straightforward. By the definition of covered, $|\calR_j|\geq L(j)$. We expect the bins in~$\calR_j$ to stay full for long enough to fill bin~$j$, and the probability that this fails is sufficiently small that we can sum it over~$j$, taking a union bound. 

The difficulty arises when $j$ is weakly exposed.
In this case, we start with the good situation where
$\dist(\bvect^X_{[j-1]}(t)) \gtrsim \Po(\epsilon \lambda W_1,\ldots,\epsilon \lambda W_{j-1})$. When this good situation fails (which it will, with high probability!) there will no longer be enough expected noise to continue filling bin~$j$. Balls will start escaping. However, the $\LL(j)$ very high-weight bins in $\calR_j$ will maintain $\Omega(\log \log j)$ expected noise for at least $\poly(j)$ steps. This noise will be enough to enable the process to re-establish the initial good situation.
The big problem (for the proof) is how to drop down to this lower noise level without losing the  Poisson domination on the population of the bins. We tackle this problem with the following ideas.

It helps to consider what could cause the noise level to fail before the process naturally fills bin~$j$ and establishes $\dist(\bvect^X_{[j]}(t)) \gtrsim \Po(\epsilon \lambda W_1,\ldots,\epsilon \lambda W_{j})$. Crudely, there are two things that could go wrong. 
The first possibility (``a send problem'') is that many balls have unusual send patterns through the process, leading to (for example) a long stretch of empty low-weight bins, despite the fact that
no more balls escape than the expected number. 
The second possibility (``an escape problem'') is that more balls escape than what is expected. 
The crucial insight is that, while the events involved in an escape problem are heavily dependent, they are unlikely enough to union bound away. Meanwhile, while send problems are likely enough to be a concern, the events that they depend on (individual ball trajectories) are truly independent. 

Instead of trying to directly analyse the backoff 
process~$X$, we couple it 
with a ``volume process'' $V$ that tracks send patterns
and with an ``escape process'' $E$ that tracks escapes.
The volume process and escape process will interact through a tightly-controlled and \emph{limited} interface, which we will design to make the escape process more independent.
The coupling will 
maintain the following important
coupling invariant.
\begin{equation*}\tag*{(Coupling invariant)}
\b_i^V(t) \setminus \b_i^E(t) \subseteq \b_i^X(t).
\end{equation*}
Once we've done this,  to show that the backoff process is not positive recurrent, 
it will suffice to 
establish the following proof property.
\begin{equation*}\tag*{(Proof property)}
\mbox{With positive probability, $\forall t \,\exists i \,\b_i^V(t) \setminus \b_i^E(t)$ is non-empty.}
\end{equation*}

In order to describe the volume process it helps to first describe a variant of the externally jammed channel called a $j$-jammed process (Definition~\ref{def:j-jammed}).
The $j$-jammed process~$Y$ with send sequence~$\p$ and birth rate~$\lambda$ is the same as the backoff process with these parameters, except that for all bins $i$ with $i<j$, $\e_i^Y(t)=\emptyset$, whereas $\e_j^Y(t) = \s_j^Y(t)$, so all balls stay in the system until they send from bin~$j$, and they escape at that point, whether or not they have a collision. The balls in a $j$-jammed process behave independently, so we can show (Lemma~\ref{lem:slowfill-j-jammed}) that 
$\dist(b_{[j]}^Y(t)) = \Po(f_1(t),\ldots,f_j(t))$ where
$f_{[j]}(t)$ follows a straightforward recurrence (Definition~\ref{def:f}) and,
for all $t = \Omega( \sum_{x\in j} W_x)$, 
$f_k(t) \geq 9\lambda W_k/10$.

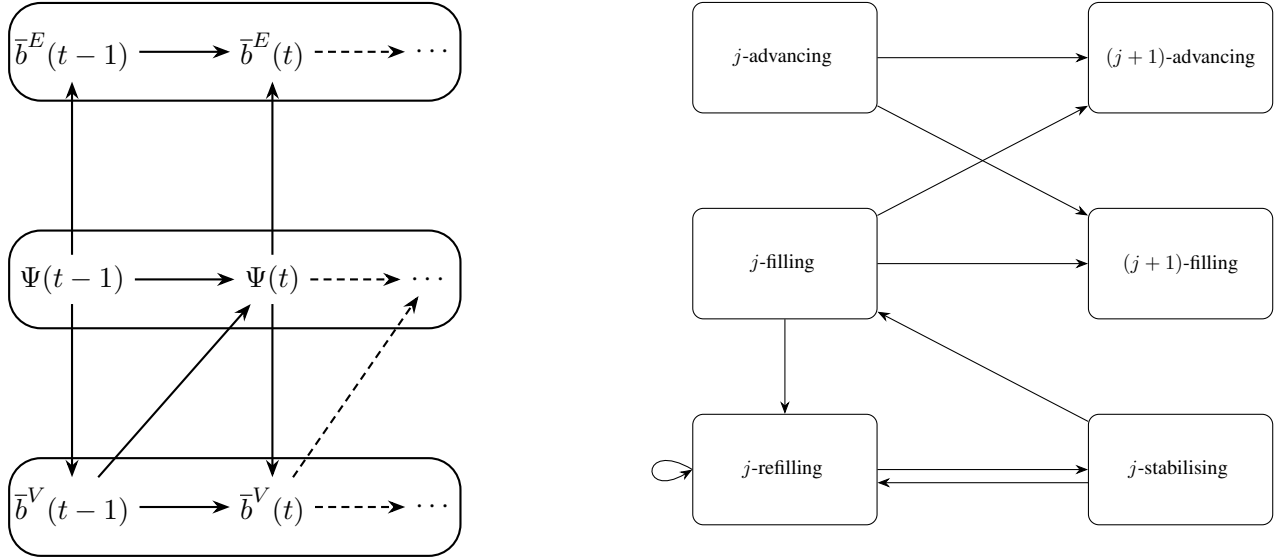
\begin{figure}[ht]
\begin{minipage}{0.45\textwidth}   
 \centering

\begin{tikzpicture}[
  >=Stealth,
  node distance=2.0cm,
  cloudbox/.style={
    rectangle,
    rounded corners=12pt,
    draw,
    thick,
    minimum width=6.0cm,   
    minimum height=1.3cm
  }
]

\node[cloudbox] (cloudV) {};

\node (bVt)  at ([xshift=0.85cm]cloudV.west) {$\overline{b}^{V}(t-1)$};
\node (bVtp) at ([xshift=0.5cm]cloudV.center) {$\overline{b}^{V}(t)$};
\draw[->, thick] (bVt) -- (bVtp);

\node[cloudbox, above=1.7cm of cloudV] (cloudPsi) {};
\node (Psit)  at ([xshift=0.85cm]cloudPsi.west) {$\Psi(t-1)$};
\node (Psitp) at ([xshift=0.5cm]cloudPsi.center) {$\Psi(t)$};
\draw[->, thick] (Psit) -- (Psitp);

\draw[->, thick] (bVt) -- (Psitp);

\node[cloudbox, above=1.7cm of cloudPsi] (cloudE) {};
\node (bEt)  at ([xshift=0.85cm]cloudE.west) {$\overline{b}^{E}(t-1)$};
\node (bEtp) at ([xshift=0.5cm]cloudE.center) {$\overline{b}^{E}(t)$};
\draw[->, thick] (bEt) -- (bEtp);

\draw[->, thick] (Psit) -- (bVt);
\draw[->, thick] (Psit) -- (bEt);

\draw[->, thick] (Psitp) -- (bVtp);
\draw[->, thick] (Psitp) -- (bEtp);

\node[right=1.2cm of bVtp]   (dotsV)   {$\cdots$};
\node[right=1.2cm of Psitp]  (dotsPsi) {$\cdots$};
\node[right=1.2cm of bEtp]   (dotsE)   {$\cdots$};

\draw[->, thick, densely dashed] (bVtp)  -- (dotsV);
\draw[->, thick, densely dashed] (Psitp) -- (dotsPsi);
\draw[->, thick, densely dashed] (bEtp)  -- (dotsE);

\draw[->, thick, densely dashed] (bVtp) -- (dotsPsi);

\end{tikzpicture}

\end{minipage}%
\hfill%
\begin{minipage}{0.55\textwidth}

\centering
\begin{tikzpicture}[
  scale=0.7,         
  transform shape,    
  node distance = 1.8cm and 4cm,
  >=Stealth,
  box/.style={
    draw,
    rounded corners,
    minimum width=3.5cm,
    minimum height=2.1cm,  
    align=center
  }
]

\node[box] (jadv)
  {$j$-advancing  };

\node[box, below=of jadv] (jfill)
  {$j$-filling  };

\node[box, below=of jfill] (jrefill)
  {$j$-refilling  };

\node[box, right=of jadv] (jp1adv)
  {$(j+1)$-advancing};

\node[box, below=of jp1adv] (jp1fill)
  {$(j+1)$-filling};

\node[box, below=of jp1fill] (jstab)
  {$j$-stabilising  };

\draw[->] (jadv)   -- (jp1adv);      
\draw[->] (jadv)   -- (jp1fill);   

\draw[->] (jfill)  -- (jp1adv);    
\draw[->] (jfill)  -- (jp1fill);    
\draw[->] (jfill)  -- (jrefill);    

\draw[->] (jrefill) -- (jstab);    
\draw[->] (jstab)  -- (jfill);    

\draw[->] ([yshift=-0.25cm] jstab.west) -- ([yshift=-0.25cm] jrefill.east); 
\draw[->]
  (jrefill.west)
    .. controls +(-1,0.6) and +(-1,-0.6)
    .. (jrefill.west);

\end{tikzpicture}

\end{minipage}
\caption{Left: The high-level state $\Psi(t)$ is the interface between $V$ and $E$ at time $t$.
Right: Our specific high-level state transition algorithm (described in Section~\ref{sketch:ourhls})}
\label{fig:hls-trans}

\end{figure}

While bin~$j$ is filling, the volume process follows a $j$-jammed process.
The escape process $E$ also acts as a $j$-jammed process, except that it has no births. 
Instead, it has independent Bernoulli arrivals in each bin at the end of each step -- these account for balls that escape in the backoff process, but are modelled with independent Bernoulli variables to obtain more independence. The arrival rates depend on the interface between~$V$ and~$E$.
We must design the interface so that these arrival rates are large enough to maintain the coupling invariant but small enough to enable us to prove the proof property. 
The interface between the volume process~$V$ and the escape process~$E$ at time~$t$ is captured by an object $\Psi(t)$ which we call a ``high-level state''.
The evolution of high-level states is intimately connected
with the evolution of~$V$.
Also, $\Psi(t)$ controls the arrival rates in the escape process~$E$ at time~$t$.  $E$~does not depend on~$V$, except through the high-level states.
A depiction of the evolution of these processes is given in Figure~\ref{fig:hls-trans} (left picture).    

In general, a high-level state 
(Definition~\ref{def:high-level-state-newer}) is a tuple
$\Psi = (g^{\Psi},\tau^{\Psi}, j^{\Psi},\zvect^{\Psi},\calS^{\Psi},\type^{\Psi})$. 
At time~$t$, the process is in high-level state~$\Psi(t)$, which may or may not differ from $\Psi(t-1)$. 
If $\Psi = (g,\tau,j,\zvect,\calS,\type)$ then
$\tau$ is a non-negative integer, $j$ is a positive integer, $\zvect$ is a $j$-tuple of non-negative reals, and $\calS$ is a set of disjoint sets of bins in~$[j]$.
(We will explain the parameters~$g$ and $\type$ later.) 
Recall that while the process is in the high-level state~$\Psi = (g,\tau,j,\zvect,\calS,\type)$, the volume process $V$ acts as a $j$-jammed channel. More precisely, let $F_{[j]}^{\Psi}$ be the function that tracks the expected bin-populations in a $j$-jammed process, starting from $\Po(\zvect)$ at time~$\tau$. Then the volume process maintains the following very important invariant:
\begin{equation*}\tag*{(Volume invariant)}
\mbox{Conditioned on $\Psi(\tau_0),\dots,\Psi(t)$ and $\Psi(t) = \Psi$, $\dist(\bbar^V_{[j^\Psi]}(t)) = \Po(F_{[j^\Psi]}^{\Psi}(t))$.}
\end{equation*}
While the combined process is in the high-level state~$\Psi$,
the invariant will be  very useful for lower-bounding the population of 
bins in $\bbbar_{[j]}^V$ (which is useful  for the proof property). In particular, it 
solves
the problem identified at the beginning of this section --- when a backoff process~$X$ is
filling a bin $j$ that is weakly exposed it starts in a good situation with $\dist(\bvect^X_{[j-1]}(t)) \gtrsim \Po(\epsilon \lambda W_1,\ldots,\epsilon \lambda W_{j-1})$
but when this good situation 
fails the expected noise drops to $\Omega(\log \log j)$. 
Nevertheless, we do not lose the helpful Poisson domination.

It is easy to see roughly why the volume invariant is true --- the volume process is essentially following a $j$-jammed process. Balls behave independently in a $j$-jammed process~$Y$, and the distribution of 
$\bbar^Y_{[j]}(t)$ 
is exactly $\Po(F_{[j]}(t))$  (Lemma~\ref{lem:slowfill-j-jammed}) where 
the function~$F_{[j]}$  follows a simple recurrence (similar to the function~$f_{[j]}$ that we have already discussed),  tracking the means of bin populations in a $j$-jammed process starting from distribution~$\Po(\zvect)$ at time~$\tau$. So the volume invariant is ``easy'' as long as the high-level state doesn't change. In fact, until the high-level state changes, $V$ will be closely coupled with a $j$-jammed process $Y$ from $\tau$ so that $Y$ can be used to determine the next high-level state.

So why does the high-level state ever change? 
The purpose of the high-level state is to provide enough information to the escape process~$E$ so that the independent end-of-step Bernoulli arrivals can be tuned with the correct parameters -- these arrival rates must be large enough to maintain the coupling invariant, but small enough to enable us to prove the proof property. 
If the volume process 
is well-populated at time~$t-1$, so that for sets~$S \in \calS^\Psi$ the expected noise 
$\calN_S^V(t):= \sum_{i\in S} p_i b_i^V(t-1)$ from bins in~$S$ is high, and if the escape process is not too overloaded at time~$t-1$ in those same sets $S$, 
then the coupling invariant at time~$t-1$ guarantees that balls are not very likely to escape in the backoff process~$X$ at time~$t$, so the end-of-step arrivals in~$E$ at time~$t$ can be small without breaking the coupling invariant at time $t$. This helps with the proof property. 
To make it work, the combined process therefore maintains the following important invariant for the escape process.
\begin{align*}
&\mbox{Condition on $\Psi(\tau_0),\ldots,\Psi(t)$
and suppose that $\Psi(t) = (g,\tau,j,\zvect,\calS,\type)$.
Then
for all $S \in \calS$,}\\\tag*{(Escape invariant)}
&\mbox{$\calN_S^Y(t+1) \ge \epsilon \lambda |S|$
where $\bbbar^{Y}$ is 
the $j$-jammed process starting from $\Po(\zvect)$ at time~$\tau$.}
\end{align*}

As we have seen,   the end-of-step arrival rates in the escape process~$E$ may be  tuned to small values 
as long as
the expected noise in~$V$ from bins in sets in~$\calS$ is high, all without breaking the coupling invariant.
However, the escape process has no independent information about noise in~$V$ --- all that it sees is the high-level state. Thus, to make things work, the volume process~$V$ must change the high-level state whenever the noise from set $S\in \calS$ falls below what is expected, and this is what we get from the escape invariant.\footnote{The fact that the escape invariant is written in terms of the underlying $j$-jammed process~$Y$ rather than in terms of the volume process itself helps with conditioning, but to think about it here, it helps to note that the volume process agrees with~$Y$ from~$\tau$ until the high-level state changes.}

At this point we can say something about $\calS$. 
If the high-level state is 
$(g,\tau,j,\zvect,\calS,\type)$
and bin~$j$ is covered, then $\calS$ contains a single set~$S$, namely the reservoir~$\calR_j$ which can be expected to provide noise while bin~$j$ is filling.
If bin~$j$ is weakly exposed, then the role of~$\calS$ is more complicated. In this case, we initially expect quite a lot of noise towards filling bin~$j$, but this noise will fail, and the combined process cycles through a sequence of high-level states capturing the available noise (sometimes lots of noise, sometimes merely noise from the $O(\log \log j)$ very high-weight bins in~$[j-1]$) and these high-level states change in order to appropriately change the end-of-step arrival rate in~$E$.

So we've seen one situation where~$V$ has to change the high-level state --- namely running out of noise. There is another situation where $V$ needs to change the high-level state --- this is when   the volume invariant implies that bin~$j$ has filled, and it is time 
for the process to move on to filling bin~$j+1$, thus entering a new high-level state where the volume process will act as a $(j+1)$-jammed process.

All of this explains why the algorithmic dependency arrows go both directions between the 
evolution of the volume process and the 
evolution of the high-level states in the depiction of the combined process (Figure~\ref{fig:hls-trans}, left). Lack of noise in~$V$, or progress on filling~$j$ in~$V$, causes  $V$ to trigger a high-level state change. However, the volume invariant that we need is conditioned on the high-level state and therefore the high-level state~$\Psi(t)$ is a factor in determining the next volume state~$\bbar^V(t)$. 
The transitions depicted by the figure
are therefore as follows, where $\hls$
is a particular high-level state transition function that we will describe in Section~\ref{sketch:ourhls}.

First,
if $\Psi(t-1)$ is the high-level state
$\Psi = (g,\tau,j,\zvect,\calS,\type)$, then 
$\Psi(t) = \hls(\Psi, \bbar^{Y}_{[j]}(t), t)$
where $Y$ is the $j$-jammed process starting from $\Po(\zvect)$ at time~$\tau$. Now that $\Psi(t)$
is determined, 
$\bbar^V_{[j^{\Psi(t)}]}(t)$ needs to be defined in a way that
establishes
the volume invariant. 
Rather than just continuing the $j$-jammed process
$Y$ for another step,  the
state $\bbbar^V_{[j]}(t)$ needs
to be adjusted using a (known) coupling between 
the distribution of 
$\bbbar^{Y}_{[j]}(t)$, which is 
$\Po(F^{\Psi}_{[j]}(t))$, and the desired
distribution 
$\Po(F^{\Psi(t)}_{[j^{\Psi(t)}]}(t))$.

Everything that we have sketched so far is infrastructure for our proof, by which we mean that it 
is our intuition for how to prove instability for any such process but it
does not yet present our  intuition about how this specific backoff process evolves.
That intuition enters into the design of our actual high-level state transition algorithm~$\hls$, which is sketched in Section~\ref{sketch:ourhls}. The ``infrastructure'' is the collection of mathematical structures that we have built to make the proof work. In particular, the notion of the volume process, the escape process, and the high-level state transition algorithm that sits between them, and the invariants that are required. It is too complicated to prove all of the invariants for our actual high-level state transition algorithm, so what we do instead, in Section~\ref{sec:hls-properties} (and in the corresponding part of the full proof) is to distill all of the properties that a high-level state transition algorithm must have in order to facilitate a coupling that maintains all of the invariants that we just described. The main issue is conditioning -- the invariants are conditioned on the sequence of high-level states, and in particular, the Poisson domination of the volume process needs to be conditioned on this to maintain the volume invariant. In order to facilitate that, the high-level state transition algorithm needs to be designed with care, and the issues are sketched in Section~\ref{sec:hls-properties}. 
After that, we will sketch the description of our actual high-level state transition algorithm in Section~\ref{sketch:ourhls}. In our full proof, we will need to show that our high-level state transition algorithm satisfies the properties that we've required  and also that it enables us    to establish the proof property.

\subsubsection{Key properties of the high-level state transition algorithm}
 \label{sec:hls-properties}

Here we briefly pick up on a few of the properties that the high-level state transition algorithm~$\hls$ must have in order to enable the invariants from Section~\ref{sec:proofsketch}.
We start with the escape process.
Escapes in backoff processes are dependent and messy, so we have simplified things  in the escape process~$E$ by dominating escapes with (carefully tuned) independent Bernoulli end-of-step arrivals. 
Even with this simplification, we are not able to get a clean Poisson domination of the growing population of the escape process (this is because changes to the high-level state are not controlled by the escape process, so we will need to expose the escape process each time the high-level state changes).
From the point of view of the escape process, the high-level state changes at unpredictable times -- its change is entirely driven by the volume process.
What we do in the analysis is to derive pseudorandomness conditions  (invariants) that the population of~$E$ satisfies while in a particular high-level state and 
we show that, if these invariants hold
when the high-level state is entered, then the invariants required by the next high-level state are likely to hold by the time that high-level state is entered. 

Making this work constrains the high-level state transition algorithm~$\hls$ in the following way.
The design of the algorithm needs to ensure that  the volume process does not change the high-level state too often (so that these pseudorandomness conditions can be re-established) though of course, the volume process must change the high-level state when it runs out of noise, so there is a tension there.
 
In Section~\ref{sec:proofsketch} we noted two situations where the volume process needs to change the high-level state. The first is due to having too few balls - it runs out of expected noise in the relevant sets $S\in \calS$. The second is essentially due to having many balls - either it finishes filling bin~$j$ and moves on to bin~$j+1$, or it merely has enough noise, and continues in its current $j$-jammed process (without changing the high-level state). Either  
way, the volume invariant demands that 
the distribution of $\bbar_{[j]}^V(t)$ be an independent Poisson tuple \emph{conditioned on the sequence of high-level states up to $\Psi(t)$}.
We know from our analysis of $j$-jammed channels that the distribution of their bin populations is a Poisson tuple, but here we need the conditioning.
In order to condition on low noise 
(or to accommodate the increase in~$j$) 
we must ensure that the desired distribution 
$\Po(F^{\Psi(t)}_{[j^{\Psi(t)}]}(t))$
is stochastically dominated above by
$\Po(F^{\Psi(t-1)}_{[j]}(t))$, the distribution that would be obtained by taking another step in the $j$-jammed process from $\tau^{\Psi(t-1)}$.  To enable this, 
we constrain the high-level state transition function
(see (V\ref{item:valid-zlower}) in the hypothesis of Lemma~\ref{lem:valid-transition})
so that, if $\Psi(t)$ differs from $\Psi(t-1)$, then 
for all $\ell \in [j^{\Psi(t-1)}]$,
$z_\ell^{\Psi(t)} \leq F_\ell^{\Psi(t-1)}(t)$.
Furthermore, we constrain the high-level state transition 
function (see (T\ref{item:trans4}) in Definition~\ref{def:transition-rule}) so that
for all $\ell \in [j^{\Psi(t)}] \setminus [j^{\Psi(t-1)}]$,
$z_\ell^{\Psi(t)} =0$.

In order to condition on high noise, which allows the high-level state to stay unchanged, we use Harris's inequality (Lemma~\ref{lem:Harris}) which shows that the probability of any two (dependent) events is dominated below by the product of their probabilities as long as the events are increasing functions of some independent set of variables. 
For this to be applicable, the high-level state transition algorithm~$\hls$ has to be constrained so that 
remaining in a high-level state is an increasing event (see (V\ref{item:valid-up}) in the hypothesis of Lemma~\ref{lem:valid-transition}).
In particular, we require that 
the high-level state transition algorithm $\hls$ be designed in such a way that
increasing the
entries of $\bbar^{Y}_{[j]}(t)$ in the transition
$\Psi(t) = \hls(\Psi(t-1), \bbar^{Y}_{[j]}(t), t)$
cannot change the output of $\hls$ from $\Psi(t-1)$ (i.e.\ remaining in the same state) to another state.
Actually, we similarly need to apply Harris's inequality 
in order to condition on high noise in certain changes to the high-level state. Fortunately, we do not need these events to be completely increasing (which would be impossible to achieve) but it suffices to have them increasing in bins~$\ell$ such that $z_\ell^{\Psi(t)}$ is not~$0$ (see (V\ref{item:valid-up-non-zero}) in the hypothesis of Lemma~\ref{lem:valid-transition}).
In particular, this means that whenever the high-level state changes due to a set of bins not containing enough balls, those bins  must be completely emptied in $V$.

We have already mentioned that the high-level state has to change when the volume process runs out of noise. This is captured by (T\ref{item:trans3}) in Definition~\ref{def:transition-rule} which requires that,
for all $S \in \calS^{\Psi(t)}$, $\sum_{x \in S} p_x b^Y_x(t) \geq \delta \lambda |S|$ for some $\delta$.

The coupling is given in Section~\ref{sec:VEB-couple} assuming a generic high-level state transition algorithm~$\hls$ 
that satisfies all of the necessary properties.
A detail that may help the reader to pattern-match against this introduction is that, in the full proof, we don't have one escape process~$E$. Instead, a new escape process is started every time the high-level state changes -- this is to smooth writing, because we expose the state of the escape process every time the high-level state changes. Thus, there is an infinite sequence of escape processes $E_1,E_2,\ldots$ These are indexed by the parameter~$g$ of the high-level state. In the coupling, the sends of the balls are synchronised so that they stay in the same bin in the various processes.

A good way to think about the escape
processes is as follows.
Moving on from the start state, within an escape process, balls evolve independently (and any arrivals of new balls after the process started are also independent) though the rate of end-of-step arrivals depends on the existing balls in the escape process. 
Thus, by defining the escape process in this way, we gain independence from the volume process, which enables us to upper-bound the number of balls in the escape processes. 
What we gave up for that is that the coupling between backoff and the volume/escape processes isn't perfect - that is, the backoff process does not precisely match the difference between the volume processes and the escape processes, but the domination is in the right direction -- every ball that is in the volume process but not the relevant escape process is in the backoff process, so we  an maintain the coupling invariant.

We note that the properties that we have identified for the high-level state transition algorithm are stringent and we cannot expect them to be maintained deterministically. Instead, one of the possibilities for the type parameter in a high-level state~$\Psi$ is~$\Failure$. The high-level state transition function~$\hls$ is designed to transition to a failure state whenever anything goes wrong, and to stay in a failure state forever after, in that case.
When we analyse the volume process (Section~\ref{sec:backoff-bounding-analysis}), we must ensure that all transitions to failure occur with summable failure probability. To do so, we analyse how long the volume process stays in a given high-level state and what high-level states it is likely to transition to from there. This information is then also used for analysing the escape processes (Section~\ref{sec:escape-analysis}) -- for example, to argue that $j$ advances sufficiently quickly so 
the escape process doesn't build up too many balls. 

\subsubsection{A brief description of our specific high-level state transition algorithm and its properties} 
\label{sketch:ourhls}

Recall that a high-level state is a tuple of the form 
$\Psi = (g,\tau,j,\zvect,\calS,\type)$.
The possibilities for the $\type$ parameter are
$\Failure$, $\initialising$, $\advancing$, $\filling$, $\refilling$, and $\stabilising$.
We've already mentioned~$\Failure$ -- the whole process transitions into a~$\Failure$ state in the unlikely event that something goes wrong.
We will not dwell on the~$\initialising$ state --- it is simply the first high-level state 
so that (conditioning on the event $\Einit$) we have  
$\dist(\bvect^X_{[j_0]}(\tau_0)) \gtrsim \Po(3 \lambda W_1/4,\ldots,3\lambda W_{j_0}/4)$  
and we start coupling the backoff process~$X$
with the volume/escape/high-level state process at time~$\tau_0$.
We refer to a high-level state~$\Psi$ 
as $j$-advancing, $j$-filling, $j$-refilling, or $j$-stabilising, otherwise (capturing $j$ in our description of the type of~$\Psi$). 
The basic operation of the high-level state transition algorithm is captured in Figure~\ref{fig:hls-trans} (right).
When $j$ is a covered bin currently being filled, a high-level state can only be $j$-advancing. 
When $j$ is weakly exposed, it can be $j$-filling, $j$-refilling, or $j$-stabilising.

The simplest case to understand is when bin~$j$ is covered.
Let $j$ be a covered bin and let $\Psi=(g,\tau,j,\zvect,\calS,\advancing)$ be a $j$-advancing state.
Recall from the volume invariant that
$\Psi(t) = \Psi$ entails
$\dist(\bbar^V_{[j]}(t)) = \Po(F_{[j]}^{\Psi}(t))$, where $F_{[j]}^{\Psi}$ is the function 
that tracks the expected bin-populations in a $j$-jammed process,
starting from
$\Po(\zvect)$ at time~$\tau$.
In particular, 
$\dist(\bbar^V_{[j]}(\tau)) = \Po( \zvect)$.
We therefore define $\zvect$ 
for the $j$-advancing state~$\Psi$
so that bins in $[j-1]$ are full and bin~$j$ is empty.
The set $\calS$ contains a single set $S=\calR_j$, which is the reservoir that can be relied on to supply noise as bin~$j$ fills.
It is likely that bin~$j$ fills in $W_j \poly(j)$ time and the high-level state then changes to $(j+1)$-advancing or $(j+1)$-filling (depending on whether or not bin~$j+1$ is covered).
In the unlikely event that this doesn't happen, the high-level state changes to a failure state (not depicted in the figure) -- the probability that this happens is small enough that we can sum it over~$j$.
We next consider what happens to the escape process while the whole process is in the $j$-advancing state~$\Psi$. Basically, since $S = \calR_j$ has size at least~$\L(j)$, 
the expected noise is high in a $j$-advancing state.
So from a $j$-advancing state an escape process that is not over-full can enjoy a very small end-of-step arrival rate (as small as $1/j^C$ for some big constant~$C$). Under this arrival rate, the stationary population of $b_k^E$ is approximately $\Po(W_k/j^C)$ so, even though its population gets exposed at time~$\tau$ when high-level state~$\Psi$ is entered, 
the escape process $E$ empties quickly. This makes it possible to establish the pseudorandomness condition that will be required when the high-level state transitions to $(j+1)$-advancing or $(j+1)$-filling.

The case where $j$ is weakly exposed is more complicated. On the step~$\tau$, 
when high-level state $\Psi=(g,\tau,j,\zvect,\calS,\filling)$ is first entered,
the volume process has filled bins~$[j-1]$ and this is reflected in~$\zvect$.
In high-level state~$\Psi$, $\calS$ contains a single set~$S_1$ containing all low-weight bins with weight above some~$\Wtilde[j]$ but it also contains singleton sets for every very high-weight bin $i\in [j-1]$. 
The point is that $\Wtilde[j]$ is chosen so that $\bbbar_{S_1}^V$ is likely to 
provide $\Omega(|S_1|) = \Omega(\log j)$ noise for $\exp(\Omega(j))$ time. During this  time, around $\epsilon$ balls are added   to $\b_j^V$ per step (in expectation, according to the $j$-jammed process), making progress towards filling bin~$j$. When the noise from~$S_1$ fails in~$V$ (as it is very likely to do, eventually) the high-level state changes to $j$-refilling.
The singleton very high-weight bins will be likely to provide $\Omega(\log\log j)$ noise until the process next enters a $j$-filling state (preventing it from transferring to a failure state instead!). These very high-weight bins are not in~$S_1$ so the high-level state transition to the $j$-refilling state, conditioned on the failure of noise in~$S_1$, does not have to reduce the population of these bins. 

When the process transitions at some step~$\tau'$ to the $j$-refilling state 
$\Psi'=(g',\tau',j,\zvect',\calS',\refilling)$ the bins in~$S_1$ have failed to provide
sufficient noise, so these have value~$0$ in~$\zvect'$.
There is now not much noise in the volume process, so 
$\calS'$ only contains singleton sets for the very high-weight bin $i\in [j-1]$. The purpose of the $j$-refilling state is to wait long enough   until it is likely that the bins in $\bbar_{[j-1]}^V$ are again full (recall that we have a Poisson domination on their expectation), so that the process may re-enter a $j$-filling state (via a $j$-stabilising state - which we discuss in a minute). If the transition into the
$j$-stabilising state fails (due to the fact that 
that the volume process has less noise than its expectation would suggest),  it simply starts a new $j$-refilling state (this is the meaning of the self-loop on ``$j$-refilling'' in Figure~\ref{fig:hls-trans} -- we do not depict transitions that do not change the high-level state, but there is a possibility of moving from one $j$-refilling state to another).
While the process is in the $j$-refilling state we do not expect much noise so this high-level state determines a high end-of-step arrival rate for the escape process~$E$.
Fortunately, the whole process is unlikely to stay in a $j$-refilling state for long, so this does not cause too much harm.

Nevertheless, by the time that the bins in $[j-1]$ have refilled (in expectation) in the volume process and the whole process is ready to transition out of $j$-refilling and back to $j$-filling, the escape process will have been enduring a high end-of-step arrival rate 
for a while (for the period of $j$-refilling) so it will not satisfy the pseudorandomness conditions necessary for a transition directly back to $j$-filling. This is the purpose of the $j$-stabilising state. When the process transitions at some step~$\tau''$ to the $j$-stabilising state 
$\Psi''=(g'',\tau'',j,\zvect'',\calS'',\stabilising)$ the volume process will have refilled bins in $[j-1]$ in expectation  and this is reflected in~$\zvect''$. From the point of view of the volume process, nothing happens in $j$-stabilising -- it is simply waiting to transition back to $j$-filling.
The purpose of the $j$-stabilising state is to enable the escape process to recover the stringent pseudorandomness condition for re-entry into a $j$-filling state (the proof relies on the fact that the probability that the whole process enters a $j$-filling state without having this condition satisfied has very small failure probability, which can be summed over~$j$). For this purpose, $\calS''$ has several sets~$S$.
These allow the end-of-step arrival rate in~$E$ to be dominated appropriately.
At first, shortly after~$\tau''$, the end-of-step arrival rate may be fairly high -- as high as $1/(\log j)^C$ for some large constant~$C$. 
(This is necessary, and is why the pseudorandomness conditions for refilling and stabilising are so much weaker than those for advancing and filling).
However, bins $1,\dots,(\log j)^2$ rapidly empty out, enabling domination with a much lower end-of-step arrival rate~$1/j^C$ for as long as these bins remain full in the volume process. To allow the escape process to rely on the fact that these bins are still full in volume, low-weight bins in this range are included as a set in $\calS$, and if they do empty out in $V$ (which they may) then the process transitions back to a $j$-refilling state. If no such transition occurs, the lower arrival rate enables $E$ to re-establish the pseudorandomness conditions needed for transitioning to the next $j$-filling state. The process may bounce back and forth between $j$-refilling and $j$-stabilising states up to $O(j)$ times before transitioning back to a $j$-filling state.

 Fortunately, the good progress that is being made every time the process is in $j$-filling 
outweighs the loss of balls in the volume process 
during $j$-refilling and $j$-stabilising, 
so bin~$j$
is likely to fill   within~$W_j \poly(j)$ steps, enabling an eventual transition to $(j+1)$-filling or $(j+1)$-refilling.

\subsubsection{Putting it all together (for the case with finitely-many strongly exposed bins)}
\label{sec:bottomline}

Recall from the beginning of Section~\ref{sec:finiteSE}
that the analysis started 
by choosing a finite time~$\tau_0$, depending on~$\p$ and~$j_0$, and a possible, but unlikely event~$\Einit$ that entails
$\dist(\bvect^X_{[j_0]}(\tau_0)) \gtrsim \Po(3 \lambda W_1/4,\ldots,3\lambda W_{j_0}/4)$ and furthermore, entails that the backoff process does not return to the empty state before step~$\tau_0+1$.
The analysis of the volume process, the escape process, and the high-level state transition algorithm
in Sections~\ref{sec:backoff-bounding-analysis} and~\ref{sec:escape-analysis}
enables us to draw the following conclusions about the time steps following~$\tau_0$.
First, from~Corollary~\ref{cor:VEB-volume}, conditioning on $\Einit$,
with probability at least~$2/3$,
\begin{itemize}
\item the high-level state transition algorithm never enters a failure state, and
\item for all  $j\geq j_0$, it takes at most $\poly(j)(W_{j_0} + \cdots + W_{j-1})$ steps to reach a $j$-advancing or $j$-filling state.
\end{itemize}

Also from Lemma~\ref{lem:VEB-escape}, conditioning on $\Einit$, with probability at least~$2/3$,
\begin{itemize}
\item The very high-weight bins are never over-full in~$E$,
\item and bins in~$\calR_j$ are never over-full in~$E$ in a $j$-advancing state.
\end{itemize}

Putting these together (via a union bound) in the proof of Theorem~\ref{thm:goaljkillerSE},
and using lower-bounds on the population of the volume process from the constraints on the high-level state transition algorithm, we obtain that, 
conditioning on $\Einit$,
with probability at least~$1/3$,
\begin{itemize}
\item the high-level state transition algorithm never enters a failure state, and
\item for all $t>\tau$ such that the high-level state is not a failure state, there a bin~$i$ such that $\b_i^V(t) \setminus \b_i^E(t)$ is non-empty.
\end{itemize}

We next apply the coupling invariant.
$\Einit$ occurs with some positive probability~$\epsilon$. Thus, 
with probability at least~$\epsilon/3$, 
the backoff process does not visit the empty state before time~$\tau_0+1$ 
and after that   there is always a bin~$i$ such that $\b_i^X(t)$ is non-empty.
We conclude that the backoff process~$X$ is not positive recurrent (in fact, in this case it is transient, though we don't get a full proof of transience, see Section~\ref{sec:discussion}).

\subsubsection{Infinitely-many strongly exposed bins}\label{sec:intro-inf-SE}

Suppose that $\p$ is a send sequence with infinitely many strongly-exposed bins. In this section we briefly discuss
how to prove that for any positive birth rate~$\lambda$,
the backoff process with send sequence~$\p$ and birth rate~$\lambda$ is not positive recurrent.
For this we define (Definition~\ref{def:bottleneck}) the notion of a ``bottleneck'' bin.
Bin~$i$ is a bottleneck if  
$W_i \ge \exp(i/\exp((\log\log i)^2))$, 
for all $\ell \in [i-1]$,  $W_\ell \le W_i^{1/\Theta(\log \log i)}$,
and there is a  
strongly exposed bin~$i^+$ in the range $i \leq i^+ \leq i \exp({(\log i)}^{1/2})$.
Lemma~\ref{lem:nearbybottleneck} shows that if $j$ is sufficiently large and bin~$j$ is strongly exposed, then 
there is a bottleneck bin~$i$ with $j=i^+$. Thus, the range on~$i^+$ in the definition of bottleneck ensures that
having infinitely many strongly exposed bins implies having infinitely many bottlenecks -- this is captured in 
Corollary~\ref{cor:infinite-KMP-bins-2}.

The point of bottlenecks is that, if bin~$i$ is a bottleneck, then it contributes 
at least $\exp(i^{1/2})$ to the expected sojourn time from the empty state back to the empty state.
The key lemma for showing this is Lemma~\ref{lem:GetIntoi}
which, roughly, shows that, with probability at least~$W_i^{-6/7}$, a newborn ball born at step~$1$ gets to bin~$i$
without escaping.
The important aspects of the argument are as follows.
\begin{itemize}
\item   The upper bound on the weights of bins in $[i-1]$ ensures it is sufficiently likely that the newborn takes only $i^{1+o(1)}$ steps to send from each of the bins in $[i-1]$, provided that it doesn't escape. This is essentially
Lemma~\ref{lem:bottleneck-geo-bound} -- the analysis is like an externally-jammed process, so it is about summing the geometric waiting times. 
\item On the other hand, noise from other newborns is sufficiently likely to block the first $O(W_0 + \cdots + W_{\log i})$ steps during this journey (Lemma~\ref{lem:slow-fillstreams}).
\item Finally (Lemma~\ref{lem:bottle-events-TwoStream}), the resulting
balls in the first $\log i$ bins are sufficiently likely to block escapes during the remaining steps when the newborn is on its way to bin~$i$.
\end{itemize}

If the newborn reaches bin~$i$ (which happens with probability at least $W_i^{-6/7}$) then it takes roughly~$W_i$ time before it sends from bin~$i$ so its contribution to the expected sojourn time is at least 
$W_i W_i^{-6/7}$, which is at least  $\exp(i^{1/2})$ given the lower bound on~$W_i$, enabling us to show that the expected sojourn time is unbounded since there are infinitely many bottlenecks.

In order to make this proof rigorous 
(in order to rigorously handle the dependence between balls) we dominate the backoff process with two processes that have slightly more independence. One of these is the two-stream process from~\cite{GL-oldcontention} 
(Definition~\ref{def:two-stream}).  In a two-stream process, escapes are only blocked by balls in the other stream, so a two-stream process has a little more independence than a backoff process. It can be coupled with a backoff process (Definition~\ref{def:standard-twostream})
so that balls in the two-stream process are also in the backoff process. The other dominating process   is called an   ``under-backoff process'' (Definition~\ref{def:under-backoff}) -- this process has three streams. Streams~$A$ and~$B$ operate together as a backoff process  (ignoring stream~$C$). However, escapes in stream~$C$ are blocked by all three streams. Our analysis follows a newborn in stream~$C$ but the blocking comes from streams~$A$ and~$B$ (and the noise there comes by domination with a two-stream process). The under-backoff process can be coupled with a real backoff process (Definition~\ref{def:standard_under}) so that balls in the under-backoff process are also in the backoff process, so the the bound on the probability that the newborn in stream~$C$ reaches bin~$i$ gives a similar bound in the backoff process.

\subsection{Context} \label{sec:context} 

In this section, we review some of the academic context of our work. In particular, we review work on contention-resolution in Kelly's model~\cite{Kelly}. Work in other models will be discussed in Section~\ref{sec:related-work}.

In the paper~\cite{kelly-macphee},
Kelly and MacPhee gave  a complete characterisation of send sequences for which, with probability~$1$, the number of balls that escape is finite. Clearly, a backoff process with such a send sequence is not positive recurrent.
The following corollary of their theorem, from~\cite{GL-oldcontention},
shows that the characterisation applies if the weight  of bin~$i$ grows too slowly as a function of~$i$.\footnote{The statement of Corollary~\ref{KM-cor} as Corollary~9 of~\cite{GL-oldcontention} restricted attention to $\lambda\in (0,1)$ and omitted the part about finitely many balls escaping, but the version stated here follows directly from Kelly and MacPhee's theorem and from the proof in~\cite{GL-oldcontention}.}
\begin{corollary}\label{KM-cor}
(This is \cite[Corollary 9]{GL-oldcontention}, a corollary of \cite[Theorem 3.10]{kelly-macphee}) Let $\p$ be a send sequence such that $\log(W_i) = o(i)$ as $i\rightarrow \infty$. Then for all $\lambda >0$, the backoff process~$X$ with send sequence~$\p$ and birth rate~$\lambda$ is not positive recurrent. Moreover, with probability~$1$, only finitely many balls escape from~$X$.
\end{corollary}

We have
have already discussed the proof (and conjecture) of Aldous~\cite{Aldous}, which was the direct inspiration for this work.
MacPhee~\cite{MacPhee} remarked on the difficulty of extending Aldous's proof:
{\it ``Aldous believes that the proof can be modified to cope with non-exponential backoff. The changes required to do this for any particular case 
are necessarily quite substantial and the possibility of finding a considerably larger class of backoff regimes than the exponential ones for which a reasonable proof (by which is meant a proof that does not degenerate into the enumeration of a long list of different cases) seems remote.''}

The authors of~\cite{GJKP} obtained
an instability result that applies to all backoff protocols
but this result is much weaker than Aldous's conjecture, showing instability for birth rate $\lambda \geq 0.42$.

One crucial feature of Aldous's proof is that bin weights grow exponentially in binary exponential backoff. This provides concentration -- the probability that things go wrong while bin~$j$ is filling is so low that it may be summed off with a union bound.
Send sequences that include low-weight bins lack this concentration, since low-weight bins will often be much emptier than their expected population sizes. Any instability proof that also applies for send sequences with many low-weight bins   therefore needs to avoid relying on such union bounds and to instead consider the (messy) joint distribution of the bin populations. This is what is done in~\cite{GL-oldcontention} - where   the distribution of $(b_1^X(t),\ldots,b_j^X(t))$ is dominated below by a $j$-tuple of independent Poisson random variables. 
This gave instability for $\lambda>0$ for a   a wide class of send sequences~$\p$, including all monotonically non-increasing sequences and, more interestingly, all sequences such that the median of $p_0,\ldots,p_k$ is $o(1)$ as a function of~$k$. Nevertheless, the proof 
relied crucially on the fact that the expected noise stays sufficiently high. Thus, the proof technique is inssuficient for proving Aldous's conjecture.

We conclude this section by briefly discussing an explicit backoff protocol that goes quiet infinitely often.
Specifically, note that
the send sequence $\p = (p,p,p,\ldots)$ 
from the slotted ALOHA protocol leads to a backoff process with ever-increasing expected noise. Similarly, the send sequence of exponential backoff with $p_k=C^{-k}$ leads to a backoff process with ever-increasing expected noise. But it turns out that these can be interleaved so that each $p_k$ is taken from one or the other of these sequences 
(so 
$p_k \in 
\{p,C^{-k}\}$) in such a way that the 
resulting backoff process goes quiet infinitely often. This example is from a family of send sequences from~\cite{GL-oldcontention}. 
Roughly, the interleaving   is done so that  
the new send sequence~~$\p$
looks like slotted ALOHA for all~$k$ in some initial interval, then it looks like exponential backoff while $k$ ranges over a (longer) interval, then it looks like slotted ALOHA  while $k$ ranges over an (even longer) interval, and so on. The behaviour of the backoff protocol with this send sequence is interesting. 
Every $k$ that marks the beginning of   an ``exponential backoff interval'' 
is a strongly-exposed bin. 
The likely behaviour of the backoff process
is that it eventually enters a state in which \emph{all messages in the system} have collided exactly~$k$ times, but the number of these messages is substantially smaller than~$C^k$, so the expected noise is near~$0$. This happens infinitely often, because it happens for every strongly exposed bin~$k$. Thus, the protocol alternates between ``quiet periods'' with very low expected noise (expected noise near~$0$) and ``noisy periods'' in which the expected noise first grows rapidly, then stays high for a very long time, before going quiet again.
For concreteness, here are some parameters that work.
Let $C$ be a sufficiently large integer (for example, $C=10$).
let $\avect = (a_0,a_1,\ldots)$ be a fast-growing sequence of integers, indexed by the natural numbers, such that $a_0=0$. For example, consider the sequence with $a_0=0$ given by 
$a_k = (2C)^{a_{k-1}}$. Then
let $\p$ be the send sequence defined 
as follows.
 
 \[
p_j = \begin{cases}
C^{-j} & \mbox{ if, for some even $k\geq 0$, }
a_{k} \le j \leq a_{k+1}-1,\\
1/2 & \mbox{ if, for some odd $k\geq 0$, }
a_{k} \le j \le a_{k+1}-1.
\end{cases}
\]  

For every sufficiently large even~$k$, bin~$a_k$ is strongly exposed and it is 
likely that the bins in $[a_k-1]$
simultaneously empty before there are any sends
from bin~$a_{k}$, so the process is silent.

\subsection{Related work}\label{sec:related-work} 

In this section we briefly discuss work on contention resolution in other models that are related to our model.
This is not meant to be a comprehensive treatment -- the purpose is just to mention some highlights.

First, we alluded in the introduction to the reasons that the Ethernet network is stable, despite the fact that binary exponential backoff is unstable. 
The first reason is that dropping messages that have collided at least~$K$ times (for any positive integer~$K$) makes any backoff process stable (at the expense of reliability since the dropped messages are just dropped!). The second reason is that, when the relevant application makes it possible to supplement backoff processes with queues, this can sometimes lead to stability.
In some applications (including Ethernet networks, but not including sharing of resources in the cloud) it makes sense to think of a multiple-access channel as having $n$ queues
for some fixed constant~$n$. Messages at the heads of the queues run backoff protocols and newborn messages arrive at the tails of queues. 
These queue-assisted backoff processes are sometimes stable, mainly because the queues ensure that there cannot be more than $n$~sends at any one time.
The original Ethernet
makes both of these adjustments, with $K=16$ and $n=1024$~\cite{kelly-macphee}. 

We briefly mentioned in the introduction the work of \Hastad, Leighton, and Rogoff in the queue-assisted model~\cite{HLR}. This is important work which should be better known. In this paper, \Hastad{} et al.~show that (in the queue-assisted model) binary exponential backoff is unstable for birth rates~$\lambda>1/2$ (assuming that $n$ is sufficiently large) but polynomial backoff is stable for any birth rate that is less than~$1$.
It seems from this result that polynomial backoff should be better known and more-often used.
It would be interesting to know the exact threshold where binary exponential backoff becomes stable in this setting. We believe that the best that is known~\cite{hesham} is that it is stable for birth rates up to $O(n^{-0.75})$.

There is quite a bit of interesting work on contention resolution in settings where messages have strictly more information than they do in the original foundational model of Kelly. 
For example, algorithmic improvements are possible in the \emph{full sensing} model, where messages have the ability to continually  ``listen'' to the channel, even when they do not send, allowing all messages at each step to distinguish between some or all of these outcomes -- silence, a successful send, or a collision. 
For historical work on variants of this model, we refer the reader to the review article of Ephremides and Hajek~\cite{EH}, the paper of Greenberg, Flajolet, and Ladner~\cite{GFL}, and the survey of Chlebus~\cite{Chlebus}.

For contention resolution with both queue-assistance and full-sensing some very powerful protocols are known.
In particular, there is a stable full-sensing protocol for the model where listening messages can distinguish silence from the other outcomes, 
even for the more general model in which some specified pairs of queues are allowed to use the channel simultaneously~\cite{ShahShin,SST}.

More recently,
Bender et al.~\cite{BKKP-adversarial}
consider a full-sensing model without collision detection, meaning that messages are allowed to continually listen to the channel, and they can distinguish success from the other outcomes, but they are unable to distinguish silence from a collision. In this model, Bender et al.{}
give a protocol which achieves constant throughput, 
even when the message arrival is adversarial rather than random. 
What this means is that there is a uniform positive lower bound on $n_t/s_t$, where $n_t$ is the number of messages that arrive up to time~$t$ and $s_t$ is the number of steps up to time~$t$ when the channel is non-empty.
Chen et al.~\cite{CJZ-jamming} 
gives an upper bound, in the same model, on the best possible throughput that can be achieved under a given amount of adversarial jamming. They give a contention-resolution algorithm that  achieves the optimal throughput in the model.

Bender et al.~\cite{Bender2025SICOMP} 
considered a ternary-feedback full-sensing model with adversarial arrivals and adversarial jamming. (Ternary feedback means that all three possibilities can be distinguished.) They gave an algorithm that achieves constant throughput  while guaranteeing (with high probability) that the number of times any message accesses the channel in the first $t$ steps (either to send or to listen) is at most a polynomial in the logarithim of $n_t + j_t$, where $j_t$ is the number of jammed slots up to time~$t$.

Continuing the theme of limiting channel access,    
\cite{ 
XZ-SPAA25} considers a 
full-sensing model without collision-detection and
with $n$ messages.  
They give a randomised algorithm which,
under  some mild additional assumptions,
enables all messages to escape in $O(n)$ time steps with high probability,  where   each message only sends or listens $O((\log n)^2)$ times. 

Finally, \cite{CCDKPP} studies contention resolution in settings where there are $n$ messages in total, with adversarial arrival times, but the feedback that these messages get is acknowledgement-based (whether or not they succeed when they send). They study latency (the amount of time that a message is in the system before it is sent successfully). In addition to deriving interesting tradeoffs on latency bounds, they show that latency can be improved substantially when messages have access to a global clock.

\subsection{Future work}\label{sec:introfuture}

Here we'd like to point out a  direction that we have not explored in this paper. In this paper, we study backoff protocols, where the probability that a message sends at each step is a function of $k$, the number of collisions that it has already had, so that, during each step following the $k$'th collision, the message sends with probability~$p_k$ and is silent with probability~$1-p_k$. Other possibilities exist. For example, it is easy to see that in a backoff protocol, the total amount of time that a message waits after its $k$'th collision is a geometric random variable with mean $W_k=1/p_k$. But we could imagine contention-resolution protocols where the total waiting time is chosen from some other (non-geometric) distribution (still depending 
on~$k$). Also, the random waiting time could potentially depend on other information, for example, the absolute time of the collision. Such contention-resolution protocols are referred to as \emph{acknowledgement-based protocols}. 
Kelly~\cite{Kelly} asked whether there is any   acknowledgement-based contention-resolution protocol that is stable  
for any positive birth rate~$\lambda$.  MacPhee~\cite{MacPhee} conjectured that there is not. See Section~\ref{sec:discussion} for a further discussion of open problems.

\subsection{Organisation of the paper (with commentary)}
\label{sec:outline}

The rest of this document is an appendix with all of the proofs. It is
 organised as follows.
 Section~\ref{sec:prelim} gives preliminary definitions, defining general notation and standard facts about probability distributions (dominations and tail bounds) -- it is best to skip this section, and to refer back to it as needed. Section~\ref{sec:prelim}  also proves that the backoff process is irreducible and aperiodic (a fact that we have already used to derive the consequences of stability). 
 
 Section~\ref{sec:backoff-defs} describes a stochastic process called a ``generalised backoff process''. Certainly a backoff process is a generalised backoff process, but in the course of our proof we will encounter various other processes that we will use for domination. The definitions in Section~\ref{sec:backoff-defs} provide a common framework. The definitions look a a little bit too general at first sight, but that turns out to be useful. In Section~\ref{sec:backoff-defs} we also introduce a very basic method (which we call the ``standard coupling'') for coupling these processes with each other. 
 Section~\ref{sec:reduction} provides a reduction to simplify our task. 
 The theorem that is stated there
 (Theorem~\ref{thm:goaljkiller}) is the same as Theorem~\ref{thm:goal} except that we are allowed to assume that $\lambda$ is sufficiently small, that $p_0=1$, and that every bin weight $W_j$ is upper-bounded by $W_j \leq (2/\lambda)^j$. In Section~\ref{sec:reduction} we prove Theorem~\ref{thm:goal} using Theorem~\ref{thm:goaljkiller}, so these simplifications are without loss of generality. The rest of the paper focuses on proving Theorem~\ref{thm:goaljkiller}.
 Section~\ref{sec:j-jammed} sets up a very important tool for our analysis, the $j$-jammed process.  This section provides detailed analysis of the movement of balls in a $j$-jammed process. We use this throughout our proof. (The proofs in this section may be skimmed, but it is good to have the definitions.)
 
 Section~\ref{sec:covered} provides the definitions of covered, weakly exposed, and strongly exposed bins that we have already discussed informally. It also derives some important consequences of these definitions. Section~\ref{sec:infinite} gives the proof for the case where there are infinitely many strongly exposed bins. The rest of the paper is therefore about the case where the number of strongly exposed bins is finite. 
 
 Section~\ref{sec:hls} gives our generic definitions of a high-level-state and derives the properties that a high-level state transition  rule must have in order to enable our coupling. Section~\ref{sec:vol} defines the volume process, including its intricate connection to the high-level state transitions. Lemma~\ref{lem:dist-volume} is important - showing that the distribution of the bin populations for the volume process are give by a Poisson tuple, even conditioning on the sequence of high-level states. Section~\ref{sec:def-escape} defines an escape process and Section~\ref{sec:VEB-couple} gives our coupling of a volume process, escape processes, and a backoff process. The description of the coupling is rather long, unfortunately -- it is just complicated.
 Given the coupling, we know that, conditioned on the sequence of high-level states, the distribution of the first $j^{\Psi(t)}$ bins of the volume process are distributed as 
 $\Po(\F^{\Psi(t)}_{[j^{\Psi(t)}]}(t))$. We also know 
 the collection $\calS^{\Psi(t)}$ of sets
 expected to have noise in the volume process at time~$t$, and these control the end-of-step arrival rates of the escape process at time~$t$. Crucially, we have the invariant that 
 $\b_i^{V}(t) \setminus \b_i^E(t) \subseteq \b_i^X(t)$. All of this is infrastructure that we have built, depending only on the assumption that the high-level state transition rule satisfies the properties that we laid down in Section~\ref{sec:hls}.
 We don't yet know anything about the likely sizes of the means $F_i^{\Psi(t)}$ or about the actual end-of-step arrival rates --- these depend on the actual high-level state transition rule, which comes next.

In Section~\ref{sec:sets}, we describe the specific types of high-level states that we will use for our proof. 
 This is followed by Section~\ref{sec:backoff-bounding-valid} which gives our high-level state transition algorithm, the backoff bounding transition rule. The idea behind this algorithm is one of the main ideas in this paper. 
 The particular high-level states and the way that transitions between them happen, according to the backoff bounding rule, is where our intuition about how the backoff process actually evolves enters the proof.
 Section~\ref{sec:backoff-bounding-valid} also proves that the backoff bounding rule satisfies all of the generic properties that we have required for a high-level state transition function in order to do the coupling. 
 
 Finally, Section~\ref{sec:backoff-bounding-analysis} analyses the volume process, Section~\ref{sec:escape-analysis} analyses the escape processes, and Section~\ref{sec:final} combines the results of those two sections to prove our theorem (for the case with finitely-many strongly exposed bins).
The bulk of the technical work of the paper is in Sections~\ref{sec:backoff-bounding-analysis} and~\ref{sec:escape-analysis}.
 
In Section~\ref{sec:discussion} we discuss open problems, and ideas for future work.
There is a symbol table at the end of the paper which may be useful for reference.

\section{Preliminaries}
\label{sec:prelim}

\subsection{General notation}

All logarithms in this paper are to the base~$e$.
We use $\newsym{Z}{integers}{\integers}$ for the set of integers and $\newsym{R}{reals}{\reals}$ for the set of reals.
We use a subscript to restrict this set, so for example $\integers_{\geq 0}$ is the set of non-negative integers and $(\integers_{\geq 0})^j$  is the set of $j$-tuples of non-negative integers.
Vectors are represented as tuples, and denoted with an overline, for example we write $\avect= (a_1,a_2,\ldots)$ to mean that $\avect$ is a vector indexed by the positive integers. 
Given a vector~$\avect$ and a subset~$S$ of its index set, $\newsym{avect-S}{sub-vector}{\avect_{S}}$ denotes the vector formed by taking the elements of~$\avect$ that are indexed by elements of~$S$.
For example, if $\avect$ is indexed by the positive integers, then  $\avect_{[j]}$  denotes the sub-vector $(a_1,\ldots,a_j)$.
We use $\newsym{0-vect}{zero vector}{\zerovect}$ to denote the all-$0$ vector -- its length will be clear from context. Similarly $\newsym{1-vect}{1 vector}{\onevect}$ denotes the all-$1$ vector.
Given two $j$-tuples $\avect$ and $\bvect$
with the same index set, 
$\avect \newsym{leq}{vector inequality}{\leq }\bvect$ means that for every index~$i$, $a_i \leq b_i$.
We often write ``let $x=y$'' in order to define a quantity $x$ but if we want to emphasise that it is a definition, we write
$x:=y$.  

\begin{definition}
\label{def:[r]}
For a real number 
$r$ we define 
$\newsym{[]}{Set of positive integers $\leq r$ - Def~\ref{def:[r]}}{[r]}$ to be the set of positive integers that are at most~$r$.
\end{definition}

The following trivial observation will be useful.

\begin{observation}
\label{obs:dumb} 
For any real numbers $A$, $B$, and $C$,
$\min\{A,B\} - \min\{C,B\} \geq \min\{0,A-C\}$
and $\max\{0,B-C\} \leq \max\{0,A-C\} + \max\{0,B-A\}$.
\end{observation}

\subsection{Probability distributions and domination}

In this paper a geometric distribution counts the number of Bernoulli trials up to and including success, so that e.g.\ the expectation of a geometric variable with parameter $p$ is $1/p$.

Given any random variable~$X$, we use $\newsym{dist}{distribution of any random variable~$X$} {\dist(X)}$ to denote its distribution.

\begin{definition}
\label{def:Poisson}
For every  non-negative real number~$\lambda$,
\newsym{Po}{Poisson distribution with mean $\lambda$}{\Po(\lambda)} is the Poisson distribution with mean~$\lambda$. 
\end{definition}

\begin{definition}\label{def:PoTuple}
For all integers $j\geq 1$ and all 
$\muvect \in  (\reals_{\geq 0})^j$, 
\newsym{Po-1}{tuple of independent Poisson r.v.s with given means}{\Po(\muvect)}
is the distribution with support~$(\integers_{\geq 0})^j$ such that
$(P_1,\dots,P_j) \sim  \Po(\muvect)$  means that $P_1,\dots,P_j$ are independent Poisson variables 
with 
 $\E(P_x) = \mu_x$.
\end{definition}

\begin{definition}\label{def:phi}
Let $A$ and $B$ be distributions on $(\integers_{\geq 0})^j$. $A$ \emph{dominates $B$ from above} 
(written $A \newsym{>}{Dominates from above}{\gtrsim} B$)
if there is a distribution $C$ on pairs $(\bar{a},\bar{b}) \in (\integers_{\geq 0})^j \times (\integers_{\geq 0})^j$ such that $\bar{a}\sim A$, $\bbar \sim B$, and $\bar{a}\geq \bbar$.
If $A$ dominates $B$ from above let \newsym{phi-AB}{Stochastic domination map from distribution $A$ to $B$ where $A \gtrsim B$}{\phi_{A,B}} be the (random) function from $(\integers_{\geq 0})^j$ to $(\integers_{\geq 0})^j$ which maps $\bar{a}$ to $\bbar$ with probability $C(\bar{a},\bbar)/A(\bar{a})$. By construction, $\phi_{A,B}(\bar{a}) \leq \bar{a}$ pointwise. Also, when $\bar{a} \sim A$, the random variable $\phi_{A,B}(\bar{a})$  has distribution~$B$.   
We naturally use the notation 
\newsym{phi-AB-minus}{}{\phi^{-1}_{A,B}} for the (random) function from $(\integers_{\geq 0})^j$ to $(\integers_{\geq 0})^j$ which maps $\bbar$ to $\bar{a}$ with probability $C(\bar{a},\bbar)/B(\bar{b})$.
\end{definition}

\begin{observation}\label{obs:dominate-suff-cond}
Let $A$ and $B$ be  
distributions on $(\integers_{\geq 0})^j$.
Let $\avect$ be drawn from~$A$ and $\bvect$ from~$B$. 
Suppose that, for 
all $\xbar \in
(\integers_{\geq 0})^j$,   $\pr(\avect \ge \xbar) \ge \pr(\bvect \ge \xbar)$. Then $A\gtrsim B$.    
\end{observation}

\begin{observation} \label{obs:dom-Poisson-tuples}
If $j$-tuples $\avect$ and $\bvect$ satisfy $\avect \geq \bvect$ then
$\Po(\avect) \gtrsim \Po(\bvect)$.
\end{observation}

\begin{lemma}\label{lem:Harris}
Let $\rho_1,\rho_2\colon \reals^j\to\reals$ be non-decreasing functions. Let $X_1,\dots,X_j$ be independent real-valued random variables and define the random vector $\Xvect = (X_1,\dots,X_j)$ taking values in $\reals^j$. Then $\E(\rho_1(\Xvect)\rho_2(\Xvect)) \ge \E(\rho_1(\Xvect))\E(\rho_2(\Xvect))$.
\end{lemma}
\begin{proof}
This is Harris's inequality, copied from~\cite[Theorem~2.15]{harris-textbook} -- with renaming.
\end{proof} 
\begin{corollary}\label{cor:Harris}
Let $X_1,\dots,X_j$ be independent real-valued random variables and define the random vector $\Xvect = (X_1,\dots,X_j)$ taking values in $\reals^j$. 
Let $\calE_1,\calE_2$ be increasing events that can be determined from the value of~$\Xvect$ with  $\Pr(\calE_2)>0$.
Then $\Pr( \calE_1 \mid \calE_2) \geq \Pr(\calE_1)$.
\end{corollary}
\begin{proof}
 Let $\rho_1$ be the indicator  r.v.{} for~$\calE_1$ and let $\rho_2$ be the indicator for~$\calE_2$. Then apply Lemma~\ref{lem:Harris} and divide through by $\Pr(\rho_2=1)$.
\end{proof} 

\subsection{Tail bounds}\label{sec:tail-bounds}

\begin{lemma}\label{lem:chernoff-large-dev} (E.g.\ \cite[Theorem~4.4(iii) and Theorem~5.4(i)]{MU})
    Let $Y$ be a real random variable with mean $\mu$ which is either Poisson or a sum of independent Bernoulli variables. Then for all $R \ge 6\mu$, we have $\pr(Y \ge R) \le 2^{-R}$.
\end{lemma}

\begin{lemma}\label{lem:chernoff-small-dev-upper}(E.g.\ \cite[Theorems~4.4(ii) and 5.4(i)]{MU})
Let $Y$ be a real random variable with mean $\mu$ which is either Poisson or a sum of independent Bernoulli variables. Then for all $0 < \delta \le 1$, we have $\pr(Y \ge (1+\delta)\mu) \le e^{-\delta^2\mu/3}$.
\end{lemma}

\begin{lemma}\label{lem:chernoff-small-dev} (E.g.\ \cite[Theorems~4.5(ii) and 5.4(ii)]{MU})
Let $Y$ be a real random variable with mean $\mu$ which is either Poisson or a sum of independent Bernoulli variables. Then for all $0 < \delta < 1$, we have $\pr(Y \le (1-\delta)\mu) \le e^{-\delta^2\mu/2}$.
\end{lemma}

\begin{lemma}\label{lem:poisson-noisy}
Fix a send sequence $\p$. Fix real numbers $c>0$ and $W\geq 1$. Let $S$ be a set of positive integers such that, for every $k\in S$, $W_k \geq W$. Let $\{ P_k : k\in S\}$ be a set of   independent Poisson variables 
such that, for each $k\in S$,
$\E(P_k) \ge cW_k$.   Then for all $\delta \in (0,1)$,
\[
\pr\big(\sum_{x \in S} P_xp_x < (1-\delta)c|S|\big) \le \exp(-\delta^2 c |S| W/2).
\]
\end{lemma}
\begin{proof}
We will apply the standard proof of a Chernoff bound with an added convexity argument. Let 
$t > 0$ (to be fixed later), let $X = \sum_{x \in S}P_x p_x$, and let $\delta \in (0,1)$. Then by Markov's inequality, since $\E(X) \ge c|S|$,
\begin{align*}
\pr(X < (1-\delta)c|S|) &= \pr(e^{-t X} > e^{-t(1-\delta)c|S|})
\leq \frac{\E[e^{-tX}]}{e^{-t(1-\delta)c|S|}} 
= \frac{\prod_{x \in S} \E[e^{-t p_xP_x}]}{e^{-t(1-\delta)c|S|}}\\
&\leq\frac{\prod_{x \in S} 
\exp(c W_x (e^{-t p_x} - 1))} {e^{-t(1-\delta)c|S|}}
=\exp\big( \sum_{x \in S} c W_x (e^{-t p_x} - 1)) - t|S| c(1-\delta)\big).
\end{align*}
Define a function $f$ by $f(z) = e^{-tz/W}$; then $f$ is convex, so for all $a \in [0,1]$, $f(a) \le (1-a)f(0) + af(1)$. On taking $a = Wp_x$ (which is in $[0,1]$ since $x \in S$), we obtain
\[
e^{-tp_x} \le 1 - Wp_x + Wp_xe^{-t/W} = 1 + Wp_x(e^{-t/W}-1);
\]
hence $cW_x(e^{-tp_x} - 1) \le cW(e^{-t/W}-1)$ and so
\[
\pr(X < (1-\delta)c|S|) \le \exp\big( cW|S|(e^{-t/W}-1) - t|S| c(1-\delta)\big).
\]
Now take $t = -W\log(1-\delta)>0$ to get
\begin{align*}
\pr(X < (1-\delta)c|S|) &\le 
\exp( -cW|S|\delta + W \log(1-\delta)|S| c(1-\delta)) \\
&= 
\exp( - c |S| W (\delta -  \log(1-\delta)(1-\delta))).
\end{align*}

As in the standard Chernoff proof, for all $0 < \delta < 1$, 
$\delta^2/2 \leq \delta -  \log(1-\delta)(1-\delta)$,  so the result follows.
\end{proof}

\begin{corollary}\label{cor:poisson-noisy}
Fix a real number $\lambda>0$.
Fix a send sequence $\p$. Fix a real number $W\geq 1$.  Let $S$ be a set of positive integers such that, for every $k\in S$, $W_k \geq W$. Let $\{ P_k : k\in S\}$ be a set of   independent Poisson variables 
such that, for each $k\in S$,
$\E(P_k) \ge  \lambda W_k/10$.   Then \[
\pr\big(\sum_{x \in S} P_xp_x < \lambda |S|/20\big) \le \exp(-\lambda |S| W /80).
    \]
\end{corollary}
\begin{proof}
Take $c=\lambda/10$ and $\delta=1/2$ in Lemma~\ref{lem:poisson-noisy}.
\end{proof}

\begin{lemma}\label{lem:janson-geo} \cite[Theorem 2.1]{Janson}
Let $X_1,\dots,X_z$ be independent geometric variables where $X_i$ has parameter $p_i$, and let $X = X_1 + \dots + X_z$. Let $\mu = \sum_{i=1}^z 1/p_i$ and let $p_* = \min_i p_i$. Then for all $C \ge 1$,
\[
\pr(X \ge C\mu) \le e^{-p_*\mu(C-1-\log C)}.
\]
\end{lemma}

\begin{definition}\label{def:weak-neg-assoc}
    We say that a family of random variables $Z_1, \dots, Z_n$ is \emph{negatively associated} if for all tuples $A_1$ and $A_2$ whose elements lie in $\{Z_1,\dots,Z_n\}$ and are mutually distinct, and all real non-decreasing functions $f_1,f_2$ on those tuples,
    \begin{equation*}
        \E\big(f_1(A_1)f_2(A_2)\big) \le \E(f_1(A_1))\E(f_2(A_2)).
    \end{equation*}
\end{definition}

\begin{lemma}\label{lem:coord-neg-assoc}
    Suppose $Z_1,\dots,Z_n$ are negatively-associated random variables. Then for all real non-decreasing functions $f_1,\dots,f_n$,
    \[
        \E\Big(\prod_{i=1}^n f_i(Z_i)\Big) \le \prod_{i=1}^n \E(f_i(Z_i)).
    \]
\end{lemma}
\begin{proof}
    This is proved as~\cite[Lemma~4]{dr-na}, with our $Z_i$'s mapping to their $X_i$'s.
\end{proof}

We next define a generalised multinomial distribution referred to in~\cite[p.~100]{dr-na} as ``the balls and bins experiment''. Essentially, it is the distribution obtained by relaxing the typical requirement of a multinomial distribution that all independent trials have the same distribution.

\begin{definition}\label{def:gen-multi}
Let $m$ and $n$ be positive integers, and let $A_1,\dots,A_m$ be independent integer random variables, each of which has its
support contained in $[n]$. 
For every $i\in [m]$ and $k\in[n]$,
let $Q_{i,k} = \pr(A_i = k)$.
For all $k \in [n]$, let 
$a_k = |\{i \in [m] \colon A_i = k\}|$. Then $\avect=(a_1,\ldots,a_n)$ is a   generalised multinomial  random variable with parameter~$Q$.
\end{definition}

\begin{lemma}\label{lem:multinomial-neg-assoc}
If $Z_1,\dots,Z_n$ are components of a  generalised multinomial random variable, then $Z_1,\ldots,Z_n$ are negatively associated.
\end{lemma}
\begin{proof}
It is immediate from Definition~\ref{def:weak-neg-assoc} that it suffices to prove that the set of all components of a generalised multinomial random variable $(Z_1',\dots,Z_{n'}')$ is negatively associated (even if it has other components besides $Z_1,\dots,Z_n$). This is proved as~\cite[Theorem~14]{dr-na}, with our $Z_i'$'s mapping to their $B_i$'s and our $n'$ mapping to their $n$.
\end{proof}

\begin{lemma}
\label{lem:excess-chernoff}
Let $n$ be an integer, and let $Z_1,\dots,Z_n$ be components of a generalised multinomial  random variable. Let $\muvect \in (\reals_{\ge 0})^n$, and suppose $\E(Z_k) \le \mu_k$ for all $k \in [n]$. Fix $\delta \ge e$, and for all $k \in [n]$, let $Z_k^+ = \max\{0, Z_k - \delta\mu_k\}$. Let $Z^+ = \sum_{k=1}^n Z_k^+$. Then for all $z\geq 0$, $\pr(Z^+ > z) \le 2^n\delta^{-z}$.
\end{lemma}

\begin{proof}
The proof is based on the standard Chernoff bound method as used in e.g.~\cite{MU}.
Fix $z\geq 0$ and let $\theta = \log \delta$; observe that $\theta > 0$ since $\delta \ge e$. Then by Markov's inequality,
\begin{align*}
\pr(Z^+ > z) &= \pr(e^{\theta Z^+}\geq e^{\theta z})
\leq e^{-\theta z} \E(e^{\theta Z^+}) = e^{-\theta z}\E\Big(\prod_{k=1}^ne^{\theta Z_k^+}\Big).
\end{align*}
By Lemma~\ref{lem:multinomial-neg-assoc}, $Z_1,\dots,Z_n$ are negatively associated. Thus by applying Lemma~\ref{lem:coord-neg-assoc} and taking $f_k(Z_k) = e^{\theta Z_k^+}$, it follows that
\[
\pr(Z^+ > z) \le e^{-\theta z} \prod_{k=1}^n \E(e^{\theta Z_k^+}).
\]
For all 
real numbers~$x$,  $e^{\theta \max\{0,x\}} 
  \le 1 + e^{\theta x}$. 
Thus, $e^{\theta Z_k^{+}} \leq 1 + e^{\theta (Z_k - \delta \mu_k)}$.
Hence
\begin{equation}\label{eq:excess-chernoff-0}
\pr(Z^+ > z) \le e^{-\theta z}\prod_{k=1}^n\big(1+\E(e^{\theta (Z_k-\delta\mu_k)})\big).
\end{equation}

Without loss of generality, suppose that $Z_1,\dots,Z_n$ are the first $n$ coordinates of a generalised multinomial random variable with parameter~$Q$.
Let $A_1,\ldots,A_m$ be the random variables in the definition of a generalised multinomial random variable.
For $i\in [m]$ let $I_{i,k}$ be the indicator for the event~$A_i=k$. 
Thus,
for every $i\in [m]$ and $k\in[n]$,
$Q_{i,k} = \pr(A_i = k) = \E[I_{i,k}]$.
Also, $Z_k = \sum_{i\in [m]} I_{i,k}$
and the variables in $\{I_{1,k},\ldots,I_{m,k}\}$ are independent.
By~\eqref{eq:excess-chernoff-0}, it follows that
\begin{equation}\label{eq:excess-chernoff-1}
\pr(Z^+ > z)
= e^{-\theta z}\prod_{k=1}^n\bigg(1 + e^{-\theta\delta \mu_k}\E\Big(\prod_{i=1}^{m} e^{\theta I_{i,k}}\Big)\bigg),
\end{equation}
where
\[
    \E\Big(\prod_{i=1}^m e^{\theta I_{i,k}}\Big) = \prod_{i=1}^m \E(e^{\theta I_{i,k}}) = \prod_{i=1}^m \big(1 + Q_{i,k}(e^\theta-1)\big).
\]
Moreover, for all $i$ and $k$, $1+Q_{i,k}(e^\theta - 1) \le e^{Q_{i,k}(e^\theta - 1)}$, so
\[
\E\Big(\prod_{i=1}^m e^{\theta I_{i,k}}\Big) \le \exp\Big(\sum_{i=1}^m Q_{i,k}(e^\theta-1)\Big) = \exp\big((e^\theta-1)\E(Z_k)\big) \le \exp\big((e^\theta-1)\mu_k\big).
\]
Hence by~\eqref{eq:excess-chernoff-1},
\[
\pr(Z^+ > z) \le e^{-\theta z}\prod_{k=1}^n\Big(1+\exp\big(\mu_k(-\theta\delta+e^\theta-1)\big)\Big).
\]
Since $\theta = \log \delta$ and $\delta \ge e$ by hypothesis, $e^\theta - 1 - \theta\delta < \delta - \delta\log\delta \le 0$. Hence
\[
\pr(Z^+ > z) \le e^{-\theta z}(1+1)^n = 2^n\delta^{-z}.
\]
\end{proof}

\subsection{The backoff Markov chain}

Theorem~\ref{thm:goal} shows that for any send sequence~$\p$ and any positive birth rate~$\lambda$, the backoff process with send sequence~$\p$ and birth rate~$\lambda$ is not positive recurrent.
Specifically, the proof of Theorem~\ref{thm:goal} shows that the empty state, with no balls waiting, is not positive recurrent.
Except in a pathological case which is easily handled,
this implies that no states are positive recurrent, since Proposition~\ref{prop:irreducible-aperiodic} shows that the Markov chain is irreducible and aperiodic. 

\begin{definition}
We say that a send sequence~$\p$ \emph{has an infinite tail of $1$s} if 
there is a positive integer $k$ such that, for all $j\geq k$, $p_j=1$.
\end{definition}

\begin{proposition}
\label{prop:irreducible-aperiodic} Let $\p$ be a send sequence that does not have an infinite tail of $1$s. Fix any birth rate $\lambda >0$. Then the backoff process with send sequence~$\p$ and birth rate~$\lambda$ is irreducible and aperiodic.
\end{proposition}
\begin{proof}
Let $X$ be the backoff process with send sequence~$\p$ and birth rate~$\lambda$.
Let $\Omega$ be the set of all states that are reachable from the initial empty state of~$X$. That is, $\Omega$ contains all states~$\omega$ 
such that there exists a $t\geq 0$ such that $\Pr(\bbbar^X(t) = \omega) > 0$.

First, from every state $\bbbar^X(t)$, it is possible to reach the empty state. To see this,
let $\beta_1,\ldots,\beta_K$ be the balls in $\bbbar^X(t)$.
For each $k\in [K]$, let $J[k]$ be the bin such that
$\beta_k \in \b^X_{J[k]}(t)$.
Let $J'[k] = \min \{j' \geq J[k]: p_{j'} < 1\}$ and let
$\ell = \max_{k\in [K]} J'[k] - J[k]$.
Consider the evolution from $\bbbar^X(t)$ where
there are no births 
during steps $t + 1,\ldots, t+\ell+K$. 
While it is still in the system, ball $\beta_k$
sends at steps $t+1,\ldots, t+ J'[k] - J[k]$ 
and then at step $t+\ell+k$ -- if it is still in the system then it is in bin~$J'[k]$ at time~$t+\ell$ and it is silent 
for all other steps in the range $t+\ell+1,\ldots, t+\ell + K$.  Thus, $\bbbar^X(t+\ell+K)$ is the empty state.
 
By definition, the empty state can reach every state in~$\Omega$, so the chain is irreducible.
It is aperiodic since the empty state can make a transition to itself.
\end{proof}

We have necessarily excluded send sequences with an infinite tail of $1$s from Proposition~\ref{prop:irreducible-aperiodic}. However, Corollary~\ref{KM-cor} (a corollary of the characterisation of Kelly and MacPhee) shows that, with probability~$1$, only finitely balls escape in  a backoff process with such a send sequence   so we do not need to treat such send sequences here.

The reader may wonder about the set~$\Omega$ in the
proof of Proposition~\ref{prop:irreducible-aperiodic}.  
The exact nature of~$\Omega$ depends on~$\p$.
For example, suppose that $p_j=1$ for some $j\geq 2$.
Then it is not possible
to reach a state such that $b_j^X(t)=1$ and, 
for all $i\neq j$, $b_j^X(t)=0$.
The single ball~$\beta$ in $\b_j^X(t)$ would have to be in $\b_{j-1}^X(t-1)$ since $p_j=1$. However, in this case it would only arrive in $\b^X_j(t)$ if at least two balls send 
at time~$t$, which would mean $|\cup_{i} b_i^X(t)| \geq 2$.
In any case, we don't need to know the exact nature of the state space~$\Omega$, 
the proof simply shows that the Markov chain~$X$ is irreducible and aperiodic with state space~$\Omega$.

\section{Generalised backoff processes}
\label{sec:backoff-defs}

In order to prove Theorem~\ref{thm:goal} we will couple a backoff process with many other stochastic processes, and analyse these processes in detail. We start by defining the notion of a \emph{generalised backoff process}. A backoff process is an example of a generalised backoff process. However, a generalised backoff process may differ from a backoff process in many ways: Instead of starting from time~$0$ with the first time step being step~$1$, it may start from an arbitrary time~$\tau$ with the first time step being step~$\tau+1$. Instead of starting from an empty state, it may start with an arbitrary initial population.  In a backoff process, the number of births just before each send-step is a Poisson random variable with mean~$\lambda$, but in general, the number may follow any distribution. Also, in a generalised backoff process, there is the possibility of end-of-step arrivals (new balls added to the bins) at the end of each step.  The decisions that balls make about whether or not to send in a generalised backoff process is determined by a send sequence, in the same way as in a backoff process. However, in a generalised backoff process there is more flexibility about which balls escape. Definition~\ref{def:genbackoff}
therefore leaves many details of a generalised backoff unspecified and all of these need to be specified to define a particular generalised backoff process (see Remark~\ref{rem:genbackoff}).

A backoff process is of course a special case of a generalised backoff process (see Definition~\ref{def:backoffasgen}). Most of the generalised backoff processes that we consider will not have end-of-step arrivals -- certainly backoff processes do not have end-of-step arrivals. The reason that we include the possibility of end-of-step arrivals in generalised backoff processes is that they are useful in escape processes. Recall from the Introduction (Section~\ref{sec:intro}) that we will couple a backoff process with a ``volume process'' that tracks its send patterns and an infinite number of ``escape processes'' that track escapes. The coupling will maintain the invariant that every ball present in the volume process but not in an escape process is also present in the original backoff process. Thus, the end-of-step arrivals (into escape processes) are designed to track escapes from the original backoff process, to enable this invariant. 

The ``names'' of balls are unimportant but it is important that balls have distinct names, so we use the following convention in Definition~\ref{def:genbackoff}. Let $\newsym{B-cal}{ball names}{\ball(t)} = \{
\ball(t',i,i') \mid \mbox{$t'$, $i$, $i'$ are non-negative integers with $t'\leq t$}
\}$.

\begin{definition}\label{def:genbackoff} 
Fix a \emph{send sequence}  
$\p$, a non-negative integer~$\tau$, and 
a distribution~\newsym{D-cal-0}{birth distribution}{\calD} on non-negative integers. 
Fix a (possibly random)
tuple $\arrvect^Z(\tau)$ of disjoint sets 
$\arrvect^Z(\tau) = (\arr^Z_{1}(\tau),\arr^Z_2(\tau),\ldots)$ where
each $\arr_i^Z(\tau) 
\subseteq \ball(\tau)$.
A \emph{generalised backoff process}~\newsym{Z-roman}{generalised backoff process, Definition~\ref{def:genbackoff}}{Z} with send sequence~$\p$, start time~$\tau$,
birth distribution~$\calD$, and 
initial population $\arrvect^Z(\tau)$ is defined as follows: First,
$\b^Z_0(\tau) = \emptyset$ and,
for all positive integers~$j$,
$\b^Z_j(\tau) = \arr^Z_j(\tau)$.
For each time step~$t>\tau$ and each bin $j\geq 0$, the process has random variables 
\newsym{Bjz-bf}{population of bin~$j$ in~$Z$ just before sends at $t$}{\B^Z_j(t)}, 
\newsym{sjz-bf}{balls that send from $j$ at $t$ in $Z$}{\s^Z_j(t)},
\newsym{ejz-bf}{balls that escape from $j$ at $t$ in $Z$}{\e^Z_j(t)}, and 
\newsym{bjz-bf}{population of bin~$j$ in~$Z$ just after step $t$}{\b^Z_j(t)} of respective sizes 
\newsym{Bjz}{$|\B^Z_j(t)|$}{B^Z_j(t)}, 
\newsym{sjz}{$|\s^Z_j(t))|$}{s^Z_j(t)}, 
\newsym{ejz}{$|\e^Z_j(t)|$}{e^Z_j(t)}, and 
\newsym{bjz}{$|\b^Z_j(t)|$}{b^Z_j(t)}.
By analogy to backoff processes, 
$\b^Z_j(t)$ is the population of bin~$j$ just after step~$t$ in~$Z$, 
$\B^Z_j(t)$ is the population of bin~$j$ just before the sends at step~$t$,
$\s^Z_j(t)$ is the set of balls that send from bin~$j$ at step~$t$, and $\e^Z_j(t)$ is the set of balls that escape from bin~$j$ at step~$t$.
The process evolves as   follows for $t>\tau$.
\begin{itemize}

\item   The quantity \newsym{n-z}{$|n^Z(t)|$}{n^Z(t)} is drawn from~$\calD$.  If $n^Z(t) = 0$ then $\n^Z(t) = \emptyset$. Otherwise,
the set 
$\newsym{nz-bf}{set of balls born at $t$ in $Z$}{\n^Z(t)} = \{\ball(t,0,1),\ldots,\ball(t,0,n^Z(t))\}$.

\item $\B_0^Z(t) = \n^Z(t) \cup \b_0^Z(t-1)$.   
For $j\geq 1$, $\B_j^Z(t) =  \b_j^Z(t-1)$.

\item    
For all $j\geq 0$,  the set $\s^Z_j(t)$ is formed
by taking each element of $\B_j^Z(t)$ independently with probability~$p_j$. 
$\newsym{sz-bf}{set of balls that send at $t$ in $Z$}{\s^Z(t)} = \cup_{j\geq 0}\ 
\s^Z_j(t)$
and $\newsym{sz}{$|\s^Z(t)|$}{s^Z(t)} = |\s^Z(t)|$.

\item For all $j\geq 0$, 
a random subset $\e^Z_j(t)\subseteq \s^Z_j(t)$ is drawn from some distribution (to fix a particular generalised backoff process this distribution will be specified).

\item $\b^Z_0(t) = \B^Z_0(t) \setminus \s^Z_0(t)$.

\item Fix a (possibly random) tuple 
$\arrvect^Z(t)$ of disjoint sets
$\newsym{arr-bf}{end-of-step arrivals}{\arrvect^Z(t)}  = ( \arr^Z_{1}(t),\arr^Z_2(t),\ldots)$
where each $\arr^Z_i(t) 
\subseteq \ball(t)
\setminus  \cup_{k\geq 0} \B^Z_k(t)$.

\item For all $j\geq 1$, 
$\b^Z_j(t) = \arr^Z_j(t) \cup (\B^Z_j(t) \setminus \s^Z_j(t)) \cup 
(\s^Z_{j-1}(t) \setminus \e^Z_{j-1}(t))$.

\end{itemize}

For $t\geq \tau$ we use the notation $\bbbar^Z(t)$ for the  tuple $(\b_0^Z(t),\b_1^Z(t),\ldots)$
and we use $\bbar^Z(t)$ for the tuple 
 $(b_0^Z(t),b_1^Z(t),\ldots)$. We define $\romanarrvect^Z(t) = (\romanarr_1^Z(t), \romanarr_2^Z(t),\ldots)$ where $\newsym{arr-roman}{$|\arr^Z_j(t)|$}{\romanarr^Z_j(t)} = |\arr^Z_j(t)|$.

\end{definition}

\begin{remark}\label{rem:genbackoff}
To specify a specific generalised backoff process~$Z$ it is necessary to specify $\p$, $\tau$, $\calD$, $\arrvect^Z(t)$ for $t\geq \tau$, and the mechanism for choosing $\e_j^Z(t)$.
\end{remark}

Given a specific send sequence~$\p$ and birth rate~$\lambda>0$ the backoff process with send sequence~$\p$ and birth rate~$\lambda$ 
(Definition~\ref{def:backoff})
can be formalised as a generalised backoff process.
In Definition~\ref{def:backoffasgen} we do this, but 
instead of fixing the birth distribution~$\calD$ to be $\Po(\lambda)$, we allow
 for a more general birth distribution -- this will be useful for domination in proofs. When we speak of a backoff process with birth rate~$\lambda$
 we mean that the birth distribution is $\Po(\lambda)$.

\begin{definition}
\label{def:backoffasgen}
Fix a send sequence~$\p$ and a 
a distribution $\calD$ on non-negative integers. The \emph{backoff process}  $X$
with send sequence~$\p$ and birth distribution~$\calD$ 
is 
the generalised backoff process with send sequence~$\p$, start time~$0$,  and birth distribution $\calD$, where for all $t\geq 0$ and $j\geq 1$,
$\arr^X_j(t) = \emptyset$ and  where, for all $t\geq 1$ and $j\geq 0$,
\begin{itemize}
    \item 
$\e^X_j(t) =\s^X_j(t)$ if 
$s_j^X(t) = s^X(t) = 1$
and equal to $\emptyset$, otherwise.
\end{itemize}
\end{definition} 

Sometimes in the paper it will be useful to divide the balls of the backoff process into two or three cohorts. 
Cohorts are different from the streams of the two-stream process  that we discussed in the introduction (which will be defined formally in Definition~\ref{def:two-stream}) -- the cohorts of a backoff process do not affect the way that it evolves -- they are purely for accounting. Definition~\ref{def:backoff-cohorts} is phrased in terms of a general set~$\calC$ of cohorts, but the two cases that we will use are $\calC = \{A,B\}$ with two cohorts and $\calC = \{A,B,C\}$ with three cohorts.

\begin{definition}
\label{def:backoff-cohorts} 
Fix a send sequence~$\p$ and 
a distribution $\calD$ on non-negative integers. 
Let $X$ be the {backoff process}  
with send sequence~$\p$ and birth distribution~$\calD$.
Let $\calC$ be a finite set of cohorts. The $\calC$-cohorts of~$X$  are defined as follows.
For every $t\geq 1$, partition the set  $\n^{X}(t)$ of balls born at time $t$ 
by choosing a cohort~$D\in \calC$ uniformly at random for each ball,
and letting $\n^{X^D}(t)$ be the balls in~$\n^X(t)$ that choose cohort~$D$.
 The cohorts have no effect on the evolution of the backoff process -- they will be  useful only for accounting in proofs.
For each $t\geq 1$, each $j\geq 0$, and each $D\in \calC$, define
\begin{align*}
\newsym{b-j-X-D}{balls in bin~$j$ after step~$t$ in cohort~$D$, Def~\ref{def:backoff-cohorts}}{\b_j^{X^D}(t)} &= \b_j^{X}(t) \cap \cup_{t'=1}^{t} \n^{X^D}(t'),\\
\newsym{B-j-X-D}{balls in bin~$j$ before sends at~$t$ in cohort~$D$, Def~\ref{def:backoff-cohorts}}{\B_j^{X^D}(t)} &= \B_j^{X}(t) \cap \cup_{t'=1}^{t} \n^{X^D}(t'),\\
\newsym{s-j-X-D}{sends from bin~$j$ at step~$t$ in cohort~$D$, Def~\ref{def:backoff-cohorts}}{\s_j^{X^D}(t)} &= \s_j^{X}(t) \cap \cup_{t'=1}^{t} \n^{X^D}(t'),\\
\newsym{e-j-X-D}{escapes from bin~$j$ at step~$t$ in cohort~$D$, Def~\ref{def:backoff-cohorts}}{\e_j^{X^D}(t)} &= \e_j^{X}(t) \cap \cup_{t'=1}^{t} \n^{X^D}(t').
\end{align*}
The respective sizes of these are $b_j^{X^D}(t)$, $B_j^{X^D}(t)$, 
$s_j^{X^D}(t)$, and $e_j^{X^D}(t)$.
\end{definition}

Definition~\ref{def:standard-couple} defines a coupling -- the ``standard coupling'' -- between two backoff processes. The definition is somewhat lengthy because of the need to specify the invariants that are maintained as the coupling evolves -- the invariants at step $t-1$ are what guarantees the coupling to be well-defined at step $t$. Thus, it is most convenient to prove the invariants as part of the definition. The idea is straightforward, and will be used throughout the paper.

\begin{definition}\label{def:standard-couple}
Fix a send sequence~$\p$ and two distributions $\calD$ and~$\calD'$ on non-negative integers such that~$\calD \gtrsim \calD'$.
Let $X$ be the backoff process with 
send sequence~$\p$ and
birth distribution~$\calD$ 
and let $X'$ be the backoff process with send sequence~$\p$ and birth distribution~$\calD'$.
The \emph{\newsym{standard coupling}{standard coupling of two backoff processes}{\textrm{standard coupling}}} of~$X'$ and~$X$ has the following invariants for every $i\geq 0$ and $t\geq 0$.
\begin{description}  
\item  [Inv-$b_i(t)$:] $\b_i^{X'}(t) \subseteq \b_i^{X}(t)$.
\item  [Inv-$B_i(t)$:] $\B_i^{X'}(t) \subseteq \B_i^{X}(t)$.
\item [Inv-$s_i(t)$:] $\s_i^{X'}(t) \subseteq \s_i^{X}(t)$.
\end{description}
The coupling is defined as follows.
For each $i\geq 0$, $\b_i^{X'}(0) = \b_i^{X}(0) = \emptyset$, establishing Inv-$b_i(0)$.
 
For each positive integer~$t$, 
the coupling evolves as follows.
First,
the pair $(n^X(t),n^{X'}(t))$ is chosen from the joint distribution~$C$ from Definition~\ref{def:phi}, so that $n^X(t) \sim \calD$, $n^{X'}(t) \sim \calD'$ and $n^X(t) \geq n^{X'}(t)$.
Then for each $Z\in\{X',X\}$, the set $\n^Z(t)$
is defined according  to Definition~\ref{def:backoffasgen}. In particular, if
 $n^{Z}(t) = 0$ then $\n^{Z}(t) = \emptyset$. Otherwise, ${\n^Z(t)} = \{\ball(t,0,1),\ldots,\ball(t,0,n^Z(t))\}$.
For all $i\geq 0$, the sets
$\B^{X'}_i(t)$, $\s^{X'}_i(t)$,  
$\e^{X'}_i(t)$, $\b^{X'}_i(t)$, and $B^X_i(t)$ are chosen according to Definition~\ref{def:backoffasgen}.

Using Inv-$b_0(t-1)$, we get
$\B_0^{X'}(t) = \n^{X'}(t) \cup \b_0^{X'}(t-1) \subseteq
\n^{X}(t) \cup \b_0^{X}(t-1) = \B_0^{X}(t)$, 
establishing Inv-$B_0(t)$.
Also, 
for $i\geq 1$ 
using Inv-$b_i(t-1)$, we get $\B_i^{X'}(t) = \b_i^{X'}(t-1) \subseteq \b_i^{X}(t-1) = \B_i^{X}(t)$, establishing Inv-$B_i(t)$

For every   $i\geq 0$, $\S^{X}_i(t)$ is formed by including each ball in 
$\B_i^{X}(t) \setminus \B_i^{X'}(t)$ 
independently with probability~$p_i$.
Then $\s^{X}_i(t) = \S^{X}_i(t) \cup \s^{X'}_i(t)$, establishing Inv-$s_i(t)$.

Finally, for every $i\geq 0$, $\e_i^{X}(t)$ and $\b_i^{X}(t)$ are chosen according to Definition~\ref{def:backoffasgen}.
The invariant Inv-$b_i(t)$ is established for all $i\geq 0$ as follows.  Fix $i$ and consider   a ball $\beta \in \b_i^{X'}(t)$. We wish to show $\beta \in \b_i^{X}(t)$.
There are two possibilities.
\begin{itemize}
\item   $\beta \in \B_i^{X'}(t) \setminus \s_i^{X'}(t)$:  
In this case, by Inv-$B_i(t)$, $\beta \in \B_i^{X}(t)$. By construction,   $\beta\notin \s_i^{X}(t)$ so $\beta \in \b_i^{X}(t)$.
\item (for $i\geq 1$) $\beta \in \s^{X'}_{i-1}(t) \setminus \e^{X'}_{i-1}(t)$:
From Definition~\ref{def:backoffasgen}, $s^{X'}(t) \geq 2$, By
Inv-$s_{i-1}(t)$, $\beta \in \s^{X}_{i-1}(t)$. By all of the constraints Inv-$s_k(t)$, $s^{X}(t) \geq s^{X'}(t) \geq 2$ so $e^{X}_{i-1}(t) = \emptyset$ and
$\beta \in \s^{X}_{i-1}(t) \setminus \e^{X}_{i-1}(t)$ which implies
$\beta \in \b_i^{X}(t)$.
\end{itemize}
\end{definition}

Definition~\ref{def:expected-noise} extends the notion of expected noise from the introduction to generalised backoff processes. 

\begin{definition}\label{def:expected-noise}
Let $Z$ be a generalised backoff process with send sequence $\p$ and start time $\tau$. For all $t \ge \tau+1$ and all sets $S$ of positive integers, let 
$\newsym{N-S}{expected noise, Def~\ref{def:expected-noise}}{\calN^Z_S(t)} = \sum_{i \in S}p_ib_i^Z(t-1)$.
\end{definition}

Observe that $\calN^Z_S(t)$ is the  expected number of sends at time $t$, conditioned on the 
(random) values in $\{ b^Z_i(t-1): i\geq 1\}$.

We conclude with Lemma~\ref{lem:single-send-bound}, which gives an upper bound on the probability that a step has at most one send 
(from non-newborns)
in a generalised backoff process.

\begin{lemma}\label{lem:single-send-bound}
Fix a \emph{send sequence}  
$\p$ and a non-negative integer~$\tau$.
Fix $\lambda \in (0,1/60)$.
Fix a positive integer~$j$ and a $j$-tuple
$\avect \in (\integers_{\geq 0})^j$.
Let $Z$ be
a generalised backoff process with  
send sequence~$\p$, start time~$\tau$,   and an initial population~$\arrvect^Z(\tau)$ such that $\romanarrvect^Y_{[j]}(\tau) = \avect$. 
Fix a  set $S \subseteq [j]$ such that 
$|S| \geq 100 /\lambda^2$ and 
$\calN^Z_{S}(\tau+1) = \sum_{i\in S} p_i a_i  \geq \lambda |S|/80$.
Let $\quiet$ be the probability that
$\sum_{i\geq 1} s_i^Z(\tau+1) \leq 1$. 
Then $\quiet \leq  
\exp(-\sum_{i \in S}p_i a_i/16)$.
\end{lemma}
\begin{proof} 
Let $\calN$ denote $\calN_S^Z(\tau+1) = \sum_{i\in S} p_i a_i$. 
By assumption, $\calN \geq \lambda |S|/80 \geq  
100/(80 \lambda) \geq 4$.
We split into two cases depending on 
$\avect_{[j]}$.

\medskip\noindent\textbf{Case 1:} There is an $i \in S$ such that $p_i a_i \geq \calN /2$. In this case, $p_i a_i \geq 2$.  
By a Chernoff bound (Lemma~\ref{lem:chernoff-small-dev} with $\delta=1/2$),
\begin{align*}
\quiet &\le \pr(s_i^Z(\tau+1) \le 1) \le \pr(s_i^Z(\tau+1) \le p_i a_i/2) \le e^{-a_ip_i/8} \le \exp(-\calN/16).
\end{align*}

\medskip\noindent\textbf{Case 2:} For all $i \in S$, $p_i a_i \leq \calN/2 $.
For this case note that we already established
$\calN/8 \geq 100/(640 \lambda) = (50/640)(2/\lambda)$, 
which is at least $\log(2/\lambda)$.
Let $R = \lambda |S|$. We
also have $\calN/8 \geq R/640$ which is at least $\log(R)$
for $R\geq 6000$ which holds since
$R \geq 100/\lambda$ and $1/\lambda \geq 60$.
Choose $i\in S$ to maximise 
$\prod_{x \in S\setminus \{i\}} (1-p_x)^{a_x} $.
Then 
 \begin{align*}
 \quiet &\leq \prod_{x\in S} (1-p_x)^{a_x} + \sum_{k\in S} \prod_{x \in S \setminus \{k\}} (1-p_x)^{a_x} 
 \le (|S|+1)\prod_{x \in S \setminus \{i\}} (1-p_x)^{a_x} \\&\le (|S|+1)\exp\big({-}\sum_{x \in S\setminus \{i\}}p_x a_x\big)
 \le (|S|+1)\exp\big({-}\sum_{x \in S}p_x a_x/2\big) = (|S|+1) \exp(-\calN/2).  
 \end{align*}
Now $|S|+1 \leq 2 |S| = R(2/\lambda)$
so $|S|+1 \leq \exp( \log(R) + \log(2/\lambda)) 
\leq \exp( \calN/8 + \calN/8) = \exp(\calN/4)$
so $q \leq \exp(-\calN/4)$. 
\end{proof}

\section{Reducing to the case where \texorpdfstring{$p_0=1$}{p\_0 = 1}
with upper bounds on~\texorpdfstring{$\lambda$}{lambda} and the bin weights \texorpdfstring{$W_i$}{W\_i} }  
\label{sec:reduction}

The purpose of this section is to state Theorem~\ref{thm:goaljkiller} (which will be proved in the rest of this paper) and to prove Theorem~\ref{thm:goal} using it.  
The statement of Theorem~\ref{thm:goaljkiller}
is similar to that of Theorem~\ref{thm:goal} except that it allows us to assume that $\lambda$ is sufficiently small, that $p_0=1$, and to assume a useful upper bound on the weight $W_i$ of each bin~$i$. (Recall that $W_i=1/p_i$.) The upper bound is inspired by Lemma~\ref{lem:killer-fact}, which is from \cite{GL-oldcontention} and was known to the authors of~\cite{GJKP}, but here we strengthen the bound by applying it to all bins rather than to all bins with sufficiently large index.

\begin{restatable} {theorem}{thmgoaljkiller}\label{thm:goaljkiller}
There exists a 
real number $\lambda_0 \in (0,1]$  
such that the following holds.
Consider any $\lambda \in (0,\lambda_0)$.
Consider any send sequence~$\p$ such that 
$p_0=1$ and, for all $j \ge 1$, $W_j \leq (2/\lambda)^j$. 
Then 
the backoff process with birth rate $\lambda$ and
send sequence $\p$ is not positive recurrent. 
\end{restatable}

We start by proving some observations and lemmas that we will use to prove Theorem~\ref{thm:goal} using Theorem~\ref{thm:goaljkiller}. First, Observation~\ref{obs:decrease-lambda} is useful for allowing the arbitrary choice of~$\lambda_0$ in Theorem~\ref{thm:goaljkiller}.

\begin{observation}\label{obs:decrease-lambda}
Fix positive reals $\lambda'$ and $\lambda$ 
such that $0 < \lambda' < \lambda$.
Fix a send sequence~$\p$.
Let $X$ be a backoff process with birth distribution~$\Po(\lambda)$ and send sequence~$\p$ and let $X'$ be a backoff process with birth distribution~$\Po(\lambda')$ and send sequence~$p$.
If $X'$ is not positive recurrent then $X$ is not positive recurrent.     
\end{observation}
\begin{proof}
The standard coupling of~$X'$ and~$X$
has the property that, for all $i\geq 0$ and $t\geq 0$, 
$\b_i^{X'}(t) \subseteq \b_i^{X}(t)$.
In the coupling, $X$ cannot return to the state in which every bin is empty before $X'$ does, so the expected return time for~$X$ is at least the expected return time for~$X'$.
\end{proof}

We next give Observation~\ref{obs:p0one} which 
is useful for allowing the restriction $p_0=1$ in 
Theorem~\ref{thm:goaljkiller}.

\begin{observation}\label{obs:p0one}
Fix $\lambda >0$. Fix a send sequence~$\p$.
Let $\p'$ be the send sequence that is identical to~$\p$ except that $p'_0=1$.
Let $X$ be the backoff process 
with send sequence~$\p$ and birth rate~$\lambda$.  Let $X'$ be the backoff process with
send sequence~$\p'$ and
birth rate $\lambda p_0$. If $X'$ is not positive recurrent then $X$ is not positive recurrent.
\end{observation}
\begin{proof}
We will give a coupling of~$X'$ and~$X$ with the property that, for all $i\geq 0$ and $t\geq 0$, 
$\b_i^{X'}(t) \subseteq \b_i^{X}(t)$.
The argument is then identical to the proof of Observation~\ref{obs:decrease-lambda}.

The coupling is essentially the standard coupling (Definition~\ref{def:standard-couple}), with 
$\calD = \Po(\lambda)$ and
$\calD' = \Po(\lambda')$ where $\lambda' = \lambda p_0$. 
The invariants are the same as those of the standard coupling. 
The coupling is also the same, 
except that the definition of $\s_0^X(t)$ 
is modified to get the correct marginals (this is necessary since   $p_0 < p'_0=1$).

Since $p'_0=1$, $\b_0^{X'}(t-1)=\emptyset$ and $\B_0^{X'}(t) = 
\n^{X'}(t)$.
Also, $\B_0^X(t) = \b_0^X(t-1) \cup \n^X(t)$.
In this coupling, 
$\S^X_0(t)$  
is formed by including each ball in $\b_0^X(t-1)$ independently with probability~$p_0$. This differs from the standard coupling in the sense that the balls in $\n^X(t) \setminus \n^{X'}(t)$ are all excluded from $\S^X_0(t)$. As in the standard coupling,  
then $\s^{X}_0(t) = \S^{X}_0(t) \cup \s^{X'}_0(t)$, establishing Inv-$s_0(t)$.

It is easy to see that the marginals are correct since the size of $\s^X_0(t) \cap \n^X(t)$ is $n^{X'}(t)$ which is a Poisson random variable with mean~$\lambda p_0$, as required. Furthermore, the size of $\n^X(t) \setminus \s_0(t)$ 
is Poisson with mean $\lambda - \lambda' = \lambda(1-p_0)$.
 \end{proof}

Our next goal is Corollary~\ref{cor:skip-j}, which says that, if we are trying to prove that a backoff process with 
send sequence $\p$ is not positive recurrent for any birth rate $\lambda\in (0,1)$, then we can freely remove any finite number of bins from the 
start of $\p$ and work with the new send sequence. The heart of the proof of Corollary~\ref{cor:skip-j} is Lemma~\ref{lem:dom1}, which we prove next.

\begin{lemma}\label{lem:dom1}
Fix $\lambda >0$ and a send sequence~$\p$.
Let $X$ be the backoff process with birth rate $\lambda$ and send sequence $\p$.
Let $J$ be a positive integer.
Let $\p'$ be the send sequence 
such that, for all $i\geq 0$, $p'_i = p_{i+J}$.
There is a real number $\lambda'\in (0,\lambda)$ and a distribution $\calD$ on non-negative integers with  $\calD \gtrsim \Po(\lambda')$ such that 
the following holds.
Let $X'$ be the backoff process with 
birth distribution $\calD$ and  send sequence $\p'$. Then 
there is a coupling of $X'$ and~$X$ such that, 
for all non-negative integers~$i$ and $t$,
$\b^{X'}_{i}(t) \subseteq \b^{X}_{i+J}(t+J)$.
\end{lemma}

\begin{proof} 
Let $\lambda' = \lambda^2 e^{-\lambda} \tfrac12 \prod_{\ell=0}^{J-1} p_\ell^2$.
We first define the birth distribution $\calD$. This will be defined in terms of the backoff process~$X$.
For   positive integers~$t$ and~$k$, let
$\mathcal{E}_{k,t+J}$ be the event that 
$n^X(t) = k+1$ (so
$\n^X(t) = \{\calB(t,0,1),\ldots,
\calB(t,0,k+1)\}$)  and
$\n^X(t) 
\subseteq \bigcap_{\ell \in \{0,\ldots,J-1\}}
\s^X_\ell(t+\ell) $. Note that $\Pr(\mathcal{E}_{k,t+J}) =  \frac{\lambda^{k+1}}{ e^{\lambda} (k+1)!} \prod_{\ell=0}^{J-1} p_\ell^{k+1}$, which does not depend on~$t$ and is
at most $\lambda^k / (e^\lambda k!)$.

We define $\calD$ to be the distribution which has value~$k$ with probability $\Pr(\mathcal{E}_{k,t+J}) $
and has value~$0$ otherwise. 
For every positive integer~$k$,
$\pr(\Po(\lambda')=k) $ is equal to
\begin{align*}
\frac{{\lambda'}^k}{ e^{\lambda'} k!} = 
\big(
\frac{\lambda^{k+1}}{ e^{\lambda} (k+1)!}
 \big)
\big(
  \prod_{\ell=0}^{J-1}  p_\ell^{k+1} \big)
\frac{k+1}{2^k} \lambda^{k-1} e^{-\lambda(k-1)}
  \big( \prod_{\ell=0}^{J-1} p_\ell^{k-1} \big)
\exp\big(-
\lambda^2 e^{-\lambda} \tfrac12 \prod_{\ell=0}^{J-1} p_\ell^2 
\big),
\end{align*}
which is at most $\Pr(\mathcal{E}_{k,t+J}) = \Pr(\calD=k)$. Thus, $\calD$ is dominated below by $\Po(\lambda')$, as required.

The coupling is similar to the standard coupling (Definition~\ref{def:standard-couple}) except that it is modified to account for the  
shifting of~$\p$ to create~$\p'$. We give the details for completeness.  
The invariants are as follows for every $i\geq 0$ and $t\geq 0$.
\begin{description}  
\item  [Inv-$b_i(t)$:] $\b_i^{X'}(t) \subseteq \b_{i+J}^X(t+J)$.
\item  [Inv-$B_i(t)$:] $\B_i^{X'}(t) \subseteq \B_{i+J}^X(t+J)$.
\item [Inv-$s_i(t)$:] $\s_i^{X'}(t) \subseteq \s_{i+J}^X(t+J)$.
\end{description}
For each $i\geq 0$, $\b_i^{X'}(0) = \b_i^{X}(0) = \emptyset$, establishing Inv-$b_i(0)$.
Before beginning step~$1$ of the coupling, the $X$ process is run
for the first $J-1$ steps (for $t\in \{1,\ldots,J-1\}$).

Then for each $t\geq 1$ the coupling evolves as follows.
We assume at this point that
$X'$ has been run through step $t-1$
and that $X$ has been run through step $t+J-1$
and that, for all $i\geq 0$, Inv-$b_i(t-1)$ holds.

We now show how  to coordinate step~$t$ of~$X'$ with step~$t+J$ of~$X$. 
The set $\n^{X}(t+J)$ is chosen independently from its distribution (its size is a Poisson random variable with parameter~$\lambda$).

The variables that are already defined include $\n^X(t)$
and, for all $\ell\in \{0,\ldots,J-1\}$, $\s_\ell^X(t+\ell)$.
Thus, it is clear whether there is a positive integer~$k$ such that $\calE_{k,t+J}$ occurs. 
\begin{itemize}
\item  If so, then we set
$\n^{X'}(t) = 
\{\calB(t,0,1),\ldots,
\calB(t,0,k)\}$. Note that all of these balls are in $\s_{J-1}^X(t+J-1)$ and 
$s_{J-1}^X(t+J-1)\geq k+1>1$ so all of them are in $\b_{J}^X(t+J-1)$ and they will also be in $\B_J^{X}(t+J)$.
\item If there is no such positive integer~$k$, then $\n^{X'}(t) = \emptyset$. 
\end{itemize}
Either way, we will have $\n^{X'}(t) \subseteq \B_J^X(t+J)$.

Then all of the sets $\B_i^{X'}(t)$ and $\B_i^{X}(t+J)$ are chosen from the definition.
Since $\B_0^{X'}(t) = \n^{X'}(t) \cup \b_0^{X'}(t-1)$ 
by definition 
and $\b_0^{X'}(t-1) \subseteq \b_{J}^X(t+J-1)$
by Inv-$b_0(t-1)$ and
$\b_J^X(t+J-1) = \B_J^X(t+J)$ by definition, we have
Inv-$B_0(t)$.
For any $i>0$, 
via Inv-$b_i(t-1)$ we get
$\B_i^{X'}(t) = \b_i^{X'}(t-1) \subseteq \b_{i+J}^{X}(t-1+J) = \B_{i+J}^X(t+J)$ so we get Inv-$B_i(t)$.

The rest of the coupling is similar to the standard coupling.  
For all $i\geq 0$, $\s_i^{X'}(t)$ 
is formed by taking each ball in
$\B_i^{X'}(t)$ independently with probability~$p'_i=p_{i+J}$.
Then $\S_{i+J}^{X}(t+J)$ is formed by taking each ball in $\B_{i+J}^X(t+J) \setminus \B_i^{X'}(t)$ independently with probability~$p_{i+J}$. Finally, $\s_{i+J}^{X}(t+J) = \s_i^{X'}(t) \cup \S_{i+J}^{X}(t+J)$, establishing Inv-$s_i(t)$. The sets $\s_0^X(t),\ldots,\s_{J-1}^X(t)$ are chosen following the correct marginals.

Finally, for every $i\geq 0$, $\e_i^X(t)$, $\e_i^{X'}(t)$, $\b_i^X(t)$ and $\b_i^{X'}(t)$ are chosen according to the definition.  

We finish by establishing 
Inv-$b_i(t)$ for all $i\geq 0$. 
Fix $i$  and consider a ball $\beta \in \b_i^{X'}(t)$. We wish to show $\beta \in \b_{i+J}^X(t+J)$.
There are two possibilities.
\begin{itemize}
\item   $\beta \in \B_i^{X'}(t) \setminus \s_i^{X'}(t)$:  
In this case, by Inv-$B_i(t)$, $\beta \in \B_{i+J}^X(t+J)$. By construction,   $\beta\notin \s_{i+J}^X(t+J)$ so $\beta \in \b_{i+J}^X(t+J)$.
\item (for $i\geq 1$) $\beta \in \s^{X'}_{i-1}(t) \setminus \e^{X'}_{i-1}(t)$:
From the definition, $s^{X'}(t) \geq 2$, By
Inv-$s_{i-1}(t)$, $\beta \in \s^X_{i+J-1}(t+J)$. By all of the constraints Inv-$s_k(t)$, $s^X(t+J) \geq s^{X'}(t) \geq 2$ so $\e^X_{i+J-1}(t+J) = \emptyset$ and
$\beta \in \s^{X}_{i+J-1}(t+J) \setminus \e^{X}_{i+J-1}(t+J)$ which implies
$\beta \in \b_{i+J}^X(t+J)$.
\end{itemize}    
\end{proof}

\begin{lemma}\label{lem:skip-j}
Fix $\lambda >0$ and a send sequence~$\p$.
Let $X$ be the backoff process with birth rate~$\lambda$ and send sequence~$\p$. 
Let $J$ be a positive integer.
Let $\p'$ be the send sequence 
such that, for all $i\geq 0$, $p'_i = p_{i+J}$.
There is a real number $\lambda''\in (0,\lambda) $  such that for any   $\lambda' \in (0, \lambda'')$  the following holds.
Let $X'$ be the backoff process with birth rate~$\lambda'$ and send sequence~$\p'$. Then there is a coupling of~$X'$ with~$X$ such
 that, for all $i\geq 0$ and $t\geq 0$,
$\b^{X'}_{i}(t) \subseteq \b^{X}_{i+J}(t+J)$. 
\end{lemma}
\begin{proof}
 
By Lemma~\ref{lem:dom1}, there is a real number $\lambda''\in (0,\lambda)$ and a distribution $\calD$ on non-negative integers 
with $\calD \gtrsim  \Po(\lambda'')$ such that the following holds.
Let $X''$ be the backoff process with birth distribution~$\calD$ and send sequence~$\p'$.
There is a coupling of~$X''$ with~$X$
such that, for all $i\geq 0$ and $t\geq 0$,  
$\b^{X''}_{i}(t) \subseteq \b^{X}_{i+J}(t+J)$.

Let $\calD' = \Po(\lambda')$ so that process~$X'$ has birth distribution~$\calD'$ and send sequence~$p'$.
Note that $\calD \gtrsim \calD'$.
The standard coupling of~$X'$ and~$X''$ (Definition~\ref{def:standard-couple})
has the property that, for every $i\geq 0$ and $t\geq 0$, $\b_i^{X'}(t) \subseteq \b_i^{X''}(t)$.

Putting the coupling of $X$ and $X''$ together with the coupling of $X''$ and $X'$ completes the proof.
   \end{proof}

\begin{corollary}\label{cor:skip-j}
Fix $\lambda >0$ and a send sequence~$\p$.
Let $X$ be the backoff process with birth rate~$\lambda$ and send sequence~$\p$. 
Let $J$ be a positive integer.
Let $\p'$ be the send sequence 
such that, for all $i\geq 0$, $p'_i = p_{i+J}$.
There is a real number $\lambda''\in (0,\lambda) $  such that for any   $\lambda' \in (0, \lambda'')$  the following holds.
Let $X'$ be the backoff process with birth rate~$\lambda'$ and send sequence~$\p'$. 
If $X'$ is not positive recurrent then $X$ is not positive recurrent.
\end{corollary}

\begin{proof}
 
Consider the coupling of~$X'$ and~$X$ from Lemma~\ref{lem:skip-j}. 
Let $T$ be the first time that $X$ returns to the empty state and let $T'$ be the first time that $X'$ returns to the empty state.
By Lemma~\ref{lem:skip-j}, for all non-negative integers~$i$ and~$t$, 
$\b^{X'}_{i}(t) \subseteq \b^{X}_{i+J}(t+J)$.  
Now, for all $t\in \{1,\ldots, T'-1\}$ there   is an index~$i_t$ such that 
$\b_{i_t}^{X'}(t)\neq \emptyset$   so $\b_{i_t+J}(t+J)\neq \emptyset$.
We conclude that $T\in \{1,\ldots,J\}$ or $T \geq T'+J$.

If $X'$ is not positive recurrent then $E[T']$ is infinite.
Since $\Pr(T>J)$ is bounded away from~$0$ (for example, 
$\Pr(T>J)$ is at least the probability that there are two births for each of the first $J$ steps, which is  ${(\lambda^2 e^{-\lambda}/2!)}^J$), we conclude that $E[T]$ is infinite, so  $X$ is not positive recurrent. 
\end{proof}

We finish this section by 
proving Theorem~\ref{thm:goal} using Theorem~\ref{thm:goaljkiller} (which will be proved in the rest of the paper).
We will use the following Lemma from~\cite{GL-oldcontention} (which was also known to the authors of~\cite{GJKP}). Its proof is straightforward and is based on showing that, under the hypothesis of the lemma, the expected lifetime of a newborn is infinite even if we dominate by letting it escape
when it sends, even when it collides, except for collisions with newborns.
\begin{lemma}\label{lem:killer-fact} \cite[Lemma 7]{GL-oldcontention}
Consider any birth rate $\lambda \in (0,1)$.
Consider any send sequence $\p$ such that, for infinitely many~$j$, $p_j \leq (\lambda p_0/2)^j$.
The backoff process with send sequence~$\p$ and birth rate~$\lambda$ is not positive recurrent.
\end{lemma}

\thmgoal*

\begin{proof}

In order to prove the Theorem, we first use Observation~\ref{obs:decrease-lambda}. By this observation, it suffices to prove the following claim.

\noindent{\bf Claim:\quad}{\sl
There exists a $\lambda_0 \in (0,1]$  such that the following holds.
Consider any $\lambda \in (0,\lambda_0)$
and any send sequence~$\p$. Then the backoff process with send sequence~$\p$ and birth rate~$\lambda$ is not positive recurrent.}

For the rest of the proof we will prove the Claim. Let $\lambda_0$ be the real number in $(0,1]$ from the statement of Theorem~\ref{thm:goaljkiller}. 
Consider any $\lambda \in (0,\lambda_0)$ 
and any send sequence~$\p$. 
Let $X$ be the backoff process with 
send sequence~$\p$ and
birth rate~$\lambda$. 
We wish to show that $X$ is not positive recurrent.

Suppose that there 
 is a positive integer $J'$ such that, for all $j\geq J'$,
$p_j \geq (\lambda p_0/2)^j$
(otherwise, we finish by  Lemma~\ref{lem:killer-fact}).
Suppose that there is a positive integer $J\geq J'$ such that $p_J \leq p_0$
(otherwise, we finish by Corollary~\ref{KM-cor}).
Note that, for all $j\geq J$, $p_j \geq (\lambda p_J/2)^j$.

Let $\p'$ be the send sequence such that, for all $j\geq 0$, $p'_j = p_{j+J}$.
By  Corollary~\ref{cor:skip-j}, there is a $\lambda'''\in (0,\lambda)$ such that for any $\lambda' \in (0,\lambda''')$ the following is true.  Let $X'$ be a backoff process with birth rate~$\lambda'$ and send sequence~$\p'$. 
If $X'$ is not positive recurrent then $X$ is not positive recurrent.

To finish, we fix this $\lambda'''\in (0,\lambda)$ from Corollary~\ref{cor:skip-j}.
Now let $\lambda' = \min(\lambda {(\lambda p_J/2)}^J, \lambda'''/2)$.
Let $X'$ be a backoff process with birth rate~$\lambda'$ and send sequence~$\p'$. Our goal is to show that $X'$ is not positive recurrent (then we are finished by Corollary~\ref{cor:skip-j} since $\lambda' \in (0,\lambda''')$).

By the definition of~$\lambda'$,  
$\lambda' \leq \lambda {(\lambda p_J/2)}^J $ so
for any positive integer~$\ell$, 
$(\lambda p_J/2)^J \geq (\lambda'/\lambda)^\ell$ 
and therefore $(\lambda p_J/2)^{\ell + J} \geq (\lambda' p_J/2)^\ell$.

Taking $\ell = j-J$ note that for all $\ell\geq 1$, $p'_\ell = p_{\ell+J} \geq (\lambda p_J/2)^{\ell+J}
\geq (\lambda' p_J/2)^\ell = (\lambda' p'_0/2)^\ell
$.

Finally, let $\p''$ be the send sequence that is identical to $\p'$ except that $p''_0=1$. Let $X''$ be a backoff process with birth rate 
 $\lambda'' = \lambda' p'_0$ and send sequence $\p''$. By Observation~\ref{obs:p0one}, we just need to show that $X''$ is not positive recurrent and this implies that $X'$ is not positive recurrent.
 
Note that 
$\lambda'' \leq \lambda ' \leq \lambda''' \leq \lambda < \lambda_0$. Also, $p''_0=1$ and
for all $\ell \geq 1$, $p''_\ell = p'_\ell \geq  
(\lambda' p'_0/2)^\ell = 
(\lambda'' /2)^{\ell}$. 
Now we  use Theorem~\ref{thm:goaljkiller}.
  \end{proof}

\section{The \texorpdfstring{$j$}{j}-jammed process}
\label{sec:j-jammed}

The $j$-jammed process from Definition~\ref{def:j-jammed} will be used throughout our proof. For example, in the case where there are finitely-many strongly exposed bins, it will be used in the analysis of the volume process that tracks the sends.

Since our goal is to prove Theorem~\ref{thm:goaljkiller} where $p_0=1$ we assume $p_0=1$ throughout the rest of the paper -- this simplifies notation since every generalised backoff process~$Z$ then has $b^Z_0(t) = 0$ (so we omit bin~$0$ from the list of bins).
Since Theorem~\ref{thm:goaljkiller} 
only requires us to find a real number $\lambda_0\in(0,1]$ and to show lack of positive recurrence for $\lambda \in (0,\lambda_0)$ we freely apply upper bounds to $\lambda$ where convenient.

\begin{definition}\label{def:j-jammed}
Fix a positive integer~$j$. Fix a send sequence~$\p$ with $p_0=1$, a non-negative integer~$\tau$, and a $\lambda\in (0,1)$. 
Fix a $j$-tuple $\avect \in (\integers_{\geq 0})^j$.
A \emph{$j$-jammed} process~$Y$ with send sequence~$\p$, start time~$\tau$, birth rate~$\lambda$, 
and initial population size tuple~$\avect$   
is a generalised backoff process~$Y$ with  
send sequence~$\p$, start time~$\tau$, birth distribution~$\Po(\lambda)$ and an initial population~$\arrvect^Y(\tau)$ that satisfies the  following constraints.  
\begin{itemize}
\item $\romanarrvect^Y_{[j]}(\tau) = \avect$ and for $i>j$, $ \romanarr^Y_i(\tau) = 0$.
\item For all $t> \tau$ and $i\geq 1$, $\romanarr_i^Y(t) = 0$.
\item For all $t> \tau$ and $i \in [j-1]$, $\e_i^Y(t) = \emptyset$.
\item For all $t>\tau$ and $i \ge j$, $\e^Y_i(t) = \s^Y_i(t)$.
\end{itemize}
\end{definition}

Note that Definition~\ref{def:j-jammed} does not specify the names of the balls in $\cup_{i\in [j]} \arr_i^Y(\tau)$ so there are many $j$-jammed processes $Y$ with send sequence~$\p$, start time~$\tau$, birth rate~$\lambda$, and initial population size tuple~$\avect$, depending on these ball names. However, the ball names won't turn out to be important. Definition~\ref{def:j-jammed-dist} specifies a distribution on the vector~$\avect$ in Definition~\ref{def:j-jammed}.

\begin{definition} \label{def:j-jammed-dist}
Fix a positive integer~$j$. Fix a send sequence~$\p$ with $p_0=1$, a non-negative integer~$\tau$, and a $\lambda\in (0,1)$. 
Let $\calP$ be a distribution with support $(\integers_{\geq 0})^j$.
A \emph{$j$-jammed} process~$Y$ with send sequence~$\p$, start time~$\tau$, birth rate~$\lambda$, and initial population size distribution~$\calP$ is a 
$j$-jammed process with send sequence~$\p$, start time~$\tau$, birth rate~$\lambda$, and initial population size tuple~$\avect$, where $\avect \sim \calP$.

\end{definition}

Fix a send sequence~$\p$ with $p_0=1$, a birth rate~$\lambda$, and a $j$-tuple $\zvect\in (\reals_{\geq 0})^j$.
We will define a function~$f^{\zvect}$, depending on these parameters,
and will consider $Y$, a $j$-jammed process with send sequence~$\p$, a start time~$\tau$, birth rate~$\lambda$, and initial population size distribution $\Po(\zvect)$.
Lemma~\ref{lem:f-track-jammed}
shows that for $t\geq \tau$,
$\b_{[j]}^Y(t) \sim \Po(f_{[j]}^{\zvect}(t-\tau))$.
In order to allow re-use of the lemmas that we build up,
we define~$f^{\zvect}$ with an additional argument - a $j$-tuple~$\Gammavect$ noted in its superscript, but taken to be~$\zerovect$ here.
We also allow an arbitrary offset~$\Delta$ (a non-negative integer) and start~$Y$ from an initial population-size tuple drawn from distribution 
$\Po(f_{[j]}^{\zerovect,\zvect}(\Delta))$, adjusting the result accordingly.
We start with the definition of~$f^{\Gammavect,\zvect}$.

\begin{definition}\label{def:f}
Fix a positive integer~$j$, a send sequence~$\p$ with $p_0=1$, a $\lambda\in (0,1)$, a $j$-tuple $\zvect\in (\reals_{\geq 0})^j$, and a
$j$-tuple 
$\Gammavect \in [0,1]^j$. 
We now define a function $f^{\Gammavect,\zvect}$ from non-negative integers to 
tuples of non-negative reals indexed by $\{0,\ldots,j\}$.
\newsym{f-x}{component $x$ of function $f^{\Gammavect, \zvect}$ (Definition~\ref{def:f})}{f_x^{\Gammavect,\zvect}(t)} denotes the $x$'th component of 
$f^{\Gammavect,\zvect}(t)$.
For all $x\in  [j]$, $f_x^{\Gammavect,\zvect}(0) := z_x$. For all $t\geq 0$,
$f_0^{\Gammavect,\zvect}(t) := \lambda$. 
Finally,
for all $t\geq 1$ and all $x\in  [j]$, 
$f_x^{\Gammavect,\zvect}(t) := (1-p_x)f_x^{\Gammavect,\zvect}(t-1) + (1-\Gamma_{x})p_{x-1}f^{\Gammavect,\zvect}_{x-1}(t-1)$.
\end{definition}

The most straightforward part of Lemma~\ref{lem:f-track-jammed}
is calculating the expectation of $b_x^Y(t)$
using Definition~\ref{def:f}.
This is done in Observation~\ref{obs:f-gen}.

\begin{observation}\label{obs:f-gen}
Fix a positive integer~$j$, a send sequence~$\p$ with $p_0=1$, a $\lambda\in (0,1)$, and a $j$-tuple $\zvect\in (\reals_{\geq 0})^j$.
Fix non-negative integers~$\tau$ and~$\Delta$.
Let $Y$ be a $j$-jammed process with send sequence~$\p$, start time~$\tau$, birth rate~$\lambda$, and initial population size distribution $\Po(f_{[j]}^{\zerovect,\zvect}(\Delta))$.  
Then for all $t\geq \tau$ and $x\in  [j]$, $\E[b_x^{Y}(t)] = f_x^{\zerovect,\zvect}(\Delta + t-\tau)$. 
\end{observation}
\begin{proof} 
Write $f$ for $f^{\zerovect,\zvect}$.
The proof is by induction on~$t$, with base case $t=\tau$.
For $t=\tau$ and $x\in [j]$, $\E[b_x^{Y}(t)] = f_x(\Delta)$. 
For the inductive step (from Definitions~\ref{def:j-jammed-dist},
\ref{def:j-jammed}, and~\ref{def:genbackoff}) for $t>\tau$,
\[\E[b_x^{Y}(t)] 
= \E[B_x^{Y}(t) - s_x^{Y}(t) + s_{x-1}^{Y}(t) ]
= \E[B_x^{Y}(t)(1-p_x)  + 
B_{x-1}^{Y}(t) p_{x-1}].\]
If $x=1$ then
this is $\E[B_1^{Y}(t)(1-p_1)  + n^{Y}(t)]
= (1-p_1)\E[b_1^{Y}(t-1)]  + \lambda$ so by the inductive hypothesis, this is 
$(1-p_1)
f_1(\Delta + t-1-\tau)  + \lambda$, which is $f_1(\Delta + t-\tau)$, as required.
For $x>1$
this is 
$\E[b_x^{Y}(t-1)(1-p_x)  + 
b_{x-1}^{Y}(t-1) p_{x-1}]$, which is again $f_x(\Delta + t-\tau)$, as required. 
\end{proof}

In order to establish the fact that
$b_1^Y(t),\ldots,b_j^Y(t)$ are independent Poisson random variables, we define the notion of \emph{ball trajectories} (Definition~\ref{def:trajectory}).

\begin{definition}
\label{def:trajectory}
Fix a positive integer~$j$. Fix a send sequence~$\p$ with $p_0=1$, a non-negative integer~$\tau$, and a $\lambda\in (0,1)$. 
Let $\calP$ be a distribution with support $(\integers_{\geq 0})^j$.
Let $Y$ be a $j$-jammed process with send sequence~$\p$, start time~$\tau$, birth rate~$\lambda$, and initial population size distribution~$\calP$.
Let $t\geq\tau$ be an integer.

A \emph{ball trajectory} up to time~$t$ is a 
tuple \newsym{tr}{ball trajectory up to time $t$}{\vectraje  := ({\traje}_\tau,\ldots,{\traje}_t)} of integers in $\{0, \ldots, j+1\}$ 
such that
$\traje_\tau \leq j$ and, for all $t' \in \{\tau+1,\ldots,t\}$,
$\traje_{t'} \in 
\{\traje_{t'-1}, \traje_{t'-1} + 1\}$.

For any $t'\in \{\tau,\ldots,t\}$ such that ${\traje}_{t'} \in [j]$, a ball following trajectory~$\vectraje$ is in  $\b_{{\traje}_{t'}}^Y(t')$.
If  ${\traje}_{t'} = 0$, then the ball is born after time~$t'$ -- that is, there is a $t''>t'$ such that the ball is in $\n^Y(t'')$.
If ${\traje}_{t'} = j+1$, then the ball sends from bin~$j$ (and leaves the system) by time~$t'$ -- that is, there is a $t'' \leq t'$ such that the ball is in $\s_j^Y(t'')$.

For every   integer~$t\geq \tau$ and every
$x\in \{0,\ldots,j+1\}$, let 
\newsym{T-x-t}{set of ball trajectories ending in bin $x$}{\calT^Y(x,t)} be the set of all ball trajectories up to time~$t$ 
with ${\traje}_t = x$. 
Let 
\newsym{T-x}{set of ball trajectories}{\calT^Y(t)} be the set of all ball trajectories up to time~$t$, that is $\calT^Y(t) = \cup_{x\in \{0,\ldots,j+1\}} \calT^Y(x,t)$.
For $\traje \in \calT^Y(t)$, let 
\newsym{b-tr}{set of balls following a trajectory}{\b^Y_\traje} be the set of balls that follow trajectory~$\traje$ and let 
\newsym{b-tr-2}{$|\b^Y_\traje|$}{b^Y_\traje} be  $|\b^Y_\traje|$.

\end{definition}

Observation~\ref{obs:latertraj} follows easily from Definition~\ref{def:trajectory}.

\begin{observation}
\label{obs:latertraj}
Fix a positive integer~$j$. Fix a send sequence~$\p$ with $p_0=1$, a non-negative integer~$\tau$, and a $\lambda\in (0,1)$. 
Let $\calP$ be a distribution with support $(\integers_{\geq 0})^j$.
Let $Y$ be a $j$-jammed process with send sequence~$\p$, start time~$\tau$, birth rate~$\lambda$, and initial population size distribution~$\calP$.
Let $t_1$ and $t_2$ be  
integers with $\tau \leq t_1 \leq t_2$.
Then, for any $\traje \in  \calT^Y(t_1)$, there is a subset $S \subseteq \calT^Y(t_2)$ 
such that $\b_\traje^Y = \cup_{\traje' \in S} \b_{\traje'}^Y$. (In particular, $S$ is the set of all trajectories up to~$t_2$ that extend $\traje$.)
\end{observation}

Observations~\ref{obs:firsttraj} and \ref{obs:secondtraj}
establish independence between trajectories and bins.

\begin{observation}\label{obs:firsttraj}
Fix a positive integer~$j$. Fix a send sequence~$\p$ with $p_0=1$, a non-negative integer~$\tau$, and a $\lambda\in (0,1)$. 
Fix $\muvect \in (\reals_{\geq 0})^j$. 
Let $Y$ be a $j$-jammed process with send sequence~$\p$, start time~$\tau$, birth rate~$\lambda$, and initial population size distribution~$\Po(\muvect)$.
Let $t\geq \tau$ be an integer.
For every $\traje \in \calT^Y(t)$, $b^Y_\traje$ is a Poisson random variable.
Also, $\{\b_{\traje}^Y\colon \traje \in \calT^Y(t)\}$ are disjoint sets whose sizes are mutually independent random variables.
\end{observation}
\begin{proof}
Let $j'$ be the maximum entry in the trajectory $\traje$
and for each $\ell \in [j']$ let $M_\ell$ be the number of entries in~${\traje}$ that are~$\ell$.
If $j'\leq j$ then let $W = W_{j'}$. Otherwise, let $W=1$.
Let $\rho = W \prod_{\ell = \max\{{\traje}_{\tau},1\}}^{\min\{j,j'\}} 
(1-p_\ell)^{M_\ell-1} p_\ell
$.

If ${\traje}_{\tau}\in [j]$ then 
$b^Y_\traje$ is Poisson with mean 
$\rho \mu_{{\traje}_{\tau}} $. 
If ${\traje}_{\tau} = 0$ then
$b^Y_\traje$ is Poisson with mean $\rho \lambda$.
Disjointness follows immediately from the definition of trajectories, and independence follows immediately from the fact that Bernoulli events (such as ball-sending) split Poisson random variables into independent Poisson random variables.
\end{proof}

\begin{observation} \label{obs:secondtraj}
Fix a positive integer~$j$. Fix a send sequence~$\p$ with $p_0=1$, a non-negative integer~$\tau$, and a $\lambda\in (0,1)$. 
Fix $\muvect \in (\reals_{\geq 0})^j$. 
Let $Y$ be a $j$-jammed process with send sequence~$\p$, start time~$\tau$, birth rate~$\lambda$, and initial population size distribution~$\Po(\muvect)$.
Let $t\geq \tau$ be an integer.
Then $b_1^Y(t),\ldots,b^Y_j(t)$ are independent Poisson random variables.
\end{observation} 
\begin{proof}
For $x\in [j]$, $b_x^Y(t) = \sum_{\traje \in \calT^Y(x,t)} b_\traje^Y$. By Observation~\ref{obs:firsttraj}, the terms in the sum are Poisson and independent, so the sum is a Poisson random variable.
Independence for different $x$-values follows from the fact that they correspond to sums over different $\traje$ values, and these are independent.
\end{proof}

 Lemma~\ref{lem:f-track-jammed} 
combines these observations and
adds information about positive correlation.

\begin{lemma}\label{lem:f-track-jammed}
Fix a positive integer~$j$, a send sequence~$\p$ with $p_0=1$, a $\lambda\in (0,1)$, and a $j$-tuple $\zvect\in (\reals_{\geq 0})^j$.
Fix non-negative integers~$\tau$ and~$\Delta$.
Let $Y$ be a $j$-jammed process with send sequence~$\p$, start time~$\tau$, birth rate~$\lambda$, and initial population size distribution $\Po(f_{[j]}^{\zerovect,\zvect}(\Delta))$.  

Then, for all $t\geq \tau$, 
$b_{[j]}^Y(t) \sim \Po(f_{[j]}^{\zerovect,\zvect}(\Delta+t-\tau))$.
Furthermore, the following positive correlation property holds.
Let $\calE_1$ and $\calE_2$ be any events 
determined by  
the random variables in 
$\{ b_\traje^Y \mid \traje \in \calT^Y(t)\}$ such that $\calE_2$ occurs with positive probability. 
If the indicator variables of  
$\calE_1$ and $\calE_2$ are increasing functions of 
these random variables   then $\Pr(\calE_2 \mid \calE_1) \geq \Pr(\calE_2)$.
\end{lemma}
\begin{proof}
The first part of the lemma 
(the distribution of $b_{[j]}^Y(t)$)
follows immediately from Observation~\ref{obs:f-gen} 
and Observation~\ref{obs:secondtraj}. 
For the positive correlation, note that the random variables in
$\{ b_\traje^Y \mid \traje \in \calT^Y(t)\}$ are independent of each other by Observation~\ref{obs:firsttraj}. 
Let $g_1$ and $g_2$ be the indicator variables of $\calE_1$ and $\calE_2$. 
Now Harris's inequality  (Lemma~\ref{lem:Harris})
guarantees that $\E[g_1 g_2] \geq \E[g_1] \E[g_2]$
which gives 
$\pr(\calE_2 \wedge \calE_1) \ge \pr(\calE_2)\pr(\calE_1)$, as desired.
\end{proof}

In order to make use of Lemma~\ref{lem:f-track-jammed} it will be useful to compute lower bounds on $f_x^{\zerovect, \zvect}(t)$,
starting from  
$\Delta=0$ and
$\zvect=\zerovect$. We will do this using Lemma~\ref{lem:f-fill-all-bins}, which is stated more generally so that it can be re-used. We first present the definitions and lemmas that we need in order to prove Lemma~\ref{lem:f-fill-all-bins}.
Definition~\ref{def:mu} defines the quantity $\mu_x^{\Gammavect}$, which will be used throughout the paper.

\begin{definition}\label{def:mu}
Fix a positive integer~$j$, a send sequence~$\p$ with $p_0=1$, a $\lambda\in (0,1)$, and a $j$-tuple $\Gammavect\in  [0,1]^j$.
For all $x \in  \{0,\ldots,j\}$
let  
\newsym{mu-Gamma}{quantity from Definition~\ref{def:mu}}{\mu^{\Gammavect}_x =  {\lambda}{W_x}\prod_{a=1}^x (1-\Gamma_a)}.  
Note that $\mu^{\Gammavect}_0 = {\lambda}\prod_{a=1}^0 (1-\Gamma_a)$.
The empty product is~$1$ by convention, so $\mu^{\Gammavect}_0=\lambda$. 
\end{definition}

Lemma~\ref{lem:f-fill-bin} tracks the lower bound on $f_x^{\zerovect,\zvect}$. The way to think of it is that the premise is that all bins up to bin~$k$ are full in expectation, at least in the sense that the lower bound on their expectation is at least $\mu_x^{\Lambdavect}$.
Item (i) says that this continues. 
Item (ii) tracks progress on filling bin~$k+1$ (in expectation) for as long as it is not already full by bounding how much gets added to the expectation during~$t$ steps.
Item (iii) uses this to calculate a bound on how long it takes until bin~$k+1$ is also full in expectation.

 \begin{lemma}\label{lem:f-fill-bin}
Fix a positive integer~$j$, a send sequence~$\p$ with $p_0=1$, a $\lambda\in (0,1)$, 
a $j$-tuple
$\zvect \in (\reals_{\geq 0})^j$, and
$j$-tuples $\Gammavect$ and $\Lambdavect$
in $[0,1)^j$ satisfying
$\Gammavect \leq \Lambdavect$.
Fix a non-negative integer~$t_0$ and
an integer $k\in [j]$. Suppose that, for all $x \in [k]$, $f_x^{\Gammavect,\zvect}(t_0) \ge \mu_x^{\Lambdavect}$. 
Then the following hold.
\begin{enumerate}[(i)]
\item  For all $t \ge 0$ and $x \in [k]$,  $f_x^{\Gammavect,\zvect}(t_0+t) \ge \mu_x^{\Lambdavect}$. 
\item If $k+1 \leq j$ then for all $t \ge 0$,
 $
f_{k+1}^{\Gammavect,\zvect}(t_0+t) \ge \min\big\{\mu_{k+1}^{\Lambdavect},\ f_{k+1}^{\Gammavect,\zvect}(t_0) + t(\Lambda_{k+1} - \Gamma_{k+1})p_k\mu_k^{\Lambdavect}\big\}$.
\item If $k+1 \leq j$ then for all $t \ge  
W_{k+1}
/(\Lambda_{k+1} - \Gamma_{k+1})$,  $f_{k+1}^{\Gammavect,\zvect}(t_0+t) \ge \mu_{k+1}^{\Lambdavect}$.
\end{enumerate}
\end{lemma}
\begin{proof}
We abbreviate $f^{\Gammavect,\zvect}$ by $f$
and $\mu^{\Lambdavect}$ by~$\mu$
throughout. 

We first prove (i) by induction on~$t$.
The base case is $t=0$ which is by supposition. For the inductive step fix $t\geq 0$ and assume (i) for~$t$. 
By Definition~\ref{def:f}, for all $x \in [j]$,
 $
f_x(t_0+t+1) = (1-p_x)f_x(t_0+t) + (1-\Gamma_x)p_{x-1}f_{x-1}(t_0+t)
$.
By the induction hypothesis, for $x-1 \leq k$, 
$f_{x-1}(t_0+t) \geq \mu_{x-1}$.
By  the definition of $\mu_x$,
for all $x\in [j]$,
$(1-\Lambda_x) p_{x-1} \mu_{x-1} = p_x \mu_x$. 
Thus, for $x\in [\min\{k+1,j\}]$, 
\begin{align*}
f_x(t_0+t+1) &\ge (1-p_x)f_x(t_0+t) + (1-\Gamma_{x})p_{x-1}\mu_{x-1}\\
 &= (1-p_x)f_x(t_0+t) + (1-\Lambda_x)p_{x-1}\mu_{x-1} + (\Lambda_x - \Gamma_x)p_{x-1}\mu_{x-1}\\
&= (1-p_x)f_x(t_0+t) + p_x\mu_x + (\Lambda_x - \Gamma_x)p_{x-1}\mu_{x-1}\\
&\geq (1-p_x)f_x(t_0+t) + p_x\mu_x,  
\end{align*}
establishing (i) for $x\in [k]$ since, by
the inductive hypothesis, 
$f_x(t_0+t) \geq  \mu_x$.

For the rest of the proof suppose $k+1 \leq j$. We will now prove (ii) by induction on~$t$. Once again, the base case is $t=0$, which is straightforward.
Fix $t\geq 0$ and assume that (ii) holds for~$t$. 
The calculation that we did in part (i) showed
that 
\[ f_{k+1} (t_0+t+1) \geq (1-p_{k+1}) f_{k+1}(t_0+t) + p_{k+1} \mu_{k+1}.\]
If $f_{k+1}(t_0+t) \geq \mu_{k+1}$, then we have established
$f_{k+1}(t_0+t+1) \geq \mu_{k+1}$ so we have established
(ii) for~$t+1$.
Otherwise,  
applying the inductive hypothesis (ii) to~$t$ 
we have 
$f_{k+1}(t_0+t) \geq f_{k+1}(t_0) + t(\Lambda_{k+1} - \Gamma_{k+1}) p_k \mu_{k}$.
The calculation that we did in part (i) showed that
\[
f_{k+1} (t_0+t+1) \geq (1-p_{k+1}) f_{k+1}(t_0+t) + p_{k+1} \mu_{k+1}
+ (\Lambda_{k+1} - \Gamma_{k+1})p_{k}\mu_{k},\]
and plugging in the lower bound that we just derived for $f_{k+1}(t_0+t)$, we get
 \[
f_{k+1}(t_0+t+1) \ge \big(f_{k+1}(t_0) + t(\Lambda_{k+1} - \Gamma_{k+1})p_k\mu_k\big) - p_{k+1}f_{k+1}(t_0+t) + p_{k+1}\mu_{k+1} + (\Lambda_{k+1} - \Gamma_{k+1})p_k\mu_k.
 \]
    
Since $f_{k+1}(t_0+t) < \mu_{k+1}$, it follows that
\begin{align*}
f_{k+1}(t_0+t+1) &\ge \big(f_{k+1}(t_0) + t(\Lambda_{k+1} - \Gamma_{k+1}) p_k \mu_k\big) + (\Lambda_{k+1} - \Gamma_{k+1})p_k\mu_k\\
&= f_{k+1}(t_0) + (t+1)(\Lambda_{k+1} - \Gamma_{k+1})p_k\mu_k,
\end{align*} establishing (ii) for $t+1$.

Finally, to prove (iii) suppose $k+1\leq j$ and consider 
$t \ge  
W_{k+1}
    /(\Lambda_{k+1} - \Gamma_{k+1})$.
Note that
 $t (\Lambda_{k+1} - \Gamma_{k+1})\geq W_{k+1}(1-\Lambda_{k+1})$ so 
\[t (\Lambda_{k+1} - \Gamma_{k+1}) p_k \mu_k 
\geq
W_{k+1}(1-\Lambda_{k+1}) 
\frac{p_{k+1}\mu_{k+1}}
{1-\Lambda_{k+1}}
= 
\mu_{k+1},\]
so by (ii)
\[f_{k+1}(t_0+t) \geq 
\min\big\{\mu_{k+1},\ f_{k+1}(t_0) + t(\Lambda_{k+1} - \Gamma_{k+1})p_k\mu_k\big\}
\geq \min\big\{\mu_{k+1},
\ f_{k+1}(t_0) + \mu_{k+1}\big\} \geq \mu_{k+1}.\]

\end{proof}

Lemma~\ref{lem:f-bound} gives an upper bound on the expected
fullness of bins.

\begin{lemma}\label{lem:f-bound}
Fix a positive integer~$j$, a send sequence~$\p$ with $p_0=1$, a $\lambda\in (0,1)$, 
a $j$-tuple
$\Gammavect \in [0,1]^j$
and a $j$-tuple
$\zvect \in (\reals_{\geq 0})^j$, 
such that 
$z_x \leq \mu_x^{\Gammavect}$ for all $x\in [j]$.
Then for all $x \in \{0,\ldots,j\}$ and all $t \ge 0$,  $0 \le f^{\Gammavect,\zvect}_x(t) \le \mu_x^{\Gammavect}$. 
\end{lemma}

\begin{proof}
We abbreviate $f^{\Gammavect,\zvect}$ by~$f$
and $\mu^{\Gammavect}$ with~$\mu$ throughout.
From Definition~\ref{def:f} and Definition~\ref{def:mu},
for all $t\geq 0$,
$f_0(t) = \lambda = \mu_0$.

We prove 
$0 \le f_x(t) \le \mu_x$
for 
all $x\in [j]$ and for
$t\geq 0$ by induction on $t$. 
The base case is $t=0$. By definition
$f_x(0)  =z_x \le \mu_x$.

For the inductive step fix $t\geq 0$ and suppose 
that 
for all $x \in \{0,\ldots,j\}$,  $0 \le f_x(t) \le \mu_x$.  Then, for $x\in [j]$,
\begin{align*}
f_x(t+1) &= (1-p_x)f_x(t) + (1-\Gamma_x)p_{x-1}f_{x-1}(t) 
\le (1-p_x)\mu_x + (1-\Gamma_x)p_{x-1}\mu_{x-1} \\
&= (1-p_x){\lambda W_x}\prod_{a=1}^x (1-\Gamma_a) + \lambda \prod_{a=1}^x (1-\Gamma_a) = {\lambda W_x}\prod_{a=1}^x (1-\Gamma_a) = \mu_x,
\end{align*}
giving  
$f_x(t+1) \leq \mu_x$.
The fact that $0\leq f_x(t+1)$ is straightforward from the first equality and the inductive hypothesis. 
\end{proof}

Observation~\ref{obs:rearrange} is a technicality that it is useful in the proof of Lemma~\ref{lem:gen-f-monotone}, which shows that the sum of bin-expected fullnesses does not shrink over time (assuming that they weren't over-full in expectation to start with).

 \begin{observation}
\label{obs:rearrange}
Fix a positive integer~$j$ and a
a $j$-tuple
$\Gammavect \in [0,1]^j$.  Fix $k\in \{0,\ldots,j-1\}$. 
Then 
$(1-\Gamma_1) - \sum_{x=1}^k \Gamma_{x+1} \prod_{a=1}^x (1-\Gamma_a) = \prod_{a=1}^{k+1} (1-\Gamma_a)$.
\end{observation}   
\begin{proof}
The proof is by induction on~$k$.
The base case is $k=0$. For the inductive step fix $k\in [j-1]$ and assume 
the result for $k-1$. Then the  left-hand side equals 
\[
(1-\Gamma_1) - \sum_{x=1}^{k-1} \Gamma_{x+1} \prod_{a=1}^x (1-\Gamma_a) - \Gamma_{k+1}\prod_{a=1}^k (1-\Gamma_a) = 
\prod_{a=1}^{k} (1-\Gamma_a) - \Gamma_{k+1} \prod_{a=1}^k (1-\Gamma_a) =  \prod_{a=1}^{k+1} (1-\Gamma_a).
\]
\end{proof}

\begin{lemma}\label{lem:gen-f-monotone}
Fix positive integers~$j'$ and~$j$ with $j' \leq j$. Fix  a send sequence~$\p$ with $p_0=1$, a $\lambda\in (0,1)$, 
a $j$-tuple
$\Gammavect \in [0,1]^j$
and a $j$-tuple
$\zvect \in (\reals_{\geq 0})^j$, 
such that 
$z_x \leq \mu_x^{\Gammavect}$ for all $x\in [j]$.
Then for all $t \ge 1$,     $\sum_{x \in [j']} f^{\Gammavect,\zvect}_x(t) \ge \sum_{x \in [j']} f^{\Gammavect,\zvect}_x(t-1)$.
\end{lemma}
\begin{proof}
We abbreviate $f^{\Gammavect,\zvect}$ by~$f$
and $\mu^{\Gammavect}$ by~$\mu$ throughout.
By Definition~\ref{def:f}, for $t\geq 1$,
\begin{align*}
\sum_{x = 1}^{j'}f_x(t) - \sum_{x=1}^{j'}f_x(t-1)
&= \sum_{x=1}^{j'} \big( (1-\Gamma_x)p_{x-1}f_{x-1}(t-1) -p_xf_x(t-1) \big)\\
&= (1-\Gamma_1) f_0(t-1) - \sum_{x=1}^{j'-1} \Gamma_{x+1}p_xf_x(t-1)  -p_{j'}f_{j'}(t-1).
\end{align*}
For all $t\geq 0$, $f_0(t)=\lambda$.
By Lemma~\ref{lem:f-bound}, 
for all $x \in \{0,\ldots,j\}$ and all $t \ge 0$,  $0 \le f_x(t) \le \mu_x$. 
Thus, for $t\geq 1$,  
\begin{align*}
\sum_{x = 1}^{j'} f_x(t) - \sum_{x=1}^{j'} f_x(t-1) &\ge \lambda(1-\Gamma_1) - \sum_{x=1}^{j'-1} \Gamma_{x+1}p_x\mu_x - p_{j'}\mu_{j'}\\
&= \lambda(1-\Gamma_1) - \sum_{x=1}^{j'-1} \Gamma_{x+1}\lambda \prod_{a=1}^x (1-\Gamma_a) - \lambda\prod_{a=1}^{j'} (1-\Gamma_a).
\end{align*}

By Observation~\ref{obs:rearrange} with $k=j'-1$,
$
        \sum_{x = 1}^{j'} f_x(t) - \sum_{x=1}^{j'} f_x(t-1) \ge 0,
    $
    as required.
\end{proof}

Lemma~\ref{lem:f-increasing} establishes a useful monotonicity with respect to the superscripts of the $f$-function.

\begin{lemma}\label{lem:f-increasing}
Fix   a send sequence~$\p$ with $p_0=1$, a $\lambda\in (0,1)$.
Fix a positive integer~$j$
and $j$-tuples $\Gammavect$ and $\Lambdavect$ 
in $[0,1]^j$ satisfying $\Gammavect \leq \Lambdavect$. Fix $j$-tuples $\avect$ 
and $\zvect$ in $(\reals_{\geq 0})^j$
satisfying $\zvect \geq \avect$.
Then for all $k \in [j]$ and $t \ge 0$, $f_k^{\Gammavect,\zvect}(t) \ge f_k^{\Lambdavect,\avect}(t)$.
\end{lemma}
\begin{proof}
The proof is by induction on $t$. 
The base case is $t=0$. For this case, for all $k \in [j]$, $f_k^{\Gammavect,\zvect}(0) = z_k \ge a_k = f_k^{\Lambdavect,\avect}(0)$ as required. 
Suppose that the result holds for some $t \ge 0$. Then by Definition~\ref{def:f}, for all $k \in [j]$,
\[
f_k^{\Gammavect,\zvect}(t+1) = (1-p_k)f_k^{\Gammavect,\zvect}(t) + (1-\Gamma_k)p_{k-1}f_{k-1}^{\Gammavect,\zvect}(t).
    \]
By the induction hypothesis, since $\Gammavect \le \Lambdavect$, it follows that
\[
f_k^{\Gammavect,\zvect}(t+1) \ge (1-p_k)f_k^{\Lambdavect,\avect}(t) + (1-\Lambda_k)p_{k-1}f_{k-1}^{\Lambdavect,\avect}(t) = f_k^{\Lambdavect,\avect}(t+1).
    \]
\end{proof}

Finally, Lemma~\ref{lem:f-fill-all-bins} gives us the desired upper bound on how long it takes to fill all $j$~bins (in expectation). The quantity~$R$ measures, collectively, how far under-full the bins are at the start. The amount of time needed to fill is linear in~$R j$.

\begin{lemma}\label{lem:f-fill-all-bins}
Fix a positive integer~$j$, a send sequence~$\p$ with $p_0=1$, a $\lambda\in (0,1)$, 
a $j$-tuple
$\zvect \in (\reals_{\geq 0})^j$, 
and
$j$-tuples $\Gammavect\in [0,1]^j$ 
and $\Lambdavect\in (0,1]^j$.
Suppose $\Gammavect \leq \tfrac12 \Lambdavect $ and, for all $k \in [j]$, $\mu_k^{\Lambdavect} \ge \lambda W_k/2$. 
Fix a non-negative real number~$R$
such that $\sum_{k \in [j]} \max\{0, \mu_k^{\Lambdavect} - z_k\} \le R$. 
Then, for all \[t \ge (4/\lambda)jR/\min\{\Lambda_k\colon k \in [j]\}\] and all $
k \in [j]$, $f_k^{\Gammavect,\zvect}(t) \ge \mu_k^{\Lambdavect}$.
\end{lemma}

\begin{proof}
For all $k\in [j]$ let $a_k = \min\{z_k, \mu^{\Lambdavect}_k\}$
and note that $0 \leq a_k \leq \mu^{\Lambdavect}_k$.
For every $\hat{j}\in\{0,\ldots,j\}$, let
$t_0(\hat{j}) = (4/\lambda) \hat{j} R/\min\{{\Lambda}_k\colon k \in [j]\}$.

We will prove by induction on~$\hat{j}$
that for all $k\in [\hat{j}]$ and all $t\geq t_0(\hat{j})$, $f^{\Gammavect,\zvect}_k(t) \ge \mu^{\Lambdavect}_k$.

The base case, $\hat{j}=0$, is vacuous, since there are no $k\in [\hat{j}
]$.
Suppose that the inductive hypothesis holds  for a $\hat{j} \in \{0,\ldots,j-1\}$. We will establish the inductive hypothesis for $\hat{j}+1$.
Note that we just need to show, 
for $t\geq t_0(\hat{j}+1)$, that 
$f^{\Gammavect,\zvect}_{\hat{j}+1}(t) \ge \mu^{\Lambdavect}_{\hat{j}+1}$ -- the other values of~$k$ follow immediately from the inductive hypothesis.

Since $0 \le a_k \le \mu_k$ for all $k\in [\hat{j}+1]$, we can apply Lemma~\ref{lem:gen-f-monotone} 
(with the $\Gammavect$ in the lemma statement equal to our $\Lambdavect$ and the $j'$ in the lemma statement equal to our~$\hat{j}+1$ and the $\zvect$ in the lemma statement equal to our~$\avect$) to show that for all $t \ge 0$, 
\[
\sum_{k=1}^{\hat{j}+1} f_k^{\Lambdavect,\avect}(t) \ge \sum_{k=1}^{\hat{j}+1} f_k^{\Lambdavect,\avect}(0) = \sum_{k=1}^{\hat{j}+1} \mu_k^{\Lambdavect} - \sum_{k=1}^{\hat{j}+1} (\mu_k^{\Lambdavect} - a_k) = \sum_{k=1}^{\hat{j}+1} \mu_k^{\Lambdavect} - \sum_{k=1}^{\hat{j}+1} \max\{0, \mu_k^{\Lambdavect}-z_k\} \ge \sum_{k=1}^{\hat{j}+1} \mu_k^{\Lambdavect} - R.
    \]

Moreover, by Lemma~\ref{lem:f-bound}, for all $k \in [\hat{j}+1]$ and all $t \ge 0$ we have $0 \le f^{\Lambdavect,\avect}_k(t) \le \mu_k^{\Lambdavect}  $.  
It therefore follows that, for all $t \ge 0$ and all $k \in [\hat{j}+1]$,  $f^{\Lambdavect,\avect}_k(t) \ge \mu_k^{\Lambdavect}-R$.  Taking $t=t_0(\hat{j})$ and applying Lemma~\ref{lem:f-increasing}, 
\begin{equation}\label{eq:f-fill-all-bins}
f_{\hat{j}+1}^{\Gammavect,\zvect}(t_0(\hat{j})) \ge \mu_{\hat{j}+1}^{\Lambdavect}-R.
\end{equation}

By the induction hypothesis,  for all $t \ge t_0(\hat{j})$ and all $k \in [\hat{j}
]$,  $f_k^{\Gammavect,\zvect}(t) \ge \mu_k^{\Lambdavect}$. By Lemma~\ref{lem:f-fill-bin}(ii)
(with the $t_0$ of Lemma~\ref{lem:f-fill-bin} equal to $t_0(\hat{j})$ and the $k$
of Lemma~\ref{lem:f-fill-bin} equal to $\hat{j}$ and the $t$ of Lemma~\ref{lem:f-fill-bin} equal to $t-t_0(\hat{j})$) it follows that for all $t\geq t_0(\hat{j})$, 
\[
f_{\hat{j}+1}^{\Gammavect,\zvect}(t) \ge \min\big\{\mu_{\hat{j}+1}^{\Lambdavect}, f_{\hat{j}+1}^{\Gammavect,\zvect}(t_0(\hat{j})) + (t-t_0(\hat{j}))(\Lambda_{\hat{j}+1} - \Gamma_{\hat{j}+1})p_{\hat{j}}\mu_{\hat{j}}^{\Lambdavect}\big\}.
\]

Now suppose $t \ge t_0(\hat{j}+1)$. By hypothesis,
$t-t_0(\hat{j}) \geq (4/\lambda) R/\min\{\Lambda_k: k\in [j]\}$. Also, $\Lambda_{\hat{j}+1} - \Gamma_{\hat{j}+1} \ge \Lambda_{\hat{j}+1}/2$, and $p_{\hat{j}}\mu_{\hat{j}}^{\Lambdavect} \ge \lambda/2$, so it follows that 
\[
f_{\hat{j}+1}^{\Gammavect,\zvect}(t) \ge \min\big\{\mu_{\hat{j}+1}^{\Lambdavect}, f_{\hat{j}+1}^{\Gammavect,\zvect}(t_0(\hat{j})) + \frac{R\Lambda_{\hat{j}+1}}{\min\{\Lambda_k\colon k \in [j]\}}\big\} \ge \min\big\{\mu_{\hat{j}+1}^{\Lambdavect}, f_{\hat{j}+1}^{\Gammavect,\zvect}(t_0(\hat{j})) + R\big\}.
\]
It therefore follows by~\eqref{eq:f-fill-all-bins} that $f_{\hat{j}+1}^{\Gammavect,\zvect}(t) \ge \mu_{\hat{j}+1}^{\Lambdavect}$, as required.
\end{proof}

Combining Lemma~\ref{lem:f-track-jammed} and Lemma~\ref{lem:f-fill-all-bins} we can now get a useful bound on how long it takes a $j$-jammed process to fill up. Lemma~\ref{lem:slowfill-j-jammed} underlies the initial filling of the volume process in our proof for the case of finitely many strongly exposed bins, and its  main idea also gets used in the case where we have infinitely many strongly exposed bins.

\begin{definition} \label{def:Tprot}
Given any positive integer~$j$ and any  send sequence~$\p$,
let 
\newsym{T-p-j}{time to fill $j$-jammed process, Definition~\ref{def:Tprot}}{T^{\p,j } = 80 j^2 \lfloor \sum_{x\in [j]} W_x \rfloor}.
\end{definition}

\begin{lemma}\label{lem:slowfill-j-jammed}
Fix a positive integer~$j$, a send sequence~$\p$ with $p_0=1$, and a $\lambda\in (0,1)$.
Fix a non-negative integer~$\tau$.
Let $Y$ be a $j$-jammed process with send sequence~$\p$, start time~$\tau$, birth rate~$\lambda$, and initial population size distribution $\Po(\zerovect)$.  
Then for all $t\geq T^{\p,j}$,
$\dist(\bvect^Y_{[j]}(t+\tau)) \gtrsim \Po(\lambda 9 W_1/10,\ldots, \lambda 9W_j/10)$.
\end{lemma}
\begin{proof}
We will show
\begin{itemize} 
\item For all $t \geq 0$, 
$b_{[j]}^Y(t + \tau) \sim \Po(f_{[j]}^{\zerovect,\zerovect}(t))$.
\item 
Also, for all $t \geq T^{\p,j}$ and all $k\in [j]$, 
$f_k^{\zerovect,\zerovect}(t) \geq 9 \lambda W_k/10$.
\end{itemize}
The lemma then follows from Observation~\ref{obs:dom-Poisson-tuples}.

The first item comes from Lemma~\ref{lem:f-track-jammed} 
with $\Delta=0$ and $\zvect = \zerovect$ (using Definition~\ref{def:f}). For the second item, we use Lemma~\ref{lem:f-fill-all-bins} with $\zvect = \zerovect$
and $\Gammavect = \zerovect$. We let $\Lambdavect$ be the vector in which every entry is $1/(10j)$.
The reason for this choice is that, for all $x\in [j]$, 
$\mu_x^{\Lambdavect} = \lambda W_x(1-1/(10j))^x$ so
$9 \lambda W_x/10 \leq \mu_x^{\Lambdavect} \leq \lambda W_x$.
We then choose $R := \lambda \sum_{x\in [j]} W_x$ so that the pre-condition 
$\sum_{k \in [j]} \max\{0, \mu_k^{\Lambdavect} - z_k\} \le R$
is met.
The lemma shows that for all 
$t \ge (4/\lambda)jR/\min\{\Lambda_k\colon k \in [j]\}$ and all $
k \in [j]$, $f_k^{\zerovect,\zerovect}(t) \ge \mu_k^{\Lambdavect}$ and we already showed that this is at least $9\lambda W_k/10$.
So to finish we need only show that
 $T^{\p,j} \geq (4/\lambda) j R/\min\{\Lambda_k: k\in [j]\}$,
 that is, taking $W = \sum_{x\in [j]} W_x$ and plugging in the definitions, we
 must show that 
$ 80 j^2 \lfloor  W \rfloor
\geq 40  j^2  W$ which holds since $W \geq W_1 \geq 1$.
\end{proof}
 
We conclude this section with some bounds 
(Lemma~\ref{lem:incf} and Corollary~\ref{cor:incf})
which follow from the definitions in this section (Definitions~\ref{def:f} and~\ref{def:mu}), and will be used later.

 \begin{lemma}\label{lem:incf}
Fix a positive integer~$j$, a send sequence~$\p$ with $p_0=1$, a $\lambda\in (0,1)$, a $j$-tuple $\zvect\in (\reals_{\geq 0})^j$, and 
$j$-tuple 
$\Gammavect \in [0,1]^j$. 
Fix $k\in [j]$.
Suppose that $p_1 z_1 \le (1-\Gamma_1)\lambda$ and that for all $x \in \{2,\ldots,k\}$,  $p_x z_x \le (1-\Gamma_x)p_{x-1}z_{x-1}$. Then, for all $t \ge 1$ and all $x \in [k]$, $f_x^{\Gammavect,\zvect}(t) \ge f_x^{\Gammavect,\zvect}(t-1)$.\end{lemma}
\begin{proof}
Write $f$ for $f^{\Gammavect,\zvect}$.
We will prove the following by induction on $t$, for all $t\geq 0$.
\begin{itemize}
\item For all $x\in [k]$, 
$p_xf_{x}(t) \le (1-\Gamma_x) p_{x-1} f_{x-1}(t)$.
\item For all $x\in \{0,\ldots k\}$, 
$f_x(t+1) \ge f_x(t)$.
\end{itemize}
The second item holds for $x=0$ by Definition~\ref{def:f}. 
To get the second item from the first one,  for  $x\in [k]$ using the definition,
 $f_x(t+1) =  
f_x(t) - p_x f_x(t) + (1-\Gamma_{x})p_{x-1}f_{x-1}(t)\geq f_x(t)$.

The first item holds for $t=0$ by hypothesis.
For the inductive step, suppose that
both items hold for some $t-1 \ge 0$. 
We will establish the first item for $x\geq 1$ as follows (the calculation uses  Definition~\ref{def:f} for the first equality and uses the two items of the inductive hypothesis for the two inequalities).
 
\begin{align*}
p_x&f_x(t) - (1-\Gamma_x)p_{x-1}f_{x-1}(t)\\
&= (1-p_x)p_x f_x(t-1) + p_xp_{x-1}(1-\Gamma_x)f_{x-1}(t-1) - (1-\Gamma_x)p_{x-1}f_{x-1}(t)\\
&= (1-p_x) p_x f_x(t-1) 
- (1-p_x) p_{x-1} (1-\Gamma_x)f_{x-1}(t-1)\\&\qquad\qquad\qquad\qquad 
+ p_{x-1} (1-\Gamma_x) f_{x-1}(t-1)
- (1-\Gamma_x)p_{x-1}f_{x-1}(t)\\
&\leq (1-p_x) ( p_xf_x(t-1) -
p_{x-1}(1-\Gamma_x)f_{x-1}(t-1)) \leq 0.
\end{align*} \end{proof}

\begin{corollary}\label{cor:incf}
Fix a positive integer~$j$,    a send sequence~$\p$ with $p_0=1$, a $\lambda\in (0,1)$,  and 
$j$-tuples 
$\Gammavect$ and $\Gammavect'$ in $[0,1]^j$
with $\Gammavect \leq \Gammavect'$.
Fix a $k\in [j]$
and a $j$-tuple $\zvect\in (\reals_{\geq 0})^j$,
where 
$z_\ell = \mu^{\Gammavect'}_\ell$  
for $\ell \in [k]$.
Then, for all $t \ge 1$ and all $x \in [k]$, 
$f_x^{\Gammavect,\zvect}(t) \ge f_x^{\Gammavect,\zvect}(t-1)$. 
\end{corollary}
\begin{proof}
To use Lemma~\ref{lem:incf}, we must show 
$p_1 \mu_1^{\Gammavect'} \le (1-\Gamma_1)\lambda$ and that for all $x \in \{2,\ldots,k\}$,  $p_x \mu^{\Gammavect'}_x \le (1-\Gamma_x)p_{x-1}\mu^{\Gamma'}_{x-1}$.
Both of these follow since for $x\in \{0,\ldots,k\}$,
$p_x \mu_x^{\Gammavect'} = \lambda \prod_{a=1}^x (1-\Gamma'_a)$.
\end{proof}

\section{Covered and exposed bins}\label{sec:covered}

Our proofs rely on various constants that depend on $\lambda\in (0,1)$. These are collected in Definition~\ref{def:constants}. Even though they depend on $\lambda$, we omit the explicit $\lambda$
from the names of these constants, in order to reduce clutter.  

\begin{definition}\label{def:constants}
Fix $\lambda \in (0,1)$.
Let 
$\newsym{C-noise}{$300$}{\Cnoise}= {300}$, 
$\newsym{C-SE}{$10\Cnoise \lceil\log(1/\lambda)\rceil$}{\CSE} = 10\Cnoise \lceil\log(1/\lambda)\rceil$,
$\newsym{C-Ups}{320}{\CUpsilon} = 320$,
$\newsym{c-NB}{$1/(30\CUpsilon)$}{\CNoiseBelow} = 1/({30}\CUpsilon)$,
$\newsym{C-F}{2}{\CFill} = 2$,
$\newsym{phi}{100}{\phi}=100$, and 
$\newsym{chi}{1000}{\chi} = 1000$.
\end{definition}

\begin{definition}\label{def:noise}
Fix $\lambda \in (0,1)$.
For every positive integer~$j$,
let $\newsym{LL(j)}{A function that is $\Theta(\log\log j)$ - Def~\ref{def:noise}}{\LL(j)} = \CSE \lceil \max\{1, \log\log j\} \rceil$
and let
$\newsym{L(j)}{A function that is $\Theta(\log j)$ - Def~\ref{def:noise}}{\L(j)} = \CNoise\lceil\max\{\log(1/\lambda),\CNoise \log j \}\rceil
$.
\end{definition}

\begin{definition}\label{def:Upsilon} Let $\p$ be a send sequence with $p_0=1$.
For all
positive integers~$j$ and all real numbers $W\geq 1$, let $\newsym{Upsilon-j-W}
{bins in $[j]$ with weight $\geq W$ - Def~\ref{def:Upsilon}}{\Upsilon_{j,\ge W}} = \{\ell \in [j-1]\colon W_\ell \ge W\}$. 
\end{definition}

From the definition~$[0]$ is the empty set, so $\Upsilon_{1,\geq W}$ is also empty.
Definition~\ref{def:Wtilde} defines a positive integer~$\Wtilde[j]$.
Note in the definition that
$ |\Upsilon_{j,\ge {1}}| \ge \min \{ j-1,C_\Upsilon \L(j)/\lambda\} $, so $\Wtilde[j]$ is well-defined, and is always at least~$1$ -- a fact that we use often, and is the basis of Observation~\ref{obs:Wtilde}.

\begin{definition}\label{def:Wtilde}
Fix $\lambda \in (0,1)$, a send sequence $\p$ with $p_0=1$, and a positive integer~$j$.
If $j=1$ then $\Wtilde[j] = 1$. Otherwise,
let
 $\newsym{W-tilde-j}{A particular $W$ - Def~\ref{def:Wtilde}}{\Wtilde[j]}$ 
be the smallest  integer~$W$ maximising
${W}|\Upsilon_{j,\ge {W}}|$
over all positive integers $W$ such that $ |\Upsilon_{j,\ge {W}}| \ge \min \{ j-1,C_\Upsilon \L(j)/\lambda\} $.  
\end{definition}

\begin{observation}\label{obs:Wtilde}
Fix $\lambda \in (0,1)$, a send sequence $\p$ with $p_0=1$, and a positive integer~$j$.
Then 
\[\Wtilde[j] \,|\Upsilon_{j,\geq \Wtilde[j]}|
\geq j-1.\]
\end{observation}

Definition~\ref{def:covered} gives the definition of covered and exposed bins. 
The idea is that if bin~$j$ is many-covered, then~$W_j$ is small compared to $\Wtilde[j] \,|\Upsilon_{j,\geq \Wtilde[j]}|$.  The
bins in $\Upsilon_{j,\geq \Wtilde[j]}$ can be expected to make noise for approximately $\exp( \Wtilde[j]\, |\Upsilon_{j,\geq \Wtilde[j]}|)$ steps -- and this  will be enough to fill bin~$j$. 
The ``many'' refers to the fact that 
$ |\Upsilon_{j,\ge {W}}| \ge \min \{ j-1,C_\Upsilon \L(j)/\lambda\} $, so many bins are providing noise while bin~$j$ fills.  
Unless $j$ is small (which won't be relevant for us), $ |\Upsilon_{j,\ge {W}}| \ge   C_\Upsilon \L(j)/\lambda$, so the expected noise will be $\Omega(\log j)$.
If bin~$j$ is
heavy-covered, then there are fairly many bins of high weight (weight at least~$j^2$) in $[j-1]$. Even though there aren't quite as many of these as the many-covered case (due to the factor of~$2$ in the denominator in (Prop~\ref{cov-prop-2}), whose importance we explain later), they are heavy (with weight at least~$j^2$), so they can still be expected to provide $\Omega(\log j)$ noise for the time that it takes bin bin~$j$ to fill. If bin~$j$ is exposed then we do expect the bins in $[j-1]$ to go quiet as bin~$j$ fills. However, if it is weakly exposed then property~(Prop~\ref{cov-prop-3}) ensures that there are enough bins in~$[j-1]$ with very high weight (weight~$j^{\chi}$) -- and these are very likely to stay full, and to allow the bins in $[j-1]$ to refill.

\begin{definition}\label{def:covered}
Fix $\lambda \in (0,1),$ a send sequence~$\p$ with $p_0=1$, and a positive integer~$j$.   
We consider three properties.
\begin{enumerate}[({Prop }1):]
\item\label{cov-prop-1} $\CFill W_j j^{\phi}\le 
\exp(\CNoiseBelow\lambda\Wtilde[j]\,|\Upsilon_{j,\ge\Wtilde[j]}|) $.
\item \label{cov-prop-2} $|\Upsilon_{j,\geq j^2}| \geq  \CUpsilon \L(j)/(2\lambda)$.
\item \label{cov-prop-3} $|\Upsilon_{j,\ge j^{\chi}}| \ge \CUpsilon \LL(j)/\lambda$.
\end{enumerate}

\begin{itemize}
\item  We say that bin $j$ is \emph{many-covered} if 
(Prop~\ref{cov-prop-1}) holds. 
\item We say that bin $j$ is \emph{heavy-covered} if (Prop~\ref{cov-prop-1}) is false and  (Prop~\ref{cov-prop-2}) holds.
\item We say that bin $j$ is \emph{weakly exposed} if
(Prop~\ref{cov-prop-1}) is false,
(Prop~\ref{cov-prop-2}) is false, and
(Prop~\ref{cov-prop-3}) holds.
\item We say that bin $j$ is \emph{strongly exposed} 
if (Prop~\ref{cov-prop-1}) is false, 
(Prop~\ref{cov-prop-2}) is false, and
(Prop~\ref{cov-prop-3}) is false.
     
\end{itemize}

We say that bin $j$ is \emph{covered} if it is many-covered or heavy-covered -- that is, 
if (Prop~\ref{cov-prop-1}) or
(Prop~\ref{cov-prop-2}) holds.
We say that bin $j$ is \emph{exposed} if it is weakly exposed or strongly exposed -- that is, if 
(Prop~\ref{cov-prop-1}) is false and
(Prop~\ref{cov-prop-2}) is false.

\end{definition}

We will use Observation~\ref{obs:SE-large-weight} which follows directly from   Definition~\ref{def:covered} -- its purpose is just to suppress some of the constants in (Prop~\ref{cov-prop-1}) to make the writing simpler.

Observation~\ref{obs:SE-large-weight}
and many subsequent lemmas have the hypothesis that $j$ is ``sufficiently large'', without specifying the details. Whenever we say this, we mean that $j$ is sufficiently large with respect to~$\lambda$ and~$\p$. Thus, it is sufficiently large with respect to the constants in Definition~\ref{def:constants}. Whenever we say ``sufficiently large'' (anywhere in this work) we always mean with respect to~$\lambda$ and~$\p$ and with respect to all constants.

\begin{observation}\label{obs:SE-large-weight}
Fix $\lambda \in (0,1)$, a send sequence~$\p$ with $p_0=1$, and a sufficiently large positive integer~$j$. If bin~$j$ is exposed then 
$W_j \geq \exp(\CNoiseBelow\lambda j/2)$.
\end{observation}
\begin{proof}
By Observation~\ref{obs:Wtilde}, $\Wtilde[j]\,|\Upsilon_{j,\ge\Wtilde[j]}|\geq j-1$. Thus, since (Prop~\ref{cov-prop-1}) is false,
\[W_j > \exp(\CNoiseBelow\lambda\Wtilde[j]\,|\Upsilon_{j,\ge\Wtilde[j]}|)/(\Cfill j^{\phi} ) \geq
\exp(\CNoiseBelow\lambda (j-1))/(\Cfill j^{\phi} ).\]
\end{proof}

\begin{remark}\label{rem:infWE}
Section~\ref{sec:context} contains an example of a send sequence with infinitely many strongly exposed bins.
Here is an example of a send sequence with infinitely many weakly exposed bins. 
As in Section~\ref{sec:context} let $C$ be sufficiently large
and consider the sequence $(a_0,a_1,\ldots)$ 
with $a_0=0$ and, for $k\geq 1$, $a_k = (2 C)^{a_{k-1}}$.
Let 
$R(k) = \lceil \CUpsilon \LL(a_{k})/\lambda \rceil $.
Let $k_0$ be a sufficiently large positive integer. 
Let $\p$ be the sequence
defined as follows.

 \[
p_j = \begin{cases}
C^{-j} & \mbox{ if $j= a_k$ for some $k\geq k_0$,}\\
a_k^{-\chi} & \mbox{
if $j\in \{a_k - R(k),\ldots,a_k-1\}$ for some $k\geq k_0$
,}\\
1/2 & \mbox{  otherwise}.
\end{cases}
\]   
We claim that for $k\geq k_0$, bin~$a_k$ is weakly exposed.
(Prop~\ref{cov-prop-3}) holds by construction.
To show that (Prop~\ref{cov-prop-2}) is false note that
the highest-weight bin in
$[a_k-R(k)-1]$ has weight $C^{a_{k-1}} < a_k < a_k^2$. 
For (Prop~\ref{cov-prop-1}) we note by Definition~\ref{def:Wtilde} that 
$\Wtilde[a_k] \leq (a_{k-1})^{\chi}$  
since there are at most $R(k)+1$ bins with weight 
larger than $(a_{k-1})^{\chi}$  in $[a_k-1]$. 
If $\Wtilde[a_k] \leq 2$
then $
\exp(\CNoiseBelow\lambda\Wtilde[a_k]\,|\Upsilon_{a_k,\ge\Wtilde[a_k]}|) 
<
\exp(\CNoiseBelow\lambda 
2 a_k )
<
C^{ a_k}
= W_{a_k} <
\CFill W_{a_k} a_k^{\phi}
$
so (Prop~\ref{cov-prop-1}) is false, as required. 
Otherwise, $2 < \Wtilde[a_k] \le (a_{k-1})^\chi$. For $k\geq k_0$, there are at most $k(R(k)+1)\leq 2kR(k)$ bins of weight greater than $2$ in $[a_k-1]$ so there are at most $2k R(k)$ bins with weight at least $
\Wtilde[a_k]$ in this range. Since $2R(k) \le (4C_{\Upsilon}\CSE\log\log a_k)/\lambda \le a_{k-1}$, $\exp(\CNoiseBelow\lambda\Wtilde[a_k]\,|\Upsilon_{a_k,\ge\Wtilde[a_k]}|) < \exp(\CNoiseBelow\lambda k(a_{k-1})^{\chi+1}) < C^{a_k} = W_{a_k} < C_FW_{a_k}a_k^\phi$ and again (Prop~\ref{cov-prop-1}) is false, as required. 
\end{remark}

The goal of
Section~\ref{sec:infinite}  is Corollary~\ref{cor:backoff-finite-SE-bins},
which proves lack of positive recurrence in
the case where there are infinitely many strongly exposed bins.
The proof in that section relies on the fact that every strongly exposed bin has a \emph{bottleneck} bin (Definition~\ref{def:bottleneck}) nearby -- this is Lemma~\ref{lem:nearbybottleneck}, which we prove here. The idea is that if bin~$j$ is a bottleneck bin, then the contribution to the sojourn time (between the starting empty state and the next empty state) when balls are in bin~$j$ is~$\omega(1)$, as a function of~$j$. This is a consequence of 
 Lemma~\ref{lem:bottleneck-geo-bound}, 
which implies a 
not-too-unlikely upper bound on the time that it takes a newborn ball to reach bin~$j$ in a $j$-jammed process when bin~$j$ is a bottleneck. 
Before formally defining a bottleneck bin,
we give (in 
Definition~\ref{def:forbottleneck}) a couple of definitions that extend Definitions~\ref{def:constants} and~\ref{def:noise}. 
Note that, since $\LL(i) = \Theta(\log \log i)$, it is also true that  $K(i) = \Theta(\log \log i)$.

\begin{definition}\label{def:forbottleneck} 
Fix $\lambda \in (0,1)$.
Let $\newsym{C-0}{
$2 \phi \CFill\CUpsilon\CNoise^2\log(6/\lambda)/(\CNoiseBelow\lambda)$}{C_0} = 2 \phi \CFill\CUpsilon\CNoise^2\log(6/\lambda)/(\CNoiseBelow\lambda)$ and, for every positive integer~$i$,
let
$\newsym{K-i}{$2\lceil \CUpsilon \LL(i) /\lambda \rceil$}{\K(i)} = 2\lceil \CUpsilon \LL(i) /\lambda \rceil$.
\end{definition}

\begin{definition}\label{def:bottleneck}
Fix $\lambda \in (0,1)$.  Let $\p$ be a send sequence with $p_0=1$ and let $i$ be a positive integer. We say that bin 
$i$ is a \emph{bottleneck} if all of the following hold.
\begin{enumerate}[(i)]
\item $W_i \ge \exp(i/\exp((\log\log i)^2))$.
\item For all $\ell \in [i-1]$,  $W_\ell \le W_i^{1/(2\K(i))}$.
\item There is a 
strongly exposed bin~$i^+$ which satisfies $i \leq i^+ \leq i \exp({(\log i)}^{1/2})$.
\end{enumerate}
\end{definition}

\begin{lemma}
\label{lem:nearbybottleneck}
Fix $\lambda \in (0,1)$.
Let $\p$ be a send sequence with $p_0=1$ such that, for all $j\geq 1$, $W_j \leq (6/\lambda)^j$.
Let $i$ be 
strongly exposed bin such that $i$ is  sufficiently large. Then there is a bottleneck bin $I$ with $I^+ = i$.
\end{lemma}

\begin{proof}
Let $\rr(i) = 4\K(i)\log(6/\lambda)/ (\lambda\CNoiseBelow)$.
We will find a sequence $i_1,\dots,i_k$ of positive integers with $i_1=i$ such that the following hold.
\begin{enumerate}
\item[(S1)] For all $\ell \in [k-1]$, $i_{\ell+1} \in [i_1/\rr(i_1)^{\ell},i_\ell)$.
\item[(S2)] For all $\ell \in [k-1]$, $W_{i_{\ell+1}} \ge W_{i_1}^{1/(2\K(i_1))^\ell}$.
\item[(S3)] Either $k = \K(i_1)$ or ($k < \K(i_1)$ 
and (S3b) holds).
\begin{enumerate}
\item [(S3b)] For all $y \in  [i_k-1]$, $W_y < W_{i_1}^{1/(2\K(i_1))^k}$.
\end{enumerate}
\end{enumerate}
 
\medskip\noindent \textbf{Proving the result from the sequence:} Suppose that we have such a sequence $i_1,\ldots,i_k$. We will show that we can take $I = i_k$ because bin~$i_k$ is a bottleneck with $i_k^+=i$. 

We first argue that $k < \K(i_1)$. Suppose 
for contradiction $k\geq K(i_1)$. Then by (S3), $k = \K(i_1)$. By (S2), for all $\ell \in [K(i_1)-1]$,
\begin{equation}\label{eq:i-ell-big}
W_{i_{\ell+1}} \ge W_{i_1}^{1/(2\K(i_1))^{\ell}} \ge W_{i_1}^{1/(2\K(i_1))^{k-1}} > W_{i_1}^{1/(2\K(i_1))^{\K(i_1)}}.
\end{equation}
Since $i_1$ is exposed, by Observation~\ref{obs:SE-large-weight}, $W_{i_1} \ge \exp({\CNoiseBelow \lambda i_1/2})$. Since $i_1$ is large, it follows 
from~\eqref{eq:i-ell-big} that for 
each $\ell \in [K(i_1)]$,
$W_{i_\ell} \geq   i_1^{\chi}$. 
However, since $i_1$ is strongly exposed, (Prop~\ref{cov-prop-3}) is false so  there are fewer than $\CUpsilon \LL(i_1)/\lambda < \K(i_1)-1$ bins in $[i_1-1]$ with weight at least $i_1^{\chi}$. This contradicts our initial assumption,  so $k < \K(i_1)$. Hence (S3b) holds.

We next show that $i_k$ satisfies all items in the definition of a bottleneck bin 
(Definition~\ref{def:bottleneck}) so bin~$i_k$ is a bottleneck.  To establish part (i) of Definition~\ref{def:bottleneck}, observe that 
(S2) with $\ell=k-1$ implies $W_{i_{k}} \ge W_{i_1}^{1/(2\K(i_1))^{k-1}}$. Since
 $k < \K(i_1)$, this implies $W_{i_k} \ge W_{i_1}^{1/(2\K(i_1))^{\K(i_1)}}$. 
Using again $W_{i_1} \ge \exp({\CNoiseBelow \lambda i_1/2})$,  it follows that
\[
    W_{i_k} \ge \exp\Big(\frac{\CNoiseBelow\lambda i_1}{2(2\K(i_1))^{\K(i_1)}}\Big) > \exp\Big(\frac{\CNoiseBelow\lambda i_1}{2\exp((\log\log i_1)^{3/2})}\Big) > \exp\Big(\frac{i_1}{\exp((\log\log i_1)^2)}\Big).
\]
Since (S1) gives $i_k < i_1$,
this implies $W_{i_k} > \exp(\frac{i_k}{\exp((\log\log i_k)^2)})$,
as required in Definition~\ref{def:bottleneck}(i).

We next show that $i_k$ satisfies   part (ii) of Definition~\ref{def:bottleneck}.
If $k=1$ then vacuously $W_{i_k}^{1/(2\K(i_k))}  \geq
W_{i_1}^{1/(2\K(i_1))^k}$. If $k>1$ we prove this re-using first
  $W_{i_k} \ge W_{i_1}^{1/(2\K(i_1))^{k-1}}$ and then  
$i_1 > i_k$, as follows. 
\[
W_{i_k}^{1/(2\K(i_k))} \ge 
W_{i_1}^{(1/(2\K(i_1))^{k-1})\cdot (1/(2 \K(i_k))}\geq
W_{i_1}^{1/(2\K(i_1))^k}.
\]  

By (S3b) it follows that for all $y \in [i_k-1]$, $W_y < W_{i_1}^{1/(2\K(i_1))^k} \leq W_{i_k}^{1/(2\K(i_k))}$, as required in part (ii) of Definition~\ref{def:bottleneck}.

Finally, to show that $i_k$ satisfies part (iii) of Definition~\ref{def:bottleneck}, we will take $i_k^+ = i_1$. Since $i_1$ is strongly exposed and $i_k < i_1$, it suffices to prove $i_1 \le i_k 
\exp({(\log i_k)}^{1/2})$. This is obvious if $k=1$. Otherwise, by (S1) with $\ell=k-1$, since $i_1$ is large, 
\[
i_k \ge i_1/\rr(i_1)^{k-1} > i_1/\rr(i_1)^{\K(i_1)} > i_1/\exp((\log\log i_1)^3).
\]
Since $i_1$ is large, it follows that $i_k \ge {i_1}^{1/2}$ and hence $\log \log i_k \ge (\log\log i_1)/2$; thus
\[
    i_1 < i_k\exp((\log\log i_1)^3) \le i_k\exp(8(\log\log i_k)^3) < i_k
  \exp({(\log i_k)}^{1/2}),
\]
as required
in part (iii) of Definition~\ref{def:bottleneck}. We have therefore shown that given a sequence $i_1,\dots,i_k$ satisfying (S1)--(S3), the result follows.

\medskip\noindent\textbf{Finding the sequence:} It remains to find such a sequence. We proceed iteratively. Suppose that we have found $i_1,\dots,i_x$ satisfying (S1) and (S2) for some $x < \K(i_1)$, where $i_1 = i$. (For $x=1$ this is immediate.) We prove that either $i_1,\dots,i_x$ satisfies (S3b) (in which case we can take $k=x$) or that we can find $i_{x+1}$ such that $i_1,\dots,i_{x+1}$ satisfies (S1) and (S2). We then find the required sequence by repeating this procedure at most $\K(i_1)-1$ times, since (S3) is satisfied automatically if $x=\K(i_1)$ and we can take $k=x$ in that case. To do this, we split into two cases.

\medskip\noindent\textbf{Case 1:} There is a positive integer $\ell < i_x$ with $W_\ell \ge W_{i_1}^{1/(2\K(i_1))^{x}}$. In this case we will take $i_{x+1} = \ell$, and we must prove that (S1) and (S2) hold for $i_1,\dots,i_{x+1}$. (S2) is immediate. To see (S1), observe that $W_\ell \le (6/\lambda)^{\ell}$ by hypothesis; thus
$(6/\lambda)^\ell \ge W_{i_1}^{1/(2\K(i_1))^x}$.
Reusing $W_{i_1} \ge \exp(\CNoiseBelow \lambda i_1/2)$, we get
$ (6/\lambda)^\ell \ge \exp({\CNoiseBelow \lambda i_1/(2 (2\K(i_1))^x)})
$.
After taking logarithms of both sides and dividing through by $\log(6/\lambda)$ 
 {using the fact that $\CNoiseBelow \leq 1$ and $2 \log(6/\lambda)/\lambda \ge 1$, }
we obtain
\[
    \ell \ge \frac{\CNoiseBelow\lambda}{2\log(6/\lambda)}\cdot \frac{i_1}{(2\K(i_1))^x} \ge \Big(\frac{\CNoiseBelow\lambda}{4\log(6/\lambda)\K(i_1)}\Big)^x i_1 = i_1/r(i_1)^x,
\]
as required by (S1).

\medskip\noindent\textbf{Case 2:} For all positive integers $\ell < i_x$, $W_\ell < W_{i_1}^{1/(2\K(i_1))^x}$. In this case we will take $k=x$; (S1) and (S2) hold by construction, and (S3b) holds by hypothesis.
\end{proof}

Corollary~\ref{cor:infinite-KMP-bins-2} follows easily from Lemma~\ref{lem:nearbybottleneck} and will be useful in Section~\ref{sec:infinite}.

\begin{corollary}\label{cor:infinite-KMP-bins-2}
Fix $\lambda \in (0,1)$.
Let $\p$ be a send sequence with $p_0=1$ such that, for all $j\geq 1$, $W_j \leq (6/\lambda)^j$. 
Suppose that there are infinitely many strongly exposed bins. Then there are infinitely many bottleneck bins.
\end{corollary}
\begin{proof}
Suppose that $i$ is sufficiently large and is strongly exposed.  Then by  Lemma~\ref{lem:nearbybottleneck} there is a bottleneck bin~$I$ with $I^+=i$.
By Definition~\ref{def:bottleneck}(iii), $I \leq i \leq I\exp((\log I)^{1/2})$.
Then if we consider another strongly exposed bin~$i'$ with 
$i' > I\exp((\log I)^{1/2})$ the same considerations will give another bottleneck bin~$I'$ that is strictly greater than~$I$
since we do not have 
$I \leq i' \leq I\exp((\log I)^{1/2})$. 
\end{proof}

We will conclude this section by proving 
Lemma~\ref{lem:bottleneck-geo-bound}, which 
which implies a 
not-too-unlikely upper bound on the time that it takes a newborn ball to reach bin~$j$ in a $j$-jammed process when bin~$j$ is a bottleneck. First we give a technical lemma, 
Lemma~\ref{lem:se-banding}, which will get used in the argument.

\begin{lemma}\label{lem:se-banding}
Fix $\lambda \in (0,1)$.
Let $\p$ be a send sequence with $p_0=1$ such that, for all $j\geq 1$, $W_j \leq (6/\lambda)^j$.  Let $x$ be sufficiently large  and suppose that bin $x$ is exposed.   Then for all positive integers~$w$,  
$|\Upsilon_{x,\ge w}| < C_0\max\{\log x, x/w\}$.
\end{lemma}
\begin{proof}
Since bin~$x$ is exposed, (Prop~\ref{cov-prop-1}) is false so, using $\CFill> 1$ and $\phi > 1$,
$ 
\exp(\CNoiseBelow\lambda\Wtilde[x]\,|\Upsilon_{x,\ge\Wtilde[x]}|) 
< \CFill W_{x} {x}^{\phi}
\leq (\CFill W_x x)^{\phi}
$. Taking logarithms,
$ 
\CNoiseBelow\lambda\Wtilde[x]\,|\Upsilon_{x,\ge\Wtilde[x]}|
<  \phi  \log(\CFill x W_x) \leq \CFill \phi \log(x W_x)
$.
From the definition of $\Wtilde[x]$ (Definition~\ref{def:Wtilde}),
if $|\Upsilon_{x,\ge w}| \geq \CUpsilon \L(x) /\lambda$, 
then $w | \Upsilon_{x,\geq w}|
\leq \Wtilde[x]\,|\Upsilon_{x,\ge\Wtilde[x]}|$.
Thus, 
\begin{equation}\label{eq:se-banding-max}
|\Upsilon_{x,\ge w}| < \max\{\CUpsilon \L(x)/\lambda, \CFill \phi\log(x W_{x})/(\CNoiseBelow\lambda w)\}.
\end{equation}

We bound each term of~\eqref{eq:se-banding-max} separately. Since $x$ is large, using Definition~\ref{def:noise},
\begin{align}\nonumber
\CUpsilon \L(x)/\lambda 
&= (\CUpsilon/\lambda) \cdot 
\CNoise\lceil\max\{\log(1/\lambda),\CNoise \log x \}\rceil 
\\ \label{eq:se-banding-term-1} 
&= (\CUpsilon\CNoise/\lambda) \lceil\CNoise\log x\rceil <
C_0\log x.
\end{align}
We now bound the second term of~\eqref{eq:se-banding-max}. Since $W_{x} \le (6/\lambda)^{x}$ by hypothesis and $x$ is large, 
\begin{equation}\label{eq:se-banding-term-2}
\CFill\phi\log(x W_{x})/(\CNoiseBelow\lambda w) \leq \CFill \phi \cdot 2x\log(6/\lambda)/(\CNoiseBelow\lambda w) < C_0x/w.
\end{equation}
The result follows by combining~\eqref{eq:se-banding-max}, \eqref{eq:se-banding-term-1} and~\eqref{eq:se-banding-term-2}.
\end{proof}

In the proof of Lemma~\ref{lem:bottleneck-geo-bound} we use the following 
generalisation of Definition~\ref{def:Upsilon}. We won't put this notation in the index, since it is only used in this proof.

\begin{definition}\label{def:gen-Upsilon} Let $\p$ be a send sequence with $p_0=1$.
For all
positive integers~$j$ and all real numbers $1\leq u_1 < u_2$,
let $ {\Upsilon_{j, [u_1,u_2) }} = \{\ell \in [j-1]\colon u_1 \leq W_\ell < u_2\}$. 
\end{definition}

\begin{lemma}\label{lem:bottleneck-geo-bound}
Fix $\lambda \in (0,1)$. Let $\p$ be a send sequence with $p_0=1$ such that, for all $x\geq 1$, $W_x \leq (6/\lambda)^x$. 
Let $i$ be a bottleneck bin such that $i$ is  sufficiently large. 
Let $G_0,\ldots, G_{i-1}$ be independent geometric variables with $\E(G_x) = W_x$. Then
$\pr\big(G_0 + \dots + G_{i-1} \le i\cdot   \exp((\log i)^{2/3}\big)   \ge 1/W_i^{1/2}$.
\end{lemma}
 
\begin{proof}
Let $i^+$ be the smallest integer satisfying (iii) in the definition of bottleneck (Definition~\ref{def:bottleneck}).
We divide the geometric variables into three classes. Let
    \begin{align*}
        \Sigma_1 &= \{j \in \{0,\ldots,i-1\} \colon W_j \le i^+/\log i^+\},\\
        \Sigma_2 &= \{j \in \{0,\ldots,i-1\} \colon i^+/\log i^+ < W_j \le (i^+)^{\chi}\},\\
\Sigma_3 &= \{j \in \{0,\ldots,i-1\} \colon W_j > (i^+)^{\chi}\},
    \end{align*}
where $\chi = 1000$ from Definition~\ref{def:constants}. 
Note that
\begin{equation}\label{eq:bottleneck-geo-bound-classes}
\sum_{j=0}^{i-1}G_j = \sum_{j \in \Sigma_1} G_j + \sum_{j \in \Sigma_2} G_j + \sum_{j \in \Sigma_3} G_j.
\end{equation}

We will bound each of these three sums above separately.

\noindent\textbf{Bounding the sum over $\boldsymbol{\Sigma_1}$:} We will prove that 
\begin{equation}\label{eq:bottleneck-geo-bound-class1-goal}
        \sum_{j \in \Sigma_1}W_j \le 8C_0i^+\log i^+
    \end{equation}
and then apply Markov's inequality. To do so, we first break these weights into bands; observe that every $j \in \Sigma_1$ satisfies $W_j \le i^+/\log i^+ < 2^{2\lceil \log i^+ \rceil}$ and hence (from Definition~\ref{def:gen-Upsilon}) $j$~lies in exactly one set $\Upsilon_{i,[2^k,2^{k+1})}$ with $0 \le k \le 2\lceil \log i^+ \rceil - 1$ and therefore (since $i\leq i^+$ by Definition~\ref{def:bottleneck}) $j$~lies in  
exactly one set $\Upsilon_{i^+,[2^k,2^{k+1})}$ with $0 \le k \le 2\lceil \log i^+ \rceil - 1$.
Thus
\begin{equation}\label{eq:bottleneck-geo-bound-class1-bands}
\sum_{j \in \Sigma_1} W_j = \sum_{k=0}^{2\lceil \log(i^+)\rceil - 1} \sum_{j \in \Sigma_1 \cap \Upsilon_{i^+,[2^k,2^{k+1})}} W_j \le \sum_{k=0}^{2\lceil \log(i^+)\rceil - 1} 2^{k+1}|\Sigma_1 \cap \Upsilon_{i^+,[2^k,2^{k+1})}|.
    \end{equation}

Since bin $i^+$ is exposed,  Lemma~\ref{lem:se-banding}  with $w=2^k$ 
and $x=i^+$ 
implies
    \[
        |\Sigma_1 \cap \Upsilon_{i^+,[2^k,2^{k+1})}| \le |\Upsilon_{i^+,\ge 2^k}| \le C_0\max\{\log i^+, i^+/2^k\}.
    \]
By the definition of $\Sigma_1$, for all $j \in \Sigma_1$, $W_j \le i^+/\log i^+$. Hence if $\Sigma_1 \cap \Upsilon_{i^+,[2^k,2^{k+1})} \ne \emptyset$, then $2^k \le i^+/\log i^+$ and hence $i^+/2^k \ge \log i^+$. It follows that for all $k$,  $
        |\Sigma_1 \cap \Upsilon_{i^+,[2^k,2^{k+1})}| \le C_0i^+/2^k$.
    Substituting this into~\eqref{eq:bottleneck-geo-bound-class1-bands} yields
    \[
        \sum_{j \in \Sigma_1}W_j \le \sum_{k=0}^{2\lceil \log i^+\rceil -1} 2^{k+1} \cdot C_0i^+/2^k \le 8C_0i^+\log i^+.
    \]
    By Markov's inequality, it follows that
    \begin{equation}\label{eq:bottleneck-geo-bound-class1-result}
        \pr\Big(\sum_{j \in \Sigma_1} G_j \le 16C_0i^+\log i^+\Big) \ge 1/2.
    \end{equation}

    \medskip\noindent \textbf{Bounding the sum over $\boldsymbol{\Sigma_2}$:} First observe that since $G_0,\dots,G_{i-1}$ are independent,
    \[
        \pr\Big(\sum_{j \in \Sigma_2}G_j = |\Sigma_2|\Big) = \prod_{j \in \Sigma_2}\pr(G_j = 1) = \prod_{j \in \Sigma_2} (1/W_j) \ge 1/(i^+)^{\chi|\Sigma_2|}.
    \]
    By Lemma~\ref{lem:se-banding} applied with $w=i^+/\log i^+$ and $x=i^+$, 
    $
        |\Sigma_2| \le |\Upsilon_{i^+,\ge i^+/\log i^+}| \le C_0\log i^+$. It follows that
    \begin{equation}\label{eq:bottleneck-geo-bound-class2-result}
        \pr\Big(\sum_{j \in \Sigma_2}G_j = |\Sigma_2|\Big) \ge \exp\big({-}\chi C_0(\log i^+)^2\big).
    \end{equation}

    \medskip\noindent \textbf{Bounding the sum over $\boldsymbol{\Sigma_3}$:} First observe that since 
    $G_0,\dots,G_{i-1}$ are independent,
    \begin{equation}\label{eq:bottleneck-geo-bound-class3-initial}
        \pr\Big(\sum_{j \in \Sigma_3}G_j = |\Sigma_3|\Big) = \prod_{j \in \Sigma_3}\pr(G_j = 1) = \prod_{j \in \Sigma_3} (1/W_j).
    \end{equation}
Since $i$ is a bottleneck (and $0\notin \Sigma_3$ since $W_0=1$) by Definition~\ref{def:bottleneck}(ii) it follows that
    \[
        \pr\Big(\sum_{j \in \Sigma_3}G_j = |\Sigma_3|\Big) \ge (1/W_i^{1/(2\K(i))})^{|\Sigma_3|}
    \]
Since $i^+$ is strongly exposed, (Prop~\ref{cov-prop-3}) is false, so 
$
        |\Sigma_3| \le |\Upsilon_{i^+, \ge (i^+)^{\chi}}| < \CUpsilon \LL(i^+)/\lambda$.
Substituting  this into the previous equation yields
    \[
        \pr\Big(\sum_{j \in \Sigma_3}G_j = |\Sigma_3|\Big) \ge 1/W_i^{\CUpsilon\LL(i^+)/(2\lambda \K(i))}.
    \]
Since $i$ is large, 
since $\K(i) = 2 \lceil C_\Upsilon \LL(i) / \lambda \rceil$ (Definition~\ref{def:forbottleneck}),
$\LL(j) = \CSE \lceil \log\log j\rceil $  for $j\geq 3$ (Definition~\ref{def:noise}) and
$i \leq i^+ \leq i \exp({(\log i)}^{1/2})$
(
Definition~\ref{def:bottleneck}(iii)),
    \[
        \frac{\CUpsilon\LL(i^+)}{2\lambda \K(i)} \le \frac{\LL(i^+)}{4\LL(i)} = \frac{\lceil \log \log i^+ \rceil}{4\lceil \log \log i\rceil} < \frac{1}{3}.
    \]
Thus,
    \begin{equation}\label{eq:bottleneck-geo-bound-class3-result}
        \pr\Big(\sum_{j \in \Sigma_3}G_j = |\Sigma_3|\Big) > 1/W_i^{1/3}.
    \end{equation}

 \medskip\noindent \textbf{Combining the bounds:}
    Since $G_0,\dots,G_{i-1}$ are independent, combining~\eqref{eq:bottleneck-geo-bound-classes}, \eqref{eq:bottleneck-geo-bound-class1-result}, \eqref{eq:bottleneck-geo-bound-class2-result}, and \eqref{eq:bottleneck-geo-bound-class3-result} yields
    \begin{align*}
\pr\Big(\sum_{j=0}^{i-1}G_j \le 16C_0i^+\log i^+ + |\Sigma_2| + |\Sigma_3|\Big) \ge        \frac{\exp(-\chi C_0(\log i^+)^2)}{2W_i^{1/3}}.
    \end{align*}
Since $i$ is sufficiently large and bin~$i$ is a bottleneck,  
$W_i \geq \exp(i/\exp((\log\log i)^2)) > \exp({i^{1/2}})$, $i^+ \le i
\exp({(\log i)}^{1/2})$ and $\log i^+ \le 2\log i$. Thus, since $|\Sigma_2| + |\Sigma_3| \leq i$,

\[
16C_0i^+\log i^+ + |\Sigma_2| + |\Sigma_3| 
\le 32C_0i(\log i)
\exp({(\log i)}^{1/2})  + i \le i \exp((\log i)^{2/3}).\]
Also, 
\[
\frac{\exp(-\chi C_0(\log i^+)^2)}{2W_i^{1/3}} \geq 
\frac{\exp(-4 \chi C_0(\log i)^2)}{2W_i^{1/3}} \geq \frac{1}{W_i^{1/2}},\]
    and the result follows.
\end{proof}

\section{Infinitely many strongly exposed bins}
\label{sec:infinite}

The purpose of this section is to prove Corollary~\ref{cor:backoff-finite-SE-bins} which shows that if $\p$ has infinitely many bins which are strongly exposed then the backoff process with send sequence~$\p$ is not positive recurrent.
In order to get this result we consider two auxiliary processes -- the two-stream process, which is defined in Section~\ref{sec:two-stream} and the under-backoff process, which is defined in Section~\ref{sec:under-backoff}. Both of these can be coupled with the backoff process (Definitions~\ref{def:standard-twostream} and \ref{def:standard_under}) so that the population of the backoff process stays above the population of the auxiliary process.
However, in each case the auxiliary process is a little bit easier to analyse.

Throughout the section we use the constants in Definition~\ref{def:bottleconstants}
and the events in Definition~\ref{def:bottle-events}. The $\bottleneck$ subscript stands for ``bottleneck''. 

\begin{definition}\label{def:bottleconstants}
Fix $\lambda \in (0,1)$. Let $\newsym{c-b}{Def \ref{def:bottleconstants}}{\cbottle = 1/(8\log(7/\lambda))}$ and let $\newsym{a-b}{Def~\ref{def:bottleconstants}}{\abottle = \cbottle\lambda/640}$.  
\end{definition}

Recall from Definition~\ref{def:expected-noise} that, given a generalised backoff process~$Z$ with send sequence~$\p$, and a set~$S$ of bins, the expected noise at time~$t$ is  
${\calN^Z_S(t)} = \sum_{i \in S}p_ib_i^Z(t-1)$

\begin{definition}\label{def:bottle-events} 
Fix $\lambda\in (0,1)$. 
Let $i$ be a sufficiently large integer and let 
$j = \lfloor \cbottle \log i \rfloor$. 
Let $Z$ be a generalised backoff process with send sequence $\p$ and start time~$\tau$.  
Let $\calE^{Z,\calN}_{\tau}$
and $\calE^{Z,\s}_{\tau}$ be trivial events that always occur.
For all $t\geq \tau+1$
define the following events.
Let $\calE^{Z,\calN}_t$ be the event that 
$\calN_{[j]}^Z(t)
\geq  \abottle\log i$ 
and 
let $\calE^{Z,\s}_t$ be the event that 
$\s^{{Z}}(t) \setminus \B_0^{{Z}}(t)$ is non-empty.   
Let $\calE^{Z,\calN}_{\leq t} = 
\bigcap_{t'\in \{\tau+1,\ldots,t\}} \calE^{Z,\calN}_{t'} $ and
let $\calE^{Z,\s}_{\leq t} = 
\bigcap_{t'\in \{\tau+1,\ldots,t\}} \calE^{Z,\s}_{t'}. $ 
\end{definition}

Section~\ref{sec:inf-j-jammed} 
gives a lower bound on the probability of the events in Definition~\ref{def:bottle-events} for a $j$-jammed process. Specifically, when bin $i$ is a bottleneck 
and $j = \Theta( \log(i))$ 
Lemma~\ref{lem:bottle-events-pair-jammed} shows that with probability at least $W_i^{-1/6}$,
a pair of independent $j$-jammed processes both 
satisfy these events 
for $i \lceil  \exp((\log i)^{2/3} ) \rceil$ steps.
Section~\ref{sec:two-stream}  (Lemma~\ref{lem:bottle-events-TwoStream}) 
exploits the nature of the events in Definition~\ref{def:bottle-events}. Specifically, when these events occur, it is possible 
to couple a two-stream process with the pair of $j$-jammed processes so that 
the two-stream population is at least as big as the $j$-jammed population (intuitively, this is because the events prevent any escapes). 
Hence, the probability bound from Lemma~\ref{lem:bottle-events-pair-jammed} also holds  for the two-stream process, showing that for a long period of time, both streams of the two-stream process are noisy. Section~\ref{sec:under-backoff} defines the under-backoff process (which has three streams) and shows how to couple   the two-stream process 
with the first two of these streams so that the two-stream  population stays below the population of the relevant streams of the under-backoff process. This guarantees (Lemma~\ref{lem:FinalLB-underbackoff}) that two of the streams of the under-backoff process are likely to be noisy for a fairly long period.  This (Lemma~\ref{lem:GetIntoi}) makes it fairly likely that a ball gets all the way to bin~$i$ in the third stream of the under-backoff process. This makes it easy to prove that the expected sojourn time (from an empty state to an empty state) in the under-backoff process is an increasing function of~$i$ (Lemma~\ref{lem:se-sojourn}) which implies, in the presence of infinitely many bottlenecks, that the under-backoff process is not positive recurrent.
A final coupling (Definition~\ref{def:standard_under}) establishes the same result (Corollary~\ref{cor:backoff-finite-SE-bins}) for backoff.

\subsection{\texorpdfstring{$j$}{j}-jammed processes and bottlenecks}
\label{sec:inf-j-jammed}

Lemma~\ref{lem:bottle-events-jammed}
gives a lower bound on the probability of the events in Definition~\ref{def:bottle-events} for a $j$-jammed process.
The lower bounds are expressed as a function of $i$ where $j = \Theta( \log(i))$.
We will use the lemma in the situation where bin~$i$ is a bottleneck and where the start state comes from Lemma~\ref{lem:slowfill-j-jammed}. Lemma~\ref{lem:bottle-events-pair-jammed} expresses the failure probability from Lemma~\ref{lem:bottle-events-jammed} as a function of $W_i$ for the case where bin~$i$ is a bottleneck.
For convenient application later, this is written for a pair of independent $j$-jammed processes.
  
\begin{lemma}
\label{lem:bottle-events-jammed}
Fix $\lambda \in (0,1/2)$. Let $\p$ be a send sequence with $p_0=1$. 
Let $i$ be a sufficiently large integer   
and let $j= \lfloor \cbottle \log i \rfloor$.  
Fix a non-negative integer~$\tau$. 
Let $\zvect\in (\reals_{\geq 0})^j$ be a $j$-tuple such that, for all $x\in [j]$, $z_x \geq \lambda W_x/4$. 
Let $Y$ be a \emph{$j$-jammed} process with send sequence~$\p$, start time~$\tau$, birth rate~$\lambda$, and initial population size distribution~$\Po(\zvect)$.
Then for all integers~$t>\tau$,
\begin{itemize}
\item $\pr(  {\calE^{Y,\s}_{t}} \mid
\calE^{Y,\calN}_{\leq t} \cap \calE^{Y,\s}_{\leq t-1} ) \geq 1- 1/i^{\abottle}  $, and

\item 
$\pr(   {\calE^{Y,\calN}_t} \mid
\calE^{Y,\calN}_{\leq t-1} \cap \calE^{Y,\s}_{\leq t-1} ) \geq 1-
1/i^{\abottle}$.

\end{itemize}

\end{lemma}
\begin{proof}
We first establish~\eqref{eq:add-t} for all $t>\tau$. Let $t>\tau$ be an integer.
First we show that for $t'\in \{\tau+1,\ldots,t\}$ the indicator variables of the events $\calE^{Y,\calN}_{t'}$ and $\calE^{Y,\s}_{t'}$ are increasing functions of the variables in 
$\{ b_\traje^Y \mid \traje \in \calT^Y(t)\}$ from Definition~\ref{def:trajectory}. 
Consider the event $\calE^{Y,\calN}_{t'}$.
This is the event that 
$\sum_{x=1}^j p_x b_x^{{Y}}(t'-1)  \geq  \abottle\log i$ 
which is the event that 
$\sum_{x=1}^j \sum_{\traje \in \calT^{Y}(x,t'-1)}b_\traje^{Y} p_x \geq \abottle \log i$. The result now follows from Observation~\ref{obs:latertraj}
which says that, for every $\traje \in \calT^Y(t'-1)$,
there is a set $S\subseteq \calT^Y(t)$
such that $\b_\traje^Y = 
\cup_{\traje'\in S} \b_{\traje'}^Y$. 
Finally, consider the event $\calE^{Y,\s}_{t'}$. 
This is
the event that 
$\s_1^{Y}(t') \cup \cdots \cup \s_j^{Y}(t') $ is non-empty. 
This is the event that there is a   trajectory $\traje \in \calT^{Y}(t')$ with 
$\traje_{t'-1} \in [j]$ and $\traje_{t'} \neq \traje_{t'-1}$
such that   $|\b_{\traje}^{Y}|>0$. The result again follows from Observation~\ref{obs:latertraj}. 

\smallskip We have showed that the event $ \calE^{Y,\calN}_{\leq t-1}  \cap \calE^{Y,\s}_{\leq t-1}$ is an increasing function of the variables in 
$\{ b_\traje^Y \mid \traje \in \calT^Y(t)\}$.
So by Lemma~\ref{lem:f-track-jammed}, for all $t>\tau$,
\begin{equation}\label{eq:add-t}
\pr(   {\calE^{Y,\calN}_t} \mid
\calE^{Y,\calN}_{\leq t-1} \cap \calE^{Y,\s}_{\leq t-1} ) 
\geq \pr( {\calE^{Y,\calN}_t}).\end{equation}
Making this substitution in the item
that we wish to prove (and re-writing
the first item), it suffices to show that for all integers $t> \tau$,

\begin{itemize}
\item $\pr(  {\calE^{Y,\s}_{t}} \mid
 {\calE^{Y,\calN}_t} \cap  \calE^{Y,\calN}_{\leq t-1}  \cap \calE^{Y,\s}_{\leq t-1}) \geq 1 - 1/i^{\abottle}  $.
\item 
$\pr(   {\calE^{Y,\calN}_t}   ) \geq 
1-
1/i^{\abottle}$.  
\end{itemize}

For the first item,
note that any outcome of steps $\tau+1,\ldots,t-1$ that satisfies the conditioning satisfies   $\calE^{Y,\calN}_t$ 
so $\sum_{x=1}^j p_x b_x^{Y}(t-1) \geq \abottle \log i$.
The probability that $\calE^{Y,\s}_t$ does not occur (subject to the conditioning) is at most $ 
\prod_{x=1}^j (1-p_x)^{b_x^{{Y}}(t-1)} $ and this probability
(subject also to the conditioning) is at most
$
 \exp(- \sum_{x=1}^j p_x b_x^{{Y}}(t-1)
 ) \leq \exp(- \abottle \log i) = 1/i^{ \abottle}$.
(Here we are essentially using the law of total probability,
conditioning on possible values of $\bbar_{[j]}^{Y}(t-1)$.)

For the 
second
item 
we want a lower bound on the probability that
$
\sum_{x=1}^j p_x b_x^{{Y}}(t-1)  \geq  \abottle\log i$ 
where (from   Lemma~\ref{lem:f-track-jammed} with $\Delta=0$) 
$b_{[j]}^{Y}(t-1) \sim 
\Po(f^{\zerovect,\zvect}_{[j]}(t-1-\tau))$.

Now let $\Gammavect'$ be the $j$-tuple that is all~$0$ except $\Gamma'_1 = 3/4$. 
We wish to apply Corollary~\ref{cor:incf} with $\Gammavect = \zerovect$ and $k=j$. This requires
$z_x = \mu^{\Gammavect'}_x$ for $x\in [j]$, which follows from the definition of~$\zvect$ and~$\Gammavect'$.
Corollary~\ref{cor:incf} shows that $f_x^{\zerovect,\zvect}(t-1-\tau) \geq f_x^{\zerovect,\zvect}(0) = z_x \geq \lambda W_x/4$.

Now apply Lemma~\ref{lem:poisson-noisy}
with $c=\lambda/4$, $W=1$, $S=[j]$
and $P_k = b_k^Y(t-1)$.
Fix $\delta = 1-1/80$.
Note that $\abottle \log i \leq 
(1-\delta) c j $ 
so 
$\Pr(\sum_{x=1}^j p_x b_x^{{Y}}(t-1)  
<  \abottle\log i)
\leq 
\Pr( \sum_{x=1}^j p_x b_x^{{Y}}(t-1)  
<  (1-\delta) c j)$.
By the lemma, this is at most
$\exp(-\delta^2 c j/2)$. 
Since $\delta^2 \geq 39/40$, this is at most
$
\exp (-\lambda (39/320) j ) \leq
\exp(- \lambda (39/640) \cbottle \log i)
= \exp(- 39 \abottle \log i)
\leq \exp(- \abottle \log i)
  $.
\end{proof}

Lemma~\ref{lem:bottle-events-pair-jammed} expresses the failure probability from Lemma~\ref{lem:bottle-events-jammed} as a function of $W_i$ for the case where bin~$i$ is a bottleneck. As noted,  this is written for a pair of independent $j$-jammed processes.

\begin{lemma}
\label{lem:bottle-events-pair-jammed}
Fix $\lambda \in (0,1/2)$. Let $\p$ be a send sequence with $p_0=1$. 
Let $i$ be a sufficiently large integer   
and suppose that bin~$i$ is a bottleneck. Let $j= \lfloor \cbottle \log i \rfloor$. 
Fix a non-negative integer~$\tau$. 
Let $\zvect\in (\reals_{\geq 0})^j$ be a $j$-tuple such that, for all $x\in [j]$, $z_x \geq \lambda W_x/4$. 
Let $Y^A$ and $Y^B$ be independent
$j$-jammed processes with send sequence~$\p$, start time~$\tau$, birth rate~$\lambda$, and initial population size distribution~$\Po(\zvect)$. 
Let $\TEnd = \tau + i \lceil  \exp((\log i)^{2/3} ) \rceil$.
Then $\Pr(
\calE^{Y^A,\calN}_{\leq \TEnd} \cap \calE^{Y^B,\calN}_{\leq \TEnd} \cap 
\calE^{Y^A,\s}_{\leq \TEnd} \cap \calE^{Y^B,\s}_{\leq \TEnd}
) \geq W_i^{-1/6}$.

\end{lemma}

\begin{proof}
Let $I = i \lceil  \exp((\log i)^{2/3} ) \rceil$.
By Lemma~\ref{lem:bottle-events-jammed},
\[\Pr(
\calE^{Y^A,\calN}_{\leq \TEnd} \cap \calE^{Y^B,\calN}_{\leq \TEnd} \cap 
\calE^{Y^A,\s}_{\leq \TEnd} \cap \calE^{Y^B,\s}_{\leq \TEnd}
)  \] is at least
 $   (1-1/i^{\abottle})^{4I}
\geq \exp(-2(4I)/i^{\abottle} )$.
Since $i$ is a bottleneck,
$W_i^{-1/6} \leq \exp(-i/(6\exp((\log\log i)^2)))$.
So we just need to show
that $8 I / i^{\abottle} \leq i/
(6\exp((\log\log i)^2))$.  Plugging in the definition of~$I$, we need to show that
\[ 
8 i^{1-\abottle}  \lceil  \exp((\log i)^{2/3} ) \rceil \leq i /
(6\exp((\log\log i)^2)).\]
This holds since $i$ is sufficiently large.
\end{proof}

\subsection{Two-stream processes}
\label{sec:two-stream}

In this section we extend the results of Lemma~\ref{lem:bottle-events-pair-jammed} to a pair of generalised backoff processes called a \emph{two-stream process}. The two-stream process was introduced in~\cite{GL-oldcontention} 
as a useful device for dominating backoff processes.

\begin{definition}\label{def:two-stream}
Fix $\lambda \in (0,1)$,
a \emph{send sequence}  
$\p$ with $p_0=1$, and a non-negative integer~$\tau$.
A \emph{two-stream process} 
with birth rate~$\lambda$, send sequence~$\p$ and start time~$\tau$ 
is a pair $Z=(Z^A,Z^B)$ of generalised backoff processes with 
send sequence~$\p$, start time~$\tau$, 
birth distribution
 $\calD =  \Po(\lambda/2)$
and some initial populations $\arrvect^{Z^A}(\tau)$ and $\arrvect^{Z^B}(\tau)$. 
We refer to $Z^A$ and $Z^B$ as the \emph{streams} of~$Z$. 
The streams
evolve independently  
except 
that $\e_j^{Z^A}(t)$ and $\e_j^{Z^B}(t)$ are defined as follows for $t> \tau$. 
\begin{itemize}
    \item  
If there are distinct $C$ and $D$ in $\{A,B\}$ such that
$s^{Z^C}(t) = 0$ and $s^{Z^{D}}(t)>0$ then for all $j$, $\e_j^{Z^C}(t)=\emptyset$. In this case, a bin~$J^D$ is chosen with probability proportional to $s^{Z^{D}}_{J^D}(t)/s^{Z^{D}}(t)$ and $\e_{J^D}^{Z^D}(t)$
contains a single uniformly-chosen ball from 
$\s^{Z^D}_{J^D}$.
Also, 
for all $j\neq J^D$, $\e_j^{Z^D}(t)=\emptyset$. 
\item Otherwise, for   all $j$, $\e_j^{Z^A}(t)= \e_j^{Z^B}(t)=\emptyset$.
\end{itemize}
There are no end-of-step arrivals after step~$\tau$. That is, for all $D\in \{A,B\}$, all $t> \tau$, and all positive integers~$j$, $\arr_j^{Z^D} = \emptyset$.
\end{definition}

The point of the two-stream process is that sends in one stream are only blocked by the other stream (although randomness is used to ensure that at most one ball escapes from a given stream in a step).  This provides a little more independence than a backoff process, but it is nevertheless easy to couple a backoff process~$X$ with a two-stream process~$Z$ so that every ball in~$Z$ is in~$X$ -- the standard coupling for this is given in Definition~\ref{def:standard-twostream}.

 \begin{definition}
\label{def:standard-twostream} 
Fix $\lambda \in (0,1)$
and a \emph{send sequence}  
$\p$ with $p_0=1$.
Let $Z=(Z^A,Z^B)$ be the two-stream process with birth rate~$\lambda$, send sequence~$\p$ and start time~$0$  
where for all $D\in \{A,B\}$ and all
positive integers~$k$,
$\arr_k^{Z^D}(\tau) = \emptyset$. 
Let $X$ be the backoff process with birth distribution~$\Po(\lambda)$, send sequence~$\p$, and cohort set $\calC=\{A,B\}$ (Definition~\ref{def:backoff-cohorts}).
The standard coupling of~$Z$ and~$X$ has the following invariants for every $i\geq 1$ and $t\geq 1$.
\begin{description}  
\item  [Inv-$b_i(t)$:] For each stream $D\in \{A,B\}$, $\b_i^{Z^D}(t-1) \subseteq \b_i^{X^D}(t-1)$.
\item  [Inv-$B_i(t)$:] For each stream $D\in \{A,B\}$, $\B_i^{Z^D}(t) \subseteq \B_i^{X^D}(t)$.
\item [Inv-$s_i(t)$:] 
For each stream $D\in \{A,B\}$, $\s_i^{Z^D}(t) \subseteq \s_i^{X^D}(t)$.
\end{description}
Initially, Inv-$b_i(1)$ holds for every $i\geq 1$ since all of these sets are empty.

For any $t$ assuming Inv-$b_i(t)$, the coupling at step~$t$ proceeds as follows.
For each $D\in \{A,B\}$, $\n^{Z^D}(t) = \n^{X^D}(t)$, so for all $i\geq 1$, 
Inv-$B_i(t)$ holds.
For every   
$D\in \{A,B\}$ and
$i\geq 1$, $\S^{X^D}_i(t)$ is formed by including each ball in 
$\B_i^{X^D}(t) \setminus \B_i^{Z^D}(t)$ 
independently with probability~$p_i$.
Then $\s^{X^D}_i(t) = \S^{X^D}_i(t) \cup \s^{Z^D}_i(t)$, establishing Inv-$s_i(t)$.

The assignments to $\e_i^{X^D}(t)$ follow deterministically, using  Definition~\ref{def:backoffasgen}
and the assignments to $\e_i^{Z^D}(t)$  are chosen  following Definition~\ref{def:two-stream}.

Since $\s_i^{X^D}(t) \cap \B_i^{Z^D}(t) = \s_i^{Z^D}(t)$, Inv-$b_i(t+1)$  follows from  
$ \s_i^{Z^D}(t) \cap \e_i^{X^D}(t) \subseteq \e_i^{Z^D}(t)$, which holds:
If a ball~$\beta$ is in $ \s_i^{Z^D}(t) \cap \e_i^{X^D}(t) $ then
$\s^{X^A}(t) \cup \s^{X^B}(t)= \{\beta\}$ so
$\s^{Z^A}(t) \cup \s^{Z^B}(t) = \{\beta\}$ so 
$\e_i^{Z^D}(t) = \{\beta\}$.

\end{definition}

In order to ease the notation for dealing with the initial populations in the two streams of a two-stream process we generalise Definition~\ref{def:PoTuple} in the following natural way. The notation in Definition~\ref{def:genPotuple} is local to this section, so it is not included in the general index.

\begin{definition}\label{def:genPotuple}.
For every $j$-tuple $\zvect$, we use the notation $\zvect,\zvect$ to denote a tuple of length $2j$ where the first $j$ entries are the entries of $\zvect$ and similarly the last $j$ entries are the entries of $\zvect$. Thus, if $A$ and $B$ are random variables with support $(\integers_{\geq 0})^j$ then $\dist(A,B) \gtrsim \Po(\zvect,\zvect)$ indicates that the distribution of the pair $(A,B)$ is dominated below by an independent product of $2j$ Poisson distributions, where, for $k\in [j]$,  the $k$'th and $(j+k)$th Poisson distributions have mean~$z_k$. 
\end{definition}

Recall from Definition~\ref{def:Tprot} 
that $T^{\p,j } = 80 j^2 \lfloor \sum_{x\in [j]} W_x \rfloor$. Lemma~\ref{lem:slow-fillstreams}
extends the filling for $j$-jammed processes (Lemma~\ref{lem:slowfill-j-jammed}) to two-stream processes.

\begin{lemma}\label{lem:slow-fillstreams}
Fix $\lambda \in (0,1/2)$ and a send sequence $\p$ with $p_0=1$. Let $Z=(Z^A,Z^B)$ be a two-stream process with birth rate~$2\lambda$, send sequence~$\p$, and start time~$\tau$,
where for all $D\in \{A,B\}$ and all
positive integers~$k$,
$\arr_k^{Z^D}(\tau) = \emptyset$. 
Fix a positive integer~$j$ and let  $\zvect$ be the $j$-tuple 
with  $z_x = \lambda 9 W_x/20$ for all $x\in [j]$.
Then there is an event $\Edom$ which depends only on the evolution of $Z$ during steps $\{\tau+1,\ldots, \tau+T^{\p,j}\}$ such that the following hold.

\begin{itemize}
\item $\Pr(\Edom) \geq (1-\exp({-\lambda/2}))^{2 T^{\p,j}}$,
\item If $\Edom$ occurs then
for each $t\in \{\tau+1, \ldots, \tau+ T^{\p,j}\}$, $\s^{Z^A}(t)$ and $\s^{Z^B}(t)$ are both non-empty,   and

\item 
$\dist( (\bbar^{Z^{A}}_{[j]}(\tau + T^{\p,j}),\bbar^{Z^{B}}_{[j]}(\tau + T^{\p,j})) \mid \Edom)
\gtrsim   \Po(\zvect,\zvect)$.  
\end{itemize}

\end{lemma}

\begin{proof}
For each stream $D\in \{A,B\}$ and
each $t>\tau$,  partition the
set $\n^{Z^D}(t)$ 
of balls that are born at time~$t$
into~$Z^D$ into two sub-streams 
$\n^{Z^{D'}}(t)$ and $\n^{X^{D''}}(t)$ 
by choosing a sub-stream independently and uniformly at random, for each ball. The sub-streams have no effect on the evolution of the two-stream process.
Let $\Edom$ be the event that, 
for each $D\in \{A,B\}$ and all $t\in \{\tau+1, \ldots, \tau+ T^{\p,j}\}$,
$n^{Z^{D''}}(t) >0$. The second item is immediate from the definition of $\Edom$.
For each $D\in \{A,B\}$, $\n^{Z^{D''}}$
is a Poisson random variable with mean~$\lambda/2$ so the   first item holds.
For the third item, we will show that
\[\dist( (\bbar^{Z^{A'}}_{[j]}(\tau + T^{\p,j}),\bbar^{Z^{B'}}_{[j]}(\tau + T^{\p,j})) \mid \Edom)
\gtrsim   \Po(\zvect,\zvect).\]

Given $\Edom$, 
each
$Z^{D'}$ behaves as an independent $j$-jammed process with birth rate $\lambda/2$ during steps $\tau+1,\ldots,\tau+T^{\p,j}$,
except that balls that send from bin~$j$ don't escape (and there may be balls in bins $i>j$). By Lemma~\ref{lem:slowfill-j-jammed}, and Observation~\ref{obs:dom-Poisson-tuples},
$\dist(\bbar^{Z^{D'}}_{[j]}(\tau + T^{\p,j})) 
\gtrsim \Po(\zvect)$.
\end{proof}

Lemma~\ref{lem:bottle-events-TwoStream} couples a two-stream process with a pair of $j$-jammed processes to extend the result of Lemma~\ref{lem:bottle-events-pair-jammed} to a two-stream process.

\begin{lemma}

\label{lem:bottle-events-TwoStream}
Fix $\lambda \in (0,1/2)$ and a send sequence $\p$ with $p_0=1$.  
Let $i$ be a sufficiently large integer and suppose that bin~$i$ is a bottleneck. 
Let $j= \lfloor \cbottle \log i \rfloor$.  
Let $\zvect\in (\reals_{\geq 0})^j$ be a $j$-tuple such that, for all $x\in [j]$, $z_x \geq \lambda W_x/4$. 
 Let $Z=(Z^A,Z^B)$ be a two-stream process with birth rate~$2\lambda$, send sequence $\p$, and start time~$\tau$ where 
$\dist(\romanarrvect^{Z^A}_{[j]}(\tau) ,\romanarrvect^{Z^B}_{[j]}(\tau))\gtrsim \Po(\zvect,\zvect)$.
(We do not make any assumptions about
$\romanarr_k^{Z^A}(\tau)$ or $\romanarr_k^{Z^B}(\tau)$ for $k>j$.)
Let $\TEnd' = \tau + i \lceil  \exp((\log i)^{2/3} ) \rceil$.
Then $\Pr(
\calE^{Z^A,\calN}_{\leq \TEnd'} \cap \calE^{Z^B,\calN}_{\leq \TEnd'} \cap 
\calE^{Z^A,\s}_{\leq \TEnd'} \cap \calE^{Z^B,\s}_{\leq \TEnd'}
) \geq W_i^{-1/6} $.
\end{lemma}

\begin{proof}

Let $\calD = \dist(\romanarrvect^{Z^A}_{[j]}(\tau) ,\romanarrvect^{Z^B}_{[j]}(\tau))$
and $\calD' = \Po(\zvect,\zvect)$ 
so that $\calD \gtrsim \calD'$.
Recall from Definition ~\ref{def:phi} 
that $\phi_{\calD,\calD'}$ maps  
$(\integers_{\geq 0})^{2j}$
to $(\integers_{\geq 0})^{2j}$
so that, for $j$-tuples $\avect$ and $\avect'$,
$\phi_{\calD,\calD'}(\avect,\avect') \leq \avect,\avect'$.
Also, when $(\avect,\avect') \sim \calD$,
the random variable $\phi_{\calD,\calD'}(\avect,\avect')$ has distribution $\Po(\zvect,\zvect)$.

Draw $\romanarrvect^{Z^A}_{[j]}(\tau) ,\romanarrvect^{Z^B}_{[j]}(\tau)$ from $\calD$. Let $\arrvect^{Z^A}_{[j]}(\tau) ,\arrvect^{Z^B}_{[j]}(\tau)$ be the initial populations of~$Z^A$ and~$Z^B$
with, for each $D\in \{A,B\}$ and $k\in [j]$,
$|\arr^{Z^D}_{k}(\tau) | = 
\romanarr^{Z^D}_{k}(\tau)$.
Let $(\yvect,\yvect') = \phi_{\calD,\calD'}(
\romanarrvect^{Z^A}_{[j]}(\tau) ,\romanarrvect^{Z^B}_{[j]}(\tau)) $ and 
let $Y^A$ be a $j$-jammed process with birth rate~$\lambda$, send sequence $\p$, start time~$\tau$, and initial population size~$\yvect$. 
Specifically, for each $k\in [j]$,
let $\arr_k^{Y^A}(\tau)$ be the lexicographically least $y_k$ 
balls in $\arr_k^{Z_A}(\tau)$. 
$Y^A$ has no initial population in bins $k>j$.
Similarly, 
let $Y^B$ be an independent $j$-jammed process with birth rate~$\lambda$, send sequence $\p$, start time~$\tau$, and initial population size~$\yvect'$. Specifically, for each $k\in [j]$,
let $\arr_k^{Y^B}(\tau)$ be the lexicographically least $y_k'$ 
balls in $\arr_k^{Z_B}(\tau)$. $Y^B$ has no initial population in bins $k>j$.
By construction, the initial population size distribution of~$Y^A$ and~$Y^B$ is~$\Po(\zvect)$.
 
For all $t \ge 0$, let $\calE_{\le t}$ be the event 
$\calE^{Y^A,\calN}_{\leq t} \cap \calE^{Y^B,\calN}_{\leq t} \cap 
\calE^{Y^A,\s}_{\leq t} \cap \calE^{Y^B,\s}_{\leq t}$.
By Lemma~\ref{lem:bottle-events-pair-jammed}, 
$\Pr(\calE_{\le \TEnd'}) \geq W_i^{-1/6}$.

We will couple $Y^A$, $Y^B$ and $Z$  
for steps $t\in \{\tau+1,\ldots,\TEnd'\}$
such that if $\calE_{\le t}$ occurs then the following invariants hold
for every   $k\geq 1$.
\begin{description}  
\item  [Inv-$b_k(t)$:] For all $D\in \{A,B\}$, $\b_k^{Y^D}(t-1) \subseteq \b_k^{Z^D}(t-1)$.
\item  [Inv-$B_k(t)$:] For all $D\in \{A,B\}$, $\B_k^{Y^D}(t) \subseteq \B_k^{Z^D}(t)$.
\item [Inv-$s_k(t)$:] For all $D\in \{A,B\}$, $\s_k^{Y^D}(t) \subseteq \s_k^{Z^D}(t)$.
\end{description}

From the definition of the events (Definition~\ref{def:bottle-events}),
these invariants and~$\calE_{\le \TEnd'}$ imply 
$\calE^{Z^A,\calN}_{\leq \TEnd'} \cap \calE^{Z^B,\calN}_{\leq \TEnd'} \cap 
\calE^{Z^A,\s}_{\leq \TEnd'} \cap \calE^{Z^B,\s}_{\leq \TEnd'}$, giving the lemma.

The coupling follows the idea of the standard coupling (Definition~\ref{def:standard-couple}).
For all $k\geq 1$, Inv-$b_k(\tau+1)$ holds by construction.

Given any $t\in \{\tau+1,\ldots,T\}$, and 
for each $D\in \{A,B\}$ 
given
$\bbbar^{Z^D}(t-1)$ and 
$\bbbar^{Y^D}(t-1)$ satisfying 
Inv-$b_k(t)$ for all $k\geq 1$, and given 
$\n^{Y^D}(t)$,
$\overline{\B}^{Y^D}(t)$, and
$\s^{Y^D}(t)$,
step~$t$ is as follows. If $\n^{Y^D}(t)$, $\overline{B}^{Y^D}(t)$ and $\s^{Y^D}(t)$ contradict the occurrence of $\calE_{\le t}$, then
step~$t$ of~$Z$ evolves independently of step~$t$ of $Y^A$ and~$Y^B$ and   there is no need to prove the invariants. 

Suppose instead that $\n^{Y^D}(t)$, $\overline{B}^{Y^D}(t)$ and $\s^{Y^D}(t)$ are consistent with the occurrence of $\calE_{\le t}$. For each $D\in \{A,B\}$, $\n^{Z^D}(t) = \n^{Y^D}(t)$.  
This, together with the deterministic definition of $\B_k^{Z^D}(t)$, gives
Inv-$B_k(t)$.
For every   $k\geq 1$, $\S^{Z^D}_k(t)$ is formed by including each ball in 
$\B_k^{Z^D}(t) \setminus \B_k^{Y^D}(t)$ 
independently with probability~$p_k$.
Then $\s^{Z^D}_k(t) = \S^{Z^D}_k(t) \cup \s^{Y^D}_k(t)$, establishing Inv-$s_k(t)$.

Since $\s^{Y^D}(t)$ is consistent with $\calE_{\le t}$ for
each $D\in \{A,B\}$, $\calE_t^{Y_D,\s}$  occurs, so $\s^{Y^D}(t)\ne \emptyset$. Since Inv-$s_k(t)$ holds for each $k\ge 1$, it follows that
$\s^{Z^A}(t)$ and $\s^{Z^B}(t)$ are both non-empty.
By the definition of the two-stream process (Definition~\ref{def:two-stream}),
for all $k$, 
$e_k^{Z^A}(t) = e_k^{Z^B}(t) = 0$. The definition of generalised backoff process (Definition~\ref{def:genbackoff}) then 
gives Inv-$b_k(t+1)$ for all $k$.
\end{proof}

Suppose that $\p$ has the property that, 
for all positive integers~$x$, $W_x \le (6/\lambda)^x$. 
Suppose also that $i$ is a bottleneck bin and that $i$ is sufficiently large.  
Lemma~\ref{lem:FinalLBTwoStream}
combines Lemmas~\ref{lem:slow-fillstreams} and~\ref{lem:bottle-events-TwoStream}
in this case
to give a lower bound of $W_i^{-1/3}$ on the probability that a two-stream process, starting from the empty state, has sends in both streams for the first 
$i \lceil  \exp((\log i)^{2/3} ) \rceil$  
steps.   
We will use the lemma to give the same probability bound for the event that there are no escapes during this period in a backoff process (though to do this we will first, in Section~\ref{sec:under-backoff}, consider an intermediate process called an ``under-backoff process'').

\begin{lemma}

\label{lem:FinalLBTwoStream}
Fix $\lambda \in (0,1/2)$ and a send sequence $\p$ with $p_0=1$
such that, for all positive integers~$x$, $W_x \le (6/\lambda)^x$. 
Let $Z=(Z^A,Z^B)$ be a two-stream process with birth rate~$2\lambda$, send sequence~$\p$, and start time~$0$
where for all $D\in \{A,B\}$ and all
positive integers~$k$,
$\arr_k^{Z^D}(0) = \emptyset$. 
Let $i$ be a sufficiently large integer and suppose that bin~$i$ is a bottleneck. 
Let $\TEnd = i \lceil  \exp((\log i)^{2/3} ) \rceil$.
Then 
with probability at least 
$W_i^{-1/3}  $,
for all $t\in [\TEnd]$, 
$\s^{Z^A}(t)$ and $\s^{Z^B}(t)$ are both non-empty.
\end{lemma}

\begin{proof}
Let $j= \lfloor \cbottle \log i \rfloor$.  
Let  $\zvect$ be the $j$-tuple 
with  $z_x = 9\lambda  W_x/20$ for all $x\in [j]$.  
By Lemma~\ref{lem:slow-fillstreams} (with $\tau=0$),
there is an event $\Edom$ which depends only on the evolution of $Z$ during steps in $ [T^{\p,j}]$ such that the following hold.

\begin{itemize}
\item $\Pr(\Edom) \geq (1-\exp({-\lambda/2}))^{2 T^{\p,j}}$,
\item If $\Edom$ occurs then
for each $t\in [T^{\p,j}]$, $\s^{Z^A}(t)$ and $\s^{Z^B}(t)$ are both non-empty,   and

\item 
$\dist( (\bbar^{Z^{A}}_{[j]}( T^{\p,j}),\bbar^{Z^{B}}_{[j]}( T^{\p,j})) \mid \Edom)
\gtrsim   \Po(\zvect,\zvect)$.  
\end{itemize}

Now condition on~$\Edom$ and consider the two-stream process~$Z$
starting from time $T^{\p,j}$.
Apply Lemma~\ref{lem:bottle-events-TwoStream} with $\tau = T^{\p,j}$ and
$\TEnd' = \tau +  \TEnd$.
Given the conditioning on~$\Edom$, 
this lemma shows that 
\[\Pr(
\calE^{Z^A,\calN}_{\leq \TEnd'} \cap \calE^{Z^B,\calN}_{\leq \TEnd'} \cap 
\calE^{Z^A,\s}_{\leq \TEnd'} \cap \calE^{Z^B,\s}_{\leq \TEnd'}
) \geq W_i^{-1/6} .\]
The event $\calE^{Z^D,\s}_{\leq \TEnd'} $
(for any $D\in \{A,B\}$) implies that $\s^{Z^D}(t)$ is non-empty 
for every $t\in \{\tau+1,\ldots, \TEnd'\}$,
so certainly for every 
$t\in \{\tau+1,\ldots,\TEnd\}$.

Thus to finish we need only show 
$(1-\exp^{-\lambda/2})^{2 T^{\p,j}} \geq W_i^{-1/6}$.

Recall that $j= \lfloor \cbottle \log i \rfloor$ where, from Definition~\ref{def:bottleconstants},
$\cbottle = 1/(8\log(7/\lambda))$. We will now use the fact that, for all positive integers~$x$, $W_x \leq (6/\lambda)^x$. Specifically,
from Definition~\ref{def:Tprot}, 
\[2T^{\p,j } = 160 j^2 \lfloor \sum_{x\in [j]} W_x \rfloor 
\le 
160 j^3 (6/\lambda)^j \leq (7/\lambda)^j  \leq i^{1/8}.\]
Thus, 
$(1-\exp({-\lambda/2}))^{2 T^{\p,j}} \geq 
(1-\exp({-\lambda/2}))^{i^{1/8}}$,
and we wish to show that this is at least $W_i^{-1/6}$.
Equivalently, we wish to show that
$W_i \geq 1 /  (1-\exp(-\lambda/2))^{6 i^{1/8}}$.
By the Definition of bottleneck 
(Definition~\ref{def:bottleneck}(i)),
$W_i \ge \exp(i/\exp((\log\log i)^2))$, so this follows since~$i$ is sufficiently large.
\end{proof}

\subsection{Under-backoff processes}
\label{sec:under-backoff}

Now we have the result that we need for two-stream processes, specifically Lemma~\ref{lem:FinalLBTwoStream}
which says that, with probability at least $W_i^{-1/3}$, both streams are blocked
for $i \lceil  \exp((\log i)^{2/3} ) \rceil$ steps by senders (assuming that bin~$i$ is a bottleneck and $i$ is sufficiently large).
Next, we define a final type of auxiliary process, a so-called \emph{under-backoff process}.
In Lemma~\ref{lem:FinalLB-underbackoff}
we apply a coupling to lift the result of Lemma~\ref{lem:FinalLBTwoStream}, so that it applies to an under-backoff process. 
Roughly, the under-backoff has three streams and 
Lemma~\ref{lem:FinalLB-underbackoff} 
ensures that the first two of these stay noisy for $i \lceil  \exp((\log i)^{2/3} ) \rceil$ steps so that the third stream (which is defined slightly differently) acts as a jammed channel.
This enables us, in 
Lemma~\ref{lem:GetIntoi}, 
to show that, given the bottleneck~$i$, the under-backoff process  
has probability at least $W_i^{-6/7}$ of avoiding the empty state until a ball is in bin~$i$ of its third stream.
This makes the expected sojourn time at least $e^{i^{1/2}}$ 
(Lemma~\ref{lem:se-sojourn}).  Since 
infinitely many strongly exposed bins lead to infinitely many bottlenecks (Corollary~\ref{cor:infinite-KMP-bins-2}),  we get  the final result for under-backoff processes, Corollary~\ref{cor:finite-SE-bins}, which states that under-backoff processes   are not positive recurrent when there are infinitely many strongly exposed bins. From there it will be a simple matter, in Section~\ref{sec:final-infinite}, to lift this result (using one final coupling) so that it applies to backoff processes.

\begin{definition}\label{def:under-backoff}
Fix $\lambda \in (0,1/3)$ and
a \emph{send sequence}  
$\p$ with $p_0=1$.
The \emph{under-backoff process} $U$
with birth rate~$3\lambda$ and send sequence~$\p$  
is a 3-tuple $U=(U^A,U^B,U^C)$ of generalised backoff processes with 
send sequence~$\p$, start time~$0$, 
birth distribution
 $\calD =  \Po(\lambda)$
where for each stream $D\in \{A,B,C\}$ 
each bin $k\geq 1$,
and each time $t\geq 0$,
$\arrvect_k^{U^D}(t) = \emptyset$. 
The streams
evolve independently  
except 
that $\e_j^{U^A}(t)$, $\e_j^{U^B}(t)$, and
$\e_j^{U^C}(t)$ are defined as follows for $t>0$ and $j\geq 1$. 
 
\begin{itemize}

\item For $D\in \{A,B\}$, if $s_j^{U^D}(t) = s^{U^A}(t)+ s^{U^B}(t)=1$ then $\e_j^{U^D}(t) = \s_j^{U^D}(t)$. 
Otherwise $\e_j^{U^D}(t) = \emptyset$.

\item If $s_j^{U^C}(t) = s^{U^A}(t)+ s^{U^B}(t) + s^{U^C}(t)=1$ then $\e_j^{U^C}(t) = \s_j^{U^C}(t)$. Otherwise $\e_j^{U^C}(t) = \emptyset$.

\end{itemize}
\end{definition}

It is important to note that, in the under-backoff process $U= (U^A,U^B,U^C)$ with birth rate~$3\lambda$ and sequence~$\p$, 
the balls in~$U^A$ and~$U^B$, taken together, are a backoff process 
(Definition~\ref{def:backoffasgen}, Definition~\ref{def:backoff-cohorts})
with birth distribution~$\Po(2\lambda)$, send sequence~$\p$, and 
cohort set~$\calC= \{A,B\}$.   
On the other hand, the balls in~$U^C$ can have their escapes blocked by balls in~$U^A$ and~$U^B$.

Next, in  Lemma~\ref{lem:FinalLB-underbackoff}
we   lift the result of Lemma~\ref{lem:FinalLBTwoStream}, so that it applies to the under-backoff process.

\begin{lemma}
\label{lem:FinalLB-underbackoff}
Fix $\lambda \in (0,1/3)$. Let $\p$ be a send sequence with $p_0=1$ such that, for all positive integers~$x$, $W_x \le (6/\lambda)^x$. 
Let $i$ be a sufficiently large integer and suppose that bin~$i$ is a bottleneck. 
Let $U = (U^A,U^B,U^C)$ be the under-backoff process with birth rate~$3\lambda$ and send sequence~$\p$.
Let $\TEnd = i \lceil  \exp((\log i)^{2/3} ) \rceil$.
Then 
with probability at least 
$W_i^{-1/3}  $,
for all $t\in [\TEnd]$, 
$\s^{U^A}(t)$ and $\s^{U^B}(t)$ are both non-empty.
\end{lemma}
\begin{proof}
We have observed that the balls in $U^A$ and $U^B$, taken together, are a backoff process~$X$ with birth distribution $\Po(2\lambda)$, send sequence~$\p$, and cohort set~$\calC=\{A,B\}$.    Let $Z=(Z^A,Z^B)$ 
be a two-stream process with birth rate~$2\lambda$, send sequence~$\p$ 
and start time~$0$,  
where for all $D\in \{A,B\}$ and all
positive integers~$k$,
$\arr_k^{Z^D}(\tau) = \emptyset$. The standard coupling of $(U^A,U^B)$ and
$(Z^A,Z^B)$ (Definition~\ref{def:standard-twostream}) has the property that, for every $i\geq 1$, $t\geq 1$, and stream $D\in \{A,B\}$, $\s_i^{Z^D}(t) \subseteq \s_i^{U^D}(t)$.
By Lemma~\ref{lem:FinalLBTwoStream}, with probability at least $W_i^{-1/3}$, for all $t\in [\TEnd]$, $\s^{Z^A}(t)$ and $\s^{Z^B}(t)$ are both non-empty, which implies that 
(under the standard coupling) $\s^{U^A}(t)$ and $\s^{U^B}(t)$ are both non-empty.
\end{proof}

\begin{lemma}
\label{lem:GetIntoi}
Fix $\lambda \in (0,1/3)$. Let $\p$ be a send sequence with $p_0=1$ such that, for all positive integers~$x$, $W_x \le (6/\lambda)^x$. 
Let $i$ be a sufficiently large integer and suppose that bin~$i$ is a bottleneck. 
Let $U = (U^A,U^B,U^C)$ be the under-backoff process with birth rate~$3\lambda$ and send sequence~$\p$.
Let $\TEnd = i \lceil  \exp((\log i)^{2/3} ) \rceil$.
Then 
with probability at least 
$  W_i^{-6/7}  $,
for all $t\in [\TEnd]$, 
$\s^{U^A}(t)$ and $\s^{U^B}(t)$ are both non-empty
and there is a $T^U\in [\TEnd]$ such that $\b_i^{U^C}(T^U)$ is non-empty
and for every $t\leq T^U$, $\bigcup_{k\in [i]} 
\b_k^{U^C}(t)  $ is non-empty.
\end{lemma}

\begin{proof}
By   Lemma~\ref{lem:FinalLB-underbackoff}, with probability at least $W_i^{-1/3}$, for all $t\in [\TEnd]$, $\s^{U^A}(t)$ and $\s^{U^B}(t)$ are both non-empty. 
We will refer to this as the event that $U^A$ and~$U^B$ are noisy.
From Definition~\ref{def:under-backoff}, 
if $U^A$ and~$U^B$ are noisy then
$U^C$ operates independently of $U^A$ and~$U^B$ except that nothing ever escapes.
Thus, conditioned on having~$U^A$ and~$U^B$ be noisy, the event in the lemma is satisfied if $\n^{U^C}(1)$ is non-empty and 
$G_0 + \cdots + G_{i-1} \leq \TEnd$, where $G_x$ is an independent geometric random variable with $\E[G_x] = W_x$. 
By Lemma~\ref{lem:bottleneck-geo-bound}, the probability of this is at least
$(1-e^{-\lambda}) W_i^{-1/2}$.
So the overall  probability is
at least $(1-e^{-\lambda}) W_i^{-5/6} \geq W_i^{-6/7}$.
\end{proof}

\begin{lemma}\label{lem:se-sojourn}
Fix $\lambda \in (0,1/3)$. Let $\p$ be a send sequence with $p_0=1$ such that, for all positive integers~$x$, $W_x \le (6/\lambda)^x$. 
Let $i$ be a sufficiently large integer and suppose that bin~$i$ is a bottleneck. 
Let $U = (U^A,U^B,U^C)$ be the under-backoff process with birth rate~$3\lambda$ and send sequence~$\p$.
Then the expected sojourn time of $U$ from the empty state back to the empty state is at least $e^{{i}^{1/2}}$.
\end{lemma}

\begin{proof}
By Lemma~\ref{lem:GetIntoi}, with probability at least $1/W_i^{6/7}$, 
$U$ reaches a state with $b_i^{U^C}(t) > 0$ before returning to the empty state. Suppose this first occurs at time $T^U$. 
Conditioned on the state of~$U$ at time~$T^U$,
let $\beta$ be an arbitrary ball in $\b_i^{U^C}(T^U)$.
Then $U$  cannot return to the empty state after $T^U$ until $\beta$ sends from bin $i$, which takes expected time $W_i$. Thus the expected sojourn time of $U$ from the empty state is at least $W_i/W_i^{6/7} = W_i^{1/7}$. 
The result follows from $W_i \geq e^{7 i^{1/2}}$, which follows from the Definition of bottleneck
 (Definition~\ref{def:bottleneck}(i)), since $i$ is sufficiently large. 
\end{proof}

\begin{corollary}\label{cor:finite-SE-bins}
Fix $\lambda \in (0,1/3)$. Let $\p$ be a send sequence with $p_0=1$ such that, for all positive integers~$x$, $W_x \le (6/\lambda)^x$. 
 If $\p$ has
 infinitely many   strongly exposed bins
 then the under-backoff process 
 with birth rate~$3\lambda$ and send sequence~$\p$
is not positive recurrent (specifically, the empty state is not positive recurrent).
\end{corollary} 
\begin{proof}
Immediate from Corollary~\ref{cor:infinite-KMP-bins-2} and Lemma~\ref{lem:se-sojourn}.
\end{proof}

\subsection{The result for backoff processes}
\label{sec:final-infinite}

Corollary~\ref{cor:finite-SE-bins} shows that the under-backoff process is not positive recurrent if $\p$ has infinitely many strongly exposed bins. Specifically, the empty state is not positive recurrent, so starting there, the expected time until the process is again empty is infinite.
In order to get the same result for backoff (Corollary~\ref{cor:backoff-finite-SE-bins}) we couple the under-backoff process with the corresponding backoff process so that the population of the under-backoff process stays under the population of the backoff process. The standard coupling, maintaining this invariant, is given in Definition~\ref{def:standard_under}.    This coupling follows the pattern of other standard couplings (for example, Definition~\ref{def:standard-twostream}) though of course the details are tuned to the specific definition of the under-backoff process.

 \begin{definition}\label{def:standard_under}
Fix $\lambda \in (0,1/3)$. Let $\p$ be a send sequence with 
$p_0=1$. 
Let $U=(U^A,U^B,U^C)$ be the under-backoff process with 
birth rate $3\lambda$ and
send sequence $\p$.
Let $X$ be the backoff process with 
birth rate~$3\lambda$, send sequence~$\p$, and cohort set~$\calC=\{A,B,C\}$ (Definition~\ref{def:backoff-cohorts}). 
The standard coupling of $U$ and $X$ has the following invariants for every $i\geq 1$ and $t\geq 1$.
\begin{description}  
\item  [Inv-$b_i(t)$:] For each stream $D\in \{A,B,C\}$, 
$\b_i^{U^D}(t-1) \subseteq \b_i^{X^D}(t-1)$.
\item  [Inv-$B_i(t)$:] For each stream $D\in \{A,B,C\}$, $\B_i^{U^D}(t) \subseteq \B_i^{X^D}(t)$.
\item [Inv-$s_i(t)$:] 
For each stream $D\in \{A,B,C\}$, $\s_i^{U^D}(t) \subseteq \s_i^{X^D}(t)$.
\end{description}
Initially, Inv-$b_i(1)$ holds for every $i\geq 1$ since all of these sets are empty.

For any $t$ assuming Inv-$b_i(t)$, the coupling at step~$t$ proceeds as follows.
For each $D\in \{A,B,C\}$, $\n^{U^D}(t) = \n^{X^D}(t)$, so for all $i\geq 1$, 
Inv-$B_i(t)$ holds.
For every   
$D\in \{A,B,C\}$ and
$i\geq 1$, $\S^{X^D}_i(t)$ is formed by including each ball in 
$\B_i^{X^D}(t) \setminus \B_i^{Z^D}(t)$ 
independently with probability~$p_i$.
Then $\s^{X^D}_i(t) = \S^{X^D}_i(t) \cup \s^{Z^D}_i(t)$, establishing Inv-$s_i(t)$.

The assignments to $\e_i^{X^D}(t)$ 
and $\e_i^{U^D}(t)$ follow deterministically, using  Definitions~\ref{def:backoffasgen} and~\ref{def:under-backoff}.

Since for all $i$ $\s_i^{X^D}(t) \cap \B_i^{U^D}(t) = \s_i^{U^D}(t)$,
the invariant Inv-$b_i(t+1)$ for all $i$ follows from  
$ \s_i^{U^D}(t) \cap \e_i^{X^D}(t) \subseteq \e_i^{U^D}(t)$, which holds: 
If a ball~$\beta$ is in $ \s_i^{U^D}(t) \cap \e_i^{X^D}(t) $ then
$\s^{X^A}(t) \cup \s^{X^B}(t) \cup s^{X^C}(t) = \{\beta\}$ so
$\s^{U^A}(t) \cup \s^{U^B}(t) \cup s^{U^C}(t) = \{\beta\}$ so 
$\e_i^{U^D}(t) = \{\beta\}$.
\end{definition}

 We can now prove the main result of this section, 
 Corollary~\ref{cor:backoff-finite-SE-bins}. Note that we don't need the flexibility of~$\lambda_0$ and we could weaken the hypothesis  $W_i\leq (2/\lambda)^i$, but we write the corollary this way for consistency where it is used.

\begin{corollary}\label{cor:backoff-finite-SE-bins}
There exists a real number $\lambda_0\in (0,1)$ such that the following holds. 
Consider any $\lambda \in (0,\lambda_0)$. Consider any send sequence~$\p$ such that $p_0=1$ and, for all $i\geq 1$, $W_i \leq (2/\lambda)^i$. Suppose that $\p$ has infinitely many bins which are strongly exposed. Then
the backoff process with birth rate~$\lambda$ and send sequence~$\p$ is not positive recurrent. 
\end{corollary}

 \begin{proof}
 Fix any $\lambda_0 \in (0,1)$. Consider any $\lambda \in (0,\lambda_0)$ 
 and let $\lambda' = \lambda/3$ so that $\lambda' \in (0,1/3)$. 
 Consider any send sequence~$\p$ such that $p_0=1$ and, for all $i\geq 1$, 
 $W_i \leq (2/\lambda)^i \leq (2/\lambda')^i$. Suppose that $\p$ has infinitely many bins which are strongly exposed. Let $X$ be the backoff process with birth rate $\lambda = 3\lambda'$ and send sequence~$\p$ and let $U$ be the under-backoff process with birth rate~$\lambda = 3\lambda'$ and send sequence~$\p$.
 The standard coupling of~$U$ and~$X$ (Definition~\ref{def:standard_under})
 ensures that for each 
 $i\geq 1$, each $t\geq 1$ and each
 stream $D\in \{A,B,C\}$,
 $\b_i^{U^D}(t-1) \subseteq \b_i^{X^D}(t-1)$.
Corollary~\ref{cor:finite-SE-bins} shows that~$U$ is not positive recurrent, so the expected return time to the empty state, after having started there, is infinite.
Since $\b_i^{U^D}(t-1) \subseteq \b_i^{X^D}(t-1)$, the same is true of~$X$.
\end{proof}

\section{High-level states for the case with finitely many strongly exposed bins}
\label{sec:hls}\label{sec:finite}

As we explained in the introduction, we analyse the backoff process by coupling it with a volume process~$V$ that tracks sends and an infinite number of escape processes that track the escapes. The invariant is that every ball present in the volume process~$V$ but not in an escape process is also present in the original backoff process. Thus, to show that the backoff process is not positive recurrent, it suffices to show an infinitely-growing population in the volume process~$V$ that is not in the escape processes.

The volume process~$V$ is designed so that, for all steps $t$, there is a $j$ (depending on $t$) so that the populations of the first $j$ bins are independent Poisson variables (the volume process functions as a $j$-jammed process).  The dependence of the escape process on the volume process is complex, so in order to keep things manageable, the relevant information about the state of the volume process is 
summarised by the so-called ``high-level state''. The high-level state is a function of the volume process, and the escape process depends on it (but does not otherwise depend on the volume process).

A high-level state 
(Definition~\ref{def:high-level-state-newer}) is a tuple
$\Psi = (g,\tau, j,\zvect,\calS,\type)$. 
At time~$t$, the volume process is in high-level state~$\Psi(t)$.
The meaning of the elements of~$\Psi$ are as follows. $j$ is a positive integer and $\zvect=(z_1,\ldots,z_j)$ is a tuple of non-negative reals. The idea is that when $\Psi(t)=\Psi$, 
the volume process~$V$ acts as a $j$-jammed process.  
The non-negative integers~$g$ and~$\tau$ are for book-keeping: $\tau$ is the first time such that $\Psi(\tau) = \Psi$, and $\Psi$ is the $g$'th unique high-level state in the sequence $\Psi(\tau_0),\Psi(\tau_0+1),\ldots$ where $\tau_0$ is the time of the start of the volume process~$V$. The $j$-tuple $\zvect$ 
has the property that, conditioned on the transition into high-level state~$\Psi$ at time~$\tau$, 
$\dist(\bbar^{V}(\tau)) =\Po(\zvect)$.
 $\calS$ is a set of disjoint sets of bins in $[j]$. The idea is that when $\Psi(t) = \Psi$, for all $S \in \calS$, most bins in $S$ provide noise in the volume process; thus having multiple sets in $\calS$ provides something stronger than a uniform bound on noise in $\cup_{S \in \calS}S$.

\begin{definition}\label{def:high-level-state-newer}
A \emph{high-level state} is a tuple $\newsym{Psi}{High-level state}{\Psi} = (g, \tau, j, \zvect, \calS,\type)$ where 
$g,\newsym{tau}{Start time of some process}{\tau}\ge 0$ and $j \ge 1$ are integers, $\newsym{z-bf}{Tuple of non-negative reals, typically expected initial ball counts}{\zvect}=(z_1,\ldots,z_j)$ is a tuple of non-negative reals, $\calS$ is a set of disjoint 
subsets of $[j]$,
and $\type$ is a type in 
$\{ \Failure, \initialising, \advancing, \filling, \refilling, \stabilising\}$. 
Given a high-level state $\Psi = (g,\tau,j,\zvect,\calS,\type)$, we use $g^\Psi$ to denote $g$, $\tau^\Psi$ to denote $\tau$, $j^\Psi$ to denote~$j$, ${\zvect}^\Psi$ to denote~$\zvect$, $\calS^\Psi$ to denote $\calS$, and $\type^{\Psi}$ to denote $\type$. 
\end{definition}

We will say more about most of
the types in Definition~\ref{def:high-level-state-newer} in Section~\ref{sec:sets}, but
we first define the simplest types of high-level state -- failure and initialising.

\begin{definition}\label{def:failure-state}
A high-level state 
$\Psi=(g,\tau,j,\zvect,\calS,\type)$ 
is a \emph{failure state} if 
$\type = \newsym{failure}{failure high-level state, Definition~\ref{def:failure-state}}{\Failure}$, $\zvect= \zerovect$  and $\calS=\emptyset$.  
\end{definition}

To emphasise~$j$, we refer to
a failure state $\Psi=(g,\tau,j,\zvect,\calS,\Failure)$ as a $j$-failure state.
We use similar terminology for the other types of high-level states.

\begin{definition}\label{def:init-state}    
A high-level state $\Psi=(g,\tau,j,\zvect,\calS,\type)$ 
is an \emph{initialising state} if 
$\type = \newsym{initialising}{initialising high-level state, Definition~\ref{def:init-state}}{\initialising}$, 
for each $k\in [j]$,
$z_k = 3 \lambda W_k/4$,  
$\calS=\emptyset$, and $g=0$.  
\end{definition}

It is important that $j$-failure states  have $\calS = \emptyset$.  Once the volume process enters a high-level state with type $\Failure$ we expect nothing from it (we have given up) -- we will show that this happens with summable failure probability. Similarly,
we don't require any noise from a $j$-initialising state so $\calS = \emptyset$. The reason  that $\calS=\emptyset$ is OK for
the initialising high-level state is that we will set up the transition rule between high-level states so that the volume process transitions immediately from its initialising state (and never visits it again). The vector $\zvect$ of an initialising state reflects the fact that the volume process starts after the first $j$ bins have (mostly) filled up in expectation (for some $j$).

 At the end of step~$t$ the volume process transitions from high-level state $\Psi(t-1)$ to high-level state $\Psi(t)$ (which may or may not differ from $\Psi(t-1)$). The transition from $\Psi(t-1)$ to $\Psi(t)$ is determined by applying a specific high-level transition rule~``$\hls$''
 that we will define in 
 Section~\ref{sec:backoff-bounding-valid}.
 The new high-level state $\Psi(t)$
 is determined by $\Psi(t) = \hls(\Psi(t-1), 
\bbar,t)$, where $\bbar$ is a $j$-tuple depending
on $\bbar^V_{[j]}(t-1)$ -- see Definition~\ref{def:volume-process} for details.  
 Definition~\ref{def:transition-rule} first gives the generic definition of a high-level transition rule, specifying some important properties which we will need, to make the coupling of the backoff process with the volume process and the escape processes feasible.

\begin{definition}
\label{def:transition-rule}
Fix a real number $\lambda \in (0,1)$ and a send sequence $\p$ with $p_0=1$.
Let $\newsym{Psi-star}{set of high-level states}{\Psi^*}$ be a set of high-level states.
The domain 
$\newsym{R-Psi}{domain of high-level transition rule}{\calR(\Psi^*)}$  is the set of triples
$(\Psi,\bbar,t)$ such that
\begin{itemize}
\item  $\Psi \in \Psi^*$, ${\bbar}\in 
(\integers_{\geq 0})^{j^\Psi}
 $,   and $t> \tau^{\Psi}$ is an integer.  
\item If $\type^\Psi=\initialising$  then $t=\tau^\Psi+1$.   
\end{itemize}

A \emph{transition rule~$\hls$ on $\Psi^*$} is a function~$\newsym{hls}{high-level state transition rule (generic) - Definition~\ref{def:transition-rule}}{\hls}\colon \calR(\Psi^*) \to \Psi^*$ with the following properties
for all tuples
$(\Psi, \bbar, t) \in \calR(\Psi^*)$, where
$\Psi= (g,\tau,j,\zvect,\calS,\type)$
and $\Psi' = (g',\tau',j',\zvect',\calS',\type') = \hls(\Psi,\bbar,t)$. 

\begin{enumerate}[(T1)]
\item If $\Psi' \ne \Psi$, then $g'=g+1$ and $\tau' = t$.\label{item:trans1}
\item If $\type=\Failure$, then  
$\Psi' = (g+1,t,j,\zerovect,\emptyset, \Failure)$.
\label{item:trans2}
\item For all $S \in \calS'$, $\sum_{x \in S} p_x b_x \ge \lambda |S|/40$.\label{item:trans3}
\item   $j' \geq j$. For all $\ell \in [j'] \setminus [j]$, $z'_\ell = 0$.  
If $\type=\initialising$ then $j'=j$.
\label{item:trans4}
\item $\type'\neq \initialising$. \label{item:trans5}
\end{enumerate}  
\end{definition}

The definition of the domain $\calR(\Psi,\bvect, t)$ captures the fact that  
the
transition from a high-level state~$\Psi$ 
depends on~$\Psi$, and on a $j$-tuple $\bvect$ (which for us will be the bin populations of the first $j$ bins in the volume process), and on the time~$t$ of the new high-level state. The second item ensures that if the initialising state starts at time~$\tau$, then the transition out of the initialising state is at the next step, step $\tau+1$
(and a new initialising state is not re-entered, by (T\ref{item:trans5})).
Item (T\ref{item:trans1}) is for book-keeping -- it ensures that if the high-level state changes then the time of this change is properly recorded in the new high-level state, and the ``sequence number'' $g$ is incremented. 

While the volume process is
in a high-level state $\Psi$,
it evolves like a $j$-jammed process with $j= j^{\Psi}$.
The evolution of the volume process only differs from the $j$-jammed process in two ways. First, if the high-level state changes to a $\Psi(t)$ that differs from $\Psi(t-1)$, then $j^{\Psi(t)}$ may be one larger than $j^{\Psi(t-1)}$ (if we are moving attention from bin $j^{\Psi(t-1)}$ to the next bin along). Second, some balls are "dropped" at each time step, whether or not the high-level state changes. This is in order to maintain the volume invariant discussed in the introduction, namely that for a volume process V with start time $\tau_0$, conditioned on $\Psi(\tau_0),\dots,\Psi(t)$ and $\Psi(t)=\Psi$, $\dist(\bvect^V_{[j^\Psi]}(t)) = \Po(F^\Psi_{[j^\Psi]}(t))$, where $F^\Psi_{[j^\Psi]}(t)$ is defined below in Definition~\ref{def:F} and tracks the expectations for a $j$-jammed process with start time $\tau^\Psi$ and initial state $\zvect^\Psi$. Note that the volume invariant implies that, given the conditioning on $\Psi(\tau_0),\ldots,\Psi(t)$ and $\Psi(t)=\Psi$, the distribution of~$\bvect_{[j^\Psi]}^V(t)$ just depends on~$\Psi$ and on~$t$. Dropping balls in the transition to~$\bvect_{[j^\Psi]}^V(t)$ has to be done carefully in a way that keeps the coupling valid. Definition~\ref{def:transition-rule}
puts some restrictions on how this can be done, as follows.

  (T\ref{item:trans2}) guarantees that if $\Psi$ is a failure state then so is $\Psi'$ (so once the volume process enters a failure state it stays in failure states). 
 
   (T\ref{item:trans3}) guarantees that  
for all $S \in \calS^{\Psi'}$, $\sum_{x\in S} p_x \bbar_x \geq \lambda |S|/40$. Thus, every $S \in \calS^{\Psi'}$ acts as a guaranteed source of expected noise for the escape process while the volume process is in state $\Psi'$.

(T\ref{item:trans4}) guarantees 
that $j$ does not decrease as high-level state changes. Also, if $j$ increases then
we do not assume that the volume process has already acquired a population in the new bin.
Finally, $j$ does not increase as the process leaves the initialising state.

Since the value of $j$ can change when the high-level state changes, it will be useful to re-index the function $f$ (Definition~\ref{def:f}) that we used to track the means of the bins of the $j$-jammed process. We do this in Definition~\ref{def:F}.

\begin{definition}\label{def:F}
Let $\Psi = (g,\tau,j,\zvect, \calS, \type)$ be a high-level state. For $k\in [j]$ 
and $t\geq \tau$,
we use $\newsym{F}{Function to track means of population of bin~$k$ - Definition~\ref{def:F}} {F_k^{\Psi}(t)}$ to denote 
$f_k^{\zerovect,\zvect}(t-\tau)$, where the $0$-vector has length~$j$.
We also use 
$\newsym{F-2}{projection of $F$ onto $[k]$ - Definition~\ref{def:F}} {F_{[k]}^{\Psi}(t)}$ 
for the tuple $(F_1^{\Psi}(t),\ldots,F_k^{\Psi}(t))$.
\end{definition}

\begin{observation}\label{obs:F}
Let $\Psi = (g,\tau,j,\zvect, \calS, \type)$ be a high-level state. For $k\in [j]$,
$F_k^{\Psi}(\tau) = z_k$.
\end{observation}
\begin{proof}
By Definition~\ref{def:F}, 
$F_k^{\Psi}(\tau) = f_k^{\zerovect,\zvect}(0)$. By Definition~\ref{def:f}, this is~$z_k$. 
\end{proof}

Lemma~\ref{lem:slowfill-j-jammed} shows that for any sufficiently large~$t$ 
and any process~$X$ that dominates a $j$-jammed process from above, 
$\dist(\bvect^X_{[j]}(t)) \gtrsim \Po(\lambda 9 W_1/10,\ldots, \lambda 9W_j/10)$.
Looking ahead, this will be used in the coupling of volume, escape, and backoff (Definition~\ref{def:VEB-couple}) to show that, from some initialising state~$\Psi_0 = (0,\tau_0,j_0,\zvect,\emptyset,\initialising)$, 
$\dist(\bvect^{V}_{[j_0]}(\tau_0)) \gtrsim \Po(\zvect^{\Psi_0})$.
Lemma~\ref{lem:dist-volume} will establish a similar, but more general, property for the volume process, 
in terms of the changing high-level state,
for all times $t\geq \tau_0$.
As the high-level state of the volume process evolves, it behaves as a $j$-jammed process, but $j$ increases over time. We will need to
ensure that this kind of  Poisson domination continues, even as the high-level state changes.
Thus, we will need the Poisson domination for the volume process, even conditioned on the current high-level state.
To enable this, we need to further restrict the high-level state transition function. We will do this in Definition~\ref{def:valid}, following Definition~\ref{def:hls-closed},
which identifies a relevant set of high-level states.

\begin{definition}\label{def:hls-closed}
Fix $\lambda \in (0,1)$
and let  $\Psi_0$ be an initialising high-level state.
Let $\Psi^*$ be a set of high-level states which includes $\Psi_0$. 
Let $\hls$ be a transition rule on $\Psi^*$.
Define the relation $R$ on $\Psi^*$ 
by $R(\Psi,\Psi')$ iff $\exists \bbar \in 
(\integers_{\geq 0})^{j^\Psi} $, $t> \tau^\Psi$ such that
$\Psi' = \hls(\Psi,\bbar,t)$.
We define the 
\emph{\newsym{Psi-closure}{closure of high-level state transition function}
{\Psi_0\mbox{-closure of $\hls$}}}  
to be the minimal subset of~$\Psi^*$
which contains~$\Psi_0$ and is closed under
the relation~$R$.
\end{definition}

\begin{definition}\label{def:valid}
Fix $\lambda \in (0,1)$ and a send sequence $\p$ with $p_0=1$. Fix an initialising state~$\Psi_0$. Let $\hls$ be a transition rule on
a set of high-level states which includes~$\Psi_0$. Let $\Psi^{*}$ be the $\Psi_0$-closure of $\hls$. The transition rule~$\hls$ is 
\emph{\newsym{Psi-valid}{validity condition for hls - Definition~\ref{def:valid}} {\Psi^*\mbox{-valid}}}  if 
for every $\psi \in \Psi^*$ and every integer $t > \tau^\psi$, the following holds.
\begin{itemize}
\item  { Let $\Xvect(\psi,t) \sim \Po(F_{[j^\psi]}^\psi(t))$.}
Let $\Psi'$ be the random high-level state 
$\Psi' = \hls(\psi, { \Xvect(\psi,t)}, t)$.
For every high-level state $\psi'$ 
such that $\Pr(\Psi' = \psi')>0$, 
let $A(\psi,\psi',t) = \dist({ \Xvect(\psi,t)} \mid \Psi'=\psi')$.  
Then for all $\psi'$, $A(\psi,\psi',t) \gtrsim \Po(F_{ [j^\psi]}^{\psi'}(t))$. 
\end{itemize}

Fix a $\Psi^*$-valid transition rule 
$\hls$. 
Recall the function~$\phi_{\cdot,\cdot}$ from Definition~\ref{def:phi}.
For any high-level states
$\psi,\psi' \in \Psi^*$ 
and $t> \tau^{\psi}$,  we define 
\newsym{phi-2}{stochastic domination 
map determined by hls $\psi$ and $\psi'$ - see Definition~\ref{def:valid}}
{\phi_{\psi,\psi',t}} as follows.
If $A(\psi,\psi',t) \gtrsim 
\Po(F_{ [j^\psi]}^{\psi'}(t))$ then
$\phi_{\psi,\psi',t} = \phi_{A(\psi,\psi',t),
\Po(F_{ [j^\psi]}^{\psi'}(t)) }$.
Otherwise $\phi_{\psi,\psi',t}$ is the identity.
\end{definition}

In practice, it is difficult to apply Definition~\ref{def:valid} so what we do instead is define some properties (V\ref{item:valid-zlower}), (V\ref{item:valid-up}), and (V\ref{item:valid-up-non-zero}) and prove in Lemma~\ref{lem:valid-transition} that any transition rule with these properties is valid  in the sense of Definition~\ref{def:valid}. When we define our specific high-level transition rule in Section~\ref{sec:backoff-bounding-valid}, we will verify the properties 
(V\ref{item:valid-zlower}), (V\ref{item:valid-up}), and (V\ref{item:valid-up-non-zero}).

(V\ref{item:valid-zlower}) guarantees  that the assumed expected starting position~$z'$ of the next high-level state $\psi'$ cannot be above the expected population sizes that would arise
if the volume process just stayed in~$\psi$. (This enables us to keep the coupling valid -- we can drop the assumed population by ignoring some balls -- that only hurts us --  but we would be cheating if we created balls out of thin air.)
The point of properties (V\ref{item:valid-up}) and (V\ref{item:valid-up-non-zero}) is to formalise our strategy for conditioning on sufficient noise (as captured by the high-level state) while maintaining independent Poisson populations -- it sets up conditions which guarantee this in Definition~\ref{def:valid}.

\begin{lemma}\label{lem:valid-transition}
Fix $\lambda \in (0,1)$ and a send sequence $\p$ with $p_0=1$. 
Fix an initialising state~$\Psi_0$. Let $\hls$ be a transition rule on a set of high-level states which includes~$\Psi_0$. Let $\Psi^*$ be the $\Psi_0$-closure of~$\hls$. The
transition rule $\hls$   is 
$\Psi^*$-valid if it has properties (V\ref{item:valid-zlower})--(V\ref{item:valid-up-non-zero}), stated below,
for all 
$(\psi= (g,\tau,j,\zvect,\calS,\type),\bbar,T)\in \calR(\Psi^*)$ where $\psi' = (g',\tau',j',z',\calS',\type') = \hls(\psi,\bbar,T)$.

In order to state (V\ref{item:valid-zlower})--(V\ref{item:valid-up-non-zero}),
let $ \newsym{over-1}{set of $j$-tuples above $\bbar$}{
\bupset} = \{\bbarup\in (\integers_{\geq 0})^j  \mid
 \bbarup \geq \bbar\}$. Let 
 \newsym{over-2}{set of $j$-tuples above $\bbar$ with specified components equal}{
\bupsomeset}  be the set of all tuples $\bbarupsome\in \bupset $ such that
for all $k\in  [j]$ with $z'_k=0$, 
$\bupsome_k = b_k$. The properties are as follows.

\begin{enumerate}[(V1)]
\item   \label{item:valid-zlower} If ${\psi'} \ne \psi$, then for all $\ell \in  [j]$, 
$z'_\ell \le F_\ell^{\psi}(T)$  
        
\item  \label{item:valid-up}   If ${\psi'}=\psi$, then for all 
$\bbarup \in \bupset$, 
$\hls(\psi,\bbarup, T)={\psi}$.  
        
\item  \label{item:valid-up-non-zero}  If ${\psi'} \ne \psi$,
then for all $\bbarupsome \in \bupsomeset$, 
 $\hls(\psi,\bbarupsome,T)=\psi'$.
\end{enumerate}
\end{lemma}

\begin{proof}

Fix $\psi= (g,\tau,j,\zvect,\calS,\type)\in \Psi^*$ and $t> \tau^{\psi}$. Let $\Xvect(\psi,t)$ be a random $j$-tuple
of non-negative integers distributed as
$\Xvect(\psi,t) \sim \Po(F_{[j]}^\psi(t))$.

Let $\Psi'$ be the random high-level state  
$\Psi' = \hls(\psi, { \Xvect(\psi,t)}, t)$.
Fix a  high-level state $\psi'=(g',\tau',j',z',\calS',\type')$ 
such that $\Pr(\Psi' = \psi')>0$. 
Let $A(\psi,\psi',t) = \dist({ \Xvect(\psi,t)} \mid \Psi'=\psi')$
We wish to show that $A(\psi,\psi',t) \gtrsim 
\Po(F_{ [j]}^{\psi'}(t))$.

{\bf Case 1: $\mathbf \psi'= \psi$.\quad}
Since $\dist(\Xvect(\psi,t)) = \Po(F_{[j]}^\psi(t))$, it suffices to show
\begin{equation}\label{eq:Ashowone}
A(\psi,\psi',t) \gtrsim \dist(\Xvect(\psi,t)),
\end{equation}
and
\begin{equation}\label{eq:Bshowone}
\Po(F_{[j]}^\psi(t))  \gtrsim   \Po(F_{ [j]}^{\psi'}(t)).
\end{equation}
Since $\psi' = \psi$, \eqref{eq:Bshowone} is immediate.
To finish, we will prove~\eqref{eq:Ashowone}.
Let $\bbar^A \sim A(\psi,\psi',t)$.
Using Observation~\ref{obs:dominate-suff-cond}
what we must prove is that, for every   $\xbar\in (\integers_{\geq 0})^j$, 
$\Pr(\bbar^A \geq \xbar) \geq \Pr( \Xvect(\psi,t) \geq \xbar)$.
Fix~$\xbar\in (\integers_{\geq 0})^j$
such that $\Pr( \Xvect(\psi,t) \geq \xbar)>0$.
Let $\Xvect$ be the tuple $\Xvect(\psi,t)$.
Note crucially that the elements of this tuple are independent.
Let $\calE_1$ be the event that $\Xvect\geq \xbar$. 
Clearly this is an increasing event in
the $\sigma$-field of~$\Xvect$, $\sigma(\Xvect)$.
Let
$\calE_2$ be the event that $\Psi' = \psi'$. 
Given the fixed~$\psi$ and~$t$, $\calE_2 \in \sigma(\Xvect)$. 
We have chosen~$\psi'$ so that $\Pr(\calE_2)>0$.
(V\ref{item:valid-up}) guarantees that $\calE_2$ is increasing in~$\Xvect$. 
The result then follows from
Corollary~\ref{cor:Harris}.

{\bf Case 2: $\mathbf \psi' \neq \psi$.\quad}
Let $W = \{ i\in [j] \mid z'_i > 0\}$
and let $ W^0 = [j]\setminus W$.
For convenience, for the rest of this proof, we will represent a $j$-tuple $\xbar$ of non-negative integers as the pair
$(\xbar_{W},\xbar_{W^0})$.  
We wish to show that $A(\psi,\psi',t) \gtrsim 
\Po(\f_{ [j]}^{\psi'}(t))$
but note that any~$\xbar$
that has positive probability as an output of
$\Po(\f_{ [j]}^{\psi'}(t))$
is of the form $(\xbar_W,\zerovect)$.
This is important because when we do the stochastic domination the 
output of 
 $\Po(\f_{ [j]}^{\psi'}(t))$
 will therefore be of the form
 $(\xbar_W,\zerovect)$.
 The output of $A(\psi,\psi',t)$ is not similarly constrained so it will be of the general form 
 $(\xbar_{W^+},\xbar_{W^0})$ -- and we will often denote $\xbar_{W^0}$
 as a vector~$\avect$.

For every tuple 
$\xbarWbar\in (\integers_{\geq 0})^j_{W^0}$ 
let $\calE_{\xbarWbar}$ be 
the event that $\Xvect_{W^0}(\psi,t) = \xbarWbar$.
Let $\Omega^0$ be the set of
tuples
$\xbarWbar\in (\integers_{\geq 0})^j_{W^0}$ 
such that $\Pr(\calE_{\xbarWbar})>0$
and $\Pr(\Psi' = \psi' \wedge \calE_{\xbarWbar})>0$.

Let $ A(\psi,\psi',t, \xbarWbar) = \dist({ \Xvect_W(\psi,t)} \mid \Psi'=\psi',  \calE_{\xbarWbar})$.
Drawing a sample from $A(\psi,\psi',t)$  is the same as sampling $\xbarWbar$ from some appropriate distribution on $\Omega^0$ (let's call that distribution $\calL$) and then sampling $\xbar_W$ from 
$ A(\psi,\psi',t, \xbarWbar) $.
It is important here that $\xbarWbar \in \Omega^0$ -- otherwise there are no outputs $(\cdot,\xbarWbar)$ 
with positive probability in~$A(\psi,\psi',t)$.

We have now reached an important conclusion. To show 
$A(\psi,\psi',t) \gtrsim 
\Po(\f_{ [j]}^{\psi'}(t))$
it suffices to show that, for every $\xbarWbar \in \Omega^0$,
$ A(\psi,\psi',t, \xbarWbar) \gtrsim 
\Po(\f_{ W}^{\psi'}(t))$.
This is because we sample from $A(\psi,\psi',t)$ by first choosing $\xbarWbar$ and then sampling from 
$ A(\psi,\psi',t, \xbarWbar)$
whereas we sample from 
$\Po(\f_{ [j]}^{\psi'}(t))$
by choosing $\zerovect$ on $W^0$
and then sampling from $\Po(\f_{ W}^{\psi'}(t))$.
(The way to think about this 
in the language of Definition~\ref{def:phi} is that, once we have a joint distribution $C_{\xbarWbar}$ for
each~$\xbarWbar$, we can combine them to get~$C$ by just sampling from $C_{\xbarWbar}$ with the probability that $\xbarWbar$ is output in~$\calL$.)

Let $\Xvect(\psi,t,\xbarWbar) \sim \dist({ \Xvect_W(\psi,t)} \mid    \calE_{\xbarWbar}) = \Po(F_W^\psi(t))$.
The final equality in the previous sentence is crucial, and comes from the definition of a Poisson tuple and specifically from the fact that the components of the tuple are independent.
To finish, we need to show, for all $\xbarWbar \in \Omega^0$, that
\begin{equation}\label{eq:Ashowtwo}
A(\psi,\psi',t,\xbarWbar) \gtrsim  \dist(\Xvect(\psi,t,\xbarWbar)),
\end{equation}
and
\begin{equation}\label{eq:Bshowtwo}
 \Po(F_W^\psi(t))  \gtrsim   \Po(\f_{ W}^{\psi'}(t)).
\end{equation}

For~\eqref{eq:Bshowtwo},  it suffices to show that
$\F_W^\psi(t) \geq \f_{ W}^{\psi'}(t)$
then the result follows by Observation~\ref{obs:dom-Poisson-tuples}.
From (T\ref{item:trans1}), $\tau'=t$.
Taking the definition of $\F$ from
Definition~\ref{def:F},
we wish to show for $k\in W$
that 
$F_k^\psi(t)  \geq 
f_k^{\zerovect,\zvect'}(t-\tau')=
f_k^{\zerovect,\zvect'}(0) = z'_k
$ where   $f$ is defined in  Definition~\ref{def:f}. By (V\ref{item:valid-zlower})  
$z'_k \leq F_k^\psi(t)$, as required.

To finish, we will prove~\eqref{eq:Ashowtwo}.
Let $\bbar^A \sim A(\psi,\psi',t,\xbarWbar)$
and let $\Xvect \sim \Xvect(\psi,t, \xbarWbar)$.
Using Observation~\ref{obs:dominate-suff-cond}
what we must prove is that, for every  
$\xbar_W$,
$\Pr(\bbar^A \geq \xbar_W) \geq \Pr(  \Xvect \geq \xbar_W)$.
Fix $\xbar_W$.
Note that both distributions are conditioned on $\calE_{\xbarWbar}$.
Let $\calE_1$ be the event that 
$\Xvect\geq \xbar_W$. 
The event~$\calE_1$ is clearly in $\sigma(\Xvect)$ and is increasing in~$\Xvect$.
Let
$\calE_2$ be the event that $\psi' = \Psi'$. 
Given the fixed~$\psi$ and~$t$ it is clear that 
$\calE_2$ is in $\sigma(\Xvect(\psi,t))$.
Conditioned on $\calE_{\xbarWbar}$, it is in $\sigma(\Xvect)$.
We have chosen~$\psi'$ and $\xbarWbar$ so that $\Pr(\calE_2\mid \calE_{\xbarWbar})>0$.
(V\ref{item:valid-up-non-zero}) guarantees that $\calE_2$ is increasing in~$\Xvect$. The result then follows from
Corollary~\ref{cor:Harris}.
\end{proof}

\section{The volume process}\label{sec:vol}

Definition~\ref{def:volume-process} 
gives the formal definition of the volume process, which we will couple with a backoff process in order to track sends.
At any point in time, the volume process mimics a $j$-jammed process for some $j$.
The value of~$j$ increases as the high-level state changes - which happens
as part of the volume process, according to the high-level state transition rule~$\hls$.

\begin{figure}[ht]   
\centering
 
\begin{tikzpicture}[
  box/.style={draw, rounded corners, align=center, minimum width=28mm, minimum height=10mm},
  >={Latex[length=2.2mm]},
  arrow/.style={->, thick, shorten <=2pt, shorten >=2pt}
]

\node[box] (psi_prev) at (0,  2.2) {$\Psi(t-1)$};
\node[box] (bV_prev)  at (0,  0.0) 
{$\bbbar^{V}(t-1)$};
\node[box] (bY_prev)  at (0, -2.2) 
{$\bbbar^{Y_{t-1}}(t-1)$};

\node[box] (bY)       at (3.8,-2.2) {$\bbbar^{Y_{t-1}}(t)$};

\node[box] (psi)      at (7.6, 0.0) {$\Psi(t)$};

\node[box] (bV)       at (11.4,0.0) {$\bbbar^{V}(t)$};

\draw[arrow] (bV_prev) -- (bY_prev);
 
\draw[arrow] (bY_prev) -- (bY);

\draw[arrow] (psi_prev.east) to[out=0,in=160] (psi.west);

\draw[arrow] (bY.east) to[out=0,in=-90] (psi.south);

\draw[arrow] (psi) -- (bV);

\draw[arrow] (bY.east) to[out=0,in=-120] (bV.south west);

\draw[arrow] (psi_prev.east) to[out=0,in=140] (bV.north west);

\end{tikzpicture}
 
\caption{Illustration of Definitions~\ref{def:volume-process-ext} and \ref{def:volume-process}.}
\label{fig:volume}
\end{figure}
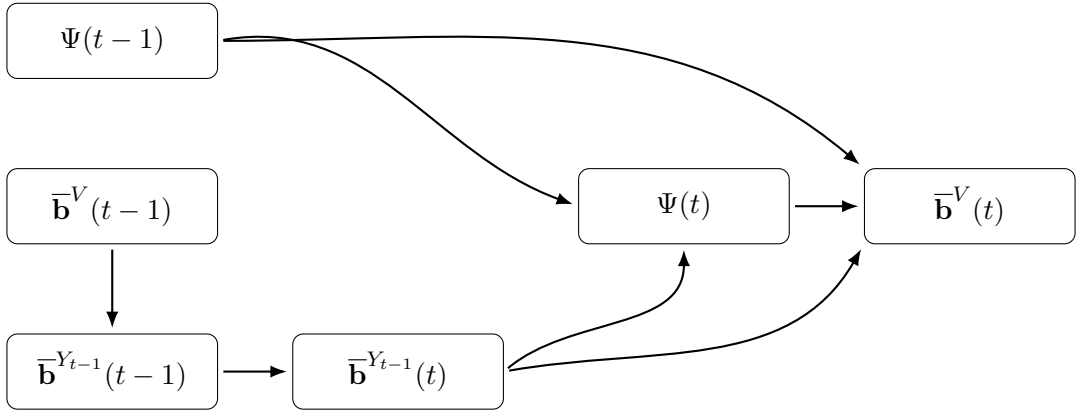

The definition is difficult to understand, so we first give a longer
definition, Definition~\ref{def:volume-process-ext}, which includes extra explanation. 
The explanation will be made rigorous later in Lemma~\ref{lem:dist-volume} --  it is just there to aid understanding.
We then give the  shorter definition, Definition~\ref{def:volume-process}, which contains only the definition (without explaining invariants), and can be referred to later. The reader is encouraged to read these in order -- first Definition~\ref{def:volume-process-ext} with the informal explanation, then
Definition~\ref{def:volume-process} (the definition, for future reference), and finally the lemma, Lemma~\ref{lem:dist-volume}.
Figure~\ref{fig:volume} illustrates the definitions.

\begin{definition}\label{def:volume-process-ext}
Fix a real number $\lambda \in (0,1)$ and a send sequence $\p$ with $p_0=1$. 
Fix an initialising state 
$\Psi_0 = (0, \tau_0,j,\zvect,\emptyset,\initialising)$.
Let $\hls$ be a transition rule on a set of high-level states which includes~$\Psi_0$. Let $\Psi^*$ be the $\Psi_0$-closure of~$\hls$. Suppose that $\hls$ is $\Psi^*$-valid.
A \emph{volume process} $V$ 
from $\Psi_0$ with transition rule~$\hls$, send sequence~$\p$ and birth rate~$\lambda$ is defined as follows. 

Choose the initial population as follows.
Choose 
a random vector 
$\romanarrvect = (\romanarr_1,\ldots,\romanarr_j)
\sim\Po(F_{[j]}^{\Psi_0}(\tau_0))$.
Let $\arrvect^V(\tau_0)$ be a tuple
of $j$~disjoint sets
$\arrvect^V(\tau_0) = (\arr^V_1(\tau_0),\ldots,\arr^V_j(\tau_0))$ 
where each $\arr^V_i(\tau_0) \subseteq \calB(\tau_0)$ and has size 
$\romanarr^V_i(\tau_0)$.
Let $\bbbar^V_{[j]}(\tau_0) = \arrvect^V(\tau_0)$.
For all $i>j$, let   ${\b^V_i}(\tau_0) = \emptyset$. 

Let the initial high-level state $\Psi(\tau_0)$ be~$\Psi_0$. 
Crucially, we will show in Lemma~\ref{lem:dist-volume} that the following invariant is maintained, conditioned on the value of $\Psi(\tau_0),\ldots,\Psi(t)$.
\begin{itemize}
\item Invariant: $
\dist( \bbar^{V}_{[j^{\Psi(t)}]}(t))= \Po(F_{[j^{\Psi(t)}]}^{\Psi(t)}(t))$. For $i> j^{\Psi(t)}$, $b_i^V(t)=0$.

\end{itemize}

For all $t> \tau_0$ the evolution of $V$ at step~$t$ is as follows, having 
already defined $\Psi(t-1)$ and  $\bbbar^V(t-1)$.

\begin{itemize}

\item Let $j = j^{\Psi(t-1)}$.

\item Let $Y_{t-1}$ be a $j$-jammed process with 
send sequence~$\p$, start time~$t-1$, 
birth rate~$\lambda$, and initial population  $\arrvect^{Y_{t-1}}(t-1) = \bbbar_{[j]}^V(t-1)$.

\item Explanation: 
Let $\zvect = \zvect^{\Psi(t-1)}$.
From the invariant, 
$\dist(\romanarrvect^{Y_{t-1}}(t-1) ) =
\Po(F_{[j]}^{\Psi(t-1)}(t-1))  
$   which, by 
Definition~\ref{def:F}, is 
$\Po( f_{[j]}^{\zerovect,  
\zvect}(t-1 - \tau^{\Psi(t-1)}))$.  
By 
Lemma~\ref{lem:f-track-jammed},
$\dist(\bbar_{[j]}^{Y_{t-1}}(t)) = 
\Po(f_{[j]}^{\zerovect, \zvect}(t-\tau^{\Psi(t-1)}))
$ which, by Definition~\ref{def:F} again, is
$
\Po(F_{[j]}^{\Psi(t-1)}(t))
$.

\item  Let $\Psi(t) = \hls(\Psi(t-1), \bbar^{Y_{t-1}}_{[j]}(t), t)$.

\item Explanation:  
Let $\Xvect \sim 
\dist(\bbar_{[j]}^{Y_{t-1}}(t)) = 
\Po(F^{\Psi(t-1)}_{[j]}(t))$
and let $A = \dist(\Xvect \mid \hls(\Psi(t-1), \Xvect, t) = \Psi(t))$.
The function $\phi_{\Psi(t-1),\Psi(t),t}$ is  
$\phi_{A(\Psi(t-1),\Psi(t),t),
\Po(F^{\Psi(t)}_{[j]}(t))}$  in Definition~\ref{def:valid} since $\hls$ is $\Psi^*$-valid.  
By construction  (from Definition~\ref{def:phi}), for all $\abar\in (\integers_{\geq 0})^j$,  $\phi_{\Psi(t-1),\Psi(t),t}(\abar) \leq \abar $.
When
$\abar \sim A$, 
the random variable 
$\phi_{\Psi(t-1),\Psi(t),t}(\avect)$
has distribution $\Po(F_{[j]}^{\Psi(t)}(t))$.

\item Let $\bbar^V_{[j]}(t) = \phi_{\Psi(t-1),\Psi(t),t}(\bbar^{Y_{t-1}}_{[j]}(t))$.

\item Explanation: It will be the case that
$\dist(\bbar^V_{[j]}(t)) = 
\Po(F_{[j]}^{\Psi(t)}(t))
$ and, for all $i \in [j]$,
$b^V_i(t) \leq b_i^{Y_{t-1}}(t)$.  

\item  For all $i\in [j]$,  let  $\b_i^{V}(t)$ be the lexicographically least
subset of $\b_i^{Y_{t-1}}(t)$ 
of size $b^V_i(t)$.

\item {For all $i > j$, let $\b^V_i(t) = \emptyset$.}
\end{itemize}
\end{definition}

We now give a short version of Definition~\ref{def:volume-process-ext},
which omits the explanation and can be used for future reference.
Definition~\ref{def:volume-process} 
is illustrated in Figure~\ref{fig:volume}.

\begin{definition}\label{def:volume-process} (short version of Definition~\ref{def:volume-process-ext}, without the explanation)
Fix a real number $\lambda \in (0,1)$ and a send sequence $\p$ with $p_0=1$. 
Fix an initialising state 
$\Psi_0 = (0, \tau_0,j,\zvect,\emptyset,\initialising)$.
Let $\hls$ be a transition rule on a set of high-level states which includes~$\Psi_0$. Let $\Psi^*$ be the $\Psi_0$-closure of~$\hls$. Suppose that $\hls$ is $\Psi^*$-valid.
A \emph{volume process} $V$ 
from $\Psi_0$ with transition rule~$\hls$, send sequence~$\p$ and birth rate~$\lambda$ is defined as follows. 

Choose the initial population as follows.
Choose 
a random vector 
$\romanarrvect = (\romanarr_1,\ldots,\romanarr_j)
\sim\Po(F_{[j]}^{\Psi_0}(\tau_0))$.
Let $\arrvect^V(\tau_0)$ be a tuple
of $j$~disjoint sets
$\arrvect^V(\tau_0) = (\arr^V_1(\tau_0),\ldots,\arr^V_j(\tau_0))$ 
where each $\arr^V_i(\tau_0) \subseteq \calB(\tau_0)$ and has size 
$\romanarr^V_i(\tau_0)$.
Let $\bbbar^V_{[j]}(\tau_0) = \arrvect^V(\tau_0)$.
For all $i>j$, let   ${\b^V_i}(\tau_0) = \emptyset$. 
Let the initial high-level state 
$\Psi(\tau_0)$ be~$\Psi_0$.

For all $t> \tau_0$ the evolution of 
$\newsym{V}{Volume process - Definition~\ref{def:volume-process}}{V}$ at step~$t$ is as follows, having 
already defined $\Psi(t-1)$ and  $\bbbar^V(t-1)$.

\begin{itemize}

\item Let $j = j^{\Psi(t-1)}$.

\item Let
$\newsym{Y-t}{part of volume process
- Definition~\ref{def:volume-process}
}{Y_{t-1}}$ be a $j$-jammed process with 
send sequence~$\p$, start time~$t-1$, 
birth rate~$\lambda$, and initial population  $\arr^{Y_{t-1}}(t-1) = \bbbar_{[j]}^V(t-1)$.

\item  Let $\newsym{Psi-t}{high-level state determined by volume process - Definition~\ref{def:volume-process} }
{\Psi(t) }
= \hls(\Psi(t-1), \bbar^{Y_{t-1}}_{[j]}(t), t)$.

\item Let $\bbar^V_{[j]}(t) = \phi_{\Psi(t-1),\Psi(t),t}(\bbar^{Y_{t-1}}_{[j]}(t))$.

\item  For all $i\in [j]$,  let  $\b_i^{V}(t)$ be the lexicographically least
subset of $\b_i^{Y_{t-1}}(t)$ 
of size $b^V_i(t)$.

\item {For all $i > j$, let $\b^V_i(t) = \emptyset$.}
\end{itemize}
\end{definition}

Lemma~\ref{lem:dist-volume} gives a rigorous version of the explanation in the extended version (Definition~\ref{def:volume-process-ext})
of Definition~\ref{def:volume-process}.
In particular, Lemma~\ref{lem:dist-volume} pins down the conditioning that we rely on.
Note that
item (i) in the lemma statement implies 
$\dist(\bbar_{[j^{\psi_t}]}^{Y_{t}}(t+1) \mid \Psivect_{T(t)} = \psivect_{T(t)})  = \dist(\bbar_{[j^{\psi_t}]}^{Y_{t}}(t+1) \mid \Psi(t)= \psi_t)$.
Similarly, item (ii) implies
$\dist(\bbar^V_{[j^{\psi_t}]}(t) \mid \Psivect_{T(t)} = \psivect_{T(t)})  =
\dist(\bbar^V_{[j^{\psi_t}]}(t) \mid  \Psi(t)=\psi_t)$.

\begin{lemma}\label{lem:dist-volume}
Fix a real number $\lambda \in (0,1)$ and a send sequence $\p$ with $p_0=1$. 
Fix an initialising state 
$\Psi_0 = (0, \tau_0,j,\zvect,\emptyset,\initialising)$.
Let $\hls$ be a transition rule on a set of high-level states which includes~$\Psi_0$. Let $\Psi^*$ be the $\Psi_0$-closure of~$\hls$. Suppose that $\hls$ is $\Psi^*$-valid.

Let $V$ be a volume process from~$\Psi_0$ with transition rule~$\hls$, send sequence~$\p$ and birth rate~$\lambda$. 
For any
$t \ge \tau_0$, 
let $Y_t$ be the process defined in Definition~\ref{def:volume-process} and
let $T(t) = \{\tau_0,\ldots,t\}$.

 Then the following statements hold for
 every $t\geq \tau_0$
 and every tuple
$ 
\psivect_{T(t)} = (\psi_{\tau_0},\dots,\psi_{t})$  of high-level states 
such that $\Failure\notin \{\type^{\psi_{t'}}
: t'\in T(t)\}$
and
$\Pr( \Psivect_{T(t)} = \psivect_{T(t)}) > 0$. 
\begin{enumerate}[(i)]
\item $\dist(\bbar_{[j^{\psi_t}]}^{Y_{t}}(t+1) \mid \Psivect_{T(t)} = \psivect_{T(t)})  = \Po(F_{[j^{\psi_t}]}^{\psi_{t}}(t+1))$.
\item $\dist(\bbar^V_{[j^{\psi_t}]}(t) \mid \Psivect_{T(t)} = \psivect_{T(t)})  = \Po(F_{[j^{\psi_t}]}^{\psi_t}(t))$.
\end{enumerate}
\end{lemma}

\begin{proof}
For brevity, let $\calE_{t}$ denote the event $\Psivect_{T(t)} = \psivect_{T(t)}$.
We will prove (ii) for all $t\geq \tau_0$ by induction on~$t$. Along the way, we will prove that (ii)   implies (i).

For the base case note that (ii) holds for $t=\tau_0$ by the definition of the volume process (Definition~\ref{def:volume-process}). 
Fix $t \ge \tau_0$, 
let $j= j^{\psi_t}$ 
and $\zvect = \zvect^{\psi_t}$
and suppose (for induction) that (ii) holds 
at time~$t$, meaning that
\[
\dist(\bbar^V_{[j]}(t) \mid \calE_{t}) = 
\Po(F_{[j]}^{\psi_{t}}(t)).
\]
From Definition~\ref{def:volume-process}, 
$Y_t$ is the $j$-jammed process with send sequence~$\p$, start time~$t$, birth rate~$\lambda$, and initial population 
$\bbbar_{[j]}^V(t)$, which, conditioned on $\calE_t$,  
has distribution $ \Po(F_{[j]}^{\psi_{t}}(t))$.
By 
Definition~\ref{def:F}, this is 
$\Po( f_{[j]}^{\zerovect,  
\zvect}(t - \tau^{\psi_t}))$. 
We now apply  
Lemma~\ref{lem:f-track-jammed}
to the $j$-jammed process~$Y_t$,
noting that the $\tau$ of Lemma~\ref{lem:f-track-jammed} is our $t$ and taking the $\Delta$ of Lemma~\ref{lem:f-track-jammed} to be~$t-\tau^{\psi_t}$.
By Lemma~\ref{lem:f-track-jammed} applied to time $t+1$,
$\dist(\bbar_{[j]}^{Y_{t}}(t+1) \mid \calE_t) = 
\Po(f_{[j]}^{\zerovect, \zvect}(t+1-\tau^{\psi_t}))
$.
It follows from Definition~\ref{def:F} that
\begin{equation}\label{eq:blabla}
\dist(\bbar^{Y_{t}}_{[j]}(t+1) \mid \calE_{t}) = \Po(F_{[j]}^{\psi_{t}}(t+1)),
\end{equation}
as required by (i).  

To finish the inductive proof of (ii), we will show (ii) at step~$t+1$. 
For brevity, let   $\Xvect = \bbar^{Y_{t}}_{[j]}(t+1)$.
From \eqref{eq:blabla},
$\dist(\Xvect \mid \calE_t) = 
\Po(F_{[j]}^{\psi_{t}}(t+1))$.
From Definition~\ref{def:volume-process},
$\Psi(t+1) = \hls(\psi_t, \Xvect,t+1)$.
So consider $\psivect_{T(t+1)}$
extending $\psivect_{T(t)}$ and suppose
$\Psi(t+1) = \psi_{t+1} $.

Since $\hls$ is valid, 
from Definition~\ref{def:valid} applied with our $\Xvect$ conditioned on $
\calE_t$,  
\[A(\psi_t,\psi_{t+1},t+1) \gtrsim \Po(F^{\psi_{t+1}
}_{[j]}(t+1))\] and
the function 
$\phi_{\psi_t,\psi_{t+1},t+1}$ from Definition~\ref{def:valid} is  
$\phi_{A(\psi_t,\psi_{t+1},t+1),\Po(F^{\psi_{t+1}}_{[j]}(t+1))}$.
By the construction of   $\bbar^V_{[j]}(t+1)$ in Definition~\ref{def:volume-process},
\begin{align}\nonumber
\dist(\bbar^V_{[j]}(t+1) \mid \calE_{t}, \Psi(t+1) = \psi_{t+1}) &= 
\dist(\phi_{\psi_{t},\psi_{t+1},t+1}(\Xvect) \mid \calE_t, \Psi(t+1) = \psi_{t+1})\\\label{eq:blabla-b}
&= \dist(\phi_{\psi_{t},\psi_{t+1},t+1}(\Xvect) \mid \calE_t, \hls(\psi_{t}, \Xvect, t+1) = \psi_{t+1}).
\end{align}
By~\eqref{eq:blabla}, $ 
\dist(\Xvect\mid \calE_t) = 
\Po(F_{[j]}^{\psi_{t}}(t+1))$, and the event $\{\hls(\psi_{t}, \Xvect, t+1) = \psi_{t+1}\}$ lies in the $\sigma$-field of~$\Xvect$, $\sigma(\Xvect)$. Recalling the definition of~$A(\psi_{t},\psi_{t+1},t+1)$ from Definition~\ref{def:valid} again applied with our $\Xvect$ conditioned on $\calE_t$,  
\begin{equation}\label{eq:blabla-c}
A(\psi_{t},\psi_{t+1},t+1) = 
\dist(\Xvect \mid \calE_t, \hls(\psi_{t}, \Xvect, t+1) = \psi_{t+1}).
\end{equation}
We've seen that
$\phi_{\psi_t,\psi_{t+1},t+1}$ 
is $\phi_{ A(\psi_t,\psi_{t+1},t+1),
\Po(F^{\psi_{t+1}}_{[j]}(t+1))}$.
By construction  (from Definition~\ref{def:phi}), for all $\abar\in (\integers_{\geq 0})^j$,  $\phi_{\psi_t,\psi_{t+1},t+1}(\abar) \leq \abar $.
When
$\abar \sim A(\psi_t,\psi_{t+1},t+1)$, 
the random variable 
$\phi_{\psi_t,\psi_{t+1},t+1}(\avect)$
has distribution $\Po(F_{[j]}^{\psi_{t+1
}}(t+1))$.
Thus by~\eqref{eq:blabla-c},
\[
\dist(\phi_{\psi_{t},\psi_{t+1},t+1}(\Xvect) \mid \calE_{t}, \hls(\psi_{t}, \Xvect, t+1) = \psi_{t+1}) =  \Po(\F^{\psi_{t+1}}_{[j]}(t+1)),
\]
so it follows from~\eqref{eq:blabla-b} that
\begin{align}\label{eq:volume-invariant-proof-1}
\dist(\bbar^V_{[j]}(t+1) \mid \calE_{t}, \Psi(t+1) = \psi_{t+1}) &= \Po(\F_{[j]}^{\psi_{t+1}}(t+1)).
\end{align}
This establishes (ii) at step~$t+1$ 
if $j^{\psi_{t+1}} = j$.

Otherwise, by (T\ref{item:trans1}), 
$\tau^{\psi_{t+1}} = t+1$ 
and by (T\ref{item:trans4}) 
$j^{\psi_{t+1}} \geq j$ and
for all $\ell \in [j^{\psi_{t+1}}] \setminus [j]$,
$z^{\psi_{t+1}}_\ell = 0$.
It remains to show that, for all $\ell \in 
[j^{\psi_{t+1}}] \setminus [j]$,
$b_\ell^V(t+1)=0$ and 
$F_\ell^{\psi_{t+1}}(t+1) = 0$.
The first of these is immediate from Definition~\ref{def:volume-process}.
To show 
$F_\ell^{\psi_{t+1}}(t+1) = 0$
we use Definition~\ref{def:F}
which states that
$F_\ell^{\psi_{t+1}}(t+1) = 
f_\ell^{\zerovect,\zvect^{\psi_{t+1}}}(t+1- \tau^{\psi_{t+1}}) = z_\ell^{\psi_{t+1}} = 0$, as required.
\end{proof}

\section{Escape processes}\label{sec:def-escape}

As explained in the introduction, our proof strategy is to couple the backoff process with the volume process~$V$ 
from Definition~\ref{def:volume-process}
(this tracks the sends of the backoff process), and with an infinite number of escape processes, which track escapes.
A new escape process starts every time the high-level state changes. The invariant is that every ball present in the volume process, but not in the relevant escape process, is also present in the original backoff process. The next step is to define the notion of escape process.

\begin{definition}
\label{def:escape}
Fix a send sequence~$\p$ with $p_0=1$, a $\lambda \in (0,1)$, and a high-level state $\Psi= (g,\tau,j,\zvect,\calS,\type)$.

Fix a (possibly random) tuple $\arrvect^E(\tau)$ of
disjoint sets
$\arrvect^E(\tau) = (\arr_1^E(\tau),\arr_2^E(\tau),\ldots)$
where each $ \arr^E_{i}(\tau)\subseteq \ball(\tau)$ and, 
for $i>j$, $\arr_i^E(\tau) = \emptyset$. 

An escape process \newsym{E}{Escape process}{E} with 
high-level state $\Psi$, send sequence~$\p$, birth rate~$\lambda$, and
initial population $\arrvect^E(\tau)$
is the generalised backoff process with 
send sequence~$\p$, start time~$\tau$, the birth distribution~$\calD$ that always returns~$0$,  initial population~$\arrvect^E(\tau)$, and which satisfies the following constraints.

For all $t\geq \tau$, define 
$\mathbb{S}^E(t) := \{S \in \calS\colon \calN^E_S(t+1) \le \lambda |S|/80\}$.
Let
$\hat{S}^E(t) := 
\cup_{S \in \mathbb{S}^E(t)} S$. Finally,
define
\[
\Xi^E(t) := \begin{cases}
\lambda |\hat{S}^E(t)|/80 & \mbox{ if }|\hat{S}^E(t)| \ge 100/\lambda^2,\\
0 & \mbox{ otherwise,}
\end{cases}
\]
and $\newsym{nu-E}{end-of-step arrival rate of escape process}{\nu^E(t)} := \exp(-\Xi^E(t)/192)$.
For $t>\tau$, $\arrvect^E(t)$ and $\evect^E(t)$ are defined as follows. 
\begin{itemize}
\item  For each $i\in [j]$, 
$\romanarr_i^E(t)$ is a Bernoulli r.v.~with parameter~$\nu^E(t-1)$. These are independent for different~$i$ and for different~$t$. For all $i>j$,  $\romanarr^E_i(t) = 0$.
\item For all $i \in [j-1]$,   $\e_i^E(t) = \emptyset$. For all $i \ge j$, $\e^E_i(t) = \s^E_i(t)$.
\end{itemize}

Note that we have not specified the name of the ball in $\arr_i^E(t)$, in the event that $\romanarr_i^E(t)=1$.

\end{definition}

Definition~\ref{def:escape} needs some explanation.
In the coupling that we will describe in Section~\ref{sec:VEB-couple}, a new escape process will start every time the volume process changes the high-level state.
The initial population of the escape process comes directly from the population of the previous escape process. The birth rate~$\lambda$ is present as a parameter because the rate of end-of-step arrivals depends on~$\lambda$, but there are no births in an escape process. Instead, balls arrive as end-of-step arrivals - these correspond to escapes in the backoff process.  
In Definition~\ref{def:escape}, $\nu^E(t)$ is defined to give an upper bound on the probability that a ball escapes from any give bin~$i\in [j]$ during step~$t$ of the backoff process. Its definition
roughly captures the probability of simultaneous silence 
from each set of bins~$S \in \calS$. The idea is that for each $S$, (T\ref{item:trans3}) of Definition~\ref{def:transition-rule} guarantees that the expected noise in the volume process 
from balls in $S$ is at least $\lambda|S|/40$. 
As long as these bins are not too overloaded in the escape process 
(so that they get included in $\mathbb{S}^E(t)$) this 
expected volume noise also corresponds to expected noise
at least $\lambda|S|/80$ in the backoff process, which is accounted for in the
upper bound $\nu^E(t)$ on the probability of escape.
If the bins in~$S$ are already too overloaded in~$E$ (captured by $\calN_S^E(t+1)>\lambda|S|/80$) then these balls just contribute the trivial upper bound $\nu^E(t) \leq 1$. Apart from the end-of-step arrivals, escape process is similar to a $j$-jammed process -- this is because there is no need to track balls after bin~$j$.

\section{The Volume-Escape-Backoff (VEB) Coupling}
\label{sec:VEB-couple}

We are now ready to describe the volume-escape-backoff coupling (Definition~\ref{def:VEB-couple}). 
In Definition~\ref{def:VEB-couple} it is helpful to recall (from Definition~\ref{def:F})
that every high-level state~$\psi=(g,\tau,j,\zvect,\calS,\type)$ has the property that  
$F_{[j]}^{\psi}(\tau) = f^{\zerovect,\zvect}_{[j]}(0) = \zvect$. In the coupling we restrict $\lambda$ to the range $(0,1/120)$. The reason for this is that it makes it easy to use Lemma~\ref{lem:single-send-bound}. It is without loss of generality since our goal (Theorem~\ref{thm:goaljkiller}) allows us to set any upper bound on~$\lambda$.
Recall from Definition~\ref{def:Tprot} that 
$T^{\p,j}  = 80 j^2 \lfloor \sum_{x\in [j]} W_x \rfloor$.

\begin{definition}
\label{def:VEB-couple}  
Fix $\lambda \in (0,1/120)$. Let $\p$ be a send sequence with $p_0=1$. 
Let $j_0$ be a positive integer.
Let $X$ be a backoff process 
with 
send sequence~$\p$, 
birth distribution $\Po(2 \lambda)$, and cohort set~$\calC=\{A,B\}$ (Definition~\ref{def:backoff-cohorts}).
Let $\Einit$ be any event determined by 
$\{n^{X^B}(t') : t'\in [T^{\p,j_0}]\}$
which entails that each 
$n^{X^B}(t')$ in this set is at least~$1$.
We will condition on $\Einit$.

Let
$\tau_0 =  T^{\p,j_0}$. 
Let
$\Psi_0$ be the  initialising high-level state   
with $\tau^{\Psi_0} = \tau_0$ and $j^{\Psi_0} = j_0$. 
Let $\hls$ be a transition rule on a set of high-level states which includes~$\Psi_0$. 
Let $\Psi^*$ be the $\Psi_0$-closure of~$\hls$. Suppose that $\hls$ is $\Psi^*$-valid.

The VEB coupling couples $X$ with a volume process~$V$ from $\Psi_0$
with transition rule~$\hls$, send sequence~$\p$ and birth rate~$\lambda$
and with a sequence $E_1,E_2,\ldots$ of escape processes with send sequence~$\p$ and birth rate~$\lambda$.

From the definition of volume process (Definition~\ref{def:volume-process}), 
$\Psi(\tau_0) = \Psi_0$. Due to event~$\Einit$, 
the balls in in $\bbar^{X^A}_{[j_0]}$ behave as a $j_0$-jammed process during steps 
in~$[\tau_0]$. By Lemma~\ref{lem:slowfill-j-jammed}, 
$\dist(\bvect^{X^A}_{[j_0]}(\tau_0)) \gtrsim \Po(\zvect^{\Psi_0})$.
Let $\varphi$ be the function
$\varphi_{ \dist(\bvect^{X^A}_{[j_0]}(\tau_0))
,\Po(z^{\Psi_0})}$ from Definition~\ref{def:phi}. 
Let  
$\romanarrvect =  
\varphi (\bvect^{X^A}_{[j]}(\tau_0))$
so that 
$\dist(\romanarrvect) = \Po(\zvect^{\Psi_0})$
and $\bvect^{X^A}_{[j]}(\tau_0)\geq
\romanarrvect$.
For each $i\in [j_0]$
let $b_{i}^V(\tau_0) = \arr^V_{i}(\tau_0)$ be a subset of $\b_i^{X^A}(\tau_0)$ of size 
$\romanarr_i$.

For each $t > \tau_0$ let $Y_{t-1}$ 
be a $j^{\Psi(t-1)}$-jammed process with send sequence~$\p$,
start time~$t-1$, 
birth rate $\lambda$, and initial population $\arr^{Y_{t-1}}(t-1) = \bbbar_{[j^{\Psi(t-1)}]}^V(t-1)$, as defined in the definition of a volume process (Definition~\ref{def:volume-process}).

For every positive integer~$g$, let
$\tau_g = \min \{t > \tau_0 \mid g^{\Psi(t)} = g\}$. Note that $\tau_g$ is a random variable, depending on the evolution of~$V$.  Also by (T\ref{item:trans5}) and (T\ref{item:trans1}), $\tau_1 = \tau_0+1$. By (T\ref{item:trans4}), $j^{\Psi(\tau_1)} = j_0$.

For $i\in [j_0]$ we will use the notation $\b_i^{E_0}(\tau_1)$ as a synonym for~$\e_i^X(\tau_1)$. 
For each~$g\geq 1$, the escape process $E_g$ has high-level state $\Psi(\tau_g)$. The end-of-step arrivals of~$E_g$ are as specified in the Definition of escape process (Definition~\ref{def:escape})
with the additional constraint that for all   $i \in [j^{\Psi(\tau_g)}]$,
$\b^{E_g}_i(\tau_g) = \arr^{E_g}_i(\tau_g)  = \b_i^{E_{g-1}}(\tau_g)$.

{\bf Invariant:}  The coupling maintains the following invariants for all $t\geq \tau_0$.

\begin{enumerate}[({I}1)]
\item  
Either 
$\type^{\Psi(t)} \in \{\initialising,\Failure\}$ or, 
for all $S\in \calS^{\Psi(t)}$,
 $
    \sum_{x\in S} p_x b_x^{Y_{t-1}}(t) \geq \lambda |S|/40$.
\item Either $\type^{\Psi(t)} \in \{\initialising,\Failure\}$ or, for all   $i\geq 1$, 
$\b_i^{Y_{t-1}}(t) \setminus \b_i^{E_{g^{\Psi(t-1)}}}(t) \subseteq 
\b_i^X(t)$.
\end{enumerate}
 
Note that (I1) and (I2) are true for $t=\tau_0$ since $\type^{\Psi(\tau_0)} = \initialising$.
Before doing the general step of the coupling we explain the coupling at step $\tau_1 = \tau_0+1$, and show how (I1) and (I2) are established. We have already defined the
$j_0$-jammed process $Y_{\tau_0}$. We  take $\n^{X^A}(\tau_1) = \n^{Y_{\tau_0}}(\tau_1)$ (and choose $\n^{X^B}(\tau_1)$ independently)
so that, for all $i\geq 0$, 
$\B_i^{Y_{\tau_0}}(\tau_1) \subseteq \B_i^X(\tau_1)$.
The sends in $Y_{\tau_0}$ and $X$ are coordinated as follows. 
For all non-negative integers~$i$ and
all balls~$\beta \in \B_i^{Y_{\tau_0}}(\tau_1)  $,
$\beta$ is either in  both of
$\s_i^{Y_{\tau_0}}(\tau_1) $ and $ \s_i^{X}(\tau_1)$
or in 
neither of them. The escapes in $Y_{\tau_0}$ at step~$\tau_1$ are deterministic (balls only escape when they send from bin~$j_0$). 
The escapes in $X$ at step~$\tau_1$ are according to the process (at most one ball escapes).
There are no end-of-step arrivals in~$X$ or $Y_{\tau_0}$ at step~$\tau_1$. Then $\Psi(\tau_1)$
is chosen via
$\Psi(\tau_1) = \hls(\Psi_0, \bbar_{[j_0]}^{Y_{\tau_0}}(\tau_1),\tau_1)$. 
As we have already noted, $j^{\Psi(\tau_1)} = j_0$.
$\bbbar^V_{[j_0]}(\tau_1)$
is chosen according to the volume process
(for each $i\in [j_0]$, 
$\b_i^V(\tau_1)$ is a subset of 
$\b^{Y_{\tau_0}}_i(\tau_1)$, and
for $i>j_0$, 
$\b_i^V(\tau_1)=\emptyset$). 
As noted, we define 
$\b_i^{E_0}(\tau_1) =  \e_i^X(\tau_1)$.
Invariant (I2) at time~$\tau_1$ is by construction since, for all 
$i\geq 1$, $\b_i^{Y_{\tau_0}}(\tau_1) \setminus \b_i^{E_0}(\tau_1) \subseteq \b_i^X(\tau_1)$.
Invariant (I1)  comes from (T\ref{item:trans3}).
Escape process $E_1$ has high-level state $\Psi(\tau_1)$. Following the general pattern that we have already described, for 
$i \in [j_0]$,
$\b^{E_1}_i(\tau_1) = \arr^{E_1}_i(\tau_1)  = \b_i^{E_{0}}(\tau_1)$.

We now do the general step of the coupling.
Before step~$t$, for some $t>\tau_1$,
we have $\Psi(t-1)$, $\bbbar^V(t-1) =\bbbar^{Y_{t-1}}(t-1)$, $\bbbar^{X}(t-1)$, and 
$\bbbar^{E_{g^{\Psi(t-1)}}}(t-1)$ which  satisfies
$\bbbar^{E_{g^{\Psi(t-1)}}}(t-1) =
\bbbar^{E_{g^{\Psi(t-2)}}}(t-1)
$.  
If $\type^{\Psi(t-1)} = \Failure$ then all processes evolve independently at step~$t$.
Otherwise, the coupling proceeds as follows at 
step~$t$.

\begin{itemize}

\item Let $\psi = \Psi(t-1)$, 
$j = j^{\psi}$, and $\nu = \nu^{E_{g^{\psi}}}(t-1)$.

\item $\n^{X^A}(t) = \n^{Y_{t-1}}(t)$. Choose $\n^{X^B}(t)$ independently.

\item 
For all $g$ such that $1 \leq g \leq g^{\psi} -1$, $E_g$ 
evolves at step~$t$ independently of the other processes.

\item The sends in $Y_{t-1}$, $X$, and $E_{g^{\psi}}$
are coordinated in the following way at step~$t$. 
If, for any non-negative integer~$i$ and any
$\Phi \subseteq \{ Y_{t-1}, X, E_{g^{\psi}}\}$,
a ball $\beta$ is in 
$\cap_{U\in \Phi} \B_i^U(t)$
then $\beta$ is either in 
$\cap_{U\in \Phi} \s_i^U(t)$
or in 
$\cap_{U\in \Phi} \B_i^U(t)\setminus \s_i^U(t)$, according to the correct marginal distribution.

\item The escapes in $Y_{t-1}$ 
and $E_{g^{\psi}}$ at step~$t$ 
are deterministic (Definition~\ref{def:j-jammed} and \ref{def:escape}) -- balls do not escape, except that they always escape when they send from bin~$j$.
There are no end-of-step arrivals in~$X$ or $Y_{t-1}$ at step~$t$.  

\item  
We next coordinate the 
end-of-step arrivals in
$E_{g^{\psi}}$ at step~$t$ (Definition~\ref{def:escape} -- for each $i\in [j]$ an independent Bernoulli arrival in bin~$i$ with parameter~$\nu$)   with the
escapes in~$X$ 
at step~$t$ (Definition~\ref{def:backoff}).

 \begin{itemize}
\item Let $\calP_{\mathrm{Ber}}(t)$ be a distribution on $j$-tuples of independent Bernoulli random variables with parameter~$\nu$. Let $\calP_{\Esc}(t)$ be the distribution of the $j$-tuple
$\ebar^X_{[j]}(t)$ conditioned on the whole coupling up to and including step~$t-1$, that is on the values of
$\Psi(t-1)$, $\bbbar^V(t-1) =\bbbar^{Y_{t-1}}(t-1)$, $\bbbar^{X}(t-1)$, and 
$\bbbar^{E_{g^{\psi}}}(t-1)$.
We will show $\calP_{\mathrm{Ber}}(t)\gtrsim \calP_{\Esc}(t)$.

\item 
Let $\barr \sim \calP_{\mathrm{Ber}}(t)$. By Observation~\ref{obs:dominate-suff-cond} it suffices to show that for all $j$-tuples $\bara$ of Boolean values,
$\Pr(\barr \geq \bara) \geq \Pr(\ebar^X_{[j]}(t) \geq \bara)$. 
Since otherwise the right-hand probability is~$0$, we may assume that $\bara$ has at most one entry that is~$1$. To avoid trivialities (where both probabilities are~$1$) we may assume that $\bara$ has exactly one entry that is~$1$.
Thus, our goal is to show, for all $i\in [j]$, that $\Pr(e_i^X(t)=1) \leq \nu$.
From the definition of 
$\nu = \nu^{E_{g^{\psi}}}(t-1)$  in Definition~\ref{def:escape}
this is vacuous unless $|\hatS^{E_{g^{\psi}}}(t-1)|\geq 100/\lambda^2$, so suppose that this is the case.
Fix $i\in [j]$.  Recall that the type of $\psi$ is not $\initialising$ or $\Failure$ (since $t>\tau_1$ and we have already specified that the processes evolve independently after a failure type is reached).

\begin{itemize}
\item From (I1) at $t-1$, for all $S\in \calS^{\psi}$, 
$ \sum_{x\in S} p_x b_x^{Y_{t-2}}(t-1) \geq \lambda |S|/40$.

\item From (I2) at $t-1$,  for all $x \geq 1$,
$
\b_x^{Y_{t-2}}(t-1) \setminus \b_x^{E_{g^{\Psi(t-2)}}}(t-1) \subseteq \b_x^X(t-1)$.

\item Recall from the information that we collected just before step~$t$ that  
$\bbbar^{E_{g^{\Psi(t-2)}}}(t-1) = \bbbar^{E_{g^{\psi}}}(t-1)$. Thus, for all $x\geq 1$, 
\begin{equation}
\label{eq:keep}
\b_x^{Y_{t-2}}(t-1) \setminus \b_x^{E_{g^{\psi}}}(t-1) \subseteq \b_x^X(t-1).
\end{equation}
\item   It follows that for all $S \in \calS^{\psi}$,
\[
 \sum_{x\in S}p_xb_x^X(t-1) \ge \sum_{x \in S}p_x \big(b_x^{Y_{t-2}}(t-1) - b_x^{E_{g^{\psi}}}(t-1)\big) \ge \lambda|S|/40 - \sum_{x \in S} p_xb_x^{E_{g^{\psi}}}(t-1).
\]

\item Let $\mathbb{S} = 
\mathbb{S}^{E_{g^{\psi}}}(t-1)
$ and $\hat{S} = 
\hat{S}^{E_{g^{\psi}}}(t-1)
$ from Definition~\ref{def:escape}. 
$\mathbb{S}$ is a subset of the sets in $\calS^{\psi}$ and these sets are disjoint by the definition of a high-level state.
$\hat{S}$ is the union of the sets in $\mathbb{S}$. 
Thus
the following holds, where the second inequality is from the definition of $\mathbb{S}$ which guarantees that
these sets have $\calN_S^{E_{g^{\psi}}}(t) = \sum_{x\in S} p_x b_x^{E_{g^\psi}}(t-1) \leq \lambda |S|/80$.

\begin{align*} \sum_{S \in \mathbb{S}}\sum_{x\in S}p_xb_x^X(t-1) &\ge \frac{\lambda |\hat{S}|}{40} - \sum_{S \in \mathbb{S}} \sum_{x \in S} p_xb_x^{E_{g^{\psi}}}(t-1)\\
 &=\frac{\lambda |\hat{S}|}{40} - 
 \sum_{S \in \mathbb{S}} 
 \calN_S^{E_{g^\psi}}(t) \\
&\geq \frac{\lambda |\hat{S}|}{80} 
=  \Xi^{E_{g^{\psi}}}(t-1).\end{align*}

\item  We apply Lemma~\ref{lem:single-send-bound} 
with the generalised backoff process $Z$ of Lemma~\ref{lem:single-send-bound} as~$X$ (started from time $t-1$). The $j$-tuple 
$\avect$ of Lemma~\ref{lem:single-send-bound} 
is  $\bbar_{[j]}^{X}(t-1)$.
We take the set $S$ from Lemma~\ref{lem:single-send-bound}
to be~$\hat{S}$.
The pre-conditions are that $|\hat{S}| \geq 100/\lambda^2$
and $\sum_{i\in \hatS} p_i b_i^X(t-1) \geq \lambda |\hatS|/80$, both of which we have established.  
The lemma shows that the probability $q$ that  
$\sum_{i\geq 1} s_i^X(t) \leq 1$
is at most 
\[\exp\big(-\sum_{x \in \hat{S}} p_x b_x^X(t-1)/16\big)
\leq 
\exp(- \Xi^{E_{g^{\psi}}}(t-1)/16) \leq
\exp(- \Xi^{E_{g^{\psi}}}(t-1)/192) =  \nu\] as desired.

\end{itemize} 

\item So we have now shown $\calP_{\mathrm{Ber}}(t)\gtrsim \calP_{\Esc}(t)$.

\item  
Let $
\romanarrvect^{E_{g^\psi}}(t)
= \phi^{-1}_{\calP_{\mathrm{Ber}}(t),\calP_{\Esc}(t)}(\bare^X_{[j]})$. By construction, $ 
\romanarrvect^{E_{g^\psi}}(t)
\geq \bare^X_{[j]}(t)$.

\item Choose the sets of balls in 
$\arrvect^{E_{g^\psi}}(t)$ as follows.
If $e_i^X(t)=1$ then
$\arr_i^{E_{g^{\psi}}}(t) = \e_i^X(t)$. 
If $e_i^X(t)=0$
and $a_i=1$ then 
$\arr_i^{E_{g^{\psi}}}(t) =  \{\calB(t,i,1)\}$.

\item In order to conform to Definition~\ref{def:genbackoff} we are obligated to demonstrate that
$\arr_i^{E_{g^{\psi}}}(t)$
is disjoint from $\n^{E_{g^{\psi}}}(t)$ (which is empty, so that is easy!) and also from every set $\b_{i'}^{E_{g^{\psi}}}(t-1)$. When $e_i^X(t)=0$ this is clear from the name $\calB(t,i,1)$. Otherwise the relevant ball is in $\e_i^X(t)$ so it can't be in any $\e_{i'}^X(t')$ for any $t' < t$ so can't be in $\b_{i'}^{E_{g^{\psi}}}(t')$.

\end{itemize}

\item From the definition of the volume process,  $\Psi(t)= \hls(\psi, \bbar^{Y_{t-1}}_{[j]}(t),t)$.

\item   The determination of $\bbbar^V(t)$ is deterministic from Definition~\ref{def:volume-process}. Finally, if
$\Psi(t)=\psi$
then it is obvious that
$\bbbar^{E_{g^{\Psi(t)}}}(t) =
\bbbar^{E_{g^{\psi}}}(t)
$. Otherwise, 
let $g = g^{\Psi(t)}$. By (T1), $g^{\psi} = g-1$. Then
$\tau_{g} = t$
so using our specified  end-of-step arrivals
we have
$\b^{E_g}_i(t)=
\arr^{E_g}_i(t)  = \b_i^{E_{g-1}}(t)$ which gives the guarantee 
$\bbbar^{E_{g}}(t) =
\bbbar^{E_{g-1}}(t)
$, as required.
\item It remains to establish (I1) and (I2) for step~$t$
for the case where $\type^{\Psi(t)}\neq \Failure$.

\item For (I1) we wish to show that 
for all $S\in \calS^{\Psi(t)}$,
 $
    \sum_{x\in S} p_x b_x^{Y_{t-1}}(t) \geq \lambda |S|/40$.
This follows from (T\ref{item:trans3}).

\item For (I2) we wish to show that, for all   $i\geq 1$, 
$\b_i^{Y_{t-1}}(t) \setminus \b_i^{E_{g^{\psi}}}(t) \subseteq 
\b_i^X(t)$.  
Since  
the sends in $Y_{t-1}$, $X$ and $E_{g^{\psi}}$ are coordinated 
(so any ball that is in multiple processes is in the same bin in these processes), it suffices to prove that
$
\b^{Y_{t-1}}(t) \setminus \b^{E_{g^{\psi}}}(t) \subseteq \b^X(t)$.
We proceed as follows.

\begin{itemize}

\item Applying the definition of $j$-jammed process,
\begin{equation}\label{eq:VEB-volume}
\b^{Y_{t-1}}(t) = (\b^{Y_{t-1}}(t-1) \cup \n^{Y_{t-1}}(t)) \setminus \s^{Y_{t-1}}_{j}(t). 
\end{equation}

\item Applying the definition of the escape process 
(including the fact that $\b_i^{E_{g^{\psi}}}(t)$
is empty for $i>j$),
\[
\b^{E_{g^{\psi}}}(t) = 
\bigcup_{i \in [j]}\arr_i^{E_{g^{\psi}}}(t)\cup 
\big(
\b^{E_{g^{\psi}}}(t-1) \setminus \s^{E_{g^{\psi}}}_{j}(t)   \big) .
\]
By the properties of our coupling that we have already established, for every $i\in [j]$, $\arr_i^{E_{g^{\psi}}}(t) \supseteq \e_i^X(t)$
so 
\begin{align}\nonumber
\b^{E_{g^{\psi}}}(t) &\supseteq
\e^X(t)\cup 
\big(
\b^{E_{g^{\psi}}}(t-1) \setminus 
\s^{E_{g^{\psi}}}_{j}(t)   \big)\\\label{eq:VEB-escape}
&\supseteq \big(\e^X(t)\cup \b^{E_{g^{\psi}}}(t-1)\big) \setminus \s^{E_{g^{\psi}}}_{j}(t).
\end{align}

\item Combining~\eqref{eq:VEB-volume} and~\eqref{eq:VEB-escape}, we obtain 
\[
\b^{Y_{t-1}}(t) \setminus 
\b^{E_{g^{\psi}}}(t) \subseteq 
(U_1 \setminus U_2 ) \setminus (U_3 \setminus U_4)\]
where 
$U_1=  \big(\b^{Y_{t-1}}(t-1) \cup 
\n^{Y_{t-1}}(t)\big) $, 
$U_2 =  \s^{Y_{t-1}}_{j}(t)$,
$U_3 =  \big(\e^X(t)\cup 
\b^{E_{g^{\psi}}}(t-1)\big)$ and 
$U_4 = 
\s^{E_{g^{\psi}}}_{j}(t)$.
We first wish to simplify the right-hand side
to $(U_1 \setminus U_2) \setminus U_3$. For this,
we wish to show that any balls in $U_1$ and $U_4$ are also in $U_2$. This holds by our synchronisation of sends.
Any ball in $U_4$ sends from bin~$j$ at time~$t$ in $E_{g^\psi}$. If it is in $U_1$ then it is present in $Y_{t-1}$ just before step~$t$ so it sends from bin~$j$ at time~$t$ in $Y_{t-1}$ meaning that it is in~$U_2$.
Thus we now have the right-hand-side
$(U_1 \setminus U_2) \setminus U_3$.
This is trivially contained in $U_1 \setminus U_3$.
Thus, we have shown 
\[
\b^{Y_{t-1}}(t) \setminus 
\b^{E_{g^{\psi}}}(t) \subseteq 
\big(\b^{Y_{t-1}}(t-1) 
\cup \n^{Y_{t-1}}(t)\big) 
\setminus \big(\e^X(t)\cup 
\b^{E_{g^{\psi}}}(t-1)\big).
\]

\item Applying the definition of the volume process, $\b^{Y_{t-1}}(t-1) = \b^V(t-1) \subseteq 
\b^{Y_{t-2}}(t-1)$.
Also, as we have set it up in the coupling,
$\n^{Y_{t-1}}(t) \subseteq \n^X(t)$. 
Therefore,
\[
\b^{Y_{t-1}}(t) \setminus \b^{E_{g^{\psi}}}(t) \subseteq \big(\b^{Y_{t-2}}(t-1) \cup \n^{X}(t)\big) \setminus \big(\e^X(t)\cup \b^{E_{g^{\psi}}}(t-1)\big).
\]

\item  
From~\eqref{eq:keep}, 
$\b^{Y_{t-2}}(t-1) \setminus \b^{E_{g^{\psi}}}(t-1) \subseteq \b^X(t-1)$.
So the right-hand-side is
\begin{align*}
&\big((\b^{Y_{t-2}}(t-1) \setminus 
\b^{E_{g^{\psi}}}(t-1)) \setminus
\e^X(t)\big)
\cup
\n^{X}(t) \setminus \big(\e^X(t)\cup \b^{E_{g^{\psi}}}(t-1)\big)
\\
\subseteq\ &
\big(\b^{X}(t-1)\setminus \e^X(t)\big)
\cup
\big(\n^{X}(t) \setminus \e^X(t)\big)
= \b^X(t),
\end{align*}
establishing  (I2) at $t$.

\end{itemize}

\end{itemize}

\end{definition}

\section{Types of high-level states}
\label{sec:sets}

We have already defined failure states and initialising states in Definitions~\ref{def:failure-state} and \ref{def:init-state}. In Section~\ref{sec:backoff-bounding-valid} we will define our high-level state transition rule~$\hls$. But first we define the other types of high-level states -- \advancing, \filling, \refilling, and \stabilising.

Recall that $\phi=100$, $\chi = 1000$, and $\CUpsilon = 320$ from  Definition~\ref{def:constants}
and that $\Upsilon_{j,\geq W} = \{ \ell \in [j-1]: W_\ell \geq W\}$ from Definition~\ref{def:Upsilon}. 
Recall the 
definition of $L(j)$ from Definition~\ref{def:noise}, the
definition of covered and exposed from Definition~\ref{def:covered}, and the definition of $\Wtilde[j]$ from Definition~\ref{def:Wtilde}.
Define $\gammavect$ as follows.
\begin{definition}\label{def:real-gamma}
For all integers $x \ge 1$, let 
$\newsym{gamma-x}{component of $\gammavect$ - Definition~\ref{def:real-gamma}}{\gamma_x = 1/(4x^\phi)}$ and from Definition~\ref{def:mu}
Let
$\gammavect = (\gamma_1,\gamma_2,\ldots)$.
\end{definition}

Recall from Definition~\ref{def:mu} that, 
for a fixed $\lambda \in (0,1)$, and a fixed send sequence $\p$, 
$\mu_x^{\gammavect}$ is the quantity defined by
$\mu_x^{\gammavect} = \lambda W_x \prod_{a=1}^x(1-\gamma_a)$.
The point of $\mu_x^{\gammavect}$ is that,
in a $j$-jammed process~$Y$
with initial size distribution 
$\Po({\zvect})$ and start time~$\tau$,
each $b_x^Y(t)$ is a Poisson random variable 
with mean $f^{\zerovect,\zvect}_x(t-\tau)$
(see Lemma~\ref{lem:f-track-jammed}).
However, (see Lemma~\ref{lem:f-bound} with $\Gammavect = \zerovect$),
this mean is at most $\lambda W_x$.
So we consider bin~$x$ to be ``full'' if it contains $\lambda W_x$ balls. 
The quantity $\mu_x^{\gammavect}$ allows
a little degradation so it is a little less than $\lambda W_x$. However, the following observation, Observation~\ref{obs:mu-gamma-bounds}, shows that $\mu_x^{\gammavect}$ is still a good measure of how full bin~$i$ might be.

\begin{observation} 
\label{obs:mu-gamma-bounds}
Fix $\lambda \in (0,1)$, a send sequence $\p$, and a positive integer~$x$. Then $2\lambda W_x/3
\leq \mu_x^{\gammavect} \leq 3\lambda W_x/4$.
\end{observation}
\begin{proof}
We wish to show that $
2/3 \leq \prod_{a=1}^x(1-1/(4a^{\phi})) \leq 3/4$.
The lower bound comes from the crude bound
\[\sum_{a=1}^x 1/a^{\phi} \leq 1 + 2^{-(\phi-2)}\sum_{a=2}^x 1/a^2 \leq  1 + 2^{-(\phi-2)} \pi^2/6 \leq 4/3
\]
so $\prod_{a=1}^x (1-1/(4a^{\phi})) \geq 1- \frac{1}{4} \sum_{a=1}^x 1/a^\phi \geq 2/3$.
The upper bound comes from 
$\prod_{a=2}^x(1-1/(4 a^{\phi})) \leq 1$.
\end{proof}

We next define  $j$-advancing, $j$-filling, $j$-refilling, and $j$-stabilising states.
We start by defining a $j$-advancing state, which is the only one of these types that applies to a covered bin~$j$ (Definition~\ref{def:covered}).

\begin{definition}\label{def:advancing-state}
Fix $\lambda \in (0,1)$ and a send sequence $\p$ with $p_0=1$. Let $j$ be a positive integer.
We say a high-level state $\Psi$ is \emph{$j$-advancing} if 
bin~$j$ is covered and the following properties hold. 
\begin{itemize}
\item $\type^{\Psi} = \advancing$. 
\item $j^{\Psi} = j$.
\item $z_{j}^{\Psi} = 0$. 
\item For all $k \in [j-1]$, $z_k^{\Psi} = \mu_k^{\gammavect}$.
\item If bin~$j$ is many-covered then
$\calS^{\Psi}= \{\Upsilon_{j, \geq \Wtilde[j]}\}$.
If bin~$j$ is heavy-covered then 
$\calS^{\Psi} = \{\Upsilon_{j,\geq j^2}\}$.
\end{itemize} 
\end{definition}

The key to  Definition~\ref{def:advancing-state} is
that it applies when $j$ is covered.
As we noted in the introduction, and in Section~\ref{sec:covered}, this means that we expect bins $1,\ldots,j-1$ to continue providing $\Omega(\log j)$ expected noise (for the escape process) as bin~$j$ fills.
The bins in the set in $\calS^{\Psi}$ are the ones that we expect to provide this noise.
We will choose our high-level state transition rule so that
the volume process goes to this state after all bins in $[j-1]$ are already reasonably full in expectation, so this is why
for all $k\in [j-1]$, $z_k^{\Psi} = \mu_k^{\gammavect}$. 
However, the definition sets $z_j^{\Psi}=0$. This is because, when the volume process enters 
a $j$-advancing state, we don't assume anything
about how full bin~$j$ already is -- the 
high-level state is ``advancing'' to this~$j$ in order to fill it. 

As we noted, 
we expect at least $\L(j)$ noise while $j$ is being filled, but some bins may go empty. Therefore, in Definition~\ref{def:advancing-state}  we define the set  $\calS^{\Psi}$  
to be the single set $S$ given in the definition.
This choice of~$S$ fulfils two requirements --
we require both that $| S
| = \Omega(\log j)$ (so this set provides $\L(j)$ noise when full) and that, subject to this, the weights of bins in 
$S$ are as large as possible -- either $j^2$ in the heavy-covered case, or at least $\Wtilde[j]$ in the many-covered case -- 
so that $S$ stays full for as long as possible. 
The choice of $\Wtilde[j]$ 
is for this reason -- it is
is chosen in Definition~\ref{def:Wtilde} to maximise $W|\Upsilon_{j,\geq W}|$ while ensuring $|\Upsilon_{j,\geq W}|= \Omega(\log j)$.  
The reason that we want to maximise $W|\Upsilon_{j,\geq W}|$ is that, at smallest, all of the bins in $\Upsilon_{j,\geq W}$ have weight~$W$ so their population is Poisson with mean $W | \Upsilon_{j,\geq W}|$ -- 
we want this quantity to be large so that they
can maintain the desired noise for as long as possible. 

The remaining types of high-level states ($j$-filling, $j$-refilling, and $j$-stabilising) all apply to bins~$j$ that are exposed.
In fact, they will be weakly exposed since we will ensure that $j_0$ is chosen so that all strongly-exposed bins are in $[j_0-1]$. Also (T\ref{item:trans4}) ensures that $j$ does not decrease as high-level state transitions occur.
We emphasise this point with Remark~\ref{rem:j}.
\begin{remark}\label{rem:j}
High-level states are used in the proof of Theorem~\ref{thm:goaljkillerSE} which concerns send sequences~$\p$ with finitely-many strongly exposed bins. The proof starts by choosing a value~$j_0$ that is sufficiently large that all strongly-exposed bins lie in~$[j_0-1]$. The initialising high-level state~$\Psi_0$ will have $j^{\Psi_0} = j_0$. Thus, we will never encounter a high-level state~$\Psi$ with $j^{\Psi} < j_0$.  This means that bin $j^{\Psi}$ is always covered, or weakly exposed.
\end{remark}

\begin{definition}\label{def:filling-state}
Fix $\lambda \in (0,1)$ and a send sequence $\p$ with $p_0=1$. Let $j$ be a positive integer.  We say a high-level state $\Psi$ is \emph{$j$-filling} if bin~$j$ is exposed and the following properties hold.
\begin{itemize}
\item $\type^\Psi = \filling$.
\item $j^\Psi = j$.
\item For all $k \in [j-1]$, $z_k^\Psi = \mu_k^{\gammavect}$.

\item  $\calS^\Psi = \{\Upsilon_{j, \ge \Wtilde[j]} \setminus \Upsilon_{j, \ge j^2}\}
\cup \{\{k\}\colon k \in \Upsilon_{j, \ge j^\chi}\}$.

\end{itemize}
\end{definition}

A $j$-filling state, as defined in Definition~\ref{def:filling-state} is   analogous to an advancing state, but for an exposed bin. We will show that when the volume process is in high-level state~$\Psi$ for a $j$-filling state~$\Psi$, there is expected noise at least $\L(j)$ for filling bin~$j$.

Since bin $j$ is exposed, we don't expect bin $j$ to fill without losing this noise. Hence, the 
high-level state transition rule that we will use will make many   high-level-state transitions (from filling to refilling to stabilising and back to filling) as $j$ fills.
Each time the volume process enters a filling state, the bins in $[j-1]$ will be reasonably full in expectation, so for all $k\in [j-1]$, 
Definition~\ref{def:filling-state}
requires
$z_k^{\Psi} = \mu_k^{\gammavect}$, 
similarly to Definition~\ref{def:advancing-state}.  However, unlike the advancing state, it is not assumed that $z_j^{\Psi}=0$. 
Crucially, Definition~\ref{def:filling-state} does not place any restrictions at all on $z_j^\Psi$.
This is 
so that, as the volume process repeatedly enters a filling state while bin~$j$ is filling,
$z_j^{\Psi}$ can be used to track the progress
that is being made (on filling bin~$j$). 

The set $\calS^{\Psi}$ (containing sets that provide reliable noise) contains two types of set. The first type is the set $S$ which contains all bins in $\Upsilon_{j, \geq \Wtilde[j]}$ that have weight at most $j^2$.   These will typically provide expected noise $\L(j)$ (to the escape process) for $\exp(\Omega(j))$ time, ensuring that the volume processes adds $\exp(\Omega(j))$ balls to bin $j$ in expectation before it transitions out of a filling state. 

The second type is a singleton set $\{k\}$ where $k$ is a very high-weight bin with $W_k \geq j^{\chi}$. 
There are at least $\CUpsilon \LL(j)/\lambda$ such bins since bin $j$ is weakly exposed. 
Like all bins in $[j-1]$ these bins are expected to be reasonably full when the filling state is entered, and they are useful because they
will 
then provide at least $\LL(j)$ expected noise 
for a long period (polynomial in~$j$),
even if all other bins are empty and even if the volume process transitions to refilling and stabilising. We will use this to guarantee at least $\LL(j)$ noise while in refilling and stabilising states, where our goal will be for the volume process to quickly refill the bins in $1,\ldots,j-1$ (to re-establish noise $\L(j)$) and  re-enter a filling state. 

Importantly, none of the very high-weight bins in singleton sets appear in $S$; thus when the noise from $S$ drops below $\L(j)$ and the volume process transitions out of filling, zeroing bins in $S$ in the process, the very high-weight bins are not zeroed.  
This is the reason for the factor of~$2$ in the denominator in (Prop~\ref{cov-prop-2}) -- we exclude high-weight bins from~$S$, but we must ensure that $|S|$ still has size $\log j$.
We next give the definition of $j$-refilling states.

\begin{definition}\label{def:refilling-state}
Fix $\lambda \in (0,1)$ and a send sequence $\p$ with $p_0=1$. Let $j$ be a positive integer.  We say a high-level state $\Psi$ is \emph{$j$-refilling} if bin~$j$ is exposed and the following properties hold.
\begin{itemize}
\item $\type^\Psi = \refilling$.
\item $j^\Psi = j$.
\item For all positive integers $k$ with $W_k < j^2$, $z_k^\Psi = 0$.
\item $\calS^\Psi = \{\{k\}\colon 
k \in \Upsilon_{j,\geq j^\chi} \}$.
\end{itemize} 
\end{definition}

Our high-level state transition rule will
have the property that
the volume process goes to a $\refilling$ state after the low-weight bins (those with weight at most $j^2$) have become too empty, so 
in Definition~\ref{def:refilling-state},
every low-weight bin $x$ has $z_x^{\Psi}=0$. 
We keep progress on higher-weight bins, so 
Definition~\ref{def:refilling-state} doesn't restrict their $z$-values. Finally, the set $\calS$ that maintains noise while refilling is the set of singleton sets containing very-high-weight bins.
These are the only source of noise in refilling.  The purpose of refilling is
to refill bins $[j-1]$  and 
we will show that
it is very unlikely that these very-high-weight bins stop being full in the time that the volume process stays in refilling, so the probability of dropping below $\LL(j)$ expected noise while in refilling is sufficiently small. 
Once refilling finishes,
bins $[j-1]$ are all refilled (in expectation) and the process will enter a stabilising state, as a means towards re-entering a filling state. 
There is a small chance the transition from a refilling state into a stabilising state could ``fail''. If this happens the volume process will enter a new refilling state -- we will say more about this when we define the high-level state transition rule. 

One property that the transition rule will have is that 
the volume process only ever enters a stabilising state from a refilling state. The point of a stabilising state is to provide a guaranteed time window, after leaving a refilling state, in which the low-weight bins in the first $(\log j)^2$ bins continue to provide $\L(j)$ noise, which will be very important for analysing the relevant escape process. As far as the volume process is concerned, the goal of the stabilising state is just for these bins to remain full for long enough (for a period which is at most a polynomial in~$j$). We now give the definition of $j$-stabilising states.

\begin{definition}\label{def:stabilising-state}
Fix $\lambda \in (0,1)$ and a send sequence $\p$ with $p_0=1$. Let $j$ be a positive integer.  We say a high-level state $\Psi$ is \emph{$j$-stabilising} if bin~$j$ is exposed and the following properties hold.

\begin{itemize}
\item $\type^\Psi = \stabilising$.
\item $j^\Psi=j$.
\item For all $k\in [j-1]$, $z_k^\Psi = \mu_k^{\gammavect}$.
\item $\calS^\Psi = \{\{ k \in [(\log j)^2] \colon W_k < j^2\}, \{k \in [j-1] \setminus [(\log j)^2] \colon W_k < j^2\}\} \cup \{\{k\}\colon k \in \Upsilon_{j, \ge j^\chi}\}$.
\end{itemize} 
\end{definition}

For all $k\in [j-1]$, 
Definition~\ref{def:stabilising-state}
assigns $z_k^{\Psi} = \mu_k^{\gammavect}$ (since our high-level state transition rule will ensure that these bins start full in expectation from the refilling state).  Definition~\ref{def:stabilising-state} does not place any restrictions on $z_j^\Psi$, so that progress in filling bin~$j$ can be kept.

In the definition of $\calS^{\Psi}$, we include 
all low-weight bins in order to provide $\L(j)$ noise for the volume process and we also include all singleton sets containing very high-weight bins in order to provide $\LL(j)$ noise in the event of a transition back to refilling.
We split the low-weight bins into two sets -- a set $S_1$ containing those that are in the first $(\log j)^2$ bins, and a set $S_2$ containing the others.

In order to bound the ball counts in the escape process from above, we will first use the singleton sets in $\calS$ -- 
since there are at least $\LL(j)$ of these,
it will turn out that they
provide enough noise that the escape process's end-of-step arrival probability 
will be 
$\exp(- \Omega(\log \log j))$. Using this, bins in $S_1$ will quickly empty in the escape process.
Since $S_1 \in \calS^{\Psi}$, this   will drop the escape   end-of-step arrival probability 
of the escape process
to $\exp(-\Omega(|S_1|))$ which is at most a small inverse-polynomial function of~$j$. Using this, bins in $[j-1]$ will then quickly empty in the escape process, and in particular bins in $S_2$ will quickly empty; this will ensure (with summable failure probability) that the 
end-of-step arrival probability of the escape process remains low  while the volume process is  in the stabilising high-level state.

In Section~\ref{sec:backoff-bounding-valid} we will give the specific high-level state transition rule that we will use.

\section{The backoff-bounding transition rule}
\label{sec:backoff-bounding-valid}

We next give a high-level state transition rule (Definition~\ref{def:transition-rule}) giving transitions between the specific high-level states that we have described (failure states -- Definition~\ref{def:failure-state}, 
initialising states -- Definition~\ref{def:init-state}, 
advancing states -- Definition~\ref{def:advancing-state}, 
filling states -- Definition~\ref{def:filling-state}, 
refilling states -- Definition~\ref{def:refilling-state}, 
and stabilising states -- Definition~\ref{def:stabilising-state}).

We will consider a positive integer~$j_0$
and a $j_0$-initialising state~$\Psi_0$. Then the transition rule operates on the set 
$\Psi^*(\Psi_0)$ 
(Definition~\ref{def:pre-backoff-bounding-state-space}) containing $\Psi_0$ and  all
non-initialising high-level states~$\Psi$ with $j^{\Psi} \geq j_0$.
We refer to the high-level state transition rule~$\hls$ (Definition~\ref{def:backoff-bounding-rule})
as the  \emph{$\Psi_0$-backoff-bounding rule} -- we include $\Psi_0$ in the name because it is the first high-level state, and we call it the ``backoff-bounding rule'' because it is the high-level transition rule that we will used to bound (from below) the population of the backoff process.
We will later prove that the backoff-bounding rule~$\hls$ 
obeys the constraints on a high-level transition rule from Definition~\ref{def:transition-rule} (Observation~\ref{obs:BB-bounding-is-transrule})
and that it
is  $\Psi^*$-valid (Definition~\ref{def:valid}) 
which we prove (Lemma~\ref{lem:backoff-bounding-valid})
via Lemma~\ref{lem:valid-transition}.

\begin{definition}\label{def:pre-backoff-bounding-state-space}
Fix $\lambda \in (0,1)$ and a send sequence $\p$ with $p_0=1$. Let $j_0 \ge 2$ be a positive integer. 
Let $\Psi_0$ be a $j_0$-initialising state.
We define $\newsym{Psi-*}{set of high-level states
for backoff bounding rule}{\Psi^*(\Psi_0)}$   to be the set 
containing $\Psi_0$ and all $j$-failure,  $j$-advancing, $j$-filling, $j$-refilling, and $j$-stabilising states with $j \ge j_0$.
\end{definition}

The backoff-bounding rule $\hls$ is given in Definition~\ref{def:backoff-bounding-rule}.
Recall from Definition~\ref{def:transition-rule}
that $\calR(\Psi^*(\Psi_0))$ is the set of triples
$(\Psi,\bbar,t)$ such that
$\Psi \in \Psi^*(\Psi_0)$, ${\bbar}\in 
(\integers_{\geq 0})^{j^\Psi}
 $,   and $t> \tau^{\Psi}$.  
 If $\type^\Psi=\initialising$  then $t=\tau^\Psi+1$.  
Also, recall the definition 
of~$\Upsilon_{j,\geq W}$ from Definition~\ref{def:Upsilon} and the definition
of~$F$ from~Definition~\ref{def:F}.
Recall from Definition~\ref{def:constants} that $\phi=100$ and $\chi=1000$. Recall from Definition~\ref{def:real-gamma}
that
${\gamma_k = 1/(4k^\phi)}$,
and from Definition~\ref{def:mu} that
$\mu_x^{\gammavect} = \lambda W_x \prod_{a=1}^x(1-\gamma_a)$.

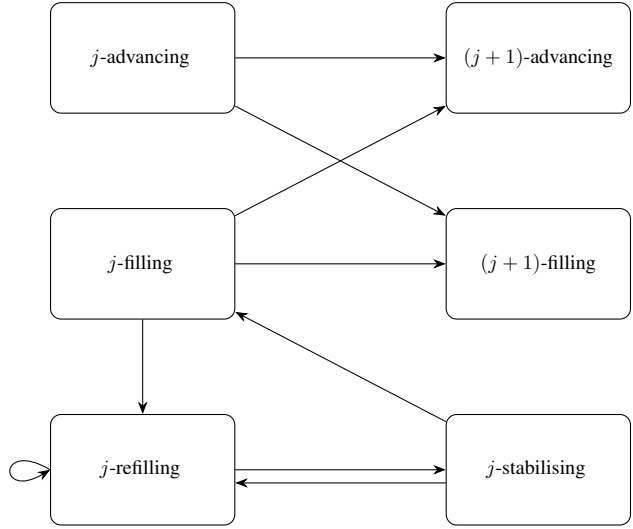
\begin{figure}[ht]   

\centering
\begin{tikzpicture}[
  scale=0.7,         
  transform shape,    
  node distance = 1.8cm and 4cm,
  >=Stealth,
  box/.style={
    draw,
    rounded corners,
    minimum width=3.5cm,
    minimum height=2.1cm,  
    align=center
  }
]

\node[box] (jadv)
  {$j$-advancing  };

\node[box, below=of jadv] (jfill)
  {$j$-filling  };

\node[box, below=of jfill] (jrefill)
  {$j$-refilling  };

\node[box, right=of jadv] (jp1adv)
  {$(j+1)$-advancing};

\node[box, below=of jp1adv] (jp1fill)
  {$(j+1)$-filling};

\node[box, below=of jp1fill] (jstab)
  {$j$-stabilising  };

\draw[->] (jadv)   -- (jp1adv);      
\draw[->] (jadv)   -- (jp1fill);   

\draw[->] (jfill)  -- (jp1adv);    
\draw[->] (jfill)  -- (jp1fill);    
\draw[->] (jfill)  -- (jrefill);    

\draw[->] (jrefill) -- (jstab);    
\draw[->] (jstab)  -- (jfill);    

\draw[->] ([yshift=-0.25cm] jstab.west) -- ([yshift=-0.25cm] jrefill.east); 
\draw[->]
  (jrefill.west)
    .. controls +(-1,0.6) and +(-1,-0.6)
    .. (jrefill.west);

\end{tikzpicture}

\caption{High-level state transitions under the backoff-bounding rule from Definition~\ref{def:backoff-bounding-rule}.  
\label{fig:backoff-bounding-rule}}

\end{figure}

Figure~\ref{fig:backoff-bounding-rule} summarises the transitions between high-level states made by the backoff-bounding rule in Definition~\ref{def:backoff-bounding-rule}.
The basic transitions follow the sketch given in the introduction to the paper.  
When the volume process starts in the initialising state~$\Psi_0$,  the bins of  $[j_0]$  are an independent Poisson tuple  by Definition~\ref{def:volume-process}. 
By Lemma~\ref{lem:dist-volume}, conditioned on each sequence of high-level states $\Psi(\tau_0),\ldots,\Psi(t)$, the bins of  $[j^{\Psi(t)}]$  are an independent Poisson tuple. 
From there, the process eventually gets to a $(j^{\Psi(t)}+1)$-advancing or $(j^{\Psi(t)}+1)$-filling state via (R\ref{item:R4}), however, in the case that 
$j^{\Psi(t)}$ is exposed this may involving going through successive refilling and stabilising states with $j=j^{\Psi(t)}$ until some filling state finally fills bin $j$. The process moves to a failure state when certain unlikely events happen, the most important of which is if a very high-weight bin is no longer full, either in expectation or in reality, 
(R\ref{item:R1})
or if a transition from $j$ to $j+1$ fails in (R\ref{item:R4}).

\begin{definition}\label{def:backoff-bounding-rule}
Fix $\lambda \in (0,1)$ and a send sequence $\p$ such that $p_0=1$.
Let $j_0\geq 2$ be a positive integer.
Let $\Psi_0$ be a $j_0$-initialising state.
The \emph{$\Psi_0$-backoff-bounding rule} is the transition rule $\hls$  
on~$\Psi^*(\Psi_0)$   with the following properties
for all 
$(\Psi= (g,\tau,j,\zvect,\calS,\type),\bbar,t)\in \calR(\Psi^*(\Psi_0))$ where $\Psi' = (g',\tau',j',z',\calS',\type') = \hls(\Psi,\bbar,t)$.

\begin{enumerate}[(R1)]

\item \label{item:R1} If any of the following hold then $\Psi'= (g+1,t,j, \zerovect, \emptyset,\Failure)$.  
\begin{enumerate}[(i)]
\item $\type=\Failure$, or
\item 
for some $k \in \Upsilon_{j, \ge j^\chi}$,   $F_k^{\Psi}(t) < \lambda W_k/2$, or
\item for some $k\in \Upsilon_{j,\geq j^2}$, $F^{\Psi}_k(t) \ge \lambda W_k/2$ and $b_k < \lambda W_k/4$.
\end{enumerate}

\item \label{item:R2} If (R\ref{item:R1}) does not apply and, for some $S\in \calS$, 
$\sum_{k \in S} p_k b_k < \lambda |S|/40$,
then 
\begin{enumerate}[(i)]
\item  
If $\type \in \{\refilling, \advancing\}$  then 
$\Psi'= (g+1,t,j, \zerovect, \emptyset,\Failure)$.  
\item Otherwise, 
$\Psi'= (g+1,t,j,\zvect', \calS', \refilling)$ is the $j$-refilling state where, for all $k \in \Upsilon_{j+1,\geq j^2}$,   $z_k' =\min\{F_k^\Psi(t),\mu_k^{\gammavect}\}$.  
\end{enumerate}

\item \label{item:R3} 
If (R\ref{item:R1}) and (R\ref{item:R2}) don't apply   and $\type=\initialising$   then
\begin{enumerate}[(i)]
\item If bin $j$ is covered then 
$\Psi'= (g+1,t,j,\zvect', \calS', \advancing)$ is 
the $j$-advancing state unless this would violate (T\ref{item:trans3}) in which case 
$\Psi'= (g+1,t,j, \zerovect, \emptyset,\Failure)$. 
\item If bin $j$ is exposed  
 then $\Psi'= (g+1,t,j,\zvect', \calS', \filling)$ is the $j$-filling state with $z_j' = 0$
 unless this would violate (T\ref{item:trans3}) in which case 
 $\Psi'= (g+1,t,j, \zerovect, \emptyset,\Failure)$. 
\end{enumerate}

\item \label{item:R4} 
If (R\ref{item:R1}) and (R\ref{item:R2}) don't apply   and $\type\in \{\advancing,\filling\}$ then
\begin{enumerate}[(i)]
\item If  $t\geq \tau + j^{24}$ and  $F_{j}^{\Psi}(t) \ge \mu^{\gammavect}_j$  
\begin{enumerate}[(a)]
\item If bin~$j+1$ is covered then 
$\Psi'= (g+1,t,j+1,z', \calS', \advancing)$ is 
the $(j+1)$-advancing state unless this would violate (T\ref{item:trans3}) in which case 
$\Psi'= (g+1,t,j, \zerovect, \emptyset,\Failure)$. 
\item If bin~$j+1$ is exposed then 
$\Psi'= (g+1,t,j+1,z', \calS', \filling)$ is the $(j+1)$-filling state with   $z'_{j+1}=0$ unless this would violate (T\ref{item:trans3}) in which case 
$\Psi'= (g+1,t,j, \zerovect, \emptyset,\Failure)$. 
\end{enumerate} 
\item Otherwise, $\Psi'=\Psi$.
\end{enumerate}

\item  \label{item:R5}
If (R\ref{item:R1}) and (R\ref{item:R2}) don't apply  and $\type=\refilling$  
then
\begin{enumerate}[(i)]
\item If  {$t \geq \tau + j^{\phi+46}$} and  
for all $k \in [j-1]$, $F_k^{\Psi}(t) \ge \mu_k^{\gammavect}$ then
$\Psi'= (g+1,t,j,\zvect', \calS', \stabilising)$ is the  $j$-stabilising state with $z_j' = F_j^{\Psi}(t)$ unless this would violate (T\ref{item:trans3}) in which case 
$\Psi'= (g+1,t,j,\zvect', \calS', \refilling)$ is the $j$-refilling state where, for all $k \in \Upsilon_{j+1,\geq j^2}$,    
$z_k' = F_k^{\Psi}(t)$.   
\item Otherwise, $\Psi' = \Psi$.
\end{enumerate}

\item \label{item:R6} 
If (R\ref{item:R1}) and (R\ref{item:R2}) don't apply   and $\type=\stabilising$ then
\begin{enumerate}[(i)]
\item If $  t \ge \tau + j^{\phi + 72}$  
then  $\Psi'= (g+1,t,j,\zvect', \calS', \filling)$ is the $j$-filling state with   $z_j' = F_j^{\Psi}(t)$
unless this would violate (T\ref{item:trans3})  
in which case 
$\Psi'= (g+1,t,j,\zvect', \calS', \refilling)$ is the $j$-refilling state where, for all $k \in \Upsilon_{j+1,\geq j^2}$,    
$z_k' = F_k^{\Psi}(t)$.  
\item Otherwise, $\Psi' = \Psi$.
\end{enumerate}
\end{enumerate}
\end{definition}

The possible transitions in Definition~\ref{def:backoff-bounding-rule} are illustrated in 
Figure~\ref{fig:backoff-bounding-rule}.
The figure doesn't include self-loops from~$\Psi$ to~$\Psi$, so the loop from $j$-refilling to $j$-refilling indicates the possible transition in (R\ref{item:R5})(i) from~$\Psi$ to a  $j$-refilling state~$\Psi'$ which differs from~$\Psi$.

In the remainder of this section, we 
first 
give a name to the closure of the backoff-bounding rule~$\hls$ 
(using the definition of closure from Definition~\ref{def:hls-closed}). We then
show that 
the backoff-bounding rule~$\hls$ is a valid high-level state transition rule.

\begin{definition}\label{def:backoff-bounding-state-space}
Fix $\lambda \in (0,1)$ and a send sequence $\p$ with $p_0=1$. Let 
$j_0 \ge 2$ be a positive integer. 
Let $\Psi_0$ be a $j_0$-initialising state.
We define $\newsym{Psi-BB}{closure of the backoff bounding rule -- Def~\ref{def:pre-backoff-bounding-state-space}}{\OurHLS}$ to be the $\Psi_0$-closure of the  $\Psi_0$-backoff bounding rule from Definition~\ref{def:backoff-bounding-rule}.  
\end{definition}

\begin{observation}\label{obs:BB-bounding-is-transrule} 
Fix $\lambda \in (0,1)$ and a send sequence $\p$ such that $p_0=1$.
Let $j_0\geq 2$ be a positive integer.
Let $\Psi_0$ be a $j_0$-initialising state.
The {$\Psi_0$-backoff-bounding rule} from  Definition~\ref{def:backoff-bounding-rule} 
satisfies (T\ref{item:trans1}) -- (T\ref{item:trans5}) of Definition~\ref{def:transition-rule}.  
\end{observation}
\begin{proof}
(T\ref{item:trans1}) and (T\ref{item:trans2}) are by inspection.
We establish (T\ref{item:trans3}) 
as follows.
\begin{itemize}
\item If $\type' = \Failure$ this follows since $\calS' = \emptyset$ from Definition~\ref{def:failure-state}.
\item If $\type' = \refilling$ this follows since 
$\calS' = \{ \{k\}: k \in \Upsilon_{j,\geq j^{\chi}}\}$ from Definition~\ref{def:refilling-state}. But by (R\ref{item:R1}), for each $k\in \Upsilon_{j,\geq j^{\chi}}$,   $p_k b_k \geq \lambda/4 \geq \lambda/40$. 
\item If $\Psi' = \Psi$ and $\type'\notin \{\Failure, \refilling \}$ this holds by (R\ref{item:R2}).
\item The other cases all have an explicit check on (T\ref{item:trans3}).
\end{itemize} 
Items
(T\ref{item:trans4}) 
and (T\ref{item:trans5})  are by inspection.
\end{proof}

\begin{lemma}\label{lem:BB-bounding-V1}
Fix $\lambda \in (0,1)$ and a send sequence $\p$ such that $p_0=1$. Let $j_0$ be a positive integer.
Let $\Psi_0$ be a $j_0$-initialising state.
The {$\Psi_0$-backoff-bounding rule} $\hls$ from  Definition~\ref{def:backoff-bounding-rule} 
satisfies
item~(V\ref{item:valid-zlower}) of Lemma~\ref{lem:valid-transition} on the state space closure~\OurHLS.
Namely, 
for all 
$(\Psi= (g,\tau,j,\zvect,\calS,\type),\bbar,t)\in \calR(\OurHLS)$ where $\Psi' = (g',\tau',j',z',\calS',\type') = \hls(\Psi,\bbar,t)$,
if ${\Psi'} \ne \Psi$, then for all $\ell \in  [j]$, 
$z'_\ell \le F_\ell^{\Psi}(t) 
$.   
\end{lemma}
\begin{proof}

From Definition~\ref{def:F}, $F_\ell^{\Psi}(t) = f_\ell^{\zerovect,\zvect}(t-\tau)$
and this quantity is non-negative by  Definition~\ref{def:f}.
We then consider cases.
\begin{itemize}
\item  If $\type' = \Failure$ then, from Definition~\ref{def:backoff-bounding-rule},
$\Psi'= (g+1,t,j, \zerovect, \emptyset,\Failure)$.
So for $\ell\in [j]$, $z'_\ell=0$ which is 
at most  the  non-negative value
$F_\ell^{\Psi}(t)$.
\item Transitions (R\ref{item:R4})(ii), (R\ref{item:R5})(ii), and (R\ref{item:R6})(ii) are not covered by the lemma since $\Psi' = \Psi$.

\item 
If $\type' = \refilling$ and $\Psi' \neq \Psi$ 
(so the transition is to a refilling state, but not via (R\ref{item:R5})(ii)) then $j' = j$ and, 
for all $k\in \Upsilon_{j+1,\geq j^2}$, $z'_k \leq F_k^{\Psi}(t)$.
The definition of refilling states -- Definition~\ref{def:refilling-state} -- ensures that for $k\in [j] \setminus \Upsilon_{j+1,\geq j^2}$, $z'_k=0$. 
So
for all $k$, 
$z'_k \leq F_k^{\Psi}(t)$.
\item Transition (R\ref{item:R3}) and $\type' \neq \Failure$:
Here $j'=j$, $\type = \initialising$, and
$\type' \in \{ \advancing, \filling\}$.
By the definitions of advancing and filling states 
(Definitions~\ref{def:advancing-state} and~\ref{def:filling-state}),
 for all $k\in [j-1]$, $z_k'  
 = \mu_k^{\gammavect} = \lambda W_k \prod_{a=1}^k(1-\gamma_a)
 $.
 $z'_j = 0$.  We have already shown that $F_k^{\Psi}(t) \geq 0$. Thus we wish to show that for $k\in [j-1]$ $\mu_k^{\gammavect} \leq F_k^{\Psi}(t)$.  By Definition~\ref{def:transition-rule}, $t= \tau+1$.
By the definition of initialising state (Definition~\ref{def:init-state}) for all $k\in [j]$, 
$z_k = 3\lambda W_k/4$.
For $k \in \{2,\ldots, j-1\}$,
Definitions~\ref{def:F} and~\ref{def:f} ensure that
$F_k^{\Psi}(t) = (1-p_k) 3 \lambda W_k/4 + p_{k-1} 3 \lambda W_{k-1}/4 = 3 \lambda W_k/4
$. By Observation~\ref{obs:mu-gamma-bounds}, this is at least $\mu_k^{\gammavect}$, as required. The case~$k=1$ is similar.
The definitions ensure that
$F_1^{\Psi}(t) = (1-p_1) 3 \lambda W_1/4 + \lambda \geq 3 \lambda W_1/4$.

\item   
Transition (R\ref{item:R4})(i) and $\type' \neq \Failure$: Here
$j'=j+1$, and both~$\type$ and $\type'$ are in $\{\advancing, \filling\}$.
By the definitions of
advancing and filling states (Definitions~\ref{def:advancing-state} and~\ref{def:filling-state}),
 for all $k\in [j-1]$, $z_k = z_k'  
 = \mu_k^{\gammavect}  
 $.  Also $z'_j = \mu_j^{\gammavect}$ and $z'_{j+1}=0$.
By Corollary~\ref{cor:incf} (with $\Gamma = \zerovect$ and $k=j-1$),
for all $\ell\in [j-1]$,
$F_\ell^{\Psi}(t) = f_\ell^{\zerovect, \zvect}(t-\tau) \geq f_\ell^{\zerovect, \zvect}(0) = z_\ell = \mu_\ell^{\gammavect} = z'_\ell$.
Also, by the guard in (R\ref{item:R4})(i),
$  
F_j^{\Psi}(t) \geq \mu_j^{\gammavect} = z'_j$.
Finally, $F_{j+1}^\Psi(t)$ is non-negative, so it is at least $z'_{j+1}=0$.

 \item  Transition (R\ref{item:R5})(i) and $\type' \neq  \refilling$: Here  
 $j'=j$ and $\type' = \stabilising$. By the definition of a stabilising state (Definition~\ref{def:stabilising-state}) for all $k\in [j-1]$, $z'_k = \mu_k^{\gammavect}$ which is at most $F_k^{\Psi}(t)$, by the guard. Also, $z'_j = F_j^{\Psi}(t)$.

 \item  Transition (R\ref{item:R6})(i) 
 and $\type' \neq  \refilling$: Here $j'=j$,
 $\type = \stabilising$, and $\type' = \filling$. By the definitions of filling state (Definition~\ref{def:filling-state}) and stabilising state (Definition~\ref{def:stabilising-state}), 
 for all $k\in [j-1]$, $z'_k = \mu_k^{\gammavect}= z_k$. 
Similar to the argument for (R\ref{item:R4})(i),
by Corollary~\ref{cor:incf} (with $\Gamma = \zerovect$ and $k=j-1$),
for all $\ell\in [j-1]$,
$F_\ell^{\Psi}(t) = f_\ell^{\zerovect, \zvect}(t-\tau) \geq f_\ell^{\zerovect, z}(0) = z_\ell = \mu_\ell^{\gammavect} = z'_\ell$.
Finally, $z'_j = F_j^{\Psi}(t)$.

\end{itemize}
\end{proof}

\begin{lemma}\label{lem:BB-bounding-V2}

Fix $\lambda \in (0,1)$ and a send sequence $\p$ such that $p_0=1$.
Let $j_0\geq 2$ be a positive integer. 
Let $\Psi_0$ be a $j_0$-initialising state.
The {$\Psi_0$-backoff-bounding rule} $\hls$ from  Definition~\ref{def:backoff-bounding-rule} 
satisfies
item~(V\ref{item:valid-up}) of Lemma~\ref{lem:valid-transition} on the state space closure~\OurHLS.
Namely, 
for all 
$(\Psi= (g,\tau,j,\zvect,\calS,\type),\bbar,t)\in \calR(\OurHLS)$ where $\Psi  = \hls(\Psi,\bbar,t)$,
for all 
$  \bbarup\in (\integers_{\geq 0})^j $
such that  $\bbarup \geq \bbar$, $\hls(\Psi, \bbarup,t) = \Psi$.
\end{lemma}
\begin{proof}
The relevant transitions are (R\ref{item:R4})(ii), (R\ref{item:R5})(ii), and (R\ref{item:R6})(ii).
We first note if (R\ref{item:R1}) and (R\ref{item:R2}) don't apply then the same is true when $\bbar$ is replaced by~$\bbarup$.
We then note that the same is true for the guards in  (R\ref{item:R4})(i), (R\ref{item:R5})(i) and (R\ref{item:R6})(i).
\end{proof}

\begin{lemma}\label{lem:BB-bounding-V3}
Fix $\lambda \in (0,1)$ and a send sequence $\p$ such that $p_0=1$.
Let $j_0\geq 2$ be a positive integer. 
Let $\Psi_0$ be a $j_0$-initialising state.
The {$\Psi_0$-backoff-bounding rule} from  Definition~\ref{def:backoff-bounding-rule} 
satisfies
item~(V\ref{item:valid-up-non-zero}) of Lemma~\ref{lem:valid-transition} on the state space closure~\OurHLS.
Namely,  
for all 
$(\Psi= (g,\tau,j,\zvect,\calS,\type),\bbar,t)\in \calR(\OurHLS)$ where $\Psi' = (g',\tau',j',z',\calS',\type') = \hls(\Psi,\bbar,t)$,
if ${\Psi'} \ne \Psi$, then
for all 
$  \bbarup\in (\integers_{\geq 0})^j $
such that  $\bbarup \geq \bbar$ 
and for all $k\in [j]$ with $z'_k=0$,
$b^+_k = b_k$,
$\hls(\Psi, \bbarup,t) = \Psi'$.
\end{lemma}
\begin{proof}
We consider the transitions that may occur. 
\begin{itemize}
\item If $\type' = \Failure$ this is vacuously true since, by 
by the definition of failure state (Definition~\ref{def:failure-state}), for all $k\in [j']$, $z'_k = 0$, so $b^+_k = b_k$. 

\item Transition by (R\ref{item:R2})(ii): Note that $j'=j$.  
Since (R\ref{item:R1})  doesn't apply, for all  
$k\in \Upsilon_{j,\geq j^\chi}$,  $b_k \geq \lambda W_k/4$.
Clearly it is also true that,
for all  
$k\in \Upsilon_{j,\geq j^\chi}$,  $b^+_k \geq \lambda W_k/4$.
Since (R\ref{item:R2})  applies,
for some $S\in \calS$, 
$\sum_{k\in S} p_k b_k < \lambda |S|/40$. 
Clearly, 
$\type\in \{\filling, \stabilising\}$
($\type$ cannot be in $\{\Failure,\initialising\}$ since those have $\calS = \emptyset$ and it cannot be in $\{\refilling, \advancing\}$ since those lead to (R\ref{item:R2})(i)). Also $\type' = \refilling$.
From the definitions of filling and stabilising states (Definitions~\ref{def:filling-state} and~\ref{def:stabilising-state}), 
$\calS$ contains some sets~$S$ such that,  for every $k\in S$, $W_k < j^2$. 
It causes no trouble if one of these sets $S$ triggers (R\ref{item:R2}) 
with 
$\sum_{k\in S} p_k b_k < \lambda |S|/40$.
The reason that it causes no trouble is that
such a set $S$ has the property that, for every   $k\in S$, $z'_k=0$ (by the definition of refilling state, Definition~\ref{def:refilling-state}),
so $\sum_{k\in S} p_k b^+_k   = 
\sum_{k\in S} p_k b_k
$, as required.  There is one more possibility to consider.
$\calS$ contains   the sets $\{k\}$ such that $k\in \Upsilon_{j,\geq j^{\chi}}$.
So we have to consider the possibility that (R\ref{item:R2})(ii) is triggered because there is a  $k\in 
\Upsilon_{j,\geq j^{\chi}}$ such that
$p_k b_k < \lambda /40$. However, we have already seen that this is ruled out by (R\ref{item:R1}).

\item For the remaining items, note that if (R\ref{item:R1}) and (R\ref{item:R2}) don't apply then the same is true when $\bbar$ is replaced with $\bbarup$. Also, if a proposed transition obeys (T\ref{item:trans3}) then the same is true when $\bbar$ is replaced with $\bbarup$. These two observations are all that is needed for (R\ref{item:R3}) and (R\ref{item:R4})(i).
These two observations are also all that is needed
for transitions in (R\ref{item:R5})(i), and (R\ref{item:R6})(i), except transitions to refilling that are caused by violation of (T\ref{item:trans3}).

 \item Transitions (R\ref{item:R4})(ii), (R\ref{item:R5})(ii) and (R\ref{item:R6})(ii) are not covered by the lemma since $\Psi' = \Psi$.
 
\item  Transitions to a refilling state $\Psi'$
in   (R\ref{item:R5})(i) and (R\ref{item:R6})(i)
caused by the fact that an attempted transition to a stabilising or filling state~$\Psi''$ would violate (T\ref{item:trans3}).
Stabilising and filling states have some sets $S\in \calS^{\Psi''}$ 
such that for all $k\in S$, $W_k< j^2$. However, these states $S$ are not a problem,  since, 
from the definition of refilling state (Definition~\ref{def:refilling-state}),
$z'_k=0$ for every $k$ with $W_k < j^2$.
(This is like the argument for (R\ref{item:R2})(ii).)
$\calS^{\Psi''}$ also contains
the sets $\{k\}$ such that $k\in \Upsilon_{j,\geq j^{\chi}}$. If the violation of (T\ref{item:trans3}) were triggered by one of these sets $\{k\}$ we would have
 $p_k b_k < \lambda / 40$. This is impossible since (R\ref{item:R1}) would have applied to the transition from~$\Psi$.
 We conclude that the violation of (T\ref{item:trans3}) gets triggered by a set $S\in \calS^{\Psi''}$ whose bins~$k$
 have $b^+_k = b_k$ so the violation of 
(T\ref{item:trans3}) still occurs when~$\bbar$ is replaced by~$\bbarup$.

 \end{itemize}
\end{proof}

\begin{lemma}
\label{lem:backoff-bounding-valid}

Fix $\lambda \in (0,1)$ and a send sequence $\p$ such that $p_0=1$. Let $j_0$  be a positive integer.
Let $\Psi_0$ be a $j_0$-initialising state.
The {$\Psi_0$-backoff-bounding rule} $\hls$ from  Definition~\ref{def:backoff-bounding-rule} 
is $\OurHLS$-valid (Definition~\ref{def:valid}).
\end{lemma}
  
\begin{proof}
Observation~\ref{obs:BB-bounding-is-transrule}
shows that the $\Psi_0$-backoff-bounding rule satisfies (T\ref{item:trans1})--(T\ref{item:trans5}) of Definition~\ref{def:transition-rule}.
To show that it is $\OurHLS$-valid (from Definition~\ref{def:valid}) we  use Lemma~\ref{lem:valid-transition}, which shows that we need only establish (V\ref{item:valid-zlower}), (V\ref{item:valid-up}) and (V\ref{item:valid-up-non-zero}). These follow from Lemmas \ref{lem:BB-bounding-V1}, \ref{lem:BB-bounding-V2}, and \ref{lem:BB-bounding-V3}.
\end{proof}

\section{Analysis of the volume process}
\label{sec:backoff-bounding-analysis}

The main purpose of this section is to prove Lemma~\ref{lem:VolumeAnalysis}
which shows that, with probability at least~$99/100$, the volume process never enters a failure state, and for all $j\geq j_0$, it quickly gets to a $j$-filling or $j$-advancing state. Lemma~\ref{lem:VolumeAnalysis} 
immediately implies
Corollary~\ref{cor:VEB-volume} (which has the same content, but is expressed in terms of the VEB-coupling).  We will use Corollary~\ref{cor:VEB-volume} later, when we use the VEB coupling to bound the population of a backoff process.

When we prove our final theorem (specifically Theorem~\ref{thm:goaljkillerSE}) using
the VEB coupling we will choose $j_0$ to be sufficiently large and the volume process will start with a $j_0$-initialising state. The backoff bounding rule (Definition~\ref{def:backoff-bounding-rule}) never decreases $j$ in high-level state transitions, so in our analysis we will always be able to assume that $j$ is sufficiently large. We don't spell out the exact requirements each time since we will be able to freely choose $j_0$, but the way to think about it is that $j_0$ is sufficiently large depending on~$\lambda$  
and~$\p$. Specific things that we assume include $j_0\geq 5$ and that (since these theorems only apply when  $\p$ has  only finitely-many strongly exposed bins),   all strongly-exposed bins lie in $[j_0-1]$. We also assume that  
$j_0-1 \geq C_\Upsilon \L(j_0)/\lambda$, 
so that we can rely on $|\Upsilon_{j_0,\geq \Wtilde[j_0]}| \geq C_\Upsilon \L(j_0)/\lambda$ (see Definition~\ref{def:Wtilde}). And we assume that
$\CUpsilon \L(j_0)/\lambda$ and 
$\CUpsilon \LL(j_0)/\lambda$ are at least~$2$, to avoid trivialities. To avoid lengthening the paper, we won't spell all of these assumptions each each time, we will just state in each relevant lemma that  ``$j$ is sufficiently large''.

Towards proving Lemma~\ref{lem:VolumeAnalysis},
we start with Section~\ref{sec:volume-M}
which gives some general lemmas about how the expected bin populations change as the volume process progresses. It introduces the key concept of ``missing'' balls.
In Section~\ref{sec:volume-transient} we will show that the volume process is not likely to get stuck for too long in any single high-level state. In Section~\ref{sec:volume-bad-trans} we will show that certain ``bad'' high-level state transitions (including transitions to failure states) are unlikely. In Section~\ref{sec:volume-theloop}
we will analyse the progress of the process through the refilling-stabilising loop. Finally, in Section~\ref{sec:volume-progress} we will show that the process is likely to progress~$j$ sufficiently quickly, and we will prove Lemma~\ref{lem:VolumeAnalysis}.

\subsection{Tracking expectations and ``missing'' balls}
\label{sec:volume-M}
 
The analysis of the volume process   focusses on expectations because the process follows a $j$-jammed  process and the balls are independent. 
Indeed, we have already shown
in Lemma~\ref{lem:dist-volume}
that the distribution of
$\bbar_{[j^{\Psi(t)}]}^{V}(t)$, conditioned on the sequence of high-level states up through time~$t$, is a tuple of
independent Poisson values whose means are given by the function~$F^{\Psi(t)}$. To be specific,
the distribution is 
$\Po(F_{[j^{\Psi(t)}]}^{\Psi(t)}(t))$.

We have already seen that the function $F^{\Psi(t)}$ follows a simple recurrence. Specifically,
$F_k^{\Psi(t)}(t) = 
f_k^{\zerovect,\zvect^{\Psi(t)}}(t-\tau^{\Psi(t)})$ (Definition~\ref{def:F}), and
the function 
$f_k^{\zerovect,\zvect^{\Psi(t)}}$  follows a simple recurrence relation  
(Definition~\ref{def:f}).

The analysis of expectations in the volume process would be straightforward if $F_k^{\Psi(t)}(t)$ were monotonically increasing between changes of the high-level state. However, we've already seen that this is not the case.
Item (V\ref{item:valid-zlower}) in the conditions for the validity of a high-level-state transition rule (Lemma~\ref{lem:valid-transition})
ensures that, if the high-level state changes to $\Psi(t+1) \neq \Psi(t)$,
then the new expectation $z_k^{\Psi(t+1)}$, is at most 
$F_k^{\Psi(t)}(t)$, but it might be smaller (and it is smaller if bin~$k$ is over-full before the transition).

Recall from Definition~\ref{def:real-gamma} that $\gamma_k = 1/(4k^{\phi})$ where (from Definition~\ref{def:constants}) $\phi=100$. 
Recall from Definition~\ref{def:mu} that, 
for a fixed $\lambda \in (0,1)$, and a fixed send sequence $\p$, 
$\mu_x^{\gammavect}$ is the quantity defined by
$\mu_x^{\gammavect} = \lambda W_x \prod_{a=1}^x(1-\gamma_a)$.
We think of bin~$k$ as being full if it has at least $\mu_k^{\gammavect}$ balls -- from Observation~\ref{obs:mu-gamma-bounds}, 
$2\lambda W_k/3
\leq \mu_k^{\gammavect} \leq 3\lambda W_k/4$.
We use Definition~\ref{def:M} to
track how far the expectation 
$F_k^{\Psi(t)}(t)$ differs from 
$\mu_k^{\gammavect}$.
We think of $\calM_S^{\Psi(t)}(t)$ as the number of ``missing'' balls expected to be in a set~$S$ of bins at time~$t$ relative to~$\mu^{\gammavect}_{S}$.

\begin{definition}
\label{def:M}
Fix a send sequence~$\p$ and a $\lambda \in (0,1)$.
Let $\Psi=(g,\tau,j,\zvect,\calS,\type)$ be a high-level state
and let $S$ be a subset of $[j]$.
Let $t\geq \tau$ be an integer.
Then
$\newsym{M-S}{missing balls - Def~\ref{def:M}}{\calM_S^\Psi(t) }:= \sum_{k\in S} \max\{0, \mu^{\gammavect}_k - F_k^\Psi(t)\}$
\end{definition}

Lemma~\ref{lem:gen-F-minus-mu-grows} establishes monotonicity for $\calM_{[j']}^\Psi(t)$ (as a function of~$t$) for all $j' \le j^\Psi$.

\begin{lemma}\label{lem:gen-F-minus-mu-grows}
Fix $\lambda \in (0,1)$. Let $\p$ be a send sequence with $p_0=1$. Let   $\Psi = (g,\tau,j,\zvect,\calS,\type)$ be a high-level state. Fix $j'\in [j]$.
Then for all $t' > \tau$,
$\calM_{[j']}^{\Psi}(t') \leq
\calM_{[j']}^{\Psi}(t'-1)$.
\end{lemma}
\begin{proof}
By Definition~\ref{def:F}, the statement is equivalent to showing that, 
for all $t>0$,
\[
\sum_{x \in [j']} \max\{0, \mu_x^{\gammavect} - f^{\zerovect,\zvect}_x(t)\} \le \sum_{x \in [j']} \max\{0, \mu_x^{\gammavect} - f^{\zerovect,\zvect}_x(t-1)\}.
\]

Fix $t>0$.
For all $x \in [j]$, let $z^+_x = f_x^{\zerovect,\zvect}(t-1)$.
Let $\zvect^+$ be the $j$-tuple
$(z^+_1,\ldots,z^+_j)$.
Observe that $f^{\zerovect,\zvect}_x(t) = f_x^{\zerovect,\zvect^+}(1)$. 
For all $x\in [j]$,  let $a_x = \min\{z^+_x, \mu_x^{\gammavect}\}$. Let $\avect$ be the $j$-tuple $(a_1,\ldots,a_{j})$.
We will apply Lemma~\ref{lem:f-increasing} with $\Gammavect = \zerovect$,
$\Lambdavect = \gammavect$, the $\zvect$ of
Lemma~\ref{lem:f-increasing}
as ${\zvect}^+$ and the $\avect$ of Lemma~\ref{lem:f-increasing} as $\avect$. Note that our $\avect \leq {\zvect}^+$ as required.  The lemma shows, in our notation, that for all $k\in [j]$ and $t''\geq 0$,
$f_k^{\zerovect,{\zvect}^+}(t'')\geq f_k^{\gammavect,\avect}(t'')$.
It follows that, for all $x\in [j]$, $f_x^{\zerovect,\zvect}(t) = f_x^{\zerovect,\zvect^+}(1) \geq f_x^{\gammavect,\avect}(1)$, and hence
\[
\sum_{x \in [j']} \max\{0, \mu_x^{\gammavect} - f^{\zerovect,\zvect}_x(t)\} \le \sum_{x \in [j']} \max\{0, \mu_x^{\gammavect} - f^{\gammavect,\avect}_x(1)\}.
\]
We next apply Lemma~\ref{lem:f-bound} with $\Gammavect = \gammavect$, and the
$\zvect$ of Lemma~\ref{lem:f-bound} as our $\avect$.
This applies since $0 \leq a_x \leq \mu_{x}^{\gammavect}$ for all $x\in [j]$. The lemma, taking $t=1$, shows that  
for all $x \in [j]$,
$0 \leq f_x^{\gammavect,\avect}(1) \leq \mu_x^{\gammavect}$.
Thus,
\[
\sum_{x \in [j']} \max\{0, \mu_x^{\gammavect} - f^{\gammavect,\avect}_x(1)\} = \sum_{x \in [j']} (\mu_x^{\gammavect} - f^{\gammavect,\avect}_x(1)) = \sum_{x \in [j']}\mu_x^{\gammavect}- \sum_{x \in [j']}f_x^{\gammavect,\avect}(1).
\]
We next apply Lemma~\ref{lem:gen-f-monotone} 
with $\Gammavect = \gammavect$, $t=1$, and the
$\zvect$ of Lemma~\ref{lem:gen-f-monotone} as 
$\avect$. The lemma shows that
\[
\sum_{x\in [j']} f_x^{\gammavect,\avect}(1) \geq \sum_{x\in [j']} f_x^{\gammavect,\avect}(0).
\]
Thus
\begin{align*}
\sum_{x \in [j']} \max\{0, \mu_x^{\gammavect} - f^{\zerovect,\zvect}_x(t)\} 
&\le \sum_{x \in [j']}\mu_x^{\gammavect}- \sum_{x \in [j']}f_x^{\gammavect,\avect}(0) = \sum_{x \in [j']} (\mu_x^{\gammavect} - a_x) = \sum_{x \in [j']} (\mu_x^{\gammavect} - \min\{z_x^+, \mu_x^{\gammavect}\})\\
&= \sum_{x \in [j']} \max\{0,\mu_x^{\gammavect} - z^+_x\}
= \sum_{x \in [j']} \max\{0,\mu_x^{\gammavect} - f_x^{\zerovect,\zvect}(t-1)\},
\end{align*}
as required.
\end{proof}

Lemma~\ref{lem:gen-F-minus-mu-grows} shows that while the high-level state stays at~$\Psi$, 
$\calM_S^{\Psi}(t)$ shrinks, or at least it doesn't grow. 
We can also give a corollary (Corollary~\ref{cor:f-fill-all-bins}) of the material in Section~\ref{sec:j-jammed} 
which  implicitly gives an upper bound on how long it takes for $\calM_{[j']}^{\Psi}(t)$ to shrink to~$0$,
making all bins in $[j']$ full in expectation: this happens in $O(\calM_{[j']}^\Psi(t)\cdot \poly(j))$ time.
This upper bound will be useful 
for analysing the volume process when it enters a refilling state. As  can be seen from the 
definition of refilling state, (Definition~\ref{def:refilling-state}), and from the
backoff bounding rule  (Definition~\ref{def:backoff-bounding-rule}), when the volume process enters a $j$-refilling state~$\Psi$, it only zeroes bins of weight less than $j^2$.
Thus,  $\calM_{[j]}^\Psi(t) = O(j^3)$. In this case,  Corollary~\ref{cor:f-fill-all-bins} gives an effective bound on how long it takes to repair the damage caused by entering the refilling state, and to refill the affected bins (in expectation).

\begin{corollary}\label{cor:f-fill-all-bins}
Fix a send sequence~$\p$ and a $\lambda \in (0,1)$.
Let $\Psi=(g,\tau,j,\zvect,\calS,\type)$ be a high-level state
and let $j'$ be an integer in $[j]$. 
Fix $R\in \reals_{\geq 0}$  and suppose that
$\calM_{[j']}^{\Psi}(\tau)  \le R$.
Then for all $t' 
\ge \tau + (4/\lambda) j' R/\min\{\gamma_k\colon k \in [j']\}$ and all $k \in [j']$, $F_k^{\Psi}(t') \ge \mu_k^{\gammavect}$ so
$\calM_{[j']}^{\Psi}(t')=0$.
\end{corollary}

\begin{proof}
Take the $j$ of
Lemma~\ref{lem:f-fill-all-bins} to be~$j'$.
Apply Lemma~\ref{lem:f-fill-all-bins}
with $\Gammavect= \zerovect$ and $\Lambdavect = \gammavect_{[j']}  $. 
Note that $\mu_k^{\Lambdavect} \geq \lambda W_k/2$ by Observation~\ref{obs:mu-gamma-bounds}.
\end{proof}

We will shortly give another lemma, Lemma~\ref{lem:bound-z-shrink-one-loop},  which bounds how $\calM_S^{\Psi(t)}(t)$ changes between certain high-level state transitions, but first we give another useful corollary of Section~\ref{sec:j-jammed}, along the lines of Corollary~\ref{cor:f-fill-all-bins}, and we define the relevant high-level state transitions.

\begin{corollary}\label{cor:f-fill-bin}
Fix a send sequence~$\p$ with $p_0=1$ and a $\lambda \in (0,1)$.
Let $\Psi=(g,\tau,j,\zvect,\calS,\type)$ be a high-level state.  
For any $k\in [j]$ and any $t' \geq \tau$,  suppose that, for all $x \in [k]$,  $F_x^{\Psi}(t') \ge \mu_x^{\gammavect}$. Then for all $x\in[k]$ and $t'' \geq t'$, $F_x^{\Psi}(t'') \ge \mu_x^{\gammavect}$. Also, if $k<j$ then for all $t'' \ge t'$,
 $
F_{k+1}^{\Psi}(t'') \ge \min\{\mu_{k+1}^{\gammavect},\ F_{k+1}^{\Psi}(t') + (t''-t') \gamma_{k+1}p_k\mu_k^{\gammavect}\}$.
In particular, for all $t'' \ge t' + W_{k+1}/\gamma_{k+1}$,  $F_{k+1}^{\Psi}(t'') \ge \mu_{k+1}^{\gammavect}$.
\end{corollary}
\begin{proof}
We will apply Lemma~\ref{lem:f-fill-bin}
with $\Gammavect = \zerovect$ and $\Lambdavect = \gammavect$.
It will be useful to take the $j$ of Lemma~\ref{lem:f-fill-bin} to be our $j+1$ -- for this take (our) $z_{j+1} = 0$.
Take $t_0 = t'-\tau$. 

To show that  
for all $x\in[k]$ and $t'' \geq t'$, $F_x^{\Psi}(t'') \ge \mu_x^{\gammavect}$
apply item (i) with
$t = t''-t'$. This shows
that $f^{\zerovect,\gammavect}_x(t_0+t''-t') \geq \mu_x^{
\gammavect}$, which is what we require since $t_0 + t'' - t' = t'' - \tau  $.

To show that, if $k<j$, then for all  $t'' \ge t'$,
 $
F_{k+1}^{\Psi}(t'') \ge \min\{\mu_{k+1}^{\gammavect},\ F_{k+1}^{\Psi}(t') + (t''-t') \gamma_{k+1}p_k\mu_k^{\gammavect}\}$, take item (ii) with $t=t''-t'$.

It remains to show that when $t'' \ge t' + W_{k+1}/\gamma_{k+1}$, $F_{k+1}^\Psi(t'') \ge \mu_{k+1}^{\gammavect}$. Given this bound on $t''$, we have
\[
(t''-t')\gamma_{k+1}p_k\mu_k^{\gammavect} \ge W_{k+1}p_k\mu_k^{\gammavect} = W_{k+1}p_k \cdot \lambda W_k\prod_{a=1}^k (1-\gamma_a) \ge \lambda W_{k+1} \prod_{a=1}^{k+1} (1-\gamma_a) = \mu_{k+1}^{\gammavect}.
\]
Thus by the second part of the result,
\[
F_{k+1}^{\Psi}(t'') \ge \min\{\mu_{k+1}^{\gammavect},\ F_{k+1}^{\Psi}(t') + (t''-t') \gamma_{k+1}p_k\mu_k^{\gammavect}\} \ge \min\{\mu_{k+1}^{\gammavect},\ F_{k+1}^{\Psi}(t') + \mu_{k+1}^{\gammavect}\} = \mu_{k+1}^{\gammavect},
\]
as required.
\end{proof}

We next single out certain high-level state transitions that we will use to bound $\calM_S^{\Psi(t)}$. From a $j$-filling state (see Figure~\ref{fig:backoff-bounding-rule}) the progress of the volume process that will be helpful 
for bounding the population of the backoff process  will be a move to a $(j+1)$-advancing or $(j+1)$-filling state (signalling that bin~$j$ has filled). As we already explained, 
it is unlikely that bin $j$ will fill without losing noise, so the backoff-bounding rule will make many high-level-state transitions from $j$-filling to $j$-refilling to $j$-stabilising and back to $j$-filling (following the cycle at the bottom of the figure) while bin~$j$ is filling up. 
Sometimes the transition from a $j$-refilling state into a $j$-stabilising state fails -- if this happens the process instead takes the self-loop to a new $j$-refilling state -- we refer to this self-loop as a ``back transition''. The transition from $j$-stabilising to $j$-filling could also fail, leading to a back-transition back to $j$-refilling (as shown by the left-pointing horizontal arrow in the figure). Definition~\ref{def:back-transition} gives the formal definition of these back transitions.
(These are the back-transitions of the refilling-stabilising loop.)

\begin{definition}
\label{def:back-transition}
Fix $\lambda \in (0,1)$, a send sequence $\p$ with $p_0=1$, and a positive integer~$j_0$. Let $\tau_0$ be a non-negative integer. Let $\Psi_0$ be the  initialising high-level state with $\tau^{\Psi_0} = \tau_0$ and $j^{\Psi_0} = j_0$. 
Let $\hls$ be the $\Psi_0$-backoff bounding rule from Definition~\ref{def:backoff-bounding-rule}.
Let $V$ be a volume process from $\Psi_0$ with transition rule $\hls$, send sequence~$\p$ and birth rate~$\lambda$. 
For any $t>\tau_0$ 
the high-level state transition 
given by 
$\Psi(t) 
= \hls(\Psi(t-1), \bbar^{Y_{t-1}}_{[j^{\Psi(t-1)}]}(t), t)$ is a \emph{back transition} if 
$\type^{\Psi(t)} = \refilling$ and
one of the following hold:
\begin{itemize}
\item $\type^{\Psi(t-1)} = \refilling$, and the transition is via (R\ref{item:R5})(i); or
\item $\type^{\Psi(t-1)} = \stabilising$, and the transition is via (R\ref{item:R2})(ii) or (R\ref{item:R6})(i).
\end{itemize}
\end{definition}

 Lemma~\ref{lem:bound-z-shrink-one-loop} says that $\calM_{[j-1]}^{\Psi(t)}(t)$ increases by $O(j^3)$ between time~$t$ (where $\Psi(t)$  is refilling or stabilising) and the next back transition before re-entering a filling state.

\begin{lemma}
\label{lem:bound-z-shrink-one-loop}
Fix $\lambda \in (0,1)$ and a send sequence $\p$ with $p_0=1$.
Let $j_0$ be a sufficiently large positive integer. Let $\tau_0$ be a non-negative integer. Let $\Psi_0$ be the  initialising high-level state with $\tau^{\Psi_0} = \tau_0$ and $j^{\Psi_0} = j_0$. 
Let $\hls$ be the $\Psi_0$-backoff bounding rule from Definition~\ref{def:backoff-bounding-rule}.
Let $V$ be a volume process from $\Psi_0$ with transition rule $\hls$, send sequence~$\p$ and birth rate~$\lambda$. For any
$t \ge \tau_0$, 
let   $T(t) = \{\tau_0,\ldots,t\}$.
Consider any integers $t' > t \geq \tau_0$ and any sequence 
$\psivect_{T(t')} = (\psi_{\tau_0},\dots,\psi_{t'})$ 
of states in~$\OurHLS$ 
such that the following hold.
\begin{itemize}
\item  
$\Pr( \Psivect_{T(t')} = \psivect_{T(t')}
)>0
$.
\item 
For all $t'''\in \{t,\ldots,t'-1\}$,
$\type^{\psi_{t'''}} \in \{\refilling,\stabilising\}$.
\item  For all $t'' \in\{t,\ldots,t'-2\}$,
the transition 
$\psi_{t''} \rightarrow \psi_{t''+1}$ is not a back transition.
\item  The
transition 
$\psi_{t'-1} \rightarrow \psi_{t'}$ is a back transition.
\end{itemize}
Let $j= j^{\psi_t}$.
Whenever $\Psivect_{T(t')} = \psivect_{T(t')} $,
$
\calM_{[j-1]}^{\psi_{t'}}(t') 
\le 
\calM_{[j-1]}^{\psi_t}(\tau^{\psi_t})
+ 
\lambda j^3$.

\end{lemma}
\begin{proof}
Note that for all $t''\in \{t,\ldots,t'\}$, 
$j^{\psi_{t''}} = j$.
The relevant transition is $\psi_{t'} = \hls(\psi_{t'-1}, \bbar^{Y_{t'-1}}_{[j]}(t'),t')$. 
Recall from Definition~\ref{def:M} that 
for any high-level state~$\psi$, any 
$S\subseteq [j^{\psi}]$, and any $t\geq \tau^{\psi}$,
$\calM_S^\psi(t) = \sum_{k\in S} \max\{0, \mu^{\gammavect}_k - F_k^\Psi(t)\}$.
Let $S_1 = \Upsilon_{j,\geq j^2}$
and let $S_2 = [j-1] \setminus S_1$. From its definition,
$\calM_{[j-1]}^{\psi_{t'}}(t') 
 = 
 \calM_{S_1}^{\psi_{t'}}(t') + 
 \calM_{S_2}^{\psi_{t'}}(t')$.
By the definition of refilling state (Definition~\ref{def:refilling-state}), for all $k \in  S_2$, $z_k^{\psi_{t'}} = 0$. 
Since the transition to~$\psi_{t'}$ is a back-transition, $\psi_{t'} \neq \psi_{t'-1}$ so $\tau^{\psi_{t'}} = t'$ 
and by Observation~\ref{obs:F},
$F^{\psi_{t'}}_{[j]}(t') = \zvect^{\psi_{t'}}$.  For
$k\in S_2$,   $\mu_k^{\gammavect} \leq \lambda W_k \leq \lambda j^2$, so
$\calM_{S_2}^{\psi_{t'}}(t') 
 \leq \lambda j^3$.
We will finish by showing
$\calM_{S_1}^{\psi_{t'}}(t')
\leq \calM_{[j-1]}^{\psi_t}(\tau^{\psi_t})$.
We split into  cases.

\medskip\noindent\textbf{Case 1: } $\type^{\psi_t} = \stabilising$. 
By the  definition
of stabilising state 
(Definition~\ref{def:stabilising-state}),
for all $k\in [j-1]$, $z_k^{\psi_t} = \mu_k^{\gammavect}$ so, from its definition, 
$ 
\calM_{[j-1]}^{\psi_t}(\tau^{\psi_t})=  0$. 
 By the backoff bounding rule (Definition~\ref{def:backoff-bounding-rule}) the transition to $\psi_{t'}$ is the first transition that does not follow (R\ref{item:R6})(ii). Thus, for all $t''\in \{t,\ldots,t'-1\}$, $\psi_{t''} = \psi_t$.
By Lemma~\ref{lem:gen-F-minus-mu-grows}    
since $t'> t \geq \tau^{\psi_t}$,
it follows that 
$
\calM_{[j-1]}^{\psi_t}(t') \leq 
\calM_{[j-1]}^{\psi_t}(\tau^{\psi_t})
=   0$.  Since 
$
\calM_{[j-1]}^{\psi_t}(t')  (\tau^{\psi_t})
=   0$,
each $k\in [j-1]$ has
$F_k^{\psi_t}(t') \geq \mu_k^{\gammavect}$.

\medskip\noindent\textbf{Case 1A: } The back transition to $\psi_{t'}$ is via (R\ref{item:R6})(i): We've just seen that $F^{\psi_{t'}}_{[j]}(t') = \zvect^{\psi_{t'}}$. 
From the backoff bounding rule   for all $k\in  S_1$, $z^{\psi_{t'}}_k = 
F_k^{\psi_{t'-1}}(t')= F_k^{\psi_t}(t')$, and we have already seen that this is at least $ \mu_k^{\gammavect}$.
Thus, from its definition, 
$\calM_{S_1}^{\psi_{t'}}(t') = 0$.

\medskip\noindent\textbf{Case 1B: } The back transition to $\psi_{t'}$ is via (R\ref{item:R2})(ii): In this case, from the backoff bounding rule   for all $k\in  S_1$, $z^{\psi_{t'}}_k = \min\{F_k^{\psi_t}(t'),\mu_k^{\gammavect}\}$.  
Since $
F_k^{\psi_t}(t')\geq \mu_k^{\gammavect}$, this is at least 
$ \mu_k^{\gammavect}$.
Thus, from its definition, 
$\calM_{S_1}^{\psi_{t'}}(t') = 0$.

\medskip\noindent\textbf{Case 2: } $\type^{\psi_t} = \refilling$. 
If there is a $t''\in \{t,\ldots,t'\}$ such that $\type^{\psi_{t''}} = \stabilising$ then we are finished by the calculation in Case~1. Suppose that this does not happen. Then by the backoff bounding rule   the transition to $\psi_{t'}$ is the first transition that does not follow (R\ref{item:R5})(ii).
Thus, for all $t''\in \{t,\ldots,t'-1\}$, $\psi_{t''} = \psi_t$. As in Case~1,
by Lemma~\ref{lem:gen-F-minus-mu-grows}   
since $t'> t \geq \tau^{\psi_t}$,
it follows that 
$
\calM_{[j-1]}^{\psi_t}(t') \leq 
\calM_{[j-1]}^{\psi_t}(\tau^{\psi_t})$.
 
Recall that $F^{\psi_{t'}}_{[j]}(t') = \zvect^{\psi_{t'}}$. 
From the backoff bounding rule (transition (R\ref{item:R5}(i))),   for all $k\in  S_1$, $z^{\psi_{t'}}_k = 
F_k^{\psi_{t'-1}}(t') = F_k^{\psi_t}(t')$.  Thus, from its definition, 
$
\calM_{S_1}^{\psi_{t'}}(t') 
\leq \calM_{[j-1]}^{\psi_{t}}(t')  \leq 
\calM_{[j-1]}^{\psi_t}(\tau^{\psi_t})$.
\end{proof}

\subsection{Quickly changing high-level state} \label{sec:volume-transient}

In this section we show that the volume process is not likely to get stuck for too long in any single 
advancing, refilling, or stabilising state. The reason for this is that the guard on the relevant transition out of this state is likely to become true quickly.
Filling states will be handled as part of the refilling-stabilising loop analysis in Section~\ref{sec:volume-theloop} -- we leave filling states until then because we expect to visit filling states multiple times in the process of filling bin~$j$, going through the loop each time, and losing a little progress towards
establishing   the guard on (R\ref{item:R4})(i).

\begin{lemma}\label{lem:advancing-transient}
Fix $\lambda \in (0,1)$ and a send sequence $\p$ with $p_0=1$.  
Let $j_0$ be a sufficiently large positive integer. Let $\tau_0$ be a  non-negative integer. Let $\Psi_0$ be the   initialising high-level state with  $\tau^{\Psi_0} = \tau_0$ and $j^{\Psi_0} = j_0$. 
Let $\hls$ be the $\Psi_0$-backoff bounding rule from Definition~\ref{def:backoff-bounding-rule}.  
Let $V$ be a volume process from $\Psi_0$ with transition rule $\hls$ send sequence~$\p$ and birth rate~$\lambda$. 
Let $t> \tau_0$ be a positive integer such that $\Psi(t)$ is $j$-advancing and $\tau^{\Psi(t)} = t$.
Then there is an integer $t'$
in the range
$t < t' \le t+ 12 j^{\phi+1}W_j$ such that $\Psi(t') \neq \Psi(t)$.
\end{lemma}

\begin{proof}
The goal is to show that we don't expect (R\ref{item:R4})(ii) to occur $12 j^{\phi+1}W_j$ times in a row starting from the transition from~$\Psi(t)$.
By the definition of advancing state (Definition \ref{def:advancing-state}), for all $k \in [j-1]$, $z_k^{\Psi(t)} = \mu_k^{\gammavect}$, and $z_j^{\Psi(t)} = 0$.  Therefore
$\calM_{[j]}^{\Psi(t)}(t) =  
\mu_j^{\gammavect}$. 
We apply Corollary~\ref{cor:f-fill-all-bins} with $j' = j$ and $R = \mu_j^{\gammavect}$
to show that for all $t' \geq 
t + (4/\lambda) j \mu_j^{\gammavect}/\min\{\gamma_k\colon k \in [j]\}$ and all $k \in [j]$, $F_k^{\Psi(t)}(t') \ge \mu_k^{\gammavect}$. 
Note that $\min\{\gamma_k: k\in [j]\} = \gamma_j = 1/(4 j^\phi)$. By Observation~\ref{obs:mu-gamma-bounds}, $W_j \geq 4 \mu_j^{\gammavect} / (3 \lambda)$.
Thus if $t'\geq t + 12 j^{\phi+1} W_j$ 
then 
\[ t' \geq 
t + 3  j W_j/
\min\{\gamma_k: k\in [j]\}
\geq 
t + (4/\lambda) j \mu_j^{\gammavect}/\min\{\gamma_k: k\in [j]\},
\]
so  $F_j^{\Psi(t)}(t') \geq \mu_j^{\gammavect}$. 
 
Transition (R\ref{item:R4})(i) triggers from the first $t'\geq t + j^{24}$ 
such that  $F_j^{\Psi(t)}(t') \geq \mu_j^{\gammavect}$, unless some other transition other than (R\ref{item:R4}(ii)) triggers first. Either way $\Psi(t') \neq \Psi(t)$. 
The lemma follows since
$j^{24} \leq 12 W_j j^{\phi+1}$.
\end{proof}

\begin{lemma}\label{lem:refilling-transient} 
Fix $\lambda \in (0,1)$ and a send sequence $\p$ with $p_0=1$.  
Let $j_0$ be a sufficiently large positive integer. Let $\tau_0$ be a  non-negative integer. Let $\Psi_0$ be the   initialising high-level state with $\tau^{\Psi_0} = \tau_0$ and $j^{\Psi_0} = j_0$. 
Let $\hls$ be the $\Psi_0$-backoff bounding rule from Definition~\ref{def:backoff-bounding-rule}.  
Let $V$ be a volume process from $\Psi_0$ with transition rule $\hls$ send sequence~$\p$ and birth rate~$\lambda$.  
Let $t> \tau_0$ be a positive
integer such that $\Psi(t)$ is $j$-refilling,
$\tau^{\Psi(t)} = t$, and
$\calM_{[j-1]}^{\Psi(t)}(t)
\leq 2j^{44}$. 
Then 
there is an integer $t'\in \{t+1,\ldots,t+j^{\phi+46}\}$  
such that $\Psi(t') \neq \Psi(t)$.   
\end{lemma}
\begin{proof} 
By the backoff bounding rule (Definition~\ref{def:backoff-bounding-rule}) we just need to rule out applying (R\ref{item:R5})(ii) $j^{\phi+46}$ times in a row starting from the transition from $\Psi(t)$.  
So we just need to show that,  for all
$k \in  [j-1]$, 
$F_k^{\Psi(t)}(t+{j^{\phi+46}}) \geq \mu_k^{\gammavect}$.

Towards this end, let $R=2j^{44}$ and note from the hypothesis of the lemma that
$\calM_{[j-1]}^{\Psi(t)}(t)
\leq R$. 
We apply Corollary~\ref{cor:f-fill-all-bins} with $j' = j-1$  
to show that for all $t' \geq 
t + (4/\lambda) (j-1) (2 j^{44})/\min\{\gamma_k\colon k \in [j-1]\}$ and all $k \in [j-1]$, $F_k^{\Psi(t)}(t') \ge \mu_k^{\gammavect}$. As in the proof of Lemma~\ref{lem:advancing-transient},
 $\min\{\gamma_k: k\in [j-1]\} = \gamma_{j-1} = 1/(4 (j-1)^\phi)$ where $\phi=100$.
Thus 
\[ t + j^{\phi+46} \geq 
t + (4/\lambda) (j-1) (2 j^{44})/\min\{\gamma_k\colon k \in [j-1]\},\]
since $j\geq j_0$ is sufficiently large. 
\end{proof}

\begin{lemma}\label{lem:stabilising-transient} 
Fix $\lambda \in (0,1)$ and a send sequence $\p$ with $p_0=1$.  
Let $j_0$ be a sufficiently large positive integer. Let $\tau_0$ be a  non-negative integer. Let $\Psi_0$ be the   initialising high-level state with $\tau^{\Psi_0} = \tau_0$ and $j^{\Psi_0} = j_0$. 
Let $\hls$ be the $\Psi_0$-backoff bounding rule from Definition~\ref{def:backoff-bounding-rule}.  
Let $V$ be a volume process from $\Psi_0$ with transition rule $\hls$ send sequence~$\p$ and birth rate~$\lambda$. 
Let $t> \tau_0$ be a positive
integer such that $\Psi(t)$ is $j$-stabilising and
$\tau^{\Psi(t)} = t$.
Then 
there is a $t'\in \{t+1,\ldots,t+j^{\phi+72}\}$ such that $\Psi(t') \neq \Psi(t)$.  
\end{lemma}
\begin{proof}
By the backoff bounding rule (Definition~\ref{def:backoff-bounding-rule}) we just need to rule out applying (R\ref{item:R6})(ii) $j^{\phi+72}$ times in a row starting from the transition from $\Psi(t)$.  This is ruled out by (R\ref{item:R6})(i).
\end{proof}
\subsection{Upper-bounding probability of ``bad'' high-level state transitions}
\label{sec:volume-bad-trans}

In this section we show that certain ``bad'' high-level state transitions (including transitions to failure states) are unlikely.  
In order to avoid repeated calculation, we start by giving the following lemmas that will be used repeatedly.

\begin{lemma}\label{lem:vol-F-lb}
Fix a real number $\lambda \in (0,1)$ and a send sequence $\p$ with $p_0=1$. Let $j_0$ be a sufficiently large integer.
Let  
$\Psi_0 
$ 
be the initialising state
with $j^{\Psi_0} = j_0$ and
$\tau^{\Psi_0} = \tau_0$.
Let $\hls$ be  
the $\Psi_0$-backoff bounding rule from Definition~\ref{def:backoff-bounding-rule}.
Let $V$ be a volume process from~$\Psi_0$ with transition rule~$\hls$, send sequence~$\p$ and birth rate~$\lambda$. 
For any
$t \ge \tau_0$, 
let   $T(t) = \{\tau_0,\ldots,t\}$.
Then   for
every $t >  \tau_0$
 and every tuple
$ 
\psivect_{T(t-1)} = (\psi_{\tau_0},\dots,\psi_{t-1})$  of high-level states such that $\type^{\psi_{t-1}} \notin\{\Failure, \refilling\}$, 
if
$\Psivect_{T(t-1)} = \psivect_{T(t-1)}$, then
for every $k\in [j^{\psi_{t-1}}-1]$,
and every $t' \geq \tau^{\psi_{t-1}}$,
$F_k^{\psi_{t-1}}(t') \geq \mu_k^{\gammavect} \geq 2 \lambda W_k/3$. In particular, $F_k^{\psi_{t-1}}(t) \geq \mu_k^{\gammavect} \geq 2 \lambda W_k/3$.
\end{lemma}

\begin{proof}
Fix $t > \tau_0$.  
Let $\psi_{t-1} = (g,\tau,j,\zvect,\calS,\type)$.
We will apply Corollary~\ref{cor:f-fill-all-bins} to $\psi_{t-1}$
with $j' = j-1$ and $R=0$.
We need that 
$\calM_{[j -1]}^{\psi_{t-1}}(\tau) = 0$, which by the definition of~$\calM$ (Definition~\ref{def:M}) is the same as requiring that, 
for $k\in [j-1]$, 
$z_k \geq \mu_k^{\gammavect}$.
We do this in cases.
\begin{itemize}
\item If   $\psi_{t-1}$ is initialising:
From the definition of initialising state
(Definition~\ref{def:init-state}) 
for $k\in [j]$, $z_k = 3 \lambda W_k/4 $. By
Observation~\ref{obs:mu-gamma-bounds},
$\mu_k^{\gammavect} \leq 3 \lambda W_k/4 = z_k$.
\item If $\psi_{t-1}$ is advancing,  filling, or stabilising: From the definitions of those states (Definition~\ref{def:advancing-state}, Definition~\ref{def:filling-state}, Definition~\ref{def:stabilising-state}),
for all $k\in [j -1]$, $z_k = \mu_k^{\gammavect}$. 
\end{itemize}

Corollary~\ref{cor:f-fill-all-bins} shows that for all $t' \geq \tau$ and all $k\in [j -1]$, 
$F_k^{\psi_{t-1}}(t') \geq \mu_k^{\gammavect}$. We use $t'=t$.
The final lower bound 
$\mu_k^{\gammavect}\geq 2 \lambda W_k/3$ comes from Observation~\ref{obs:mu-gamma-bounds}.
   \end{proof}

We now proceed to show that certain high-level state transitions are unlikely.
We will do this by considering the transitions of the backoff bounding rule in order, starting from (R\ref{item:R1}).
There is no need to bound the probability of (R\ref{item:R1})(i) since this transition is from a state that is already failing. Observation~\ref{obs:no-R1-from-adv-fill-stab} rules out (R\ref{item:R1})(ii), except from refilling states which will be considered later, as part of the analysis of the refilling-stabilising loop.

\begin{observation}\label{obs:no-R1-from-adv-fill-stab}
Fix a real number $\lambda \in (0,1)$ and a send sequence $\p$ with $p_0=1$. Let $j_0$ be a sufficiently large integer.
Let  
$\Psi_0 
$ 
be the initialising state
with $j^{\Psi_0} = j_0$ and
$\tau^{\Psi_0} = \tau_0$.
Let $\hls$ be  
the $\Psi_0$-backoff bounding rule from Definition~\ref{def:backoff-bounding-rule}.
Let $V$ be a volume process from~$\Psi_0$ with transition rule~$\hls$, send sequence~$\p$ and birth rate~$\lambda$. 
For any
$t \ge \tau_0$, 
let   $T(t) = \{\tau_0,\ldots,t\}$.
Then   for
every $t >  \tau_0$
 and every tuple
$ 
\psivect_{T(t-1)} = (\psi_{\tau_0},\dots,\psi_{t-1})$  of high-level states such that $\type^{\psi_{t-1}} \neq \refilling$, 
if
$\Psivect_{T(t-1)} = \psivect_{T(t-1)}$, then
the transition to $\Psi(t)$ is not via (R\ref{item:R1})(ii).
\end{observation}
\begin{proof}
Let $\psi_{t-1} = (g,\tau,j,\zvect,\calS,\type)$. 
The transition is defined by 
${\Psi(t) }
= \hls(\psi_{t-1}, \bbar^{Y_{t-1}}_{[j]}(t), t)$   where $Y_{t-1}$ is the process defined in the definition of the volume process (Definition~\ref{def:volume-process}).
It is not possible 
for (R\ref{item:R1})(ii) to trigger if
$\type=\Failure$, since that would trigger (R\ref{item:R1})(i) first.
Otherwise, Lemma~\ref{lem:vol-F-lb} shows
that for every $k\in [j-1]$, $F_k^{\psi_{t-1}}(t) \geq   2 \lambda W_k/3$ so (R\ref{item:R1})(ii) doesn't trigger. 
\end{proof}

Lemma~\ref{lem:volume-easyChernoffs} upper-bounds the probability of a transition via (R\ref{item:R1})(iii).

\begin{lemma}
\label{lem:volume-easyChernoffs}
Fix a real number $\lambda \in (0,1)$ and a send sequence $\p$ with $p_0=1$. Let $j_0$ be a sufficiently large integer.
Let  
$\Psi_0 
$ 
be the initialising state
with 
$j^{\Psi_0} = j_0$ and
$\tau^{\Psi_0} = \tau_0$.
Let $\hls$ be  
the $\Psi_0$-backoff bounding rule from Definition~\ref{def:backoff-bounding-rule}.
Let $V$ be a volume process from~$\Psi_0$ with transition rule~$\hls$, send sequence~$\p$ and birth rate~$\lambda$. 
For any
$t \ge \tau_0$, 
let   $T(t) = \{\tau_0,\ldots,t\}$.
Then   for
every $t >  \tau_0$
 and every tuple
$ 
\psivect_{T(t-1)} = (\psi_{\tau_0},\dots,\psi_{t-1})$  of high-level states 
such that  
$\Pr( \Psivect_{T(t-1)} = \psivect_{T(t-1)}) > 0$, 
the probability that the transition to $\Psi(t)$ is via 
(R\ref{item:R1})(iii), conditioned on
$\Psivect_{T(t-1)} = \psivect_{T(t-1)}$,
is at most $\exp(- (j^{\psi_{t-1}})^{1.9})$.
\end{lemma}
\begin{proof}

Let $\calE_{t-1}$ be the event that 
$\Psivect_{T(t-1)} = \psivect_{T(t-1)}$. Let $j = j^{\psi_{t-1}}$.
The transition is defined by 
${\Psi(t) }
= \hls(\Psi(t-1), \bbar^{Y_{t-1}}_{[j]}(t), t)$   where $Y_{t-1}$ is the process defined in the definition of the volume process (Definition~\ref{def:volume-process}).
We may assume that
$\Failure\not\in \{\type^{\psi_{\tau_0}},\ldots,\type^{\psi_{t-1}}\}$
since otherwise the transition would take (R\ref{item:R1})(i).

By Lemma~\ref{lem:backoff-bounding-valid},
$\hls$ is $\OurHLS$-valid. 
Since $\Failure\not\in \{\type^{\psi_{\tau_0}},\ldots,\type^{\psi_{t-1}}\}$, Lemma~\ref{lem:dist-volume} ensures the following.
\begin{equation}\label{eq:use}
\dist(\bbar_{[j]}^{Y_{t-1}}(t) \mid  \calE_{t-1}) = \Po(F_{[j]}^{\psi_{t-1}}(t)).\end{equation}

We may assume that, for all $k\in \Upsilon_{j,\geq j^{\chi}}$, $F_k^{\psi_{t-1}}(t) \geq \lambda W_k/2$ since otherwise (R\ref{item:R1})(ii) would 
occur.

To upper-bound the probability that (R\ref{item:R1})(iii) occurs,   we wish to upper bound the probability that,
for some $k\in \Upsilon_{j,\geq j^2}$,
$b^{Y_{t-1}}_k(t)< \lambda W_k/4$. 
By~\eqref{eq:use} this is upper-bounding the probability that a variable that is distributed as 
$\Po(F_k^{\psi_{t-1}}(t))$
is at most $\lambda W_k/4$ where $
F_k^{\psi_{t-1}}(t) \geq \lambda W_k/2$ since $\chi\geq 2$.

By Lemma~\ref{lem:chernoff-small-dev} with $\delta=1/2$, this
probability is at most $\exp(-\lambda W_k/2^4)\leq \exp(- \lambda j^2 2^{-4})$, so summing over~$k$, 
the upper bound   is 
$j\exp(- \lambda j^2 2^{-4})
\leq \exp(- j^{1.9})$,
where the final inequality comes from $j \geq j_0$ 
where $j_0$  is sufficiently large.  
 \end{proof}

We next consider transition (R\ref{item:R2}). We start with an easy observation, Observation~\ref{obs:no-R2-fromrefilling}, showing that that (R\ref{item:R2})(i) cannot occur from a refilling state.

\begin{observation}
\label{obs:no-R2-fromrefilling}
Fix a real number $\lambda \in (0,1)$ and a send sequence $\p$ with $p_0=1$. Let $j_0$ be a sufficiently large integer.
Let  
$\Psi_0 
$ 
be the initialising state
with 
$j^{\Psi_0} = j_0$ and
$\tau^{\Psi_0} = \tau_0$.
Let $\hls$ be  
the $\Psi_0$-backoff bounding rule from Definition~\ref{def:backoff-bounding-rule}.
Let $V$ be a volume process from~$\Psi_0$ with transition rule~$\hls$, send sequence~$\p$ and birth rate~$\lambda$. 
For any
$t \ge \tau_0$, 
let   $T(t) = \{\tau_0,\ldots,t\}$.
Then   for
every $t >  \tau_0$
 and every tuple
$ 
\psivect_{T(t-1)} = (\psi_{\tau_0},\dots,\psi_{t-1})$  of high-level states such that $\type^{\psi_{t-1}} = \refilling$, 
if
$\Psivect_{T(t-1)} = \psivect_{T(t-1)}$, then
the transition to $\Psi(t)$ is not via (R\ref{item:R2})(i).
\end{observation}

\begin{proof}
The transition is defined by 
${\Psi(t) }
= \hls(\psi_{t-1}, \bbar^{Y_{t-1}}_{[j]}(t), t)$   where $Y_{t-1}$ is the process defined in the definition of the volume process (Definition~\ref{def:volume-process}) and $j = j^{\psi_{t-1}}$. 
For any set $S$ of bins, 
let $\calE(S)$ be the event that
$\sum_{x\in S} p_x b^{Y_{t-1}}_x(t) < \lambda |S|/40$.

Transition (R\ref{item:R2})(i) is taken from the refilling state~$\psi_{t-1}$  if
$\calE(S)$ occurs for some $S\in \calS^{\psi_{t-1}}$.
From the definition of refilling state (Definition~\ref{def:refilling-state})
the only possibility is that $S = \{k\}$ where $k\in [j-1]$
has weight at least~$j^{\chi}$ 
but $\calE(S)$ would trigger (R\ref{item:R1}).
\end{proof}

Lemma~\ref{lem:AdvancingNoiseFail} upper-bounds the probability of  (R\ref{item:R2})(i) 
for the remaining case where $\type^{\psi_{t-1}}$ is advancing.  
It also rules out
transition (R\ref{item:R3})(i) 
to non-advancing states since the proof technique is the same. Lemma~\ref{lem:AdvancingIncj-NoiseFail} is almost identical 
to Lemma~\ref{lem:AdvancingNoiseFail} except that it rules out transition
(R\ref{item:R4}(i))(a) to non-advancing states  
and this transition would increase~$j$, so the statement of the lemma is adjusted accordingly.

\begin{lemma}
\label{lem:AdvancingNoiseFail}

Fix a real number $\lambda \in (0,1)$ and a send sequence $\p$ with $p_0=1$. Let $j_0$ be a
sufficiently large integer.  
Let  
$\Psi_0 
$ 
be the initialising state
with 
$j^{\Psi_0} = j_0$ and $\tau^{\Psi_0} = \tau_0$.
Let $\hls$ be  
the $\Psi_0$-backoff bounding rule from Definition~\ref{def:backoff-bounding-rule}.
Let $V$ be a volume process from~$\Psi_0$ with transition rule~$\hls$, send sequence~$\p$ and birth rate~$\lambda$. 
For any
$t \geq \tau_0$, 
let   $T(t) = \{\tau_0,\ldots,t\}$.
Then  the following statements hold for
every $t >  \tau_0$ and
every tuple
$ 
\psivect_{T(t-1)} = (\psi_{\tau_0},\dots,\psi_{t-1})$  of high-level states 
such that  
$\Pr( \Psivect_{T(t-1)} = \psivect_{T(t-1)}) > 0$,
where $j=j^{\psi_{t-1}}$.
Let $\calE[1]$ be the event that 
$\Psi(t-1)$ is advancing and
the transition to $\Psi(t)$ is
via (R\ref{item:R2})(i).
Let $\calE[2]$ be the event that 
the transition to $\Psi(t)$ is via  (R\ref{item:R3})(i)  and
$\Psi(t)$ is not an advancing state (so it is a failure state).
\begin{enumerate}[(i)]
\item  If  bin $j$ is many-covered
then $\Pr(\calE[1] \mid \Psivect_{T(t-1)} = \psivect_{T(t-1)}) \leq 
\exp(-  \lambda |\Upsilon_{j,\geq \Wtilde[j]}| \, \Wtilde[j] /80)$.

\item If  bin $j$ is heavy-covered 
then $\Pr(\calE[1] \mid \Psivect_{T(t-1)} = \psivect_{T(t-1)}) \leq 
\exp(-  j^{1.9})$.

\item  If  bin $j$ is many-covered
then $\Pr(\calE[2] \mid \Psivect_{T(t-1)} = \psivect_{T(t-1)}) \leq 
\exp(-  \lambda |\Upsilon_{j,\geq \Wtilde[j]}| \, \Wtilde[j] /80)$.

\item If  bin $j$ is heavy-covered 
then $\Pr(\calE[2] \mid \Psivect_{T(t-1)} = \psivect_{T(t-1)}) \leq 
\exp(-  j^{1.9})$.

\end{enumerate}
\end{lemma}
\begin{proof}

Fix $t > \tau_0$. Let $\calE_{t-1}$ be the event that 
$\Psivect_{T(t-1)} = \psivect_{T(t-1)}$.  
The transition is defined by 
${\Psi(t) }
= \hls(\psi_{t-1}, \bbar^{Y_{t-1}}_{[j]}(t), t)$   where $Y_{t-1}$ is the process defined in the definition of the volume process (Definition~\ref{def:volume-process}). 
For any set $S$ of bins, 
let $\calE(S)$ be the event that
$\sum_{x\in S} p_x b^{Y_{t-1}}_x(t) < \lambda |S|/40$.

If  
$\calE[1]$ occurs 
then $\calE(S)$ occurs for some $S \in \calS^{\psi_{t-1}}$ 
for the $j$-advancing state~$\psi_{t-1}$.
If $\calE[2]$ occurs then 
(T\ref{item:trans3}) fails
while attempting to move to a $j$-advancing state~$\Psi'$ via (R\ref{item:R3})(i).
From Definition~\ref{def:transition-rule} this means that $\calE(S)$ occurs for some $S\in \calS^{\Psi'}$.
Let $\Psi = \psi_{t-1}$ for upper-bounding the probability of $\calE[1]$ and let $\Psi = \Psi'$ for upper-bounding the probability of~$\calE[2]$.
From the definition of advancing state (Definition~\ref{def:advancing-state}), 
$\calS^{\Psi}$ contains a single set~$S\subseteq [j-1]$ (the exact definition of~$S$ depends on whether bin~$j$ is many-covered or heavy-covered).
We will prove the bounds in the lemma statement by upper-bounding $\Pr(\calE(S) \mid \calE_{t-1})$ in each case.

We may assume that none of the high-level states in $\psivect_{T(t-1)}$ are failure states since all transitions from failure states are via (R\ref{item:R1}).
By Lemma~\ref{lem:dist-volume}(i), 
$\dist(\bbar_{[j]}^{Y_{t-1}}(t) 
\mid   \calE_{t-1}   ) =  \Po(F_{[j]}^{\psi_{t-1}}(t))$.
Given the transitions that we are considering, $ \type^{\psi_{t-1}} \in \{ \initialising, \advancing\}$.
By Lemma~\ref{lem:vol-F-lb}, 
for every $k\in [j -1]$,
$F_k^{\psi_{t-1}}(t) \geq \mu_k^{
\gammavect}\geq 2 \lambda W_k/3 > \lambda W_k/10$.
For each $k\in S$, 
$b^{Y_{t-1}}_k(t) $ is an independent Poisson variable with mean at least $\lambda W_k/10$.
If, for some real number $W\geq 1$,
all $k\in S$ have $W_k\geq W$ then
Corollary~\ref{cor:poisson-noisy}
gives $\Pr(\calE(S) \mid \calE_{t-1}) \leq \exp(- \lambda |S|\, W/80)$.

Items (i) and (iii) follow from the fact that, when bin~$j$ is many-covered,
$S = \Upsilon_{j, \geq \Wtilde[j]}$.
By definition, all $k\in S$ have $W_k \geq \Wtilde[j]$.

Items (ii) and (iv) follow from the fact that, when bin~$j$ is heavy-covered,
$S = \Upsilon_{j,\geq j^2}$. By the definition of heavy-covered (Definition~\ref{def:covered})
$|S| \geq \CUpsilon L(j)/(2 \lambda)$, so taking $W=j^2$ and using the fact that $j\geq j_0$ is sufficiently large, 
$\lambda |S| W /80 \geq 
j^2 \CUpsilon L(j)/160 \geq j^{1.9}$. 
\end{proof}

Lemma~\ref{lem:AdvancingIncj-NoiseFail} is similar to Lemma~\ref{lem:AdvancingNoiseFail} except that the relevant transition increases~$j$, so the case analysis in the statement is on 
whether $j^{\Psi(t)}$ is many-covered or heavy-covered and not on whether $j^{\psi_{t-1}}$ is many-covered or heavy-covered. The proof is otherwise the same as the $\calE_2$ case of the proof of Lemma~\ref{lem:AdvancingNoiseFail}.

\begin{lemma}
\label{lem:AdvancingIncj-NoiseFail}

Fix a real number $\lambda \in (0,1)$ and a send sequence $\p$ with $p_0=1$. Let $j_0$ be a
sufficiently large integer.  
Let  
$\Psi_0 
$ 
be the initialising state
with 
$j^{\Psi_0} = j_0$ and $\tau^{\Psi_0} = \tau_0$.
Let $\hls$ be  
the $\Psi_0$-backoff bounding rule from Definition~\ref{def:backoff-bounding-rule}.
Let $V$ be a volume process from~$\Psi_0$ with transition rule~$\hls$, send sequence~$\p$ and birth rate~$\lambda$. 
For any
$t \geq \tau_0$, 
let   $T(t) = \{\tau_0,\ldots,t\}$.
Then  the following statements hold for
every $t >  \tau_0$ and
every tuple
$ 
\psivect_{T(t-1)} = (\psi_{\tau_0},\dots,\psi_{t-1})$  of high-level states 
such that  
$\Pr( \Psivect_{T(t-1)} = \psivect_{T(t-1)}) > 0$.
Let  $j = j^{\psi_{t-1}}+1$ and let
$\calE$ be the event that 
the transition to $\Psi(t)$ is via  (R\ref{item:R4})(i)(a) and
$\Psi(t)$ is not an advancing state (so it is a failure state).
\begin{enumerate}[(i)]

\item  If  bin $j$ is many-covered
then $\Pr(\calE \mid \Psivect_{T(t-1)} = \psivect_{T(t-1)}) \leq 
\exp(-  \lambda |\Upsilon_{j,\geq \Wtilde[j]}| \, \Wtilde[j] /80)$.

\item If  bin $j$ is heavy-covered 
then $\Pr(\calE \mid \Psivect_{T(t-1)} = \psivect_{T(t-1)}) \leq 
\exp(-  j^{1.9})$.

\end{enumerate}
\end{lemma}
\begin{proof}

Fix $t > \tau_0$. Let $\calE_{t-1}$ be the event that 
$\Psivect_{T(t-1)} = \psivect_{T(t-1)}$.  
The transition is defined by 
${\Psi(t) }
= \hls(\psi_{t-1}, \bbar^{Y_{t-1}}_{[j-1]}(t), t)$   where $Y_{t-1}$ is the process defined in the definition of the volume process (Definition~\ref{def:volume-process}). 
For any set $S$ of bins, 
let $\calE(S)$ be the event that
$\sum_{x\in S} p_x b^{Y_{t-1}}_x(t) < \lambda |S|/40$.

If $\calE$ occurs then 
(T\ref{item:trans3}) fails
while attempting to move to a $j$-advancing state~$\Psi$ via (R\ref{item:R4})(i)(a).
From Definition~\ref{def:transition-rule} this means that $\calE(S)$ occurs for some $S\in \calS^{\Psi}$.
From the definition of advancing state, 
$\calS^{\Psi}$ contains a single set~$S\subseteq [j-1]$.
We will prove the bounds in the lemma statement by upper-bounding $\Pr(\calE(S) \mid \calE_{t-1})$ in each case.

We may assume that none of the high-level states in $\psivect_{T(t-1)}$ are failure states since all transitions from failure states are via (R\ref{item:R1}).
By Lemma~\ref{lem:dist-volume}(i), 
$\dist(\bbar_{[j]}^{Y_{t-1}}(t) 
\mid   \calE_{t-1}   ) =  \Po(F_{[j]}^{\psi_{t-1}}(t))$.
Given the transition that we are considering, $ \type^{\psi_{t-1}} \in \{   \advancing,\filling\}$.
By Lemma~\ref{lem:vol-F-lb}, 
for every $k\in [j -1]$,
$F_k^{\psi_{t-1}}(t) \geq \mu_k^{
\gammavect}\geq 2 \lambda W_k/3 > \lambda W_k/10$.
For each $k\in S$, 
$b^{Y_{t-1}}_k(t) $ is an independent Poisson variable with mean at least $\lambda W_k/10$.
If, for some real number $W\geq 1$,
all $k\in S$ have $W_k\geq W$ then
Corollary~\ref{cor:poisson-noisy}
gives $\Pr(\calE(S) \mid \calE_{t-1}) \leq \exp(- \lambda |S|\, W/80)$.

Items (i) follows from the fact that, when bin~$j$ is many-covered,
$S = \Upsilon_{j, \geq \Wtilde[j]}$.
By definition, all $k\in S$ have $W_k \geq \Wtilde[j]$.

Items (ii) follows from the fact that, when bin~$j$ is heavy-covered,
$S = \Upsilon_{j,\geq j^2}$. By the definition of heavy-covered (Definition~\ref{def:covered})
$|S| \geq \CUpsilon L(j)/(2 \lambda)$, so taking $W=j^2$ and using the fact that $j\geq j_0$ is sufficiently large, 
$\lambda |S| W /80 \geq 
j^2 \CUpsilon L(j)/160 \geq j^{1.9}$. 
\end{proof}

The next lemma, Lemma~\ref{lem:T3tofill}, is similar to Lemmas~\ref{lem:AdvancingNoiseFail} and~\ref{lem:AdvancingIncj-NoiseFail} except that it applies to filling states rather than to advancing states.
It upper-bounds the probability of (R\ref{item:R2})(ii) when $\type^{\psi_{t-1}}$ is filling. It also rules out 
transitions (R\ref{item:R3})(ii), 
(R\ref{item:R4})(i)(b) 
and (R\ref{item:R6})(i)
to non-filling states   since the proof technique is the same.
The proof technique is broadly similar to that of Lemma~\ref{lem:AdvancingNoiseFail} 
and~\ref{lem:AdvancingIncj-NoiseFail}
except that singleton high-weight bins in the set $\calS$ corresponding to the filling state must be considered, and high-weight bins are deleted from the set $\Upsilon_{j,\geq \Wtilde[j]}$ 
in~$\calS$ (which complicates the proof a little, but only reduces the size of~$S$ by a constant factor).

\begin{lemma}\label{lem:T3tofill}
Fix a real number $\lambda \in (0,1)$ and a send sequence $\p$ with $p_0=1$. Let $j_0$ be a sufficiently large  integer.  
Let  $\Psi_0$ be the initialising state
with $j^{\Psi_0} = j_0$ and
$\tau^{\Psi_0} = \tau_0$. 
Let $\hls$ be  
the $\Psi_0$-backoff bounding rule from Definition~\ref{def:backoff-bounding-rule}.
Let $V$ be a volume process from~$\Psi_0$ with transition rule~$\hls$, send sequence~$\p$ and birth rate~$\lambda$. 
For any
$t \geq \tau_0$, 
let   $T(t) = \{\tau_0,\ldots,t\}$.
Then  the following statements hold  for
every $t >  \tau_0$ and
every tuple
$ 
\psivect_{T(t-1)} = (\psi_{\tau_0},\dots,\psi_{t-1})$  of high-level states 
such that  
$\Pr( \Psivect_{T(t-1)} = \psivect_{T(t-1)}) > 0$.
\begin{itemize}
\item  If $\type^{\psi_{t-1}} = \filling$ then the probability 
that the transition to $\Psi(t)$ is
via (R\ref{item:R2})(ii),
conditioned on
$\Psivect_{T(t-1)} = \psivect_{T(t-1)}$,
is   at most  $\exp(- \lambda(j^{\psi_{t-1}}-1)/400)$. 
\item 
The probability that 
$\Psi(t)$ is not a filling state and
the transition is via  
(R\ref{item:R3})(ii), 
(R\ref{item:R4})(i)(b), 
or (R\ref{item:R6})(i), conditioned on
$\Psivect_{T(t-1)} = \psivect_{T(t-1)}$,
is at most $\exp(- \lambda(j^{\psi_{t-1}}-1)/400)$. 
\end{itemize} 
\end{lemma}

\begin{proof}
Fix $t > \tau_0$. Let $\calE_{t-1}$ be the event that 
$\Psivect_{T(t-1)} = \psivect_{T(t-1)}$. Let $j^- = j^{\psi_{t-1}}$.
The transition is defined by 
${\Psi(t) }
= \hls(\psi_{t-1}, \bbar^{Y_{t-1}}_{[j^-]}(t), t)$   where $Y_{t-1}$ is the process defined in the definition of the volume process (Definition~\ref{def:volume-process}). 
For any set $S$ of bins, 
let $\calE(S)$ be the event that
$\sum_{x\in S} p_x b^{Y_{t-1}}_x(t) < \lambda |S|/40$.

The first item in the lemma statement occurs if (R\ref{item:R2})(ii) is taken from a filling state~$\psi_{t-1}$. Thus, it occurs if
$\calE(S)$ occurs for some $S\in \calS^{\psi_{t-1}}$ for the filling state $\psi_{t-1}$.
The second item in the lemma statement occurs if  
transitions  
(R\ref{item:R3})(ii), 
(R\ref{item:R4})(i)(b), 
or (R\ref{item:R6})(i)  violate (T\ref{item:trans3}) while attempting to move to a filling state~$\Psi$. 
From Definition~\ref{def:transition-rule} this means that $\calE(S)$ occurs for some $S\in \calS^{\Psi}$.

From the definition of 
a filling state (Definition~\ref{def:filling-state}), 
$\calS^{\psi_{t-1}}$ and $\calS^{\Psi}$ are both of the form 
$\calS = \{\Upsilon_{j, \ge \Wtilde[j]} \setminus \Upsilon_{j, \ge j^2}\}
\cup \{\{k\}\colon k \in \Upsilon_{j, \ge j^\chi}\}$ where $j\in \{j^-,j^-+1\}$ (depending on the transition).
We will complete the proof by showing 
that  the sum of $\Pr( \calE(S) \mid \calE_{t-1}
  )$, summed over $S\in \calS$, is at most $\exp(- \lambda(j^{-}-1)/400)$.

First consider $S = \{ k\}$ where $k\in [j-1]$ has weight at least~$j^{\chi}$. 
If $k< j^-$ then $\calE(\{k\})$ cannot trigger any of these transitions
because   $\calE(\{k\})$ would trigger (R\ref{item:R1}) first. 
One more case arises if the transition is (R\ref{item:R4})(i)(b) 
so that
 $j = j^-+1$ and $W_{j^-} \geq {j}^\chi$. In this case, we use the same argument as Lemma~\ref{lem:volume-easyChernoffs} to show
$\Pr(\calE(\{j^-\}) \mid \calE_{t-1}) \leq 
\exp(- \lambda j^{\chi} 2^{-4})$.
Note from the backoff bounding rule that in this case (by the guard of the transition) $F_{j^-}^{\psi_{t-1}}(t) \geq \mu_{j^-}^{\gammavect}$ which is at least 
$2 \lambda W_{j^-}/3$ by Observation~\ref{obs:mu-gamma-bounds}.  
We may assume that
$\Failure\not\in \{\type^{\psi_{\tau_0}},\ldots,\type^{\psi_{t-1}}\}$
since otherwise the transitions that we are interested in would not occur
in the transition to $\Psi(t)$.   Lemma~\ref{lem:dist-volume} ensures  
that
$\dist(\bbar_{j^-}^{Y_{t-1}}(t) \mid  \calE_{t-1}) = \Po(F_{j^-}^{\psi_{t-1}}(t))$.
By Lemma~\ref{lem:chernoff-small-dev} with $\delta=1/2$, 
$\Pr(\calE(\{j^-\}) \mid \calE_{t-1}) \leq  \exp(-\lambda W_{j^-}/2^4)\leq 
\exp(- \lambda j^{\chi} 2^{-4})$, which is much less than the failure probability in the lemma statement.

Finally, we consider the set
$S = \Upsilon_{j, \ge \Wtilde[j]} \setminus \Upsilon_{j, \ge j^2}$.
We will show that $\Pr(\calE(S) \mid \calE_{t-1}) \leq  \exp(-\lambda (j^{-}-1)/320)$. To get the lemma, we use a union bound, adding this to the (much smaller) failure probability from the previous part.

Bin $j$ is exposed 
from the high-level states involved in all of these transitions. 
By the definition of covered/exposed (Definition~\ref{def:covered}),
Property~(Prop~\ref{cov-prop-2}) is false so 
$|\Upsilon_{j,\geq j^2}| <  \CUpsilon \L(j)/(2\lambda)$.
By the definition of $\Wtilde[j]$  
(Definition~\ref{def:Wtilde}),
$|\Upsilon_{j,\geq \Wtilde[j]}| \geq \CUpsilon \L(j)/\lambda$
so 
$|\Upsilon_{j,\geq j^2}| \leq \tfrac12 |\Upsilon_{j,\geq \Wtilde[j] }|$ which means that $|S| \geq \tfrac12
|\Upsilon_{j,\geq \Wtilde[j]}|$.
At this point, the proof becomes similar to the proof of Lemma~\ref{lem:AdvancingNoiseFail}.

By Lemma~\ref{lem:dist-volume}(i), 
$\dist(\bbar_{[j^-]}^{Y_{t-1}}(t) 
\mid \Psivect_{T(t-1)} = \psivect_{T(t-1)} ) =  \Po(F_{[j^-]}^{\psi_{t-1}}(t))$.
Given the transitions that we are considering, $ \type^- \in \{ \initialising, \advancing, \filling, \stabilising \}$.
By Lemma~\ref{lem:vol-F-lb}, 
for every $k\in [j^- -1]$,
$F_k^{\psi_{t-1}}(t) \geq \mu_k^{
\gammavect}\geq 2 \lambda W_k/3 > \lambda W_k/10$.

Now let $S^- = S \setminus \{j^-\}$ and note that $|S^-| \geq |S|-1 \geq 
\tfrac14
|\Upsilon_{j,\geq \Wtilde[j]}|$, using the fact that $j\geq j_0$ is sufficiently large.
We will now apply Corollary~\ref{cor:poisson-noisy} 
to~$\p$ and~$S^-$ with 
$W = \Wtilde[j]$. 
For each $k\in S^-$, the Poisson variable~$P_k$ from
Corollary~\ref{cor:poisson-noisy} is $b_k^{Y_{t-1}}(t)$ whose mean we have just shown is at least $\mu_k^{\gammavect}$ which is at least 
$ \lambda W_k/10$.
The corollary shows that 

\begin{align*} \Pr( \calE(S) \mid 
  \calE_{t-1})  &\leq
\Pr(\sum_{x\in S^-} p_x b^{Y_{t-1}}_x(t) < \lambda |S|/40
\mid \calE_{t-1}
)\\ &\leq
\Pr(\sum_{x\in S^-} p_x b^{Y_{t-1}}_x(t) < \lambda |S^-|/20 \mid \calE_{t-1})\\
&\leq\exp(-\lambda|S^-| \Wtilde[j]/80).\end{align*}

The bound 
$\Pr(\calE(S) \mid \calE_{t-1}) \leq    \exp(-\lambda (j^{-}-1)/320)
$
follows since 
$\Wtilde[j]|S^-| \geq 
\Wtilde[j] |\Upsilon_{j,\geq \Wtilde[j]}|/4
\geq
(j-1)/4 \geq (j^{-}-1)/4$.
\end{proof}

Transition (R\ref{item:R2}) cannot be taken from initialising and failure states since 
they have an empty set $\calS$. 
Lemma~\ref{lem:T3tostab} upper-bounds the probability of the remaining transition from (R\ref{item:R2}), which is from a stabilising state. It is significant that the failure probability in Lemma~\ref{lem:T3tostab} is higher than the ones that we have encountered so far, but it will still be small enough for us since the outcome is a refilling state, rather than a failing state. The lemma also upper-bounds the probability of transition (R\ref{item:R5})(i) to non-stabilising states.

\begin{lemma}
\label{lem:T3tostab}
Fix a real number $\lambda \in (0,1)$ and a send sequence $\p$ with $p_0=1$. Let $j_0$ be a
sufficiently large integer.  
Let  
$\Psi_0$ be the initialising state with
$j^{\Psi_0} = j_0$ and $\tau^{\Psi_0} = \tau_0$.  
Let $\hls$ be  
the $\Psi_0$-backoff bounding rule from Definition~\ref{def:backoff-bounding-rule}.
Let $V$ be a volume process from~$\Psi_0$ with transition rule~$\hls$, send sequence~$\p$ and birth rate~$\lambda$. 
For any
$t \geq \tau_0$, 
let   $T(t) = \{\tau_0,\ldots,t\}$.
Then  the following statements hold  for
every $t >  \tau_0$ and
every tuple
$ 
\psivect_{T(t-1)} = (\psi_{\tau_0},\dots,\psi_{t-1})$  of high-level states 
such that  
$\Pr( \Psivect_{T(t-1)} = \psivect_{T(t-1)}) > 0$.
\begin{itemize}
\item  If $\type^{\psi_{t-1}} = \stabilising$ then the probability 
that the transition to $\Psi(t)$ is
via (R\ref{item:R2})(ii),
conditioned on
$\Psivect_{T(t-1)} = \psivect_{T(t-1)}$,
is   at most 
$\exp(-\lambda (\log j^{\psi_{t-1}})^2/200)$. 
\item 
The probability that 
$\Psi(t)$ is not a stabilising state and
the transition is via  
(R\ref{item:R5})(i), conditioned on
$\Psivect_{T(t-1)} = \psivect_{T(t-1)}$,
is at most 
$\exp(-\lambda (\log j^{\psi_{t-1}})^2/200)$. 
\end{itemize} 
\end{lemma}
\begin{proof}  
Fix $t > \tau_0$. Let $\calE_{t-1}$ be the event that 
$\Psivect_{T(t-1)} = \psivect_{T(t-1)}$. 
Let $j = j^{\psi_{t-1}}$.
The transition is defined by 
${\Psi(t) }
= \hls(\psi_{t-1}, \bbar^{Y_{t-1}}_{[j]}(t), t)$   where $Y_{t-1}$ is the process defined in the definition of the volume process (Definition~\ref{def:volume-process}). 
For any set $S$ of bins, 
let $\calE(S)$ be the event that
$\sum_{x\in S} p_x b^{Y_{t-1}}_x(t) < \lambda |S|/40$.

The first item in the lemma statement occurs if (R\ref{item:R2})(ii) is taken from~$\psi_{t-1}$. Thus, it occurs if
$\calE(S)$ occurs for some $S\in \calS^{\psi_{t-1}}$ for the stabilising state $\psi_{t-1}$.
The second item in the lemma statement occurs if  
transition  
(R\ref{item:R5})(i)   violates (T\ref{item:trans3}) while attempting to move to a stabilising state~$\Psi$. 
From Definition~\ref{def:transition-rule} this means that $\calE(S)$ occurs for some $S\in \calS^{\Psi}$.

From the definition of stabilising state (Definition~\ref{def:stabilising-state})
$\calS^{\psi_{t-1}}$ and $\calS^{\Psi}$ are both of the form
\[\calS= \{\{ k \in [(\log j)^2] \colon W_k < j^2\}, \{k \in [j-1] \setminus [(\log j)^2] \colon W_k < j^2\}\} \cup \{\{k\}\colon k \in \Upsilon_{j, \ge j^\chi}\}.\]

We will complete the proof by showing 
that  the sum of $\Pr( \calE(S) \mid \calE_{t-1}
  )$, summed over $S\in \calS$, is at most
$\exp(-\lambda (\log j)^2/200)$.

First consider $S = \{ k\}$ where $k\in [j-1]$ has weight at least~$j^{\chi}$. 
These are straightforward -- $\calE(\{k\})$ cannot trigger   the transitions that we are considering
because    it would trigger (R\ref{item:R1}) first. 
The remaining two sets are 
$S_1 = \{ k \in [(\log j)^2] \colon W_k < j^2\}$ and
$S_2 = \{k \in [j-1] \setminus [(\log j)^2] \colon W_k < j^2\}\}$.

By Lemma~\ref{lem:dist-volume}(i), 
$\dist(\bbar_{[j]}^{Y_{t-1}}(t) 
\mid \Psivect_{T(t-1)} = \psivect_{T(t-1)} ) =  \Po(F_{[j]}^{\psi_{t-1}}(t))$. 
For transition (R\ref{item:R2})(ii),
$\type= \stabilising$.  In this case,
by Lemma~\ref{lem:vol-F-lb}, 
for every $k\in [j -1]$,
$F_k^{\psi_{t-1}}(t) \geq \mu_k^{
\gammavect}\geq 2 \lambda W_k/3 > \lambda W_k/10$. 
For transition (R\ref{item:R5})(i), the guard guarantees that for all 
$k\in [j -1]$,
$F_k^{\psi_{t-1}}(t) \geq \mu_k^{
\gammavect}  > \lambda W_k/10$. 

Applying
 Corollary~\ref{cor:poisson-noisy}
 to each $S\in \{S_1,S_2\}$, we get 
 \[\Pr( \calE(S) \mid 
  \calE_{t-1})  
\leq\exp(-\lambda|S| \Wtilde[j]/80)
\leq\exp(-\lambda|S|/80). \] 
To finish, we require lower bounds on~$|S_1|$ and $|S_2|$.
Since $\psi_{t-1}$ is 
a refilling or stabilising state,  bin $j$ is exposed. 
By the definition of covered/exposed (Definition~\ref{def:covered}), (Prop~\ref{cov-prop-2}) is false, so 
$|\Upsilon_{j,\geq j^2}| <  \CUpsilon \L(j)/(2\lambda)$, so $|S_1| \geq
\lfloor(\log j)^2\rfloor - \CUpsilon \L(j)/(2\lambda)$ which is at least $(\log j)^2/2$ since $j\geq j_0$ is sufficiently large.
Finally, $|S_2| \ge j-1 - \lceil (\log j)^2\rceil - C_\Upsilon L(j)/(2\lambda) > (\log j)^2/2$.
The failure probability is 
therefore at most
$2\exp(-\lambda (\log j)^2/160)
\leq \exp(-\lambda (\log j)^2/200)$.
\end{proof}

\subsection{Analysis of the refilling-stabilising loop}\label{sec:volume-theloop}

In this section we analyse the progress of the volume process through the refilling-stabilising loop. 

The main lemma in the section is Lemma~\ref{lem:main-loop-bounds} that
says, starting from a transition from a $j$-filling high-level state~$\Psi$
to a $j$-refilling high-level state, the volume process is likely to (fairly) quickly escape from the refilling-stabilising loop, transitioning back to a $j$-filling
high level state that has not lost many balls from bin~$j$ (in expectation), relative to~$\Psi$.

Recall from Definition~\ref{def:back-transition}
that  \emph{back transitions} are
the back-transitions of the refilling-stabilising loop. These are
the self-loop transition in Figure~\ref{fig:backoff-bounding-rule}
from a refilling state to another refilling state via (R\ref{item:R5})(i) and the left-arrow transition in the figure from a stabilising state to a refilling state via (R\ref{item:R2})(ii) or (R\ref{item:R6})(i). 

The stopping times in Definition~\ref{def:stopping-times} will be useful in this section.
\begin{definition}
\label{def:stopping-times}
Fix a real number $\lambda \in (0,1)$ and a send sequence $\p$ with $p_0=1$. Let $j_0$ be a
sufficiently large integer.  
Let  
$\Psi_0$ be the initialising state with
$j^{\Psi_0} = j_0$ and $\tau^{\Psi_0} = \tau_0$.  
Let $\hls$ be  
the $\Psi_0$-backoff bounding rule from Definition~\ref{def:backoff-bounding-rule}.
Let $V$ be a volume process from~$\Psi_0$ with transition rule~$\hls$, send sequence~$\p$ and birth rate~$\lambda$. 
For any
$t \geq \tau_0$,  
let   $\newsym{T-t}{sequence of times - Def~\ref{def:stopping-times}}{T(t) = \{\tau_0,\ldots,t\}}$ and define the following stopping times (which could possibly have the value~$\infty$).
\begin{align*}
\newsym{trans}{next change of state - Def~\ref{def:stopping-times}}{\trans({t})} &= \min\{t' > {t} \mid \Psi(t') \neq \Psi({t})\},\\
\newsym{back}{next back transition - Def~\ref{def:stopping-times}}{\back(t)} &= \min \{ t' > {t} \mid  
\Psi(t' - 1) \rightarrow 
\Psi(t') \mbox{ is a back transition} 
\},\\
\newsym{exit}{next exit from loop - Def~\ref{def:stopping-times}}{\exit({t})} &= \min\{t' > {t} \mid \type^{\Psi(t')}  \notin \{\refilling,\stabilising\}\},\\
\newsym{backexit}{min of $\back(t)$ and $\exit(t)$ - Def~\ref{def:stopping-times}}{\backexit(t)} &= \min \{\back(t),\exit(t)\}.
\end{align*}

\end{definition}

We start with Lemma~\ref{lem:nonback}, which  says that, starting from a
state in the refilling-stabilising loop  without too many missing balls, it does not take too long before the next back transition, or the next transition out of the refilling-stabilising loop. Lemma~\ref{lem:tautilde-unlikely-back} then gives an upper bound on the probability that this transition is a back transition (rather than an exit from the loop).

\begin{lemma}
\label{lem:nonback} 
Fix a real number $\lambda \in (0,1)$ and a send sequence $\p$ with $p_0=1$. Let $j_0$ be a
sufficiently large integer.  
Let  
$\Psi_0$ be the initialising state with
$j^{\Psi_0} = j_0$ and $\tau^{\Psi_0} = \tau_0$.  
Let $\hls$ be  
the $\Psi_0$-backoff bounding rule from Definition~\ref{def:backoff-bounding-rule}.
Let $V$ be a volume process from~$\Psi_0$ with transition rule~$\hls$, send sequence~$\p$ and birth rate~$\lambda$.   For any
$t \ge \tau_0$, 
let   $T(t) = \{\tau_0,\ldots,t\}$.
Then  the following statement holds  for
every $t >  \tau_0$ and
every tuple
$ 
\psivect_{T(t-1)} = (\psi_{\tau_0},\dots,\psi_{t-1})$  of high-level states 
such that  
$\Pr( \Psivect_{T(t-1)} = \psivect_{T(t-1)}) > 0$ and 
$\type^{\psi_{t-1}} \in \{ \refilling,\stabilising\}$, where
$j=j^{\psi_{t-1}}$.  
If
$\Psivect_{T(t-1)} = \psivect_{T(t-1)}
$ 
then 
\begin{enumerate}[(i)]
\item for all $t'$ in the range $t-1 \leq t' \leq \backexit(t-1)$, $j^{\Psi(t')} = j$,
\item for all $\hat{t}$ in  the range $t-1 \le \hat{t} < \backexit(t-1)$ with  $\type^{\Psi(\hat{t})} = \stabilising$,  
$\trans(\hat{t}) \leq \hat{t} + j^{\phi+72}$, 
\item 
If $\calM_{[j-1]}^{\psi_{t-1}}(\tau^{\psi_{t-1}}) \leq 
2j^{44} 
$ then
for all 
$\hat{t}$ in the range $t-1 \le \hat{t} < \backexit(t-1)$
with $\type^{\Psi(\hat{t})} = \refilling$,   
$\trans(\hat{t})  \leq \hat{t} + j^{\phi+{46}}$, and
\item 
If $\calM_{[j-1]}^{\psi_{t-1}}(\tau^{\psi_{t-1}}) \leq 
2j^{44} 
$ then
$\backexit(t-1) \leq t-1 + j^{\phi+73}$. 
\end{enumerate}
\end{lemma}

\begin{proof}  Item~(i) follows from the fact that, according to the backoff bounding rule,  transitions from refilling and stabilising states don't change~$j$. 
Also from the backoff bounding rule,
there can only be one $\hat{t}$
in the range $t-1 < \hat{t} < \backexit(t-1)$  
such that $\Psi(\hat{t}) \neq \Psi(\hat{t}-1)$ -- if there is such 
a~$\hat{t}$ then $\Psi(t-1) = \cdots = \Psi(\hat{t}-1)$ 
is a $j$-refilling state
and
$\Psi(\hat{t})= \cdots = \Psi(\backexit(t-1)-1) $ is a $j$-stabilising state. Thus $\hat{t} = \trans(t-1)$
and $\backexit(t-1) = \trans(\hat{t})$.
Thus, item~(iv) follows from items~(ii) and~(iii)  since
$j\geq j_0$ is sufficiently large so $j^{\phi+46} + j^{\phi+72} \leq j^{\phi+73}$. 

Item~(ii) follows immediately from Lemma~\ref{lem:stabilising-transient}.
For  item~(iii), we've observed that for all 
$\hat{t}$ in the range $t-1 \le \hat{t} < \backexit(t-1)$
with $\type^{\Psi(\hat{t})} = \refilling$, $\Psi(\hat{t}) = \psi_{t-1}$, so the hypothesis of the lemma guarantees that 
$\calM_{[j-1]}^{\Psi(\hat{t})}(\tau^{\Psi(\hat{t})}) \leq 
2j^{44}  
$. By Lemma~\ref{lem:refilling-transient} 
there is an integer $t''\in \{\tau^{\Psi(\hat{t})}+1,\ldots,\tau^{\Psi(\hat{t})}+j^{\phi+46}\}$  
such that $\Psi(t'') \neq \Psi(\hat{t})$, giving the result.
\end{proof}

\begin{lemma}
\label{lem:tautilde-unlikely-back}
Fix a real number $\lambda \in (0,1)$ and a send sequence $\p$ with $p_0=1$. Let $j_0$ be a
sufficiently large integer.  
Let  
$\Psi_0$ be the initialising state with
$j^{\Psi_0} = j_0$ and $\tau^{\Psi_0} = \tau_0$.  
Let $\hls$ be  
the $\Psi_0$-backoff bounding rule from Definition~\ref{def:backoff-bounding-rule}.
Let $V$ be a volume process from~$\Psi_0$ with transition rule~$\hls$, send sequence~$\p$ and birth rate~$\lambda$. For any
$t \ge \tau_0$, 
let   $T(t) = \{\tau_0,\ldots,t\}$.
Then  the following  holds for
every $t >  \tau_0$ and
every tuple
$ 
\psivect_{T(t-1)} = (\psi_{\tau_0},\dots,\psi_{t-1})$  of high-level states 
such that  
$\Pr( \Psivect_{T(t-1)} = \psivect_{T(t-1)}) > 0$, 
$\type^{\psi_{t-1}} \in \{ \refilling,\stabilising\}$, 
and 
$\calM_{[j-1]}^{\psi_{t-1}}(\tau^{\psi_{t-1}}) \leq 
2j^{44} 
$, where $j=j^{\psi_{t-1}}$.  
\[
\pr(\mbox{the transition into 
$\Psi(\backexit(t-1))$ is a back transition }  \mid 
\Psivect_{T(t-1)} = \psivect_{T(t-1)}
 )\leq 
 \exp(-\lambda (\log j )^2/202).
\]

\end{lemma}

\begin{proof} 
For any $t' \geq t$ and
 every tuple
$\psivect_{T(t'-1)}$
of high-level states that extends
$\psivect_{T(t-1)}$, has
$\Pr( \Psivect_{T(t'-1)} = \psivect_{T(t'-1)}) > 0$,  
and entails $\backexit(t-1) \geq t'$, we will show 
\begin{equation}\label{eq:todo}
\pr( 
\mbox{the transition into  $
\Psi(t')$  is a back transition }  
\mid 
 \Psivect_{T(t'-1)} = \psivect_{T(t'-1)}
 ) \leq 
 \exp(-\lambda (\log j )^2/201).
\end{equation}
By Lemma~\ref{lem:nonback}(i), 
$j^{\psi_{t'-1}} = j$ so
the relevant transition is $\Psi(t') = \hls(\psi_{t'-1}, \bbar^{Y_{t'-1}}_{[j]}(t'),t')$.
We consider cases based on $\type^{\psi_{t'-1}}$.

If $\type^{\psi_{t'-1}} = \refilling$ then the only possible back transition is via (R\ref{item:R5})(i), and Lemma~\ref{lem:T3tostab}
shows that the probability of taking this transition is at most 
$\exp(-\lambda (\log j )^2/200)$.  

If $\type^{\psi_{t'-1}} = \stabilising$ then the two possible back transitions are via (R\ref{item:R2})(ii) or (R\ref{item:R6})(i). Lemma~\ref{lem:T3tostab} shows that the probability of taking an (R\ref{item:R2})(ii) transition is at most $\exp(-\lambda (\log j )^2/200)$,   and Lemma~\ref{lem:T3tofill} shows that the probability of taking an (R\ref{item:R6})(i) transition is at most 
$\exp(- \lambda(j-1)/400))$. 
Using the fact that $j\geq j_0$ is sufficiently large 
so $\exp(-\lambda (\log j )^2/200)
+ \exp(- \lambda(j-1)/400))
\leq \exp(-\lambda (\log j )^2/201)$,
in both cases we have proved~\eqref{eq:todo}.

By Lemma~\ref{lem:nonback}(iv), which applies because 
$\calM_{[j-1]}^{\psi_{t-1}}(\tau^{\psi_{t-1}}) \leq 
2j^{44} 
$ by hypothesis, 
assuming $\Psivect_{T(t-1)} = \psivect_{T(t-1)}$,
we obtain $\backexit(t-1) \le t-1 + j^{\phi+73}$. Thus we obtain
\begin{align*}
&\pr(
\mbox{the transition into  $
\Psi(\backexit(t-1))$  is a back transition }  
 \mid 
\Psivect_{T(t-1)} = \psivect_{T(t-1)}
 )\\
=&\sum_{t' = t}^{t-1+j^{\phi+73}} \pr(\backexit(t-1) = t' 
\mbox{and the transition into  $
\Psi(t')$  is a back transition }  
 \mid 
\Psivect_{T(t-1)} = \psivect_{T(t-1)}
 )\\
\leq&\sum_{t' = t}^{t-1+j^{\phi+73}} 
\max_{\psivect_{T(t'-1)}} 
 \pr(
 \mbox{the transition into  $
\Psi(t')$  is a back transition }     
 \mid 
 \Psivect_{T(t'-1)} = \psivect_{T(t'-1)}
  ),
\end{align*}
where the maximum is over 
$\psivect_{T(t'-1)}$ which extend
$\psivect_{T(t-1)}$ and have
$\Pr( \Psivect_{T(t'-1)} = \psivect_{T(t'-1)}) > 0$  
and entail $\backexit(t-1) \geq t'$. 
We have already proved each term of this sum is at most 
$\exp(-\lambda (\log j )^2/201)$,  so the result follows since $j\geq j_0$ is sufficiently large.
\end{proof}

Lemma~\ref{lem:leave-loop-quickly} says that, after starting in a $j$-filling state and moving to a $j$-refilling state (entering the loop), with high probability it takes only $\poly(j)$ steps before the process leaves the loop. Moreover, 
at most $2j$   distinct high-level states   are entered before this happens. The point of Lemma~\ref{lem:leave-loop-quickly} is that the failure probability is very low. By Lemma~\ref{lem:tautilde-unlikely-back}, the probability that a back transition happens before exiting the loop is $\exp(-\Omega( (\log j)^2))$. So Lemma~\ref{lem:leave-loop-quickly} shows that the probability that this  happens $j$~times in a row is $\exp(-\Omega( j(\log j)^2))$ -- this failure probability is small enough that we will later be able to apply a union bound.

\begin{lemma}
\label{lem:leave-loop-quickly}
Fix a real number $\lambda \in (0,1)$ and a send sequence $\p$ with $p_0=1$. Let $j_0$ be a
sufficiently large integer.  
Let  
$\Psi_0$ be the initialising state with
$j^{\Psi_0} = j_0$ and $\tau^{\Psi_0} = \tau_0$.  
Let $\hls$ be  
the $\Psi_0$-backoff bounding rule from Definition~\ref{def:backoff-bounding-rule}.
Let $V$ be a volume process from~$\Psi_0$ with transition rule~$\hls$, send sequence~$\p$ and birth rate~$\lambda$. For any
$t \ge \tau_0$, 
let   $T(t) = \{\tau_0,\ldots,t\}$.
Then  the following statement holds  for
every $t \geq  \tau_0+2$ and
every tuple
$ 
\psivect_{T(t-1)} = (\psi_{\tau_0},\dots,\psi_{t-1})$  of high-level states 
such that  
$\Pr( \Psivect_{T(t-1)} = \psivect_{T(t-1)}) > 0$, 
$\type^{\psi_{t-2}} = \filling$, 
$\type^{\psi_{t-1}} = \refilling$, and $j = j^{\psi_{t-1}}$.
Conditioned on $\Psivect_{T(t-1)} = \psivect_{T(t-1)} $, with 
probability at least $1-\exp(-\lambda j(\log j )^2/202)$,  
the following events occur.
\begin{enumerate}[(i)]
\item $\exit(t-1) \le t-1+j^{\phi+74}$,
\item for all $\hat{t}$ 
in the range $t-1 \leq \hat{t} < \exit(t-1)$  with $\type^{\Psi(t)} = \stabilising$, $\trans(\hat{t}) \le \hat{t} + j^{\phi+72}$,
\item for all $\hat{t}$ 
in the range $t-1 \leq \hat{t} < \exit(t-1)$ 
with $\type^{\Psi(t)} = \refilling$, $\trans(\hat{t}) \le \hat{t} + j^{\phi+46}$, and
\item $|\{\Psi(t') \colon t-1 \le t' \le \exit(t-1)-1\}| \le 2j$.
\end{enumerate}
\end{lemma}

\begin{proof}  
Fix $t$ and 
$\psivect_{(t-1)}$ as in the lemma statement.
Suppose that
$\Psivect_{T(t-1)} = \psivect_{T(t-1)}$.

Let $\tilde{\tau}_0 = t-1$ and for all positive integers~$x$, let $\tilde\tau_x = \min \{ \exit(t-1), \back(\tilde\tau_{x-1})\}$.
If  $\tilde\tau_x< \exit(t-1) $ then
$\tilde \tau_x = \back(\tilde\tau_{x-1})$ so 
$\tau^{\Psi(\tilde \tau_x)} = \tilde\tau_x$ and $\type^{\Psi(\tilde \tau_x)} = \refilling$.

By
the backoff bounding rule
(Definition~\ref{def:backoff-bounding-rule}), the transition from~$\psi_{t-2}$ to~$\psi_{t-1}$ is via (R\ref{item:R2})(ii) so $j^{\psi_{t-2}} = j$ and $\tau^{\psi_{t-1}} = t-1$.
Applying Lemma~\ref{lem:vol-F-lb} to the $j$-filling state $\psi_{t-2}$,
we get that for every $k\in [j-1]$, 
$F_k^{\psi_{t-2}}(t-1) \geq \mu_k^{\gammavect}$.
By transition (R\ref{item:R2})(ii) of the backoff bounding rule, 
for all $k\in [j-1]$ with $W_k\geq j^2$, 
$z_k^{\psi_{t-1}} = \min\{ F_k^{\psi_{t-2}}(t-1),\mu_k^{\gammavect}\} = \mu_k^{\gammavect}$.  
Thus, \[\calM_{[j-1]}^{\psi_{t-1}}(t-1) \leq \sum_{k\in [j-1]: W_k \leq j^2} \mu_k^{\gammavect} \leq \sum_{k\in [j-1]: W_k \leq j^2} \lambda W_k \leq \lambda j^3.\]

Using the notation from the start of the proof, we have shown that $\calM_{[j-1]}^{\Psi(\tilde\tau_0)}(\tilde\tau_0) \leq \lambda j^3$.
If $x\geq 1$ and $\tilde{\tau}_x < \exit(t-1)$
then Lemma~\ref{lem:bound-z-shrink-one-loop} 
applied with 
the $t'$ of the lemma as~$\tilde\tau_x$ and the $t$ of the lemma as  
$\tilde\tau_{x-1}$
guarantees that 
$
\calM_{[j-1]}^{\Psi(\tilde\tau_x)}(\tilde\tau_x) 
\le 
\calM_{[j-1]}^{ \Psi(\tilde\tau_{x-1})}(\tilde\tau_{x-1})
+ 
\lambda j^3$
so 
\begin{equation}\label{eq:precond}
\calM_{[j-1]}^{\Psi(\tilde\tau_x)}(\tilde\tau_x) 
\le 
\calM_{[j-1]}^{\Psi(\tilde\tau_{0})}(\tilde\tau_{0})
+ 
\lambda x j^3
\leq  (x+1) \lambda j^3.\end{equation}

For all $x \in \{0,\dots,j-1\}$, let $H(x)$ be the set of all pairs $(t_x,\psivect^x)$ where $t_x \ge t-1$ is an integer and $\psivect^x$ is a vector $(\psivect_{t-1}^x,\dots,\psivect_{t_x-1}^x)$ such~that:
\begin{itemize}
\item $\psivect^x_{T(t-1)} = \psivect_{T(t-1)}$;
\item $\pr(\Psivect_{T(t_x-1)} = \psivect^x) > 0$; and 
\item $\Psivect_{T(t_x-1)} = \psivect^x$ entails $\tilde\tau_x = t_x-1 < \exit(t-1)$. 
\end{itemize}
(Thus $(t_x, \psivect^x) \in H(x)$ means that $t_x$ is a possible value of $\tilde\tau_x+1$ subject to $\tilde\tau_x < \exit(t-1)$, and that $\psivect^x$ is a possible history of high-level states from $t-1$ up to $t_x-1$ with this property.)
For any $x \in \{0,\dots,j-1\}$ and $(t_x,\psivect^x) \in H(x)$
we can apply Lemma~\ref{lem:tautilde-unlikely-back} with the~$t$ of the lemma as~$t_{x}$,
conditioning on 
$\Psivect_{T(t_{x}-1)} = \psivect^{x}$.
 By Equation~\eqref{eq:precond}, since $x\leq j$ and $j\geq j_0$ is sufficiently large, the precondition of
Lemma~\ref{lem:tautilde-unlikely-back} is satisfied.
The lemma shows that 
\begin{equation} \label{eq:firstfact}
\Pr( \tilde\tau_{x+1} \neq \exit(t-1) \mid
\Psivect_{T(t_{x}-1)} = \psivect^{x}  )
 \leq 
  \exp(-\lambda (\log j )^2/202). 
\end{equation}

We can also apply Lemma~\ref{lem:nonback} with the $t$ of  the lemma as $t_x = \tilde\tau_{x}+1$. 
Part~(iv) shows that 
$\Psivect_{T(t_{x}-1)} = \psivect^{x}$ 
 implies
\begin{equation}
\label{eq:secondfact}
\tilde\tau_{x+1} \leq \tilde\tau_{x}+ j^{\phi+73}.
\end{equation}
Parts (ii) and (iii), together with \eqref{eq:precond} show that 
$\Psivect_{T(t_{x}-1)} = \psivect^{x}$ 
 implies
\begin{equation}
\label{eq:thirdfact}
\mbox{for all $\hat{t}$ in the range $\tilde\tau_x \leq \hat{t} < \tilde\tau_{x+1}$, 
$\trans(\hat{t}) \leq \hat{t} + j^{\phi+72}$}
\end{equation}
and 
\begin{equation}\label{eq:fourthfact}
    \mbox{for all $\hat{t}$ in the range $\tilde\tau_x \leq \hat{t} < \tilde\tau_{x+1}$ with $\type^{\psi_{\hat{t}}} = \refilling$, 
$\trans(\hat{t}) \leq \hat{t} + j^{\phi+46}$.}
\end{equation}

We will next repeatedly apply~\eqref{eq:firstfact}
and~\eqref{eq:secondfact}. 
The first and third equalities in the following calculation are from the definition of $\tilde\tau_{x+1} = \min\{\exit(t-1),\back(\tilde{\tau}_x)\}$ together with the fact that the sequence $\back(\tilde{\tau}_0),\back(\tilde\tau_1),\dots$ is non-decreasing. The second inequality is from~\eqref{eq:firstfact}.
\begin{align*}
&\Pr( \tilde\tau_j \neq \exit(t-1) \mid \Psivect_{T(t-1)} = \psivect_{T(t-1)} )\\
&\qquad= \Pr\big( \bigwedge_{x=0}^{j-1} \tilde\tau_{x+1} = \back(\tilde\tau_x) < \exit(t-1) \mid
\Psivect_{T(t-1)} = \psivect_{T(t-1)}\big)\\
&\qquad= \prod_{x=0}^{j-1}\pr\big(\tilde\tau_{x+1} = \back(\tilde\tau_x) < \exit(t-1) \mid \Psivect_{T(t-1)} = \psivect_{T(t-1)} \wedge \bigwedge_{y=0}^{x-1} (\tilde\tau_{y+1} = \back(\tilde\tau_y) < \exit(t-1)) \big)\\
&\qquad= \prod_{x=0}^{j-1}\pr\big(\tilde\tau_{x+1} \ne \exit(t-1) \mid \Psivect_{T(t-1)} = \psivect_{T(t-1)} \wedge \tilde\tau_{x} < \exit(t-1)) \big)\\
&\qquad\le \prod_{x=0}^{j-1} 
\max_{(t_x,\psivect^x) \in H(x)} \pr\big(\tilde\tau_{x+1} \ne \exit(t-1) \mid \Psivect_{T(t_x-1)}=\psi^x\big)\\
&\qquad\leq\prod_{x=0}^{j-1} \exp(-\lambda (\log j )^2/202) \leq \exp(-\lambda j(\log j )^2/202).
\end{align*}

To finish, we show that items (i)--(iv) in the lemma statement follow from 
$\Psivect_{T(t-1)} = \psivect_{T(t-1)}$ and
$\tilde \tau_j = \exit(t-1)$. By applying~\eqref{eq:secondfact}  for all $x\in  [j]$,
we find that  
$\Psivect_{T(t-1)} = \psivect_{T(t-1)}$ 
 implies
$\tilde\tau_{j} \leq t-1 + j^{\phi+74}$. giving (i).  
It also implies (iv) 
since 
$|\{\Psi(t') : t-1 \leq t' \leq \tilde\tau_j - 1\}| \leq 2j$ from the backoff bounding rule and the definition of~$\tilde\tau_j$.
Finally, items (ii) and (iii) follow by applying~\eqref{eq:thirdfact} and~\eqref{eq:fourthfact}
for all $x\in \{0,\ldots,j-1\}$.
\end{proof}

Lemma~\ref{lem:F-step-in-loop} gives an upper bound on how much the expected bin population of high-weight bins that are not over-full can decrease in certain high-level state transitions.
Lemma~\ref{lem:loop-F-bound} then uses the bound from Lemma~\ref{lem:F-step-in-loop} to upper bound how much these expectations can decrease over the course of the loop.

\begin{lemma}\label{lem:F-step-in-loop}
Fix a real number $\lambda \in (0,1)$ and a send sequence $\p$ with $p_0=1$. Let $j_0$ be a
sufficiently large integer.  
Let  
$\Psi_0$ be the initialising state with
$j^{\Psi_0} = j_0$ and $\tau^{\Psi_0} = \tau_0$.  
Let $\hls$ be  
the $\Psi_0$-backoff bounding rule from Definition~\ref{def:backoff-bounding-rule}.
Let $V$ be a volume process from~$\Psi_0$ with transition rule~$\hls$, send sequence~$\p$ and birth rate~$\lambda$. 
Fix any $t> \tau_0$ and any
$\psi_{t-1} \in \OurHLS$ with $j=j^{\psi_{t-1}}$.
Let $\psi_t = \hls(\psi_{t-1},\bbar, t)$ for some $\bbar$.
Suppose that the transition $\psi_{t-1}\rightarrow \psi_t$ is via (R\ref{item:R2})(ii), (R\ref{item:R5}), or (R\ref{item:R6}).
Then for all $k \in \Upsilon_{j+1,\ge j^2}$,
$\min \{F_k^{\psi_t}(t+1),\mu_k^{\gammavect}\} \geq \min \{F_k^{\psi_{t-1}}(t),\mu_k^{\gammavect}\} - 2\lambda$.
\end{lemma}
\begin{proof}
 
Fix $k\in \Upsilon_{j+1,\geq j^2}$.
The lemma follows by adding the inequalities~\eqref{eq:AAA}
and \eqref{eq:BBB}. The rest of the proof establishes these two inequalities.
\begin{align}
\label{eq:AAA} 
\min\{F_k^{\psi_t}(t+1),\mu_k^{\gammavect}\} - \min\{F_k^{\psi_{t}}(t),\mu_k^{\gammavect}\} &\ge -2\lambda.\\\label{eq:BBB}
\min\{F_k^{\psi_t}(t),\mu_k^{\gammavect}\} - \min\{F_k^{\psi_{t-1}}(t),\mu_k^{\gammavect}\} &\ge 0.
\end{align}

It is clear from the backoff bounding rule that $j^{\psi_t} = j$.
Let $\psi_t = (g',\tau', j,\zvect', \calS', \type')$.
Recall from Definition~\ref{def:F} that $F_k^{\psi_t}(t+1) = f_k^{\zerovect,\zvect'}(t+1-\tau')$,
so since $t+1 > \tau'$, 
by Definition~\ref{def:f},
\begin{equation}\label{eq:pbound}
F_k^{\psi_t}(t+1)=
f_k^{\zerovect,\zvect'}(t+1-\tau')  
\ge (1-p_k)f_k^{\zerovect,\zvect'}(t-\tau') = (1-p_k)F_k^{\psi_t}(t).
\end{equation}

To prove~\eqref{eq:AAA}, we split into three cases depending on $F_k^{\psi_t}(t)$.

\medskip\noindent\textbf{Case 1:} $F_k^{\psi_t}(t) \ge 2\lambda W_k$ and $p_k < 1/2$. By Observation~\ref{obs:mu-gamma-bounds},  
$\mu_k^{\gammavect} \leq \lambda W_k$, so  
$F_k^{\psi_t}(t) \geq 2 \lambda W_k \ge \mu_k^{\gammavect}$ and, using~\eqref{eq:pbound}, $F_k^{\psi_t}(t+1) \ge F_k^{\psi_t}(t)/2 \ge \lambda W_k \ge \mu_k^{\gammavect}$. Thus $\min\{F_k^{\psi_t}(t+1),\mu_k^{\gammavect}\} = \min\{F_k^{\psi_t}(t),\mu_k^{\gammavect}\} = \mu_k^{\gammavect}$ and~\eqref{eq:AAA} follows immediately.

\medskip\noindent\textbf{Case 2:} $F_k^{\psi_t}(t) \ge 2\lambda W_k$ and $p_k \ge 1/2$. Then $\min\{F_k^{\psi_t}(t+1),\mu_k^{\gammavect}\} - \min\{F_k^{\psi_t}(t),\mu_k^{\gammavect}\} \ge -\mu_k^{\gammavect} \ge -\lambda W_k = -\lambda/p_k \ge -2\lambda$, as required by~\eqref{eq:AAA}.

\medskip\noindent\textbf{Case 3:} $F_k^{\psi_t}(t) < 2\lambda W_k$. 
Then from~\eqref{eq:pbound},
\[
F_k^{\psi_t}(t+1) -  F_k^{\psi_t}(t) \ge -p_kF_k^{\psi_t}(t) > -2p_k\lambda W_k = -2\lambda,
\]
so~\eqref{eq:AAA} follows by Observation~\ref{obs:dumb}.
In all three cases, we have proved~\eqref{eq:AAA} as required.

We now prove~\eqref{eq:BBB},
namely that 
$\min\{F_k^{\psi_t}(t),\mu_k^{\gammavect}\} \ge \min\{F_k^{\psi_{t-1}}(t),\mu_k^{\gammavect}\}$.
We will check each of the relevant transitions (R\ref{item:R2})(ii), (R\ref{item:R5})  and (R\ref{item:R6}). If $\psi_t = \psi_{t-1}$, the statement is obvious --  this holds in (R\ref{item:R5})(ii) and (R\ref{item:R6})(ii). Otherwise, $\min\{F_k^{\psi_t}(t),\mu_k^{\gammavect}\} = \min\{z_k',\mu_k^{\gammavect}\}$
so we wish to show
$\min\{z_k',\mu_k^{\gammavect}\} \ge \min\{F_k^{\psi_{t-1}}(t),\mu_k^{\gammavect}\}$.

We consider the relevant transitions, recalling that 
$k\in \Upsilon_{j+1,\geq j^2}$.
\begin{itemize}
\item (R\ref{item:R2})(ii):
$z_k' =\min\{F_k^{\psi_{t-1}}(t),\mu_k^{\gammavect}\}$.

\item (R\ref{item:R5})(i) to stabilising: For $k\in [j-1]$, $z_k' =\mu_k^{\gammavect}$, $z_j' = F_j^{\psi_{t-1}}(t)$.
\item (R\ref{item:R5})(i) to refilling: $z_k' = F_k^{\psi_{t-1}}(t)$.
\item (R\ref{item:R6})(i) to filling: For $k\in [j-1]$, $z_k' =\mu_k^{\gammavect}$, $z_j' = F_j^{\psi_{t-1}}(t)$.
\item (R\ref{item:R6})(i) to refilling: $z_k' = F_k^{\psi_{t-1}}(t)$.
\end{itemize}

In  all cases,
$z'_k   \geq \min\{F_k^{\psi_{t-1}}(t),\mu_k^{\gammavect}\}$, so $ 
\min\{z_k', \mu_k^{\gammavect}\} \ge \min\{F_k^{\psi_{t-1}}(t),\mu_k^{\gammavect}\}$. 
Hence 
$
\min\{F_k^{\psi_t}(t),\mu_k^{\gammavect}\} = 
\min\{z_k', \mu_k^{\gammavect}\} \ge \min\{F_k^{\psi_{t-1}}(t),\mu_k^{\gammavect}\}$, as required.
\end{proof}

\begin{lemma}
\label{lem:loop-F-bound}
Fix a real number $\lambda \in (0,1)$ and a send sequence $\p$ with $p_0=1$. Let $j_0$ be a
sufficiently large integer.  
Let  
$\Psi_0$ be the initialising state with
$j^{\Psi_0} = j_0$ and $\tau^{\Psi_0} = \tau_0$.  
Let $\hls$ be  
the $\Psi_0$-backoff bounding rule from Definition~\ref{def:backoff-bounding-rule}.
Let $V$ be a volume process from~$\Psi_0$ with transition rule~$\hls$, send sequence~$\p$ and birth rate~$\lambda$. For any
$t \ge \tau_0$, 
let   $T(t) = \{\tau_0,\ldots,t\}$.
Then  the following statement holds  for
every $t \geq  \tau_0+2$ and
every tuple
$ 
\psivect_{T(t-1)} = (\psi_{\tau_0},\dots,\psi_{t-1})$  of high-level states 
such that  
$\Pr( \Psivect_{T(t-1)} = \psivect_{T(t-1)}) > 0$, 
$\type^{\psi_{t-2}} = \filling$, 
$\type^{\psi_{t-1}} = \refilling$, and $j = j^{\psi_{t-1}}$.
Conditioned on $\Psivect_{T(t-1)} = \psivect_{T(t-1)} $,
for all $k\in \Upsilon_{j+1,\geq j^2}$
and all 
$t'$ in the range 
$t-2 \leq t' < \exit(t-1)$, 
$F_k^{\psi_{t'}}(t'+1) \geq \min\{F_k^{\psi_{t-2}}(t-1),\mu_k^{\gammavect}\} - 2 \lambda(t'-t+2)$.
\end{lemma}
\begin{proof} 
Consider any $t'$ in the given range. Then
\[
\min\{F_k^{\psi_{t'}}(t'+1),\mu_k^{\gammavect}\} - \min\{F_k^{\psi_{t-2}}(t-1),\mu_k^{\gammavect}\} = \sum_{t''=t-1}^{t'} \big(\min\{F_k^{\psi_{t''}}(t''+1),\mu_k^{\gammavect}\} - \min\{F_k^{\psi_{t''-1}}(t''),\mu_k^{\gammavect}\}\big).
\]

We will apply Lemma~\ref{lem:F-step-in-loop} to show that each term
$\min\{F_k^{\psi_{t''}}(t''+1),\mu_k^{\gammavect}\} - \min\{F_k^{\psi_{t''-1}}(t''),\mu_k^{\gammavect}\})$ is at least $-2\lambda$. The $t$ in Lemma~\ref{lem:F-step-in-loop} is our $t''$ so the application requires 
$t''  > \tau_0$, which is satisfied, and that the transition $\psi_{t''-1} \rightarrow \psi_{t''}$ is via (R\ref{item:R2})(ii), (R\ref{item:R5}), or (R\ref{item:R6}). This follows because 
any other transition would lead to
$\type^{\psi_{t''}} \not\in \{\refilling,\stabilising\}$ which would imply $\exit(t-1) \leq t''$, and we have $t'' \le t'  < \exit(t-1)$. The result therefore follows.
\end{proof}

We are finally ready to prove Lemma~\ref{lem:main-loop-bounds}, which says  that, starting from a transition from a $j$-filling high-level state~$\Psi$
to a $j$-refilling high-level state, the volume process is likely to (fairly) quickly escape from the refilling-stabilising loop, transitioning back to a $j$-filling
high level state that has not lost many balls from bin~$j$ (in expectation), relative to~$\Psi$.

\begin{lemma}\label{lem:main-loop-bounds}
Fix a real number $\lambda \in (0,1)$ and a send sequence $\p$ with $p_0=1$. Let $j_0$ be a
sufficiently large integer.  
Let  
$\Psi_0$ be the initialising state with
$j^{\Psi_0} = j_0$ and $\tau^{\Psi_0} = \tau_0$.  
Let $\hls$ be  
the $\Psi_0$-backoff bounding rule from Definition~\ref{def:backoff-bounding-rule}.
Let $V$ be a volume process from~$\Psi_0$ with transition rule~$\hls$, send sequence~$\p$ and birth rate~$\lambda$. For any
$t \ge \tau_0$, 
let   $T(t) = \{\tau_0,\ldots,t\}$.
Then  the following statement holds  for
every $t \geq  \tau_0+2$ and
every tuple
$ 
\psivect_{T(t-1)} = (\psi_{\tau_0},\dots,\psi_{t-1})$  of high-level states 
such that  
$\Pr( \Psivect_{T(t-1)} = \psivect_{T(t-1)}) > 0$, 
$\type^{\psi_{t-2}} = \filling$, 
$\type^{\psi_{t-1}} = \refilling$, and $j = j^{\psi_{t-1}}$.
Conditioned on $\Psivect_{T(t-1)} = \psivect_{T(t-1)} $,
with 
probability at least $1-\exp(-\lambda j(\log j )^2/203)$,  
the following events occur. 
\begin{enumerate}[(i)]
\item $\exit(t-1) \leq t-1+ j^{\phi+74}$.
\item The transition 
into $\Psi(\exit(t-1))$ is via (R\ref{item:R6})(i) 
into a $j$-filling state. 
\item $F_j^{\Psi(\exit(t-1))}(\exit(t-1)) \geq \min\{F_j^{\psi_{t-2}}(t-1),\mu_j^{\gammavect}\} - 2\lambda j^{\phi+74}$.
\end{enumerate}
\end{lemma}
\begin{proof}
Let $\calE_1$ be the event that $\exit(t-1) \leq t -1 + j^{\phi+74}$.
Let \[A = \{t' 
\mid t -1\leq t' \leq  
\min \{ \exit(t-1), t-1+j^{\phi+74}\}\}.\]
For $t'\in A$, 
let $\calE_2(t')$ be the event
that the transition 
into $\Psi(t')$   is not via (R\ref{item:R1})(iii).  
 
We will first prove $\pr(\neg \calE_1 \cup \cup_{t'\in A} \neg\calE_2(t') \mid 
\Psivect_{T(t-1)} = \psivect_{T(t-1)} ) \le 
\exp(-\lambda j (\log j )^2/204)$, then prove that
if $\calE_1 \cap \cap_{t'\in A} \calE_2(t')$ holds, then so do events
(i)--(iii).

We start with the failure probability. By Lemma~\ref{lem:leave-loop-quickly}, \[\pr(\neg \calE_1 \mid \Psivect_{T(t-1)} = \psivect_{T(t-1)}
 ) \leq 
 \exp(-\lambda j (\log j )^2/202).\]  

Consider $t'\in A$. Since $t' \le \exit(t-1)$, $j^{\psi_{t'}} = j$. By Lemma~\ref{lem:volume-easyChernoffs} and the law of total probability, 
\[ \pr(\Psi(t')\rightarrow\Psi(t'+1) \mbox{is via (R\ref{item:R1})(iii)} \mid  \Psivect_{T(t-1)} = \psivect_{T(t-1)})   \le \exp(-j^{1.9})\]
Since $j\geq j_0$ is sufficiently large, the failure probability bound follows by a union bound  
over the events $\neg \calE_1$ and, summed over all $t'\in A$, $\neg \calE_2(t')$.

It remains to prove that 
if $\calE_1 \cap \cap_{t'\in A} \calE_2(t')$ holds, then so do events
(i)--(iii).
Event (i) is immediate from $\calE_1$. For (ii), note 
from the backoff bounding rule that the only possible transitions 
into $\Psi(\exit(t-1))$ are via  (R\ref{item:R1})(ii), (R\ref{item:R1})(iii), (R\ref{item:R2})(i) from a $j$-refilling state, and (R\ref{item:R6})(i) 
to a $j$-filling state.  Since $\calE_2(t')$ holds for all $t'\in A$, and by $\calE_1$, $\exit(t-1) \in A$,
$\calE_2(\exit(t-1))$ guarantees that  
 the transition into $\Psi(\exit(t-1))$ is not via (R\ref{item:R1})(iii).
To demonstrate that event (ii) holds, we will rule out 
(R\ref{item:R1})(ii) and (R\ref{item:R2})(i) from a $j$-refilling state. Observation~\ref{obs:no-R2-fromrefilling} rules out 
(R\ref{item:R2})(i). 
So for proving (ii) it suffices to show that the transition is not via (R\ref{item:R1})(ii). 
Lemma~\ref{lem:vol-F-lb} 
guarantees that for all   $k\in [j-1]$,
$F_k^{\psi_{t-2}}(t-1) \geq \mu_k^{\gammavect}$.

Thus by Lemma~\ref{lem:loop-F-bound}, for all $k \in \Upsilon_{j, \geq j^\chi}$ 
and all $t'$ in the range $t-2 \leq t' < \exit(t-1)$,
$F_k^{\psi_{t'}}(t'+1) \geq 
\mu_k^{\gammavect} - 2\lambda (t'-t+2) \geq
\mu_k^{\gammavect} - 2\lambda (\exit(t-1)-t+2) 
\geq \mu_k^{\gammavect} - 2\lambda (j^{\phi+74} +1) 
$.
To rule out (R\ref{item:R1})(ii) we will show that $ 
\mu_k^{\gammavect} - 2\lambda (j^{\phi+74} +1) 
\geq \lambda W_k/2$. 
By Observation~\ref{obs:mu-gamma-bounds}, $\mu_k^{\gammavect} \ge 2\lambda W_k/3$.
So this holds since $W_k \geq j^{\chi}\geq j^{\phi+75}$ and $j\geq j_0$ is sufficiently large.

It remains to show (iii), namely that
$F_j^{\Psi(\exit(t-1))}(\exit(t-1)) \geq \min\{F_j^{\psi_{t-2}}(t-1),\mu_j^{\gammavect}\} - 2\lambda j^{\phi+74}$
By Lemma~\ref{lem:loop-F-bound}, $F_j^{\Psi(\exit(t-1)-1)}(\exit(t-1)) \ge \min\{F_j^{\psi_{t-2}}(t-1),\mu_j^{\gammavect}\} - 2\lambda(\exit(t-1)-(t-1))$. 
Since $\calE_1$ occurs,  
this is at least $
\min\{F_j^{\psi_{t-2}}(t-1),\mu_j^{\gammavect}\} - 2\lambda j^{\phi+74}$.

Since (ii) holds, the transition 
into $\Psi(\exit(t-1))$ is 
via (R\ref{item:R6})(i) into a $j$-filling state, so from the backoff bounding rule
$ z_j^{\Psi(\exit(t-1))} = F_j^{\Psi(\exit(t-1)-1)}(\exit(t-1))$.

Also, by the definition of~$F$,
$F_j^{\Psi(\exit(t-1))}(\exit(t-1)) = z_j^{\Psi(\exit(t-1))}$. Thus
\[
F_j^{\Psi(\exit(t-1))}(\exit(t-1)) = F_j^{\Psi(\exit(t-1)-1)}(\exit(t-1)) \ge \min\{F_j^{\psi_{t-2}}(t-1),\mu_j^{\gammavect}\} - 2\lambda j^{\phi+74},
\]
as required.
\end{proof}

\subsection{Likely progress towards increasing \texorpdfstring{$j$}{j}}
\label{sec:volume-progress}

In this section we show that the volume process is likely to progress~$j$ sufficiently quickly. The main goal is Lemma~\ref{lem:VolumeAnalysis},
which shows that, with probability at least~$99/100$, the volume process never enters a failure state, and for all $j\geq j_0$, it quickly gets to a $j$-filling or $j$-advancing state. 
The main parts of the proof are Lemmas~\ref{lem:OnejVolumeAnalysis}  
and~\ref{lem:OnejVolumeAnalysisAdvancing}  
which show  this for a single bin -- namely
that from a $j$-filling high-level state (or a $j$-advancing one in the case of Lemma~\ref{lem:OnejVolumeAnalysisAdvancing}) the process is very likely to reach a $(j+1)$-advancing or $(j+1)$-filling state before failing, and one of these is likely to happen in $\poly(j) W_j$ steps.

Lemma~\ref{lem:OnejVolumeAnalysis} follows the progress of the volume process from a $j$-filling state. In order to facilitate the proof, which is based on quite a few lemmas, Definition~\ref{def:fromfilling} 
captures the initial setting (so we don't have to keep repeating it) and
defines various sets and stopping times. Some explanation of these stopping times  follows the definition.

\begin{definition}\label{def:fromfilling}
The initial setting is defined as follows:
Fix $\lambda \in (0,1)$ and a send sequence $\p$ with $p_0=1$ such that,
for all $i \ge 1$, $W_i \le (4/\lambda)^i$.
Let $j_0$ be a sufficiently large positive integer
where all strongly-exposed bins are in $[j_0-1]$.
Fix an initialising high-level state $\Psi_0$  with  $\tau^{\Psi_0} = \tau_0$ and $j^{\Psi_0} = j_0$. 
Let $\hls$ be the $\Psi_0$-backoff bounding rule from Definition~\ref{def:backoff-bounding-rule}.  
Let $V$ be a volume process from $\Psi_0$ with transition rule $\hls$, 
 send sequence~$\p$, and birth rate~$\lambda$  (Definition~\ref{def:volume-process}).
For any
$t \ge \tau_0$, 
let   $T(t) = \{\tau_0,\ldots,t\}$.
Fix $t>\tau_0$, 
 and fix a tuple  
$ 
\psivect_{T(t-1)} = (\psi_{\tau_0},\dots,\psi_{t-1})$  of high-level states 
such that  
$\Pr( \Psivect_{T(t-1)} = \psivect_{T(t-1)}) > 0$, 
$\type^{\psi_{t-1}} = \filling$, 
$\tau^{\psi_{t-1}} = t-1$, and $z_j^{\psi_{t-1}} = 0$ where $j=j^{\psi_{t-1}}$.
Consider the evolution of~$V$, conditioned on
$\Psivect_{T(t-1)} = \psivect_{T(t-1)}$.

Define the following quantities, sets, and stopping times.
\[
\tau_{j,\full} = \min\{t'  \geq t \colon \type^{\Psi(t')} = \Failure \mbox{ or }j^{\Psi(t')} = j+1\}.\]
Let $\Out_0 = t-1$. For all positive integers $x$, let
\begin{align*}
\In_x &= \min(\{\tau_{j,\full}\} \cup \{t > \Out_{x-1} \colon 
\mbox{$\Psi(t-1) \rightarrow \Psi(t)$
is via (R\ref{item:R2})(ii) from $j$-filling to $j$-refilling
} \}).\\
\Out_x &= \min(\{\tau_{j,\full}\} \cup \{t > \In_{x} \colon 
\mbox{$\Psi(t-1) \rightarrow \Psi(t)$
is via  (R\ref{item:R6})(i) from $j$-stabilising to $j$-filling
} \}).
\end{align*}
Let 
\begin{align*} 
A &= \lceil 4(4 j)^{\phi} W_j\rceil.\\
\Ready &  
= \{ T'\geq t-1 \colon  
\mbox{ $\Psi(T')$ is a $j$-filling state,  
$T' + j^{24} \leq t+A$, and
$
F_j^{\Psi(T')}(T'+1) \ge \mu_j^{\gammavect} $ }
\}. 
\\
\tmin &= \min(\{t+A\} \cup \Ready).
\end{align*} 
\end{definition}

\begin{remark}\label{rem:fromfilling}
The setting of Definition~\ref{def:fromfilling}
starts at time~$t-1$ when the volume process enters a $j$-filling state.
The stopping times apply from then until time~$\tau_{j,\full}$ when the process has either reached a failure state, or has reached a $(j+1)$-filling state or $(j+1)$-advancing state.
Before $\tau_{j,\full}$ the process proceeds in a very structured way (from the backoff bounding rule). We view time $t-1$ as the $0$'th time that the process exits the $j$-refilling/$j$-stabilising loop (into a $j$-filling state). After the process exits 
the loop for the $x$'th time (unless $\tau_{j,\full}$ happens), it stays in the same $j$-filling state (via transition (R\ref{item:R4})(ii)) for some number of steps and then, at time $\In_x$, it goes via (R\ref{item:R2})(ii) into the $j$-filling/$j$-stabilising loop for the $x$'th time. Unless $\tau_{j,\full}$ happens, it stays in this loop until time $\Out_x$, when it goes out via (R\ref{item:R6})(i) into a $j$-filling state.

Thus, if the
volume process
goes into the refilling-stabilising loop
at least $x$ times while filling  bin~$j$, then
$\In_x$ is the  $x$'th time that it goes in. $\Out_x$ is the $x$'th time that it comes out of the loop. However, these are capped at $\tau_{j,\full}$, since from time~$\tau_{j,\full}$, the volume process stops filling bin~$j$ -- it either moves on to bin~$j+1$ or it fails.

We will use the fact that if $\In_x < \tau_{j,\full}$ then $\Out_x = \exit(\In_x)$, where $\exit(\In_x)$ is the stopping time defined in 
 Section~\ref{sec:volume-theloop} -- it is the first time after~$\In_x$ that the high-level state is not refilling or stabilising.  

 For every non-negative integer~$x$,
 $t-1 \leq \Out_{x} \leq \In_{x+1} \leq \Out_{x+1} \leq \tau_{j,\full}$. Also, if $\Out_{x+1} < \tau_{j,\full}$ then these inequalities are strict so 
 $\Out_{x} < \In_{x+1} < \Out_{x+1}$. In fact, if $\Out_{x}< \tau_{j,\full}$ then
 the definition of $\In_{x+1}$ gives 
 $\Out_x < \In_{x+1}$. Similarly, if $\In_{x+1} < \tau_{j,\full}$ then $\In_{x+1} < \Out_{x+1}$.

The exact definition of~$A$ is not very important but it is important that it is a large enough polynomial times $W_j$. This is adequate time for the process to be likely to fill bin~$j$. Thus in Lemma~\ref{lem:OnejVolumeAnalysis} we will show that it is likely that 
$\tau_{j,\full} \leq t+A$.

$\Ready$ is the set of steps~$T'$ such that the transition from $\Psi(T')$ to $\Psi(T'+1)$ could potentially be via (R\ref{item:R4})(i), which would trigger~$\tau_{j,\full}$, and hopefully move to a $(j+1)$-advancing or $(j+1)$-filling state. The main point is that the value $F_j^{\Psi(T')}(T'+1)$ is at least $\mu_j^{\gammavect}$, which is required by the guard on (R\ref{item:R4})(i) in the backoff bounding rule.
In fact, the backoff bounding rule also requires $T'+1 \geq \tau^{\Psi(T')} + j^{24}$ which may not be true of~$T'$. 
Nevertheless, the value of $F_j^{\Psi(T')}(T'+a)$ does not decrease
below $\mu_j^{\gammavect}$ as $a$ increases -- we will prove this  in Lemma~\ref{lem:filling-F-monotonoe} -- so if the transitions from $\Psi(T'),\Psi(T'+1),\ldots,\Psi(T'+a-1)$ are via
(R\ref{item:R4})(ii)
and $T' + a \geq \tau^{\Psi(T')} + j^{24}$, the process may then transition via (R\ref{item:R4})(i). 
For this $a \leq  j^{24}$ suffices, and
the potential that $a$ might be this large is why we demand $T' + j^{24} \leq t+A$ in the definition of~$\Ready$. 
 
The stopping time $\tmin$ is the first 
time~$T' \in \Ready$, or it is $t+A$ 
if $\Ready=\emptyset$ (this will be the low-probability  situation where we assume that the volume process may fail to reach~$j+1$).
\end{remark}

Lemma~\ref{lem:filling-F-monotonoe}
is the monotonicity lemma discussed in Remark~\ref{rem:fromfilling}.

\begin{lemma}\label{lem:filling-F-monotonoe}
Consider the setting of Definition~\ref{def:fromfilling}. 
Suppose that $t'\in \Ready$.
Then for all $a\in \integers_{\geq 1}$,
$F_j^{\Psi(t')}(t'+a) \geq \mu_j^{\gammavect}$.
\end{lemma}
\begin{proof}
Let $\Psi(t') = (g,\tau, j,\zvect, \calS, \filling)$. From the definition of filling state (Definition~\ref{def:filling-state}),
for all $x\in [j-1]$, $z_x = \mu_x^{\gammavect}$. 
Recall from Definition~\ref{def:F} that
for $x\in [j]$ and $t\geq \tau$,
$F_x^{\Psi(t')}(t) =f_x^{\zerovect,\zvect}(t-\tau) $.
Note that $t' \geq \tau$.
Applying  Lemma~\ref{lem:f-fill-bin}(i)  with the $k$ of the lemma as $j-1$, $\Gammavect = \zerovect$,
$\Lambdavect = \gammavect$, and $t_0 = 0$,
for all $x\in [j-1]$,
$f_x^{\zerovect,\zvect}(t'+1-\tau) \geq \mu_x^{\gammavect}$.
Since 
$t' \in \Ready$, 
$f_j^{\zerovect,\zvect}(t'+1-\tau) \geq \mu_j^{\gammavect}$.
Now apply Lemma~\ref{lem:f-fill-bin}(i)
with $k=j$, $\Gammavect = \zerovect$,
$\Lambdavect = \gammavect$, $t_0 = t'+1-\tau$,
and $t=a-1$
to get 
for all $x\in [j]$,
$f_x^{\zerovect,\zvect}(t_0+t) \geq \mu_x^{\gammavect}$, giving the lemma statement. 
\end{proof}

Definition~\ref{def:morefromfilling}
extends the setting of Definition~\ref{def:fromfilling}
by defining the particular value of ``$a$'' that is discussed in Remark~\ref{rem:fromfilling}, and defining some useful sets. The notation in Definitions~\ref{def:fromfilling} and~\ref{def:morefromfilling} is local to this section, so it is omitted  from the index. 

\begin{definition}\label{def:morefromfilling}
We work in the setting of Definition~\ref{def:fromfilling}. Also,
let 
\begin{align*} 
\amin &= \min\{ a\in  \integers_{\geq 1}
: \tau^{\Psi(\tmin)} + j^{24}\leq
\tmin+ a  
\}.\\
\LoopStarts &= \{{t'} : 
t \leq t' \leq 
\min\{
t+A+1,
\tau_{j,\full}\}
-1 \mbox{ and for some $x\geq 1$, $t' = \In_x$} 
 \}.\\
\Fills &= \{ t' : 
t-1 \leq t' < t+A 
\mbox{ and $\Psi(t')$ is a $j$-filling state 
}
 \}.\\
\end{align*}  
\end{definition}

\begin{remark}\label{rem:morefromfilling}
The intention is that the setting of 
Definition~\ref{def:morefromfilling}
includes, and extends, the setting of
Definition~\ref{def:fromfilling}.
Thus, when we consider the setting of 
Definition~\ref{def:morefromfilling} we implicitly use both definitions.

Note that 
$\tau^{\Psi(\tmin)}\leq \tmin$,
so $\amin \leq j^{24}$.
The purpose of $\amin$ is as follows.
Suppose that  $\tmin  \in \Ready$. Suppose   that for all $a\in [\amin-1]$, the
transition into $\Psi(\tmin+a)$ is 
not via (R\ref{item:R1}) or (R\ref{item:R2}). Then it is  via (R\ref{item:R4})(ii).
Also, by Lemma~\ref{lem:filling-F-monotonoe}, as long as 
the transition into~$\Psi(\tmin+\amin)$ 
is not via (R\ref{item:R1}) or (R\ref{item:R2}) it will be via 
(R\ref{item:R4})(i), triggering~$\tau_{j,\full}$. Since $\tmin + j^{24} \leq t+A$,
$\tmin + \amin \leq t+A$.

$\LoopStarts$ is the set of times $t' \leq t+A$ that the
volume process goes into the refilling-stabilising loop   while filling bin~$j$.

$\Fills$ is the 
set of times in the interval $t-1\leq t' < t+A$ 
such that $\Psi(t')$ is $j$-filling.
From the definition of~$\tau_{j,\full}$
and the backoff bounding rule, there are no times~$t'\geq \tau_{j,\full}$ in $\Fills$. From Definition~\ref{def:fromfilling}, for
each $t'\in \Fills$ there is a non-negative integer~$x$ such that 
$\Out_x \leq t' < \In_{x+1}$.
\end{remark}

We have just seen, in  Remark~\ref{rem:morefromfilling},
that for each $t'\in \Fills$, there is a 
positive integer $x$   such that 
$\Out_{x-1} \leq t' < \In_{x}$. Observation~\ref{obs:stay-jfill-state} expands 
on this, capturing the fact that the high-level state does not change in this period. 

\begin{observation}\label{obs:stay-jfill-state}
Consider the setting of Definition~\ref{def:morefromfilling}. For all 
positive integers~$x$  
and all $t'$ such that 
$\Out_{x-1} \leq t' \leq \In_x-1$,
$\Psi(t') = \Psi(\Out_{x-1})$ and this is a $j$-filling state.
\end{observation}
\begin{proof}
Since $\Out_{x-1} < \In_x$, $\Out_{x-1}$ is not $\tau_{j,\full}$ so 
$\Psi(\Out_{x-1})$ is a $j$-filling state reached via (R\ref{item:R6})(i) (or, in the case where $x=1$, it is state $\psi_{t-1}$).
From the backoff bounding rule, 
any sequence of transitions from $\Out_{x-1}$ which repeatedly use (R\ref{item:R4})(ii) keep the high-level state unchanged.
Since $\In_x \le \tau_{j,\full}$, 
the backoff bounding rule ensures that, for all 
$t'$ with $\Out_{x-1} \le t' \leq \In_x-2$, 
the transition $\Psi(t')\rightarrow\Psi(t'+1)$ is
not via (R\ref{item:R1}) (which would trigger $\tau_{j,\full}$), (R\ref{item:R2})(ii) (which would trigger $\In_x$) or (R\ref{item:R4})(i) (which would trigger $\tau_{j,\full}$), so it must be via (R\ref{item:R4})(ii). Thus $\Psi(t'+1) = \Psi(\Out_{x-1})$.   
\end{proof} 

Observation~\ref{obs:readysoon} shows that transition (R\ref{item:R4})(i) is not taken before a time step in~$\Ready$.

\begin{observation}
\label{obs:readysoon}\label{obs:RS}
Consider the setting of Definition~\ref{def:morefromfilling}.  
Consider any $T'$ in the range $t-1 \leq T' \leq t+A-j^{24}$ 
such that 
$[T' ] \cap \Ready = \emptyset$.
Then for any $T''\in \{t-1,\ldots,T'\}$
such that
$\Psi(T'')$ is $j$-filling,
$F_j^{\Psi(T'')}(T''+1) < \mu_j^{\gammavect}$ and the transition from $\Psi(T'')$ is not via 
(R\ref{item:R4})(i).
\end{observation}
\begin{proof}
If $F_j^{\Psi(T'')}(T''+1) \geq \mu_j^{\gammavect}$ then $T''\in \Ready$, contrary to the hypothesis of the observation.
Since $F_j^{\Psi(T'')}(T''+1) < \mu_j^{\gammavect}$ the guard on (R\ref{item:R4})(i) fails in the transition from~$\Psi(T'')$. 
\end{proof}

Definition~\ref{def:fillingevents}
defines some events from  
the setting of these definitions. Once again, these events are local to this section, and are not in the global index.

\begin{definition}\label{def:fillingevents}
Consider the situation of Definition~\ref{def:morefromfilling}. We define the following events.

\begin{itemize}
\item 
$\calE_{{\badRone}}$ is the event that,
for all $t'\in \Fills$, 
the transition $\Psi(t')\rightarrow\Psi(t'+1)$ is not via (R\ref{item:R1}).
\item $\calE_{{\badRtwo}}$ is the event that, for all $\ell \in [j^{24}]$, 
the transition into $\Psi(\tmin + \ell)$ 
is not from a $j$-filling state via (R\ref{item:R2}).
\item 
$\calE_{{\badRfour}}$ is the event that
the transition into $\Psi(\tmin + \amin)$ 
is not 
from a $j$-filling state 
via (R\ref{item:R4})(i) 
into a failure state.
\item 
$\calE_{\badLoopStarts}$ is the event that 
$|\LoopStarts|
\leq  (A+1) \exp(- \lambda(j-1)/800)$.

\item 
$\calE_{\badLongLoop}$ is the event that, for all $x\in [A]$ with $\In_{x} < \tau_{j,\full}$, the following hold: 
$\Out_{x} \leq \In_{x} + j^{\phi+74}$, 
$\Out_{x} < \tau_{j,\full}$,
 and
$F_j^{\Psi(\Out_{x})}(\Out_{x}) \ge \min\{ F_j^{\Psi(\In_{x}-1)}(\In_{x}),\mu_j^{\gammavect}\} - 2\lambda j^{\phi+74}$.

\end{itemize}
 
\end{definition}

Lemmas~\ref{lem:badRone}, \ref{lem:badRtwo}, \ref{lem:badR4},
\ref{lem:badLoopStarts}, and \ref{lem:badLongLoop}  give upper bounds on the probability of
$\neg \calE_{\badRone}$,
$\neg \calE_{\badRtwo}$, 
$\neg \calE_{\badRfour}$, 
$\neg \calE_{\badLoopStarts}$, and 
$\neg \calE_{\badLongLoop}$.
These are combined in Corollary~\ref{cor:fromfillingEs} to upper bound the probability  
of 
$\neg \calE_{\badRone} \vee \neg \calE_{\badRtwo} \vee \neg \calE_{\badRfour} \vee \neg \calE_{\badLoopStarts} \vee \neg \calE_{\badLongLoop}$.
After that we give some lemmas 
which derive some consequences of $
\calE_{\badRone} \wedge \calE_{\badRtwo} \wedge \calE_{\badRfour} \wedge \calE_{\badLoopStarts} \wedge \calE_{\badLongLoop}$.   Finally, the proof of Lemma~\ref{lem:OnejVolumeAnalysis} uses these consequences to show that, from a $j$-filling high-level state, the volume process is very likely to reach a $(j+1)$-advancing or $(j+1)$-filling state before failing, and one of these is likely to happen in $\poly(j)W_j$ steps.

\begin{lemma}\label{lem:badRone}
Consider the setting of Definition~\ref{def:morefromfilling} and the events of Definition~\ref{def:fillingevents}.  Then
\[\Pr(\neg \calE_{{\badRone}}\mid \Psivect_{T(t-1)} = \psivect_{T(t-1)}) \leq  (A+2) \exp(-j^{1.9}).\]
\end{lemma}
\begin{proof}
We wish to upper bound the probability that there is a $t'\in \Fills$ such that the transition 
from $\Psi(t')$ is via 
(R\ref{item:R1}). Since $\Psi(t')$ 
is a $j$-filling state (from the definition of~$\Fills$), it is not via
(R\ref{item:R1})(i) and, by Observation~\ref{obs:no-R1-from-adv-fill-stab}, it is not via (R\ref{item:R1})(ii). 

We will therefore upper bound the probability 
that there is a $t'\in \Fills$ such that the transition 
from $\Psi(t')$ is via 
(R\ref{item:R1})(iii).
From the definition, $\Fills$ 
is the set of times $t'$  in the interval $t-1\leq t' < t+A$ such that $\Psi(t')$ is $j$-filling. By Lemma~\ref{lem:volume-easyChernoffs}, the probability of this transition, conditioned on any state $\Psi(t')$, is at most  $\exp(-j^{1.9})$. The lemma follows by a union bound.
\end{proof}

\begin{lemma}\label{lem:badRtwo}
Consider the setting of Definition~\ref{def:morefromfilling} and the events of Definition~\ref{def:fillingevents}. 
 Then
\[\Pr(\neg \calE_{{\badRtwo}}\mid \Psivect_{T(t-1)} = \psivect_{T(t-1)}) \leq j^{24}
\exp(-\lambda (j-1)/400).\]
\end{lemma}
\begin{proof} 
We wish to upper bound the probability that, for some $\ell \in [j^{24}]$,
the transition into $\Psi(\tmin + \ell)$ 
is from a $j$-filling state via (R\ref{item:R2}). Note, from the backoff bounding rule, that in this case the transition is via (R\ref{item:R2})(ii).
 
Consider any fixed time $t' \geq t-1$ and any vector  $\psivect_{T(t')}$ of high-level states
such that $\Pr(\Psivect_{T(t')} = \psivect_{T(t')}) > 0$ and $\psivect_{t'}$ is a $j$-filling state.
Let $\calE$ be the event that 
the transition $\Psi(t') \rightarrow \Psi(t' + 1)$ is via (R\ref{item:R2})(ii).  By Lemma~\ref{lem:T3tofill},
$\Pr(\calE \mid \Psivect_{T(t')} = \psivect_{T(t')}) \leq 
\exp(-  \lambda  (j-1)/400)$.

We can conclude, for any fixed~$\ell$, 
that the probability that the transition into $\Psi(\tmin + \ell)$ is from a $j$-filling state via (R\ref{item:R2}) is at most $\exp(-  \lambda  (j-1)/400)$. This follows,
since the probability of the given transition,  assuming 
$\Psivect_{T(t-1)} = \psivect_{T(t-1)}$
is at most 
\[ 
\sum_{t' = t-1}^{t+A+\ell-1} 
\Pr(\tmin = t'- \ell+1
\mid \Psivect_{T(t-1)} = \psivect_{T(t-1)}
)
\max_{\psivect_{T(t')}}
\Pr(\calE \mid \Psivect_{T(t')} = \psivect_{T(t')}),
\]
where   the maximum is over 
vectors  $\psivect_{T(t')}$ of high-level states which extend $\psivect_{T(t-1)}$,
have positive probability as a value of 
$\Psivect_{T(t')}$, have $\psi_{t'}$ as a $j$-filling state,
and entail $\tmin = t'-\ell+1$.
The conclusion follows since
$\Pr(\calE \mid \Psivect_{T(t')} = \psivect_{T(t')})$
is uniformly at most $\exp(-\lambda (j-1)/400)$.

Finally, the lemma follows by a union bound over~$\ell$.  
\end{proof}

\begin{lemma}\label{lem:badR4}
Consider the setting of Definition~\ref{def:morefromfilling} and the events of Definition~\ref{def:fillingevents}.  Then
\[\Pr(\neg\calE_{{\badRfour}}\mid \Psivect_{T(t-1)} = \psivect_{T(t-1)}) \leq \exp(- \lambda j/500)
  .\]
\end{lemma}
\begin{proof}

We wish to upper-bound the probability 
that the transition into $\Psi(\tmin + \amin) $ is from a $j$-filling state via 
(R\ref{item:R4})(i) into a failure state.
 
The proof is similar to the proof of Lemma~\ref{lem:badRtwo} except that we use a different lemma, and we don't need the union bound.
Consider any fixed time $t' \geq t-1$ and any vector  $\psivect_{T(t')}$ of high-level states
such that $\Pr(\Psivect_{T(t')} = \psivect_{T(t')}) > 0$ and $\psivect_{t'}$ is a $j$-filling state.
Let $\calE$ be the event that 
the transition $\Psi(t') \rightarrow \Psi(t' + 1)$ is via (R\ref{item:R4})(i) with $\type^{\Psi(t' + 1)} = \Failure$. 
We will show that 
$\Pr(\calE \mid \Psivect_{T(t')} = \psivect_{T(t')}) \leq 
\exp(-  \lambda  j /500)$.

If bin~$j+1$ is many-covered, 
then Lemma~\ref{lem:AdvancingIncj-NoiseFail}
upper-bounds this probability with 
\[\exp(-  \lambda |\Upsilon_{j+1,\geq \Wtilde[j+1]}| \, \Wtilde[j+1] /80),\] which is at most
$\exp(-  \lambda  j /500)$ since  $j\geq j_0$ is sufficiently large
and  $\Wtilde[j+1] \, |\Upsilon_{j+1,\ge {\Wtilde[j+1]}}|\geq j$ (by Observation~\ref{obs:Wtilde}).
If $j+1$ is heavy-covered then the same lemma upper-bounds this probability with
$\exp(-(j+1)^{1.9}) \leq \exp(-  \lambda  j /500)$.
Finally, if $j+1$ is exposed then Lemma~\ref{lem:T3tofill}
upper bounds the probability with
$\exp(- \lambda(j-1)/400)) \leq \exp(-  \lambda  j /500)$.

The lemma follows in a similar way to the proof of Lemma~\ref{lem:badRtwo} 
since the probability of the given transition, assuming 
$\Psivect_{T(t-1)} = \psivect_{T(t-1)}$
is at most 
\[ 
\sum_{t' = t-1}^{t+A+ j^{24}-1} 
\sum_{a =1}^{j^{24}}
\Pr(\tmin = t'  -a + 1 \wedge \amin = a
\mid \Psivect_{T(t-1)} = \psivect_{T(t-1)}
)
\max_{\psivect_{T(t')}}
\Pr(\calE \mid \Psivect_{T(t')} = \psivect_{T(t')}),
\]
where   the maximum is over 
vectors  $\psivect_{T(t')}$ of high-level states which extend $\psivect_{T(t-1)}$,
have positive probability as a value of 
$\Psivect_{T(t')}$, have $\psi_{t'}$ as a $j$-filling state,
and entail $\tmin = t' - a + 1$ and $\amin = a$.
The result follows since
$\Pr(\calE \mid \Psivect_{T(t')} = \psivect_{T(t')})$
is uniformly at least $\exp(-\lambda j/500)$.
\end{proof}

\begin{lemma}\label{lem:badLoopStarts}
Consider the setting of Definition~\ref{def:morefromfilling} and the events of Definition~\ref{def:fillingevents}.  Then
\[\Pr(\neg\calE_{{\badLoopStarts}}\mid \Psivect_{T(t-1)} = \psivect_{T(t-1)}) \leq  \exp(-\lambda(j-1)/800) .\]
\end{lemma}
\begin{proof}
We wish to upper-bound the probability that
$|\LoopStarts| >  (A+1) \exp(- \lambda(j-1)/800)$.
For any $t'$ satisfying 
$t\leq t' \leq t+A$,
we will upper bound the probability that $t' \in \LoopStarts$. 
From the definition, $t'$ can only be in $\LoopStarts$ if $\Psi(t'-1)$ is a $j$-filling state and the transition to $\Psi(t')$ is via (R\ref{item:R2})(ii). By Lemma~\ref{lem:T3tofill}, conditioning on $\Psivect_{T(t'-1)} = \psivect_{T(t'-1)}$ for any $\psivect_{T(t'-1)}$ with $\type^{\psi_{t'-1}} = \filling$,
the probability of this transition is at most 
$\exp(- \lambda(j-1)/400)$.

Let $\mu := (A+1) \exp(- \lambda(j-1)/400)$.
Then $\E[|\LoopStarts| :  \Psivect_{T(t-1)} = \psivect_{T(t-1)}] \leq  \mu$.
By Markov's inequality,
$\Pr[|\LoopStarts| \geq 
\exp(\lambda(j-1)/800)   \mu 
\leq \exp(-\lambda(j-1)/800)$.
\end{proof}

\begin{lemma}\label{lem:badLongLoop}
Consider the setting of Definition~\ref{def:morefromfilling} and the events of Definition~\ref{def:fillingevents}.  Then
\[\Pr(\neg\calE_{{\badLongLoop}}\mid \Psivect_{T(t-1)} = \psivect_{T(t-1)}) \leq  A \exp(-\lambda j(\log j )^2/203).\]
\end{lemma}
\begin{proof}
We wish to upper bound the probability of the
event that, for some  $x\in [A]$ with $\In_{x} < \tau_{j,\full}$, 
\begin{enumerate}[(I)]
    \item  $\Out_{x} > \In_{x} + j^{\phi+74}$, or
    \item 
$\Out_{x} \geq \tau_{j,\full}$,
 or
 \item 
$F_j^{\Psi(\Out_{x})}(\Out_{x}) < \min\{ F_j^{\Psi(\In_{x}-1)}(\In_{x}),\mu_j^{\gammavect}\} - 2\lambda j^{\phi+74}$. \end{enumerate}

Similar to the proofs of the other lemmas, 
we do this 
for a fixed~$x$
by first considering a fixed time~$t'\geq t-1$ and a fixed vector $\psivect_{T(t')}$ of high-level states which extends $\psivect_{T(t)}$, 
such that $\Pr(\Psivect_{T(t')} = \psivect_{T(t')}) > 0$ 
and entails  
that $I_x = t' < \tau_{j,\full}$. 
From the definition of~$I_x$, the transition from $\psi_{t'-1}$ to $\psi_{t'}$ is via (R\ref{item:R2})(ii)
and $\psi_{t'}$ is $j$-refilling.

We are going to use Lemma~\ref{lem:main-loop-bounds}.  The $t-1$ of 
Lemma~\ref{lem:main-loop-bounds}
is our $t'$. The lemma shows that,
conditioned on $\Psivect_{T(t')} = \psivect_{T(t')}$,
with 
probability at least $1-\exp(-\lambda j(\log j )^2/203)$,  
the following events occur. 
\begin{enumerate}[(i)]
\item $\exit(\In_x) \leq \In_x+ j^{\phi+74}$.
\item The transition 
into $\Psi(\exit(\In_x))$ is via (R\ref{item:R6})(i) 
into a $j$-filling state. 
\item $F_j^{\Psi(\exit(\In_x))}(\exit(\In_x)) \geq \min\{F_j^{\Psi({\In_x-1})}(\In_x),\mu_j^{\gammavect}\} - 2\lambda j^{\phi+74}$.
\end{enumerate}
From Remark~\ref{rem:fromfilling}, $\Out_x = \exit(\In_x)$.
Thus, (i)  
rules out (I) and (ii) rules out (II).
Finally, (iii)   is 
$F_j^{\Psi(\Out_{x})}(\Out_{x}) \ge \min\{ F_j^{\Psi(\In_{x}-1)}(\In_{x}),\mu_j^{\gammavect}\} - 2\lambda j^{\phi+74}$,
which  rules out (III) for this~$\psivect_{T(t')}$.

Similarly to the other lemmas, since this bound on the failure probability is uniform over the possible vectors 
$\psivect_{T(t')}$, we get the same upper bound on the failure probability for this~$x$. The lemma follows by a union bound over~$x$. 
\end{proof}

\begin{corollary}\label{cor:fromfillingEs}
Consider the setting of Definition~\ref{def:morefromfilling} and the events of Definition~\ref{def:fillingevents}.  Then
\[\Pr( \neg \calE_{\badRone} \vee \neg \calE_{\badRtwo} \vee \neg \calE_{\badRfour} \vee \neg \calE_{\badLoopStarts} \vee \neg \calE_{\badLongLoop} )  \leq   
\exp(-\lambda j/1000).\]
\end{corollary}
\begin{proof}
We apply a union bound to the 5 events, using
Lemmas~\ref{lem:badRone}, \ref{lem:badRtwo}, \ref{lem:badR4},
\ref{lem:badLoopStarts}, and \ref{lem:badLongLoop} to get bounds on the individual failure probabilities, and adding those failure probabilities. 

We will show that, since $j\geq j_0$ is sufficiently large, the sum is at most the bound given in the lemma statement.

For this we use the hypothesis from the initial setting in Definition~\ref{def:fromfilling} that 
for all $i \ge 1$, $W_i \le (4/\lambda)^i$.
Note that the quantity $A$ is at most  
$8(4j)^{\phi} W_j \leq 8 (4j)^{\phi} (4/\lambda)^j$. In the failure probabilities in Lemmas~\ref{lem:badRone} and~\ref{lem:badLongLoop}, this quantity is multiplied by a quantity that is smaller than any inverse exponential function of~$j$
so, since $j$ is sufficiently large,  these two bounds are each at most $\exp(-\lambda j/900)$. This quantity similarly upper bounds the failure probabilities of the other three lemmas. Thus, very crudely,
adding 5~of these is at most $\exp(-\lambda j/1000)$.
\end{proof}

We next give some lemmas 
which derive some consequences of 
$
\calE_{\badRone} \wedge \calE_{\badRtwo} \wedge \calE_{\badRfour} \wedge \calE_{\badLoopStarts} \wedge \calE_{\badLongLoop}$. Finally, we use these to prove Lemma~\ref{lem:OnejVolumeAnalysis}, which uses
the fact (from Corollary~\ref{cor:fromfillingEs})
that the events are likely, so it is likely that these consequences hold. Also, given 
the consequences, the volume process transitions from  a $j$-filling high-level state into a $(j+1)$-advancing or $(j+1)$-filling state before failing, and   in $\poly(j)W_j$ steps.

Lemma~\ref{lem:fill-Fj-bound} is an important part of our analysis. It considers a time $t'$ which is after the process entered a $j$-filling state at time $\Out_{x-1}$ but is before this $j$-filling state is left.
It gives a lower-bound on 
$F_j^{\Psi(t'-1)}(t')$ -- once this quantity is at least~$\mu_j^{\gammavect}$ the process has the possibility of moving to a  $(j+1)$-filling or $(j+1)$-advancing state via (R\ref{item:R4})(i). At time~$\Out_0=t-1$
we have $F_j^{\Psi(\Out_0)}(\Out_0) = 0$ (from  the setting of Definition~\ref{def:fromfilling}, which arises when the $j$-filling state is reached via (R\ref{item:R3}) or (R\ref{item:R4})). The bound shows that roughly $2\lambda/3$ is added to~$F_j$ with each step. In fact this is how much is gained during steps when the process is in a $j$-filling state, and not in the refilling-stabilising loop -- our precise bound is that the amount added for each such step is at least $2\lambda/3$ divided by a polynomial in~$j$.
The (small) multiple of~$W_j$ that is subtracted off  in the lemma statement arises in following way.
The number of steps spent  
in the refilling-stabilising loop is roughly at most $A j^{\phi+74} \exp(-\lambda(-(j-1)/800)$  and in each of these steps $F_j$ decreases by around at most~$\lambda$.

\begin{lemma}\label{lem:fill-Fj-bound}
Consider the setting of Definition~\ref{def:morefromfilling} and the events of Definition~\ref{def:fillingevents}.
Suppose that $\calE_{\badLoopStarts} \wedge \calE_{\badLongLoop} $ occurs.
Fix  
any $t'$ in the range 
$t\leq t' \leq t+A-j^{24}+1$ 
such that 
$\Psi(t'-1)$ is $j$-filling
and $[t'-1 ] \cap \Ready = \emptyset$.
Let $x$ be the positive integer so that
$\Out_{x-1} \leq t'-1 < \In_{x}$.
Then $F_j^{\Psi(\Out_{x-1})}(t') \geq
\frac{2\lambda(t'-(t-1))}{3(4j)^\phi} -   W_j \exp(-\lambda(j-1)/1600)$. 
\end{lemma}

\begin{proof}
Consider any $x$ and $t'$ that satisfy the conditions in the lemma statement. 
Since $\Out_{x-1} \leq t'-1 < \In_x$, $\Out_{x-1}$ is  finite and $\Out_{x-1} \leq t'-1 < \tau_{j,\full}$.
 First write $F_j^{\Psi(\Out_{x-1})}(\Out_{x-1})$ as a telescoping sum to
 observe that
\begin{align}\nonumber
F_j^{\Psi(\Out_{x-1})}(t') &= F_j^{\Psi(\Out_{x-1})}(t') - F_j^{\Psi(\Out_{x-1})}(\Out_{x-1}) + \sum_{y=1}^{x-1} \big(F_j^{\Psi(\Out_y)}(\Out_y) - F_j^{\Psi(\Out_{y-1})}(\Out_{y-1}) \big) + F_j^{\Psi(\Out_0)}(\Out_0)\\\label{eq:one-j-vol-claim-1-main}
&\ge F_j^{\Psi(\Out_{x-1})}(t') - F_j^{\Psi(\Out_{x-1})}(\Out_{x-1}) + \sum_{y=1}^{x-1} \big(F_j^{\Psi(\Out_y)}(\Out_y) - F_j^{\Psi(\Out_{y-1})}(\Out_{y-1}) \big).
\end{align}
The reason for writing this expression in this curious way is that the expression in the sum,
\[(F_j^{\Psi(\Out_y)}(\Out_y) - F_j^{\Psi(\Out_{y-1})}(\Out_{y-1})),\] measures the change in  
the $F_j$ value (for the appropriate high-level state) between $\Out_{y-1}$ (when the loop is left for the $(y-1)$'th time), and $\Out_{y}$, when it is left for the $y$'th time. In order to study this quantity, we separate it, in~\eqref{eq:one-full-loop}, into 
the portion between $\Out_{y-1}$ and $\In_y$ (when the volume process is in a $j$-filling state), and  the portion between $\In_y$ and $\Out_y$ (when the volume process is in the $y$'th refilling-stabilising loop). 
\begin{equation}
\label{eq:one-full-loop}  
F_j^{\Psi(\Out_y)}(\Out_y) - F_j^{\Psi(\Out_{y-1})}(\Out_{y-1}) = \big(F_j^{\Psi(\Out_y)}(\Out_y) - F_j^{\Psi(\In_y-1)}(\In_y)\big) + \big(F_j^{\Psi(\In_y-1)}(\In_y) - F_j^{\Psi(\Out_{y-1})}(\Out_{y-1})\big).
\end{equation}

We wish to bound the quantity on the right-hand-side of~\eqref{eq:one-full-loop} below.
Our first goal is to bound the term from~$\In_y$ to~$\Out_y$ when the volume process is in the refilling-stabilising loop. We wish to show that the amount lost during the loop is at most~$2\lambda j^{\phi+74}$ 
(essentially losing~$\lambda$ for each step in the loop, according to Lemma~\ref{lem:main-loop-bounds} -- the factor of~$2$ is an over-estimate). Thus, we will lower-bound the 
first quantity on the right-hand-side of
\eqref{eq:one-full-loop}  by $-2\lambda j^{\phi+74}$.

While we are doing this, we will use the fact that $\Psi(\In_y - 1) = \Psi(\Out_{y-1})$ and that this is a $j$-filling state (these follow from Observation~\ref{obs:stay-jfill-state} using the fact that $\Out_{y-1} \leq \In_y -1$ which follows from Remark~\ref{rem:fromfilling} since $\Out_{y-1} < \tau_{j,\full}$).
We will derive~\eqref{eq:filling-todo-fullloop} 
from~\eqref{eq:one-full-loop} both by substituting the lower bound for the period in the loop, and making this state translation for the other term (ready for the next part of our analysis).

In order to derive our lower bound on how $F_j$ changes during the refilling-stabilising loop, we first  derive a crude upper bound on~$x$.
Let $\calI = \{ \hat{t}: 
t-1 \leq \hat{t} \leq t'-1 \mbox{ and, for some $k\geq 0$, $\hat{t} \in \{\In_k,\Out_k\}$}
\}$. Since $\Out_{x-1} \leq t'-1$,
there are at least $2(x-1)+1$ elements in~$\calI$ so  $t'-(t-1) \geq 2(x-1)+1$ and since $t'-1 \leq t+A$, we have $x-1\leq (A+1)/2 \leq A$. 

We have just shown that every $y\in [x-1]$ satisfies the pre-conditions of $\calE_{\badLongLoop}$:  
$y\in [A]$ and $\In_{y} < \tau_{j,\full}$. By $\calE_{\badLongLoop}$, $F_j^{\Psi(\Out_{y})}(\Out_{y}) \geq \min\{ F_j^{\Psi(\In_{y}-1)}(\In_{y}),\mu_j^{\gammavect}\} - 2\lambda j^{\phi+74}$. We remove the minimisation by using  Observation~\ref{obs:readysoon} 
with $T' = t'-1$ and $T'' = \In_y-1$
to show $F_j^{\Psi(\In_{y}-1)}(\In_{y}) < \mu_j^{\gammavect}$.  Thus, $F_j^{\Psi(\Out_{y})}(\Out_{y}) \geq   F_j^{\Psi(\In_{y}-1)}(\In_{y})  - 2\lambda j^{\phi+74}$.
Plugging this into the first expression on the right-hand-side of~\eqref{eq:one-full-loop} while simultaneously substituting $\Psi(\Out_{y-1})$ 
for $\Psi(\In_y-1)$ in the second expression on the right-hand-side,
it follows that for all $y \in [x-1]$,
\begin{equation}\label{eq:filling-todo-fullloop}
F_j^{\Psi(\Out_y)}(\Out_y) - F_j^{\Psi(\Out_{y-1})}(\Out_{y-1}) \ge \big(F_j^{\Psi(\Out_{y-1})}(\In_y) - F_j^{\Psi(\Out_{y-1})}(\Out_{y-1})\big) - 2\lambda j^{\phi+74}.
\end{equation} 

Our next goal is to bound the 
change in~$F_j$ between $\Out_{y-1}$ when a $j$-filling state is entered and~$\In_y$, when it is left again. We will show that the amount gained is roughly at least $2\lambda/3$ per step (our precise bound is this divided by a polynomial). Thus,  we will bound the first 
term on the right-hand-side of~\eqref{eq:filling-todo-fullloop} below by 
$(2\lambda/3) (\In_y - \Out_{y-1} )/(4j)^\phi$
to get~\eqref{eq:bounded-full-loop} for all $y\in [x-1]$.
We start by lower-bounding $F_j^{\Psi(\Out_{y-1})}(\In_y)$. By the definition of $\Out_{y-1}$ (using the fact that it is less than $\tau_{j,\full}$),
$\Psi(\Out_{y-1})$ is a $j$-filling state, so by the definition of filling state (Definition~\ref{def:filling-state}), for all $\ell \in [j-1]$, $F_\ell^{\Psi(\Out_{y-1})}(\Out_{y-1}) = \mu_\ell^{\gammavect}$. Applying Corollary~\ref{cor:f-fill-bin} to state~$\Psi(\Out_{y-1})$ with the $k$ of the corollary as~$j-1$, the $t'$ of the corollary as $\Out_{y-1}$, and the $t''$ in the penultimate sentence of the corollary statement as~$\In_y$, we obtain 
 \[
F_{j}^{\Psi(\Out_{y-1})}(\In_y) \ge \min\{\mu_{j}^{\gammavect},\ F_{j}^{\Psi(\Out_{y-1})}(\Out_{y-1}) + (\In_y-\Out_{y-1}) \gamma_{j}p_{j-1}\mu_{j-1}^{\gammavect}\}.\] We have already 
shown that $F_{j}^{\Psi(\Out_{y-1})}(\In_y)=
F_j^{\Psi(\In_{y}-1)}(\In_{y}) < \mu_j^{\gammavect}$, so the minimum cannot be attained at $\mu_j^{\gammavect}$ and we have
 $
F_{j}^{\Psi(\Out_{y-1})}(\In_y) -  \ F_{j}^{\Psi(\Out_{y-1})}(\Out_{y-1}) \geq (\In_y-\Out_{y-1}) \gamma_{j}p_{j-1}\mu_{j-1}^{\gammavect}$.

We finish by showing 
$\gamma_j p_{j-1} \mu_{j-1}^{\gammavect} \geq (2\lambda/3)/(4j)^{\phi}$ -- a fact that we will re-use it later. It follows because 
Observation~\ref{obs:mu-gamma-bounds} gives $p_{j-1}\mu_{j-1}^{\gammavect} \geq 2\lambda/3$ and $\gamma_j = 1/(4 j)^\phi$ (Definition~\ref{def:real-gamma}). Thus, we have the desired lower bound for the first term on the right-hand-side of \eqref{eq:filling-todo-fullloop} so we achieve the following for all $y\in [x-1]$.

\begin{equation}\label{eq:bounded-full-loop}
F_j^{\Psi(\Out_y)}(\Out_y) - F_j^{\Psi(\Out_{y-1})}(\Out_{y-1}) \geq (2\lambda/3) (\In_y - \Out_{y-1} )/(4j)^\phi - 2\lambda j^{\phi+74}.
\end{equation}
With~\eqref{eq:bounded-full-loop}, we have achieved the goal of lower-bounding each term in the sum in~\eqref{eq:one-j-vol-claim-1-main}, representing the change in~$F_j$ over a  full iteration from $\Out_{y-1}$ to $\Out_y$. Plugging this in, we get 

\begin{equation}\label{eq:one-j-vol-bound}
F_j^{\Psi(\Out_{x-1})}(t') \ge F_j^{\Psi(\Out_{x-1})}(t') - F_j^{\Psi(\Out_{x-1})}(\Out_{x-1}) + \Big(\frac{2\lambda}{3(4j)^\phi}\sum_{y=1}^{x-1} (\In_y - \Out_{y-1} ) \Big)- 2\lambda(x-1) j^{\phi+74}.
\end{equation}

Recall that our goal is to lower-bound the right-hand-side of~\eqref{eq:one-j-vol-bound} with the quantity in the lemma statement, which is $(2\lambda/(3 (4j)^{\phi})$ for each time step, with all of the error that gets subtracted off upper-bounded in terms of~$W_j$.
Our first goal will be to lower-bound the 
change between $\Out_{x-1}$ and $t'$ via 
$F_j^{\Psi(\Out_{x-1})}(t') - F_j^{\Psi(\Out_{x-1})}(\Out_{x-1})\geq (2\lambda/3)(t' - \Out_{x-1} )/(4j)^\phi$ to achieve~\eqref{eq:ready-to-simplify}. This is similar to what we just did for the period between refilling-stabilising loops.
To  lower-bound $F_j^{\Psi(\Out_{x-1})}(t')$
We apply Corollary~\ref{cor:f-fill-bin}, using the fact that $\Psi(\Out_{x-1})$ is a $j$-filling state so 
for all $\ell \in [j-1]$, $F_\ell^{\Psi(\Out_{x-1})}(\Out_{x-1}) = \mu_\ell^{\gammavect}$. Applying Corollary~\ref{cor:f-fill-bin} to state~$\Psi(\Out_{x-1})$ with the $k$ of the corollary as~$j-1$, the $t'$ of the corollary as $\Out_{x-1}$, and the $t''$ in the penultimate sentence of the corollary statement as our~$t'$, we obtain 
 $
F_{j}^{\Psi(\Out_{x-1})}(t') \ge \min\{\mu_{j}^{\gammavect},\ F_{j}^{\Psi(\Out_{x-1})}(\Out_{x-1}) + (t'-\Out_{x-1}) \gamma_{j}p_{j-1}\mu_{j-1}^{\gammavect}\}$.
We remove the minimisation by showing that the left-hand-side, $F_{j}^{\Psi(\Out_{x-1})}(t')$ is less than $\mu_j^{\gammavect}$, so
the minimum can't be achieved at $\mu_j^{\gammavect}$. To do this, we show $F_{j}^{\Psi(\Out_{x-1})}(t')
< \mu_j^{\gammavect}
$ by applying Observation~\ref{obs:readysoon} with $T' = T''=t'-1$ noting that $t'-1< \In_x$ so $\Psi(t'-1) = \Psi(\Out_{x-1})$.
Substituting without the minimisation, we have  $
F_{j}^{\Psi(\Out_{x-1})}(t') -  \ F_{j}^{\Psi(\Out_{x-1})}(\Out_{x-1}) \geq (t'-\Out_{x-1}) \gamma_{j}p_{j-1}\mu_{j-1}^{\gammavect}$. We have already shown that $\gamma_j p_{j-1} \mu_{j-1}^{\gammavect} \geq (2\lambda/3)/(4j)^{\phi}$, so plugging this into~\eqref{eq:one-j-vol-bound}, we get the following.
\begin{equation}\label{eq:ready-to-simplify}
F_j^{\Psi(\Out_{x-1})}(t') \ge  \frac{2\lambda}{3(4j)^\phi}\Big(t' - \Out_{x-1}  + \sum_{y=1}^{x-1} (\In_y - \Out_{y-1} ) \Big)- 2\lambda(x-1) j^{\phi+74}.
\end{equation}

Our goal is now is to lower-bound the right-hand-side of~\eqref{eq:ready-to-simplify} with the quantity in the lemma statement. 
The expression in parantheses   in~\eqref{eq:ready-to-simplify} is the amount of time spent in $j$-filling, and we wish to bound this in terms of the total time. We start with a telescoping sum, 
\[
t' = t' - \Out_{x-1} + \sum_{y=1}^{x-1} (\Out_y - \In_y) + \sum_{y=1}^{x-1}(\In_y - \Out_{y-1}) + \Out_0.
\]
We have previously shown that every $y\in [x-1]$ satisfies the pre-conditions of $\calE_{\badLongLoop}$:
$y\in [A]$ and $\In_{y} < \tau_{j,\full}$. By $\calE_{\badLongLoop}$, $\Out_{y} - \In_{y} \leq j^{\phi+74}$.
Substituting this in, along with $\Out_0=t-1$, we get
\[
t' \le t' - \Out_{x-1} + (x-1)j^{\phi+74} + \sum_{y=1}^{x-1}(\In_y - \Out_{y-1}) + t-1.
\]

Rearranging, we obtain
\[
t' - \Out_{x-1} + \sum_{y=1}^{x-1} (\In_y-\Out_{y-1}) \ge t' - (t-1) -(x-1)j^{\phi+74}.
\]

Plugging this into~\eqref{eq:ready-to-simplify} yields
\begin{align*}
F_j^{\Psi(\Out_{x-1})}(t') &\ge  \frac{2\lambda}{3(4j)^\phi}(t' - (t-1)  -(x-1)j^{\phi+74}) - 2\lambda(x-1) j^{\phi+74}\\
&\ge \frac{2\lambda(t'-(t-1))}{3(4j)^\phi} - 3\lambda(x-1) j^{\phi+74}.
\end{align*}

It is clear from the lemma statement that to finish, we need only show 
$W_j \exp(-\lambda(j-1)/1600) \geq 
3\lambda(x-1) j^{\phi+74}$.
For this will obtain a slightly tighter upper bound on $x-1$
noting from the lemma statement that 
$\Out_{x-1}\leq t'-1 \leq \In_x$ 
and  (as we have already concluded) $\Out_{x-1} < \tau_{j,\full}$, so
$\In_{x-1} < \Out_{x-1}$. Thus,
there are at least $x-1$ loop starts, namely $\In_1, \ldots, \In_{x-1}$,  
before time $t'-1$. Since $t'-1 \leq t+A$ by the lemma statement, all of these loop starts are in the set $\LoopStarts$.
By $\calE_{\badLoopStarts}$, 
$|\LoopStarts|
\leq  (A+1) \exp(- \lambda(j-1)/800))$, 
so $x-1 \leq (A+1) \exp(- \lambda(j-1)/800))$. 
Thus, the bound follows 
since $A$ is $W_j$ times a polynomial in~$j$, and $j\geq j_0$ is sufficiently large. 
\end{proof}

Lemmas~\ref{lem:ReadySoonHappens} and~\ref{lem:ReadySoonGoodEnough} are the main building blocks for proving Lemma~\ref{lem:OnejVolumeAnalysis}.

\begin{lemma}\label{lem:ReadySoonHappens}
Consider the setting of Definition~\ref{def:morefromfilling} and the events of Definition~\ref{def:fillingevents}. 
Suppose that  $\calE_{\badRone}\wedge 
\calE_{\badLoopStarts}\wedge
\calE_{\badLongLoop}$ occurs.
Then
$\Ready$ is nonempty.
\end{lemma}
\begin{proof}
Suppose for contradiction that $\Ready = \emptyset$. 
We will show that there is a $t'$ in the range
$t-1 \leq t' \leq t+A -j^{24}$ such that
\begin{enumerate}[(i)]
\item $ \Psi(t')$ is a $j$-filling state, and
\item $F_j^{\Psi(t')}
(t' + 1) \geq \mu_j^{\gammavect}$.
\end{enumerate}
From the definition of~$\Ready$, these imply that $t' \in \Ready$, giving a contradiction.

In the next part of the proof we will establish that (due to the fact that
$\calE_{\badRone}\wedge 
\calE_{\badLoopStarts}\wedge
\calE_{\badLongLoop}$ occurs and that $\Ready=\emptyset$, which we have assumed for contradiction) the volume process
has the following predictable behaviour, starting from $\Out_0=t-1$, 
and continuing until (at least) time $t+A-j^{24}$.
After it is in a $j$-filling state $\Psi(\Out_x)$,
for some non-negative integer~$x$,
it stays in that high-level state for a while (for one or more steps), and
then enters the refilling-stabilising loop at~$I_{x+1}$, staying in the loop for at most $j^{\phi+74}$ steps, and then exits into a $j$-filling state $\Out_{x+1}$, repeating this behaviour for~$x+1$.

To see this, we show that 
$\calE_{\badRone}\wedge 
\calE_{\badLoopStarts}\wedge
\calE_{\badLongLoop}$ and the fact that $\Ready=\emptyset$ rules out any other transitions that the process could take.

Start from the $j$-filling state $\Psi(\Out_x)$.
From this high-level state, $\calE_{\badRone}$ rules out a transition via (R\ref{item:R1}). Transition (R\ref{item:R2}) is no problem - it just triggers~$\In_{x+1}$, which is part of the ``predictable behaviour''.  
Since we have assumed $\Ready=\emptyset$,
Observation~\ref{obs:readysoon} rules out 
a transition via  (R\ref{item:R4})(i). 
Transition (R\ref{item:R4})(ii) is also part of the ``predictable behaviour''  -- it keeps the high-level state unchanged. No other transitions are applicable from filling states.

So now consider what happens when the refilling-stabilising loop is entered at $\In_{x+1}$. First, we can assume $x+1 <A$ since $\calE_{\badLoopStarts}$ guarantees 
$|\LoopStarts|
\leq  (A+1) \exp(- \lambda(j-1)/800)) < A$. But then $\calE_{\badLongLoop}$ guarantees that 
$\Out_{x+1} \leq \In_{x+1} + j^{\phi+74}$ and 
$\Out_{x+1} < \tau_{j,\full}$. So the predictable behaviour continues, and a $j$-filling state is entered at $\Out_{x+1}$.
Thus, we have established the predictable behaviour 
starting from time~$t-1$ and continuing until time $t+A-j^{24}$.

Now consider any $t''$ in the range $t \leq t'' 
\leq t+A-j^{24}+1$ such that $\Psi(t''-1)$ is $j$-filling.
Let $x$ be the  positive integer so that $\Out_{x-1} \leq t''-1 <\In_{x}$.
By Lemma~\ref{lem:fill-Fj-bound},
\[F_j^{\Psi(\Out_{x-1})}(t'') \geq
\frac{2\lambda(t''-(t-1))}{3(4j)^\phi} -   W_j \exp(-\lambda(j-1)/1600).\] 
 
Also, by Observation~\ref{obs:stay-jfill-state}, $\Psi(\Out_{x-1}) = \Psi(t''-1)$.
So $t'=t''-1$ is in the required range and satisfies (i) It will satisfy~(ii) as long as
\[\frac{2\lambda(t''-(t-1))}{3(4j)^\phi} \geq \mu_j^{\gammavect}
 +
  W_j \exp(-\lambda(j-1)/1600).\]

Since $j\geq j_0$ is sufficiently large
and $\mu_j^{\gammavect} \leq 3 \lambda W_j/4$ (by Observation~\ref{obs:mu-gamma-bounds}),
it suffices to show 
\[
2\lambda(t''-(t-1))/(3(4j)^\phi) \geq \lambda W_j,\]
so plugging in the definition of~$A$, it suffices to show $t'' - (t-1) \geq A/2$.

So we have to be able to a choose a $t''$ in the
range $(t-1)+A/2 \leq t'' 
\leq t+A-j^{24}+1$ such that $\Psi(t''-1)$ is a $j$-filling state.  
The ``predictable behaviour'' is useful, since there are only $j$-filling states (apart from states in the refilling-stabilising loop).
This range has at least $A/2-j^{24} \geq A/3$ steps.
But by $\calE_{\badLoopStarts}$ 
there are at most 
$(A+1) \exp(- \lambda(j-1)/800))$ loop starts,
each of which takes at most $j^{\phi+74}$ time steps by $\calE_{\badLongLoop}$, and
$(A+1) \exp(- \lambda(j-1)/800))j^{\phi+74} < A/3$, so there is a suitable~$t''$.
\end{proof}
\begin{lemma}\label{lem:ReadySoonGoodEnough}
Consider the setting of Definition~\ref{def:morefromfilling} and the events of Definition~\ref{def:fillingevents}. 
Suppose that  $\calE_{\badRone} \wedge 
\calE_{\badRtwo} \wedge \calE_{\badRfour}$ 
 occurs.
If $\Ready$ is non-empty then $\tau_{j,\full} \leq t+A$ and 
$\Psi(\tau_{j,\full})$ is 
$(j+1)$-filling or $(j+1)$-advancing. Also, $z_{j+1}^{\Psi(\tau_{j,\full})}=0$.

\end{lemma} 
\begin{proof}
From the definition of~$\Ready$, any step~$t'$
that is in~$\Ready$ is in the range $t-1 \leq t' \leq t+A-j^{24}$.

Fix any $t'$ in this range  
Let $\psivect_{T(t')}$ be a vector of high-level states such that
$\Pr(\Psivect_{T(t')} = \psivect_{T(t')})>0$
and $\psivect_{T(t')}$ entails
that 
$\Ready$ is non-empty and $t' = \tmin$ (so $t'$ is the first time step in $\Ready$). 
From the definition
of $\Ready$,  $\psi_{t'}$ is a $j$-filling state. 
Now consider the volume process, conditioned on 
$\Psivect_{T(t')} = \psivect_{T(t')}$.
Note that the value of~$\amin$ can be determined from $\psivect_{T(t')}$, and (see Remark~\ref{rem:morefromfilling}), it is  in~$[j^{24}]$.

We will show that for all $a\in [\amin -1]$,
the transition into~$\Psi(t'+a)$ is via (R\ref{item:R4})(ii), so the high-level state does not change.
To see this, note that the transition is not via (R\ref{item:R1}) since we have $\calE_{\badRone}$.
It is not via (R\ref{item:R2}) since we have $\calE_{\badRtwo}$
 and $\amin \leq j^{24}$.
It is not via (R\ref{item:R4})(i) from the definition of~$\amin$ since $a<\amin$ so 
$\tau^{\Psi(t'+a-1)} + j^{24} > \tmin + a$, so the guard  
on transition (R\ref{item:R4})(i) in the backoff bounding rule fails.
No other transitions apply.

The same considerations rule out 
(R\ref{item:R1}) and (R\ref{item:R2}) for the transition into $\Psi(t'+\amin)$. 
There are no other applicable transitions except (R\ref{item:R4}). Note that this transition is from the high-level state
$\Psi(t'+\amin-1)$ which is the same as $\Psi(t')$.
By the definition of $\amin$, 
$\tau^{\Psi(t'+\amin-1)} + j^{24} \leq \tmin + \amin$
so the first guard on (R\ref{item:R4})(i) succeeds.
Lemma~\ref{lem:filling-F-monotonoe} shows that
$F_j^{\Psi(t')}(t'+\amin) \geq \mu_j^{\gamma}$, so the second guard also succeeds and this transition is taken.
The transition is not to a failure state since we have $\calE_{\badRfour}$, so it is to a $(j+1)$-filling or $(j+1)$-advancing state $\Psi(t'+\amin)$.
From the backoff bounding rule, and the definition of advancing state, $z_{j+1}^{\Psi(t'+\amin)} = 0$.
The transition triggers~$\tau_{j,\full}$, so $\tau_{j,\full} = t'+\amin \leq t'+ j^{24} \leq t+A$.
 \end{proof}

We can finally prove Lemma~\ref{lem:OnejVolumeAnalysis}. The setting of the lemma is the initial setting from Definition~\ref{def:morefromfilling} with 
an additional constraint which we don't use here, but we add it to the statement for easy combination with other lemmas. Namely, we assume that $\lambda \leq \lambda_0$ for some~$\lambda_0$. We repeat the setting in the lemma statement in order to make it easier to apply the lemma.

 \begin{lemma}
\label{lem:OnejVolumeAnalysis}
There exists a real number $\lambda_0 \in (0,1/2)$ such that the following holds. Consider any $\lambda \in (0,\lambda_0)$ and a send sequence $\p$ such that $p_0=1$ and, for all $i \ge 1$, $W_i \le (4/\lambda)^i$. Suppose that $\p$ has finitely many bins which are strongly exposed. 
Let $j_0$ be a sufficiently large positive integer
where  all strongly-exposed bins are in $[j_0-1]$.
Let $\tau_0$ be a non-negative integer. Let $\Psi_0$ be the   initialising high-level state with  $\tau^{\Psi_0} = \tau_0$ and $j^{\Psi_0} = j_0$. 
Let $\hls$ be the $\Psi_0$-backoff bounding rule from Definition~\ref{def:backoff-bounding-rule}.  
Let $V$ be a volume process from $\Psi_0$ with transition rule $\hls$, 
 send sequence~$\p$, and birth rate~$\lambda$  (Definition~\ref{def:volume-process}).
For any
$t \ge \tau_0$, 
let   $T(t) = \{\tau_0,\ldots,t\}$.
Then  the following statement holds  for
every $t > \tau_0$ and
every tuple
$ 
\psivect_{T(t-1)} = (\psi_{\tau_0},\dots,\psi_{t-1})$  of high-level states 
such that  
$\Pr( \Psivect_{T(t-1)} = \psivect_{T(t-1)}) > 0$, 
$\psi_{t-1}$ is a $j$-filling state, 
$\tau^{\psi_{t-1}} = t-1$, and $z_j^{\psi_{t-1}} = 0$.
Conditioned on $\Psivect_{T(t-1)} = \psivect_{T(t-1)} $,
with 
probability at least $1-
\exp(-\lambda j/1000) $,  
the following events occur, where  
\[
\tau_{j,\full} = \min\{t'  \geq t \colon \type^{\Psi(t')} = \Failure \mbox{ or }j^{\Psi(t')} = j+1\}.\]
\begin{enumerate}[(i)]
\item $\Psi(\tau_{j,\full})$ is $(j+1)$-filling or $(j+1)$-advancing
and $z_{j+1}^{\Psi(\tau_{j,\full})}=0$.

\item $\tau_{j,\full} \leq t+
\lceil 4(4 j)^{\phi} W_j\rceil
$.
\end{enumerate}
\end{lemma}

\begin{proof}
We use the notation from the initial setting in Definition~\ref{def:morefromfilling} and the events from Definition~\ref{def:fillingevents}. By Corollary~\ref{cor:fromfillingEs},
\[\Pr( \calE_{\badRone} \wedge \calE_{\badRtwo} \wedge \calE_{\badRfour} \wedge \calE_{\badLoopStarts} \wedge \calE_{\badLongLoop}  )  \geq 1-   
\exp(-\lambda j/1000).\]

Suppose that $\calE_{\badRone} \wedge \calE_{\badRtwo} \wedge \calE_{\badRfour} \wedge \calE_{\badLoopStarts} \wedge \calE_{\badLongLoop}  $ occurs.
Then by Lemma~\ref{lem:ReadySoonHappens}, the set $\Ready$ is non-empty. By Lemma~\ref{lem:ReadySoonGoodEnough}, (i) and (ii) hold.
\end{proof}

Lemma~\ref{lem:OnejVolumeAnalysis}  
showed
that, from a $j$-filling high-level state,   the volume process is very likely to reach a $(j+1)$-advancing or $(j+1)$-filling state before failing, and one of these is likely to happen in $\poly(j) W_j$ steps.
Lemma~\ref{lem:OnejVolumeAnalysisAdvancing}
proves a similar result, starting from a $j$-advancing high-level state.

 \begin{lemma}
\label{lem:OnejVolumeAnalysisAdvancing}
There exists a real number $\lambda_0 \in (0,1/2)$ such that the following holds. Consider any $\lambda \in (0,\lambda_0)$ and a send sequence $\p$ such that $p_0=1$ and, for all $i \ge 1$, $W_i \le (4/\lambda)^i$. Suppose that $\p$ has finitely many bins which are strongly exposed. 
Let $j_0$ be a sufficiently large positive integer
where  all strongly-exposed bins are in $[j_0-1]$.
Let $\tau_0$ be a non-negative integer. Let $\Psi_0$ be the   initialising high-level state with  $\tau^{\Psi_0} = \tau_0$ and $j^{\Psi_0} = j_0$. 
Let $\hls$ be the $\Psi_0$-backoff bounding rule from Definition~\ref{def:backoff-bounding-rule}.  
Let $V$ be a volume process from $\Psi_0$ with transition rule $\hls$, 
 send sequence~$\p$, and birth rate~$\lambda$  (Definition~\ref{def:volume-process}).
For any
$t \ge \tau_0$, 
let   $T(t) = \{\tau_0,\ldots,t\}$.
Then  the following statement holds  for
every $t > \tau_0$ and
every tuple
$ 
\psivect_{T(t-1)} = (\psi_{\tau_0},\dots,\psi_{t-1})$  of high-level states 
such that  
$\Pr( \Psivect_{T(t-1)} = \psivect_{T(t-1)}) > 0$, 
$\psi_{t-1}$ is a $j$-advancing state, and 
$\tau^{\psi_{t-1}} = t-1$.
Conditioned on $\Psivect_{T(t-1)} = \psivect_{T(t-1)} $,
with 
probability at least $1-\exp(-\lambda j/500)
 $,  
the following events occur, where  
\[
\tau_{j,\full} = \min\{t'  \geq t \colon \type^{\Psi(t')} = \Failure \mbox{ or }j^{\Psi(t')} = j+1\}.\]
\begin{enumerate}[(i)]
\item $\Psi(\tau_{j,\full})$ is $(j+1)$-filling or $(j+1)$-advancing
and $z_{j+1}^{\Psi(\tau_{j,\full})}=0$.
\item $\tau_{j,\full} \leq t+ 12 j^{\phi+1} W_j 
$.
\end{enumerate}

\end{lemma}
\begin{proof}
Let  $A = 12 j^{\phi+1} W_j$.
By the backoff bounding rule (Definition~\ref{def:backoff-bounding-rule}) and
Lemma~\ref{lem:advancing-transient}, $\tau_{j,\full }\leq t+A$, which implies (ii). 
Also, $\tau_{j,\full}$ is the first time $t'\geq t$ such that the transition $\Psi(t'-1)\rightarrow\Psi(t')$ is not via (R\ref{item:R4})(ii). Such a transition must be one of (R\ref{item:R1})(ii), (R\ref{item:R1})(iii), (R\ref{item:R2})(i), and (R\ref{item:R4})(i).

For each of these transitions that does not lead to a $(j+1)$-advancing or $(j+1)$-filling 
state with $z_{j+1}^{\Psi(\tau_{j,\full})}=0$
we will show that  
it is unlikely that this is the first transition 
during steps $t,\ldots,t+A$
that is not via (R\ref{item:R4})(ii).

 \begin{itemize}
 \item
{(R\ref{item:R1})(ii):} By Observation~\ref{obs:no-R1-from-adv-fill-stab}, such a transition cannot occur.

 \item
{(R\ref{item:R1})(iii):} 
By Lemma~\ref{lem:volume-easyChernoffs} the probability of such a transition into any given~$\Psi(t')$ is at most $\exp(-j^{1.9})$. The probability that this ever happens for $t'\in \{t,\ldots, t+A\}$ is therefore, by a union bound, at most  $(A+1) \exp(-j^{1.9})$. By hypothesis, $W_j \le (4/\lambda)^j$ and $j\geq j_0$ is sufficiently large so this is at most $\exp(-j)$.

 \item
{(R\ref{item:R2})(i):} 
There are two cases. If bin $j$ is heavy-covered then 
by Lemma~\ref{lem:AdvancingNoiseFail} the probability of such a transition into any given~$\Psi(t')$ is at most $\exp(-j^{1.9})$. 
Similar to the argument for (R\ref{item:R1})(iii), 
taking a union bound over all~$t'$ gives a failure probability that is at most $\exp(-j)$.

If $j$ is many-covered then 
by Lemma~\ref{lem:AdvancingNoiseFail}
the probability of such a transition into any given~$\Psi(t')$ is at most
$\exp(-  \lambda |\Upsilon_{j,\geq \Wtilde[j]}| \, \Wtilde[j] /80)$. By a union bound the probability that this ever happens is at most 
\[(A+1) \exp(-  \lambda |\Upsilon_{j,\geq \Wtilde[j]}| \, \Wtilde[j] /80)
\leq 13 j^{\phi+1}W_j \exp(-  \lambda |\Upsilon_{j,\geq \Wtilde[j]}| \, \Wtilde[j] /80).\]
By the definition of many-covered (Definition~\ref{def:covered}), (Prop~\ref{cov-prop-1}) holds so 
\[ 13 j j^{\phi} W_j \le (13/\CFill) j
\exp(\CNoiseBelow\lambda\Wtilde[j]\,|\Upsilon_{j,\ge\Wtilde[j]}|) \] 
where (from Definition~\ref{def:constants}), $\CFill = 2$ and
$\CNoiseBelow = 1/(30 \times 320)$.
The failure probability is therefore at most
\[(13/\CFill) j
\exp(-\lambda\Wtilde[j]\,|\Upsilon_{j,\ge\Wtilde[j]}| (1/80-\CNoiseBelow)) \]
Since $\Wtilde[j]\,|\Upsilon_{j,\ge\Wtilde[j]}|\geq j-1$
(from Observation~\ref{obs:Wtilde}) 
and $(1/80 - \CNoiseBelow) \geq 1/90$
the failure probability is at most 
\[(13/2) j
\exp(-\lambda (j-1)/90  ) \leq \exp(-\lambda j/100).\]

\item
{(R\ref{item:R4})(i):}
Note that if (R\ref{item:R4})(i) results in a $(j+1)$-advancing or $(j+1)$-filling state $\Psi(t')$ then $z_{j+1}^{\Psi(t')}=0$ -- this is by definition for an advancing state and is from the backoff bounding rule for a filling state. Thus, we just need to upper-bound the probability that   $\type^{\Psi(t')}= \Failure$.
We do this as follows.
Let $\tmin = \min\{t'' \geq t -1 + j^{24} \colon  
F^{\psi_{t-1}}_j(t'') \ge \mu_j^{\gammavect}\}$. 
If the first transition that is not via (R\ref{item:R4})(ii) is via (R\ref{item:R4})(i) then (from the backoff bounding rule) it will be the transition from 
$\Psi(t'-1)$ where $t'=\tmin$.  
If bin~$j+1$ is many-covered then by Lemma~\ref{lem:AdvancingIncj-NoiseFail}(i) the probability that this transition goes to a failure state (rather than to a $(j+1)$-advancing state) is at most 
$\exp(-  \lambda |\Upsilon_{j+1,\geq \Wtilde[j+1]}| \, \Wtilde[j+1] /80)$.
Similarly, if bin~$j+1$ is heavy-covered then by 
Lemma~\ref{lem:AdvancingIncj-NoiseFail}(ii) the probability that this transition goes to a failure state  is at most 
$\exp(-  (j+1)^{1.9})$.
If bin~$j+1$ is exposed then by Lemma~\ref{lem:T3tofill},
the probability that this transition goes to a failure state (rather than to a $(j+1)$-filling state) is at most 
$\exp(- \lambda(j-1)/400)$ (which is the largest of the three failure probabilities in the three cases where (R\ref{item:R4})(i) applies, since $j\geq j_0$ is sufficiently large).
\end{itemize} 

The result follows by a union bound, 
since the sum of the failure probabilities for (R\ref{item:R1})(iii), (R\ref{item:R2})(i) and (R\ref{item:R4})(i) is at most $\exp(-\lambda j/500)$.
\end{proof}

We can now prove Lemma~\ref{lem:VolumeAnalysis}, the main result of this section.

\begin{lemma}
\label{lem:VolumeAnalysis}
There exists a real number $\lambda_0 \in (0,1/2)$ such that the following holds. Consider any $\lambda \in (0,\lambda_0)$ and a send sequence $\p$ such that $p_0=1$ and, for all $i \ge 1$, $W_i \le (4/\lambda)^i$. Suppose that $\p$ has finitely many bins which are strongly exposed.  
Let $j_0$ be a sufficiently large positive integer
where  all strongly-exposed bins are in $[j_0-1]$.
Let $\tau_0$ be a non-negative integer. Let $\Psi_0$ be the   initialising high-level state with  $\tau^{\Psi_0} = \tau_0$ and $j^{\Psi_0} = j_0$. 
Let $\hls$ be the $\Psi_0$-backoff bounding rule from Definition~\ref{def:backoff-bounding-rule} and let $\OurHLS$  (Definition~\ref{def:backoff-bounding-state-space}) be the  $\Psi_0$-closure (Definition~\ref{def:hls-closed}) of~$\hls$. The $\Psi_0$-backoff bounding rule~$\hls$ is a transition rule (Definition~\ref{def:transition-rule}) by Observation~\ref{obs:BB-bounding-is-transrule} and 
it is $\OurHLS$-valid (Definition~\ref{def:valid}) by
Lemma~\ref{lem:backoff-bounding-valid}.

 Let $V$ be a volume process from $\Psi_0$ with transition rule $\hls$, 
 send sequence~$\p$, and birth rate~$\lambda$  (Definition~\ref{def:volume-process}). 
With probability at least~$99/100$, both of the following hold.
\begin{itemize}
\item for all $t \ge \tau_0$, $\type^{\Psi(t)} \ne \Failure$. 
\item for all $j \ge j_0$, there is
a positive integer $t \le \tau_0 + 1+ {(8 j)}^{\phi+1}(W_{j_0}+\dots + W_{j-1})$ such that 
$\Psi(t)$ is $j$-filling or $j$-advancing.
\end{itemize}
\end{lemma}
\begin{proof} 
To make the proof more succinct, we use the terminology ``$j$-good state'' for a high-level state that is either 
$j$-advancing or
$j$-filling with $z_j=0$.

From the definition of initialising state (Definition~\ref{def:init-state}),
$\Psi_0 = (0,\tau_0,j_0,\zvect,\calS,\initialising)$ where
$\calS = \emptyset$.  
From the definition of the volume process (Definition~\ref{def:volume-process}), the first transition is given by
$\Psi(\tau_0+1) = \hls(\Psi_0,\bbar_{[j_0]}^{Y_{\tau_0}}(\tau_0+1),\tau_0+1)$.
From the backoff bounding rule (Definition~\ref{def:backoff-bounding-rule}),  using the fact that $\calS=\emptyset$,
the only possibilities for this transition  are   (R\ref{item:R1})(ii), (R\ref{item:R1})(iii), or (R\ref{item:R3}). (R\ref{item:R1})(ii) is ruled out by Observation~\ref{obs:no-R1-from-adv-fill-stab}.
By Lemma~\ref{lem:volume-easyChernoffs}, the probability that the transition is via (R\ref{item:R1})(iii) is at most $\exp(-j_0^{1.9})$. 

A transition via (R\ref{item:R3}) either leads to a $j_0$-good state or to a failure state.  
There are three cases.
\begin{itemize}
\item If bin~$j_0$ is heavy-covered: 
By Lemma~\ref{lem:AdvancingNoiseFail} the probability  that the transition leads to a failure state is at most $\exp(-{j_0}^{1.9})$.

\item If bin~$j_0$ is many-covered:
By Lemma~\ref{lem:AdvancingNoiseFail}
the probability 
that the transition leads to a failure state is at most 
$\exp(-  \lambda |\Upsilon_{j_0,\geq \Wtilde[j_0]}| \, \Wtilde[j_0] /80)$.
Since $\Wtilde[j_0]\,|\Upsilon_{j_0,\ge\Wtilde[j_0]}|\geq j_0-1$
(from Observation~\ref{obs:Wtilde}), this is at most
$\exp(- \lambda (j_0-1)/80)$.

\item If bin~$j_0$ is exposed: By Lemma~\ref{lem:T3tofill}
the probability 
that the transition leads to a failure state is at most 
$\exp(- \lambda(j_0-1)/400)$. Since $j_0$ is sufficiently large, this is the largest failure probability for the (R\ref{item:R3}) cases.
\end{itemize}
Since the sum of
the failure probabilities for (R\ref{item:R1})(iii) and (R\ref{item:R3}) is at most $-\exp(-\lambda j_0/500)$, 
a union bound tells us that the probability that $\Psi(\tau_0+1)$ is a $j_0$-good state is at least
$1- \exp(-\lambda j_0/500)$.
For any
$t \ge \tau_0$, 
let   $T(t) = \{\tau_0,\ldots,t\}$.
We now combine Lemmas~\ref{lem:OnejVolumeAnalysis} and \ref{lem:OnejVolumeAnalysisAdvancing} to get the following claim, for any $j\geq j_0$.

\noindent {\bf Claim:\quad} 
For any $t> \tau_0$ 
and any vector $ 
\psivect_{T(t-1)} = (\psi_{\tau_0},\dots,\psi_{t-1})$  of high-level states 
such that  
$\Pr( \Psivect_{T(t-1)} = \psivect_{T(t-1)}) > 0$
and
$\psi_{t-1}$ is a $j$-good state with $\tau^{\psi_{t-1}}=t-1$, 
conditioned on $\Psivect_{T(t-1)} = \psivect_{T(t-1)} $,
with 
probability at least $1-
\exp(-\lambda j/1000) $, the process reaches a $(j+1)$-good state by time 
$t+ (8j)^{\phi+1}W_{j}$ (without entering any failure states along the way).
 
Lemma~\ref{lem:OnejVolumeAnalysis} establishes the claim for the case that the $j$-good state is $j$-filling. Here we use the fact that $\lceil 4(4 j)^{\phi} W_j\rceil
\leq (8j)^{\phi+1}W_{j}$. The failure probability is at most $\exp(-\lambda j/1000)$.  
Lemma~\ref{lem:OnejVolumeAnalysisAdvancing} establishes the claim for the case that the $j$-good state is $j$-advancing. Here we use the fact that $ 12 j^{\phi+1} W_j
\leq (8j)^{\phi+1}W_{j}$. The failure probability is at most $\exp(-\lambda j/500) \leq \exp(-\lambda j/1000)$.

To prove the lemma we just need to
sum the failure probabilities, 
noting that
\[ \exp(- \lambda j_0/500) + \sum_{j\geq j_0} \exp(- \lambda j/1000) \leq 1/100,\] which holds as long as $j_0$ is sufficiently large (with respect to~$\lambda$).   
\end{proof}

\begin{corollary}
\label{cor:VEB-volume}
There exists a real number $\lambda_0 \in (0,1/2)$ such that the following holds. Consider any $\lambda \in (0,\lambda_0)$ and a send sequence $\p$ such $p_0=1$ and, for all $i \ge 1$, $W_i \le (4/\lambda)^i$. Suppose that $\p$ has finitely many bins which are strongly exposed. 
Let $X$ be the backoff process with
birth rate~$2 \lambda$, send sequence~$\p$, and cohort set~$\calC=\{A,B\}$ (Definition~\ref{def:backoff-cohorts}).
Let $j_0\geq 5$ be a sufficiently large positive integer
where  all strongly-exposed bins are in $[j_0-1]$.
Let $\Einit$ be any event determined by 
$\{n^{X^B}(t') : t'\in [T^{\p,j_0}]\}$
which entails that each 
$n^{X^B}(t')$ in this set is at least~$1$.
We will condition on $\Einit$.

Let
$\tau_0 =  T^{\p,j_0}$. 
Let
$\Psi_0$ be the  initialising high-level state   
with $\tau^{\Psi_0} = \tau_0$ and $j^{\Psi_0} = j_0$. 
 
Let $\hls$ be   the $\Psi_0$-backoff bounding rule from Definition~\ref{def:backoff-bounding-rule} and let $\OurHLS$  (Definition~\ref{def:backoff-bounding-state-space}) be the  $\Psi_0$-closure (Definition~\ref{def:hls-closed}) of~$\hls$.  
The $\Psi_0$-backoff bounding rule~$\hls$ is a transition rule (Definition~\ref{def:transition-rule}) by Observation~\ref{obs:BB-bounding-is-transrule} and 
it is $\OurHLS$-valid (Definition~\ref{def:valid}) by
Lemma~\ref{lem:backoff-bounding-valid}.

Conditioned on $\Einit$, the VEB-coupling (Definition~\ref{def:VEB-couple})   couples $X$ with a volume process $V$ from $\Psi(\tau_0)=\Psi_0$ with  transition rule $\hls$, send sequence~$\p$ and birth rate~$\lambda$
and with a sequence $E_1, E_2,\ldots $ of escape processes with send sequence~$\p$ and birth rate~$\lambda$. 
With probability at least~$99/100$,  
both of the following hold.

\begin{itemize}
\item for all $t \ge \tau_0$, $\type^{\Psi(t)} \ne \Failure$. 
\item for all $j \ge j_0$, there is
a positive integer $t \le \tau_0 + 1+ {(8j)}^{\phi+1}(W_{j_0}+\dots + W_{j-1})$ such that 
$\Psi(t)$ is $j$-filling or $j$-advancing.
\end{itemize}
\end{corollary}
\begin{proof}
Immediate from Lemma~\ref{lem:VolumeAnalysis}
and the definition of the VEB coupling (Definition~\ref{def:VEB-couple}).
\end{proof}

\section{Analysis of escape processes}\label{sec:escape-analysis}

The goal of this section is to prove Lemma~\ref{lem:VEB-escape}, which says that with constant probability the escape processes in a VEB coupling never fill up. Our approach will be as follows. 

Recall from Definition~\ref{def:escape} that an escape process $E$ takes a fixed high-level state $\Psi$ as a parameter, and that the expected number $\overline{\romanarr}(t)$ of end-of-step arrivals in $E$ at time $t$ is a function $\nu^E(t)$ of the vectors $\bbar_S^E(t)$ with $S \in \calS^\Psi$. This lack of independence would complicate the analysis substantially, so in Section~\ref{sec:constant-escape} we first introduce a ``constant-escape process'' whose end-of-step arrival rate $\nu$ is constant. As one might expect, an escape process can be dominated above by a constant-escape process for as long as $\nu^E(t) \le \nu$; we give this coupling in Definition~\ref{def:constant-escape-coupling}. The expected number of balls in any bin $k$ of a constant-escape process at time $t$ is then a simple function of its initial state, which we denote by $h_k(t)$ in Definition~\ref{def:h}. The role of $h$ in analysing constant-escape processes is analogous to the role of $f$ in analysing $j$-jammed processes. Like $f$, $h$ is determined by a simple recurrence relation, and we end Section~\ref{sec:constant-escape} by proving some useful upper bounds on $h_k(t)$. We then extend these properties into Chernoff-based concentration bounds on $b_k^E(t)$.

Recall from Definition~\ref{def:VEB-couple} that high-level state transitions in the VEB coupling are not determined by the  state of the (current) escape process -- instead, the high-level states evolve together with the volume process. For this reason, it is difficult to maintain any kind of Poisson domination of the states of the escape process -- from the point of view of any escape process, the transitions happen at unpredictable times. For this reason, when the high-level state changes,
the VEB coupling passes from one escape process to the next, exposing 
the position of every ball in the initial state of the new escape process. 
The analysis of the escape process has to be able to deal with these exposed states. Moreover, different high-level states contribute
different end-of-step arrival rates to the escape processes.  For example, an escape process in a filling state will gradually empty to a low level even if it is initially full, while an escape process in a refilling state will fill rapidly to a high level even if it is initially empty. We manage this complexity with a set of parameterised pseudorandomness conditions (invariants) which we set out in Section~\ref{sec:escape-pseudorandomness}; we will use these to show that while we cannot maintain a Poisson domination for escape processes, we can nevertheless recover key desirable properties of a Poisson domination. Looked at another way, these conditions will enable a complex induction invariant on the sequence of escape processes in the VEB coupling, where each step of the induction will hold with summable failure probability. 

In Section~\ref{sec:constant-escape-invariants} we show how the invariants of Section~\ref{sec:escape-pseudorandomness} evolve over the course of a constant-escape process; in Section~\ref{sec:escape-invariants} we apply the results of Section~\ref{sec:constant-escape-invariants} to show how these invariants evolve over the course of an escape process; in Section~\ref{sec:veb-transition-invariants} we transfer the results of Section~\ref{sec:escape-invariants} to the context of the VEB coupling; and finally in Section~\ref{sec:escape-analysis-final} we apply the results of Section~\ref{sec:veb-transition-invariants} to prove Lemma~\ref{lem:VEB-escape} by induction.

\subsection{Constant-escape processes}\label{sec:constant-escape}

We first define constant-escape processes and show (in Definition~\ref{def:constant-escape-coupling}) how they can dominate escape processes from above.

\begin{definition}
\label{def:constant-escape}
Fix a send sequence~$\p$ with $p_0=1$,   
a positive integer~$j$, a non-negative integer~$\tau$,
and a real number $\nu \ge 0$. 
A \emph{constant-escape process} \newsym{E+}{constant-escape process}{E^+} on $[j]$ with 
send sequence~$\p$,
start time $\tau$, and end-of-step arrival rate \newsym{nu}{end-of-step arrival rate of constant-escape process}{\nu} is a generalised backoff process with send sequence $\p$, 
start time~$\tau$, the birth distribution~$\calD$ that always returns~$0$, and an initial population $\arrvect^{E^+}(\tau)$, that satisfies the following constraints.
\begin{itemize}
\item For all $i>j$ and $t \ge \tau$, $\romanarr_i^{E^+}(t) = 0$.
\item  For all $i\in [j]$ and $t>\tau$,
$\romanarr_i^{E^+}(t)$ is a Bernoulli r.v.\ with parameter~$\nu$. These are independent for different~$i$ and for different~$t$.
\item For all $i \in [j-1]$ and $t \ge \tau$, $\e_i^{E^+}(t) = \emptyset$. 
\item For all $i \ge j$ and $t \ge \tau$, $\e_i^{E^+}(t) = \s_i^{E^+}(t)$.
\end{itemize}
\end{definition}

\begin{definition}\label{def:constant-escape-coupling}
Fix a send sequence~$\p$ with $p_0=1$, $\lambda\in (0,1)$, 
and a high-level state $\Psi=(g,\tau',j,\zvect,\calS,\type)$.
Fix a (possibly random) tuple $\arrvect^E(\tau')$ of
disjoint sets
$\arrvect^E(\tau') = (\arr_1^E(\tau'),\arr_2^E(\tau'),\ldots)$
where each $ \arr^E_{i}(\tau')\subseteq \ball(\tau')$ and, 
for $i>j$, $\arr_i^E(\tau') = \emptyset$. 
Let $E$ be an escape process
with 
high-level state $\Psi$, send sequence~$\p$, birth rate~$\lambda$, and
initial population $\arrvect^E(\tau')$.
Fix an integer  $\tau \ge \tau'$ and 
a real number $\nu \ge 0$. 
Let $\newsym{T-E-nu}{Decoupling time of a constant-escape process with start time $\tau$ and arrival rate $\nu$ from its parent escape process $E$}{T^E_{\nu,\tau}}$ be the first time $t \ge \tau$ such that $\nu^E(t) > \nu$. 
For all $i\geq 1$, let $\arr_i^{E^+}(\tau) = \b_i^E(\tau)$.
There is a constant-escape process~$E^+$ on~$[j]$ with
send sequence~$\p$, start time~$\tau$, end-of-step arrival rate~$\nu$  
and and initial population $\arrvect^{E^+}(\tau)$
such that there is a coupling satisfying the following invariants for every $i\geq 1$.
\begin{description}  
\item  [Inv-$\b_i(t)$:] For all $t$ with $\tau \le t \le T_{\nu,\tau}^E$, $\b_i^{E}(t) \subseteq \b_i^{E^+}(t)$.
\item  [Inv-$\B_i(t)$:] For all $t$ with $\tau+1 \le t \le T_{\nu,\tau}^E$, $\B_i^{E}(t) \subseteq \B_i^{E^+}(t)$.
\item [Inv-$\s_i(t)$:] For all $t$ with $\tau+1 \le t \le T_{\nu,\tau}^E$, $\s_i^{E}(t) \subseteq \s_i^{E^+}(t)$.
\end{description}
We refer to this coupling as the \emph{constant-escape coupling} of~$E$ with~$E^+$. 
The coupling is defined as follows. By construction, for every $i\geq 1$,
$\b_i^{E}(\tau) = \b_i^{E^+}(\tau)$ (these are empty for $i>j$). Thus, 
$\mathrm{Inv}$-$\b_i(\tau)$ is satisfied for every $i\geq 1$.

For integer~$t>\tau$, the coupling evolves as follows.
If $t > T^E_{\nu,\tau}$ then $E$ and $E^+$ evolve independently at time $t$.
If $t\leq T^E_{\nu,\tau}$ then
$\nu^E(t-1) \leq \nu$ and the coupling proceeds as follows.
For all $i\geq 0$, 
$\B_i^{E^+}(t)= \b_i^{E^+}(t-1)$ and 
$\B_i^{E}(t) = \b_i^{E}(t-1)$. By $\mathrm{Inv}$-$\b_i(t-1)$, we get $\mathrm{Inv}$-$\B_i(t)$.
For every $i\geq 0$, $\s^{E}_i(t)$ is determined 
by taking each element of $\B_i^{E}(t)$ independently with probability~$p_i$.
We define $\S^{E^+}_i(t)$ to be the random set formed by including each ball in 
$\B_i^{E^+}(t) \setminus \B_i^{E}(t)$ 
independently with probability~$p_i$, then take $\s^{E^+}_i(t) = \S^{E^+}_i(t) \cup \s^{E}_i(t)$ to establish $\mathrm{Inv}$-$\s_i(t)$.
By the definitions, for all $i\in [j-1]$, $\e_i^{E}(t) = \e_i^{E^+}(t) = \emptyset$ and for all $i\geq j$, $\e_i^{E^+}(t) = \s_i^{E^+}(t)$ and
$\e_i^E(t) = \s_i^E(t)$.

We next define $\overline{\arr}^E(t)$ and $\overline{\arr}^{E^+}(t)$. For each $i\in [j]$ independently,
\begin{itemize}
\item  with probability $\nu^E(t-1)$,
$\romanarr_i^E(t) = \romanarr_i^{E^+}(t) = 1$.
In this case, let $\arr_i^{E^+}(t) = \arr_i^E(t)$.
\item with probability $\nu - \nu^E(t-1)$,
$\arr_i^E(t)=\emptyset$ and $\romanarr_i^{E^+}(t)=1$.
\item with probability $1-\nu$, 
$\arr_i^E(t) = \arr_i^{E^+}(t) = \emptyset$.
\end{itemize}

The invariant $\mathrm{Inv}$-$\b_i(t)$ is established for all $i\geq 0$ as follows.  Fix $i$ and consider   a ball $\beta \in \b_i^{E}(t)$. We wish to show $\beta \in \b_i^{E^+}(t)$. Recall from the definition of a generalised backoff process (Definition~\ref{def:genbackoff}) that $\b_i^E(t) = \arr_i^E(t) \cup (\B_i^E(t) \setminus \s_i^E(t)) \cup (\s^{E}_{i-1}(t) \setminus \e^{E}_{i-1}(t))$; we split into cases accordingly.

\begin{itemize}
\item $\beta\in \arr_i^{E}(t)$: By construction this implies $\arr_i^{E^+}(t) = \arr_i^E(t)$, so $\beta \in \b_i^E(t)$.
\item   $\beta \in \B_i^{E}(t) \setminus \s_i^{E}(t)$: By $\mathrm{Inv}$-$\B_i(t)$, $\beta \in \B_i^{E^+}(t)$, and by construction, $\beta \notin \s_i^{E^+}(t)$. Hence $\beta \in \B_i^{E^+}(t) \setminus \s_i^{E^+}(t)$, so $\beta \in \b_i^{E^+}(t)$.
\item $\beta \in \s^{E}_{i-1}(t) \setminus \e^{E}_{i-1}(t)$:
By 
$\mathrm{Inv}$-$\s_{i-1}(t)$, $\beta \in \s^{E^+}_{i-1}(t)$. 
Since $\beta \notin \e^E_{i-1}(t)$, $i-1<j$,
so $\e^{E^+}_{i-1}(t)=\e^{E^+}_{i-1}(t) = \emptyset$ and in particular $\beta \notin \e^{E^+}_{i-1}(t)$. Thus 
$\beta \in \s^{E^+}_{i-1}(t) \setminus \e^{E^+}_{i-1}(t)$ which implies
$\beta \in \b_i^{E^+}(t)$.
\end{itemize}
\end{definition}

As we will see in Observation~\ref{obs:constant-escape}, the evolution of expected ball counts in a constant-escape process is governed by the following functions, which are very similar to the $f$'s of Definition~\ref{def:f}.

\begin{definition}\label{def:h}
Let $\p$ be a send sequence with $p_0=1$. Fix a positive integer~$j$, a tuple $\newsym{a-vect}{Tuple of non-negative integers, typically initial ball counts}{\avect} \in (\reals_{\ge 0})^j$, and a real number $\nu \ge 0$. We now define a function \newsym{h-nu-a}{function to track population means in a constant-escape process on $[j]$ with initial state~$\avect$ (Definition~\ref{def:h})}{\hvect^{\nu,\avect}} from non-negative integers to tuples of non-negative reals indexed by $\{0,\dots,j\}$. $h_x^{\nu,\avect}$ denotes the $x$'th component of $\hvect^{\nu,\avect}$.
For all $x\in [j]$, let $h_x^{\nu,\avect}(0) = a_x$.  For all $t \ge 0$, let $h_0^{\nu, \avect}(t) = 0$. For all $t \ge 1$ and all $x \in [j]$, let
\begin{equation}\label{eq:h-def}
h_x^{\nu, \avect}(t) = (1-p_x)h_x^{\nu, \avect}(t-1) + p_{x-1}h_{x-1}^{\nu, \avect}(t-1) + \nu.
\end{equation}
\end{definition}

Observation~\ref{obs:h-nu} will be useful.
We will also use Observation~\ref{obs:constant-escape}, whose proof is similar to   the proof of Observation~\ref{obs:f-gen}.

\begin{observation}\label{obs:h-nu}
Let $\p$ be a send sequence with $p_0=1$. Fix a positive integer~$j$, a tuple $\newsym{a-vect}{Tuple of non-negative integers, typically initial ball counts}{\avect} \in (\reals_{\ge 0})^j$, and 
real numbers $0 \leq \nu^{-} \leq \nu^{+}$.
Then for all $t\geq 0$ and all $x\in [j]$,
$h_x^{\nu^{-},\avect}(t) \leq
h_x^{\nu^{+},\avect}(t)$.
\end{observation}
\begin{proof}
The proof is by induciton on~$t$ using the definition of~$h$. The base case is $t=0$
since $h_x^{\nu^{-},\avect}(0) = a_x = 
h_x^{\nu^{+},\avect}(0)$.
\end{proof}

\begin{observation} 
\label{obs:constant-escape}
 Fix a send sequence~$\p$ with $p_0=1$. Fix integers $j\ge 1$ and $\tau \ge 0$ and a real number $\nu\ge 0$. 
 Fix $\mvect \in (\reals_{\geq 0})^j$.
 For each $i\in [j]$. let
 $A_i$ be a random variable  
 over $\integers_{\geq 0}$
 with mean~$m_i$
 and let $\Avect=(A_1,\ldots,A_j)$. We don't make any assumptions about the joint distribution of $A_1,\ldots,A_j$.
Let~$E^+$ be a constant-escape process on $[j]$ with send sequence~$\p$, start time $\tau$, end-of-step arrival rate $\nu$, and $\romanarrvect_{[j]}^{E^+}(\tau) = \Avect$.   
Then for all $i \in [j]$
and all $t\geq \tau$,
$\E(b_i^{E^+}(t)) = h_i^{\nu,\mvect}(t-\tau)$.
\end{observation}
\begin{proof}

The proof is by induction on $t$.
For the base case, $t=\tau$, 
for every $i\in [j]$,
$\E(b_i^{E^+}(\tau)) = 
\E[\romanarr_i^{E^+}(\tau)] = 
\E[A_i] = m_i$.
For the
inductive step consider $t>\tau$.
Then from the definition of generalised backoff process (Definition~\ref{def:genbackoff}
and the definition of constant escape process (Definition~\ref{def:constant-escape}),
\begin{align*}
\E[b_i^{E^+})(t)] &= 
\E[\romanarr_i^{E^+}(t)] +
\E[B_i^{E^+}(t)\setminus \s_i^{E^+}(t)] 
+ \E[s_{i-1}^{E^+}(t)]\\
&= 
 \nu +
\E[b_i^{E^+}(t-1) ] (1-p_i)
+ \E[b_{i-1}^{E^+}(t-1)]p_{i-1}.
\end{align*}
By the induction hypothesis, right-hand side is 
$ \nu + h_i^{\nu,\mvect}(t-1-\tau)  (1-p_i)
+ h_{i-1}^{\nu,\mvect}(t-1-\tau)
p_{i-1}$. By Equation~\eqref{eq:h-def}, this
is $h_i^{\nu,\mvect}(t-\tau)$.
\end{proof}

Note that if the vector $\mvect$
from Observation~\ref{obs:constant-escape} is in $(\integers_{\geq 0})^j$ then the observation also applies to the case where 
the random variable $A_i$ is deterministically equal to~$m_i$.
Recall the definition of generalised multinomial distribution from Definition~\ref{def:gen-multi}.
Lemma~\ref{lem:constant-escape-gen-multi} shows that this distribution is relevant to the distribution of bin populations of a constant escape process.

\begin{lemma} 
\label{lem:constant-escape-gen-multi}
 Fix a send sequence~$\p$ with $p_0=1$. Fix integers $j\ge 1$ and $\tau \ge 0$ and a real number $\nu\ge 0$. 
Let $R_1,\ldots,R_j$ be independent binomial random variables
where, for all $i\in [j]$, $R_i = {\mathrm Bin}(n_i,\rho_i)$.
Let $\Rvect = (R_1,\ldots,R_j)$.
Let~$E^+$ be a constant-escape process on $[j]$ with send sequence~$\p$, start time $\tau$, end-of-step arrival rate $\nu$, and $\romanarrvect_{[j]}^{E^+}(\tau) = \Rvect$.   
Then for all $t \ge \tau$, $b_1^{E^+}(t),\dots,b_j^{E^+}(t)$ can be expressed as coordinates of a generalised multinomial random variable. 
\end{lemma}
\begin{proof}
Fix $t\geq \tau$.
Let $n= \sum_{i=1}^j n_i$
and let $m= n+ j(t-\tau)$.
Let $\beta_1,\dots,\beta_{n}$ be the set of all balls which might
potentially lie in  
$\cup_{i\in [j]} \arr_i^{E^+}(\tau)$   
and let $\beta_{n+1},\dots,\beta_m$ be the set of all balls which might potentially lie in 
$\cup_{i\in [j]} \arr_{i}^{E^+}(t')$ for some $t'\in \{\tau+1,\ldots, t\}$.   
For all $i \in [m]$, let $A_i = j+1$ if $\beta_i \notin \bigcup_{k=1}^j \b_k^{E^+}(t)$, and otherwise let $A_i$ be the unique $x$ such that $\beta_i \in \b_x^{E^+}(t)$. Then the variables $A_1,\dots,A_m$ are independent, and for all $k \in [j]$, $b_k^{E^+}(t) = |\{i \in [m] \colon A_i = k\}$. Thus $(b_1^{E^+}(t), \dots, b_j^{E^+}(t), A_{j+1})$ is a generalised multinomial distribution by 
    Definition~\ref{def:gen-multi}.
\end{proof}

Clearly Lemma~\ref{lem:constant-escape-gen-multi}
also applies to the deterministic case where each $\rho_i=1$ so that $R_i$ is deterministically equal to~$n_i$.
We next note the stationary distribution of a constant-escape process. While we will not use this directly, it forms important motivation for the lemmas that follow.

\begin{definition}\label{def:kappa-x-nu}
Fix a send sequence~$\p$ with $p_0=1$, a positive integer $j$,
and a real number $\nu \ge 0$.
For all $x \in [j]$, define $\newsym{kappa-x-nu}{$x\nu W_x$ (with respect to a fixed send sequence)}{\kMean_{x,\nu}} = x\nu W_x$.
\end{definition}

\begin{remark} \label{rem:escapemean}
Fix a send sequence~$\p$ with $p_0=1$, integers $j \ge 1$ and $\tau \ge 0$,
and a real number $\nu \ge 0$. Observe that $\overline{\kMean}(\nu) := (\kMean_{1,\nu},\ldots,\kMean_{j,\nu})$ is a fixed point in the recurrence of~\eqref{eq:h-def}.
Thus a constant-escape process~$E^+$ on~$[j]$ with 
send sequence~$\p$, start time~$\tau$, end-of-step arrival rate~$\nu$ and $
\romanarrvect_{[j]}^{E^+}(\tau)
\sim \Po(\overline{\kMean}(\nu))$ 
has  
$\bbar_{[j]}^{E^+}(t) \sim \Po(\overline{\kMean}(\nu))$ for every $t\geq \tau$.
\end{remark}

As one might expect, if the
expectation of bin populations of the initial state of a constant-escape process is at most stationary then this situation will persist as the process evolves. The following lemma is a more general version of this statement, which applies even if some bins are initially above the stationary distribution in expectation.

\begin{lemma}\label{lem:h-below-stationary} 
Let $\p$ be a send sequence with $p_0=1$. 
Let $j$ be a positive integer, let $\avect \in (\reals_{\ge 0})^j$,  let $\nu\ge 0$ be a real number, and let $T\ge 0$ be an integer. Then for all 
$t \ge T$ and all $x \in [j]$,
\[
p_xh_x^{\nu,\avect}(t) \le \max(
\{x \nu\}
\cup \{p_y h_y^{\nu,\avect}(T) + (x-y)\nu \colon y \in [x]\}).
\]
\end{lemma}
\begin{proof}
Throughout, we suppress $\nu$ and $\avect$ from $h_x^{\nu,\avect}(t)$. 
For all $x\in [j]$ and 
$\ell \in [x]$, let $\zeta_\ell = \max(
\{\ell \nu\}
\cup \{p_yh_y(T) + (\ell-y)\nu \colon y \in [\ell]\})$.
Thus we wish to show that for all $t\geq T$ and all $x\in [j]$, 
$p_x h_x(t) \leq  
\zeta_x$.

First we note that
\begin{align*} \zeta_{x-1} 
 &= 
\max( 
\{(x-1) \nu\}
\cup \{p_yh_y(T) + (x-1-y)\nu \colon y \in [x-1]\})\\
&= 
\max( 
\{x \nu\}
\cup \{p_yh_y(T) + (x-y)\nu \colon y \in [x-1]\}) - \nu\\
&\leq 
\max( 
\{x \nu\}
\cup \{p_yh_y(T) + (x-y)\nu \colon y \in [x]\}) - \nu = \zeta_x - \nu.  
\end{align*}

We now proceed by induction on $t$. For $t=T$, the result is immediate. Suppose the result holds for $t-1$ for some $t \ge T+1$ and fix $x \in [j]$.  
By the inductive  hypothesis
$p_x h_x(t-1) \leq \zeta_x$ and
 $p_{x-1} h_{x-1}(t-1) \leq 
\zeta_{x-1} \leq \zeta_x - \nu$.
By the definition of~$h$ (Definition~\ref{def:h}),
\begin{align*} 
p_xh_x(t) &= p_x\Big((1-p_x)h_x(t-1) + p_{x-1}h_{x-1}(t-1) + \nu\Big)\\
&= (1-p_x)p_xh_x(t-1) + p_x\big(p_{x-1}h_{x-1}(t-1)+\nu\big)
\le (1-p_x)\zeta_x + p_x\zeta_x
= \zeta_x,
\end{align*}
as required.
\end{proof}

Informally, Lemma~\ref{lem:h-soft-nonincreasing} establishes conditions such that $h_k^{\nu,\avect}$ is ``non-increasing in $t$ above some value~$M_k$''.

\begin{lemma}\label{lem:h-soft-nonincreasing}
Let $\p$ be a send sequence with $p_0=1$. 
Let $j$ be a positive integer, let $\avect \in (\reals_{\ge 0})^j$,  and let $\nu\ge 0$ be a real number. Let
$T$, $T'$, and $t$ be integers with 
$0\leq T \leq T' \leq t$.  
For any $k\in [j]$.
Let 
$M_k= \max(\{k \nu\}\cup \{p_y h_y^{\nu,\avect}(T) + (k-y)\nu\colon y \in [k-1]\})$. Then, for all $k\in [j]$, $h_k^{\nu,\avect}(t) \le \max\{h_k^{\nu,\avect}(T'), M_k W_k\}$.
\end{lemma}

\begin{proof}
We suppress $\nu$ and $\avect$ from the notation $h^{\nu,\avect}$ throughout. By Lemma~\ref{lem:h-below-stationary} applied with the $T$ of Lemma~\ref{lem:h-below-stationary} as our $T'$, for all $k\in [j]$,
$
p_kh_k(t) \le \max(\{k \nu\}\cup \{p_y h_y(T') + (k-y)\nu\colon y \in [k]\})$.
Since $k\nu \leq M_k$ 
and (taking $y=k$)
$p_y h_y(T') + (k-y) \nu = p_k h_k(T')$,
it suffices to show
that, for all $y\in [k-1]$,
$p_y  h_y(T') + (k-y) \nu \leq 
\max\{p_k h_k(T'),M_k\}$. This is vacuously true for $k=1$.
By Lemma~\ref{lem:h-below-stationary} applied again with 
the $T$ of 
Lemma~\ref{lem:h-below-stationary} as our~$T$, for all $k\in \{2,\ldots,j\}$ and $y \in [k-1]$,
\begin{align*}
p_yh_y(T') + (k-y)\nu &\le \max\Big(\{y \nu\} \cup \big\{p_\ell h_\ell(T) + (y-\ell)\nu\colon \ell \in [y]\big\}\Big) + (k-y)\nu\\
&= \max\Big(\{k \nu\} \cup \{p_\ell h_\ell(T) + (k-\ell)\nu\colon \ell \in [y]\}\Big) \leq M_k,
    \end{align*}
completing the proof. 
\end{proof}

Our next goal is to prove Lemma~\ref{lem:escape-low-weight-h-drops}, which further generalises Lemma~\ref{lem:h-below-stationary}. (We still state both lemmas because we will use Lemma~\ref{lem:h-below-stationary} in the proof of Lemma~\ref{lem:escape-low-weight-h-drops}.) Essentially, in the max term of Lemma~\ref{lem:h-below-stationary}, Lemma~\ref{lem:escape-low-weight-h-drops} allows us to mostly ignore a set $S$ of low-weight bins that begin above the stationary distribution by first waiting $O(j\sum_{k \in S}W_k)$ time steps. It is not immediately clear that this should be enough, as we should only expect a ball $\beta$ in bin $k \in S$ to escape the process altogether by sending from bin $j$ after $\Omega(\sum_{\ell =k}^j W_\ell)$ steps have passed. The key insight is that once $\beta$ has entered a high-weight bin, it is then unlikely to ever be in a specific low-weight bin on a specific time step, so its future contribution to $\overline{h}$ will be small. We formalise this idea in the following lemma.

\begin{lemma} \label{lem:escape-high-weight-kills-balls}
Let $\p$ be a send sequence with $p_0=1$. Let $j \ge 1$ and $\tau \ge 0$ be integers and let $\nu\ge 0$ be a real number. Let $E^+$ be a constant-escape process on $[j]$ with send sequence~$\p$, start time $\tau$, and end-of-step arrival rate $\nu$. Fix $t' \ge \tau$ and a $j$-tuple $\overline{\b}$ of ball sets such that $\pr(\b_{[j]}^{E^+}(t') = \overline{\b}) > 0$. Fix $\ell \in [j]$ and $\beta \in \b_\ell$. Then for all $t'' \ge t'$ and all $k \in [j]$,
$
    \pr(\beta \in \b_k^{E^+}(t'') \mid \b_{[j]}^{E^+}(t') = \overline{\b}) \le W_k/W_\ell$.
\end{lemma}
\begin{proof}
For all $x \in [j]$ and all $t'' \ge t'$, let 
$
g_x(t'') = \pr(\beta \in \b_x^{E^+}(t'') \mid \b_{[j]}^{E^+}(t') = \overline{\b})$.
Thus we must prove that for all $k \in \{\ell,\dots,j\}$ and all $t'' \ge t$, $g_k(t'') \le W_k/W_\ell$. We will prove this by induction on $t''$. 

For the base case of $t'' = t'$, $g_\ell(t') = 1 = W_\ell/W_\ell$ and, for all $k \in [j] \setminus \{\ell\}$, $g_k(t') = 0 < W_k/W_\ell$ as required. Now suppose the result holds up to $t''-1$ for some $t'' > t'$. Fix a ball $\beta \in \b_\ell$ and $x \in [j]$. By the definition of a constant-escape process (Definition~\ref{def:constant-escape}), $\beta$ lies in $\b^{E^+}_x(t'')$ if and only if it lies in $\b^{E^+}_x(t''-1) \setminus \s^{E^+}_x(t'')$ or $\s^{E^+}_{x-1}(t'')$; hence
$
    g_x(t'') = (1-p_x)g_x(t''-1) + p_{x-1}g_{x-1}(t''-1)$.
By the induction hypothesis, it follows that
\[
    g_x(t'') \le (1-p_x)W_x/W_\ell + p_{x-1}W_{x-1}/W_\ell = W_x/W_\ell,
\]
as required.
\end{proof}

We next use Lemma~\ref{lem:escape-high-weight-kills-balls} to prove a similar bound which applies even to balls starting in bins of low weight, by arguing that if we wait long enough then they are very likely to enter a bin of high weight (from which point Lemma~\ref{lem:escape-high-weight-kills-balls} will apply).

\begin{lemma}\label{lem:tool-2-single-ball}
Let $\p$ be a send sequence with $p_0=1$, and suppose $\lambda > 0$ is sufficiently small. Let $j\ge 1$ and $\tau \ge 0$ be integers and let $\nu\ge 0$ be a real number. Let $E^+$ be a constant-escape process on $[j]$ with send sequence~$\p$, start time $\tau$, and end-of-step arrival rate $\nu$. Consider a ball $\beta \in \cup_{i\in [j]}\b^{E^+}_i(\tau) $. Then, for all integers~$t$ satisfying $\tau < t \le \tau + (4/\lambda)^{2j}$ and all $k \in [j]$,  
\[
\pr(\beta \in b_k^{E^+}(t)) \le \frac{jW_k}{\lambda(t-\tau)}.
\]
\end{lemma}
\begin{proof}
\def\calW{\mathcal{W}}

Let $i$ be the bin such that $\beta \in \b_i^{E^+}(\tau)$, and fix $t$ in the range $\tau < t \le (\lambda/4)^{2j}$. Let $\alpha = 16\log(4/\lambda)>1$, and let $M = (t-\tau)/(\alpha j)$.

Since $\lambda$ is sufficiently small, $\alpha < 1/(2\lambda)$.
First, suppose that
$W_i \geq M$. In this case, by Lemma~\ref{lem:escape-high-weight-kills-balls}, for all $t'' \geq \tau$ and all $k\in [j]$, 
\[\Pr(\beta \in \b_k^{E^+}(t'')) \leq W_k/W_i \leq W_k/M = \alpha j W_k/(t-\tau) \leq j W_k/(\lambda(t-\tau)),\]
as required. For the rest of the proof, suppose $W_i< M$. Let
\[
i^+ = \max\{x \in \{i,\dots,j\} \colon W_i, W_{i+1}, \dots, W_x < M\}.
\]
Thus when $\beta$ sends from bin $i^+$, either $\beta$ sends from bin $j$ and escapes $E^+$ or $\beta$ enters a bin $\ell \in [j]$ of weight at least $M$. Let $T$ be the unique time such that $\beta \in \s_{i^+}^{E^+}(T)$ (observing that $T$ is finite almost surely), and let $\calE$ be the event that $T \le t$. For all $k \in [j]$, 
\begin{equation}\label{eq:tool-2-single-ball-1}
    \pr(\beta \in \b_k^{E^+}(t)) \le \pr(\beta \in \b_k^{E^+}(t) \mid \calE) + \pr(\overline{\calE}).
\end{equation}
If $i^+=j$ then $\pr(\beta \in \b_k^{E^+}(t) \mid \calE)=0$; otherwise, by Lemma~\ref{lem:escape-high-weight-kills-balls} applied with $t' = T$,
\begin{equation}\label{eq:tool-2-single-ball-2}
    \pr(\beta \in \b_k^{E^+}(t) \mid \calE) \le W_k/W_\ell \le W_k/M \le \alpha jW_k/(t-\tau)
\leq j W_k/(2 \lambda(t-\tau)).
\end{equation}
Hence by~\eqref{eq:tool-2-single-ball-1} and~\eqref{eq:tool-2-single-ball-2}, it suffices to prove that
\begin{equation}\label{eq:tool-2-single-ball-goal}
    \pr(\overline{\calE}) \le \alpha jW_k/(t-\tau).
\end{equation}

Observe that for a bin $x \in \{i,\dots,i^+\}$, the number of time steps between $\beta$ arriving at $x$ (i.e.\ $\min\{\hat{t} \colon \beta \in \b_x^{E^+}(\hat{t})\}$) and $\beta$ sending from $x$ (i.e.\ $\min\{\hat{t}\colon \beta \in \s_x^{E^+}(\hat{t})\}$) is a geometric variable with mean $W_x$. Thus $T-\tau$ is a sum of independent geometric variables with means $W_i,\dots,W_{i^+} < M$, so $\E(T-\tau) < (i^+-i+1)M$. We will apply Lemma~\ref{lem:janson-geo} to this sum with $p_* = 1/\max\{W_x\colon i \le x \le i^+\}$, $\mu = \sum_{x=i}^{i^+}W_i$, and $C=(t-\tau)/\mu$. Observe that $\mu \le jM$, so $C \ge (t-\tau)/(jM) = \alpha \ge 1$ as required by Lemma~\ref{lem:janson-geo}; hence
\[
    \pr(\overline{\calE}) = \pr(T > t) = \pr(T-\tau > C\mu) \le \exp(-p_*\mu(C-1-\log C)).
\]
As above, $C \ge \alpha \ge 2$, so $C-1-\log C \ge C/8$ and $\pr(\overline{\calE}) \le \exp(-p_*\mu C/8)$. Moreover, $p_*\mu C = p_*(t-\tau) \ge (t-\tau)/M = \alpha j$, so 
\[
    \pr(\overline{\calE}) \le e^{-\alpha j/8} = (\lambda/4)^{2j} \le 1/(t-\tau) < \alpha j W_k/(t-\tau)
\]
as required by~\eqref{eq:tool-2-single-ball-goal}; thus the result follows.
\end{proof}

We are now able to generalise Lemma~\ref{lem:h-below-stationary}, minimising the impact of a set $S$ of low-weight bins by waiting a sufficiently long time.

\begin{lemma} \label{lem:escape-low-weight-h-drops}
Let $\p$ be a send sequence with $p_0=1$, and fix a sufficiently small $\lambda \in (0,1/3)$. Let $j$ be a positive integer, let $\avect \in (\integers_{\ge 0})^j$, let $\nu\ge 0$ be a real number, and let $S \subseteq [j]$. Then for all positive integers $t \le (4/\lambda)^{2j}$ and all $k \in [j]$,
    \[
    p_kh_k^{\nu,\avect}(t) \le \frac{j}{\lambda t}\sum_{\ell \in S} a_\ell + \max\Big(\{ k \nu \} \cup \{p_y a_y + (k-y)\nu \colon y \in [k] \setminus S\}\Big).\]
\end{lemma}
\begin{proof}
Let $E^+$ be a constant-escape process on $[j]$ with send sequence~$\p$, start time $0$, end-of-step arrival rate~$\nu$, and satisfying $\romanarrvect^{E^+}_{[j]}(0) = \avect$; thus by Observation~\ref{obs:constant-escape}, 
\begin{equation}\label{eq:annoying-2b-step1}
p_kh_k^{\nu,\avect}(t) = p_k\E(b_k^{E^+}(t)).
\end{equation}
We will analyse $E^+$ by splitting its balls into two categories, which we track by coupling~$E^+$ with two constant-escape processes $E^+_1$ and $E^+_2$. The first category (tracked by $E_1^+$) consists of all balls initially present in some bin in $S$, while the second category (tracked by $E_2^+$) consists of all balls initially present in some bin in $[j] \setminus S$ together with all end-of-step arrivals after time~$0$.

To this end, define vectors $\overline{a}^S$ and $\overline{a}^{\bar{S}}$ by
\[
a^S_x = \begin{cases}
a_x & \mbox{ if }x \in S,\\
0 & \mbox{ otherwise},
\end{cases}
\qquad\qquad
a^{\bar{S}}_x = \begin{cases}
a_x & \mbox{ if }x \in [j] \setminus S,\\
0 & \mbox{ otherwise}.
\end{cases}
\]
Let $E^+_1$ be a constant-escape process on $j$ with send sequence~$\p$, start time~$0$, and end-of-step arrival rate~$0$, satisfying $\romanarrvect_{[j]}^{E_1^+}(0) = \avect^S$. Let $E_2^+$ be a constant-escape process on $j$ with send sequence~$\p$, start time~$0$, and end-of-step arrival rate~$\nu$, satisfying $\romanarrvect_{[j]}^{E_2^+}(0) = \avect^{\bar{S}}$. We couple $E^+$, $E^+_1$, and $E^+_2$ in the natural way, synchronising sends and escapes and copying all end-of-step arrivals of $E^+$ after time~$0$ in $E_2^+$. Thus, for all $\ell \in [j]$ and all $t \ge 0$, $\b_\ell^{E^+}(t) = \b_\ell^{E_1^+}(t) \cup \b_\ell^{E_2^+}(t)$ and $\b_\ell^{E_1^+}(t) \cap \b_\ell^{E_2^+}(t) = \emptyset$. In particular, by~\eqref{eq:annoying-2b-step1},
\begin{equation}\label{eq:annoying-2b-step2}
    p_kh_k^{\nu,\avect}(t) = p_k\E(b_k^{E_1^+}(t)) + p_k\E(b_k^{E_2^+}(t)).
\end{equation}

We bound each term on the right-hand side of~\eqref{eq:annoying-2b-step2} separately. Since $E_1^+$ has arrival rate zero,
\[
\bigcup_{x\in [j]} \b_x^{E_1^+}(t) \subseteq \bigcup_{x\in [j]} \b_x^{E_1^+}(0) = \bigcup_{x \in S} \arr_x^{E^+}(0);
\]
thus by linearity of expectation, 
\begin{equation*}
    \E(b_k^{E_1^+}(t)) = \sum_{x \in S}\sum_{\beta\in\arr_x^{E^+}(0)} \pr(\beta \in \b_k^{E_1^+}(t)).
\end{equation*}
By Lemma~\ref{lem:tool-2-single-ball} applied to $E_1^+$ (with the same $k$ and $t$), it follows that
\begin{equation}\label{eq:annoying-2b-step3a}
p_k\E(b_k^{E_1^+}(t)) \le p_k\sum_{x \in S} \sum_{\beta \in \arr_x^{E^+}(0)} \frac{jW_k}{\lambda t} = \frac{j}{\lambda t}\sum_{x \in S}a_x.    
\end{equation}

We next bound the other term on the right-hand side of~\eqref{eq:annoying-2b-step2}. By Observation~\ref{obs:constant-escape}, $p_k\E(b_k^{E_2^+}(t)) = p_kh_k^{\nu,\avect^{\bar{S}}}(t)$. Applying Lemma~\ref{lem:h-below-stationary}, taking $T=0$, $x=k$, and the $\avect$ of Lemma~\ref{lem:h-below-stationary} to be $\avect^{\bar{S}}$, we obtain
\begin{align*}
p_k\E(b_k^{E_2^+}(t)) 
&\le \max\Big(\{k\nu\} \cup \{p_yh_y^{\nu,\avect^{\bar{S}}}(0) + (k-y)\nu \colon y \in [k]\}\Big)\\
&= \max\Big(\{k\nu\} \cup \{p_ya_y + (k-y)\nu\colon y \in [k] \setminus S\} \cup \{(k-y)\nu\colon y \in [k] \cap S\}\Big).
\end{align*}
Since for all $y \in [k] \cap S$, $k-y \le k$, it follows that
\begin{equation}\label{eq:annoying-2b-step3b}
p_k\E(b_k^{E_2^+}(t)) \le \max\Big(\{k\nu\} \cup \{p_ya_y + (k-y)\nu\colon y \in [k] \setminus S\}\Big).
\end{equation}
The result now follows from combining equations~\eqref{eq:annoying-2b-step2}, \eqref{eq:annoying-2b-step3a} and~\eqref{eq:annoying-2b-step3b}.
\end{proof}
We conclude the section with three concentration bounds for constant-escape processes. The first two rely on the same application of standard Chernoff bounds, which we now state as an ancillary lemma.

\begin{lemma}\label{lem:h-collected-chernoff}
Let $\alpha > 0$, and let $j\ge1$ be a sufficiently large integer. Let $Y$ be a sum of independent Bernoulli random variables with mean $\mu$. If $\mu \leq \alpha^2 j(\log j)^2/6$, then with probability at least $1-j^{-3j}$, 
\begin{equation}\label{eq:h-collected-chernoff-small-mu}
    Y \le \alpha^2 j(\log j)^2.
    \end{equation}
    Otherwise, with probability at least $1-j^{-3j}$,
    \begin{equation}\label{eq:h-collected-chernoff-large-mu}
        Y \le \Big(1 + \frac{\alpha j^{1/2}\log j}{4\mu^{1/2}}\Big)\mu < 2\mu.
    \end{equation}
\end{lemma}
\begin{proof}
    Let $R = \alpha^2 j(\log j)^2$; we split into cases depending on $\mu$.
    
    \medskip\noindent\textbf{Case 1:} $\mu \le R/6$. Then by Lemma~\ref{lem:chernoff-large-dev}, since $j$ is large and $R \ge 6\mu$,
    \[
        \pr(Y > \alpha^2j(\log j)^2) = \pr(Y > R) \le 2^{-R} < j^{-3j},
    \]
    and so~\eqref{eq:h-collected-chernoff-small-mu} holds as required. 

    \medskip\noindent\textbf{Case 2: }$\mu > R/6$. Let $\delta = (R/(16\mu))^{1/2} < 1$. Observe that 
    \[
    \delta = \alpha j^{1/2}\log j/(4\mu^{1/2});
    \]
    hence by Lemma~\ref{lem:chernoff-small-dev-upper},
    \begin{align*}
        \pr\bigg(b_k^{E^+}(t) \ge \Big(1 + \frac{\alpha j^{1/2}\log j}{4\mu^{1/2}}\Big)\mu\bigg) = \pr\big(Y \ge (1+\delta)\mu\big) \le e^{-\delta^2 \mu/3} = e^{-R/48}.
    \end{align*}
    Since $j$ is large, $e^{-R/48} < j^{-3j}$ 
and the first inequality in ~\eqref{eq:h-collected-chernoff-large-mu} holds. The second inequality follows from $\delta < 1$.  
\end{proof}

The first concentration bound for constant-escape processes simply applies Lemma~\ref{lem:h-collected-chernoff} directly to the number of balls in a bin.

\begin{lemma}\label{lem:h-basic-chernoff}
Fix a send sequence $\p$ with $p_0=1$. Let $\alpha$ be a real number with $0 < \alpha \le 1$, let $j\ge 1$ be an integer, and suppose $j$ is sufficiently large relative to $\alpha$. Fix an integer $\tau \ge 0$, a real number $\nu\ge 0$, and $\avect \in (\integers_{\ge 0})^j$. Let $E^+$ be a constant-escape process on $[j]$ with send sequence~$\p$, start time $\tau$, end-of-step arrival rate $\nu$, and $\romanarrvect_{[j]}^{E^+}(\tau) = \avect$. Fix $k \in [j]$ and $t \ge \tau$, and let $\mu = h_k^{\nu,\avect}(t-\tau)$. If $\mu \leq \alpha^2j(\log j)^2/6$, then with probability at least $1-j^{-3j}$,
    \begin{equation}\label{eq:h-basic-chernoff-1}
        b_k^{E^+}(t) \le \alpha^2j(\log j)^2.
    \end{equation}
    Otherwise, with probability at least $1-j^{-3j}$,
    \begin{equation}\label{eq:h-basic-chernoff-2}
        b_k^{E^+}(t) \le \Big(1 + \frac{\alpha j^{1/2}\log j}{4\mu^{1/2}}\Big)\mu < 2\mu.
    \end{equation}
\end{lemma}
\begin{proof}
Let $Y = b_k^{E^+}(t)$, and observe that $Y$ is a sum of independent Bernoulli variables (namely indicator variables for individual balls) with $\E(Y) = h_k^{\nu,\avect}(t-\tau)$ by Observation~\ref{obs:constant-escape}. The result is now immediate from Lemma~\ref{lem:h-collected-chernoff}.
\end{proof}

When $W_k$ is large, the   additive error of Lemma~\ref{lem:h-basic-chernoff} will often be unacceptably high. We therefore provide two alternative Chernoff bounds which are stronger when $E^+$ has only been running for a short time relative to $W_k$. Lemma~\ref{lem:h-conc-sends} lower-bounds total sends from bin $k$ since the start time, and Lemma~\ref{lem:h-conc-entries} upper-bounds total entries into bin $k$ since the start time (via sends and end-of-step arrivals). Both proofs rely on the following ancillary lemma.

\begin{lemma}\label{lem:h-bernoulli-decomposition}
Fix a send sequence $\p$ with $p_0=1$. Let $j\ge 1$ be an integer, and suppose $j$ is sufficiently large. Fix an integer $\tau \ge 0$, a real number $\nu\ge 0$, and $\avect \in (\integers_{\ge 0})^j$. Let $E^+$ be a constant-escape process on $[j]$ with send sequence~$\p$, start time $\tau$, end-of-step arrival rate $\nu$, and $\romanarrvect_{[j]}^{E^+}(\tau) = \avect$. Then for all integers $t_1$ and $t_2$ with 
$\tau < t_1 \leq t_2 $ and all $\ell \in [j]$, $\sum_{t=t_1}^{t_2}s_\ell^{E^+}(t)$ can be expressed as a sum of independent Bernoulli variables. Moreover, these variables and the family of variables $\{\romanarr_{\ell+1}^{E^+}(t)\colon t_1 \le t \le t_2\}$ are mutually independent.
\end{lemma}
\begin{proof}
Let $\mathbf{S} = \bigcup_{t=t_1}^{t_2}\s_\ell^{E^+}(t)$; we must decompose $S:=|\mathbf{S}|$ into a sum of independent Bernoulli variables. For all integers $x \in [\ell]$ and $t'$ with $\tau \le t' \le t_2$, let $N_{x,t'}$ be the number of balls in $\arr^{E^+}_x(t')$ that lie in $\mathbf{S}$. Since every ball in $\mathbf{S}$ entered $E^+$ at a unique bin and time, it lies in a unique such set $\arr^{E^+}_x(t')$; thus
\[
    S = \sum_{x=1}^\ell\sum_{t'=\tau}^{t_2} N_{x,t'}.
\]
For all $x \in [\ell]$ and $y \in [\romanarr_x^{E^+}(\tau)]$, let $I_{x,y}$ be the indicator variable for the event that the (lexicographically) $y$'th ball in $\arr_x^{E^+}(\tau)$ lies in $S$. Then
\begin{equation}\label{eq:bernoulli-decomposition-1}
    S = \sum_{x=1}^\ell \sum_{y=1}^{\romanarr_x^{E^+}(\tau)} I_{x,y} + \sum_{x=1}^{\ell} \sum_{t'=\tau+1}^{t_2} N_{x,t'}.
\end{equation}
Recall from Definition~\ref{def:constant-escape} that for all $x \in [\ell]$ and $t > \tau$, the variables $\romanarr_x^{E^+}(t)$ are independent Bernoullis, and that balls move through $E^+$ independently and hence also lie in $\mathbf{S}$ independently. Moreover, all variables $I_{x,y}$ and $N_{x,t'}$ are independent from all variables $\arr_{\ell+1}^{E^+}(t')$. Thus the right-hand side of~\eqref{eq:bernoulli-decomposition-1} is the desired decomposition of $S$.
\end{proof}

\begin{lemma}\label{lem:h-conc-entries}
Fix a send sequence $\p$ with $p_0=1$. Let $j\ge 1$ be an integer, and suppose $j$ is sufficiently large. Fix an integer $\tau \ge 0$, a real number $\nu\ge 0$, and $\avect \in (\integers_{\ge 0})^j$. Let $E^+$ be a constant-escape process on $[j]$ with send sequence~$\p$, start time $\tau$, end-of-step arrival rate $\nu$, and $\romanarrvect_{[j]}^{E^+}(\tau) = \avect$. Fix integers $t > \tau$ and $k \in [j]$, and let
\[
\mu = (t-\tau)\nu + \sum_{t'=0}^{t-\tau-1} p_{k-1}h_{k-1}^{\nu,\avect}(t').
\]
Then with probability at least $1-j^{-3j}$,
\[
\Big|\bigcup_{t'=\tau+1}^t\big(\s_{k-1}^{E^+}(t') \cup \arr_k^{E^+}(t')\big)\Big| \le \max\Big\{j(\log j)^2, \Big(1+\frac{j}{\mu^{1/2}}\Big) \mu \Big\}.
\]
\end{lemma}
\begin{proof}
    Fix $k \in [j]$, and let
    \[
        \mathbf{S} = \bigcup_{t'=\tau+1}^t \s_{k-1}^{E^+}(t'),\qquad Y = \Big|\mathbf{S} \cup \bigcup_{t'=\tau+1}^t \arr_k^{E^+}(t')\Big| = |\mathbf{S}| + \sum_{t'=\tau+1}^t \romanarr_k^{E^+}(t').
    \] 
Our goal will be to apply Lemma~\ref{lem:h-collected-chernoff} to $Y$. 
From the definition of a constant-escape process (Definition~\ref{def:constant-escape}), it is immediate that all variables $\romanarr_k^{E^+}(t')$ are independent and Bernoulli. By Lemma~\ref{lem:h-bernoulli-decomposition} applied with $\ell=k-1$, $|\mathbf{S}|$ can also be expressed as a sum of independent Bernoulli variables which are independent of the $\romanarr_k^{E^+}(t')$ terms. Thus $Y$ is a sum of independent Bernoulli variables as required by Lemma~\ref{lem:h-collected-chernoff}. 
By Observation~\ref{obs:constant-escape}, 
for all $t'\geq \tau$,
$\E(b_{k-1}^{E^+}(t')) = h_{k-1}^{\nu,\avect}(t'-\tau)$.
Thus, for all $t' > \tau$,
$\E(s_{k-1}^{E^+}(t')) = p_{k-1} h^{\nu,\avect}(t'-\tau-1)$. 
By Definition~\ref{def:constant-escape}, 
for all $t' > \tau$,
$\E(\romanarr_k^{E^+}(t')) = \nu$. Thus $\E(Y) = \mu$. By applying Lemma~\ref{lem:h-collected-chernoff} to $Y$ with $\alpha=1$ and taking the maximum of the two possible upper bounds~\eqref{eq:h-collected-chernoff-small-mu} and~\eqref{eq:h-collected-chernoff-large-mu}, we obtain that with probability at least $1-j^{-3j}$,
\[
Y \le \max\Big\{j(\log j)^2,\Big(1 + \frac{j^{1/2}\log j}{4\mu^{1/2}}\Big)\mu \Big\} \le \max\Big\{j(\log j)^2,\Big(1 + \frac{j}{\mu^{1/2}}\Big)\mu \Big\},
\]
    as required.
\end{proof}

\begin{lemma}\label{lem:h-conc-sends} 
Fix a send sequence $\p$ with $p_0=1$. Let $j\ge 1$ be an integer, and suppose $j$ is sufficiently large. Fix an integer $\tau \ge 0$, a real number $\nu\ge 0$, and $\avect \in (\integers_{\ge 0})^j$. Let $E^+$ be a constant-escape process on $[j]$ with send sequence~$\p$, start time $\tau$, end-of-step arrival rate $\nu$, and $\romanarrvect_{[j]}^{E^+}(\tau) = \avect$. Fix integers $t > \tau$ and $k \in [j]$, and let
\[
\mu = p_k\sum_{t'=0}^{t-\tau-1}h_k^{\nu,\avect}(t').
\]
Suppose $\mu \ge j(\log j)^2$. Then with probability at least $1-j^{-3j}$,
\[
\Big|\bigcup_{t'=\tau+1}^t \s_k^{E^+}(t')\Big| \ge \Big(1 - \frac{j}{\mu^{1/2}}\Big)\mu.
\]
\end{lemma}
\begin{proof}
Fix $k \in [j]$, and let $Y = \sum_{t'=\tau+1}^{t} s_k^{E^+}(t')$. Observe that $Y$ is a sum of independent Bernoulli variables by Lemma~\ref{lem:h-bernoulli-decomposition}.
By Observation~\ref{obs:constant-escape}, 
for all $t'\geq \tau$,
$\E(b_{k}^{E^+}(t')) = h_{k}^{\nu,\avect}(t'-\tau)$.
Thus, for all $t' > \tau$,
$\E(s_{k}^{E^+}(t')) = p_{k} h^{\nu,\avect}(t'-\tau-1)$.
Thus, $\E(Y) = \mu$. Let $\delta = j/\mu^{1/2}$. If $\delta \ge 1$ then $Y \ge 0 \ge (1-\delta)\mu$ with certainty; otherwise, by applying Lemma~\ref{lem:chernoff-small-dev} to $Y$ with $\delta = j/\mu^{1/2}$, we obtain that
\[
        \pr(Y \leq (1-\delta)\mu) \le e^{-j^2/2} < j^{-3j}
    \]
    and so the result follows immediately.
\end{proof}

\subsection{Pseudorandomness properties of escape processes}\label{sec:escape-pseudorandomness}

We first define some constants.

\begin{definition}\label{def:escape-constants}  
Let $\newsym{xi}{$\tfrac{1}{160} 
(\lambda/4)^{10^7/\lambda^3}$ - Def~\ref{def:escape-constants}}{\xi} = \tfrac{1}{160} 
(\lambda/4)^{10^7/\lambda^3}$. We define $\newsym{P}{$P_k = \xi\prod_{\ell \le k} (1+1/\ell^2)$ - Def~\ref{def:escape-constants}}{\Pvect} = (P_0, P_1, P_2,\dots)$ where $P_k = \xi\prod_{\ell \in [k]} (1+1/\ell^2)$ (so that $P_0 = \xi$). For all integers $j \ge 0$, we define $\newsym{nu-high}{$1/(\log j)^{100}$ - Def~\ref{def:escape-constants}}{\nu_j^{\high}} = 1/(\log j)^{100}$ and $\newsym{nu-low}{$1/j^{100}$ - Def~\ref{def:escape-constants}}{\nu_j^{\low}} = 1/j^{100}$.
\end{definition}

\begin{remark}\label{rem:P-bounded}
    For all integers $k \ge 0$, $P_k \le \xi \prod_{\ell=1}^\infty (1 + 1/\ell^2) = \xi \sinh(\pi)/\pi < 4\xi$.
\end{remark}

Informally, we will prove the following long-term behaviour of escape processes in the VEB coupling. Consider the high-level state $\Psi = (g,\tau,j,\avect,\calS,\type)$ (where $\type \neq \Failure$) 
and the escape process~$E_g$ in the VEB coupling.  The ``relevant'' time-steps~$t$ for the study of~$E_g$ are  
the time-steps (starting with $t=\tau$) when $\Psi(t) = \Psi$.
The high-level state transitions do not depend on the evolution of~$E_g$, but we have do have bounds on the timing of high-level state transitions.
We will show that with summable failure probability:
\begin{itemize}
    \item If $\Psi$ is a $j$-filling or $j$-advancing state, $\nu^{E_g}(t) \le \nu_j^\low < \nu_j^\high$ for all relevant $t$. For all $k \in [j]$, $b_k^{E_g}(t)$ will be approximately bounded above by $P_kW_k$ with error that drops over time.
    \item If $\Psi$ is a $j$-refilling or $j$-stabilising state, $\nu^{E_g}(t) \le \nu_j^\high$ for all relevant $t$.
    \item If $\Psi$ is a $j$-refilling state, for all $k \in [j]$ with $W_k$ suitably large, $b_k^{E_g}(t)$ will be approximately bounded above by $P_kW_k$ with error that grows over time but remains manageable for all relevant $t$.  Other bins $k \in [j]$ may reach $b_k^{E_g}(t) \approx k\nu_k^{\low}W_k < j^2W_k$.
    \item If $\Psi$ is a $j$-stabilising state, we 
will show that $E_g$ behaves the  same way as if $\Psi$ were a $j$-refilling state for a short setup period. After this setup period, we will show that $E_g$  gradually converges to the initial conditions 
that would be satisfied
if $\Psi$ were a $j$-filling state.
\end{itemize}

In short: if $\Psi$ is $j$-filling or $j$-advancing, then $E_g$ empties to a low level and remains there. If $\Psi$ is $j$-refilling, then $E_g$ gradually fills. If $\Psi$ is $j$-stabilising, then $E_g$ fills for a short period and then empties. Motivated by this informal description, we now define the following invariants.

\begin{definition}\label{def:excess}
Fix $\lambda \in (0,1)$, a send sequence~$\p$ with $p_0=1$, an integer $j \ge 1$, and $\xvect \in (\reals_{\ge 0})^j$. For all $\delta > 0$, we define 
\[
\excess(\xvect, \delta) = \sum_{k=1}^j\max\{0, x_k - \delta P_kW_k\}.
\]
\end{definition}

We can view excess as measuring the total extent to which a standard uniform bound of the form $x_k \le \delta P_k W_k$ fails over $k \in [j]$. We will not in general be able to maintain such bounds uniformly over all bins (due to the presence of low-weight bins), but we will be able to maintain low excess. 

\begin{definition}\label{def:conditions}
Fix $\lambda \in (0,1)$ and a send sequence~$\p$ with $p_0=1$, and an integer $j \ge 1$. We define the following conditions. 
\begin{enumerate}[({Con}1)]
\item 
Let $d \ge 0$ be real. A tuple $\mvect \in (\reals_{\ge 0})^j$ \emph{satisfies (Con\ref{item:Con1}) for $d$} if $\excess(\mvect, 1) \le d$.\label{item:Con1}
\item A tuple $\mvect \in (\reals_{\ge 0})^j$ \emph{satisfies 
(Con\ref{item:Con2})} if
$\excess(\mvect, \xi^{-1/3}) \le 
\lambda j/500
$.\label{item:Con2} 
\item Let $\varpi,\theta,\nu \ge 0$ be real. A tuple $\mvect \in (\reals_{\ge 0})^j$ \emph{satisfies
(Con\ref{item:Con3}) for $(\varpi,\theta,\nu)$} if,
for all $k \in [j]$ with $W_k \ge \varpi$, $p_km_k \le \nu k + (1+\theta)P_k$.
\label{item:Con3}
\item Let $G$ be a positive real.
A tuple $\mvect \in (\reals_{\ge 0})^j$ \emph{satisfies (Con\ref{item:Con4}) for~$G$} if for all $k \in [j]$, $p_km_k \le G j^2$.\label{item:Con4}
\end{enumerate}
\end{definition}

Note that we do not give $j$ as an explicit parameter in Definition~\ref{def:conditions}, as its value is specified by the length of $\mvect$.

\begin{definition}\label{def:good-start}
Fix $\lambda \in (0,1)$, a send sequence~$\p$ with $p_0=1$, 
a positive integer~$j$, and 
a tuple $\mvect \in (\reals_{\ge 0})^j$.
For all $d \geq 0$, we say $\mvect$ is a $d$-AF-good beginning 
if it has the following properties:
\begin{enumerate}[({A}FB1)]
\item $\mvect$ satisfies (Con\ref{item:Con1}) for~$d$.
 \label{item:AFB1}
\item $\mvect$ satisfies (Con\ref{item:Con2}).\label{item:AFB2}
\end{enumerate}
For all $d,\varpi,\theta,G \ge 0$, we say $\mvect$ is a $(d, \varpi, \theta, G)$-RS-good beginning if it has the following properties:
\begin{enumerate}[(RSB1)]
\item $\mvect$ satisfies (Con\ref{item:Con1}) for~$d$.\label{item:RSB1}
\item $\mvect$ satisfies (Con\ref{item:Con3}) for $(\varpi,\theta,\nu_j^\high)$.\label{item:RSB2}
\item $\mvect$ satisfies (Con\ref{item:Con4}) for $G$.\label{item:RSB3}
\end{enumerate}
\end{definition}

In general, we will require that an escape process for an Advancing or Filling high-level state has an AF-good beginning (see Lemmas~\ref{lem:advancing-escape-final} and~\ref{lem:filling-escape-final}) and that an escape process for a Refilling or Stabilising high-level state has an RS-good beginning (see Lemmas~\ref{lem:veb-refilling-transition} and~\ref{lem:veb-escape-stab}). In the remainder of this section, we will prove some generally-useful ancillary lemmas. In Section~\ref{sec:constant-escape-invariants}, we will analyse the behaviour of constant-escape processes with arrival rates $\nu_j^\low$ and $\nu_j^\high$ in terms of (Con\ref{item:Con1})--(Con\ref{item:Con4}). For example, we prove that a constant-escape process with an RS-good arrival vector remains RS-good for a short period (with slightly increased error). In Section~\ref{sec:escape-invariants}, we extend these results to escape processes by coupling them with constant-escape processes. In Section~\ref{sec:veb-transition-invariants}, we extend these results into the context of the VEB coupling, and in Section~\ref{sec:escape-analysis-final} we analyse the entire coupling to prove our goal of Lemma~\ref{lem:VEB-escape}.

\begin{lemma}\label{lem:af-good-implies-rs-good}
Fix $\lambda \in (0,1)$, a send sequence~$\p$ with $p_0=1$, 
and a positive integer~$j$.
Let $d \ge 0$, and suppose $\avect$ is a $d$-AF-good beginning. Then $\avect$ also satisfies (Con\ref{item:Con3}) for $(dj^3,1/(\xi j^3),0)$, satisfies (Con\ref{item:Con4}) for~$1$, and is a $(d,dj^3,1/(\xi j^3),1)$-RS-good beginning.
\end{lemma}
\begin{proof}
First, consider $k \in [j]$ with $W_k \ge dj^3$. Then since $\avect$ satisfies (Con\ref{item:Con1}) for $d$, 
$\excess(\avect, 1) \le d$.
Recall that 
$
\excess(\avect, 1) = \sum_{\ell=1}^j\max\{0, a_\ell - P_\ell W_\ell\}$.
So $a_k \leq   P_k W_k + d$.
Using the fact that $W_k \geq d j^3$ and   $P_k \ge \xi$,
    \[
        a_k \le  P_kW_k+d \le P_kW_k(1 + 1/(j^3P_k)) \le (1+1/(\xi j^3))P_kW_k.
    \]
    Thus $\avect$ satisfies (Con\ref{item:Con3}) for $(dj^3,1/(\xi j^3),0)$ as claimed. 
    
We next prove that $\avect$ satisfies (Con\ref{item:Con4}) for $1$. 
We will use the fact that, 
since $\avect$ satisfies (Con\ref{item:Con2}),
\[
\excess(\avect, \xi^{-1/3}) = \sum_{k=1}^j\max\{0, a_k - \xi^{-1/3} P_kW_k\} \leq \lambda j/500,
\]
so for all $k\in [j]$, 
$a_k \leq \lambda j/500 + 
\xi^{-1/3} P_kW_k$.
We wish to show that $a_k \leq j^2 W_k$. First,
$\lambda j/500 \leq j^2/2 \leq j^2 W_k/2$.
Second, 
by Remark~\ref{rem:P-bounded},
$\xi^{-1/3} P_kW_k
\leq 4 \xi^{2/3} W_k \leq j^2 W_k/2$. Thus, 
$\avect$ satisfies (Con\ref{item:Con4}) for $1$. 
 It is then immediate that $\avect$ is a $(d,dj^3,1/(\xi j^3),1)$-RS-good beginning.
\end{proof}

\begin{lemma}\label{lem:applied-excess-conc}
Fix a send sequence~$\p$ with $p_0=1$. Let $j$ be a sufficiently large positive integer.  Fix an integer $\tau \ge 0$ and a real number $\nu\ge 0$. Let $A_1,\dots,A_j$ be independent binomial random variables, 
and for each $i\in [j]$
let $m_i = \E(A_i)$. 
Let $\Avect=(A_1,\ldots,A_j)$.
Let~$E^+$ be a constant-escape process on $[j]$ with send sequence~$\p$, start time $\tau$, end-of-step arrival rate $\nu$, and $\romanarrvect_{[j]}^{E^+}(\tau) = \Avect$. Let $\zeta$ and $\eta$ be real numbers with $\zeta \ge 1$ and $\eta \ge \xi^{-1/3}\zeta$. Fix $t \ge \tau$ such that, for all $k \in [j]$, $h_k^{\nu,\mvect}(t-\tau) \le \zeta P_kW_k$. Then with probability at least $1-e^{-500j/\lambda^2}$,
$\excess(\bvect_{[j]}^{E^+}(t),\eta) \le \lambda j/500$.
\end{lemma}
\begin{proof}
By Observation~\ref{obs:constant-escape}, for all $i\in [j]$,
$\E(b_i^{E^+}(t)) = h^{\nu,\mvect}(t-\tau)$.
By Lemma~\ref{lem:constant-escape-gen-multi},
$b_1^{E^+}(t),\dots,b_j^{E^+}(t)$ can be expressed as coordinates of a generalised multinomial random variable.
We will apply Lemma~\ref{lem:excess-chernoff}, taking the $Z_1,\dots,Z_n$ of that lemma to be $b_1^{E^+}(t),\dots,b_j^{E^+}(t)$, the $\muvect$ of that lemma to 
be $\muvect = (\mu_1,\ldots,\mu_j)$ where
$\mu_k := \zeta P_kW_k$, the $\delta$ of that lemma to be $\eta/\zeta$, and the $z$ of that lemma to be $\lambda j/500$. 
By definition,
\[
\excess(\bvect_{[j]}^{E^+}(t), \eta) = \sum_{k=1}^j \max\{0, b_k^{E^+}(t) - \eta P_kW_k\},
\]
and this is equal to the~$Z^+$ of
Lemma~\ref{lem:excess-chernoff}.
Thus by Lemma~\ref{lem:excess-chernoff},
\[ 
\pr\big(\excess(\bvect_{[j]}^{E^+}(t), \eta) > \lambda j/500\big) \le 2^j (\eta/\zeta)^{-\lambda j/500}.
\]
To finish, we need only 
show $2^j (\zeta/\eta)^{\lambda j/500} \leq \exp(-500j/\lambda^2)$.This follows since $\zeta/\eta \leq \xi^{1/3}$ and from the definition of~$\xi$ (Definition~\ref{def:escape-constants}).
\end{proof}
\begin{corollary}\label{cor:applied-excess-conc}
Fix a send sequence~$\p$ with $p_0=1$. Let $j$ be a sufficiently large positive integer.  Fix an integer $\tau \ge 0$ and a real number $\nu\ge 0$. Let $A_1,\dots,A_j$ be independent binomial random variables, 
and for each $i\in [j]$
let $m_i = \E(A_i)$. 
Let $\Avect=(A_1,\ldots,A_j)$.
Let~$E^+$ be a constant-escape process on $[j]$ with send sequence~$\p$, start time $\tau$, end-of-step arrival rate $\nu$, and $\romanarrvect_{[j]}^{E^+}(\tau) = \Avect$.
Let $t \ge \tau$, and suppose that for all $k \in [j]$, $h_k^{\nu,\mvect}(t-\tau) \le P_kW_k$. Then with probability at least $1-e^{-500j/\lambda^2}$, $\bvect_{[j]}^{E^+}(t)$ satisfies (Con\ref{item:Con2}).
\end{corollary}
\begin{proof}
Immediate from Lemma~\ref{lem:applied-excess-conc}, taking   $\zeta=1$ and $\eta=\xi^{-1/3}$.
\end{proof}

\begin{lemma}\label{lem:escape-means-stabilise}
Let $\p$ be a send sequence with $p_0=1$. Fix a sufficiently small $\lambda \in (0,1/3)$. Let $j$ be a positive integer, let $\avect \in (\integers_{\ge 0})^j$, and let $\nu,\theta\ge 0$ and $\varpi \ge 1$ be real numbers. Suppose that $\avect$ satisfies (Con\ref{item:Con3}) for $(\varpi,\theta,\nu)$ and (Con\ref{item:Con4}) for~$2$. For all $x\in [j]$ let $\ell(x) = \max(\{0\} \cup \{y \in [x]\colon W_y \ge \varpi\})$. Then for all $x \in [j]$ and all positive integers $t \le (4/\lambda)^{2j}$,
\[
p_xh_x^{\nu,\avect}(t) \le x\nu + \Big(1 + \theta + \frac{2\varpi j^4}{\xi \lambda t}\Big)P_{\ell(x)}.
\]
In particular, $\hvect_{[j]}^{\nu,\avect}(t)$ satisfies (Con\ref{item:Con3}) for $(\varpi, \theta + 2\varpi j^4/(\xi \lambda t), \nu)$.
\end{lemma}
\begin{proof}
For all $y \in [j]$ and all $t$, we suppress $\nu$ and $\avect$ from $h_y^{\nu,\avect}(t)$. Let $S = \{y \in [j]\colon W_y < \varpi \}$. By Lemma~\ref{lem:escape-low-weight-h-drops} applied to this choice of $S$, for all $x\in [j]$,
\begin{equation}\label{eq:escape-means-stabilise-1}
p_xh_x(t) \le \frac{j}{\lambda t}\sum_{\ell \in S}a_\ell + \max\Big(\{x\nu\} \cup \{p_ya_y + (x-y)\nu\colon y \in [x], W_y \ge \varpi \}\Big).
\end{equation}
Observe that $|S| \le j$. Moreover, $\avect$ satisfies (Con\ref{item:Con4}) for~$2$ by hypothesis, so for $\ell \in S$, $a_\ell \le 2j^2W_\ell < 2\varpi j^2$. Hence
    \[
        \frac{j}{\lambda t}\sum_{\ell \in S}a_\ell < \frac{2 \varpi j^4}{\lambda t}.
    \]
Moreover, since $\avect$ satisfies (Con\ref{item:Con3}) for $(\varpi,\theta,\nu)$
by hypothesis, for all $y \in [x]$ with $W_y \ge \varpi $, \[p_ya_y + (x-y)\nu \le x\nu + (1+\theta)P_y.\] 
Hence
    \begin{align*}
        &\max\Big(\{x\nu\} \cup \{p_ya_y + (x-y)\nu\colon y \in [x], W_y \ge \varpi \}\Big) \\
        \qquad\qquad\qquad\qquad&\le \max\Big(\{x\nu\} \cup \Big\{x\nu + (1+\theta)P_y \colon y \in [x], W_y \ge \varpi \Big\}\Big) \le x\nu + (1+\theta)P_{\ell(x)}.
    \end{align*}
    Together with~\eqref{eq:escape-means-stabilise-1}, we have proved
    \begin{align*}
        p_xh_x(t) \le 2 \varpi j^4/(\lambda t) + x\nu + (1+\theta)P_{\ell(x)}.
    \end{align*}
    Since for all $x \in [j]$ we have $P_x \ge \xi$, it follows that $p_xh_x(t) \le x\nu + (1+\theta+2\varpi j^4/(\xi \lambda t))P_{\ell(x)}$ as required.
\end{proof}

\subsection{Analysing constant-escape processes with pseudorandom starts}\label{sec:constant-escape-invariants}

Let $\Psi$ be a high-level state, and let $j = j^\Psi$. In analysing an escape process with high-level state $E$, we will be concerned with constant-escape processes with two arrival rates: $\nu_j^\high$ and $\nu_j^\low$.
\begin{itemize}
    \item If $\Psi$ is an advancing or filling state, we will couple $E$ to a single constant-escape process of rate $\nu_j^\low$.
    \item If $\Psi$ is a refilling state, we will couple $E$ to a single constant-escape process of rate $\nu_j^\high$.
    \item If $\Psi$ is a stabilising state, our analysis will be more complicated. We will first couple $E$ to a constant-escape process of rate $\nu_j^\high$ until its behaviour has stabilised somewhat from its initial conditions, then couple $E$ to a constant-escape process of rate $\nu_j^\low$. We will show that with failure probability $\exp(-(\log j)^2)$, this second constant-escape process will reach an AF-good beginning. If it doesn't, we will repeat the process until success, thereby reducing the overall failure probability to $\exp(-300j/\lambda^2)$. After reaching an AF-good beginning, we couple to a final constant-escape process of rate $\nu_j^\low$ to maintain the AF-good beginning until the next high-level state transition.
\end{itemize}
In this section, we will state and prove the results we need on these constant-escape processes. Most of these lemmas will typically be used in multiple contexts, so we state them in a general fashion.

We now show that a constant-escape process $E^+$ with arrival rate $\nu_j^\low$ and reasonable starting conditions stabilises rapidly, with long-term ball populations bounded by $b_k^{E^+}(t) \le P_kW_K$ (Lemma~\ref{lem:af-heavy-bins-empty}) and with $\bvect_{[j]}^{E^+}(t)$ acting as an AF-good beginning (Lemma~\ref{lem:af-constant-stabilise}).

\begin{lemma}\label{lem:af-heavy-bins-empty}
Fix a send sequence $\p$ with $p_0=1$. 
Fix a sufficiently small $\lambda \in (0,1/3)$.
Let $j$ be a sufficiently large integer. Fix an integer $\tau \ge 0$ and a tuple $\avect \in (\integers_{\ge 0})^j$. Let $E^+$ be a constant-escape process on $[j]$ with send sequence~$\p$, start time $\tau$, end-of-step arrival rate $\nu_j^\low$, and $\romanarrvect_{[j]}^{E^+}(\tau) = \avect$. Let $d$, $\varpi$, and $\theta$ be real numbers with $d\ge 0$, $\varpi\ge 1$, and $0 \le \theta < 1/(64j^2)$. Suppose that $\avect$ satisfies (Con\ref{item:Con1}) for $d$, (Con\ref{item:Con3}) for $(\varpi,\theta,\nu_j^\low)$, and (Con\ref{item:Con4}) for~$2$. 
Fix $t$ in the range 
$\max\{\varpi j^{12}, d j^3\} \le t-\tau \le (4/\lambda)^{2j}$
Then, 
with probability at least $1- 2j^{-3j+1}$,
for all $k \in [j]$ with $W_k \ge j^{10}$,
$b_k^{E^+}(t) < P_kW_k$ and $h_k^{\nu_j^\low,\avect}(t-\tau) < P_kW_k$.
\end{lemma}
\begin{proof}
For all $x \in [j]$ and all $t'$, we suppress $\nu_j^\low$ and $\avect$ from $h_x^{\nu_j^\low,\avect}(t')$, and we write $\nu = \nu_j^\low$.  
By a union bound, it suffices to prove that, for any fixed bin $k \in [j]$ with $W_k \ge j^{10}$,
with probability at least $1-2j^{-3j}$,
$ b_k^{E^+}(t) > P_kW_k  $ and
$h_k^{\nu_j^\low,\avect}(t-\tau) < P_kW_k$.
As such, fix $k \in [j]$ with $W_k \ge j^{10}$. We split into two cases depending on $h_k(t-\tau)$.
    
    \medskip\noindent \textbf{Case 1:} $h_k(t-\tau) \le (1+1/(4k^2))P_{k-1}W_k$. By Lemma~\ref{lem:h-basic-chernoff} applied with $\alpha=1$, with probability at least $1-j^{-3j}$,
    \begin{align*}
        b_k^{E^+}(t) &\le \max\Big\{j(\log j)^2, \Big(1 + \frac{j^{1/2}\log j}{4h_k(t-\tau)^{1/2}}\Big)h_k(t-\tau)\Big\} < j(\log j)^2 + \Big(1 + \frac{j}{h_k(t-\tau)^{1/2}}\Big)h_k(t-\tau).
    \end{align*}
The right-hand-side is increasing in $h_k(t-\tau)$, so
since we are in Case~1,
\[ b_k^{E^+}(t) < j(\log j)^2 + \Big(1 + \frac{j}{((1+1/(4k^2))P_{k-1}W_k)^{1/2}}\Big)\Big(1+\frac{1}{4k^2}\Big)P_{k-1}W_k.
\]

Observe that  
$(1+1/(4k^2))P_{k-1}W_k \geq \xi j^{10}$,  so
\[ b_k^{E^+}(t) < j(\log j)^2 + 
\Big(1 + 
\frac{j}{ \xi^{1/2} j^{5}}
\Big)\Big(1+\frac{1}{4k^2}\Big)P_{k-1}W_k  < \Big(1 + \frac{1}{k^2}\Big)P_{k-1}W_k = P_kW_k.
\]  From the definition of~$P_k$ (Definition~\ref{def:escape-constants}) 
we also have $h_k(t-\tau) < P_kW_k$, so the result follows for Case~1.

\medskip\noindent \textbf{Case 2:} $h_k(t-\tau) > (1+1/(4k^2))P_{k-1}W_k$.  Every ball in $\b_k^{E^+}(t)$ satisfies one of the following:
it is in $\b_k^{E^+}(\tau)$, it is sent into bin $k$ from bin $k-1$
during steps $\tau+1,\ldots,t$, or it arrived in bin $k$ 
as an end-of-step arrival during steps $\tau+1,\ldots,t$. 
Moreover, any ball which is sent from bin $k$ 
during steps $\tau+1,\ldots,t$
is no longer in bin $k$ after it is sent. 
Thus 
\begin{equation}\label{eq:af-heavy-bins-empty-0}
        b_k^{E^+}(t) = a_k + \Big|\bigcup_{t'=\tau+1}^t \big(\s_{k-1}^{E^+}(t') \cup \arr_k^{E^+}(t')\big)\Big| - \Big|\bigcup_{t'=\tau+1}^t \s_k^{E^+}(t') \Big|.
    \end{equation}
We will bound the terms of the right-hand side of~\eqref{eq:af-heavy-bins-empty-0} (other than $a_k$) separately using Lemmas~\ref{lem:h-conc-entries} and~\ref{lem:h-conc-sends}. In order to usefully apply these lemmas, we bound $h_k(t'-\tau)$ for large $t'$. Let $T_0 = \tau + \varpi j^7$.  For all $x\in [j]$, let
$\ell(x) = \max(\{0\} \cup \{y \in [x]\colon W_y \ge \varpi\})$.
We can apply Lemma~\ref{lem:escape-means-stabilise}, taking the $t$ of that lemma to be any positive
integer~$t'  \leq (4/\lambda)^{2j}$
to obtain that  
\[
p_xh_x(t') \le x\nu + \Big(1 + \theta + \frac{2\varpi j^4}{\xi \lambda t'}\Big)P_{\ell(x)}.\]
Since  $ \theta + 2\varpi j^4/(\xi \lambda t') \leq 1/(32j^2)$ for every
$t' \geq T_0-\tau$, we get
\begin{equation}\label{eq:upperbound}
\forall x\in [j] ,\forall t' \in \{T_0-\tau,\ldots, (4/\lambda)^{2j}\},
p_{x} h_{x} (t') \le x\nu + 
\big(1 +  1/(32j^2) \big)P_{\ell(x)}.
\end{equation}

    \medskip\noindent \textbf{Bounding the middle term above:} For brevity, let
    \begin{align*}
        Z_1 = \Big|\bigcup_{t'=\tau+1}^t \big(\s_{k-1}^{E^+}(t') \cup \arr_k^{E^+}(t')\big)\Big|,\qquad \mu_1 = (t-\tau)\nu + \sum_{t'=0}^{t-\tau-1}p_{k-1}h_{k-1}(t').
    \end{align*}
    Then by Lemma~\ref{lem:h-conc-entries}, with probability at least $1-j^{-3j}$,
    \begin{equation}\label{eq:af-heavy-bins-middle}
        Z_1 \le \max\Big\{j(\log j)^2, \Big(1+\frac{j}{\mu_1^{1/2}}\Big)\mu_1\Big).
    \end{equation}
    
We next bound $\mu_1$ above. First observe that for all $t'\geq 0$, by Lemma~\ref{lem:h-below-stationary} (taking the $T$ of that lemma to be $0$ and the $t$ of that lemma to be  $t'$), 
    \[
        p_{k-1}h_{k-1}(t') \le \max\big(\{(k-1)\nu\} \cup \{p_ya_y + (k-1-y)\nu \colon y \in [k-1]\} \big) \le j\nu + \max\{p_ya_y \colon y \in [k-1]\}.
    \]
Since $\avect$ satisfies (Con\ref{item:Con4}) for~$2$ by hypothesis, 
each $p_y a_y \leq 2 j^2$, so
it follows that $p_{k-1}h_{k-1}(t') \le j\nu + 2j^2 < 3j^2$, and hence
    \begin{equation}\label{eq:af-heavy-bins-middle-mu-1}
        \sum_{t'=0}^{T_0-\tau-1} p_{k-1}h_{k-1}(t') \le 3j^2(T_0-\tau).
    \end{equation}

Moreover, by~\eqref{eq:upperbound} applied with $x=k-1$ (using $P_{\ell(k-1)} \leq P_{k-1}$),  
    \begin{equation}\label{eq:af-heavy-bins-middle-mu-2}
        \sum_{t'=T_0-\tau}^{t-\tau-1} p_{k-1}h_{k-1}(t') \le (t-T_0)\Big(\nu(k-1) + \Big(1+\frac{1}{32j^2}\Big)P_{k-1}\Big) < (t-\tau)\Big(1 + \frac{1}{16j^2}\Big)P_{k-1}.
    \end{equation}
    Combining~\eqref{eq:af-heavy-bins-middle-mu-1} and~\eqref{eq:af-heavy-bins-middle-mu-2} yields
    \[
        \mu_1 \le (t-\tau)\nu + 3j^2(T_0-\tau) + (t-\tau)\Big(1 + \frac{1}{16j^2}\Big)P_{k-1}.
    \]
    Since $t - \tau \ge \varpi j^{12} = (T_0-\tau)j^5$ by hypothesis, and $\nu = 1/j^{100}$, it follows that 
    \begin{equation}\label{eq:af-heavy-bins-mu-1-final}
    \mu_1 \le (t-\tau)(1 + 1/(12j^2))P_{k-1}.
    \end{equation}
    
    We now substitute~\eqref{eq:af-heavy-bins-mu-1-final} into~\eqref{eq:af-heavy-bins-middle} to obtain
    \[
Z_1 \le j(\log j)^2 + \Big(1 + \frac{j}{((t-\tau)(1 + 1/(12j^2))P_{k-1})^{1/2}}\Big)
(t-\tau)
\Big(1 + \frac{1}{12j^2}\Big)P_{k-1}
    \]
Since $t-\tau \ge  \varpi j^{12}$ and $j$ is large, it follows that, with probability at least $1-j^{-3j}$ (where the failure probability is inherited from~\eqref{eq:af-heavy-bins-middle}),
    \begin{equation}\label{eq:af-heavy-bins-middle-goal}
Z_1 \le (t-\tau)\Big(1 + \frac{1}{10j^2}\Big)P_{k-1}.
    \end{equation}

    \medskip\noindent\textbf{Bounding the right term below:} For brevity, let
    \begin{align*}
        Z_2 = \Big|\bigcup_{t'=\tau+1}^t \s_k^{E^+}(t') \Big|,\qquad \mu_2 = p_k\sum_{t'=0}^{t-\tau-1} h_k(t').
    \end{align*}
Then by Lemma~\ref{lem:h-conc-sends}, 
as long as $\mu_2 \geq j(\log j)^2$,
with probability at least $1-j^{-3j}$,
\begin{equation}\label{eq:af-heavy-bins-right}
        Z_2 \ge \Big(1 - \frac{j}{\mu_2^{1/2}}\Big)\mu_2.
    \end{equation}

In order to use Equation~\eqref{eq:af-heavy-bins-right} we must
establish $\mu_2 \geq j(\log j)^2$. In any case, we will need a tighter lower bound on~$\mu_2$. The following claim will be useful.

\medskip\noindent
{\bf Claim.} For all $t'$ with 
$T_0 -\tau  \leq t' \leq t-\tau$,
$h_k(t') > (1+1/(4k^2)) P_{k-1} W_k$.   

\medskip\noindent
{\bf Proof of Claim.}  
We first apply Lemma~\ref{lem:h-soft-nonincreasing}, taking the $t$ of Lemma~\ref{lem:h-soft-nonincreasing} to be our $t-\tau$, $T' = t'$, and the $T$ of the lemma to be our~$T_0-\tau$. 
Note that $0 \leq T_0-\tau \leq t' \leq t-\tau$
so the lemma applies.
Defining $M_k = \max(\{k\nu\} \cup \{p_yh_y(T_0-\tau) + (k-y)\nu\colon y \in [k-1]\})$ as in the lemma statement, it follows from the lemma that $h_k(t-\tau) \le \max\{h_k(t'),M_k W_k\}$. By the hypothesis of Case 2, 
$h_k(t-\tau) > (1+1/(4k^2))P_{k-1}W_k$, so it suffices to prove that
$M_k < (1+1/(4k^2))P_{k-1}$.
This would imply 
that the maximum is taken at $h_k(t')$,
so $h_k(t-\tau) \leq h_k(t')$.
We will therefore finish the proof of the claim 
by showing that $M_k < (1+1/(4k^2))P_{k-1}$.
It is clear that $k  \nu = k\nu_j^\low = k/j^{100} \leq \xi \leq P_{k-1}$ 
so, from the definition of~$M_k$, it suffices to prove that for all $y\in [k-1]$, 
$p_y h_y(T_0-\tau) + (k-y) \nu \leq (1+1/(4k^2))P_{k-1}$.
By~\eqref{eq:upperbound},
$p_{y} h_{y} (T_0-\tau) + (k-y)\nu \le k \nu + 
(1 +  1/(32j^2) )P_{\ell(y)}$, and now the claim is easy since $P_{\ell(y)} \leq   P_{k-1}$.
For example, $(1 +  1/(32j^2) )P_{\ell(y)} \leq (1+1/(32k^2))P_{k-1}$
and $k\nu \leq P_{k-1}/(32 k^2)$.\qed
 
We now use the claim to derive a lower bound on~$\mu_2$. Given the claim, 

\[
\mu_2 \ge p_k\sum_{t'=T_0-\tau}^{t-\tau-1} h_k(t')   > p_k(t-T_0)\Big(1+\frac{1}{4k^2}\Big)P_{k-1}W_k = (t-T_0)\Big(1+\frac{1}{4k^2}\Big)P_{k-1}.
    \]
Since $t-\tau \ge \varpi j^{12} = (T_0-\tau)j^5$ by hypothesis, we have $t-T_0 = t-\tau - (T_0-\tau) \ge (1-1/j^3)(t-\tau)$, and hence
    \begin{equation}\label{eq:af-heavy-bins-mu-2-final}
        \mu_2 \ge (1-1/j^3)(t-\tau)(1+1/(4k^2))P_{k-1} \ge (t-\tau)(1 + 1/(5k^2))P_{k-1}.
    \end{equation}
Note that Equation~\eqref{eq:af-heavy-bins-mu-2-final} implies $\mu_2 \geq j(\log j)^2$ which is required for~\eqref{eq:af-heavy-bins-right}.
Note also that the
bound in~\eqref{eq:af-heavy-bins-right}
is increasing in~$\mu_2$.
Substituting the lower bound from 
Equation~\eqref{eq:af-heavy-bins-mu-2-final} into~\eqref{eq:af-heavy-bins-right} and applying 
$(t-\tau)^{1/2} \geq  j^{6}$, we obtain that with probability at least $1-j^{-3j}$,
\begin{align}\nonumber
Z_2 &\geq \bigg(1 - \frac{j}{{(t-\tau)}^{1/2}(1+1/(5k^2))^{1/2}P_{k-1}^{1/2}}\bigg)(t-\tau) (1+1/(5k^2))P_{k-1}\\\label{eq:af-heavy-bins-right-goal}
&\ge (1-1/j^4) (t-\tau) (1 + 1/(5k^2))P_{k-1} \ge (t-\tau)(1 + 1/(6k^2))P_{k-1}.
\end{align}

\medskip\noindent \textbf{Combining the bounds:} 
Since $\avect$ satisfies (Con\ref{item:Con1}) for $d$, $a_k \le P_kW_k + d$. Combining this with~\eqref{eq:af-heavy-bins-empty-0}, \eqref{eq:af-heavy-bins-middle-goal} and~\eqref{eq:af-heavy-bins-right-goal} and a union bound, we obtain that with probability $1 - 2j^{-3j}$,
\begin{align*}
b_k^{E^+}(t) &\le P_kW_k + d + (t-\tau)\Big(1 + \frac{1}{10j^2}\Big)P_{k-1} 
- (t-\tau)\Big(1 + \frac{1}{6k^2}\Big)P_{k-1}\\
&\le P_kW_k + d - (t-\tau)P_{k-1}/(15j^2).
    \end{align*}

    Since $t-\tau \ge dj^3$, it follows that $b_k^{E^+}(t) \le P_kW_k$ as required.
    
It now remains only to prove that $h_k(t-\tau) < P_kW_k$. Observe from the definition of $h$  (Definition~\ref{def:h}) that 
$h_k(0)=a_k$  
and 
for all $t'\geq 1$,
$h_k(t') - h_k(t'-1) = -p_k h_k(t'-1) + p_{k-1}h_{k-1}(t'-1) + \nu$.
By a telescoping sum,
\begin{align*}
h_k(t-\tau) - h_k(0)&= \sum_{t'=1}^{t-\tau} h_k(t') - h_k(t'-1)\\ 
&= \sum_{t'=1}^{t-\tau} (-p_k h_k(t'-1) + p_{k-1}h_{k-1}(t'-1) + \nu)\\
&=  -\sum_{t'=1}^{t-\tau} p_k h_k(t'-1) 
+ (t-\tau) \nu + \sum_{t'=1}^{t-\tau} 
+ p_{k-1}h_{k-1}(t'-1)  \\
&= -\mu_2 + \mu_1,
\end{align*}
So $h_k(t-\tau) = a_k + \mu_1 - \mu_2$. Then applying our existing bounds~\eqref{eq:af-heavy-bins-mu-1-final} for $\mu_1$ and~\eqref{eq:af-heavy-bins-mu-2-final} for $\mu_2$,  we obtain
\begin{align*}
h_k(t-\tau) &\le a_k + (t-\tau)\Big(1+\frac{1}{12j^2}\Big)P_{k-1} - (t-\tau)\Big(1+\frac{1}{5k^2}\Big)P_{k-1}\\
        &< a_k - (t-\tau)P_{k-1}/(12j^2) < P_kW_k,
    \end{align*}
    as required.
\end{proof}

\begin{lemma}\label{lem:af-constant-stabilise}
Fix a send sequence $\p$ with $p_0=1$. 
Fix a sufficiently small $\lambda \in (0,1/3)$.
Let $j$ be a sufficiently large integer. Fix an integer $\tau \ge 0$ and a tuple $\avect \in (\integers_{\ge 0})^j$. Let $E^+$ be a constant-escape process on $[j]$ with send sequence~$\p$, start time $\tau$, end-of-step arrival rate $\nu_j^\low$, and $\romanarrvect_{[j]}^{E^+}(\tau) = \avect$. Let $d$, $\varpi$, and $\theta$ be real numbers with $d\ge 0$, $\varpi\ge 1$, and $0 \le \theta < 1/(64j^2)$. Suppose that $\avect$ satisfies (Con\ref{item:Con1}) for $d$, (Con\ref{item:Con3}) for $(\varpi,\theta,\nu_j^\low)$, and (Con\ref{item:Con4}) for~$2$. 
Fix $t$ in the range 
$\max\{\varpi j^{12}, j^{17},dj^3\} \le t -\tau \le  (4/\lambda)^{2j}$. Then, with probability at least $1-e^{-450j/\lambda^2}$, $\bvect_{[j]}^{E^+}(t)$ is a $j^7$-AF-good beginning.
\end{lemma}
\begin{proof}
For all $x \in [j]$ and all $t'$, we suppress $\nu_j^\low$ and $\avect$ from $h_x^{\nu_j^\low,\avect}(t')$, and we write $\nu = \nu_j^\low$.  
Fix $t$ in the range 
$\max\{\varpi j^{12}, j^{17},dj^3\} \le t -\tau \le  (4/\lambda)^{2j}$.
 
 By Lemma~\ref{lem:af-heavy-bins-empty}, with failure probability at most $2j^{-3j+1}$, 
 the following holds.
 \begin{equation}\label{eq:firstone} 
\mbox{For all $k \in [j]$ with $W_k \ge j^{10}$,  $b_k^{E^+}(t) \le P_kW_k$ and $h_k(t-\tau) < P_k W_k$.} 
    \end{equation}

We next claim that, for all bins $k \in [j]$ with $W_k < j^{10}$,
$h_k(t-\tau) \le P_kW_k$. 
Let $\varpi' = \max\{\varpi, j^{10}\}$. From Definition~\ref{def:conditions}
the definition of (Con\ref{item:Con3}) is increasing in~$\varpi$, so
$\avect$ satisfies (Con\ref{item:Con3}) for 
$(\varpi',\theta,\nu)$.
For all $x\in [j]$ let $\ell(x) = \max(\{0\} \cup \{y \in [x]\colon W_y \ge \varpi'\})$. 
For any $k\in [j]$ with $W_k < j^{10}$, 
$\ell(k) \leq k-1$.  
Now
by Lemma~\ref{lem:escape-means-stabilise}, for all $k\in [j]$ with $W_k < j^{10}$,  
\[
p_kh_k(t-\tau) \le k\nu + \Big(1 + \theta + \frac{2\varpi' j^4}{\xi\lambda (t-\tau)}\Big)P_{\ell(k)} 
\le j^{-99} + \Big(1 + \frac{1}{64j^2} + \frac{2}{\xi \lambda j^3}\Big)P_{k-1} \le P_k.
    \]
After multiplying through by $W_k$, we obtain that for all $k \in [j]$ with $W_k < j^{10}$, $h_k(t-\tau) \le P_kW_k$. Combining this with Equation~\eqref{eq:firstone}, we obtain,
 with failure probability at most $2j^{-3j+1}$, 
    \begin{equation}\label{eq:af-constant-stabilise-1}
\mbox{for all bins $k \in [j]$, $h_k(t-\tau) \le P_kW_k$.}
    \end{equation}

We next lower-bound the probability that $\bvect_{[j]}^{E^+}(t)$ satisfies (Con\ref{item:Con1}) for $j^7$, as in Definition~\ref{def:conditions}. By~\eqref{eq:firstone}, with probability at least $1-2j^{-3j+1}$,  
    \[
        \excess(\bvect_{[j]}^{E^+}(t),1) = \sum_{k=1}^j \max\{0,b_k^{E^+}(t) - P_kW_k\} \le \sum_{k\in [j] \colon W_k < j^{10}} \big(b_k^{E^+}(t) - P_kW_k\big)
    \]
    By~\eqref{eq:af-constant-stabilise-1},       for all $k \in [j]$,
    $h_k(t-\tau) \le P_kW_k$. Hence applying Lemma~\ref{lem:h-basic-chernoff}   to each term of the sum with $\alpha=1$, taking a union bound,
    we find that, with probability at least $1-j^{-3j+1}$,
    \begin{align*}
\excess(\bvect_{[j]}^{E^+}(t),1) &\le \sum_{k\in [j] \colon W_k < j^{10}} \max\Big\{j(\log j)^2, \frac{j^{1/2}\log j}{4(P_kW_k)^{1/2}}P_kW_k\Big\}\\
&\le \sum_{k\in [j] \colon W_k < j^{10}} \max\Big\{j(\log j)^2, j^{1/2}(\log j) W_k^{1/2} \Big\} \leq j^7,
    \end{align*}
    as required. By a union bound, we have shown that with probability at least $1-3j^{-3j+1}$, $\bvect_{[j]}^{E^+}(t)$ satisfies (Con\ref{item:Con1}) for $j^7$.

    By~\eqref{eq:af-constant-stabilise-1}, we may apply Corollary~\ref{cor:applied-excess-conc} to prove that $\bvect_{[j]}^{E^+}(t)$ satisfies (Con\ref{item:Con2}) with probability at least $1-e^{-500j/\lambda^2}$. The result therefore follows by a union bound and by Definition~\ref{def:good-start}.     
\end{proof}

The following lemma will be useful when we are coupling an escape process $E$ whose high-level state $\Psi$ is advancing or filling to a constant-escape process $E^+$ with arrival rate $\nu_j^\low$. In order to maintain the coupling, we will need to prove that some set in $\calS^\Psi$ maintains low expected noise in $E^+$. Initially, if $E^+$ starts at an AF-good beginning,  this will be guaranteed by (Con2). The following lemma will allow us to guarantee it at all time steps.
 
\begin{lemma}\label{lem:af-constant-strong-survive}
Fix a send sequence $\p$ with $p_0=1$. 
Fix a sufficiently small $\lambda \in (0,1/3)$.
Let $j$ be a sufficiently large integer. Fix an integer $\tau \ge 0$ and a tuple $\avect \in (\integers_{\ge 0})^j$. Let $E^+$ be a constant-escape process on $[j]$ with send sequence~$\p$, start time $\tau$, end-of-step arrival rate $\nu_j^\low$, and $\romanarrvect_{[j]}^{E^+}(\tau) = \avect$.
Let $d \ge 0$ and suppose that $\avect$ is a $d$-AF-good beginning. Then, with probability at least $1-e^{-450j/\lambda^2}$, the following holds: For all $t$ satisfying $\tau \le t \le \tau + (4/\lambda)^{2j}$, 
$\excess(\bvect_{[j]}^{E^+}(t),3\xi^{-2/3}) \le \lambda j/250$.
\end{lemma}
\begin{proof}
For brevity, let $\nu = \nu_j^\low$. 
Since $\lambda$ is
small we can take a union bound over~$t$,
showing that for every fixed $t$ in the range  $\tau \le t \le \tau + (4/\lambda)^{2j}$,
\begin{equation}\label{eq:af-constant-strong-survive-0}
\pr\big(\excess(\bvect_{[j]}^{E^+}(t),3\xi^{-2/3}) > \lambda j/250\big) < e^{-500j/\lambda^2}.
\end{equation}
    
As such, fix $t$ in the range $\tau \le t \le \tau + (4/\lambda)^{2j}$. 
Our aim will be to prove~\eqref{eq:af-constant-strong-survive-0}.
We first divide the balls of $E^+$ into two cohorts, coupling $E^+$ with two constant-escape processes $E_1^+$ and $E_2^+$ (one for each cohort). We then analyse $E_1^+$ and $E_2^+$ separately.

For each $k \in [j]$, let $\hat{a}_k = \min\{a_k, \lceil \xi^{-1/3}P_kW_k \rceil\}$, and let $\tilde{a}_k = a_k - \hat{a}_k$. Let $E_1^+$ be a constant-escape process on $[j]$ with send sequence~$\p$, start time $\tau$, end-of-step arrival rate $\nu_j^\low$, and $\romanarr_k^{E_1^+}(\tau) = \hat{a}_k$ for all $k \in [j]$. Let $E_2^+$ be a constant-escape process on $[j]$ with send sequence~$\p$, start time $\tau$, end-of-step arrival rate $0$, and $\romanarr_k^{E_2^+}(\tau) = \tilde{a}_k$ for all $k \in [j]$. We couple $E^+$, $E_1^+$ and $E_2^+$ in the natural way, synchronising sends and directing all end-of-step arrivals in $E^+$ to $E_1^+$, so that
    \[
        \mbox{for all $t' \ge \tau$ and $k \in [j]$, $\b_k^{E^+}(t') = \b_k^{E_1^+}(t') \cup \b_k^{E_2^+}(t')$ and $\b_k^{E_1^+}(t') \cap \b_k^{E_2^+}(t') = \emptyset$.}
    \]
    Observe that for all $t' \ge \tau$,
    \begin{align}\nonumber
        \excess(\bvect_{[j]}^{E^+}(t'),3\xi^{-2/3}) 
        &= \sum_{k=1}^j \max\{0, b_k^{E_1^+}(t') + b_k^{E_2^+}(t') -3\xi^{-2/3} P_kW_k\} \\\nonumber
        &\le \sum_{k=1}^j \Big(\max\{0, b_k^{E_1^+}(t')-3\xi^{-2/3} P_kW_k\} + b_k^{E_2^+}(t')\Big)\\\label{eq:af-constant-strong-survive-1}
        &= \excess(\bvect_{[j]}^{E_1^+}(t'), 3\xi^{-2/3}) + \sum_{k=1}^j b_k^{E_2^+}(t').
    \end{align}

We next bound $\excess(\bvect_{[j]}^{E_1^+}(t), 3\xi^{-2/3})$ above using Lemma~\ref{lem:applied-excess-conc}.
For each $i\in [j]$, the binomial random variable $A_i$ from that lemma
is deterministic with value $a_i$
(so its mean $m_i$ is also equal to $a_i$).
We take the parameter
$\zeta$ of that lemma to be $3\xi^{-1/3}$, and the $\eta$ of that lemma to be $3\xi^{-2/3}$. In order to apply Lemma~\ref{lem:applied-excess-conc}, we must prove that for 
our fixed $t\geq \tau$ and for all $k\in [j]$,  $h_k^{\nu,\avect}(t-\tau) \le 3\xi^{-1/3}P_kW_k$. By Lemma~\ref{lem:h-below-stationary} applied with $T=0$ and $x=k$,
\begin{align*}
p_kh_k^{\nu,\avect}(t-\tau) &\le \max(\{k\nu\} \cup \{p_y\hat{a}_y + (k-y)\nu \colon y \in [k]\})\\
&\le k\nu + \max\{p_y\hat{a}_y \colon y \in [k]\}\\
&= k\nu + \max\{p_y  \min\{a_y,\lceil\xi^{-1/3}P_yW_y\rceil\}\colon y \in [k]\}\\
&\le 
k\nu+ 
\max\{p_y\lceil\xi^{-1/3}P_yW_y\rceil\colon y \in [k]\} 
\\
&\le k\nu +  \max\{2\xi^{-1/3}P_y\colon y \in [k]\} = 
k\nu + 2\xi^{-1/3}P_k \leq 3 \xi^{-1/3} P_k,
\end{align*}
and hence $h_k^{\nu,\avect}(t-\tau) \le 3\xi^{-1/3}P_kW_k$ as required.
Thus by Lemma~\ref{lem:applied-excess-conc}, with probability at least $1-e^{-500j/\lambda^2}$,
    \begin{equation}\label{eq:af-constant-strong-survive-2}
        \excess(\bvect_{[j]}^{E_1^+}(t), 3\xi^{-2/3}) \le \lambda j/500.
    \end{equation}

We next bound $\sum_{k=1}^j b_k^{E_2^+}(t)$  from~\eqref{eq:af-constant-strong-survive-1} above. Since $E_2^+$ has end-of-step arrival rate zero,
    \begin{align*}
    \sum_{k=1}^j b_k^{E_2^+}(t) &\le \sum_{k=1}^j b_k^{E_2^+}(\tau) = \sum_{k=1}^j \tilde{a}_k = \sum_{k=1}^j \max\{0, a_k - \lceil \xi^{-1/3}P_kW_k\rceil \}\\
    &\le \sum_{k=1}^j \max\{0, a_k - \xi^{-1/3}P_kW_k\} = \excess(\avect, \xi^{-1/3}).
    \end{align*}
    Since $\avect$ is a $d$-AF-good beginning (Definition~\ref{def:good-start}), it follows from (Con\ref{item:Con2}) that $\sum_{k=1}^j b_k^{E_2^+}(t) \le \lambda j/500$. Combining this with~\eqref{eq:af-constant-strong-survive-1} and~\eqref{eq:af-constant-strong-survive-2} yields $\excess(\bvect_{[j]}^{E^+}(t),2\xi^{-2/3}) < \lambda j/250$, as required by~\eqref{eq:af-constant-strong-survive-0}.
\end{proof}

The following lemma's purpose is analogous to that of Lemma~\ref{lem:af-constant-strong-survive}, but will be used when we are coupling a constant-escape process with arrival rate $\nu_j^\low$ to an escape process with a $j$-stabilising state rather than a $j$-advancing or -filling state.

\begin{lemma}\label{lem:stab-constant-weak-survive}
Fix a send sequence $\p$ with $p_0=1$. 
Fix a sufficiently small $\lambda \in (0,1/3)$.
Let $j$ be a sufficiently large integer. Fix an integer $\tau \ge 0$ and a tuple $\avect \in (\integers_{\ge 0})^j$.
Let $E^+$ be a constant-escape process on $[j]$ with send sequence~$\p$, start time $\tau$, end-of-step arrival rate $\nu_j^\high$, and $\romanarrvect_{[j]}^{E^+}(\tau) = \avect$. 
Fix real numbers~$d$ and~$\theta$ 
satisfying  $d \ge 0$ and $0 \le \theta \leq  1/(64j^2)$. Suppose that $\avect$ is a $(d,j^{10},\theta,1)$-RS-good-start. Let $S = [(\log j)^2]$. 
Fix $t$ in the range $\tau + j^{17} \leq t \leq \tau + (4/\lambda)^{2j}$. Then,  with probability at least $1-\exp(-(\log j)^2)$, $\calN_S^{E^+}(t+1) \le \lambda |S|/160$.
\end{lemma}
\begin{proof}
For brevity, let $\nu = \nu_j^\high$, and suppress $\avect$ and $\nu$ from $h_k^{\avect,\nu}(t)$. 
Since $\avect$ is a $(d,j^{10},\theta,1)$-RS-good-start it satisfies (Con\ref{item:Con3}) for $(j^{10},\theta,\nu_j^{\high})$ and (Con\ref{item:Con4}) for~$1$ (hence, also for~$2$). Thus, we may apply 
Lemma~\ref{lem:escape-means-stabilise}.
For all $x\in [j]$ 
let $\ell(x) = \max(\{0\} \cup \{y \in [x]\colon W_y \ge \varpi\})$. Applying the lemma to~$t-\tau$ with $\beta=j^{10}$, for all $x\in S$,  
    \[
p_xh_x(t-\tau) \le x\nu + \Big(1 + \theta + \frac{2j^{14}}{\xi \lambda (t-\tau)}\Big)P_{\ell(x)} \le \nu(\log j)^2 + \Big(1 + \frac{1}{64j^2} + \frac{2}{\xi\lambda j^3}\Big)P_x.
    \]
Since $\nu(\log j)^2 = (\log j)^{-98}$ and $P_x \ge \xi$, it follows that $p_xh_x(t-\tau) \le 2P_x$ for all $x \in S$. Hence by Observation~\ref{obs:constant-escape}, for all $k \in S$, $\E(b_k^{E^+}(t)) \le 2P_kW_k$. By Lemma~\ref{lem:constant-escape-gen-multi}, $(b_k^{E^+}(t)\colon x \in S)$ can be expressed as coordinates of a generalised multinomial random variable. 
We can therefore apply Lemma~\ref{lem:excess-chernoff}, taking the $Z_1,\dots,Z_n$ of that lemma to be $b_1^{E^+}(t),\ldots,b_{|S|}^{E^+}(t)$, the $\muvect$ of that lemma to be given by $\mu_k = 2P_kW_k$, the $\delta$ of that lemma to be $\xi^{-1/2} > e$, and the $z$ of that lemma to be $\lambda |S|/500$. 
    Let $Z^+ = \sum_{x \in S}\max\{0,b_k^{E^+}(t)-2\xi^{-1/2}P_kW_k\}$. Then by Lemma~\ref{lem:excess-chernoff} and Definition~\ref{def:escape-constants},
    \[
        \pr(Z^+ > \lambda |S|/500) \le 2^{|S|}\xi^{\lambda |S|/500} = \exp\Big(\Big(\log 2 - \frac{\lambda}{500}\log(1/\xi)\Big)|S|\Big) < \exp(-2|S|) < \exp(-(\log j)^2).
    \]
If $Z^+ \leq \lambda |S|/500$,  then by Remark~\ref{rem:P-bounded} and Definition~\ref{def:escape-constants},
\begin{align*}
\calN_S^{E^+}(t+1) 
&= \sum_{x \in S} p_kb_k^{E^+}(t) 
\le \sum_{x \in S}p_k\Big(2\xi^{-1/2}P_kW_k + \max\{0, b_k^{E^+}(t) - 2\xi^{-1/2}P_kW_k\}\Big)\\
&\le 2\xi^{-1/2}\sum_{x \in S}P_k + Z^+ 
\le 8\xi^{1/2}|S| + \lambda|S|/500 
< \lambda |S|/160.
\end{align*}
The result therefore follows.
\end{proof}

 The following two lemmas are upper bounds on how much the parameters of an RS-good beginning can degrade as a constant-escape process $E^+$ runs. Lemma~\ref{lem:rs-constant} shows that the parameter $d$ of (Con1) degrades at worst linearly as time goes on, the parameters of (Con3) degrade by at most a small amount, and the paramter of (Con4) improves over time. Lemma~\ref{lem:stab-strong-con1} shows further that if the arrival rate of $E^+$ is $\nu_j^\low$, then the total degradation of (Con1) (i.e.\ the amount $d$ increases by as $E^+$ runs) is bounded.

\begin{lemma}\label{lem:rs-constant}
Fix a send sequence $\p$ with $p_0=1$. 
Fix a sufficiently small $\lambda \in (0,1/3)$.
Let $j$ be a sufficiently large integer. Fix an integer $\tau \ge 0$ and a tuple $\avect \in (\integers_{\ge 0})^j$.
Fix a real number $\nu$ in the range $0 \le \nu \le \nu_j^\high$. 
Let $E^+$ be a constant-escape process on $[j]$ with send sequence~$\p$, start time $\tau$, end-of-step arrival rate $\nu$, and $\romanarrvect_{[j]}^{E^+}(\tau) = \avect$. 
Let $d$, $\varpi$, and $\theta$ be real numbers with $d\ge 0$, $\varpi\ge 1$, and $0 \le \theta < 1/(64j^2)$.
Suppose that $\avect$ is a $(d,\varpi,\theta, G)$-RS-good beginning
for a positive integer~$G\leq 2$. Then with probability at least $1-j^{-2j}$, the following properties all hold.
\begin{enumerate}[(i)]
\item For all $t$ in the range $\tau \le t \leq \tau + \varpi j^9$, $\bvect_{[j]}^{E^+}(t)$ is a $(d+5(t-\tau)j^3, \varpi j^{17}, \theta+3/j^4, G+1)$-RS-good beginning.
\item For all $t$ in the range $\tau+\varpi j^9 \le t \le \tau + (4/\lambda)^{2j}$, $\bvect_{[j]}^{E^+}(t)$ is a $(d+5(t-\tau)j^3, j^{10}, \theta+3/j^4, 1)$-RS-good beginning.
\end{enumerate}
\end{lemma}
\begin{proof}
Since $\avect$ is a  $(d,\varpi,\theta, G)$-RS-good beginning (Definition~\ref{def:good-start}),
it satisfies (Con\ref{item:Con1}) for~$d$,
(Con\ref{item:Con3}) for~$(\varpi, \theta,\nu_j^{\high})$, and (Con\ref{item:Con4}) for~$G$. We wish to re-establish these properties with different parameters. The following events will be helpful (the subscripts in the names of the events match 
(Con\ref{item:Con1})--(Con\ref{item:Con4}).) 
For all $k \in [j]$ and $t \ge \tau$, let $\calE_1(k,t)$ be the event that $b_k^{E^+}(t) \le a_k + 5(t-\tau)j^2$. Let $\calE_3(k,t)$ be the event that  $W_k < j^{10}$ or $b_k^{E^+}(t) \le k\nu_j^{\high} + (1+\theta+3/j^4)P_kW_k$. 
If $t < \tau+\varpi j^9$ then
let $\calE_4(k,t)$ be the event that $p_kb_k^{E^+}(t) \le (G+1)j^2$. Otherwise, let
$\calE_4(k,t)$ be the event that $p_kb_k^{E^+}(t) \le j^2$.  

We first bound $\pr(\calE_1(k,t))$ below. 
$\calE_1(k,t)$ holds by construction for $t=\tau$. Observe that for any $k \in [j]$ and $t \ge \tau+1$, every ball in $\b_k^{E^+}(t) \setminus \b_k^{E^+}(\tau)$ was either sent into bin $k$ from bin $k-1$, or it arrived in bin $k$ as an end-of-step arrival after step $\tau$. Thus for any $k \in [j]$ and $t \ge \tau+1$,
    \[
        b_k^{E^+}(t) - a_k \le \Big|\bigcup_{t'=\tau+1}^t \big(\s_{k-1}^{E^+}(t') \cup \arr_k^{E^+}(t')\big)\Big|.
    \]
    By Lemma~\ref{lem:h-conc-entries}, it follows for any fixed $k\in [j]$ and $t>\tau$ that, with probability at least $1-j^{-3j}$,
    \begin{equation}\label{eq:rs-constant-1}
        b_k^{E^+}(t) - a_k \le \max\Big\{j(\log j)^2,\Big(1 + \frac{j}{\mu^{1/2}}\Big)\mu\Big\} \mbox{, where } \mu = (t-\tau)\nu + \sum_{t'=0}^{t-\tau-1} p_{k-1}h^{\nu,\avect}_{k-1}(t').
    \end{equation}
    We next bound $\mu$ above. By Lemma~\ref{lem:h-below-stationary} applied with $T=0$ and $x=k-1$, for all $t' \ge 0$,
    \begin{align*}
        p_{k-1}h^{\nu,\avect}_{k-1}(t') &\le \max\big(\{(k-1)\nu\} \cup \{p_ya_y + (k-1-y)\nu\colon y \in [k-1]\}\big)\\
        &\le j\nu + \max\{p_ya_y\colon y \in [k-1]\}\\
        &\le j + \max\{p_ya_y\colon y \in [k-1]\}.
    \end{align*}
    Since $\avect$ is a $(d,\varpi,\theta,G)$-RS-good beginning (Definition~\ref{def:good-start}), it satisfies (Con\ref{item:Con4}) for $G$ 
so $p_y a_y \leq G j^2$ 
    and thus $p_{k-1}h^{\nu,\avect}_{k-1}(t') \le j + G j^2$.  
        Hence $\mu \le (t-\tau)\nu + (t-\tau)(j+G j^2) \le (G+1/10)(t-\tau)j^2$. Thus by~\eqref{eq:rs-constant-1}, for any fixed $k \in [j]$ and  $t \ge \tau+1$, 
 with probability at least $1-j^{-3j}$,   
    $b_k^{E^+}(t) - a_k \le j(\log j)^2 + 2 ((G+1/10)(t-\tau)j^2) \le 5(t-\tau)j^2$. We have therefore proved the following.
    \begin{equation}\label{eq:rs-constant-E1}
        \mbox{for all $k \in [j]$ and $t \geq \tau$, }\pr(\calE_1(k,t)) \ge 1 - j^{-3j}.
    \end{equation}

For the next two bounds we consider $k\in [j]$ and
$t$ in the range $\tau+\varpi j^9 \le t \le \tau + (4/\lambda)^{2j}$. We
will rely on
observation~\ref{obs:h-nu}, which shows that, for $t\geq \tau$, 
$h_k^{\nu,\avect}(t-\tau) \leq
h_k^{\nu_j^{\high},\avect}(t-\tau)$. 
We will apply Lemma~\ref{lem:escape-means-stabilise}. For this we need that $\avect$ satisfies (Con\ref{item:Con3}) for $(\varpi, \theta,\nu_j^{\high})$ and (Con\ref{item:Con4}) for~$2$  which holds since (Con\ref{item:Con4}) is increasing in~$G$ and $G\leq 2$.  For all $x\in [j]$ let $\ell(x) = \max(\{0\} \cup \{y \in [x]\colon W_y \ge \varpi\})$. 
By Lemma~\ref{lem:escape-means-stabilise} applied with $x=k$ and taking the $t$ of Lemma~\ref{lem:escape-means-stabilise} to be our $t-\tau$, which is at most $(4/\lambda)^{2j}$,
\begin{equation}\label{eq:usetwice}
p_k h_k^{\nu,\avect}(t-\tau) \leq
p_kh_k^{\nu_j^{\high},\avect}(t-\tau) \le k\nu_j^{\high} + \Big(1 + \theta + \frac{2\varpi j^4}{\xi\lambda (t-\tau)}\Big)P_{\ell(k)} \le k\nu_j^{\high} + \Big(1 + \theta + \frac{1}{j^4}\Big)P_k.
\end{equation}

For the first of the two bounds, we bound $\pr(\calE_3(k,t))$ below. We may assume $W_k \ge j^{10}$, or $\calE_3(k,t)$ occurs for all $t\geq t+1$. By~\eqref{eq:usetwice} and
Lemma~\ref{lem:h-basic-chernoff} applied with $\alpha = 1$, it follows that with probability at least $1-j^{-3j}$, 
\begin{align*}
b_k^E{^+}(t) &\le j(\log j)^2 + \Big(1 + \frac{j^{1/2}\log j}{4h^{\nu_j^{\high},\avect}_k(t-\tau)^{1/2}}\Big)h^{\nu_j^{\high},\avect}_k(t-\tau)\\
&\le j(\log j)^2 + \bigg(1 + \frac{j^{1/2}\log j}
{4 {\big( (k\nu_j^{\high} + (1 + \theta + 1/j^4)P_k)W_k\big)}^{1/2}}\bigg) \bigg(k\nu_j^{\high} + \Big(1 + \theta + \frac{1}{j^4}\Big)P_k\bigg)W_k.
    \end{align*}

Since $W_k \ge j^{10}$ and $P_k \ge \xi$, $(k\nu_j^{\high} + (1+\theta+1/j^4)P_k)W_k \ge \xi j^{10}$, so
    \begin{align*}
        b_k^E{^+}(t) &\le j(\log j)^2 + \Big(1 + \frac{j^{1/2}\log j}{4\xi^{1/2}j^5}\Big) \big(k\nu_j^{\high} + (1 + \theta + 1/j^4)P_k\big)W_k \\
        &< j(\log j)^2 + (1 + 1/j^4)\big(k\nu_j^{\high} + (1 + \theta + 1/j^4)P_k\big)W_k\\
        &= j(\log j)^2 + \bigg(k\nu_j^{\high} + \frac{k\nu_j^{\high}}{j^4} + \Big(1+ \theta + \frac{2}{j^4} + \frac{\theta + 1/j^4}{j^4}\Big)P_k\bigg)W_k \le k\nu_j^{\high} + (1+\theta + 3/j^4)P_kW_k.
    \end{align*}
    We have proved that
\begin{equation}\label{eq:rs-constant-E3}
\mbox{for all $k \in [j]$ and 
$t$ in the range $\tau+\varpi j^9 \le t \le \tau + (4/\lambda)^{2j}$, }\pr(\calE_3(k,t)) \ge 1 - j^{-3j}.
\end{equation}

For the second of the two bounds, we bound $\pr(\calE_4(k,t))$ below. Applying
Lemma~\ref{lem:h-basic-chernoff}   with $\alpha=1$, we get that, with probability at least $1-j^{-3j}$, $b_k^{E^+}(t) \le \max\{j(\log j)^2, 2h^{\nu,\avect}_k(t-\tau)\} $. Clearly $j(\log j)^2 \leq j^2 W_k$.
For the other term, 
by~\eqref{eq:usetwice}, 
$2h^{\nu,\avect}_k(t-\tau)
\leq 2(
k\nu_j^{\high} + (1 + \theta + {1}/{j^4})P_k
)\leq j^2 W_k$. So we have proved the following.

\begin{equation}\label{eq:rs-constant-E4-1}
\mbox{for all $k \in [j]$ and 
$t$ in the range $\tau+\varpi j^9 \le t \le \tau + (4/\lambda)^{2j}$, }
 \pr(\calE_4(k,t)) \ge 1-j^{-3j}.
\end{equation}

We next bound $\pr(\calE_4(k,t))$ below for $t$ in the range $\tau < t \le \tau + \varpi j^9$.  Recall that in this range the definition of $\calE_4(k,t)$ allows an extra factor of~$(G+1)$.
For all $k \in [j]$, by Lemma~\ref{lem:h-below-stationary} applied with $T=0$ and $x=k$ and taking the $t$ of Lemma~\ref{lem:h-below-stationary} to be our $t-\tau$,
 \[
p_kh^{\nu,\avect}_k(t-\tau) \le \max(\{k\nu\} \cup \{p_ya_y + (k-y)\nu \colon y \in [k]\}) \le \max\{p_ya_y\colon y \in [k]\} + k\nu.
    \]
Since $\avect$ satisfies (Con\ref{item:Con4}) for $G$, $p_y a_y \leq G j^2$ and
it follows that $p_kh^{\nu,\avect}_k(t-\tau) \le G j^2 + k\nu \le (G+1/2)j^2$ and hence $h^{\nu,\avect}_k(t-\tau) \le (G+1/2)j^2W_k$. By Lemma~\ref{lem:h-basic-chernoff} applied with $\alpha=1$, it follows that with probability at least $1-j^{-3j}$, $b_k^{E^+}(t) \le (1+1/j^{1/4}) (G+1/2)j^2W_k < (G+1)j^2W_k$. Thus
\begin{equation}\label{eq:rs-constant-E4-2}
\mbox{for all $k\in [j]$ and $t$ in the range $\tau \le t \le \tau + \varpi j^9$, }\pr(\calE_4(k,t)) \ge 1-j^{-3j}.
    \end{equation}
    
Now suppose that:
\begin{itemize}
\item For all $k \in [j]$ and all $t$ in the range $\tau \le t \le \tau+(4/\lambda)^{2j}$, $\calE_1(k,t)$ occurs.
\item For all $k \in [j]$ and all $t$ in the range $\tau+\varpi j^9 \le t \le \tau+(4/\lambda)^{2j}$, $\calE_3(k,t)$ occurs.
\item For all $k \in [j]$ and all $t$ in the range $\tau \le t \le \tau+(4/\lambda)^{2j}$, $\calE_4(k,t)$ occurs.
\end{itemize}
    
By~\eqref{eq:rs-constant-E1}--\eqref{eq:rs-constant-E4-2} and a union bound, these events occur simultaneously with probability at least \[
1 - 3j((4/\lambda)^{2j}+1)j^{-3j} \geq
1-j^{-2j}. \]

It now suffices to prove that (i) and (ii) of the lemma statement hold.
We first prove that for all $t$ in the range $\tau \le t \le (4/\lambda)^{2j}$, $\bvect_{[j]}^{E^+}(t)$ satisfies (Con\ref{item:Con1})
for $d+5(t-\tau)j^3$. 
Consider any $t$ in the given range.   Since $\calE_1(k,t)$ occurs for all $k\in [j]$, $b_k^{E^+}(t) \le a_k + 5(t-\tau)j^2$, so by Observation~\ref{obs:dumb},
    \begin{align*}
        \excess(\bvect_{[j]}^{E^+}(t), 1) &= \sum_{k=1}^j \max\{0, b_k^{E^+}(t) - P_kW_k\}\\
        &\le \sum_{k=1}^j \max\{0, a_k - P_kW_k\} + \sum_{k=1}^j \max\{0,b_k^{E^+}(t) - a_k\} \le \excess(\avect, 1) + 5(t-\tau)j^3.
    \end{align*}
    Since $\avect$ satisfies (Con\ref{item:Con1}) for $d$ by hypothesis, it follows that $\bvect_{[j]}^{E^+}(t)$ satisfies (Con\ref{item:Con1}) for $d + 5(t-\tau)j^3$.

    We next prove that for all $t$ in the range $\tau \le t \leq \tau+\varpi j^9$, $\bvect_{[j]}^{E^+}(t)$ satisfies (Con\ref{item:Con3}) for $(\varpi j^{17}, \theta+3/j^4, \nu_j^\high)$. 
    We wish to show that, for any $k\in j$ with $W_k \geq \varpi j^{17}$, $p_k b_k^{E^+}(t) \leq \nu_j^{\high} k + (1+ \theta+ 3/j^4) P_k$.
Consider any
$k \in [j]$ with $W_k \ge \varpi j^{17}$ and 
any $t$ in the range. Since $\calE_1(k,t)$ occurs, 
\[
p_k b_k^{E^+}(t) \le p_k a_k + p_k 5(t-\tau)j^2 \le p_k a_k + p_k 5\varpi j^{11}.
\]
Since $\avect$ satisfies (Con\ref{item:Con3}) for $(\varpi,\theta,\nu)$,
$p_k a_k \leq \nu_j^{\high} k + (1+\theta) P_k$.
Moreover,   $P_k \ge \xi$. It follows that
\[
p_k b_k^{E^+}(t) \le \nu_j^{\high} k  + (1+\theta)P_k  + p_k 5\varpi j^{11} \le \nu_j^\high k + (1+\theta+3/j^4)P_k.
\]
Thus $\bvect_{[j]}^{E^+}(t)$ satisfies (Con\ref{item:Con3}) for $(\varpi j^{17}, \theta+3/j^4, \nu_j^\high)$ as claimed.

We next prove (i) of the lemma statement. 
Consider any $t$ in the range $\tau \le t \leq \tau + \varpi j^9$. We have shown that 
$\bvect_{[j]}^{E^+}(t)$ satisfies (Con\ref{item:Con1}) for $d + 5(t-\tau)j^3$
and (Con\ref{item:Con3}) for $(\varpi j^{17},\theta+3/j^4,\nu_j^\high)$. Since $\calE_4(k,t)$ occurs for all $k \in [j]$, 
$p_kb_k^{E^+}(t) \le (G+1)j^2$. Thus, $\bvect_{[j]}^{E^+}(t)$ satisfies (Con\ref{item:Con4}) for $G+1$. 
We conclude that  $\bvect_{[j]}^{E^+}(t)$ is a $(d+5(t-\tau)j^3, \varpi j^{17}, \theta+3/j^4, G+1)$-RS-good beginning, and we have proved part (i) of the lemma statement.

Now fix $t$ with $\tau+\varpi j^9 \le t \le \tau + (4/\lambda)^{2j}$. We have already observed that $\bvect_{[j]}^{E^+}(t)$ satisfies (Con\ref{item:Con1}) for $d+5(t-\tau)j^3$. 
Since $\calE_3(k,t)$ 
occurs for all $k\in [j]$, every $k$ with $W_k \geq j^{10}$ has $b_k^{E^+}(t) \le k\nu_j^{\high} + (1+\theta+3/j^4)P_kW_k$. 
Thus, $\bvect_{[j]}^{E^+}(t)$ satisfies 
(Con\ref{item:Con3}) for $(j^{10}, \theta+3/j^4, \nu_j^\high)$. Since $\calE_4(k,t)$ occurs for all $k \in [j]$, 
$p_kb_k^{E^+}(t) \le j^2$ so
$\bvect_{[j]}^{E^+}(t)$ also satisfies  (Con\ref{item:Con4}) for $1$. Thus $\bvect_{[j]}^{E^+}(t)$ is a $(d+5(t-\tau)j^3, j^{10}, \theta+3/j^4, 1)$-RS-good beginning, and we have proved part (ii) of the lemma statement.
\end{proof}

The proof of Lemma~\ref{lem:stab-strong-con1} is similar to (but easier than) the proof of  Lemma~\ref{lem:af-constant-strong-survive}.

\begin{lemma}\label{lem:stab-strong-con1}
Fix a send sequence $\p$ with $p_0=1$. 
Fix a sufficiently small $\lambda \in (0,1/3)$.
Let $j$ be a sufficiently large integer. Fix an integer $\tau \ge 0$ and a tuple $\avect \in (\integers_{\ge 0})^j$.
Let $E^+$ be a constant-escape process on $[j]$ with send sequence~$\p$, start time $\tau$, end-of-step arrival rate $\nu_j^\low$, and $\romanarrvect_{[j]}^{E^+}(\tau) = \avect$. 
Let $d \ge 0$ be real, and suppose $\avect$ satisfies (Con\ref{item:Con1}) for $d$. 
Fix $t$ in the range $\tau \leq t \leq \tau + (4/\lambda)^{2j}$. Then, with probability at least $1-\exp(-450j/\lambda^2)$, $\bvect_{[j]}^{E^+}(t)$ satisfies (Con\ref{item:Con1}) for $d+j^{21}$.
\end{lemma}
\begin{proof}
For brevity, let $\nu = \nu_j^\low$, and suppress $\avect$ and $\nu$ from $h_k^{\avect,\nu}(t)$. Fix $t$ in the range $\tau \le t \le \tau + (4/\lambda)^{2j}$. Our aim will be to prove
\begin{equation}\label{eq:stab-strong-con1-0}
\pr\big(\excess(\bvect_{[j]}^{E^+}(t),1) > d+j^{21}\big) < e^{-450j/\lambda^2}.
\end{equation}
We first divide the balls of $E^+$ into two cohorts, coupling $E^+$ with two constant-escape processes $E_1^+$ and $E_2^+$ (one for each cohort). We then analyse $E_1^+$ and $E_2^+$ separately.

For each $k \in [j]$, let $\hat{a}_k = \min\{a_k, \lfloor P_kW_k \rfloor\}$, and let $\tilde{a}_k = a_k - \hat{a}_k$. Let $E_1^+$ be a constant-escape process on $[j]$ with send sequence~$\p$, start time $\tau$, end-of-step arrival rate $\nu_j^\low$, and $\romanarr_k^{E_1^+}(\tau) = \hat{a}_k$ for all $k \in [j]$. Let $E_2^+$ be a constant-escape process on $[j]$ with send sequence~$\p$, start time $\tau$, end-of-step arrival rate $0$, and $\romanarr_k^{E_2^+}(\tau) = \tilde{a}_k$ for all $k \in [j]$. We couple $E^+$, $E_1^+$ and $E_2^+$ in the natural way, synchronising sends and directing all end-of-step arrivals in $E^+$ to $E_1^+$, so that
\[
    \mbox{for all $t' \ge \tau$ and $k \in [j]$, $\b_k^{E^+}(t') = \b_k^{E_1^+}(t') \cup \b_k^{E_2^+}(t')$ and $\b_k^{E_1^+}(t') \cap \b_k^{E_2^+}(t') = \emptyset$.}
\]
Observe that for all $t' \ge \tau$,
\begin{align}\nonumber
    \excess(\bvect_{[j]}^{E^+}(t'),1) 
    &= \sum_{k=1}^j \max\{0, b_k^{E_1^+}(t') + b_k^{E_2^+}(t') - P_kW_k\} \\\nonumber
    &\le \sum_{k=1}^j \Big(\max\{0, b_k^{E_1^+}(t')- P_kW_k\} + b_k^{E_2^+}(t')\Big)\\\label{eq:stab-strong-con1-1}
    &= \excess(\bvect_{[j]}^{E_1^+}(t'), 1) + \sum_{k=1}^j b_k^{E_2^+}(t').
\end{align}

We next bound $\excess(\bvect_{[j]}^{E_1^+}(t), 1)$ above. Observe that $\bvect_{[j]}^{E_1^+}(\tau) = \overline{\hat{a}}$ satisfies (Con\ref{item:Con1}) for $0$, (Con\ref{item:Con3}) for $(1,0,\nu_j^\low)$, and (Con\ref{item:Con4}) for $1$. If $t \ge \tau+j^{17}$, then by Lemma~\ref{lem:af-constant-stabilise}, it follows that with probability at least $1-\exp(-450j/\lambda^2)$, $\bvect_{[j]}^{E^+}(t)$ is a $j^7$-AF-good beginning and therefore satisfies (Con\ref{item:Con1}) for $j^7$. If instead $t < \tau+j^{17}$, then by Lemma~\ref{lem:rs-constant}, with probability at least $1-j^{-2j}$, $b_{[j]}^{E^+}(t)$ is a $(5(t-\tau)j^3,j^{17},3/j^4,2)$-RS-good beginning and therefore satisfies (Con\ref{item:Con1}) for $5j^{20}$. Either way, with probability at least $1-\exp(-450j/\lambda^2)$,
\begin{equation}\label{eq:stab-strong-con1-2}
    \excess(\bvect_{[j]}^{E_1^+}(t),1) \le 5j^{20}.
\end{equation}

We next bound $\sum_{k=1}^j b_k^{E_2^+}(t)$  from~\eqref{eq:stab-strong-con1-1} above. Since $E_2^+$ has end-of-step arrival rate zero,
    \begin{align*}
    \sum_{k=1}^j b_k^{E_2^+}(t) &\le \sum_{k=1}^j b_k^{E_2^+}(\tau) = \sum_{k=1}^j \tilde{a}_k = \sum_{k=1}^j \max\{0, a_k - \lfloor P_kW_k\rfloor \}\\
    &\le \sum_{k=1}^j \max\{0, a_k - P_kW_k + 1\} \le \excess(\avect, 1) + j.
    \end{align*}
    Since $\avect$ satisfies (Con\ref{item:Con1}) for $d$, it follows that $\sum_{k=1}^j b_k^{E_2^+}(t) \le d+j$. Combining this with~\eqref{eq:stab-strong-con1-1} and~\eqref{eq:stab-strong-con1-2} yields $\excess(\bvect_{[j]}^{E^+}(t),1) \le d+j+5j^{20} \le d+j^{21}$, as required by~\eqref{eq:stab-strong-con1-0}.
\end{proof}

\subsection{Analysing escape processes with pseudorandom starts}\label{sec:escape-invariants}

We first give two lemmas that bound the amount of time that an escape process can be coupled with a constant-escape process. We will apply Lemma~\ref{lem:escape-nu-high-couple} to couple constant-escape processes with arrival rate $\nu_j^\high$ with escape processes with $j$-refilling or $j$-stabilising high-level states. Recall from Definition~\ref{def:constant-escape-coupling}
(which gives the details of the coupling) 
that, given an escape process~$E$ 
with high-level state~$\Psi$,
a time $\tau' \geq \tau^{\Psi}$, and a positive real number~$\nu$,
$ {T^E_{\nu,\tau'}}$ is the first time $t \ge \tau'$ such that $\nu^E(t) > \nu$. 
Definition~\ref{def:constant-escape-coupling} establishes useful invariants of the coupling up until~${T^E_{\nu,\tau'}}$.

\begin{lemma}\label{lem:escape-nu-high-couple}
Fix a send sequence $\p$ with $p_0 = 1$, and fix $\lambda \in (0,1)$ with $\lambda$ sufficiently small. Let $\Psi = (g,\tau,j,\zvect,\calS,\type)$ be a high-level state with $j$ sufficiently large and $\type \in \{\refilling,\stabilising\}$. Let $E$ be an escape process with high-level state $\Psi$, send sequence $\p$, and birth rate $\lambda$.  Fix $\tau' \ge \tau$. Then 
    \[
        T_{\nu_j^\high,\tau'} \ge \min\{t \ge \tau'\colon \mbox{for some $k \in [j-1]$ with $W_k \ge j^\chi$, }b_k^E(t) > \lambda W_k/80\}.
    \]
\end{lemma}
\begin{proof}
 Let $\calS' = \{\{k\} \colon k \in [j-1], W_k \ge j^\chi\}$
Since $\Psi$ is a refilling or stabilising state (Definitions~\ref{def:refilling-state} and~\ref{def:stabilising-state}),
bin~$j$ is weakly exposed and  $\calS' \subseteq \calS$.

Consider any $t \ge \tau'$, and suppose $b_k^E(t) \le \lambda W_k/80$ for all $k \in [j-1]$ with $W_k \ge j^\chi$; then to prove the lemma, it suffices to prove that $\nu^E(t) \le \nu_j^\high$.

For all $S=\{k\} \in \calS'$, $\calN_S^E(t+1) = p_kb_k^E(t) \le \lambda/80$ by hypothesis, so
 from Definition~\ref{def:escape}, 
$\calS' \subseteq \mathbb{S}^E(t)$  and   
$\{k \in [j-1] \colon W_k \ge j^\chi\}
\subseteq \hat{S}^E(t)$. 
    Since bin $j$ is weakly exposed, by Definitions~\ref{def:covered} and~\ref{def:noise} it follows that
    \[
        |\hat{S}^E(t)| \ge |\Upsilon_{j,\ge j^\chi}| \ge \CUpsilon \LL(j)/\lambda \ge (\CUpsilon\CSE/\lambda)\log\log j > 100/\lambda^2,
    \]
    so $\Xi^E(t) = \lambda |\hat{S}^E(t)|/80 \ge (\CUpsilon\CSE/80)\log\log j$. Hence by Definitions~\ref{def:escape},
    \ref{def:constants}, and~\ref{def:escape-constants},
    \[
        \nu^E(t) = \exp(-\Xi^E(t)/192) < \exp(-100\log\log j) = \nu_j^\high,
    \]
    as required.
\end{proof}

We will apply Lemma~\ref{lem:escape-nu-low-couple} to couple constant-escape processes with arrival rate $\nu_j^\low$ with escape processes with $j$-advancing, $j$-filling, or $j$-stabilising high-level states.

\begin{lemma}\label{lem:escape-nu-low-couple}
    Fix a send sequence $\p$ with $p_0 = 1$, and fix $\lambda \in (0,1)$ with $\lambda$ sufficiently small. Let $\Psi = (g,\tau,j,\zvect,\calS,\type)$ be a high-level state with $j$ sufficiently large. Let $E$ be an escape process with high-level state $\Psi$, send sequence $\p$, and birth rate $\lambda$. Suppose that $\calS$ contains a set $S$ with $|S| \ge C_\Upsilon \L(j)/(2\lambda)$. Fix $\tau' \ge \tau$. Then 
    \[
        T_{\nu_j^\low,\tau'} \ge \min\Big\{t \ge \tau'\colon \sum_{k \in S}p_kb_k^E(t) > \lambda|S|/80 \Big\}.
    \]
\end{lemma}
\begin{proof}
    Consider any $t \ge \tau'$, and suppose $\sum_{k \in S}p_kb_k^E(t) \le \lambda |S|/80$; then to prove the lemma, it suffices to prove that $\nu^E(t) \le \nu_j^\low$.

    By hypothesis, $\calN_S^E(t+1) = \sum_{k \in S}p_kb_k^E(t) \le \lambda |S|/80$, so from Definition~\ref{def:escape}, $S \in \mathbb{S}^E(t)$ and $S \subseteq \hat{S}^E(t)$. By Definition~\ref{def:noise},
    \[
        |S| \ge C_\Upsilon \L(j)/(2\lambda) \ge (C_\Upsilon \CNoise^2/(2\lambda)) \log j > 100/\lambda^2,
    \]
    so $\Xi^E(t) \ge \lambda |S|/80$. Thus by Definitions~\ref{def:constants}, \ref{def:noise} and~\ref{def:escape-constants},
    \[
        \nu^E(t) = \exp(-\Xi^E(t)/192) \le \exp(-\lambda |S|/(192 \times 80)) \le j^{-C_\Upsilon \CNoise^2/(192 \times 160)} < j^{-100} = \nu_j^\low,
    \]
    as required.
\end{proof}

We are now able to state and prove our main lemmas for analysing escape processes whose high-level state is refilling (Lemma~\ref{lem:refilling-escape-final}), advancing (Lemma~\ref{lem:af-escape-final}) or filling (Lemma~\ref{lem:af-escape-final}). Both proofs will be by coupling to a single constant-escape process followed by straightforward applications of results from Section~\ref{sec:constant-escape-invariants}; the only somewhat non-trivial part is making sure the couplings don't break down. (In other words, recalling the definition of $T_{\nu,\tau}$ from Definition~\ref{def:VEB-couple}, we will need to bound $T_{\nu,\tau}$ below for an appropriate choice of $\nu$ and $\tau$.) Note that we will also Lemma~\ref{lem:af-escape-final} in analysing escape processes whose high-level state is stabilising, hence the more general statement.

\begin{lemma}\label{lem:refilling-escape-final}
Fix a send sequence $\p$ with $p_0 = 1$ and finitely many strongly-exposed bins, and fix $\lambda \in (0,1)$ with $\lambda$ sufficiently small. Let $\Psi = (g,\tau,j,\zvect,\calS,\refilling)$ be a high-level state with $j$ sufficiently large. Fix a tuple $\avect \in (\integers_{\ge 0})^j$. Let $E$ be an escape process with high-level state $\Psi$, send sequence $\p$, birth rate $\lambda$, and initial population satisfying $\romanarrvect^E(\tau) = \avect$. Let $d$ and $\theta$ be real numbers with $d\geq 0$ and $0 \leq  \theta \leq  1/(64j^2)$, and suppose $\avect\in (\integers_\geq 0)^j$ is a $(d,j^{27},\theta,2)$-RS-good beginning. Then with probability at least $1-j^{-2j}$, the following properties hold.
    \begin{enumerate}[(i)]
        \item $\bvect_{[j]}^E(\tau+j^{\phi+46})$ is a $(d+5j^{\phi+49}, j^{10}, \theta+3/j^4, 1)$-RS-good beginning.
        \item For all $k\in [j]$ with $W_k \ge j^\chi$ and all $t$ in the range $\tau \le t \le \tau + j^{\phi+46}$, $p_kb_k^E(t) \le \lambda/80$.
    \end{enumerate}
\end{lemma}
\begin{proof}
Using Definition~\ref{def:constant-escape-coupling}, we couple $E$ with a constant-escape process $E^+$ on $[j]$ with send sequence~$\p$, start time $\tau$, end-of-step arrival rate $\nu_j^\high$
and initial population $\arrvect^E(\tau)$.  
By Lemma~\ref{lem:rs-constant} applied to $E^+$ with $\varpi = j^{27}$, with probability at least $1-j^{-2j}$:
\begin{enumerate}[(R1)]
\item For all $t$ in the range $\tau \le t \le \tau+j^{\phi+46}$, $\bvect_{[j]}^{E^+}(t)$ satisfies (Con\ref{item:Con1}) for $d+5j^{\phi+49}$ (this follows from (i) and (ii) in the lemma statement 
since $d+ 5(t-\tau)j^3 \leq d+5j^{\phi+49}$). 
\item $\bvect_{[j]}^{E^+}(\tau+j^{\phi+46})$ is a $(d+5j^{\phi+49}, j^{10}, \theta+3/j^4,1)$-RS-good beginning (by Lemma~\ref{lem:rs-constant}(ii)).
    \end{enumerate}
    For the rest of the proof, suppose that this event occurs.

We next claim that for all $t$ in the range $\tau \le t \le \tau+j^{\phi+46}$ and all $k \in [j]$ with $W_k \ge j^\chi$, $p_kb_k^{E^+}(t) \le \lambda/80$. By (R1), 
$\excess(\bvect_{[j]}^{E^+}(t),1) \leq d+ 5j^{\phi+49}$ so
for all $k\in [j]$ 
$b_k^{E^+}(t) \le P_kW_k + j^{\phi+49}$.
Since $\phi = 100$ and $\chi=1000$ by Definition~\ref{def:constants}, it follows by Remark~\ref{rem:P-bounded} that $b_k^{E^+}(t) \le 2P_kW_k \le 8\xi W_k < \lambda W_k/80$, so as claimed,
    \begin{equation}\label{eq:refilling-escape-final-1}
        \mbox{ for all $k \in [j]$ with $W_k \ge j^\chi$ and all $t$ in the range $\tau \le t \le \tau+j^{\phi+46}$, }p_kb_k^{E^+}(t) \le \lambda/80.
    \end{equation}

Recall from the coupling of~$E$ and~$E^+$ in
Definition~\ref{def:constant-escape-coupling}
that Inv-$b_i$ holds. Namely, for all $t$ with 
$\tau \leq t \leq T^E_{\nu_{j,\high},\tau}$ and all $i\geq 1$, $\b_i^{E}(t) \subseteq \b_i^{E^+}(t)$.
Given this invariant, it  suffices to prove that $ \tau + j^{\phi+46}\leq
T^E_{\nu_j^\high,\tau}$ --- given this, 
        part (i) of the statement follows from (R2) and part (ii) follows from~\eqref{eq:refilling-escape-final-1}. 
    
    We will bound $T^E_{\nu_j^\high,\tau}$ using Lemma~\ref{lem:escape-nu-high-couple}, which shows that 
    \[
        T_{\nu_j^\high,\tau} \ge \min\{t \ge \tau\colon \mbox{for some $k \in [j-1]$ with $W_k \ge j^\chi$, }b_k^E(t) > \lambda W_k/80\}.
    \]
From the coupling invariant,
    \[
        T_{\nu_j^\high,\tau} \ge \min\{t \ge \tau\colon \mbox{for some $k \in [j-1]$ with $W_k \ge j^\chi$, }b_k^{E^+}(t) > \lambda W_k/80\}.
    \]
 Given~\eqref{eq:refilling-escape-final-1},
  $\tau + j^{\phi+46} < T^E_{\nu_j^\high,\tau}$, as required.
\end{proof}

\begin{lemma}\label{lem:af-escape-final}
Fix a send sequence $\p$ with $p_0 = 1$ and finitely many strongly-exposed bins, and fix $\lambda \in (0,1)$ with $\lambda$ sufficiently small. Let $\Psi = (g,\tau,j,\zvect,\calS,\type)$ be a high-level state with $j$ sufficiently large. 
Suppose that $\calS$ contains a set $S$ such that $|S| \ge \CUpsilon \L(j)/(2\lambda)$ and $|S|\cdot \min\{W_k\colon k \in S\} \ge (j-1)/2$. 
Fix a tuple $\avect \in (\integers_{\ge 0})^j$
such that $\avect$ is a $j^7$-AF-good beginning. Let $E$ be an escape process with high-level state $\Psi$, send sequence $\p$, and birth rate $\lambda$.
Fix $\tau' \ge \tau$. Then with probability at least 
$1 - \exp({-400 j/\lambda^2})$, conditioned on the event $\bvect_{[j]}^E(\tau') = \avect$, the following properties hold.
\begin{enumerate}[(i)]
\item For all $t$ in the range $\tau' \le t \le \tau' + (4/\lambda)^{2j}$, $\bvect_{[j]}^E(t)$ is a $(6j^{25},j^{27},2/(\xi j^3),2)$-RS-good beginning.
\item For all $t$ in the range $\tau' + j^{22} \le t \le \tau' + (4/\lambda)^{2j}$, $\bvect_{[j]}^E(t)$ is a $j^7$-AF-good beginning.
\item For all $t$ in the range $\tau' \le t \le \tau' + (4/\lambda)^{2j}$, $\calN_S^E(t+1) \le \lambda|S|/80$.
\item For all $k\in [j]$ with $W_k \ge j^\chi$ and all $t$ in the range $\tau' \le t \le \tau' + (4/\lambda)^{2j}$, $p_kb_k^E(t) \le \lambda/80$.           
\end{enumerate}
\end{lemma}

\begin{proof}
Using Definition~\ref{def:constant-escape-coupling}, we couple $E$ with a constant-escape process $E^+$ on $[j]$ with send sequence~$\p$, start time $\tau'$, end-of-step arrival rate $\nu_j^\low$
and initial population $\b^E_{[j]}(\tau')$.

By Lemma~\ref{lem:af-good-implies-rs-good},
$\avect$ satisfies (Con\ref{item:Con3}) for $(j^{10},1/(\xi j^3),0)$, satisfies (Con\ref{item:Con4}) for~$1$, and is a $(j^7,j^{10},1/(\xi j^3),1)$-RS-good 
beginning,
so it satisfies (Con\ref{item:Con1}) for~$j^7$.

We will apply three lemmas.
\begin{itemize}
    \item  
We will be able to apply Lemma~\ref{lem:af-constant-stabilise} with $d=j^7$, $\varpi = j^{10}$ and $\theta = 1/(\xi j^3)$ since (Con\ref{item:Con3}) is increasing in~$\nu$
and (Con\ref{item:Con4}) in~$G$. This guarantees that, for any fixed
 $t$ in the range 
$j^{22} \le t -\tau' \le  (4/\lambda)^{2j}$, with probability at least $1-e^{-450j/\lambda^2}$, $\bvect_{[j]}^{E^+}(t)$ is a $j^7$-AF-good beginning.
Applying a union bound over~$t$, using the fact that $\lambda$ is sufficiently small, the probability that this fails for any $t$
is at most $(4/\lambda)^{2j} e^{-450 j/\lambda^2} \leq e^{-449 j/\lambda^2}$.

\item 
We will be able to apply Lemma~\ref{lem:af-constant-strong-survive} with $d=j^7$. This guarantees that, with probability at least $1-e^{-450j/\lambda^2}$, the following holds: For all $t$ satisfying $\tau' \le t \le \tau' + (4/\lambda)^{2j}$, 
$\excess(\bvect_{[j]}^{E^+}(t),3\xi^{-2/3}) \le \lambda j/250$.  
\item We will be able to apply Lemma~\ref{lem:rs-constant} with $(d,\varpi,\theta,G) = (j^7,j^{10},1/(\xi j^3),1)$.  This lemma guarantees that, with probability at least $1-j^{-2j}$, 
\begin{enumerate}[(I)]
    \item  
for all $t$ in the range $\tau' \le t \leq \tau' +  j^{19}$, $\bvect_{[j]}^{E^+}(t)$ is a $(j^7+5(t-\tau')j^3, j^{27}, \theta+3/j^4, 2)$-RS-good beginning. An RS-good beginning is increasing in all of its arguments, so this is a 
$(6j^{25}, j^{27}, 2/(\xi j^3), 2)$-RS-good beginning.
\item for all $t$ in the range
$\tau' + j^{19} \leq t \leq \tau'+j^{22}$,
$\bvect_{[j]}^{E^+}(t)$ is a $(j^7+5(t-\tau')j^3, j^{10}, \theta+3/j^4, 1)$-RS-good beginning. Again, this is a 
$(6j^{25}, j^{27}, 2/(\xi j^3), 2)$-RS-good beginning.

\end{enumerate}
 
\end{itemize}

Combining all three of these with a union bound,
we get that, with probability at least $1 - e^{-400 j/\lambda^2}$, all of the following hold.
    \begin{enumerate}[({A}1)]
        \item For all $t$ with $\tau' \le t \le \tau'+j^{22}$, $\bvect_{[j]}^{E^+}(t)$ is a $(6j^{25}, j^{27}, 2/(\xi j^3), 2)$-RS-good beginning (by the last item).
        \item For all $t$ with $\tau' + j^{22} \le t \le \tau' + (4/\lambda)^{2j}$, $\bvect_{[j]}^{E^+}(t)$ is a $j^7$-AF-good beginning (by the first item). 
        \item For all $t$ with $\tau' \le t \le \tau' + (4/\lambda)^{2j}$, $\excess(\bvect_{[j]}^{E^+}(t), 3\xi^{-2/3}) \le \lambda j/250$ (by  the middle item). 
    \end{enumerate}
    For the rest of the proof, suppose that this event occurs.

    We next claim that for all $t$ with $\tau' \le t \le \tau' + (4/\lambda)^{2j}$, $\calN_S^{E^+}(t+1) \le \lambda/80$. Observe that by Remark~\ref{rem:P-bounded} and (A3),
\begin{align*}
\calN_S^{E^+}(t+1) &= \sum_{k \in S}p_kb_k^{E^+}(t) \le \sum_{k \in S}p_k\Big(3\xi^{-2/3}P_kW_k + \max\{0, b_k^{E^+}(t) - 3\xi^{-2/3}P_kW_k\}\Big)\\
&\le \sum_{k \in S}3\xi^{-2/3}P_k + \max\{p_k\colon k \in S\}\cdot \sum_{k=1}^j\max\{0, b_k^{E^+}(t) - 3\xi^{-2/3}P_kW_k\}\\
&\le 12|S|\xi^{1/3} + \max\{p_k\colon k \in S\}\cdot \excess(\bvect_{[j]}^{E^+}(t), 3\xi^{-2/3})\\
&\le \lambda |S|/1000 + \max\{p_k\colon k \in S\}\cdot \lambda j/250.
\end{align*}

By hypothesis, for all $k \in S$, $|S|W_k \ge (j-1)/2$ and hence $p_k \le 2|S|/(j-1)$; hence
    \begin{equation}\label{eq:af-escape-final-noise}
        \calN_S^{E^+}(t+1) \le \lambda |S|/1000 + \lambda |S|(1+1/(j-1))/125 < \lambda |S|/80
    \end{equation}
    as claimed.

Similarly, we claim that for all $t$ in the range $\tau' \le t \le \tau'+ (4/\lambda)^{2j}$ and all $k \in [j]$ with $W_k \ge j^\chi$, $p_kb_k^{E^+}(t) \le \lambda/80$. 
By (A3), 
$\excess(\bvect_{[j]}^{E^+}(t), 3\xi^{-2/3}) \le \lambda j/250$
so
for any such $k$,
$p_k b_k^{E^+}(t) \leq  3\xi^{-2/3} P_k  + \lambda j p_k/250$.
Using Remark~\ref{rem:P-bounded},
$p_k b_k^{E^+}(t) \leq  12\xi^{1/3}   + \lambda j^{-\chi+1}/250 \leq \lambda/80$.
So we have
    \begin{equation}\label{eq:AF-escape-vhw}
        \mbox{ for all $k \in [j]$ with $W_k \ge j^\chi$ and all $t$ in the range $\tau' \le t \le \tau'+ (4/\lambda)^{2j}$, }p_kb_k^{E^+}(t) \le \lambda/80.
    \end{equation}

Recall from the coupling of~$E$ and~$E^+$ in
Definition~\ref{def:constant-escape-coupling}
that Inv-$b_i$ holds. Namely, for all $t$ with 
$\tau' \leq t \leq T^E_{\nu_{j,\low},\tau'}$ and all $i\geq 1$, $\b_i^{E}(t) \subseteq \b_i^{E^+}(t)$.

We claim that it now suffices to prove that $\tau' + (4/\lambda)^{2j} \leq T^E_{\nu_j^\low, \tau'}  $. 
Given this and the coupling invariant,  part (ii) of the statement follows from (A2), part (iii) of the statement follows from~\eqref{eq:af-escape-final-noise}, 
part (iv) of the statement follows from~\eqref{eq:AF-escape-vhw},
and (A1) implies that for all $t$ with $\tau' \le t \le \tau'+j^{22}$, $\bvect_{[j]}^E(t)$ is a $(6j^{25},j^{27},2/(\xi j^3),2)$-RS-good beginning. Finally, by (ii) of the statement and Lemma~\ref{lem:af-good-implies-rs-good}, for all $t \ge \tau'+j^{22}$, $b_{[j]}^E(t)$ is a $(j^7,j^{10},1/(\xi j^3),1)$-RS-good beginning and hence also a $(6j^{25},j^{27},2/(\xi j^3),2)$-RS-good beginning, so (i) of the statement holds. It remains to prove 
that $\tau' + (4/\lambda)^{2j} \leq T^E_{\nu_j^\low, \tau'}  $.  

    We will bound $T^E_{\nu_j^\low,\tau'}$ using Lemma~\ref{lem:escape-nu-low-couple}; observe that $S$ is the set in $\calS$ required by the hypothesis, so
    \[
        T_{\nu_j^\low,\tau'} \ge \min\Big\{t \ge \tau' \colon \sum_{k \in S}p_kb_k^E(t) > \lambda |S|/80\Big\}.
    \]
    Since $b_k^E(t) \le b_k^{E^+}(t)$ for all $t \le T^E_{\nu_j^\low,\tau'}$, it follows from~\eqref{eq:af-escape-final-noise} that $T^E_{\nu_j^\low,\tau'} \ge \tau'+(4/\lambda)^{2j}$ as required.
\end{proof}

We are now able to start analysing escape processes $E$ whose high-level state is stabilising. Our final lemma statement is Lemma~\ref{lem:escape-stab-final}. This proof will be more complicated than Lemmas~\ref{lem:refilling-escape-final} and~\ref{lem:af-escape-final}, so we break it into two parts. We first prove Lemma~\ref{lem:escape-stab-middle}, which says that with non-trivial failure probability, $E$ will reach an AF-good beginning. Even if this doesn't happen, Lemma~\ref{lem:escape-stab-middle} ensure that with trivial failure probability the pseudorandomness conditions don't degrade too much, so we can try again. The proof of Lemma~\ref{lem:escape-stab-final} will then be by repeated application of Lemma~\ref{lem:escape-stab-middle}.

For the statement of Lemma~\ref{lem:escape-stab-middle}, recall the definition of 
$T_{\nu,\tau'}^E$ from the constant-escape coupling (Definition~\ref{def:constant-escape-coupling}) -- namely, for an escape process~$E$ whose start time is at most~$\tau'$, $T_{\nu,\tau'}^E$ is first time $t\geq \tau'$ such that $\nu^E(t) > \nu$.

\begin{lemma}\label{lem:escape-stab-middle}
Fix a send sequence $\p$ with $p_0 = 1$ and finitely many strongly-exposed bins, and fix $\lambda \in (0,1)$ with $\lambda$ sufficiently small. Let $\Psi = (g,\tau,j,\zvect,\calS,\stabilising)$ be a high-level state with $j$ sufficiently large. Let $E$ be an escape process with high-level state $\Psi$, send sequence $\p$, birth rate $\lambda$. Fix a tuple $\avect \in (\integers_{\ge 0})^j$, let $d$ and $\theta$ be real numbers with $0 \le d \le j^{\phi+51}$ and $0 \leq  \theta \leq  1/(64j^2)$, and suppose that $\avect\in (\integers_{\geq 0})^j$ is a $(d,j^{10},\theta,1)$-RS-good beginning. Fix $\tau' \ge \tau$
and let $T = \min\{\tau' + j^{\phi+70}, T_{\nu_j^\low, \tau' + j^{19}}\}$.
Then, conditioned on $\bvect_{[j]}^E(\tau') = \avect$, the following properties hold.
\begin{enumerate}[(i)]
\item $\tau' + j^{19} \le T \le \tau' + j^{\phi+70}$.
\item With probability at least $1-\exp(-350j/\lambda^2)$, for all $t$ in the range $\tau' \le t \le T$, $\bvect_{[j]}^E(t)$ is a $(d+6j^{22},j^{27},\theta+3/j^4,2)$-RS-good beginning.
\item With probability at least $1-\exp(-350j/\lambda^2)$, $\bvect_{[j]}^E(T)$ is a $(d+6j^{22},j^{10},\theta+3/j^4,1)$-RS-good beginning.
\item With probability at least $1-\exp(-(\log j)^2/2)$, $\bvect_{[j]}^E(T)$ is a $j^7$-AF-good beginning.
\end{enumerate}
\end{lemma}
\begin{proof}
Part (i) of the lemma statement is immediate from the definition of $T$. To prove the rest of the lemma, we apply Definition~\ref{def:constant-escape-coupling} to couple $E$ with two constant-escape processes $E_1^+$ and $E_2^+$ on~$[j]$ with send sequence~$\p$. 
$E_1^+$ has  start time $\tau'$, end-of-step arrival rate $\nu_j^\high$ and initial population $\bvect_{[j]}^E(\tau') = \avect$; $E_2^+$ has   start time $\tau' + j^{19}$, end-of-step arrival rate $\nu_j^\low$ and initial population $\bvect_{[j]}^E(\tau'+ j^{19})$. 
Let $S = [(\log j)^2]$. We define the following properties of $E_1^+$ and $E_2^+$.
\begin{enumerate}[(SM1)]
\item For all $t$ in the range $\tau' \le t \le \tau' + j^{19}$, $\bvect_{[j]}^{E_1^+}(t)$ 
is a $(d+5j^{22},j^{27},\theta+3/j^4,2)$-RS-good beginning.

\item For all $t$ in the range $\tau' + j^{19} \le t \le \tau' + j^{\phi+70}$, $\bvect_{[j]}^{E_1^+}(t)$  is a
$(d+5j^{\phi+73},j^{10}, \theta+3/j^4,1)$-RS-good beginning.

\item $\bvect_{[j]}^E(\tau'+j^{19})$ is a $(d+5j^{22},j^{10},\theta+3/j^4,1)$-RS-good beginning.
\item For all $t$ in the range $\tau' + j^{19} \le t \le j^{\phi+70}$, $\bvect_{[j]}^{E_2^+}(t)$ satisfies (Con\ref{item:Con1}) for $d+5j^{22} + j^{21}$.
\item $\bvect_{[j]}^{E_2^+}(\tau' + j^{\phi+70})$ is a $j^7$-AF-good beginning.
\item For all $t$ in the range $\tau' + j^{19} \le t \le \tau' + j^{\phi+70}$, $\calN_S^{E_1^+}(t+1) \le \lambda |S|/160$.
    \end{enumerate}
We will first lower-bound the probability that (SM1)--(SM6) occur, and then use them to derive items  (i)--(iv) of the lemma statement.

The initial population of~$E_1^{+}$ is given by $\avect$ which is a $(d,j^{10},\theta,1)$-RS-good beginning.   
By Lemma~\ref{lem:rs-constant} applied to $E_1^+$ with $\varpi = j^{10}$ and $G=1$,
    \begin{equation}\label{eq:escape-stab-middle-A}
        \pr(\mbox{(SM1) and (SM2) hold} \mid \bvect_{[j]}^E(\tau') = \avect) \ge 1-j^{-2j}.
    \end{equation}
    
In order to show that (SM1) and (SM2) imply (SM3), we first show that (SM1) and (SM2) imply lower bounds on $T_{\nu_j^\high,\tau'}$. 
\begin{itemize}
\item  Together, (SM1) and (SM2)
imply that, for all $t$ in the range 
$\tau' \leq t \leq \tau' + j^{\phi+70}$
$b_{[j]}^{E_1^+}$ satisfies (Con\ref{item:Con1}) for $d+5 j^{\phi + 73}$ 
so for all $k\in [j]$ with $W_k\geq j^{\chi}$, 
$b_k^{E_1^+}(t) \leq P_k W_k + d + 5j^{\phi+73} \leq 2 P_k W_k < \lambda W_k/80$.

  \item By Lemma~\ref{lem:escape-nu-high-couple},
    \[
    T_{\nu_j^\high,\tau'} \ge \min\{t \ge \tau' \colon \mbox{for some $k \in [j-1]$ with $W_k \ge j^\chi$, $b_k^E(t) > \lambda W_k/80$}\}.
    \]
    
\item From  Inv-$b_i$ of the constant-escape coupling (Definition~\ref{def:constant-escape-coupling}),  for all $t$ 
with $\tau' \le t \le T_{\nu_j^\high,\tau'}$ and all $k \in [j]$, $b_k^E(t) \le b_k^{E_1^+}(t)$. 
    \end{itemize}
    
Putting these three facts together, we obtain the following. 
\begin{equation}\label{eq:escape-stab-middle-1}
\mbox{(SM1) and (SM2) imply $T_{\nu_j^\high,\tau'} \ge \tau' + j^{\phi + 70}$.}
\end{equation}

(SM1) and (SM2) imply that
$\bvect_{[j]}^{E_1^+}(\tau'+j^{19})$ is a $(d+5j^{22},j^{10},\theta+3/j^4,1)$-good beginning. Given \eqref{eq:escape-stab-middle-1} and Inv-$b_i$,
this implies (SM3).

By Lemma~\ref{lem:stab-strong-con1} applied to $E_2^+$ (taking a union bound over $t$ in the range $\tau' + j^{19} \leq t \leq j^{\phi+70}$),
\begin{align}\nonumber
&\pr\big(\mbox{(SM4) occurs}\mid \mbox{(SM3) occurs and }\bvect_{[j]}^E(\tau') = \avect\big)\\\nonumber
&\qquad\qquad\ge\min\Big\{\pr\Big(\mbox{(SM4) occurs}\mid \bvect^E_{[j]}(\tau' + j^{19}) = \avect'\Big) \colon \mbox{$\avect'$ is a $(d+5j^{22},j^{10},\theta+3/j^4,1)$-RS-good beginning}\Big\}\\\label{eq:escape-stab-middle-C}
        &\qquad\qquad\ge 1 - \exp(-400j/\lambda^2).
    \end{align}

Suppose that (SM3) occurs. Then $\bvect_{[j]}^E(\tau'+j^{19})$ satisfies (Con\ref{item:Con1}) for $d+5j^{22}$ and (Con\ref{item:Con4}) for $1$. Moreover, for all $k \in [j]$ with $W_k \ge j^{\phi+54}$,
    \[
b_k^E(\tau'+j^{19}) \le P_kW_k + d + 5j^{22} \le P_kW_k + 2j^{\phi+51} \le (1 + 1/(128j^2))P_kW_k,
    \]
so $\bvect_{[j]}^E(\tau'+j^{19})$ satisfies (Con\ref{item:Con3}) for $(j^{\phi+54},1/(128j^2),\nu_j^\low)$. Thus by Lemma~\ref{lem:af-constant-stabilise} applied to $E_2^+$,
with $\varpi = j^{\phi+54}$
and the $d$ of 
 Lemma~\ref{lem:af-constant-stabilise}
 as our $d+5j^{22} \leq j^{\phi+51} + 5j^{22}$,
 for any $t$ in the 
 range 
 $j^{\phi+54} j^{12} \leq t-(\tau' + j^{19}) \leq (4/\lambda)^{2j}$,
 with probability at least $1-\exp(-450 j/\lambda^2)$, $b_{[j]}^{E_2^+}(t)$ is a $j^7$-AF-good beginning. Taking $t=\tau' + j^{\phi+70}$,

\begin{align}\nonumber
&\pr\big(\mbox{(SM5) occurs}\mid \mbox{(SM3) occurs and }\bvect_{[j]}^E(\tau') = \avect\big)\\\nonumber
&\qquad\qquad\ge\min\Big\{\pr\Big(\mbox{(SM5) occurs}\mid \bvect^E_{[j]}(\tau' + j^{19}) = \avect'\Big) \colon \mbox{$\avect'$ is a $(d+5j^{22},j^{10},\theta+3/j^4,1)$-RS-good beginning}\Big\}\\\label{eq:escape-stab-middle-D}
        &\qquad\qquad\ge 1 - \exp(-450j/\lambda^2).
    \end{align}
Finally, by Lemma~\ref{lem:stab-constant-weak-survive} applied to $E_1^+$ together with a union bound over the values of~$t$,
\begin{equation}\label{eq:escape-stab-middle-E}
\pr(\mbox{(SM6) occurs} \mid \bvect_{[j]}^E(\tau') = \avect) \ge 1 - j^{\phi+70}\exp(-(\log j)^2).
    \end{equation}

Recalling that (SM1) and (SM2) imply (SM3) we next wish to
combine equations~\eqref{eq:escape-stab-middle-A},  \eqref{eq:escape-stab-middle-C}, \eqref{eq:escape-stab-middle-D} and \eqref{eq:escape-stab-middle-E} with  union bounds.

First 
\begin{align*}
\pr(\vee_{k\in [4]} \overline{(\text{SM} k)}  \mid \bvect_{[j]}^E(\tau') = \avect) 
&\leq  
\pr(\vee_{k\in \{1,2\}} \overline{(\text{SM} k)}  \mid \bvect_{[j]}^E(\tau') = \avect) +
\pr(\vee_{k\in \{3,4\}} \overline{(\text{SM} k)}  \mid \bvect_{[j]}^E(\tau') = \avect)
\\
&\leq  
\pr(\vee_{k\in \{1,2\}} \overline{(\text{SM} k)}  \mid \bvect_{[j]}^E(\tau') = \avect) +
\pr( \overline{(\text{SM}3)}  \mid \bvect_{[j]}^E(\tau') = \avect)\\
&\quad\quad\quad + \pr( \overline{(\text{SM}4)}  \mid 
(\text{SM}3)\wedge \bvect_{[j]}^E(\tau') = \avect)
\end{align*}
Equation~\eqref{eq:escape-stab-middle-A} shows that the first term in the final sum is at most $j^{-2j}$ and, since
(SM1) and (SM2) imply (SM3),
that the second term is at most 
$j^{-2j}$.
Equation~\eqref{eq:escape-stab-middle-C} 
shows that the final term is at most~$\exp(-400 j/\lambda^2)$.
Thus, 
\[  \pr(\mbox{(SM1)--(SM4) occur} \mid \bvect_{[j]}^E(\tau') = \avect) \ge 1 - \exp(-350j/\lambda^2).\]
Similar analysis shows
\[  
\pr(\mbox{(SM1)--(SM6) occur} \mid \bvect_{[j]}^E(\tau') = \avect) \ge 1 - \exp(-(\log j)^2/2).\]
  
It therefore suffices to prove that the simultaneous occurrence of (SM1)--(SM4) implies (ii) and (iii) of the lemma statement, and that the simultaneous occurrence of (SM1)--(SM6) implies (iv) of the lemma statement.

Suppose (SM1)--(SM4) occur, so that (by~\eqref{eq:escape-stab-middle-1}) for all $t$ with $\tau' \le t \le \tau' + j^{\phi+70}$, $b_k^E(t) \le b_k^{E_1^+}(t)$. By the definition of $T_{\nu_j^\low,\tau'+j^{19}}$, for all $t$ with $\tau' + j^{19} \le t \le T$, $b_k^E(t) \le b_k^{E_2^+}(t)$. 
We make the followiong conclusions.
\begin{itemize}
\item By (SM1) and (SM4), for all $t$ with $\tau' \le t \le T$, $\bvect_{[j]}^E(t)$ satisfies (Con\ref{item:Con1}) for $d+5j^{22}+j^{21}$.
\item By (SM2) and since $T \ge \tau'+j^{19}$, $\bvect_{[j]}^E(T)$ satisfies (Con\ref{item:Con3}) for $(j^{10},\theta+3/j^4,\nu_j^\high)$ and (Con\ref{item:Con4}) for $1$.
\item By (SM1) and (SM2), for all $t$ with $\tau' \le t \le T$, $\bvect_{[j]}^E(t)$ satisfies (Con\ref{item:Con3}) for $(j^{27},\theta+3/j^4,\nu_j^\high)$ and (Con\ref{item:Con4}) for $2$.
\end{itemize}
It follows that (ii) and (iii) of the lemma statement occur.

Now suppose (SM1)--(SM6) occur. We claim that in order to prove (iv) of the lemma statement, it suffices to prove that $T_{\nu_j^\low,\tau'+j^{19}} \ge \tau' + j^{\phi+70}$. From the definition of~$T$, this would imply $T = \tau' + j^{\phi+70}$. Moreover, it would imply $b_k^E(\tau'+j^{\phi+70}) \le b_k^{E_2^+}(\tau'+j^{\phi+70})$, so (SM5) would imply (iv). 

For the rest of the proof, our goal is to prove that $T_{\nu_j^\low,\tau'+j^{19}} \ge \tau' + j^{\phi+70}$. 
Let $S' = S \setminus \Upsilon_{j, \ge j^2}$.
From the definition of stabilising state (Definition~\ref{def:stabilising-state}), $S' \in \calS$. Since $\Psi$ is a stabilising state
and $j$ is sufficiently large, bin~$j$ is weakly exposed. Thus by Definitions~\ref{def:covered}, \ref{def:constants} and~\ref{def:noise}, $|\Upsilon_{j, \ge j^2}| < \CUpsilon \L(j)/(2\lambda) \le (\log j)^2/8$, so $|S'| \ge |S| - (\log j)^2/8 \ge |S|/2$. 
By (SM6), for all $t$ in the range 
$\tau' + j^{19} \le t \le \tau' + j^{\phi+70}$, $\calN_S^{E_1^+}(t+1) \leq \lambda |S|/160$. By (SM1), (SM2), and~\eqref{eq:escape-stab-middle-1} 
this implies $\calN_S^E(t+1) \leq \lambda |S|/160$.
Then from our bound on $|S'|$,
$\calN_{S'}^E(t+1) \leq \calN_S^E(t+1) \leq \lambda |S|/160 \leq \lambda |S'|/80$. So  
from Definition~\ref{def:escape},
$S' \in \mathbb{S}^E(t)$ and $S' \subseteq \hat{S}^E(t)$. This implies $\hat{S}^E(t) \ge (\log j)^2/4 > 100/\lambda^2$, so $\Xi^E(t) \ge \lambda |S'|/80$ and
    \[
        \nu^E(t) = \exp(-\Xi^E(t)/192) \le \exp(-\lambda(\log j)^2/(192 \times 80 \times 4)) < \nu_j^\low.
    \]
    Thus $T_{\nu_j^\low,\tau'+j^{19}} \ge j^{\phi+70}$, as required.
\end{proof}

\begin{lemma}\label{lem:escape-stab-final}
Fix a send sequence $\p$ with $p_0 = 1$ and finitely many strongly-exposed bins, and fix $\lambda \in (0,1)$ with $\lambda$ sufficiently small. Let $\Psi = (g,\tau,j,\zvect,\calS,\stabilising)$ be a high-level state with $j$ sufficiently large. Let $d,\theta \ge 0$ satisfy $6j^{25} \le d \le j^{\phi+51}/2$  
and $ 2/(\xi j^3) \le \theta < 1/(100j^2)$.  Suppose that $\avect \in (\integers_{\ge 0})^j$ is a   $(d,j^{10},\theta,1)$-RS-good beginning. Let $E$ be an escape process with high-level state $\Psi$, send sequence $\p$, birth rate $\lambda$, and initial population satisfying $\romanarrvect^E(\tau) = \avect$. Then with probability at least $1-\exp(-300j/\lambda^2)$, the following properties hold.
\begin{enumerate}[(i)]
\item For all $t$ with $\tau \le t \le \tau+j^{\phi+72}$, $\bvect_{[j]}^E(t)$ is a $(d+j^{23},j^{27},\theta+1/(j^3 (\log j)^{1/4}),2)$-RS-good beginning.
\item $\bvect_{[j]}^E(\tau + j^{\phi+72})$ is a $j^7$-AF-good beginning.

\end{enumerate}
\end{lemma}
\begin{proof}

For every non-negative integer~$i$, let $d_i = d+6  i j^{22}$ 
and let $\theta_i = \theta + 3 i /j^4$.
Define a sequence of stopping times by $T_0 = \tau$ and, for all $i \ge 1$, 
$
T_i = \min\{T_{i-1} + j^{\phi+70}, T_{\nu_j^\low, T_{i-1} + j^{19}}\}$.
For all $i \ge 1$, let $\calE_i^{\mathrm{succ}}$ be the event that $\bvect_{[j]}^E(T_i)$ is a $j^7$-AF-good beginning, and let $\calE_i^{\mathrm{fail}}$ be the event that either of the following hold.
\begin{itemize}
\item for some $t$ in the range $T_{i-1} \le t \le T_i$, $\bvect_{[j]}^E(t)$ is not a $(d_i,j^{27}, \theta_i, 2)$-RS-good beginning.
\item $\bvect_{[j]}^E(T_i)$ is not a $( d_i, j^{10}, \theta_i, 1)$-good beginning.
\end{itemize}

Consider any $i\geq 0$ and suppose that $\calE_i^{\fail}$ does not occur.
Then $\bvect^E_{[j]}(T_i)$ 
is some tuple $\avect[i]$ 
that is a
$( d_i,  j^{10},  \theta_i, 1)$-good beginning.
If $ d_i \leq j^{\phi+51}$ and $\theta_i  \leq 1/(64 j^2)$ 
then 
we can apply Lemma~\ref{lem:escape-stab-middle} with $d_i$ and $\theta_i$ 
and with $\tau' = T_i$
and with the  $\avect$ of the lemma as $\avect[i]$.
Then $T$ from the statement of Lemma~\ref{lem:escape-stab-middle}
is $T_{i+1}$. 
The lemma shows that 
the probability of $\calE_{i+1}^{\fail}$
is at most $2\exp(-350 j/\lambda^2)$.
Also, 
the probability that $\calE_{i+1}^{\succ}$ does not occur is at most
$\exp(-(\log j)^2/2)$. Putting this together (with the law of total probability over~$\avect[i]$ taking $\avect[0] = \avect$),
for all $i \ge 0$ such that
$ d_i \leq j^{\phi+51}$ and $\theta_i  \leq 1/(64 j^2)$,

    \begin{align*}
        \pr\Big(\overline{\calE_i^{\mathrm{succ}}} \mid \bigcap_{x=1}^{i-1}(\overline{\calE_x^{\mathrm{succ}}} \cap \overline{\calE_x^{\mathrm{fail}}}) \Big) &\le \exp(-(\log j)^2/2),\\
        \pr\Big(\calE_i^{\mathrm{fail}} \mid \bigcap_{x=1}^{i-1}(\overline{\calE_x^{\mathrm{succ}}} \cap \overline{\calE_x^{\mathrm{fail}}}) \Big) &\le 2\exp(-350j/\lambda^2).
    \end{align*}

Let $A = \min\{i \geq 1\colon \calE_i^{\mathrm{succ}}\mbox{ or }\calE_i^{\mathrm{fail}}\mbox{ occurs}\}$, and let $a = \lfloor j/(\log j)^{1/2} \rfloor$. 
Note that $6 a j^{22} < j^{\phi+51}/2$ 
and $3 a/j^4 \leq (1/j^2)(1/64-1/100)$ so
we may take~$i$ up to~$a$.
We then get the following.
    \begin{align*}
        \pr(A > a) &\le \prod_{i=1}^a \pr\Big(\overline{\calE_i^{\mathrm{succ}}} \mid \bigcap_{x=1}^{i-1}(\overline{\calE_x^{\mathrm{succ}}} \cap \overline{\calE_x^{\mathrm{fail}}}) \Big) \le \exp(-a(\log j)^2/2) \le \exp(-350j/\lambda^2),
    \end{align*}
    and
    \begin{align*}
        \pr\Big(A \le a \mbox{ and }\calE_A^\mathrm{fail}\mbox{ occurs}\Big) \le \sum_{i=1}^a \pr\Big(\calE_i^{\mathrm{fail}} \mid \bigcap_{x=1}^{i-1}(\overline{\calE_x^{\mathrm{succ}}} \cap \overline{\calE_x^{\mathrm{fail}}}) \Big) &\le 2j\exp(-350j/\lambda^2).
    \end{align*}
By a union bound, it follows that with probability at least 
$1 - \exp(-349j/\lambda^2)$, $A \le a$, $\calE_A^\mathrm{succ}$ occurs, and $\calE_A^{\mathrm{fail}}$ does not occur. 
When this happens,
$b_{[j]}^E(T_A)$ is some tuple $\avect[A]$ that is a $j^7$-AF-good beginning.
We next wish to apply Lemma~\ref{lem:af-escape-final}, starting from the (fixed) time~$T_A$ (again, applying the law of total probability)
and conditioning on $b_{[j]}^E(T_A) =\avect[A]$ for the fixed tuple $\avect[A]$. In order to this, we first establish the conditions of Lemma~\ref{lem:af-escape-final}.

From the definition of stabilising state (Definition~\ref{def:stabilising-state}),
$\calS$ contains the set
$S = \{k \in [j-1] \setminus [(\log j)^2] \colon W_k < j^2\}\}$.
We wish to establish that $S$ satisfies the two conditions of  
Lemma~\ref{lem:af-escape-final},
namely that $|S| \ge \CUpsilon \L(j)/(2\lambda)$ and $|S|\cdot \min\{W_k\colon k \in S\} \ge (j-1)/2$. 
We will satisfy both of these by showing that $|S|\geq (j-1)/2$.
This follows from the fact that $j$ is weakly exposed (since $\Psi$ is a stabilising state and $j$ is large) so (Prop~\ref{cov-prop-2}) of Definition~\ref{def:covered} is false 
and $|\Upsilon_{j,\geq j^2}| < \CUpsilon \L(j)/(2 \lambda)$.

Applying Lemma~\ref{lem:af-escape-final}
with the $\tau'$ of the lemma as $T_A$ and the $\avect$ of the lemma as $\avect[A]$, we find, with probability at least   $1-\exp(-400j/\lambda^2)$,
the following hold.

\begin{enumerate}[(i)]
\item For all $t$ in the range $T_A \le t \le T_A + (4/\lambda)^{2j}$, $\bvect_{[j]}^E(t)$ is a $(6j^{25},j^{27},2/(\xi j^3),2)$-RS-good beginning.
\item For all $t$ in the range $T_A + j^{22} \le t \le T_A + (4/\lambda)^{2j}$, $\bvect_{[j]}^E(t)$ is a $j^7$-AF-good beginning.         
\end{enumerate}

The definition of~$A$, together with the fact that 
$\calE^{\fail}_A$ does not occur also
establishes
\begin{enumerate}[(iii)]
\item  for all $t$ in the range $\tau \le t \le T_A$, $\bvect_{[j]}^E(t)$ is a $(d_A,j^{27}, \theta_A, 2)$-RS-good beginning.   
\end{enumerate} 

The lemma statement follows from (i), (ii), and (iii),
together with the facts
that 
$d_A \leq d_a \leq d+j^{23}$, $6j^{25} \leq d+j^{23}$, $\theta_ A \leq \theta_a \leq \theta + 1/(j^3(\log j)^{1/4})$, $2/(\xi j^3) \leq \theta + 1/(j^3 (\log j)^{1/4})$
and $T_A + j^{22} \leq 
\tau+ j^{\phi+72} \leq T_A + (4/\lambda)^{2j}
$.
 
\end{proof}

\subsection{Analysing high-level state transitions in the VEB coupling}\label{sec:veb-transition-invariants}

The following notation will be used throughout Section~\ref{sec:veb-transition-invariants}
and Section~\ref{sec:escape-analysis-final}.
It extends the setting of the VEB coupling
(Definition~\ref{def:VEB-couple}) and matches the statement of Lemma~\ref{lem:VEB-escape}.

\begin{definition}\label{def:veb-setting}
The \emph{VEB setting} is described as follows. Fix a send sequence $\p$ with $p_0=1$ and finitely many strongly-exposed bins
such that, for all $i \ge 1$, $W_i \le (4/\lambda)^i$. Fix  $\lambda \in (0,1)$ with $\lambda$ sufficiently small.   Let $X$ be the backoff process with send sequence~$\p$, 
birth rate $2\lambda$,   and cohort set $\calC=\{A,B\}$ (Definition~\ref{def:backoff-cohorts}). 
Let $j_0 \ge 5$ be a sufficiently large positive integer such that all strongly-exposed bins are in $[j_0-1]$. Let $\tau_0 = T^{\p,j_0}$ (Definition~\ref{def:Tprot}). Let $\Psi_0$ be the initialising high-level state with $\tau^{\Psi_0} = \tau_0$ and $j^{\Psi_0} = j_0$. Let $\hls$ be the $\Psi_0$-backoff bounding rule (Definition~\ref{def:backoff-bounding-rule}), and let $\OurHLS$ be the $\Psi_0$-closure of $\hls$ (Definition~\ref{def:hls-closed}). Observe that $\hls$ is a transition rule by Observation~\ref{obs:BB-bounding-is-transrule} and is $\OurHLS$-valid by Lemma~\ref{lem:backoff-bounding-valid}.
    
Let $\Einit$ be the event that, for every step $t' \in [T^{\p,j_0}]$, $n^{X^B} \ge 2$. Let $\eps$ be the (positive) probability that $\Einit$ occurs. Conditioned on $\Einit$, the VEB-coupling (Definition~\ref{def:VEB-couple}) couples $X$ with a volume process $V$ from $\Psi_0$ with transition rule $\hls$, send sequence $\p$ and birth rate $\lambda$ and with a sequence $E_1,E_2,\dots$ of escape processes with send sequence $\p$ and birth rate $\lambda$. The coupling uses $b_i^{E_0}(\tau_0+1)$ as a synonym for $e_i^X(\tau_0+1)$.
\end{definition}

In order to analyse escape processes in the VEB setting, we will need some bounds on possible transition times and possible sequences of high-level states. These have already been proved (with constant failure probability) in Section~\ref{sec:backoff-bounding-analysis}, but we collect them here for easy reference.

\begin{definition}
In the VEB setting, we say that $V$ is \emph{well-behaved} if the following properties hold:
\begin{enumerate}[(WB1)]
\item For all $t \ge \tau_0$, $\type^{\Psi(t)} \ne \Failure$.\label{item:WB1}
\item For all $j \ge j_0$, there is a positive integer $t \le \tau_0 + (4/\lambda)^{2(j-1)}$ such that $\Psi(t)$ is $j$-filling or $j$-advancing.\label{item:WB2}
\item For all $g \ge 0$ with $\type^{\Psi(\tau_g)} \in \{\refilling, \stabilising\}$, there is a positive integer $g'$ with $g+1 \le g' \le g+2j$ such that $\type^{\Psi(\tau_{g'})} = \filling$.\label{item:WB3}
\item For all $g \ge 0$ with $\type^{\Psi(\tau_g)} = \refilling$, $\tau_{g+1} \le \tau_g + (j^{\Psi(\tau_g)})^{\phi+46}$.\label{item:WB4}
\item For all $g \ge 0$ with $\type^{\Psi(\tau_g)} = \stabilising$, $\tau_{g+1} \le \tau_g + (j^{\Psi(\tau_g)})^{\phi+72}$.\label{item:WB5}
\end{enumerate}
\end{definition}

\begin{lemma}\label{lem:well-behaved}
    In the VEB setting, with probability at least $98/100$, $V$ is well-behaved.
\end{lemma}
\begin{proof}
For every $j\geq j_0$, let 
$t_j = \tau_0+ 1 + (8j)^{\phi+1}\sum_{\ell=j_0}^{j-1} W_\ell $ and let
$T_j = \tau_0 + (4/\lambda)^{2(j-1)}$.
In the VEB setting, 
for all $i\geq 1$, $W_i \leq (4/\lambda)^i$.
Also, $j_0$ is sufficiently large and $\lambda$ is sufficiently small.
Thus, for all $j\geq j_0$,
\[
t_j  \le \tau_0+ j^{\phi+2}\sum_{i=1}^{j-1} (4/\lambda)^i \le \tau_0+j^{\phi+2} (4/\lambda)^{j-1} \leq T_j.
\]

Let $\calE_1$ be the event that (WB\ref{item:WB1}) occurs and for all $j \ge j_0$, there is an integer $t$ in the range $\tau_0+1 \leq t \leq  t_j$ such that $\Psi(t)$ is $j$-filling or $j$-advancing.
By Lemma~\ref{lem:VolumeAnalysis}, $\pr(\calE_1) \ge 99/100$.
Recall the following stopping times 
that are defined for every $t\geq \tau_0$
in Definition~\ref{def:stopping-times}.
\begin{align*}
 {\trans({t})} &= \min\{t' > {t} \mid \Psi(t') \neq \Psi({t})\},\\
 {\exit({t})} &= \min\{t' > {t} \mid \type^{\Psi(t')}  \notin \{\refilling,\stabilising\}\}.
\end{align*}
For all $j\geq j_0$ and $t$ in the range $\tau_0+1 \leq t \leq T_j$,
let $\calE_2(j,t)$
be the event that
\begin{itemize}
\item bin $j$ is covered, or
\item $\Psi(t-1)$ is not $j$-filling, or
\item $\Psi(t)$ is not $j$-refilling, or
\item all of the following occur.
\begin{enumerate}[(i)] 
\item $\exit(t) \le t+j^{\phi+74}$,
\item for all $\hat{t}$ 
in the range $t \leq \hat{t} < \exit(t)$  with $\type^{\Psi(t)} = \stabilising$, $\trans(\hat{t}) \le \hat{t} + j^{\phi+72}$,
\item for all $\hat{t}$ 
in the range $t \leq \hat{t} < \exit(t)$ 
with $\type^{\Psi(t)} = \refilling$, $\trans(\hat{t}) \le \hat{t} + j^{\phi+46}$, and
\item $|\{\Psi(t') \colon t \le t' \le \exit(t)-1\}| \le 2j$.
\end{enumerate} 
\end{itemize}
Lemma~\ref{lem:leave-loop-quickly}
shows that, conditioned
on $(\Psi_{\tau_0},\ldots,\Psi_t) = (\psi_{\tau_0},\ldots,\psi_t)$
where $\psi_{t-1}$ is $j$-filling and $\psi_t$ is $j$-refilling,
the probability that (i)--(iv) all occur is at least $1-\exp(-\lambda j(\log j )^2/202)$.
We define $\calE_2(j) = \cap_{t=\tau_0+1}^{T_j} \calE_2(j,t)$ and $\calE_2 = \cap_{j\geq j_0} \calE_2(j)$.
By a union bound,
\[
\pr(\calE_2) \ge 1 - \sum_{j\geq j_0} (T_j - \tau_0) \exp(-\lambda j(\log j)^2/202) \ge 99/100.
\]
    
Suppose that $\calE_1$ and $\calE_2$ occur. We will prove that (WB\ref{item:WB1})--(WB\ref{item:WB5}) occur, completing the proof. (WB\ref{item:WB1}) is part of $\calE_1$. (WB\ref{item:WB2}) follows from $\calE_1$ and $t_j \leq T_j$.  
(WB\ref{item:WB3}) follows from the backoff bounding rule (Definition~\ref{def:backoff-bounding-rule}) using
(WB\ref{item:WB1}) and
(WB\ref{item:WB2}) and 
the fact that $\calE_2$ holds (specifically, that item (iv) holds
for all $j\geq j_0$ and $t\in \{\tau_0+1,\ldots,T_j\}$).
(WB\ref{item:WB4}) follows from item (iii) and (WB\ref{item:WB5}) follows from item~(ii). 
\end{proof}

We are now in a position to translate the results of Section~\ref{sec:escape-invariants} to the VEB setting. We will work on the assumption that $V$ is well-behaved (since we will be able to union bound away the failure probability). We have one lemma for each non-trivial high-level state type: Lemma~\ref{lem:veb-refilling-transition} for refilling states; Lemma~\ref{lem:advancing-escape-final} for advancing states; Lemma~\ref{lem:filling-escape-final} for filling states; and Lemma~\ref{lem:veb-escape-stab} for stabilising states. None of these lemmas are particularly difficult, but stating them explicitly in this section will be very convenient for the proof of Lemma~\ref{lem:VEB-escape}.

\begin{lemma}\label{lem:veb-refilling-transition}
Consider the VEB setting, and let $\Psi = (g,\tau_g,j,\zvect,\calS,\refilling) \in \OurHLS$ be a high-level state.  Let $d$ and $\theta$ be real numbers with $d \ge 0$ and $0 \leq  \theta \leq 1/(64j^2)$, and suppose that $\avect\in (\integers_{\ge 0})^j$ is a $(d,j^{27},\theta,2)$-RS-good beginning. Conditioned on the events $\Psi(\tau_g) = \Psi$ and $\bvect_{[j]}^{E_g}(\tau_g) = \avect$, with probability at least $1-j^{-2j}$,  
the following events occur.
\begin{enumerate}[(i)]
\item $\bvect_{[j]}^{E_g}(\tau_g+j^{\phi+46})$ is a $(d+5j^{\phi+49}, j^{10}, \theta+3/j^4, 1)$-RS-good beginning.
\item For all $k\in [j]$ with $W_k \ge j^\chi$ and all $t$ in the range $\tau_g \le t \le \tau_g + j^{\phi+46}$, $p_kb_k^{E_g}(t) \le \lambda/80$.
\end{enumerate}
Moreover, if these events occur and $V$ is well-behaved, then the following hold.
\begin{enumerate}[(I)]
\item $\Psi(\tau_{g+1})$ is a refilling or stabilising state.
\item $\bvect_{[j]}^{E_{g+1}}(\tau_{g+1})$ is a $(d+5j^{\phi+49},j^{10},\theta + 3/j^4, 1)$-RS-good beginning.
\item For all $t$ with $\tau_g \le t < \tau_{g+1}$, 
and all $k\in [j-1]$ with $W_k\geq j^{\chi}$, 
$\{k\} \in \calS$  and $\calN^{E_{g}}_{\{k\}}(t+1) \le \lambda /80$.
Moreover, there is at least one such $k\in [j-1]$.
\end{enumerate}
\end{lemma}
\begin{proof}
The probability bound for items (i) and (ii) follows immediately from Lemma~\ref{lem:refilling-escape-final}. Suppose that
these events occur and that $V$ is well-behaved.

By the backoff-bounding rule (Definition~\ref{def:backoff-bounding-rule}), the only possible high-level state transitions from $\Psi$ are via (R\ref{item:R1}), (R\ref{item:R2})(i), (R\ref{item:R5})(i) or (R\ref{item:R5})(ii). A transition via (R\ref{item:R5})(ii) cannot occur at time $\tau_{g+1}$ as $\Psi(\tau_{g+1}) \ne \Psi$, and transitions via (R\ref{item:R1}) and (R\ref{item:R2})(i) cannot occur by (WB\ref{item:WB1}). Thus the high-level state transition from $\Psi$ to $\Psi(\tau_{g+1})$ is via (R\ref{item:R5})(i). Part (I) of the lemma follows immediately.

By the statement of (R\ref{item:R5})(i), $\tau_{g+1} \ge \tau_g + j^{\phi+46}$, so by (WB\ref{item:WB4}) it follows that $\tau_{g+1} = \tau_g + j^{\phi+46}$. Thus since (i) occurs, 
and (by the VEB coupling, Definition~\ref{def:VEB-couple}), $\bbbar_{[j]}^{E_{g+1}}(\tau_{g+1}) = 
\bbbar_{[j]}^{E_g}(\tau_{g+1})$,
item~(II) of the lemma statement follows. 

Since  (ii) occurs, for all $t$ with $\tau_g \le t \leq \tau_{g+1}$ and all $k \in [j]$ with $W_k \ge j^\chi$, $p_kb_k^{E_g}(t) \le \lambda/80$. 
By the definition of a refilling state (Definition~\ref{def:refilling-state}), 
for every $k\in [j-1]$ with $W_k\geq j^{\chi}$, 
$\{k\} \in \calS$. 
Since $\Psi$ is a refilling state, bin $j$ is exposed, and since $j \ge j_0$ bin $j$ is weakly exposed; thus by (Prop~\ref{cov-prop-3}) in Definition~\ref{def:covered}, there is at least one $k \in [j-1]$ with $W_k \ge j^\chi$. Item (III) follows.
\end{proof}

\begin{lemma}\label{lem:advancing-escape-final}
Consider the VEB setting, and let $\Psi = (g,\tau_g,j,\zvect,\calS,\advancing) \in \OurHLS$ be a high-level state. 
Let $S$ be the (sole) set in $\calS$ (Definition~\ref{def:advancing-state}).
Suppose that $\avect \in (\integers_{\ge 0})^j$  is a $j^7$-AF-good beginning. Conditioned on the events $\Psi(\tau_g) = \Psi$ and $\bvect_{[j]}^{E_g}(\tau_g) = \avect$, with probability at least $1-\exp(-400j/\lambda^2)$, 
the following events occur. 
\begin{enumerate}[(i)]
\item For all $t$ in the range $\tau_g \le t \le \tau_g + (4/\lambda)^{2j}$, $\bvect_{[j]}^{E_g}(t)$ is a $(6j^{25},j^{27},2/(\xi j^3),2)$-RS-good beginning.
\item For all $t$ in the range $\tau_g + j^{22} \le t \le \tau_g + (4/\lambda)^{2j}$, $\bvect_{[j]}^{E_g}(t)$ is a $j^7$-AF-good beginning.
\item For all $t$ in the range $\tau_g \le t \le \tau_g + (4/\lambda)^{2j}$, $\calN_S^{E_g}(t+1) \le \lambda|S|/80$.
\end{enumerate} 
Moreover, if these events do occur and $V$ is well-behaved, then the following hold:
\begin{enumerate}[(I)]
\item $\Psi(\tau_{g+1})$ is a $(j+1)$-advancing or $(j+1)$-filling state.
\item $\bvect_{[j]}^{E_{g+1}}(\tau_{g+1})$ is a $(j+1)^7$-AF-good beginning.
\item For all $t$ with $\tau_g \le t < \tau_{g+1}$,   $\calN^E_S(t+1) \le \lambda |S|/80$. Moreover, $S$ is non-empty.
\end{enumerate}
\end{lemma}
\begin{proof}
By Definition~\ref{def:advancing-state}, 
$S =\Upsilon_{j,\ge\Wtilde[j]}$ if bin $j$ is many-covered
and $S =  \Upsilon_{j,\ge j^2}$ if bin $j$ is heavy-covered.  
We will show that, in either case, 
$|S| \ge \CUpsilon \L(j)/(2\lambda)$ and $|S|\cdot \min\{W_k\colon k \in S\} \ge (j-1)/2$. 
\begin{itemize}
\item  
If bin $j$ is many-covered, then by Definition~\ref{def:Wtilde}, using the fact that $j$ is sufficiently large,
$ |\Upsilon_{j,\ge \Wtilde[j]}| \geq \CUpsilon\L(j)/\lambda$.
By Observation  Observation~\ref{obs:Wtilde},   $|\Upsilon_{j,\ge \Wtilde[j]}| \cdot \Wtilde[j] \ge j-1$.
\item 
  If   bin $j$ is heavy-covered, then by (Prop~\ref{cov-prop-2}) of Definition~\ref{def:covered}, $|\Upsilon_{j,\ge j^2}| \ge \CUpsilon\L(j)/(2\lambda)$. Thus,  $|\Upsilon_{j,\ge j^2}| \cdot j^2 > j-1$.

\end{itemize}
The probability bound for items (i)--(iii) follows from  Lemma~\ref{lem:af-escape-final}. Suppose that these events occur and that $V$ is well-behaved.

By the backoff bounding rule (Definition~\ref{def:backoff-bounding-rule}), the only possible high-level state transitions from $\Psi$ are via (R\ref{item:R1}), (R\ref{item:R2})(i), (R\ref{item:R4})(i) or (R\ref{item:R4})(ii). A transition via (R\ref{item:R4})(ii) cannot occur at time $\tau_{g+1}$ as $\Psi(\tau_{g+1}) \ne \Psi$, and transitions via (R\ref{item:R1}) or (R\ref{item:R2})(i) cannot occur by (WB\ref{item:WB1}). Thus the high-level state transition from $\Psi$ to $\Psi(\tau_{g+1})$ is via (R\ref{item:R4})(i). By (WB\ref{item:WB1}) it is to a $(j+1)$-advancing or $(j+1)$-filling state. Thus (I) of the lemma statement holds. Moreover, such a transition can only occur if $\tau_{g+1} \ge \tau_g+j^{24}$, and $\tau_{g+1} \le \tau_g + (4/\lambda)^{2j}$ by (WB\ref{item:WB2}),  so 
(ii) implies that $\bvect{[j]}^{E_g}(\tau_{g+1})$ is a $j^7$-AF-good beginning and hence a $(j+1)^7$-AF-good-begining.
By the VEB coupling,   Definition~\ref{def:VEB-couple}), $\bbbar_{[j]}^{E_{g+1}}(\tau_{g+1}) = 
\bbbar_{[j]}^{E_g}(\tau_{g+1})$, so item (II) of the lemma statement holds.  Finally, (III) of the lemma statement is immediate from (iii) and from the fact that we have a non-zero lower bound for~$|S|$.
\end{proof}

\begin{lemma}\label{lem:filling-escape-final}
Consider the VEB setting, and let $\Psi = (g,\tau_g,j,\zvect,\calS,\filling) \in \OurHLS$ be a high-level state. Let
$S = \Upsilon_{j, \ge\Wtilde[j]} \setminus \Upsilon_{j, \ge j^2}$.
Suppose that $\avect \in (\integers_{\ge 0})^j$ is a $j^7$-AF-good beginning. Conditioned on the events $\Psi(\tau_g) = \Psi$ and $\bvect_{[j]}^{E_g}(\tau) = \avect$, with probability at least $1-\exp(-400j/\lambda^2)$, the 
following events occur.
\begin{enumerate}[(i)]
\item For all $t$ in the range $\tau_g \le t \le \tau_g + (4/\lambda)^{2j}$, $\bvect_{[j]}^{E_g}(t)$ is a $(6j^{25},j^{27},2/(\xi j^3),2)$-RS-good beginning.
\item For all $t$ in the range $\tau_g + j^{22} \le t \le \tau_g + (4/\lambda)^{2j}$, $\bvect_{[j]}^{E_g}(t)$ is a $j^7$-AF-good beginning.
\item For all $t$ in the range $\tau_g \le t \le \tau_g + (4/\lambda)^{2j}$, $\calN_S^{E_g}(t+1) \le \lambda|S|/80$.

    \item For all $k\in [j]$ with $W_k \ge j^\chi$ and all $t$ in the range $\tau_g \le t \le \tau_g + (4/\lambda)^{2j}$, $p_kb_k^{E_g}(t) \le \lambda/80$.

\end{enumerate}   Moreover, if these events do occur and $V$ is well-behaved, then the following hold:
\begin{enumerate}[(I)]
\item $\Psi(\tau_{g+1})$ is a $j$-refilling, $(j+1)$-advancing or $(j+1)$-filling state.
\item If $\Psi(\tau_{g+1})$ is a $j$-refilling state, $\bvect_{[j]}^{E_{g+1}}(\tau_{g+1})$ is a $(6j^{25}, j^{27}, 2/(\xi j^3), 2)$-RS-good beginning.
\item If $\Psi(\tau_{g+1})$ is a $(j+1)$-advancing or $(j+1)$-filling state, $\bvect_{[j]}^{E_{g+1}}(\tau_{g+1})$ is a $(j+1)^7$-AF-good beginning.

\item For all $t$ with $\tau_g \le t < \tau_{g+1}$,   $\calN^{E_g}_S(t+1) \le \lambda |S|/80$. Moreover, $S$ is non-empty.

\item For all $t$ with $\tau_g \le t < \tau_{g+1}$, 
and all $k\in [j-1]$ with $W_k\geq j^{\chi}$, 
$\{k\} \in \calS$  and $\calN^{E_{g}}_{\{k\}}(t+1) \le \lambda /80$.
Moreover, there is at least one such $k\in [j-1]$.

\end{enumerate}
\end{lemma}
\begin{proof}
By Definition~\ref{def:filling-state}, $S\in \calS$. 
by Definition~\ref{def:Wtilde}, using the fact that $j$ is sufficiently large,
$  |\Upsilon_{j,\geq \Wtilde[j]}| \geq \CUpsilon\L(j)/\lambda$.
Since $j$ is not covered, (Prop~\ref{cov-prop-2}) gives that 
$|\Upsilon_{j,\geq j^2}| < \CUpsilon \L(j)/(2 \lambda)$.
Hence, $|S| \geq \CUpsilon \L(j)/(2 \lambda)$.
By Observation~\ref{obs:Wtilde},   $|\Upsilon_{j,\ge \Wtilde[j]}| \cdot \Wtilde[j] \ge j-1$
so $|S| \cdot \Wtilde[j] \geq (j-1)/2$.
The probability bound for items (i)--(iv) follows from Lemma~\ref{lem:af-escape-final}. Suppose that these events occur and that $V$ is well-behaved.

By the backoff bounding rule (Definition~\ref{def:backoff-bounding-rule}), the only possible high-level state transitions from $\Psi$ are via (R\ref{item:R1}), (R\ref{item:R2})(ii), (R\ref{item:R4})(i) or (R\ref{item:R4})(ii). A transition via (R\ref{item:R4})(ii) cannot occur at time $\tau_{g+1}$ as $\Psi(\tau_{g+1}) \ne \Psi$, and a transition via (R\ref{item:R1}) cannot occur by (WB\ref{item:WB1}). 
A transition via (R\ref{item:R2})(ii) is to a $j$-refilling state
and a  transition via  (R\ref{item:R4})(i)  is to a $(j+1)$-advancing or $(j+1)$-filling state by (WB\ref{item:WB1}). Thus (I) of the lemma statement holds. 
By (WB\ref{item:WB2}), $\tau_{g+1} \le \tau + (4/\lambda)^{2j}$, so part (II) of the lemma statement follows from~(i). 
If $\Psi(\tau_{g+1})$ is a $(j+1)$-advancing or $(j+1)$-filling state, then the transition  is via (R\ref{item:R4})(i) so $\tau_{g+1} \geq \tau_g + j^{24}$. Thus part (III) of the lemma statement follows from~(ii)
and from the fact that a $j^7$-AF-good beginning is a $(j+1)^7$-AF-good beginning. Part (IV) of the lemma statement follows from~(iii) and from our lower bound on~$|S|$. Finally, part~(V) follows from the definition of filling state (Definition~\ref{def:filling-state}) which guarantees that $\{k\} \in \calS$, from~(iv), which gives the upper bound on the noise,  and from the fact that $j$ is weakly exposed (since the high-level state is filling), so $j$~satisfies (Prop~\ref{cov-prop-3}) in Definition~\ref{def:covered}, so there are at least~$\CUpsilon \LL(j)/\lambda$ such~$k$.
\end{proof}

\begin{lemma}\label{lem:veb-escape-stab}
Consider the VEB setting, and let $\Psi = (g,\tau_g,j,\zvect,\calS,\stabilising) \in \OurHLS$ be a high-level state.  let $d$ and $\theta$ be real numbers satisfying $6j^{25} \le d \le j^{\phi+51}/2$ and $2/(\xi j^3) \le \theta < 1/(100j^2)$, and suppose that $\avect \in (\integers_{\ge 0})^j$ is a $(d,j^{10},\theta,1)$-RS-good beginning. Conditioned on the events $\Psi(\tau_g) = \Psi$ and $\bvect_{[j]}^{E_g}(\tau_g) = \avect$, with probability at least $1-\exp(-300j/\lambda^2)$, the 
following events occur.
\begin{enumerate}[(i)]
\item For all $t$ with $\tau_g \le t \le \tau_g+j^{\phi+72}$, $\bvect_{[j]}^{E_g}(t)$ is a $(d+j^{23},j^{27},\theta+1/(j^3(\log j)^{1/4}),2)$-RS-good beginning.
\item $\bvect_{[j]}^{E_g}(\tau_g + j^{\phi+72})$ is a $j^7$-AF-good beginning.

\end{enumerate}
Moreover, if these events occur and $V$ is well-behaved, then the following hold:
\begin{enumerate}[(I)]
\item $\Psi(\tau_{g+1})$ is a $j$-refilling or $j$-filling state.
\item If $\Psi(\tau_{g+1})$ is a $j$-refilling state, $\bvect_{[j]}^{E_{g+1}}(\tau_{g+1})$ is a $(d+j^{23}, j^{27}, \theta+1/(j^3(\log j)^{1/4}), 2)$-RS-good beginning.
\item If $\Psi(\tau_{g+1})$ is a $j$-filling state, $\bvect_{[j]}^{E_{g+1}}(\tau_{g+1})$ is a $j^7$-AF-good beginning.
\item For all $t$ with $\tau_g \le t < \tau_{g+1}$, 
and all $k\in [j-1]$ with $W_k\geq j^{\chi}$, 
$\{k\} \in \calS$  and $\calN^{E_{g}}_{\{k\}}(t+1) \le \lambda /80$.
Moreover, there is at least one such $k\in [j-1]$.  
    \end{enumerate}
\end{lemma}
\begin{proof}
The probability bound follows immediately from Lemma~\ref{lem:escape-stab-final}. Suppose the these events  occur and $V$ is well-behaved.

By the backoff bounding rule (Definition~\ref{def:backoff-bounding-rule}), the only possible high-level state transitions from $\Psi$ are via (R\ref{item:R1}), (R\ref{item:R2})(ii), (R\ref{item:R6})(i) or (R\ref{item:R6})(ii). A transition via (R\ref{item:R6})(ii) cannot occur at time $\tau_{g+1}$ as $\Psi(\tau_{g+1}) \ne \Psi$, and a transition via (R\ref{item:R1}) cannot occur by (WB\ref{item:WB1}). In either of the remaining cases, $\Psi(\tau_{g+1})$ is a $j$-refilling or $j$-filling state, so (I) of the lemma statement holds. 
By (WB\ref{item:WB5}), $\tau_{g+1} \le \tau_g + j^{\phi+72}$, so part (II) of the lemma statement follows from~(i) and
from the fact (from the VEB coupling, Definition~\ref{def:VEB-couple}), 
that $\bbbar_{[j]}^{E_{g+1}}(\tau_{g+1}) = 
\bbbar_{[j]}^{E_g}(\tau_{g+1})$.

If $\Psi(\tau_{g+1})$ is a $j$-filling state, then the transition is via (R\ref{item:R6})(i) so $\tau_{g+1} \ge \tau_g + j^{\phi+72}$ which means that
$\tau_{g+1} = \tau_g + j^{\phi+72}$. Thus part (III) of the lemma statement  follos from~(ii) and (WB\ref{item:WB5}).

Finally, we establish Part (IV) of the lemma. From (i), for all $t$ 
in the range $\tau_g \leq t \leq \tau_{g+1}$, $\excess(\bvect_{[j]}^{E_g}(t),1) \leq d + j^{23}$,
so for all $k\in [j-1]$ with $W_k \geq j^{\chi}$,
$b_k^{E_g}(t) \leq P_k W_k + d + j^{23} \leq 2 P_k W_k$. Therefore, 
$p_k b_k^{E_g}(t) \leq 2 P_k \leq \lambda /80$. The  definition of 
stabilising state (Definition~\ref{def:stabilising-state})  guarantees that $\{k\} \in \calS$. Since $j$ is weakly exposed
(which follows from the fact that the high-level state is stabilising),   $j$ satisfies (Prop~\ref{cov-prop-3}) in Definition~\ref{def:covered}, so there
are at least~$\CUpsilon \LL(j)/\lambda$ such~$k$.
\end{proof}

\subsection{Analysing the VEB coupling}\label{sec:escape-analysis-final}

We are now at last in a position to prove Lemma~\ref{lem:VEB-escape}, the overarching goal of Section~\ref{sec:escape-analysis-final}.

\begin{lemma}
\label{lem:VEB-escape}
There exists a real number $\lambda_0 \in (0,1/2)$ such that the following holds. Consider any $\lambda \in (0,\lambda_0)$ and a send sequence $\p$ such $p_0=1$ and, for all $i \ge 1$, $W_i \le (4/\lambda)^i$. Suppose that $\p$ has finitely many bins which are strongly exposed. 
Let $X$ be the backoff process with
birth rate~$2 \lambda$, send sequence~$\p$, and cohort set~$\calC=\{A,B\}$ (Definition~\ref{def:backoff-cohorts}).
Let $j_0\geq 5$ be a sufficiently large positive integer such that all strongly-exposed bins are in $[j_0-1]$.
Let $\Einit$ be the event that, for every step $t' \in [T^{\p,j_0}]$, $n^{X^B} \ge 2$. 
We will condition on $\Einit$.
Let
$\tau_0 =  T^{\p,j_0}$. 
Let
$\Psi_0$ be the  initialising high-level state   
with $\tau^{\Psi_0} = \tau_0$ and $j^{\Psi_0} = j_0$. 
Let $\hls$ be   the $\Psi_0$-backoff bounding rule from Definition~\ref{def:backoff-bounding-rule} and let $\OurHLS$  (Definition~\ref{def:backoff-bounding-state-space}) be the  $\Psi_0$-closure (Definition~\ref{def:hls-closed}) of~$\hls$.  
The $\Psi_0$-backoff bounding rule~$\hls$ is a transition rule (Definition~\ref{def:transition-rule}) by Observation~\ref{obs:BB-bounding-is-transrule} and 
it is $\OurHLS$-valid (Definition~\ref{def:valid}) by
Lemma~\ref{lem:backoff-bounding-valid}.

Conditioned on $\Einit$,
the VEB-coupling (Definition~\eqref{def:VEB-couple}   couples $X$ with a volume process $V$ from $\Psi_0$ with  transition rule $\hls$, send sequence~$\p$ and birth rate~$\lambda$ and with a sequence $E_1, E_2,\ldots $ of escape processes with send sequence~$\p$ and birth rate~$\lambda$. 
The coupling uses $b_i^{E_0}(\tau_0+1)$ as a synonym for $\e_i^X(\tau_0+1)$.

With probability at least~$2/3$,  
both of the following hold.
 
\begin{enumerate}[(i)]
\item for all $t> \tau_0$ such that 
$\type^{\Psi(t)}  = \advancing$
and the $S\in \calS^{\Psi(t)}$,
$\sum_{k \in  S } p_kb_k^{E_{g^{\Psi(t-1)}}}(t) \leq \lambda | S  |/80$. Moreover, $S$ is non-empty.
\item for all $t> \tau_0$ such that 
$\type^{\Psi(t)} \notin \{\Failure,\advancing\}$ and all $k \in [j^{\Psi(t-1)}-1]$ such that $W_k \ge (j^{\Psi(t-1)})^\chi$, $p_kb_k^{E_{g^{\Psi(t-1)}}}(t) \le \lambda/80$. Moreover, there is at least one such $k\in [j-1]$.
\end{enumerate}
\end{lemma}

\begin{proof} 
Observe that the hypothesis of Lemma~\ref{lem:VEB-escape} is precisely the VEB setting of Definition~\ref{def:veb-setting}.  Each non-initialising non-failure high-level state type has a corresponding lemma in Section~\ref{sec:veb-transition-invariants}, and we now define corresponding events. For brevity, for all $g \ge 1$, let $\Psi[g] = \Psi(\tau_g)$, $j[g] = j^{\Psi[g]}$, and $\avect[g] = \bvect_{[j]}^{E_g}(\tau_g)$. Then for all $g \ge 1$:
\begin{itemize}
    \item If $\type^{\Psi[g]} = \refilling$, $E_g$ \emph{starts well} if for some real $d$ and $\theta$ with $d \ge 0$ and $0 \le \theta \le 1/(64j[g]^2)$, $\avect[g]$ is a $(d,j[g]^{27},\theta,2)$-RS-good beginning. $E_g$ \emph{derails} if $E_g$ starts well but at least one event of Lemma~\ref{lem:veb-refilling-transition}(i)--(ii) does not occur.
    \item If $\type^{\Psi[g]} = \advancing$, $E_g$ \emph{starts well} if $\avect[g]$ is a $j[g]^7$-AF-good beginning. $E_g$ \emph{derails} if $E_g$ starts well but at least one event of Lemma~\ref{lem:advancing-escape-final}(i)--(iii) does not occur.
    \item If $\type^{\Psi[g]} = \filling$, $E_g$ \emph{starts well} if $\avect[g]$ is a $j[g]^7$-AF-good beginning. $E_g$ \emph{derails} if $E_g$ starts well but at least one event of Lemma~\ref{lem:filling-escape-final}(i)--(iv) does not occur.
\item If $\type^{\Psi[g]} = \stabilising$, $E_g$ \emph{starts well} if 
for some real numbers
$d$ and $\theta$ satisfying $6j[g]^{25} 
\le d \le j[g]^{\phi+51}/2$ and $2/(\xi j[g]^3) \le \theta < 1/(100j[g]^2)$,
$\avect[g]$ is a $(d,j[g]^{10},\theta,1)$-RS-good beginning. $E_g$ \emph{derails} if $E_g$ starts well but at least one event of Lemma~\ref{lem:veb-escape-stab}(i)--(ii) does not occur.
    \item If $\type^{\Psi[g]} \in \{\Failure,\initialising\}$, then $E_g$ always \emph{starts well} and never \emph{derails}.
\end{itemize}
Observe that for all $g \ge 1$, conditioned on $\Psi[g]$ and $\avect[g]$, the event that $E_g$ starts well is fully determined and the event that $E_g$ derails is determined entirely by $E_g$. We will first prove that conditioned on $\Einit$, with probability at least $2/3$, $V$ is well-behaved and no escape process $E_g$ derails. We will then prove that deterministically, when these events occur, every escape process $E_g$ starts well. Finally, we will use this fact to prove the lemma statement.

\medskip\noindent\textbf{Bounding failure probabilities:}
Each of Lemmas~\ref{lem:veb-refilling-transition}, \ref{lem:advancing-escape-final}, \ref{lem:filling-escape-final}
and~\ref{lem:veb-escape-stab} has a failure probability of at most $\exp(-300j/\lambda^2)$, and their hypothesis is precisely that $E_g$ starts well. Thus by 
Lemmas~\ref{lem:veb-refilling-transition}, \ref{lem:advancing-escape-final}, \ref{lem:filling-escape-final}
and~\ref{lem:veb-escape-stab}, for all $g \ge 1$,
\begin{equation}\label{eq:foo}
    \pr\big(E_g\mbox{ derails} \mid \avect[g]\big) \le \exp(-300j[g]/\lambda^2).
\end{equation}
Let $\calF$ be the event that $V$ is well-behaved; then $\pr(\calF) \ge 98/100$ by Lemma~\ref{lem:well-behaved}. Then by the law of total probability then a union bound, 
\begin{align}\nonumber
    \pr\Big(\bigcup_{g=1}^\infty \{E_g\mbox{ derails}\}\mid\Einit \Big) &\le \pr\Big(\calF \cap \bigcup_{g=1}^\infty \{E_g\mbox{ derails}\}\mid \Einit\Big) + \pr(\overline{\calF})\\\label{eq:bar}
    &\le \sum_{g=1}^\infty \pr(\calF\mbox{ occurs and $E_g$ derails}\mid \Einit) + 2/100.
\end{align}

We next bound $j[g]$ below, using the assumption that $\calF$ occurs. For all $g \ge 1$, let $\calF_g$ be the event that $j[g] \ge (\log g)/(2\log (4/\lambda))$; we will prove $\calF \subseteq \calF_g$. Suppose that $\calF$ occurs. Observe from the backoff-bounding rule (Definition~\ref{def:backoff-bounding-rule}) that $j^{\Psi[g]}$ is non-decreasing in $g$. Thus by (WB2), for all $j \ge j_0$ and all $t \ge \tau_0 + (4/\lambda)^{2(j-1)}$, $j^{\Psi(t)} \ge j$. Since each high-level state transition takes at least one time step, it follows that, for all $g \ge (4/\lambda)^{2(j-1)}$, $  j[g]   \ge j$, and hence that $j[g] \ge 1 + (\log_{4/\lambda} g)/2 > (\log g)/(2\log (4/\lambda))$. Thus $\calF_g$ occurs as claimed. 
 
For all $g \ge 1$, since $\calF \subseteq \calF_g$,
\[
\pr(\mbox{$\calF$ occurs and $E_g$ derails}\mid\Einit ) \le \pr(\mbox{$\calF_g$ occurs and $E_g$ derails}\mid\Einit ).
\]
Integrating over possible values $\avect$ of $\avect[g]$ which determine that $j[g] \ge (\log g)/(2\log (4/\lambda))$, it follows that
\begin{align*}
&\pr(\mbox{$\calF$ occurs and $E_g$ derails}\mid\Einit )
\le \min_{\avect}\Big(\pr\big(E_g\mbox{ derails} \mid \avect[g] = \avect, \Einit\big)\Big).
\end{align*}

By~\eqref{eq:foo} and since $\lambda$ is small, it follows that 
\[
\pr(\mbox{$\calF$ occurs and $E_g$ derails}\mid\Einit ) \le \exp\Big({-}\frac{300\log g}{\lambda^2\cdot 2\log( 4/\lambda)}\Big) \le g^{-300/\lambda}.
\]
Moreover, deterministically $j[g] \ge j_0$ for all $g \ge 1$, and so it is also true that for all $g \ge 1$,
\[
\pr(\mbox{$\calF$ occurs and $E_g$ derails}\mid\Einit) \le \exp\Big({-}\frac{300j_0}{\lambda^2}\Big) < 1/100.
\]
Thus by~\eqref{eq:bar},
\[
\pr\Big(\bigcup_{g=1}^\infty \{E_g\mbox{ derails}\}\mid\Einit\Big) \le 1/100 + \sum_{g=1}^\infty g^{-300/\lambda} +  2 /100 \le  4/100.
\]
It now follows immediately by a union bound that
\begin{equation}\label{eq:veb-escape-prob-goal}
    \pr\big(\mbox{$V$ is well-behaved and no $E_g$ derails} \mid \Einit\big) \ge 2/3.
\end{equation}

\medskip\noindent\textbf{Escape processes start well:}
We next prove that if $V$ is well-behaved and no escape process $E_g$ derails, then every escape process $E_g$ starts well. In fact, we will prove the following stronger statement. For all $g \ge 1$, let $\gamma(g) = \max\{g' \le g\colon \type^{\Psi[g']} \in \{\filling, \advancing\}\}$. Then we will prove the following.
\begin{enumerate}[(VEB1)]
    \item If $\type^{\Psi[g]} = \refilling$, then 
    \[
    \avect[g]\mbox{ is a }\Big(j[g]^{26} + 5(g-\gamma(g))j[g]^{\phi+49}, j^{27}, \frac{2}{\xi j[g]^3} + (g-\gamma(g))\cdot\frac{3}{j[g]^3(\log j[g])^{1/4}}, 2\Big)\mbox{-RS-good beginning.}
    \]
    \item If $\type^{\Psi[g]} \in \{\advancing, \filling\}$, then $E_g$ starts well.
    \item If $\type^{\Psi[g]} = \stabilising$, then 
    \[
    \avect[g]\mbox{ is a }\Big(j[g]^{26} + 5(g-\gamma(g))j[g]^{\phi+49}, j[g]^{10}, \frac{2}{\xi j^3} + (g-\gamma(g))\cdot\frac{3}{j[g]^3(\log j[g])^{1/4}}, 1\Big)\mbox{-RS-good beginning.}
    \]
    \item $\type^{\Psi[g]} \notin \{\initialising, \Failure\}$.
\end{enumerate}
Observe by (WB\ref{item:WB3}) that $g-\gamma(g) \le 2j[g]$; thus $j[g]^{26} + 5(g-\gamma(g))j[g]^{\phi+49} \le j[g]^{\phi+51}/2$ and $2/(\xi j[g]^3) + (g-\gamma(g))\cdot 3/(j[g]^3(\log j[g])^{1/4}) \le 1/(200j[g]^2)$,  so 
(VEB1), (VEB2), and (VEB3) imply that every escape process $E_g$ starts well.

We now prove (VEB1)--(VEB4) by induction on $g$. For the base case $g=1$, observe since $\Psi[0]$ is an initialising state, $\tau_1 = \tau_0 + 1$ (as there will be a transition via one of (R\ref{item:R1})--(R\ref{item:R3})).  By (WB\ref{item:WB1}), $\Psi[1]$ is not a failure state, so the transition to $\Psi[1]$ is 
not via (R\ref{item:R1}). The set $\calS$ is empty for an initialising state, so it is also not via (R\ref{item:R2}). Thus, it is 
via (R\ref{item:R3}) and $\Psi[1]$ is either an advancing or filling state. Moreover, by the definition of the VEB coupling (Definition~\ref{def:VEB-couple}), $\avect[1]=\bvect_{[j]}^{E_1}(\tau_1) = e_i^X(\tau_1)$ and in particular $\sum_{k=1}^{j[1]} \avect[1]_k \le 1$. Thus $\avect[1]$ is a $j[1]^7$-good beginning. Since $\Psi[1]$ is either an advancing or filling state, it follows that $E_1$ starts well as required.

Now fix $g \ge 1$ and suppose that (VEB1)--(VEB4) hold for $E_g$; we will prove that if $V$ is well-behaved and $E_g$ does not derail, then (VEB1)--(VEB4) hold for $E_{g+1}$ as well. By (VEB4), $\Psi[g]$ is either an advancing, filling, refilling, or stabilising state. Thus:
\begin{itemize}
\item If $\Psi[g]$ is a refilling state, 
we will use Lemma~\ref{lem:veb-refilling-transition} with $d = j[g]^{26} + 5(g-\gamma(g))j[g]^{\phi+49}$ and $\theta = 2/(\xi j[g])^3) + (g-\gamma(g))\cdot 3/(j[g]^3(\log j[g])^{1/4})$ to show that
(VEB1)--(VEB4) hold for $E_{g+1}$.
The conditions of the lemma are satisfied by (VEB1). Since $E_g$ does not derail,
the events of Lemma~\ref{lem:veb-refilling-transition}(i) and (ii) occur. Since $V$ is well-behaved, this implies the events (I), (II) and (III) of Lemma~\ref{lem:veb-refilling-transition} also occur. (I) guarantees that
$\Psi[g+1]$ is a refilling or stabilising state, so $\gamma(g+1) = \gamma(g)$ and $j[g+1] = j[g]$.
We wish to establish (VEB1)--(VEB4) for $g+1$. (VEB2) and (VEB4) follow from (I), and (VEB1) and (VEB3) follow from (II).

\item If $\Psi[g]$ is an advancing state, 
we will use Lemma~\ref{lem:advancing-escape-final}
to show that 
(VEB1)--(VEB4) hold for $E_{g+1}$. The conditions of the lemma are satisfied since (VEB2)   implies that $E_g$ starts well.
Since $E_g$ does not derail and $V$ is well-behaved, items (i), (ii), (iii), and (I), (II), and (III) in the lemma statement occur.
Item (I) shows that (VEB1), (VEB3), and (VEB4) are vacuous. Item (II) gives (VEB2) for~$g+1$.
 
\item If $\Psi[g]$ is a filling state, we will use 
Lemma~\ref{lem:filling-escape-final}. The conditions of the lemma are satisfied since (VEB2) implies that $E_g$ starts well. Since $E_g$ does not derail and $V$ is well behaved, the items in the lemma statement occur. Item (I) shows that (VEB3) and (VEB4) are vacuous.  Item (III) gives (VEB2). We will see that item (II) gives (VEB1) for~$g+1$.
To see this, note that in the case
that $\Psi[g+1]$ is a refilling state, $j[g+1]=j[g]$ (from I) and 
$\gamma(g+1)=\gamma(g)$ so $g+1- \gamma(g+1) \geq 1$.
Now to derive (VEB1) from (II)
we note that $6j[g]^{25} \leq j[g]^{\phi+49}$.

\item If $\Psi[g]$ is a stabilising state, we will use Lemma~\ref{lem:veb-escape-stab} 
with $d = j[g]^{26} + 5(g-\gamma(g))j[g]^{\phi+49}$ and $\theta = 2/(\xi j[g])^3) + (g-\gamma(g))\cdot 3/(j[g]^3(\log j[g])^{1/4})$. The conditions of the lemma follow from (VEB3) for~$g$. Item (I) shows that $j[g+1] = j[g]$ and (VEB3) and (VEB4) are vacuous. (III) gives (VEB2) for~$g+1$. If $\Psi[g+1]$ is a refilling state then $\gamma(g+1) = \gamma(g)$
so $g+1 - \gamma(g+1) > g-\gamma(g)$.
Thus (VEB1) for~$g+1$ comes from (II).

\end{itemize}
We conclude that $V$ is well-behaved and no escape process $E_g$ derails, then every escape process $E_g$ starts well.

\medskip\noindent\textbf{Starting well is sufficient:}
It now remains only to prove that if $V$ is well-behaved, no escape process $E_g$ derails, and every escape process $E_g$ starts well, then (i) and (ii) of the lemma statement are satisfied. Suppose these events occur, and fix $t > \tau_0$. Let $g = \max\{g' \colon \tau_{g'} \le t\}$, so that $\Psi(t) = \Psi[g]$. Then we split into cases depending on $\Psi[g]$. Observe that since $t > \tau_0$ and $\tau_1 = \tau_0+1$, $\Psi[g]$ is not an initialising state. Moreover, $\bvect_{[j]}^{E_{g^{\Psi(t-1)}}}(t) = \bvect_{[j]}^{E_{g^{\Psi(t)}}}(t) = \bvect_{[j]}^{E_g}(t)$; indeed, either $\Psi(t-1) = \Psi(t)$, in which case equality is immediate, or $t = \tau_g$, in which case the equality holds from the construction of $E_g$ in Definition~\ref{def:VEB-couple}.
\begin{itemize}
    \item If $\Psi[g]$ is a failure state, then (i) and (ii) are both vacuous.
    \item If $\Psi[g]$ is a refilling state, then (i) is vacuous. Moreover, (ii) holds by Lemma~\ref{lem:veb-refilling-transition}(III) since $E_g$ starts well, $E_g$ does not derail, and $V$ is well-behaved.
    \item If $\Psi[g]$ is an advancing state, then (i) holds by Lemma~\ref{lem:advancing-escape-final}(III) since $E_g$ starts well, $E_g$ does not derail, and $V$ is well-behaved. Moreover, (ii) is vacuous.
    \item If $\Psi[g]$ is a filling state, then (i) is vacuous. Moreover, (ii) holds by Lemma~\ref{lem:filling-escape-final}(V) since $E_g$ starts well, $E_g$ does not derail, and $V$ is well-behaved.
    \item If $\Psi[g]$ is a stabilising state, then (i) is vacuous. Moreover, (ii) holds by Lemma~\ref{lem:veb-escape-stab}(IV) since $E_g$ starts well, $E_g$ does not derail, and $V$ is well-behaved.
\end{itemize}
\end{proof}

\section{Main result for   finitely-many strongly exposed bins}
\label{sec:final}

Our goal in this section is to prove Theorem~\ref{thm:goaljkillerSE}, which combines the volume and escape analysis. Using it, we then prove Theorem~\ref{thm:goaljkiller}, which finishes our proof since we have already proved Theorem~\ref{thm:goal} using  
Theorem~\ref{thm:goaljkiller}.

\begin{theorem}\label{thm:goaljkillerSE}
There exists a real number $\lambda_0\in (0,1/2)$ such that the following holds. 
Consider any $\lambda \in (0,\lambda_0)$. Consider any send sequence~$\p$ such that $p_0=1$ and, for all $i\geq 1$, $W_i \leq (1/\lambda)^i$. Suppose that $\p$ has finitely many bins which are strongly exposed. Then
the backoff process with birth rate~$2\lambda$ and send sequence~$\p$ is transient. 
\end{theorem}

\begin{proof}
Choose an integer $j_0\geq 5$ large enough to apply Lemma~\ref{lem:VEB-escape} and Corollary~\ref{cor:VEB-volume}, and large enough that all strongly-exposed bins lie in $[j_0-1]$. $j_0$ should also be large enough that  
$\CUpsilon \L(j_0)/\lambda \geq 2$ and
$\CUpsilon \LL(j_0)/\lambda \geq 2$. 

Let $X$ be a backoff process with birth rate $2\lambda$, send sequence~$\p$, and cohort set~$\calC=\{A,B\}$ (Definition~\ref{def:backoff-cohorts}).
Recall from Definition~\ref{def:Tprot} that 
$T^{\p,j}  = 80 j^2 \lfloor \sum_{x\in [j]} W_x \rfloor$.
Let $\Einit$ be the event that, for every step 
$t' \in [T^{\p,j_0}]$, $n^{X^B}(t') \geq 2$.
Let $\epsilon$ be the (small, but positive) probability that~$\Einit$ occurs.

For the next portion of the proof, we condition on~$\Einit$.
Let
$\tau_0 =  T^{\p,j_0}$. 
Let
$\Psi_0$ be the  initialising high-level state   
with $\tau^{\Psi_0} = \tau_0$ and $j^{\Psi_0} = j_0$.  
Note from Definition~\ref{def:backoff-bounding-state-space} that $\OurHLS$ is 
is the $\Psi_0$-closure (Definition~\ref{def:hls-closed}) of the $\Psi_0$-backoff bounding rule~$\hls$ (Definition~\ref{def:backoff-bounding-rule}).  The $\Psi_0$-backoff bounding rule~$\hls$ is a transition rule (Definition~\ref{def:transition-rule}) by Observation~\ref{obs:BB-bounding-is-transrule} and 
it is $\OurHLS$-valid (Definition~\ref{def:valid}) by
Lemma~\ref{lem:backoff-bounding-valid}.

Conditioned on $\Einit$,
the VEB-coupling (Definition~\eqref{def:VEB-couple}   couples $X$  
with a volume process $V$ from~$\Psi_0$ with  transition rule $\hls$, send sequence~$\p$ and birth rate~$\lambda$ and with a sequence $E_1,E_2,\ldots $ of escape processes with send sequence~$\p$ and birth rate~$\lambda$.
The coupling uses $b_i^{E_0}(\tau_0+1)$ as a synonym for $\e_i^X(\tau_0+1)$.
      
By Corollary~\ref{cor:VEB-volume} and Lemma~\ref{lem:VEB-escape} and a union bound,  with probability at least $1/3$, (conditioned on~$\Einit$), items (i), (ii), and (iii) below are true. We will prove below that these items imply 
that for all $t>\tau_0$, $\sum_{i\geq 1} b_i^X(t) > 0$. Combining this with that fact
that $\Einit$ prohibits the backoff process from returning to the empty state before time~$\tau_0$, we obtain that, with at least the positive probability~$\epsilon/3$,
the backoff process does not return to the empty state.

\begin{enumerate}[(i)]
\item for all $t \ge \tau_0$, $\type^{\Psi(t)} \ne \Failure$.

\item for all $t> \tau_0$ such that 
$\type^{\Psi(t)}  = \advancing$
and the $S\in \calS^{\Psi(t)}$,
$\sum_{k \in  S } p_kb_k^{E_{g^{\Psi(t-1)}}}(t) \leq \lambda | S  |/80$. Moreover, $S$ is non-empty.
\item for all $t> \tau_0$ such that 
$\type^{\Psi(t)} \notin \{\Failure,\advancing\}$ and all $k \in [j^{\Psi(t-1)}-1]$ such that $W_k \ge (j^{\Psi(t-1)})^\chi$, $p_kb_k^{E_{g^{\Psi(t-1)}}}(t) \le \lambda/80$. Moreover, there is at least one such $k\in [j-1]$.

\end{enumerate}

Suppose that  (i), (ii), and (iii) occur. 
By (i) and 
the fact that no states after~$\tau_0$ are initialising (which follows from (T\ref{item:trans5})) and
the Coupling invariant (I2),
for all $t>\tau_0$ and all $i\geq 1$, 
\begin{equation}\label{eq:invweuse}
\b_i^{Y_{t-1}}(t) \setminus \b_i^{E_{g^{\Psi(t-1)}}}(t) \subseteq \b_i^X(t).\end{equation}
So for all $t>\tau_0$ and all $i\geq 1$,
$b^X_i(t) \ge b^{Y_{t-1}}_i(t) - b^{E_{g^{\Psi(t-1)}}}_i(t)$.
We will prove that  
$\sum_{i\geq 1} b^X_i(t) >0$. 
We split into cases depending on $\type^{\Psi(t)}$. The cases are exhaustive since $\type^{\Psi(t)}=\Failure$ is ruled out by~(i) and $\type^{\Psi(t)}=\initialising$ is ruled out by~(T5). 
To simplify the notation, let $j= j^{\Psi(t)}$.

\medskip\noindent\textbf{Case 1:} $\type^{\Psi(t)}  =\advancing $.   
By the definition of the volume process (Definition~\ref{def:volume-process}), 
the high-level state $\Psi(t)$ is reached by
$\Psi(t) = \hls(\Psi(t-1), \bbar^{Y_{t-1}}_{[j^{\Psi(t-1)}]}(t),t)$.
We consider two cases.
\begin{itemize}

\item If $\Psi(t-1)\neq \Psi(t)$: In this case, the transition is via (R\ref{item:R3})(i) or (R\ref{item:R4})(i)(a). In either case, (T\ref{item:trans3}) is satisfied  
and guarantees that $\sum_{x\in S} p_x b_x^{Y_{t-1}} \geq \lambda |S|/40$.

\item  If $\Psi(t-1) = \Psi(t)$: In this case, the transition is via (R\ref{item:R4})(ii) and (R\ref{item:R2}) does not apply.
$\calS^{\Psi(t-1)}$ is the same as $\calS^{\Psi(t)}$,
so the guard on (R\ref{item:R2}) implies that
$\sum_{k\in S} p_k b_k^{Y_{t-1}} \geq \lambda |S|/40$.

\end{itemize}
By (ii), it follows that $\sum_{i \in S} b^{Y_{t-1}}_i(t) - b^{E_{g^{\Psi(t-1)}}}_i(t) \geq \lambda |S|/80 > 0$ 
so 
\[\sum_i b^X_i(t) \geq
\sum_{i\in S} b^X_i(t) \geq 
\sum_{i\in S} b^{Y_{t-1}}_i(t) - b^{E_{g^{\Psi(t-1)}}}_i(t)
  > 0
,\] as required.

\medskip\noindent\textbf{Case 2:} $\type^{\Psi(t)} \in \{\filling,\refilling, \stabilising\}$. 
From the backoff-bounding rule (Definition~\ref{def:backoff-bounding-rule}),
$j^{\Psi(t-1)} \in \{j-1,j\}$.

Since $\type^{\Psi(t)} \in \{\filling, \refilling, \stabilising\}$, bin $j$ is exposed. By our choice of $j_0$, it follows that bin $j$ is weakly exposed. By Definition~\ref{def:covered}, this implies 
that (Prop~\ref{cov-prop-3}) holds so
$|\Upsilon_{j,\geq j^{\chi}}| \geq \CUpsilon \LL(j) /\lambda \geq 2$,
so there are at least two bins in 
$[j-1]$ with weight at least $j^{\chi}$ and hence at least one bin~$k$ in $[j-2]$ with weight at least $j^{\chi}$. 
By~\eqref{eq:invweuse} and~(iii),  
\[
\sum_i b_i^X(t) \ge b_k^X(t) \ge b_k^{Y_{t-1}}(t) - b_k^{E_{g^{\Psi(t-1)}}}(t) \ge b_k^{Y_{t-1}}(t) - W_k \lambda/80.
\]
It therefore suffices to prove that $b_k^{Y_{t-1}}(t) \ge W_k \lambda/40$. 

Since $\type^{\Psi(t)} \in \{\filling,\refilling, \stabilising\}$ the transition to $\Psi(t)$ from the backoff-bounding rule (Definition~\ref{def:backoff-bounding-rule})
occurs when (R\ref{item:R1}) 
does not apply.
Since (R\ref{item:R1}) didn't apply item (R\ref{item:R1})(ii) guarantees $F_k^{\Psi(t-1)}(t) \geq \lambda W_k/2$ and item (R\ref{item:R1})(iii) 
guarantees  $p_k b_k^{Y_{t-1}}(t)\geq \lambda/4 > \lambda/40$ as required.
\end{proof}

\begin{remark}
\label{rem:meta} 
The proof of Theorem~\ref{thm:goaljkillerSE} establishes that, with positive probability, 
no messages escape before time~$\tau_0$
and,
for all $t>\tau_0$,
$\sum_{i\geq 1} b_i^X(t) > 0$.
In fact, it proves something stronger.
For any $t>\tau_0$ 
there are two cases, depending on whether $j^{\Psi(t)}$ is covered, or weakly exposed.
If it is covered (Case 1 of the proof), the population is
actually  $\Omega(|S|) = \Omega( \log(j^{\Psi(t)}))$. Otherwise (Case 2 of the proof), all of the very-high weight bins in $[j-1]$ are non-empty so the population is   $\Omega(\log \log j^{\Psi(t)})$. However, Corollary~\ref{cor:VEB-volume}
establishes that for all $j \ge j_0$, there is
a positive integer $t \le \tau_0 + 1+ {(8j)}^{\phi+1}(W_{j_0}+\dots + W_{j-1})$ such that 
$\Psi(t)$ is $j$-filling or $j$-advancing. Also, for all $i\geq 1$,
$W_i \leq (1/\lambda)^i$ so $j^{\Psi(t)} = \Omega(\log t)$.
Thus, the population at time~$t$
is $\Omega(\log\log\log t)$.

\end{remark}
 
We can now prove Theorem~\ref{thm:goaljkiller}.

\thmgoaljkiller*
\begin{proof}

Let $\lambda'_0$ be the $\lambda_0$
of Corollary~\ref{cor:backoff-finite-SE-bins} and let $\lambda''_0$ be the $\lambda_0$ of 
Theorem~\ref{thm:goaljkillerSE}.
Let 
$\lambda_0 = \min\{\lambda'_0,2\lambda''_0\}$.
Consider any $\lambda \in (0, \lambda_0)$. Let $\lambda' = \lambda$ and let $\lambda'' = \lambda/2$. 
Consider any send sequence~$\p$ such that $p_0=1$ and, for all $j\geq 1$, $W_j \leq (2/\lambda)^j$.
If $\p$ has infinitely many bins which are strongly exposed then 
note that $\lambda' \in (0, \lambda'_0)$ and for all $j\geq 1$, 
$W_j \leq (2/\lambda')^j$.
By Corollary~\ref{cor:backoff-finite-SE-bins} the backoff process with birth rate~$\lambda = \lambda'$ and send sequence~$\p$ is not positive recurrent. If $\p$ has finitely many bins which are strongly exposed then note that 
$\lambda'' \in (0, \lambda''_0)$ and for all $j\geq 1$, 
$W_j \leq (1/\lambda'')^j$. By Theorem~\ref{thm:goaljkillerSE},
the backoff process with birth rate $\lambda = 2 \lambda''$ and send sequence~$\p$ is not positive recurrent (in fact, it is transient). 
\end{proof}

\section{Discussion}
\label{sec:discussion}

The result of this paper is 
Theorem~\ref{thm:goal}
which proves Aldous's conjecture 
by showing that every backoff protocol is unstable for every positive birth rate~$\lambda$, in the sense that the underlying Markov process is not positive recurrent.

We briefly discuss three interesting directions for future work.

First, is   every backoff protocol  unstable in the stronger sense that the corresponding Markov process is \emph{transient} for every positive birth rate? Aldous's proof of instability for binary exponential backoff does establish transience. Also, the charcterisation of Kelly and MacPhee (with probability~$1$, only finitely many balls escape) implies transience.
The most difficult part of this paper was proving Theorem~\ref{thm:goaljkillerSE}
which shows instability for send sequences~$\p$ with finitely-many strongly exposed bins. The instability result of Theorem~\ref{thm:goaljkillerSE} is actually transience. So that is a good start.
However, the corresponding result for the case with infinitely many strongly exposed bins (Corollary~\ref{cor:backoff-finite-SE-bins}) is truly a proof of lack of positive recurrence, so that result would have to be strengthened. Moreover, the statement of Theorem~\ref{thm:goaljkillerSE}
relies on several simplifications 
which are enabled via reductions in Section~\ref{sec:reduction}.
The simplification that $\lambda$ is sufficiently small (Observation~\ref{obs:decrease-lambda}) would work just as well with transience. The same holds for the simplification that $p_0=1$. However, the upper bound on bin weights relies on Lemma~\ref{lem:killer-fact}
which is~\cite[Lemma 7]{GL-oldcontention}. This is again a proof of lack of positive recurrence, so it would have to be strengthened.
  
Kelly~\cite{Kelly} asked whether there is any acknowledgement-based contention-resolution protocol that is stable for any positive birth rate. 
Acknowledgement-based protocols are a larger class of protocols than backoff protocols, so our result does not answer this question.
Kelly proved that stability of acknowledgement-based protocols is impossible for $\lambda > 0.567$. This was later improved by~\cite{GJKP}, who showed that
there is no stable acknowledgement-based protocol for any $\lambda > 0.531$. We agree with MacPhee~\cite{MacPhee} that there is probably no stable acknowledgement-based protocol for any positive~$\lambda$, and this is a very interesting open question. A particularly natural family of acknowledgement-based protocols that differs from backoff is the family of \emph{age-based protocols}.  An age-based protocol is associated with a send sequence~$\p$. A message sends with probability~$p_k$ when it has been in the system for $k$~steps. The results of Kelly and MacPhee~\cite{kelly-macphee} apply to age-based protocols as well as to backoff protocols. Ingenoso~\cite{Ingenoso} showed instability for the case that $\p$ is monotonically decreasing. It is possible that the method of this paper could prove instability for every age-based protocol for every positive birth rate (though we have not considered this).

We have shown that every backoff protocol is unstable for every positive birth rate $\lambda>0$. 
Now that we know this, it would be very interesting to investigate \emph{meta-stability}. How long does it take before the channel is expected to become congested? What is the behaviour of the channel before this point? What can we say precisely about the behaviour after that point? (For example, see Remark~\ref{rem:meta}.)

\section{Acknowledgements}

We thank Frank Kelly for useful conversations about this topic many years ago. His work inspired a lot of our work on this problem. Thanks also to Frank for getting Iain MacPhee's thesis scanned for us. We also thank Mike Paterson, Mark Jerrum, and Sampath Kannan for very many useful (and enjoyable) conservations about this problem. These conversations made a big contribution - beyond what is captured by~\cite{GJKP}.

\newpage
\section{List of Symbols}

\printglossary[title={}]
 
\newpage
 
\bibliographystyle{plain}
\bibliography{references}

\end{document}